\documentclass[12pt]{amsart}
\usepackage[margin=1in]{geometry}
\usepackage{amsmath,amssymb,amsthm,mathtools,bm,booktabs,enumitem,mathrsfs,array,longtable,tabularx,adjustbox}
\usepackage{iftex}
\ifXeTeX
  \usepackage[expansion=false]{microtype}
\else
  \usepackage[expansion=basictext]{microtype}
\fi
\usepackage{needspace,etoolbox}
\usepackage{tikz}
\usepackage[hidelinks]{hyperref}
\hypersetup{
 hypertexnames=false,
 pdftitle={Collapse of white dwarfs above the Chandrasekhar mass},
 pdfauthor={Gong Chen and Zhiwu Lin},
 pdfsubject={Finite-time collapse with exact Chandrasekhar pressure and
 physical vacuum},
 pdfkeywords={Euler--Poisson system, white dwarf, Chandrasekhar mass,
 physical vacuum, gravitational collapse}}
\allowdisplaybreaks
\numberwithin{equation}{section}
\theoremstyle{plain}
\newtheorem{theorem}{Theorem}[section]
\newtheorem{proposition}[theorem]{Proposition}
\newtheorem{lemma}[theorem]{Lemma}
\newtheorem{corollary}[theorem]{Corollary}

\theoremstyle{definition}
\newtheorem{definition}[theorem]{Definition}

\theoremstyle{remark}
\newtheorem{remark}[theorem]{Remark}

\BeforeBeginEnvironment{theorem}{\Needspace{6\baselineskip}}
\BeforeBeginEnvironment{proposition}{\Needspace{6\baselineskip}}
\BeforeBeginEnvironment{lemma}{\Needspace{6\baselineskip}}
\BeforeBeginEnvironment{corollary}{\Needspace{6\baselineskip}}
\BeforeBeginEnvironment{definition}{\Needspace{6\baselineskip}}
\BeforeBeginEnvironment{remark}{\Needspace{6\baselineskip}}

\newcommand{\dd}{\,\mathrm{d}}
\newcommand{\pa}{\partial}

\newcommand{\abs}[1]{\left\lvert #1\right\rvert}

\newcommand{\la}{\lambda}

\makeatletter
\newcommand{\equationalias}[2]{%
  \@ifundefined{r@#2}{}{%
    \protected@write\@auxout{}{%
      \string\newlabel{#1}{\csname r@#2\endcsname}}}}
\makeatother

\newcommand{\introheading}[1]{%
  \par\smallskip\noindent\emph{#1}\enspace}

\title[Collapse above the Chandrasekhar mass]{Collapse of white dwarfs\\above the Chandrasekhar mass}
\author[G. Chen]{Gong Chen}
\address{School of Mathematics, Georgia Institute of Technology,
Atlanta, Georgia 30332, USA}
\email{gc@math.gatech.edu}
\author[Z. Lin]{Zhiwu Lin}
\address{School of Mathematical Sciences, Fudan University,
Shanghai 200433, China}
\email{zwlin@fudan.edu.cn}
\date{}

\begin{document}

\begin{abstract}
We give a dynamical characterization of the Chandrasekhar mass for
cold white dwarfs governed by the radial Euler--Poisson equations with
the exact Chandrasekhar equation of state. This mass is the infimum
of masses for which a finite-energy classical solution with a physical
vacuum can contract its entire support to a point in finite time.
Below it, conservation of energy prevents complete collapse. For every
mass in a sufficiently small interval above it, we construct a
solution whose support shrinks to the origin and whose density
concentrates there with the full mass of the star. The interior
approaches a collapsing Goldreich--Weber profile. Near the vacuum
boundary, a different rescaling gives a limiting profile that connects
the high- and low-density regimes of the pressure law. We determine
this profile and show that, after adding the electron rest energy,
the conserved energy equals that of the Goldreich--Weber reference
solution. The proof combines matched asymptotic expansions with local
existence and energy estimates that extend solutions backward from
times approaching collapse to a common earlier time.
\end{abstract}

\maketitle
\tableofcontents

\section{Introduction}
\label{sec:introduction}

\subsection{The model and the collapse result}

A cold white dwarf is supported against gravity by electron degeneracy
pressure. In Chandrasekhar's model, a spherical star at rest has a
mass below a finite limiting value: as the density at its center
increases, its radius decreases and its mass approaches the
\emph{Chandrasekhar mass} $M_{\rm Ch}$
\cite{Chandrasekhar1931,Chandrasekhar1935}. This equilibrium theory
identifies the mass that pressure can support. The time-dependent
problem asks whether a star with mass just above $M_{\rm Ch}$ can
collapse completely, with all its matter reaching the origin at a
single finite time.

We construct such solutions for every mass in a sufficiently small
interval above $M_{\rm Ch}$. Below $M_{\rm Ch}$, the sharp inequality
between internal and gravitational energy prevents the whole support
from shrinking to a point. These two conclusions identify
$M_{\rm Ch}$ as the infimum of masses admitting complete collapse
in the radial classical solution class considered here. The main
task is the existence construction for the full Chandrasekhar
pressure law. Its high-density approximation describes the interior
of the collapsing star, but cannot describe its surface, where the
density must fall to zero.

We first specify the model. For a noninteracting electron gas at zero
temperature and fixed chemical composition, a choice of units gives
the pressure law
\begin{equation}
 P(0)=0,\qquad P'(\rho)=\frac{4\rho^{2/3}}{3\sqrt{1+\rho^{2/3}}}.
 \label{in:exact-pressure}
\end{equation}
Its low- and high-density limits are
\begin{equation}
 \begin{aligned}
 P(\rho)&=\frac45\rho^{5/3}(1+o(1))&&\quad(\rho\downarrow0),\\
 P(\rho)&=\rho^{4/3}(1+o(1))&&\quad(\rho\to\infty).
 \end{aligned}
 \label{in:two-regimes}
\end{equation}
The change of exponent comes from the transition from nonrelativistic
to relativistic electron momenta under compression. The fluid motion
and gravity are Newtonian. Throughout the paper we use the full law
\eqref{in:exact-pressure}; its dimensional normalization is given in
Subsection~\ref{sec:normalization}.

Let $\rho(t,X)$, $v(t,X)$, and $\Phi(t,X)$ denote density, velocity,
and gravitational potential, respectively. The Euler--Poisson equations
are
\begin{equation}
 \begin{aligned}
 \partial_t\rho+\nabla_X\!\cdot(\rho v)&=0,\\
 \rho(\partial_tv+v\!\cdot\nabla_Xv)+\nabla_XP(\rho)
      &=-\rho\nabla_X\Phi,\\
 \Delta_X\Phi&=4\pi\rho,\qquad
 \Phi(t,X)\to0\quad(|X|\to\infty).
 \end{aligned}
 \label{in:eulerian-system}
\end{equation}
We consider spherical stars in $\mathbb R^3$. Their density is positive
in the ball $B_{R_b(t)}$ and zero outside, and the surface moves with
the fluid. At each time before collapse, we impose the
\emph{physical-vacuum condition}: the squared sound speed $P'(\rho)$
vanishes linearly with the inward distance to the surface. For
\eqref{in:exact-pressure}, this means that the density vanishes like
$(R_b(t)-|X|)^{3/2}$. The total mass is
\[
 M=\int_{\mathbb R^3}\rho(t,X)\dd X,
\]
which is independent of time.

The first theorem states the mass conclusion without the asymptotic
details. We use the finite-energy radial classical solutions of
Definition~\ref{def:global-solution}, with regularity index $N=22$.

\begin{theorem}
\label{thm:principal-collapse}
For the Euler--Poisson system \eqref{in:eulerian-system} with pressure
\eqref{in:exact-pressure}, the following hold in this solution class.
\begin{enumerate}[label=(\roman*),leftmargin=2.4em]
\item If $0<M<M_{\rm Ch}$, conservation of finite energy bounds the
support radius away from zero throughout the interval of existence.
\item There is $\varepsilon>0$ such that, for every
$M\in(M_{\rm Ch},M_{\rm Ch}+\varepsilon)$, there is a solution of
mass $M$ on $[-T,0)$, for some $T>0$, satisfying
\[
 R_b(t)\longrightarrow0,\qquad
 \rho(t,X)\dd X\stackrel{*}{\rightharpoonup}M\delta_0
 \quad\text{as }t\uparrow0.
\]
\end{enumerate}
Consequently, the infimum of masses admitting complete finite-time
collapse in this class is $M_{\rm Ch}$.
\end{theorem}

Here $\delta_0$ is the unit point mass at the origin. Part~(i) is the
energy obstruction established in \cite{ChengChengLin}; we recall its
short proof below. Part~(ii) is the construction proved in this paper.
The initial density and inward velocity are chosen together in the
proof.

For the homogeneous law $P=\rho^{4/3}$, the Goldreich--Weber solutions
contract by a common scale factor \cite{GoldreichWeber}. To construct
collapse for \eqref{in:exact-pressure}, we must connect this
high-density motion to the different low-density behavior at the
surface. We do so by matching an interior expansion to a boundary
layer, and then solving the evolution problem with an error small
enough to preserve both descriptions. The local construction and the
estimates for backward continuation must use the same boundary
conditions and the same fixed mass. Subsection~\ref{in:proof-strategy}
explains how these requirements determine the proof.

\subsection{The lower bound on the mass}
\label{in:energy-bound}

The mass bound follows from the competition between pressure and
gravity under compression. Define the enthalpy and internal energy
density by
\[
 h(\rho)=\int_0^\rho\frac{P'(s)}s\dd s
       =4\bigl(\sqrt{1+\rho^{2/3}}-1\bigr),\qquad
 F(\rho)=\int_0^\rho h(s)\dd s.
\]
For a nonnegative density $\rho\in L^1(\mathbb R^3)\cap
L^{4/3}(\mathbb R^3)$ of mass $M=\int_{\mathbb R^3}\rho\dd X$, write
\[
 \mathcal W(\rho)=\frac12\iint_{\mathbb R^3\times\mathbb R^3}
             \frac{\rho(X)\rho(Y)}{|X-Y|}\dd X\dd Y.
\]
Thus $\mathcal W\ge0$ is the magnitude of the gravitational potential
energy; attraction contributes $-\mathcal W$ to the conserved total
energy
\[
 \mathscr E=\frac12\int_{\mathbb R^3}\rho|v|^2\dd X
                 +\int_{\mathbb R^3}F(\rho)\dd X-\mathcal W(\rho).
\]

Compress a density by a factor $a>0$, keeping its mass fixed:
$\rho_a(X)=a^{-3}\rho(X/a)$. The regime $a\downarrow0$ describes
concentration into a smaller ball. Since $F(\rho)=3\rho^{4/3}(1+o(1))$
at high density, the two leading potential energies scale in the
same way:
\[
 \int\rho_a^{4/3}\dd X=a^{-1}\int\rho^{4/3}\dd X,
 \qquad \mathcal W(\rho_a)=a^{-1}\mathcal W(\rho).
\]
Their sharp comparison is
\begin{equation}
 \mathcal W(\rho)\le
 3\left(\frac{M}{M_{\rm Ch}}\right)^{2/3}
                 \int_{\mathbb R^3}\rho^{4/3}\dd X.
 \label{in:threshold-inequality}
\end{equation}
The best constant is $3M_{\rm Ch}^{-2/3}$. Equality is attained by
the Lane--Emden equilibrium: the spherical, compactly supported
stationary solution for the homogeneous pressure $P=\rho^{4/3}$.
Its mass is $M_{\rm Ch}$, the limiting mass of the equilibria for the
exact pressure. Its density is described in
Subsection~\ref{in:goldreich-weber-main}; see also
\cite[Theorem~3.1 and Lemma~5.1]{ChengChengLin} and
Subsection~\ref{mat:seed} for the normalization.

The explicit enthalpy gives $F(\rho)\ge3\rho^{4/3}-4\rho$. If
$0<M<M_{\rm Ch}$, set
\[
 c_M=3\left[1-\left(\frac{M}{M_{\rm Ch}}\right)^{2/3}\right]>0.
\]
H\"older's inequality on the support gives
\[
 M=\int_{B_{R_b(t)}}\rho\dd X
 \le\left(\int_{B_{R_b(t)}}\rho^{4/3}\dd X\right)^{3/4}
       |B_{R_b(t)}|^{1/4}.
\]
Since $|B_{R_b(t)}|=4\pi R_b(t)^3/3$, \eqref{in:threshold-inequality}
and the nonnegativity of the kinetic energy imply
\[
 \mathscr E+4M\ge c_M\int\rho^{4/3}\dd X
 \ge c_M\left(\frac3{4\pi}\right)^{1/3}
                      \frac{M^{4/3}}{R_b(t)}.
\]
Mass and energy are constant along the classical flow. In particular,
$\mathscr E+4M>0$ and
\[
 R_b(t)\ge
 \frac{c_M(3/(4\pi))^{1/3}M^{4/3}}{\mathscr E+4M}>0.
\]
This proves Theorem~\ref{thm:principal-collapse}(i). It excludes
contraction of the whole support; it makes no assertion about other
possible singularities. Above $M_{\rm Ch}$ the coefficient $c_M$
changes sign, so this argument no longer prevents complete collapse.
Proving that collapse actually occurs then requires a solution of
the evolution equations.

\subsection{The collapsing interior and the vacuum surface}
\label{in:goldreich-weber-main}
\label{in:boundary-layer}

We describe the two parts of the solution before stating their precise
asymptotics. For radial fields, write
\[
 \rho(t,X)=\rho_E(t,|X|),\qquad
 v(t,X)=u_E(t,|X|)\frac{X}{|X|}\quad(X\ne0),
 \qquad u_E(t,0)=0.
\]
The subscript $E$ denotes Eulerian variables. Let $r(t,x)$ be the
radius enclosing the fixed mass $x$:
\[
 x=4\pi\int_0^{r(t,x)}\rho_E(t,s)s^2\dd s,
 \qquad 0\le x\le M.
\]
Thus keeping $x$ fixed follows the same matter as it moves. In these
mass coordinates, set
\[
 \rho_L(t,x)=\rho_E(t,r(t,x)),\qquad
 u_L(t,x)=u_E(t,r(t,x))=r_t(t,x).
\]
Differentiation of the mass identity gives
$1=4\pi\rho_Lr^2r_x$. The radial momentum equation becomes
\begin{equation}
 r_{tt}+4\pi r^2\partial_xP(\rho_L)+\frac{x}{r^2}=0,
 \qquad \rho_L=(4\pi r^2r_x)^{-1}.
 \label{in:mass-equation}
\end{equation}
The moving boundary is now the fixed endpoint $x=M$, with
$R_b(t)=r(t,M)$; see \cite[(4)]{Makino2015} and
\cite[Section~2.1]{HadzicJangExpansion}.

\introheading{The Goldreich--Weber profile.}
For $P=\rho^{4/3}$, a common scale factor $\lambda(t)$ separates
the time evolution from the density profile. The Goldreich--Weber
solutions have the form
\[
 \begin{aligned}
 \rho_{\rm hom}(t,\lambda(t)z)&=\lambda(t)^{-3}w_\delta(z)^3,
 &u_{\rm hom}(t,\lambda(t)z)&=\dot\lambda(t)z,\\
 r_{\rm hom}(t,m(z))&=\lambda(t)z,
 &m(z)&=4\pi\int_0^z w_\delta(s)^3s^2\dd s.
 \end{aligned}
\]
Here $z$ is the rescaled radius, and $w_\delta^3$ is the density at
unit scale. Normalize the first zero of $w_\delta$ to be $1$;
the mass of this reference profile is $M=m(1)$.
Substitution into \eqref{in:mass-equation} with the homogeneous
pressure gives
\[
 \begin{gathered}
 w_\delta''+\frac2z w_\delta'+\pi w_\delta^3=-\frac34\delta,
 \qquad w_\delta'(0)=0,\quad w_\delta(1)=0,\\
 \ddot\lambda=\frac{\delta}{\lambda^2}.
 \end{gathered}
\]
At $\delta=0$, taking $\lambda\equiv1$ gives the stationary density
$w_0^3$ and zero velocity. This is the Lane--Emden equilibrium used in
\eqref{in:threshold-inequality}, with $m(1)=M_{\rm Ch}$.
We use nearby profiles with
$\delta=-\beta^2/2<0$. Multiplication of the scale equation by
$2\dot\lambda$ gives
\[
 \frac{\dd}{\dd t}\left(\dot\lambda^2+\frac{2\delta}{\lambda}\right)=0,
 \qquad \dot\lambda^2=\frac{\beta^2}{\lambda}+e.
\]
The constant $e$ specifies the reference velocity once the profile
and radius are fixed. We take $e\ge0$, choose the inward branch
$\dot\lambda<0$, and set the collapse time to zero. Then
$\lambda(t)=(3\beta(-t)/2)^{2/3}(1+o(1))$ as $t\uparrow0$. The homogeneous solutions and
their conserved quantities are discussed in
\cite{GoldreichWeber,Makino1992,FuLin1998,HadzicJangExpansion}.

\introheading{The boundary layer.}
At fixed $z<1$, the homogeneous density tends to infinity as
$\lambda\downarrow0$, and the exact pressure approaches $\rho^{4/3}$.
This approximation is not uniform up to $z=1$. Indeed,
\begin{equation}
 P'(\lambda^{-3}w_\delta^3)=\frac{4w_\delta}{3\lambda}
       \left(1+\frac{\lambda^2}{w_\delta^2}\right)^{-1/2}.
 \label{in:exact-coefficient}
\end{equation}
The second factor approaches one in the interior but differs by
order one when $w_\delta$ is comparable to $\lambda$. Moreover,
the homogeneous density vanishes cubically at its surface, whereas
the physical vacuum for the exact law requires the power $3/2$.
Shrinking the star does not remove this difference: the density
still passes through the low-density regime before reaching zero.

For nonnegative quantities, we write $f\lesssim g$ if $f\le Cg$
and $f\sim g$ if $c g\le f\le Cg$ for positive constants on the
range under consideration. Parameter dependence is specified with
the estimates.

This transition determines the boundary-layer scale. Put $q=1-z$
and $\varkappa_\delta=-w_\delta'(1)>0$. The simple zero of the
profile gives, as $q\downarrow0$,
\[
 \begin{gathered}
 w_\delta(1-q)=\varkappa_\delta q(1+o(1)),\\
 \rho_{\rm hom}\sim(q/\lambda)^3,\qquad
 M-m(1-q)\sim q^4.
 \end{gathered}
\]
The density is of order one at $q\sim\lambda$. The corresponding
region has mass of order $\lambda^4$ and physical width
$\lambda q\sim\lambda^2$. We therefore introduce
\begin{equation}
 D=\frac{M-x}{\lambda^4},\qquad
 r=\lambda-\lambda^2 Z.
 \label{in:layer-scale}
\end{equation}
The variable $D$ measures the mass between a point and the surface,
in units of $\lambda^4$; $Z$ measures its inward displacement from
the reference radius $\lambda$, in units of $\lambda^2$.

The leading layer profiles $\rho_0(D)$ and $Z_0(D)$ satisfy
\[
 \rho_0=(4\pi Z_0')^{-1},\qquad
 M+\delta-4\pi\partial_D P(\rho_0)=0.
\]
The term $\delta$ comes from the acceleration of the contracting
scale. Integrating the Goldreich--Weber equation over $0<z<1$ gives
$M+\delta=4\varkappa_\delta$. The vacuum condition and matching
with the interior then yield
\begin{equation}
 P(\rho_0(D))=\frac{\varkappa_\delta}{\pi}D,\qquad
 Z_0(D)=\frac{h(\rho_0(D))}{4\varkappa_\delta}
                         +\frac1{\varkappa_\delta}.
 \label{in:leading-layer}
\end{equation}
Thus one profile contains both density regimes:
$\rho_0(D)\sim D^{3/5}$ near $D=0$ and
$\rho_0(D)\sim D^{3/4}$ as $D\to\infty$.
The additive constant in $Z_0$ is fixed by matching, as explained in
Subsection~\ref{in:proof-strategy}.

To state regularity at the vacuum endpoint, use $y=D^{1/5}$. This
choice follows from the physical-vacuum powers, at each fixed time:
\begin{equation}
 \rho_L\sim(R_b-r)^{3/2},\qquad
 M-x\sim(R_b-r)^{5/2},\qquad \xi=(M-x)^{1/5}.
 \label{in:vacuum-coordinate}
\end{equation}
In the coordinate $\xi$, the density vanishes to order three and
the distance to the surface to order two. The layer coordinate
$y=\xi/\lambda^{4/5}$ resolves these same powers on the shrinking
region. In particular, $Z_0(y^5)-Z_0(0)$ vanishes to order two,
and $\rho_0(y^5)$ to order three.
\subsection{Asymptotics of the constructed solutions}

We now state the precise form of the collapse construction. The parameter
$e\ge0$, the interior profile $w_\delta$, and the scale $\lambda$
are as above. Define the second moment of the reference density by
\[
 I_2(\delta)=\int_0^{m(1)}z(x)^2\dd x
            =4\pi\int_0^1w_\delta(z)^3z^4\dd z,
\]
where $z(x)$ is the inverse of $m(z)$. This quantity enters the
energy of the homologous motion and that of the solution below.

\begin{theorem}
\label{thm:physical-main}
Fix $e\ge0$. There exists $\varepsilon_M(e)>0$ such that, for every
prescribed mass
\[
 M_{\rm Ch}<M<M_{\rm Ch}+\varepsilon_M(e),
\]
equation \eqref{in:mass-equation} with pressure law
\eqref{in:exact-pressure} has a spherically symmetric classical
physical-vacuum solution on $[-T,0)$ for some $T>0$. Its density is
positive in the interior, its support is compact, and its mass is $M$.
The solution has the following properties.
\begin{enumerate}[label=(\roman*),leftmargin=2.4em]
\item There are $\beta>0$, $\delta=-\beta^2/2$, and a
Goldreich--Weber profile $w_\delta$, normalized as above, such that
\[
 M=4\pi\int_0^1w_\delta(z)^3z^2\dd z.
\]
The reference scale satisfies
\begin{equation}
 \dot\lambda=-\sqrt{e+\beta^2/\lambda},
 \qquad
 \lambda(t)=\left(\frac{3\beta}{2}(-t)\right)^{2/3}
      \bigl(1+O((-t)^{2/3})\bigr).
 \label{in:collapse-rate}
\end{equation}
For every $0<z_0<1$, as $t\uparrow0$,
\[
 \lambda(t)^3\rho_E(t,r)
 =w_\delta\!\left(\frac r{\lambda(t)}\right)^3(1+o(1))
 \quad\text{uniformly for }0\le r\le z_0\lambda(t),
\]
and
\[
 \sup_{0<r\le z_0\lambda(t)}
 \left|\frac{\lambda(t)u_E(t,r)}{\dot\lambda(t)r}-1\right|
 \longrightarrow0.
\]

\item The support radius tends to zero and the entire mass concentrates
at the origin:
\[
 R_b(t)\longrightarrow0,\qquad
 \rho_E(t,|X|)\dd X\stackrel{*}{\rightharpoonup}M\delta_0
 \qquad(t\uparrow0).
\]
The conserved energy satisfies
\begin{equation}
                       \mathscr E+4M=\frac e2I_2(\delta).
 \label{in:energy-identity}
\end{equation}

\item Fix $Y>0$. For $t$ sufficiently close to zero, define the rescaled
boundary radius by
\[
 \mathcal Z_\lambda(y)
 =\frac{\lambda(t)-r(t,M-\lambda(t)^4y^5)}{\lambda(t)^2},
 \qquad 0\le y\le Y.
\]
As $t\uparrow0$, $\mathcal Z_\lambda$ converges in $C^3([0,Y])$.
The three factors
\[
 \frac{\mathcal Z_\lambda'(y)}y,\qquad
 \frac{\mathcal Z_\lambda(y)-\mathcal Z_\lambda(0)}{y^2},\qquad
 \frac{\rho_L(t,M-\lambda(t)^4y^5)}{y^3}
\]
converge in $C^1([0,Y])$, with their removable endpoint values, to
positive limits. This region contains mass $\lambda^4Y^5$ and has
physical thickness of order $\lambda^2$. The limiting profiles are
specified in Subsection~\ref{in:boundary-layer}.
\end{enumerate}
The solution has the regularity of Definition~\ref{def:global-solution}
with $N=22$. The upper bound for $\lambda$ and the error constants may depend on
the selected mass and on $e$.
\end{theorem}
For each fixed $x\in(0,M]$, the radius $r(t,x)$ remains positive
before collapse and is bounded by $R_b(t)\to0$. Thus the matter
enclosing any fixed mass reaches the origin at the same time. For a
continuous test function $\varphi$, conservation of mass gives
\[
 \left|\int\varphi(X)\rho(t,X)\dd X-M\varphi(0)\right|
 \le M\sup_{|X|\le R_b(t)}|\varphi(X)-\varphi(0)|\longrightarrow0.
\]
This proves the asserted concentration once contraction of the support
is known. The energy identity requires more information, as we explain
below.

The range of masses comes from the nearby Goldreich--Weber profiles:
\begin{equation}
 M_\beta=M_{\rm Ch}
       +\frac{I_2(0)}{2M_{\rm Ch}}\beta^2+O(\beta^4),\qquad
 I_2(0)>0.
 \label{in:mass-expansion}
\end{equation}
The mass depends continuously on $\beta$, and the construction works
for all sufficiently small $\beta>0$ after choosing a sufficiently
small upper bound for $\lambda$. It therefore realizes every mass in an interval
immediately above $M_{\rm Ch}$. The detailed choices of parameters and
weighted norms are given in Theorem~\ref{thm:main-collapse}.

Part~(iii) gives asymptotics for the actual density and radius near
the moving vacuum boundary $x=M$. For each fixed $Y>0$, the change
of variables $x=M-\lambda^4y^5$, $0\le y\le Y$, describes a region
containing mass $\lambda^4Y^5$. In these coordinates,
\[
 \mathcal Z_\lambda(y)\longrightarrow Z_0(y^5)
                        \quad\text{in }C^3,
\]
and the quotient limits determine the leading coefficients in the
vacuum asymptotics
\[
 \begin{aligned}
 \frac{R_b(t)-r(t,M-\lambda^4y^5)}{\lambda^2y^2}
   &\longrightarrow \frac{Z_0(y^5)-Z_0(0)}{y^2},\\
 \frac{\rho_L(t,M-\lambda^4y^5)}{y^3}
   &\longrightarrow \frac{\rho_0(y^5)}{y^3}.
 \end{aligned}
\]
Both limits hold in $C^1([0,Y])$, with positive continuous values
at $y=0$. In particular,
\begin{equation}
 R_b(t)=\lambda(t)-\frac{\lambda(t)^2}{\varkappa_\delta}
                                     +o(\lambda(t)^2).
 \label{in:boundary-radius}
\end{equation}
With $h_E(t,r)=h(\rho_E(t,r))$, the physical-vacuum slope satisfies
\[
 \lambda(t)^2[-\partial_rh_E(t,R_b(t))]
                       \longrightarrow4\varkappa_\delta.
\]
Thus the theorem determines both the first correction to the support
radius and the leading density near vacuum. The measure limit
$M\delta_0$ alone gives neither conclusion, since the mass in this
region tends to zero. The convergence above is on each fixed
$y$-interval; see Corollary~\ref{cor:actual-layer}.

\introheading{The conserved energy.}
The moment $I_2(\delta)$ measures the mass distribution of the unit-scale
profile. For the homogeneous motion, the profile equation gives
the sum of internal and gravitational energies as
$\delta I_2(\delta)/\lambda$. Its total energy is therefore
\[
 \mathscr E_{\rm GW}
 =\frac12 I_2(\delta)\dot\lambda^2
                         +\frac{\delta I_2(\delta)}{\lambda}
 =\frac e2 I_2(\delta),
\]
since $\dot\lambda^2=\beta^2/\lambda+e$ and
$\delta=-\beta^2/2$; see \eqref{exlim:homologous-energy}.
The kinetic, internal, and gravitational energies individually have
size $\lambda^{-1}$. Their singular terms cancel. Interior convergence
of the density and velocity alone does not identify the finite
remainder of this cancellation: it neither controls the region near
the surface nor ensures an absolute error tending to zero in the
three energy terms.

The exact pressure contributes a computable change to internal
energy. At every density,
\[
 0\le F(\rho)+4\rho-3\rho^{4/3}\le3\rho^{2/3}.
\]
For a density of mass $M$ supported in $B_{R_b(t)}$, integration gives
\[
 0\le\int\bigl(F(\rho)+4\rho-3\rho^{4/3}\bigr)\dd X
 \le3M^{2/3}|B_{R_b(t)}|^{1/3}=O(\lambda).
\]
The term $-4\rho$ thus contributes exactly $-4M$, and the remaining
change of internal energy tends to zero. Adding $4\rho$ to $F$ changes
the zero of energy, but not the pressure, since
\[
 \rho(F+4\rho)'-(F+4\rho)=\rho F'-F=P(\rho).
\]
In dimensional variables this addition restores the electron rest
energy; see Subsection~\ref{sec:normalization}.

To complete the comparison, let $\mathscr E_A$ be the energy of the
matched approximate solution constructed in Section~\ref{mat:chapter}.
We estimate all three energy terms over the full mass interval,
including the region where the two expansions are joined. This gives
\[
 \left|\mathscr E_A+4M-\mathscr E_{\rm GW}\right|
       \le C\lambda(1+|\log\lambda|)^{C_{\log}},\qquad
 \mathscr E-\mathscr E_A\longrightarrow0.
\]
Here $C_{\log}$ is finite and the constants may depend on the chosen
profile and $e$. The two estimates are proved in
Lemma~\ref{exlim:global-profile-energy} and
Proposition~\ref{exlim:clock-energy}, respectively. Conservation of
energy then yields \eqref{in:energy-identity}. Thus no additional
finite contribution to the total energy remains from the surface
region or the error between the approximate and exact solutions.

For a selected profile, the identity can be written as
\[
 e=\frac{2(\mathscr E+4M)}{I_2(\delta)}.
\]
Thus the constant in the scale equation is determined by the conserved
energy and the second moment of the reference density. In particular,
$e=0$ gives $\mathscr E=-4M$ in our normalization.

The construction also gives collapse with positive energy in this
same normalization. Indeed, fix $e>8M_{\rm Ch}/I_2(0)$. As
$\beta\downarrow0$, the mass and second moment satisfy
\[
 M_\beta\longrightarrow M_{\rm Ch},\qquad
 I_2(-\beta^2/2)\longrightarrow I_2(0),
\]
and hence
\[
 \frac e2 I_2(-\beta^2/2)-4M_\beta
 \longrightarrow \frac e2 I_2(0)-4M_{\rm Ch}>0.
\]
Taking $\beta>0$ sufficiently small then gives a solution in
Theorem~\ref{thm:physical-main} with $\mathscr E>0$. Here $e$ is
chosen before the profile, as required by that theorem.
\subsection{Related work on gravitational collapse}
\label{in:critical-blowup}

Collapse solutions differ in their balance of forces, the mass
entering the singular region, and the arrival times of material
shells. These distinctions guide the comparison below. For a broader
account of collapse and expansion in the Euler--Poisson system,
see \cite{HadzicSurvey}.

\introheading{The white-dwarf model.}
Lieb--Yau derive the equilibrium theory
from a many-body quantum model \cite{LiebYauCMP,LiebYauApJ};
Hamada--Salpeter study corrections to the equation of state
\cite{HamadaSalpeter}, and Nauenberg gives approximations to the
mass--radius and energy relations \cite{Nauenberg}. In an actual
white dwarf, accretion or mergers can increase the mass; carbon
burning can cause thermonuclear disruption, while electron captures
can reduce pressure support and initiate collapse
\cite{HillebrandtNiemeyer,MaozMannucciNelemans,NomotoKondo}.
Burning, thermal evolution, and mass loss also affect merger remnants
\cite{SchwabQuataertKasen}. These processes are outside the
zero-temperature, fixed-composition model considered here.

\introheading{Weak solutions and stability of equilibria.}
For the white-dwarf law, Chen--Huang--Li--Wang--Wang
\cite[Theorem~2.3]{ChenHuangLiWangWang2024} construct radial
global finite-energy weak solutions from finite-energy initial data
with mass $M<M_{\rm Ch}$. This does not imply global classical
regularity.

The variational stability theory of gaseous equilibria
\cite{Rein,LuoSmoller2008} concerns motions near equilibrium.
For linearly stable nonrotating stars satisfying their nondegeneracy
assumption, Lin--Wang--Zhu \cite{LinWangZhu2025} prove orbital stability
under general perturbations, conditional on global weak existence and
continuity of the distance functional. Under stronger pressure
hypotheses, which include the white-dwarf law, they also construct
radial global finite-energy weak solutions that remain close to the
equilibrium. These results concern a different regime from the
trajectories constructed here, whose entire support concentrates in
finite time.

\introheading{Homologous collapse and expanding stars.}
For $P=\rho^{4/3}$, the collapsing Goldreich--Weber branches considered
here already have compact support, a physical vacuum, and simultaneous
collapse of all material shells
\cite{GoldreichWeber,Makino1992,FuLin1998}. Their mass measures therefore converge to the same point mass as in
Theorem~\ref{thm:physical-main}(ii). Homology means that every shell
moves by one scale factor; that factor need not be an exact power of
time. The inward branches with $e>0$, for example, have the leading
$(-t)^{2/3}$ rate but are not exactly self-similar
\cite[Section~2.2]{HadzicJangExpansion}.

For the outward branches, Had\v{z}i\'c--Jang proved radial nonlinear
stability of expanding Goldreich--Weber stars
\cite{HadzicJangExpansion}. Had\v{z}i\'c--Jang--Lam established
nonradial stability of the linearly expanding stars. They also proved
codimension-four stability of self-similarly expanding stars for
$\delta<0$ sufficiently close to zero, among irrotational perturbations
with the same energy and momentum as the reference star
\cite{HadzicJangLamExpansion}.

In a different expansion regime,
Had\v{z}i\'c--Jang constructed global nonsymmetric Euler--Poisson
solutions with small initial density, a physical vacuum, and an
asymptotically linear expansion \cite{HadzicJangGlobal}. Their pressure
is $P=\rho^\gamma$, with $1<\gamma<14/13$ or
$\gamma=1+1/n$ for an integer $n\ge2$. These results use decay
associated with expansion. For the self-similar collapsing reference,
infinitely many acoustic modes grow relative to its radius; the linear
calculation is given in Subsection~\ref{in:proof-strategy}.

\introheading{Isothermal collapse.}
The classical Larson--Penston profiles describe isothermal collapse
\cite{Larson,Penston}. Their rigorous construction, including passage
through the sonic point, is due to Guo--Had\v{z}i\'c--Jang
\cite{GuoHadzicJangLP}. Let $\tau=T-t$ be the remaining physical
time. For normalization constants $A,c_s>0$, their density has the form
\[
 \rho_{\rm LP}(t,r)=A\tau^{-2}D_{\rm LP}(r/(c_s\tau)).
\]
For fixed $Y>0$,
\[
 4\pi\int_0^{c_s\tau Y}\rho_{\rm LP}(t,r)r^2\dd r
 =4\pi A c_s^3\tau\int_0^Y D_{\rm LP}(y)y^2\dd y.
\]
Thus the mass in a bounded similarity region tends to zero. In our
construction the central density has the same order $\tau^{-2}$,
but the radius is of order $\tau^{2/3}$. Each fixed interior ball in
the rescaled Goldreich--Weber coordinates contains a positive limiting
fraction of the mass; contraction of the full support gives the
concentration of all the mass. The density blowup rate alone does
not distinguish these two mass distributions.

The exact Larson--Penston profile has an infinite-mass tail, but its
local collapse is compatible with finite mass and energy.
Guo--Had\v{z}i\'c--Jang--Schrecker prove radial nonlinear stability,
modulo the collapse time, in weighted spaces that allow suitable
modifications of the exterior
\cite[Theorems~2.16 and~2.22, Remark~2.24]{GuoHadzicJangSchreckerLPStability}.
Hunter's numerical study identified a discrete family of isothermal
self-similar profiles beyond the Larson--Penston solution \cite{Hunter}.
Sandine constructed infinitely many of the higher Hunter profiles,
with analytic passage through a single sonic point \cite{SandineHunter}.
These are whole-space profiles with infinite mass, whereas the
solutions constructed here have a material vacuum boundary.

\introheading{Continued collapse and polytropic profiles.}
For $1<\gamma<4/3$, Guo--Had\v{z}i\'c--Jang construct
physical-vacuum solutions undergoing continued gravitational collapse
\cite{GuoHadzicJangContinued}. The center becomes singular first,
and successive material shells reach the origin at different times.
The mass is continuously absorbed until the support disappears.
Their construction also uses a high-order approximate solution and
backward evolution. In our solutions all shells remain in the classical
fluid up to the same collapse time. For the pressure law
considered here, the leading approximation must retain pressure as
well as gravity, as the following scaling comparison shows.

Under a fixed-mass compression
$\rho_\lambda(X)=\lambda^{-3}\bar\rho(X/\lambda)$, the pressure
acceleration for $P=K\rho^\gamma$ has size
$\lambda^{2-3\gamma}$, while gravitational acceleration has size
$\lambda^{-2}$. The relative factor is
\[
                         \lambda^{4-3\gamma}.
\]
It tends to zero for $\gamma<4/3$, as in the dust approximation used
for continued collapse. At $\gamma=4/3$, both forces remain of the
same order as the acceleration for a $\tau^{2/3}$ collapse scale.
The pressure must therefore be retained in the leading equation.

Other scalings give different balances. Yahil's numerical polytropic
models connect an approximately homologous inner core to a
supersonically infalling outer region \cite{Yahil}.
Guo--Had\v{z}i\'c--Jang--Schrecker constructed smooth Yahil-type
self-similar profiles for every $1<\gamma<4/3$, retaining inertia,
pressure, and gravity \cite{GuoHadzicJangSchrecker}. These profiles
have positive density at every finite radius and hence no material
vacuum boundary. They pass through a single sonic point and have
infinite total mass; the mass in a fixed bounded similarity region
is of order $\tau^{4-3\gamma}$. The fixed-mass comparison above
therefore does not describe their force balance.

Nonconstant entropy permits another collapse regime.
Alexander--Had\v{z}i\'c--Schrecker constructed smooth spherical
self-similar solutions of the full Euler--Poisson system on
$\mathbb R^3$ for $19/12<\gamma<11/6$, with central density blowup
and flow that is supersonic relative to the similarity coordinates
away from the origin \cite{AlexanderHadzicSchrecker}. Their entropy
is transported by the flow, and the pressure is not a fixed function
of density alone. This differs from the barotropic law
\eqref{in:exact-pressure}.

\introheading{Euler implosion.}
For polytropic Euler without self-gravity,
Merle--Rapha\"el--Rodnianski--Szeftel construct smooth radial
self-similar profiles \cite{MerleRaphaelRodnianskiSzeftelI} and
finite-energy implosions \cite{MerleRaphaelRodnianskiSzeftelII}.
For the admissible pressure exponents and discrete similarity speeds,
the implosion theorem gives a finite-codimensional family of smooth
radial initial data. The density is positive on $\mathbb R^3$
before blowup, and the exterior of the self-similar core is modified
to obtain finite energy. The singularity forms only at the origin;
exterior matter remains at the blowup time.

For the nonisentropic Euler equations, Chen--Shkoller--Vicol
construct smooth self-similar implosions and prove nonlinear stability
of the ground-state profile under small smooth radial perturbations whose
initial pressure vanishes at the origin \cite{ChenShkollerVicol2026}.
The density is strictly positive, and there is no material vacuum
boundary or self-gravity. In the solutions constructed here, gravity
remains a leading force, and the entire support contracts to the origin.

\introheading{Related blowup constructions for critical equations.} Our construction is related to a number of earlier works in which an approximate solution is first constructed and then corrected to an exact solution. We recall a few examples relevant to the present argument.  Our construction draws in particular
on the method of Krieger--Schlag--Tataru
\cite{KriegerSchlagTataru} for critical wave maps: first construct
an accurate approximate solution by successive corrections, then
solve for the remainder. Bahouri--Marachli--Perelman
\cite{BahouriMarachliPerelman} adapted this approach to the
quasilinear wave equation governing hyperbolic vanishing mean
curvature flow using of matched asymptotic expansions, nonlinear energy estimates, and compactness. 
Energy estimates and compactness also underlie the constructions
of Rapha\"el--Szeftel and Le Coz--Martel--Rapha\"el for inhomogeneous
and double-power Schr\"odinger equations
\cite{RaphaelSzeftel,LeCozMartelRaphael}.
 Further examples include the $5$d Type~II blowup construction for the 
energy-critical wave equation by  Jendrej \cite{Jendrej2017}, and
compactnessthe matched expansions of
Ortoleva--Perelman \cite{OrtolevaPerelman} for the energy-critical
Schr\"odinger equation and the inner--outer gluing construction of
del Pino--Musso--Wei \cite{DelPinoMussoWei} for the energy-critical
heat equation.

For comparison with Type~II blowup, the homogeneous $4/3$
Euler--Poisson equation has the mass-preserving scaling
\[
 r(t,x)\longmapsto a\,r(a^{-3/2}t,x),\qquad a>0,
\]
whose self-similar length is $|t|^{2/3}$. Our interior retains this
rate. The smaller scale $\lambda^2$ describes the surface transition
between the two pressure regimes, not a faster temporal concentration
rate. We therefore do not describe these solutions as Type~II blowup.

The material vacuum boundary imposes constraints absent from the
whole-space constructions just discussed. The interior and surface
expansions must describe a single Lagrangian flow with prescribed
mass. Since $x$ is enclosed mass, density and radius satisfy
\[
 \dd x=4\pi r^2\rho_L\dd r,\qquad
 \rho_L=(4\pi r^2r_x)^{-1}.
\]
We therefore match the radius on the fixed interval $0\le x\le M$
and recover the density from this identity. At the surface,
$P'(\rho_L)$ vanishes linearly with the physical distance to the
boundary. Both the residual bounds and the evolution estimates must
use weights adapted to this degeneracy.

A small residual does not by itself make the approximate data
admissible for the free-boundary evolution. Successive time derivatives
of the equation impose boundary compatibility conditions on the
density and velocity. We correct the data to satisfy these conditions
while preserving the mass and the physical vacuum. Our bounds for
this correction lose powers of the collapse scale. The order of
approximation must compensate for these losses as well as for growth
during backward evolution. These requirements couple the asymptotic
construction to the local existence and continuation arguments.

\subsection{Outline of the proof}
\label{in:proof-strategy}

We first construct an approximate solution by matching the interior
and surface expansions. We then solve the evolution equation backward
from times approaching collapse, with data close to this approximation.
The residual must be small in the weighted norms used for the evolution. Local existence at each such time is insufficient: both
the lifespan and the estimates could deteriorate as the initial
radius tends to zero. We need a common earlier time, with bounds
that retain the interior and vacuum asymptotics. This requirement
determines the accuracy of the expansion, the construction of the
data, and the energy used in the evolution argument.

\introheading{Matched asymptotic expansions.}
Write $w=w_\delta$, $q=1-z$, and $\varkappa=-w'(1)$. Formula
\eqref{in:exact-coefficient} shows that the first relative pressure
correction has size $\lambda^2/w^2$. In the interior region
$q\gg\lambda$, we therefore seek an approximate Lagrangian radius
of the form
\[
 r_c(\lambda,m(z))=\lambda z
       \left[1+\sum_{n=2}^K\lambda^n\Phi_n(\log\lambda,z)\right],
\]
where $K$ is finite and each $\Phi_n$ is a polynomial in
$\log\lambda$. Substitution into the equation determines linear
equations for its coefficients. Their singular behavior as $z\to1$
must agree with the large-$D$ expansion of the boundary layer.

Even at leading order, the vacuum condition does not determine the
position of the boundary layer. Its equations give
\[
 P(\rho_0(D))=\frac{\varkappa}{\pi}D,\qquad
 Z_0(D)=\frac{h(\rho_0(D))}{4\varkappa}+C,
\]
with an undetermined constant $C$. The first equation determines
$\rho_0$; the second determines $Z_0$ only up to a translation.
To fix $C$, expand toward the interior. Since
\[
 P(\rho)=\rho^{4/3}-\rho^{2/3}+O(\log\rho),\qquad
 h(\rho)=4\rho^{1/3}-4+O(\rho^{-1/3})
 \quad(\rho\to\infty),
\]
we obtain
\[
 Z_0(D)=\frac1{\varkappa}
            \left(\frac{\varkappa D}{\pi}\right)^{1/4}
          +\left(C-\frac1{\varkappa}\right)+O(D^{-1/4})
 \quad(D\to\infty).
\]
The leading term agrees with the boundary expansion of the
homogeneous radius, because
$M-m(1-q)=\pi\varkappa^3q^4(1+O(q))$.
The constant term contributes
$-\lambda(C-1/\varkappa)$ to the relative radius $r/\lambda=1-\lambda Z$.
The interior expansion has no term of order $\lambda$, so matching
requires $C=1/\varkappa$. At the surface $\rho_0(0)=0$ and $h(0)=0$;
hence $Z_0(0)=1/\varkappa$. This is the coefficient of the boundary
displacement in \eqref{in:boundary-radius}; the full coefficient
calculation is given in \eqref{mat:connection-28}.

At higher orders, regularity at the center leaves one free coefficient
in the homogeneous solution of each core equation. The vacuum condition
leaves one additive constant in the corresponding layer equation.
We prove that the first coefficient changes the singular boundary term
by a nonzero factor. It therefore determines the singular part of the
matching; the layer constant determines the regular part.

We solve the coefficient equations successively in powers of
$\lambda$ and $\log\lambda$; see Subsection~\ref{mat:connection}.
The singular terms forced by the core equations are included in this
matching.

We join the expansions at $q\sim\sqrt\lambda$, where both $q$
and $\lambda/q$ are small, and denote the resulting radius by $A$.
Its density is defined by $(4\pi A^2A_x)^{-1}$, so the mass remains
exactly $M$. The joining region is wider than the physical layer
$q\sim\lambda$ described above. Increasing the finite order $K$
improves the error in the equation, including the spatial and time
derivatives needed below; see Proposition~\ref{mat:matched-source}.
The construction constants may depend on $\beta$ and $e$. After fixing
these parameters, we choose a sufficiently small upper bound for $\lambda$.

\introheading{Backward evolution.}
The obstruction to a forward perturbation argument is already
visible in the linearization about the homogeneous collapse with
$e=0$. Write $r=\lambda z(1+g)$, so that $g$ is the relative
radial displacement, and introduce the forward rescaled time
$\theta$ by $\dd t/\dd\theta=\lambda^{3/2}$. Then
$\lambda_\theta=-\beta\lambda$.

The linearized equation separates in a spatial eigenbasis
$\{\phi_j\}$ with eigenvalues $\mu_j\to\infty$;
see Lemma~\ref{en:homogeneous-acoustics}. These eigenfunctions
describe radial acoustic oscillations. Writing
$g(\theta,z)=\sum_j a_j(\theta)\phi_j(z)$, the linearized
equation gives
\[
 a_j''-\frac\beta2 a_j'+(\mu_j+3\delta)a_j=0,
 \qquad \delta=-\frac{\beta^2}{2}.
\]
Here primes denote $\theta$-derivatives. Substitution of
$a_j(\theta)=e^{\sigma\theta}$ gives the characteristic roots
\begin{equation}
 \sigma_j^\pm=\frac\beta4
       \pm\sqrt{\frac{25\beta^2}{16}-\mu_j},\qquad
 \operatorname{Re}\sigma_j^\pm=\frac\beta4
             \quad\left(\mu_j>\frac{25\beta^2}{16}\right).
 \label{in:acoustic-roots}
\end{equation}
For these eigenvalues the amplitudes oscillate with envelope
$e^{\beta\theta/4}$. Since $\partial_\theta\log\lambda=-\beta$,
this envelope is proportional to $\lambda^{-1/4}$.
Thus infinitely many components of the relative displacement grow
as collapse is approached. Removing finitely many growing modes
cannot give a uniform bound for the remaining linear perturbation.
This calculation concerns the homogeneous reference solution; see
Lemma~\ref{en:homogeneous-acoustics} and
\cite[Section~IV]{GoldreichWeber}.

We instead prescribe increasingly accurate data at $\lambda=a>0$
and solve toward earlier physical times. Put $\tau=-t$, and define
$s$ by $\dd\tau/\dd s=\lambda^{3/2}$. Then
\[
 \lambda_s=b\lambda,\qquad b=\sqrt{\beta^2+e\lambda}>0.
\]
We seek an interval $a\le\lambda\le A_0$ with $A_0$ independent
of $a$. The radius increases along this reversed evolution; in
physical time the resulting star is contracting.

\introheading{Compatible data and local solutions.}
The approximate solution cannot simply be used as initial data.
At a physical vacuum, successive time derivatives of the equation
impose boundary conditions on the density and velocity. We correct
finitely many boundary coefficients in the order in which they enter
these conditions. A separate perturbation supported in the interior
adjusts the mass without changing the boundary coefficients. The Jacobian matrix of this system is triangular and invertible.
A small residual therefore gives small corrections satisfying all
the required conditions; see Lemma~\ref{loc:finite-map}.

For the local construction, let $R(\tau,x)=r(-\tau,x)$ and define
the pressure and gravitational force by
\[
 \rho_R=(4\pi R^2R_x)^{-1},\qquad
 \mathcal M[R]=4\pi R^2\partial_xP(\rho_R)+\frac{x}{R^2}.
\]
We first solve the Kelvin--Voigt regularization
\[
 R_{\tau\tau}+\mathcal M[R]
       +\kappa\bigl(D\mathcal M[R]+\Lambda I\bigr)R_\tau=0,
 \qquad \kappa>0.
\]
Here $D\mathcal M[R]$ is the linearization of the force, and
$\Lambda\ge0$ makes the viscous quadratic form positive. This
choice preserves the identity
$D\mathcal M[R]R_\tau=\partial_\tau\mathcal M[R]$.
For each fixed $a>0$ and $\Lambda$, a further small correction
imposes the compatibility conditions for all sufficiently small
$\kappa>0$. For each such viscosity, the regularized equation has
a unique local solution with continuous dependence. Its continuation
criterion requires bounds on the radius and velocity, together with
quantitative preservation of the center and physical-vacuum conditions;
see Proposition~\ref{loc:high-local}. The estimates below
provide a common interval, first as $\kappa\downarrow0$ at fixed $a$,
and then as $a\downarrow0$.

\introheading{Energy estimates and continuation.}
The difficulty is to control the nonlinear force without losing a
spatial derivative in estimates uniform in viscosity. Let $Y=R-A$
be the error, expressed in rescaled time $s$, and put
$Y_j=\partial_s^jY$. Define the commutator $C_j$ by separating the
linearized term in the differentiated force difference:
\[
 \lambda^3\partial_s^j\bigl(\mathcal M[R]-\mathcal M[A]\bigr)
       =\lambda^3D\mathcal M[R]Y_j+C_j.
\]
The energy identity pairs $C_j$ with $Y_{j+1}$.
At the highest order, the uniform energy controls $Y_{j+1}$ only
in $L^2$; a direct bound of the pairing would require one more
weighted spatial derivative.

We write $C_j[f]$ for the dual pairing of $C_j$ with $f$.
The correction to the energy is suggested by the product rule
\[
 C_j[Y_{j+1}]
   =\partial_s\bigl(C_j[Y_j]\bigr)-C_{j,s}[Y_j].
\]
It replaces the problematic pairing with one involving $Y_j$, for
which the spatial derivative is controlled. The actual correction
includes the powers of $\lambda$ required by the rescaling, and
all terms from their differentiation are kept in the identity
\eqref{en:primitive-cancellation}. Lemma~\ref{en:primitive-bounds}
estimates the remaining commutators. When the lower-order errors
are small, Proposition~\ref{en:energy-coercivity} shows that the
corrected energy is equivalent to the positive quadratic energy.

We obtain the additional spatial bounds from the equation, solving
for spatial derivatives in terms of time derivatives already controlled
by the energy. These estimates supply the pointwise bounds on the
relative errors in $R$ and $R_x$ needed in the nonlinear coefficients.
The energy and spatial estimates close together: each improves one
of the two assumed bounds on the same interval. For the viscous
solutions, the uniform spatial estimate is used only through order
eight. Higher spatial bounds are recovered after removal of viscosity.

The required expansion order is determined by the allowed growth of
the error. Let $\widetilde E$ denote the modified energy and
$\Gamma_0$ its growth exponent. If the modified energy at $\lambda=a$ is $O(a^p)$
and the differentiated residual has order $\lambda^{p/2}$,
the energy estimate gives, for $p>\Gamma_0$,
\[
 \begin{aligned}
 \widetilde E(\lambda)
 &\le \left(\frac\lambda a\right)^{\Gamma_0}\widetilde E(a)
   +C\lambda^{\Gamma_0}
        \int_a^\lambda\frac{u^{p-\Gamma_0-1}}{\beta^2+eu}\dd u\\
 &\le \left(\frac\lambda a\right)^{\Gamma_0}\widetilde E(a)
   +\frac{C}{\beta^2(p-\Gamma_0)}\lambda^p.
 \end{aligned}
\]
Since $a\le\lambda$, the initial contribution satisfies
$\lambda^{\Gamma_0}a^{p-\Gamma_0}\le\lambda^p$.
High accuracy near collapse therefore compensates for growth during
backward evolution. In the proof we write $p=2J+1$: a bound of
order $\lambda^{2J+1}$ for the quadratic energy gives a bound of
order $\lambda^{J+1/2}$ for the corresponding error norm.
The integer $J$ specifies the
required decay of the error; the truncation order $K$ of the matched
expansion is chosen larger to allow for the derivatives used in
the estimates.

The constants governing energy growth and spatial recovery are fixed
before $J$ and $K$. We then choose $J$ so that $2J+1>\Gamma_0$
and the spatial estimates also close, and choose $K$ so that the
residual and the initial errors have additional positive powers of
$\lambda$; the precise choices are in \eqref{ct:orders}.
Decreasing $A_0$ makes these errors small. The continuation criterion
then extends every regularized solution to $\lambda=A_0$;
see Sections~\ref{en:section} and~\ref{ct:section}.

\introheading{Passage to the limit.}
We first let $\kappa\downarrow0$ with $a>0$ fixed. The uniform
energy and lower-order spatial estimates identify the nonlinear
force in the limit. Before applying the full inviscid spatial
estimate, we use the limiting equation to prove that its derivatives
belong to the domains of the linearized operators. This step may depend on
$a$; the resulting high-order bounds are uniform in $a$. We then let
$a\downarrow0$, taking limits on nested intervals bounded away
from $\lambda=0$, and obtain a classical Euler--Poisson solution. This procedure proves
existence of the constructed trajectory, without asserting a general
inviscid local well-posedness or stability theorem.

The estimates also imply convergence to the boundary-layer profile. For each fixed $Y_*>0$, with $J$ the order chosen in the
energy argument,
\begin{equation}
 \left\|\frac{(R-A)(s,M-\lambda^4y^5)}{\lambda^2}
           \right\|_{C^3([0,Y_*])}
     \le C_{\beta,e,Y_*}\lambda^{J+1/2}\longrightarrow0.
 \label{in:actual-layer-transfer}
\end{equation}
The approximate radius satisfies
\[
 \frac{\lambda-A(s,M-\lambda^4y^5)}{\lambda^2}
       \longrightarrow Z_0(y^5)\quad\text{in }C^3([0,Y_*]).
\]
Together with the exact density formula, these two estimates prove the radius and
density limits in Theorem~\ref{thm:physical-main}(iii). The core
limits follow from the interior bounds, and the separate global
energy comparison gives \eqref{in:energy-identity}.

The proof is organized according to these dependencies:
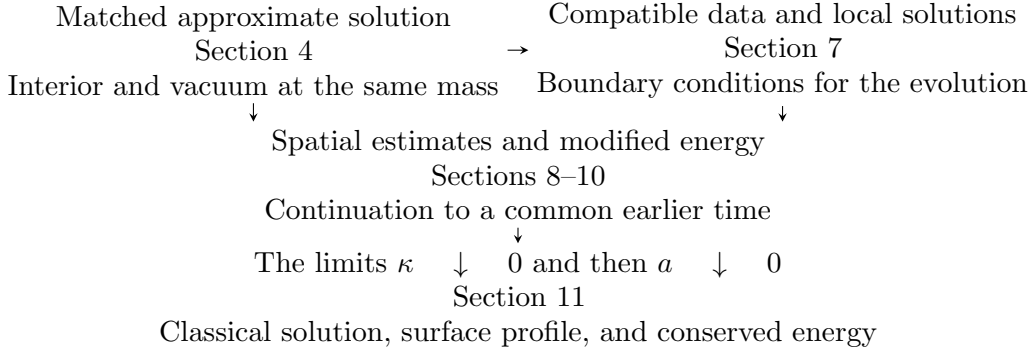
\begin{center}
\begin{tikzpicture}[>=stealth,
 every node/.style={align=center,font=\small,inner sep=3pt}]
 \node[text width=6.5cm] (profile) at (-3.5,0)
 {Matched approximate solution\\Section~\ref{mat:chapter}\\
 Interior and vacuum at the same mass};
 \node[text width=6.5cm] (local) at (3.5,0)
 {Compatible data and local solutions\\Section~\ref{loc:section}\\
 Boundary conditions for the evolution};
 \node[text width=13cm] (bounds) at (0,-1.65)
 {Spatial estimates and modified energy\\
 Sections~\ref{rec:chapter}--\ref{ct:section}\\
 Continuation to a common earlier time};
 \node[text width=13cm] (limits) at (0,-3.3)
 {The limits $\kappa\downarrow0$ and then $a\downarrow0$\\
 Section~\ref{exlim:section}\\
 Classical solution, surface profile, and conserved energy};
 \draw[->] (profile.east) -- (local.west);
 \draw[->] (profile.south) -- (profile.south |- bounds.north);
 \draw[->] (local.south) -- (local.south |- bounds.north);
 \draw[->] (bounds.south) -- (limits.north);
\end{tikzpicture}
\end{center}
The weighted spaces and linear estimates used in both the preparation
and evolution are developed in Sections~\ref{at:section},
\ref{at:calculus-section}, and~\ref{pr:section}. Auxiliary weighted
estimates, profile calculations, and local existence arguments are
collected in Appendices~\ref{app:weighted-estimates},
\ref{app:profile-estimates}, and~\ref{app:local-theory}, respectively.
An index of notation follows the bibliography.

\subsection{Acknowledgement}
G.C. was partially supported by NSF Grants DMS-2350301 and CAREER-DMS-2540992, by the Simons Foundation MP-TSM-00002258, and by the AMS Stefan Bergman Fellowship.
Z.L.  is supported in part by the National Natural Science Foundation of China (No. 12494544). We would like to thank Chongchun Zeng for many enlightening discussions.

\subsection*{AI Use Disclosure} 
The authors used OpenAI's ChatGPT as an assistive tool for language editing, improving the exposition, and checking intermediate calculations and arguments for possible errors or omissions. All mathematical statements, proofs, calculations, and references in the manuscript were independently verified by the authors, who take full responsibility for the content and correctness of the paper.
\section{The equation and main theorem}
\label{sec:module-formulation}

\subsection{The model and its normalization}
\label{sec:normalization}

We derive the dimensionless pressure law and the equation in mass
coordinates. The pressure expansions enter the profile construction
in Section~\ref{sec:EOS}. Let $\varrho$ and $p_{\rm phys}$ denote the
physical mass density and pressure. The zero-temperature Chandrasekhar
law is
\begin{equation}
 p_{\rm phys}(\varrho)=A_{\rm Ch}
 \left\{q(2q^2-3)\sqrt{1+q^2}+3\operatorname{arsinh}q\right\},
 \qquad q=\left(\frac{\varrho}{B_{\rm Ch}}\right)^{1/3},
 \label{eq:dimensional-Chandrasekhar-law}
\end{equation}
Here $m_e$ is the electron mass, $m_u$ the atomic mass unit, and
$\mu_em_u$ the mean mass per electron; $c$ and $\hbar$ denote the
speed of light and the reduced Planck constant. The constants are
\begin{equation}
 A_{\rm Ch}=\frac{m_e^4c^5}{24\pi^2\hbar^3},\qquad
 B_{\rm Ch}=\frac{\mu_em_um_e^3c^3}{3\pi^2\hbar^3}.
 \label{eq:dimensional-constants}
\end{equation}
Let $n_e$ be the electron number density and $p_F=m_ecq$ the Fermi
momentum. Filling the two electron spin states up to $p_F$ gives
\[
 \begin{aligned}
 n_e&=\frac{p_F^3}{3\pi^2\hbar^3},\\
 p_{\rm phys}&=\frac1{3\pi^2\hbar^3}
       \int_0^{p_F}\frac{p^4c^2}{\sqrt{m_e^2c^4+p^2c^2}}\dd p\\
 &=\frac{m_e^4c^5}{3\pi^2\hbar^3}
       \int_0^q\frac{s^4}{\sqrt{1+s^2}}\dd s.
 \end{aligned}
\]
The antiderivative in \eqref{eq:dimensional-Chandrasekhar-law} and
$\varrho=\mu_em_un_e$ yield \eqref{eq:dimensional-constants}.
This fixes both the pressure normalization and its dependence on composition.
If $G$ is Newton's constant, set
\begin{equation}
 R_{\rm Ch}^2=\frac{2A_{\rm Ch}}{GB_{\rm Ch}^2},\qquad
 M_{\rm Ch,*}=B_{\rm Ch}R_{\rm Ch}^3,\qquad
 T_{\rm Ch}=(GB_{\rm Ch})^{-1/2}.
 \label{eq:dimensional-scales}
\end{equation}
Under
\[
 r_{\rm phys}=R_{\rm Ch}r,\quad
 m_{\rm phys}=M_{\rm Ch,*}x,\quad
 t_{\rm phys}=T_{\rm Ch}t,\quad
 \varrho=B_{\rm Ch}\rho,\quad
 p_{\rm phys}=2A_{\rm Ch}P,
\]
the density identity becomes $\rho=(4\pi r^2r_x)^{-1}$, and
the two dimensionless coefficients are
\[
 \frac{2A_{\rm Ch}T_{\rm Ch}^2}{B_{\rm Ch}R_{\rm Ch}^2}=1,
 \qquad GB_{\rm Ch}T_{\rm Ch}^2=1.
\]
Thus the pressure and gravitational terms are $4\pi r^2P_x$ and
$x/r^2$, respectively, as in \eqref{eq:EP}. The mass unit $M_{\rm Ch,*}$
converts the dimensionless critical mass to
\[
 M_{{\rm Ch},{\rm phys}}=M_{\rm Ch,*}M_{\rm Ch}
 =\left(\frac{2A_{\rm Ch}}{G B_{\rm Ch}^{4/3}}\right)^{3/2}M_{\rm Ch}.
\]
In particular, $M_{{\rm Ch},{\rm phys}}\propto\mu_e^{-2}$ at fixed
fundamental constants. These dimensional units are fixed independently of the evolving boundary
$R_b(t)$ and of the reference scales used in the construction.

Let $x\in[0,M]$ denote enclosed mass, and let $r(\tau,x)$ be the
particle radius in reversed time $\tau=-t$. Then
\begin{equation}
 \rho=(4\pi r^2 r_x)^{-1},\qquad u=r_\tau,
 \label{eq:rho-u}
\end{equation}
and the radial Euler--Poisson equation is
\begin{equation}
 r_{\tau\tau}+4\pi r^2\pa_x P(\rho)+\frac{x}{r^2}=0.
 \label{eq:EP}
\end{equation}
Indeed, $\pa_r=r_x^{-1}\pa_x$, $m(\tau,r(\tau,x))=x$, and
$\rho^{-1}r_x^{-1}=4\pi r^2$.
The radius equation is invariant under time reversal.  Thus
\eqref{eq:EP} also holds with the collapsing time $t=-\tau$; the
corresponding collapsing velocity is $r_t=-r_\tau$.

We use the exact normalized Chandrasekhar pressure.  With $q=\rho^{1/3}$,
\begin{equation}
 P(\rho)=\frac12\left[q(2q^2-3)\sqrt{1+q^2}+3\operatorname{arsinh}q\right].
 \label{eq:ChP}
\end{equation}
Writing $S=(1+q^2)^{1/2}$, one computes
\[
 \frac{\dd}{\dd q}\{q(2q^2-3)S+3\operatorname{arsinh}q\}
 =\frac{8q^4}{S}.
\]
For $\rho>0$, the identity $\dd\rho/\dd q=3q^2$ gives
\begin{equation}
 P'(\rho)=\frac43\frac{\rho^{2/3}}{\sqrt{1+\rho^{2/3}}}>0.
 \label{eq:Pprime}
\end{equation}
Define the sound speed and the specific enthalpy by
\begin{equation}
 c(\rho)=\sqrt{P'(\rho)},\qquad
 h(\rho)=\int_0^\rho\frac{P'(s)}s\dd s.
 \label{eq:ch}
\end{equation}
The substitution $s=y^3$ yields the exact identities
\begin{align}
 h(\rho)&=4\left(\sqrt{1+\rho^{2/3}}-1\right),
 \label{eq:h-exact}\\
 \rho(h)&=\left(\frac h2+\frac{h^2}{16}\right)^{3/2},
 \label{eq:rho-h}\\
 P'(\rho(h))&=h\,a_{\rm eos}(h),\qquad
 a_{\rm eos}(h)=\frac{8+h}{3(4+h)}.
 \label{eq:a-exact}
\end{align}
Consequently,
\begin{equation}
 \frac13\le a_{\rm eos}(h)\le\frac23,
 \qquad
 \sup_{h\ge0}\abs{\pa_h^m a_{\rm eos}(h)}\le C_m
 \quad(m\ge0).
 \label{eq:a-bounds}
\end{equation}
\paragraph{Internal energy.}
We normalize the specific internal energy to vanish at vacuum.
For $\rho>0$ set
\begin{equation}
 \begin{gathered}
 \mathfrak e(\rho)=h(\rho)-\frac{P(\rho)}{\rho},\qquad
 \mathfrak e(0)=0,\qquad
 \mathfrak e'(\rho)=\frac{P(\rho)}{\rho^2},\\
 F(\rho)=\rho\mathfrak e(\rho)=\int_0^\rho h(s)\dd s.
 \end{gathered}
 \label{eq:internal-energy-normalization}
\end{equation}
The derivative identity follows from $h'=P'/\rho$; the value at zero
follows from the $5/3$ pressure law. Differentiating $\rho h-P$
and using its zero value at vacuum proves the formula for $F$.
The exact enthalpy therefore yields
\begin{equation}
 \begin{split}
 F(\rho)
 &=4\int_0^\rho\bigl(\sqrt{1+s^{2/3}}-1\bigr)\dd s\\
 &\ge3\rho^{4/3}-4\rho\qquad(\rho\ge0).
 \end{split}
 \label{eq:threshold-internal-energy}
\end{equation}
More precisely,
\begin{equation}
 0\le F(\rho)+4\rho-3\rho^{4/3}\le3\rho^{2/3}.
 \label{eq:threshold-exact-law-defect}
\end{equation}
Indeed, subtract $4s^{1/3}$ inside the integral and use
$0<\sqrt{1+s^{2/3}}-s^{1/3}\le\tfrac12s^{-1/3}$ for $s>0$.
The upper bound is integrable at zero, so the estimate holds at all
densities. Combined with the sharp Newtonian inequality and energy
conservation, the lower bound excludes subcritical collapse of the full
support. The upper bound will be used to compare the exact and polytropic
energies.

The energy shift $4M$ has a dimensional interpretation. By
\eqref{eq:dimensional-constants}--\eqref{eq:dimensional-scales}, the
energy unit is $2A_{\rm Ch}R_{\rm Ch}^3$ and
$8A_{\rm Ch}/B_{\rm Ch}=m_ec^2/(\mu_em_u)$. Write $M_{\rm phys}$
for the dimensional mass and $N_e=M_{\rm phys}/(\mu_em_u)$ for the
electron number. Adding $4\rho$ to $F$ changes the total energy by
\[
 4M\bigl(2A_{\rm Ch}R_{\rm Ch}^3\bigr)
 =\frac{m_ec^2}{\mu_em_u}M_{\rm phys}
 =N_em_ec^2.
\]
The normalization $\mathfrak e(0)=0$ therefore measures the electron
kinetic energy. Since
\[
 \rho(F+4\rho)'-(F+4\rho)=\rho F'-F=P,
\]
adding the electron rest energy changes neither the pressure nor the
evolution.

Hereafter $\beta>0$ is selected from the Goldreich--Weber family in
Subsection~\ref{sec:GW}. Fix $e\ge0$. The scale increases in the
reversed time $\tau$; define it and the rescaled time by
\begin{equation}
 \la_\tau=\sqrt{e+\frac{\beta^2}{\la}},\qquad
 \frac{\dd s}{\dd\tau}=\la^{-3/2},\qquad
 \la_s=b(\la)\la,
 \label{eq:scale}
\end{equation}
where
\begin{equation}
 b(\la)=\sqrt{\beta^2+e\la},\qquad b_s=\frac{e\la}{2},\qquad b\ge\beta>0.
 \label{eq:b}
\end{equation}

After the construction parameters have been fixed, we reduce the upper
bound for the scale so that
\begin{equation}
 0<\beta\le b(\lambda)\le1\qquad(0<\lambda\le A_{\rm out}).
 \label{eq:b-upper}
\end{equation}
This normalization is used only to replace $b^q$ by $b$ for $q\ge1$ in
the finite commutator sums.
Unless a weight is displayed, we use
\[
 (f,g)_2=\int_0^M f(x)g(x)\dd x,\qquad
 \|f\|_2=(f,f)_2^{1/2}.
\]
The operators act in the mass variable $x$.
The profile and evolution estimates below are taken on
$0<\la\le A_{\rm out}$, with $A_{\rm out}>0$ sufficiently small.

In the asymptotic expansions, $p_{\log}$ is a nonnegative integer
bounding the powers of the logarithm in the estimate under consideration.
It may be increased when finitely many such estimates are combined.
For each family constructed in the theorem, $C_{\log}$ denotes the
resulting finite logarithmic degree.

The letter $C$ denotes a positive constant that may change from line to line.
Its fixed dependencies are stated at the beginning of each major
construction. Constants used locally within a section are distinguished
from those retained in subsequent arguments.

Return to the collapsing time $t=-\tau$.  Then
\begin{equation}
 \dot\la(t)=-\sqrt{e+\frac{\beta^2}{\la(t)}},\qquad \la(t)\downarrow0\quad(t\uparrow0).
 \label{eq:collapse-scale}
\end{equation}
For $A>0$, put
\begin{equation}
 T_A=\int_0^A\frac{\dd\mu}{\sqrt{e+\beta^2/\mu}}.
 \label{eq:TA}
\end{equation}

\subsection{The classical solution class and the main theorem}
\label{sec:class-and-main}

We state the classical solution class in coordinates adapted to the
center and the vacuum boundary. The weighted spaces used to construct
such solutions are defined in Section~\ref{at:section}.

\begin{definition}
\label{def:global-solution}
Fix an even integer $N\ge22$.  A radius
$r(t,x)$ on an interval $I$ is an order-$N$ classical radial
physical-vacuum solution if the following conditions hold.
Here $N$ is the weighted regularity index used in the construction;
the ordinary endpoint regularity in (ii)--(iii) has order $N-6$.
\begin{enumerate}[label=(\roman*)]
\item \emph{Interior equation.} $r\in C^2(I\times(0,M))$, $r(t,0)=0$, and
$r_x(t,x)>0$ for $0<x<M$.  The density and radial velocity are
\[
 \rho_L=(4\pi r^2 r_x)^{-1},\qquad u_L=r_t
\]
and $r$ satisfies \eqref{eq:EP} pointwise on $I\times(0,M)$.
\item \emph{Regularity at the center.} With $z_c=x^{1/3}$, the quotient $r(t,x)/z_c$ and the scalar
functions $\rho_L(t,x)$ and $u_L(t,x)/z_c$ have even one-sided extensions in
$z_c$ near $z_c=0$.  For every compact $I'\Subset I$ and every small fixed
$z_0>0$,
\[
 \sup_{t\in I'}\sum_{j+k\le N-6}
 \bigl\|\partial_t^j\partial_{z_c}^k
 (r/z_c,\rho_L,u_L/z_c)\bigr\|_{L^\infty(0,z_0)}<\infty.
\]
This condition includes the center parity.
\item \emph{Regularity at the vacuum boundary.} There is a continuous
vacuum radius $R_b(t)=r(t,M)$ such that, with
$\xi=(M-x)^{1/5}$,
\begin{equation}
 \rho_L(t,x)=\xi^3\Theta(t,\xi),\qquad
 r(t,x)=R_b(t)-\xi^2G(t,\xi),
 \label{eq:global-solution-edge-class}
\end{equation}
where $\Theta,G$ are positive near $\xi=0$.  For every compact
$I'\Subset I$ and every sufficiently small fixed $\xi_0>0$,
\begin{equation}
 \sup_{t\in I'}
 \sum_{j+k\le N-6}
 \left(
 \|\partial_t^j\partial_\xi^k(\Theta,G)\|_{L^\infty(0,\xi_0)}
 +\|\partial_t^j\partial_\xi^k
       (\Theta^{-1},G^{-1})\|_{L^\infty(0,\xi_0)}
 \right)<\infty.
 \label{eq:explicit-edge-regularity}
\end{equation}
These bounds are taken in the coordinate $\xi$.
\item \emph{Physical-vacuum condition.} Let $x_E(t,\cdot)$ be the inverse of
$x\mapsto r(t,x)$ and put
\[
 \begin{aligned}
  \rho_E(t,r)&=\rho_L(t,x_E(t,r)),&
  u_E(t,r)&=u_L(t,x_E(t,r)),\\
  h_L(t,x)&=h(\rho_L(t,x)),&
  h_E(t,r)&=h_L(t,x_E(t,r)).
 \end{aligned}
\]
The Eulerian enthalpy extends to a $C^1$ function of the physical radius up
to the boundary and
\begin{equation}
 h_E(t,R_b(t))=0,\qquad
 -\partial_r h_E(t,R_b(t))>0.
 \label{eq:global-solution-pv}
\end{equation}
For every $j\le N-6$, the boundary identity
$\partial_t^j h_L(t,M)=0$ is understood through the factorization
$\partial_t^j h_L(t,x)=(M-x)^{2/5}\widehat h_j(t,(M-x)^{1/5})$
supplied by (iii).
Indeed, if $\rho_L=\xi^3\Theta$, the exact identity
\[
 h(\rho_L)=\xi^2
 \frac{4\Theta^{2/3}}{\sqrt{1+\xi^2\Theta^{2/3}}+1}
\]
makes this endpoint factor explicit; the stated regularity gives the
corresponding factorization for its time traces.
\end{enumerate}
Equations at the vacuum boundary are understood as the continuous limits
of the identities on $0<x<M$.
\end{definition}

For the constructed order-$22$ solution, Proposition~\ref{exlim:classical}
also gives local $C^{16}$ regularity of the material fields in the
interior. The class of collapsing solutions defined below requires
only the interior regularity in Definition~\ref{def:global-solution}.

\begin{theorem}
\label{thm:main-collapse}
\label{thm:main}
Fix $e\geq0$. There are a number $\beta_*(e)>0$ and finite
construction orders
\[
 J\ge3,\qquad K=2J+3512,
\]
independent of the profile parameter $\beta$, with the following
property. For every $0<\beta\le\beta_*(e)$, there exist
\[
 A_{\rm out}=A_{\rm out}(\beta,e)>0,\qquad C_{\log}\in\mathbb N_0
\]
for which the conclusions below hold. These choices are justified
in Lemma~\ref{ct:parameter-admissibility}.

Put
\[
 \delta=-\frac{\beta^2}{2}<0,
 \qquad
 M(\delta)=4\pi\int_0^1w_\delta(z)^3z^2\dd z,
\]
where $w_\delta$ is the Goldreich--Weber profile normalized to have its
first zero at one, and $\varkappa_\delta=-w_\delta'(1)>0$.  The constants
$A_{\rm out}$ and $C_{\log}$ depend only on $(\beta,e)$ and on the
finite construction orders fixed in the proof.

The exact Chandrasekhar Euler--Poisson system has an order-$22$
spherically symmetric classical
physical-vacuum solution in the sense of
Definition~\ref{def:global-solution} on
$[-T_{A_{\rm out}},0)$, with compact support $B_{R_b(t)}$ and
conserved mass $M(\delta)$.

For $0<A\leq A_{\rm out}$, define
\begin{equation}
 T_A=\int_0^A\frac{\mu^{1/2}}
 {\sqrt{\beta^2+e\mu}}\dd\mu,
 \qquad
 -t=\int_0^{\lambda(t)}
 \frac{\mu^{1/2}}{\sqrt{\beta^2+e\mu}}\dd\mu.
 \label{eq:main-collapse-clock}
\end{equation}
Then $\lambda(-T_A)=A$, and the restriction of this solution to
$[-T_A,0)$ has the following properties.
\begin{enumerate}[label=(\roman*),leftmargin=2em]
\item \emph{Mass and collapse scale.} The mass is strictly supercritical:
\begin{equation}
 M(\delta)>M_{\rm Ch}.
\label{eq:main-collapse-mass}
\end{equation}
The collapse scale obeys
\begin{equation}
 \lambda(t)=
 \left(\frac{3\beta}{2}(-t)\right)^{2/3}
 \left(1+O((-t)^{2/3})\right).
 \label{eq:main-collapse-scale}
\end{equation}

\item \emph{Core asymptotics.} Let $\rho_E,u_E$ be the Eulerian fields in
Definition~\ref{def:global-solution}.  Then, for every $0<z_0<1$,
uniformly on $0\leq r/\lambda(t)\leq z_0$,
\begin{equation}
 \rho_E(t,r)
 =\lambda(t)^{-3}
 w_\delta\!\left(\frac r{\lambda(t)}\right)^3(1+o(1)).
 \label{eq:main-collapse-density}
\end{equation}
The center parity gives $u_E(t,0)=0$. For positive radii,
\begin{equation}
 \sup_{0<r\le z_0\lambda(t)}
 \left|
 \frac{u_E(t,r)}
      {(\dot\lambda(t)/\lambda(t))r}-1
 \right|\longrightarrow0
 \qquad(t\uparrow0).
 \label{eq:main-collapse-velocity}
\end{equation}
In particular,
\begin{equation}
 \rho_E(t,0)=\frac{4w_\delta(0)^3}{9\beta^2}
 (-t)^{-2}(1+o(1)).
 \label{eq:main-collapse-center}
\end{equation}

\item \emph{Support and concentration.} The vacuum radius satisfies
\begin{equation}
 R_b(t)=\lambda(t)\left[1-
 \frac{\lambda(t)}{\varkappa_\delta}
 +O\!\left(\lambda(t)^2
 (1+|\log\lambda(t)|)^{C_{\log}}\right)\right]
 \longrightarrow0.
 \label{eq:main-collapse-support}
\end{equation}
If $\rho_E(t,\cdot)$ is extended by zero and
$\dd\mu_t(y)=\rho_E(t,|y|)\dd y$, then
\begin{equation}
 \mu_t\stackrel{*}{\rightharpoonup}
 M(\delta)\delta_0
 \qquad(t\uparrow0)
 \label{eq:main-collapse-measure}
\end{equation}
as finite Radon measures on $\mathbb R^3$.
\end{enumerate}

The parameter $e$ determines the reference scale in \eqref{eq:main-collapse-clock}.
Corollary~\ref{cor:physical-energy} identifies its relation to the
conserved physical energy, with the specified internal-energy normalization.
There is no uniform positive lower bound asserted for $A_{\rm out}$
as $\beta\downarrow0$.
The solution is asserted only for $t<0$.  No solution value or
physical-vacuum trace at $t=0$, and no uniqueness among all collapsing
solutions, is asserted.
\end{theorem}

Theorem~\ref{thm:main-collapse} and its consequence
Theorem~\ref{thm:physical-main} are proved in
Subsection~\ref{exlim:main-theorem-proof}, after the construction and
the asymptotic estimates.

The support estimate gives $r(t,x)\to0$ for every $0<x\le M$.
Together with mass conservation it implies the measure limit in
\eqref{eq:main-collapse-measure}. The local core estimate alone would
only identify the mass on compact subsets of the rescaled interior.

\begin{remark}
The mass of the leading profile satisfies $M_\beta\to M_{\rm Ch}$ as $\beta\downarrow0$.
Lemma~\ref{ct:parameter-admissibility} gives one choice of $J,K$ for
all sufficiently small positive profile parameters. The inverses in
the matching equations may be unbounded at zero. We therefore choose
the upper bound for the scale separately for each positive parameter,
in terms of these inverse bounds.

The mass of the leading profile has the quantitative expansion
\[
 M_\beta=M_{\rm Ch}+\frac{I_2(0)}{2M_{\rm Ch}}\beta^2+O(\beta^4),
 \qquad I_2(0)=4\pi\int_0^1w_0^3z^4\dd z>0,
\]
proved in Lemma~\ref{mat:mass-tangent}. In particular it is strictly
increasing on a possibly smaller positive parameter interval. The mass
corollary below uses only continuity, so does not require this refinement.
The logarithmic exponent is a finite degree supplied by this construction.
\end{remark}

\subsection{Prescribed masses and the boundary-layer limit}
\label{sec:mass-and-layer}

We choose the leading profile to have the prescribed mass.
Preparation and both limiting procedures preserve this mass.
Continuity in the lower-order weighted spaces then identifies the
boundary layer of the resulting solution.

\begin{corollary}
\label{cor:prescribed-supercritical-mass}\label{cor:prescribed-mass}
For each fixed $e\ge0$ there is $\varepsilon_M(e)>0$ such that every
prescribed mass
\[
 M_{\rm Ch}<M<M_{\rm Ch}+\varepsilon_M(e)
\]
is the conserved mass of a collapse solution with all the properties
of Theorem~\ref{thm:main-collapse}. The profile parameter and the upper bound for the scale may depend on
$M$ and $e$.
\end{corollary}

\begin{proof}
Use the orders and interval from
Lemma~\ref{ct:parameter-admissibility}, fix $0<\beta_1\le\beta_*$,
and put $\varepsilon_M(e)=M_{\beta_1}-M_{\rm Ch}>0$.
Subsubsection~\ref{mat:seed-part-7} proves that $M_\beta$ extends
continuously to $\beta=0$, with $M_0=M_{\rm Ch}$, and that
$M_\beta>M_{\rm Ch}$ for every positive $\beta$.
The intermediate value theorem on $[0,\beta_1]$ supplies a
$\beta\in(0,\beta_1)$ with $M_\beta=M$ for every mass in the displayed
open interval.

The orders in Lemma~\ref{ct:parameter-admissibility} apply to this
$\beta$, so Theorem~\ref{thm:main-collapse} gives the required solution. Preparation and both limits take place on the
fixed mass interval $(0,M)$. The argument uses continuity of
$M_\beta$, without requiring monotonicity.
\end{proof}

\begin{remark}
\label{rem:mass-parametrization}
After reducing $\varepsilon_M(e)$ if necessary, the leading profile realizing a
prescribed mass in Corollary~\ref{cor:prescribed-mass} is unique within
the small monotone branch of leading profiles. With $\Delta M=M-M_{\rm Ch}$,
\begin{equation}
 \begin{split}
  \beta(M)&=\left(\frac{2M_{\rm Ch}}{I_2(0)}\Delta M\right)^{1/2}
                     \bigl(1+O(\Delta M)\bigr),\\
  \delta(M)&=-\frac{M_{\rm Ch}}{I_2(0)}\Delta M+O(\Delta M^2).
 \end{split}
 \label{eq:mass-parametrization}
\end{equation}
Here the derivative in Lemma~\ref{mat:mass-tangent} is taken with
respect to $\delta$:
\[
 M'(0)=-I_2(0)/M_{\rm Ch}\ne0.
\]
The inverse-function theorem applies to the smooth function
$M(\delta)$; its inverse
$\delta(M)$ is smooth through $M=M_{\rm Ch}$ and has the displayed
Taylor expansion. On the supercritical side $\beta(M)=\sqrt{-2\delta(M)}$,
which gives the first formula. Thus $\beta$ is smooth for
$M>M_{\rm Ch}$ and has the stated square-root behavior at the endpoint.

For each fixed mass, substituting $\beta(M)$ in
\eqref{eq:main-collapse-scale} gives the collapse scale; its leading
coefficient can then be expanded as $M\downarrow M_{\rm Ch}$.
This parametrization determines the leading profile. It does not assert
uniqueness of the nonlinear solution or uniform bounds on the interval
of existence and the evolution error as the mass approaches the critical value.
\end{remark}

For the solution in Theorem~\ref{thm:main-collapse}, fix a finite
$Y>0$, put $M=M(\delta)$, and write $\lambda=\lambda(t)$.
For all sufficiently small positive $\lambda$, define
\begin{equation}
 \mathcal Z_\lambda(y)
 =\frac{\lambda-r(t,M-\lambda^4y^5)}{\lambda^2},
 \qquad 0\le y\le Y.
 \label{eq:actual-layer-Z}
\end{equation}
Let $\rho_0(D)$ be the nonnegative solution of
\begin{equation}
 P(\rho_0(D))=\frac{\varkappa_\delta}{\pi}D,
 \qquad
 Z_0(D)=\frac{h(\rho_0(D))}{4\varkappa_\delta}
                         +\frac1{\varkappa_\delta},
 \qquad \mathcal Z_0(y)=Z_0(y^5).
 \label{eq:actual-layer-leading}
\end{equation}
Define, by continuous extension at zero,
\begin{equation}
 \begin{aligned}
 \mathcal Q_\lambda(y)&=\frac{\partial_y\mathcal Z_\lambda(y)}y,
 &\mathcal G_\lambda(y)
   &=\frac{\mathcal Z_\lambda(y)-\mathcal Z_\lambda(0)}{y^2},\\
 \Theta_\lambda(y)
   &=\frac{\rho_L(t,M-\lambda^4y^5)}{y^3},
 &\Theta_0(y)&=\frac{\rho_0(y^5)}{y^3},
 \end{aligned}
 \label{eq:actual-layer-factors}
\end{equation}
and define $\mathcal Q_0,\mathcal G_0$ from $\mathcal Z_0$ in
the same way.

\begin{corollary}
\label{cor:actual-layer}
With this notation, the following conclusions hold.
\begin{enumerate}[label=(\roman*),leftmargin=2em]
\item \emph{Convergence.} As $t\uparrow0$ along the solution,
\begin{equation}
 \mathcal Z_\lambda\longrightarrow\mathcal Z_0
       \quad\hbox{in }C^3([0,Y]),\qquad
 (\mathcal Q_\lambda,\mathcal G_\lambda,\Theta_\lambda)
 \longrightarrow(\mathcal Q_0,\mathcal G_0,\Theta_0)
       \quad\hbox{in }C^1([0,Y]).
 \label{eq:actual-layer-convergence}
\end{equation}
More quantitatively, there are a finite integer $p_{\log}\ge0$ and a constant
$C_{\beta,e,Y}$ for the selected family such that
\begin{equation}
 \begin{split}
 &\|\mathcal Z_\lambda-\mathcal Z_0\|_{C^3([0,Y])}
 +\sum_{\mathcal F\in\{\mathcal Q,\mathcal G,\Theta\}}
          \|\mathcal F_\lambda-\mathcal F_0\|_{C^1([0,Y])}\\
 &\hspace{2em}\le C_{\beta,e,Y}
       \{\lambda(1+|\log\lambda|)^{p_{\log}}+\lambda^{J+1/2}\}.
 \end{split}
 \label{eq:actual-layer-rate}
\end{equation}
The first error comes from the finite matched profile and the second
from the nonlinear correction to that profile.
\item \emph{Positive factors.} The three regular factors have positive lower bounds on $[0,Y]$,
uniformly at sufficiently small scale. In particular,
\begin{equation}
 \Theta_0(0)=\left(\frac{5\varkappa_\delta}{4\pi}\right)^{3/5},
 \qquad
 \mathcal Q_0(0)=\frac{5}{4\pi\Theta_0(0)},
 \qquad
 \mathcal G_0(0)=\frac{\mathcal Q_0(0)}2.
 \label{eq:actual-layer-endpoint-values}
\end{equation}
\item \emph{Mass, thickness, and vacuum slope.} This layer contains exactly the mass $\lambda^4Y^5$, and, with
$R_b(t)=r(t,M)$, its thickness is
\begin{equation}
 R_b(t)-r(t,M-\lambda^4y^5)
       =\lambda^2y^2\mathcal G_\lambda(y).
 \label{eq:actual-layer-thickness}
\end{equation}
The physical-vacuum enthalpy slope satisfies
\begin{equation}
 \lambda(t)^2\bigl[-\partial_rh_E(t,R_b(t))\bigr]
       \longrightarrow4\varkappa_\delta=M+\delta.
 \label{eq:actual-layer-slope}
\end{equation}
\end{enumerate}
\end{corollary}

The proof in Subsection~\ref{exlim:actual-layer} uses the continuous
lower-order bounds of the solution and requires no further
subsequence. The convergence holds on each fixed $[0,Y]$, including
zero. No convergence is asserted on growing $y$-intervals, up to the
joining shell, or for material derivatives at fixed mass.

For fixed $D>0$, the density tends to the finite positive value
$\rho_0(D)$, without the core factor $\lambda^{-3}$. The leading
profile obeys the exact pressure law: its small-$D$ asymptotic
gives the $5/3$ vacuum behavior, and its large-$D$ asymptotic
matches the $4/3$ core. This last statement concerns the leading
profile, not convergence on an expanding interval.

The value $\mathcal Z_0(0)=1/\varkappa_\delta$ gives the first
support-radius correction in \eqref{eq:main-collapse-support};
the positive quotient factors give the $3/2$ density exponent and
the enthalpy slope in \eqref{eq:actual-layer-slope}.

\subsection{Conserved energy and the critical mass}
\label{sec:energy-and-threshold}

The local flux identity gives conservation of energy in the class
of Definition~\ref{def:global-solution}; see
Proposition~\ref{exlim:energy-conservation}. Comparison with the
matched profile determines the energy of the constructed solution.
Combined with the internal-energy bound
\eqref{eq:threshold-internal-energy}, the sharp Newtonian inequality
excludes collapse of the full support below the critical mass.

\begin{corollary}
\label{cor:physical-energy}
Use the specific internal energy fixed by
\eqref{eq:internal-energy-normalization}.
For the solution of Theorem~\ref{thm:main-collapse}, its conserved
physical energy satisfies
\begin{equation}
 \begin{gathered}
 \mathscr E+4M=\frac e2 I_2(\delta),\qquad
 I_2(\delta)=4\pi\int_0^1w_\delta^3z^4\dd z,\\
 \mathscr E=\int_0^M
      \left\{\tfrac12u_L^2+\mathfrak e(\rho_L)-x/r\right\}\dd x.
 \end{gathered}
 \label{eq:physical-energy-clock}
\end{equation}
Thus $e=0$ corresponds to energy $-4M$ in this normalization, and
$e\ge0$ prescribes the shifted energy $\mathscr E+4M$.
\end{corollary}

\begin{proof}
Proposition~\ref{exlim:energy-conservation} gives
$\mathscr E(t)=\mathscr E(t_0)$ for $t_0,t<0$.
The comparison in Proposition~\ref{exlim:clock-energy} yields
\[
 \lim_{t\uparrow0}\{\mathscr E(t)+4M\}
          =\frac e2 I_2(\delta).
\]
Combining the two identities proves the assertion.
\end{proof}

For comparison, the homogeneous $P=\rho^{4/3}$ Goldreich--Weber
solution with the same $(e,\delta)$ has energy
$\mathscr E_{\rm GW}(e,\delta)=eI_2(\delta)/2$ by
\eqref{exlim:homologous-energy}. Hence the corollary identifies
$\mathscr E+4M=\mathscr E_{\rm GW}(e,\delta)$ for the constructed
orbit. The shift comes from the linear term $-4\rho$ in
$F(\rho)$, whose spatial integral is $-4M$, with the vacuum
normalization fixed in \eqref{eq:internal-energy-normalization}.
Adding a constant to the specific internal energy changes the total
energy by that constant times the conserved mass and leaves the
pressure unchanged.

This identification uses the global comparison in
Lemma~\ref{exlim:global-profile-energy} and
Proposition~\ref{exlim:clock-energy}, including the joining shell
and the nonlinear correction. Those estimates exclude an additional
finite energy defect after the leading kinetic, internal, and
gravitational terms cancel. Neither compact-core convergence nor the
small mass of the physical layer alone gives that conclusion.

Positive-energy examples follow by fixing $e>8M_{\rm Ch}/I_2(0)$
and then taking sufficiently small admissible positive profile
parameters. Smooth dependence gives
\begin{equation}
 \mathscr E(\beta)
 =\frac e2 I_2(-\beta^2/2)-4M_\beta
 \longrightarrow\frac e2 I_2(0)-4M_{\rm Ch}>0.
 \label{eq:positive-energy-members}
\end{equation}
The admissible mass interval and the upper bound for the scale depend
on the fixed $e$; mass and energy are not prescribed independently. The inward
Goldreich--Weber branches also exhibit positive-energy collapse;
here the vacuum surface obeys the exact pressure law.

\paragraph{The next support-radius coefficient.}
For a fixed pair $(\beta,e)$, write the leading term of the
reference scale as
\begin{equation}
 \lambda_{\rm lead}(t)=\left(\frac{3\beta}{2}(-t)\right)^{2/3}.
 \label{eq:leading-power-scale}
\end{equation}
Expanding the integrand in \eqref{eq:main-collapse-clock} gives
\[
 -t=\frac{2}{3\beta}\lambda^{3/2}
       \left(1-\frac{3e}{10\beta^2}\lambda+O(\lambda^2)\right).
\]
Taking the $2/3$ power and inverting yields
\[
 \lambda=\lambda_{\rm lead}
 +\frac{e}{5\beta^2}\lambda_{\rm lead}^2
 +O(\lambda_{\rm lead}^3).
\]
Substitution in \eqref{eq:main-collapse-support} therefore gives
\begin{equation}
 \begin{split}
 R_b(t)={}&\lambda_{\rm lead}(t)
 +\left(\frac{e}{5\beta^2}-\frac1{\varkappa_\delta}\right)
                         \lambda_{\rm lead}(t)^2\\
 &+O\!\left(\lambda_{\rm lead}(t)^3
     (1+|\log\lambda_{\rm lead}(t)|)^{C_{\log}}\right).
 \end{split}
 \label{eq:second-support-coefficient}
\end{equation}
The coefficient of $\lambda_{\rm lead}^2$ contains the contribution
$e/(5\beta^2)$ from the scale equation and the contribution
$-1/\varkappa_\delta$ from the surface layer. Here
$e=2(\mathscr E+4M)/I_2(\delta)$, and the remainder is for the fixed
pair $(\beta,e)$. This expansion asserts no estimate for derivatives
of the remainder.

\begin{corollary}
\label{cor:critical-mass-infimum}
Let $\mathscr S_{22}$ be the set of masses $M>0$ for which the exact
system \eqref{eq:EP}--\eqref{eq:ChP} admits a finite-energy radial
order-22 classical physical-vacuum solution in the sense of
Definition~\ref{def:global-solution} with $N=22$, on some interval $[t_0,T)$
with $T<\infty$, such that its full support radius satisfies
$R_b(t)\to0$ as $t\uparrow T$. Then
\begin{equation}
 \inf\mathscr S_{22}=M_{\rm Ch}.
 \label{eq:critical-mass-infimum}
\end{equation}
Here finite energy uses \eqref{eq:internal-energy-normalization}.
The lower bound follows from the subcritical energy estimate of
Cheng--Cheng--Lin \cite[Section~5]{ChengChengLin}; the upper bound
follows from Corollary~\ref{cor:prescribed-mass}.
\end{corollary}

\begin{proof}
We first exclude collapse for $M<M_{\rm Ch}$ by deriving a positive
lower bound for the support radius. The constructed solutions then
give the reverse inequality for the infimum.

Let a solution in the stated class have mass $M<M_{\rm Ch}$.
Use the normalization \eqref{eq:internal-energy-normalization}.
By Proposition~\ref{exlim:energy-conservation}, the energy
\[
 \mathscr E=\int_0^M
   \left\{\tfrac12u_L^2+\mathfrak e(\rho_L)-\frac{x}{r}\right\}\dd x
\]
is conserved on every compact interval before collapse.
The radial identity \eqref{mat:seed-55}, valid for any radial density,
identifies the magnitude of the gravitational term with the positive
interaction energy $\mathcal W(\rho_E)$ of \eqref{mat:seed-53}.
Combining \eqref{eq:threshold-internal-energy} with the sharp
inequality \eqref{mat:seed-57}, applied to the density extended by
zero, gives
\begin{equation}
 \mathscr E+4M\ge
 c_M\int_{\mathbb R^3}\rho_E(t,|y|)^{4/3}\dd y,
 \qquad
 c_M=3\left[1-\left(\frac{M}{M_{\rm Ch}}\right)^{2/3}\right]>0.
 \label{eq:threshold-coercivity}
\end{equation}
The hypotheses $\rho_E\in L^1\cap L^{4/3}$ follow at each
time before collapse from the classical endpoint regularity and compact
support. The omitted kinetic term is nonnegative. H\"older's
inequality on $B_{R_b(t)}$ gives
\[
 \int_{\mathbb R^3}\rho_E^{4/3}\dd y
 \ge M^{4/3}|B_{R_b(t)}|^{-1/3}
 =\left(\frac3{4\pi}\right)^{1/3}\frac{M^{4/3}}{R_b(t)}.
\]
In particular $\mathscr E+4M>0$, and
\begin{equation}
 R_b(t)\ge
 \left(\frac3{4\pi}\right)^{1/3}
 \frac{c_M M^{4/3}}{\mathscr E+4M}>0
 \qquad(t_0\le t<T).
 \label{eq:threshold-radius-lower-bound}
\end{equation}
The positive lower bound holds throughout the classical interval
of this solution. Hence $\mathscr S_{22}\subset[M_{\rm Ch},\infty)$.
With $e=0$, Corollary~\ref{cor:prescribed-mass} supplies an interval
$(M_{\rm Ch},M_{\rm Ch}+\varepsilon_M(0))\subset\mathscr S_{22}$.
By Proposition~\ref{exlim:energy-conservation}, these solutions belong
to the stated finite-energy class. Consequently,
\[
 M_{\rm Ch}\le\inf\mathscr S_{22}
 \le\inf(M_{\rm Ch},M_{\rm Ch}+\varepsilon_M(0))=M_{\rm Ch},
\]
which proves \eqref{eq:critical-mass-infimum}.
\end{proof}

The corollary leaves open whether $M_{\rm Ch}\in\mathscr S_{22}$
and whether every larger mass belongs to this set. Below the threshold,
the support bound excludes collapse of the whole star but does not
exclude a local loss of classical regularity.

We next define the weighted norms in which we estimate the profile
residual and construct the solution.
\section{Weighted Sobolev spaces}
\label{at:section}

We introduce weighted Sobolev norms for the nonlinear correction and
its source. At the center, a radial displacement is represented by a
Cartesian vector field; its Sobolev norm includes the angular
derivatives. At the vacuum boundary, we use ordinary one-sided
derivatives of the regular factors in
\eqref{eq:global-solution-edge-class}, with no even extension.

We localize to four regions: the center, a fixed interior interval,
a bounded rescaled neighborhood of the vacuum, and an interval
joining the latter to the interior. The last two coordinates are
chosen according to the mass distribution near the boundary.

\subsection{Coordinates and localization}

Write $d=M-x$. Choose nonnegative smooth cutoffs
$\psi_c,\psi_i,\psi_e$ supported near the center, in a fixed compact
interior interval, and in $d<2d_0$, respectively. A cutoff whose
support reaches an endpoint is constant near that endpoint.
For a fixed integer $Q\ge24$, require
\[
 \psi_c^Q+\psi_i^Q+\psi_e^Q\ge c_0>0.
\]
Set
\begin{equation}
 h=(\psi_c^Q+\psi_i^Q+\psi_e^Q)^{1/Q},\qquad
 \theta_\nu=\psi_\nu/h,\qquad \chi_\nu=\theta_\nu^Q,
 \quad \nu=c,i,e.
 \label{at:power-partition}
\end{equation}
Choose these cutoffs so that $\chi_e=1$ for $d\le d_0$.
Split the boundary region once more, using $D=d/\lambda^4$.
Normalize two smooth nonnegative functions of $D$ as above, and
denote the first normalized $Q$-th power by $\zeta$. Choose
$\zeta=1$ for $D\le1/2$ and $\zeta=0$ for $D\ge2$, and put
\[
 \chi_5=\chi_e\zeta,\qquad \chi_8=\chi_e(1-\zeta).
\]
Thus
\[
 \chi_c+\chi_i+\chi_5+\chi_8=1.
\]
Each cutoff is a $Q$-th power of a smooth function with values in
$[0,1]$. This vanishing order is used in the localization estimates.
The cutoff power $Q$ is fixed; the derivative order $q$ below may vary.
The width $d_0$ is fixed. We write it as $d_0^{\rm edge}$ when
necessary to distinguish it from the energy amplitude $d_0$.

At the center, set $z_c=x^{1/3}$ and associate to a displacement $f$
the radial vector field
\[
 U_f(X)=f(|X|^3)\frac{X}{|X|}.
\]
We take its Cartesian Sobolev norm on an enlarged center ball $B$.
On the fixed interior interval $I_i$ we use the mass coordinate $x$.
The two boundary coordinates and their mass measures are
\begin{equation}
 y_n=(d\lambda^{-4})^{1/n},\qquad
 |\dd x|=n\lambda^4y_n^{n-1}\dd y_n,\qquad n=5,8.
 \label{at:edge-coordinates}
\end{equation}
The $y_5$ interval is bounded. The $y_8$ interval is
\[
 2^{-1/8}<y_8<C_e\lambda^{-1/2},\qquad C_e=(2d_0)^{1/8}.
\]
Beyond the matching shell, the radius is given by the core expansion.

\subsection{Weighted norms}

For an interval $I\subset[0,\infty)$ and a scalar or finite tensor
field $F$, define
\[
 \|F\|_{\mathcal B_n^q(I)}^2
   =\sum_{k=0}^q\int_I y^{n-1}|\partial_y^kF|^2\dd y.
\]
We use $B_5^q$ for this norm on the bounded $y_5$ interval and
$C_8^q$ for the corresponding norm on the longer $y_8$ interval.
The latter interval is bounded away from zero. All derivatives
in these norms are ordinary derivatives in the displayed coordinate.

\Needspace{22\baselineskip}
\begin{definition}
\label{at:atlas}
\label{en:input-covering-atlas}
\leavevmode
\begin{enumerate}[label=(\roman*),leftmargin=2em]
\item For a radial displacement, set
\begin{equation}
 \begin{split}
 \|f\|_{X^q(s)}^2={}&
  \|\chi_cU_f\|_{H^q(B;\mathbb R^3)}^2
   +\|\chi_if\|_{H^q(I_i)}^2\\
 &+\lambda^4\|\chi_5f\|_{B_5^q}^2
   +\lambda^4\|\chi_8f\|_{C_8^q}^2 .
 \end{split}
 \label{at:graph-norm}
\end{equation}
The space $X^q(s)$ is the completion in this norm of functions with
smooth equivariant center vectors and smooth one-sided vacuum
representatives satisfying $f_{y_5}(0)=0$.
All derivatives act on the localized functions, including their
cutoffs. For scalar and tensor fields, we use the corresponding
Cartesian components at the center.

\item For a field $g$ without a prescribed vacuum trace, define
$\|g\|_{q,a}$ as the square root of the derivative sum in
\eqref{at:graph-norm}, using the appropriate scalar, vector, or
tensor norm at the center. Thus
\[
 g\in X^q\quad\Longrightarrow\quad
 \|g\|_{q,a}=\|g\|_{X^q}.
\]
Unlike the definition of $X^q$, this norm imposes no approximation
by functions satisfying the vacuum condition. Its finiteness alone
does not imply $g\in X^q$.

\item For a force $g$, define the source norm
\begin{equation}
 \begin{split}
 \mathcal S_m(g)^2={}&
   \|\chi_cU_g\|_{H^m}^2+\|\chi_ig\|_{H^m}^2\\
  &+\lambda^4\|\chi_5\lambda g\|_{B_5^m}^2
   +\lambda^4\|\chi_8\lambda g\|_{C_8^m}^2.
 \end{split}
 \label{at:source-norm}
\end{equation}
\end{enumerate}
\end{definition}

The boundary source norm includes a factor $\lambda$ because the
linearized force has one additional inverse power of the scale there.
Indeed, with $L_R=\lambda^3D\mathcal M[R]$,
\[
 \lambda L_R=-A_nD_n-B_n\frac{\partial_y}{y}+V_n,
 \qquad D_n=\partial_y^2+\frac{n-1}{y}\partial_y,
 \quad n=5,8;
\]
see \eqref{at:edge-hessian}. Under the positivity and normalized
coefficient bounds of Subsection~\ref{at:current-coefficients},
$A_n,B_n,V_n$ are bounded. On the center and interior charts,
the coefficients of $L_R$ itself are bounded; see
\eqref{at:interior-hessian}.

Write $H=L^2((0,M),\dd x)$. The mass measures in
\eqref{at:edge-coordinates} and the positive lower bound for
$\sum_\nu\chi_\nu^2$ give
\[
 \|f\|_{X^0}\sim\|f\|_{L^2((0,M),\dd x)},
\]
uniformly in the scale. Higher derivatives also differentiate the
cutoffs. Their $y_8$ derivatives are uniformly bounded: on the outer
transition $d\sim d_0$,
\[
 \partial_y^k\theta_e(\lambda^4y^8)=O(\lambda^{k/2}),\qquad k\ge1;
\]
the inner transition stays on $D\sim1$. These bounds concern
$y_8$ derivatives and do not give equivalence with higher
$x$-derivative norms.

To compare different leading profiles, use their common normalized
radius $z\in[0,1]$. The mass measure and the uniform profile bound are
\[
 \dd x=4\pi w_\beta(z)^3z^2\dd z,\qquad
 0<c\le\frac{w_\beta(z)}{1-z}\le C
 \quad(0\le z<1),
\]
by \eqref{mat:seed-25}. The center and fourth-root endpoint mass
inverses have the common positive factors of
Subsubsection~\ref{mat:seed-part-5}.

Equivalently, prescribe the cutoffs in the mass fraction
$\widehat x=x/M_\beta$. Then
\[
 x^{1/3}=M_\beta^{1/3}\widehat x^{1/3},\qquad
 y_n=M_\beta^{1/n}
       \bigl((1-\widehat x)\lambda^{-4}\bigr)^{1/n},
 \quad n=5,8.
\]
Since $M_\beta$ is bounded above and away from zero, these changes
of coordinates have bounded factors at every fixed finite order.
These comparisons use the fixed cutoffs on the enlarged coordinate
intervals. At the mass endpoints, the inverse bounds hold in the
coordinates specified above; they do not give unweighted inverse
bounds or higher-order norm equivalence for arbitrary partitions.

Lemma~\ref{rec:transport} controls the time dependence of these
norms. Among the localization cutoffs, only the inner split
$\zeta((M-x)/\lambda^4)$ depends on time. The constants in that
estimate can therefore be fixed before selecting the individual profile.

\subsection{Approximation in the operator domains}
\label{at:notation-dependencies}

The mass coordinate remains fixed during the evolution and both
limits. By Definition~\ref{at:atlas}, approximation in $X^q$
requires smooth functions satisfying the center symmetry and
vacuum condition. The ordinary derivative norm $\|\cdot\|_{q,a}$
and the source norm impose no such approximation condition.

For the elliptic estimates we also require convergence of the
linearized forces. The strong operator domain is the closure of
the same compatible functions in the graph norm
\[
 \|f\|_{X^q}+\mathcal S_{q-2}(L_Rf),\qquad 2\le q\le23;
\]
that is, the approximants must satisfy
\[
 \|f_\ell-f\|_{X^q}
 +\mathcal S_{q-2}(L_Rf_\ell-L_Rf)\longrightarrow0.
\]
Propositions~\ref{at:profile-inverse} and~\ref{at:current-inverse}
prove the corresponding elliptic estimates.

The prepared data satisfy these approximation conditions. At a fixed
positive terminal scale, the inviscid equation gives the spatial
regularity and the compatible approximations needed to apply the
elliptic estimate with its uniform constant; see
Subsection~\ref{exlim:realization}.

After the terminal scale tends to zero, the resulting solution has
the lower-order continuity and highest-order essential bounds in
Proposition~\ref{exlim:terminal-diagonal}. The fixed-scale
approximation does not imply strong continuity at every retained order.

The \hyperref[notation:index]{notation index} lists the coordinates,
scales, radii, and norms, with their definitions and page references.
\section{The finite matched profile and its differentiated residual}
\label{mat:chapter}\label{sec:EOS}

The high-density expansion describes the homologous core on compact
subsets of $x<M$, but fails near vacuum, where the exact pressure has
the $5/3$ law. We construct core and boundary-layer expansions on the
same mass interval and match them in an overlap region. The matching
determines the core amplitudes and additive layer constants, preserves
the positive normalized radius and Jacobian factors, and fixes the
leading correction to the homologous surface radius.

We then estimate the residual of the resulting reference trajectory.
Its core and layer asymptotics are collected in
Corollary~\ref{mat:matched-limits}. Section~\ref{exlim:section}
passes these asymptotics to the Euler--Poisson solution;
Subsection~\ref{exlim:actual-layer} includes the vacuum endpoint on
every fixed bounded rescaled layer interval.

We use the equation of state and enthalpy
\begin{equation}
 \begin{split}
 P(\rho)&=\frac12\{\rho^{1/3}(2\rho^{2/3}-3)
                 \sqrt{1+\rho^{2/3}}
                 +3\operatorname{arsinh}(\rho^{1/3})\},\\
 P'(\rho)&=\frac{4\rho^{2/3}}{3\sqrt{1+\rho^{2/3}}},
 \qquad h(\rho)=4(\sqrt{1+\rho^{2/3}}-1).
 \end{split}
 \label{mat:eos}
\end{equation}
For a particle radius $R$ that increases with enclosed mass, write
\[
 \rho[R]=(4\pi R^2R_x)^{-1},\qquad
 \mathcal{M}[R]=4\pi R^2\partial_xP(\rho[R])+\frac{x}{R^2}.
\]
We work in reversed time $\tau=-t$, with $\dd\tau/\dd s=\lambda^{3/2}$.
The scale satisfies
\begin{equation}
 \lambda_s=b\lambda,\qquad b=\sqrt{\beta^2+e\lambda},
 \qquad b_s=\frac12e\lambda,\qquad e\ge0.
 \label{mat:clock}
\end{equation}
The boundary slope of the leading profile will be denoted by $\varkappa$; the symbol
$\kappa$ is reserved for viscosity.

The leading profiles form a fixed compact family with uniform
bounds for the positive radius and density factors. At a fixed
finite order and a fixed $\beta>0$, the higher coefficients have
bounds that may depend on both choices. After fixing these bounds,
we choose a sufficiently small upper bound for $\lambda$ in $(0,1]$.

The coefficient equations cancel successive powers of the scale in
the residual. To estimate the remainder, we prove bounds for the
coefficient functions and their derivatives. Spatial recovery uses a
mixed source estimate of total order $21$; the energy argument also
requires an $L^2(\dd x)$ estimate through material order $22$.

\medskip
\noindent\textbf{Coordinates and profile coefficients.}\phantomsection\label{mat:dictionary}
The following coordinates distinguish a fixed mass shell from the
shrinking boundary layer and their overlap. They will be derived
from the mass-coordinate map of the leading profile in Subsubsection~\ref{mat:seed-part-5}.
In the core, $x=m(z)$ is enclosed mass, $z$ is the leading radius,
and $q=1-z$ is the distance to the vacuum boundary of the leading
profile. Fixing $x$ fixes both $z$ and $q$.

The rescaled mass deficit is $D=(M-x)/\lambda^4$. Its coordinates are
$y_5=D^{1/5}$ near vacuum and $y_8=D^{1/8}$ on the long chart.
Bounded $D$ describes the boundary layer. For its high-density
expansion we use $Q=D^{1/4}=y_8^2$. The layer formula is used through
$Q=O(\lambda^{-1/2})$; the remainder of the long chart uses the core
radius.

The common overlap coordinates are $X=(M-x)^{1/4}$ and
$Y=\lambda/X=Q^{-1}$. Thus $XY=\lambda$, and the joining shell has
$X,Y\sim\sqrt\lambda$, equivalently $q\sim\sqrt\lambda$.
During coefficient extraction, $L=\log\lambda$ and $\ell=\log X$
are independent polynomial variables; after conversion,
$\log Y=L-\ell$.

The coefficients $\Phi_n$ and $Z_{n-1}$ give the corrections to the
relative core radius and to the inward displacement in the layer.
Both contribute at order $\lambda^n$ to $A/\lambda$.
Time differentiation is always at fixed mass; in the two coordinates,
\[
 \partial_\tau=\lambda^{-3/2}\partial_s,\qquad
 \partial_s=b\lambda\partial_\lambda\big|_z
 =b\left(\lambda\partial_\lambda\big|_D-4D\partial_D\right).
\]
The reference radius is $A$, and $u=A/\lambda$ is its normalized
radius. Later $R$ denotes the solution radius; in local coefficient
calculations it denotes a normalized trial radius only when this is
explicitly stated.

Three symbols have additional local uses. The layer flux coefficient
$A_0(D)=16\pi^2\rho_0^2P'(\rho_0)$ differs from the upper bound
$A_0$ for $\lambda$ in Section~\ref{ct:section}.
Subsections~\ref{mat:cutoff}--\ref{mat:material} use $D$ for a scale
derivative, and Subsection~\ref{mat:endpoint} uses $Q$ for an endpoint
quotient. Each is defined where it is used.

\medskip
\noindent\textbf{Outline of the construction.}\phantomsection\label{mat:reading-guide}
The radius $A$ in Proposition~\ref{mat:matched-family} is constructed
in four stages. The later correction in Section~\ref{loc:section}
changes the terminal data and leaves this reference trajectory fixed.
\begin{enumerate}[label=(\roman*),leftmargin=2.4em]
\item Subsection~\ref{mat:seed} constructs the leading profile, its
mass, and the positive factors in its coordinate maps.
\item Subsection~\ref{mat:connection} solves and matches the core and
layer equations. The leading layer formula \eqref{mat:connection-10}
determines its density and surface correction. At the first two
correction orders, the singular homogeneous core branch is removed;
at later orders it is matched to the layer tail. The additive layer
constant matches the regular branch. Subsubsection~\ref{mat:connection-part-8}
illustrates the first core pole and both logarithmic resonances.
Subsubsections~\ref{mat:connection-part-7}--\ref{mat:connection-part-9}
prove the recurrence and overlap estimate, using the differentiated
remainders of Subsection~\ref{app:profile-finite}.
\item Subsections~\ref{mat:native}--\ref{mat:face} give profile
regularity, bounds for the normalized coefficients, and differentiated
endpoint, core, and layer residual estimates. Their proofs are in
Appendix~\ref{app:profile-estimates}.
\item Subsections~\ref{mat:cutoff}--\ref{mat:exports} estimate the
moving-cutoff error and the material source.
Propositions~\ref{mat:matched-family}--\ref{mat:matched-source}
collect the positivity and source bounds used in the evolution.
\end{enumerate}

In the following summary, $K$ is the truncation order, $a>0$ is
the terminal scale, and $\Lambda\ge0$ is the constant shift in the
source. Write
\[
 L_\lambda=1+|\log\lambda|,\qquad
 F_j=\lambda^{-4}\partial_s^j(\lambda^4F),\qquad
 T=\lambda\partial_\lambda-\frac45y_5\partial_{y_5},
\]
where $F$ is defined in Proposition~\ref{mat:matched-source}.
The constants may depend on the fixed finite family.
Throughout this section, $p_{\log}$ is a finite nonnegative exponent
that may increase finitely many times; $p_{\log,n}$ denotes its
choice at degree $n$.

We shall use the residual bounds in the following derivative ranges.
\begin{itemize}[leftmargin=1.6em]
\item Near vacuum, on $[0,Y_0]$, the residual satisfies
\eqref{mat:endpoint-3}: its derivatives $T^q$, $q\le22$, are bounded
in $C^{28}_{y_5}$ by $C\lambda^{K-2}L_\lambda^{p_{\log}}$.
This gives the physical source estimates \eqref{mat:endpoint-20}
used in exact compatibility.

\item On the four-chart cover, the mixed source bound for $F_j$ in
\eqref{mat:strong-source-export} holds for $j+m\le21$, with size
$C\lambda^{(K-5)/2}L_\lambda^{p_{\log}}$ at $\kappa=0$.
It is used for spatial recovery at $q=m+2$, $j+q\le23$.

\item In $L^2(0,M)$, the material source bound for $F_j$ in
\eqref{mat:material-source-export} holds for $j\le22$, with
size $C\lambda^{(K-7)/2}L_\lambda^{p_{\log}}$ at $\kappa=0$.
It supplies the energy estimates through order $22$ and the kinetic
term at order $23$.

\item In the same respective derivative ranges, the viscosity
contributions to both source bounds have size
$C(\kappa\lambda^{-9/2}+\kappa\Lambda\lambda^{-3/2})$.
The viscosity bounds imposed at the terminal scale make these
$Cc\lambda^{J+1/2}$ for $\lambda\ge a$.
\end{itemize}
The mixed derivative estimates, the pure material derivative estimate,
and the endpoint estimates have separate proofs. Their stated derivative
ranges are retained throughout the evolution.

\Needspace{6\baselineskip}
\subsection{The normalized leading star}\label{mat:seed}\label{sec:GW}

The first-zero normalization places the family of leading profiles on $[0,1]$.
We prove smooth dependence on the force parameter, derive the mass
identity, and obtain uniform bounds for the center and edge coordinate
maps. These bounds are fixed before the matching order is chosen.

\subsubsection{The leading family}\label{mat:seed-part-1}

There is a number $\delta_*>0$ and a smooth family $w_\delta$, defined
for $\delta$ in an open neighborhood of $[-\delta_*,\delta_*]$, such that
\begin{equation}
 w_\delta''+\frac2z w_\delta'+\pi w_\delta^3=-\frac34\delta,\qquad
 w_\delta'(0)=0,\qquad w_\delta(1)=0,\qquad
 w_\delta>0\quad(0\le z<1).
 \label{mat:seed-1}
\end{equation}
The function is smooth at the center as a function of $z^2$. Put
\begin{equation}
 \beta_0=(2\delta_*)^{1/2},\qquad
 w_\beta=w_{-\beta^2/2},\qquad
 a_\delta=w_\delta(0),\qquad \varkappa_\delta=-w_\delta'(1).
 \label{mat:seed-2}
\end{equation}
Then $a_\delta,\varkappa_\delta>0$, uniformly on this compact parameter
interval. Every prescribed finite derivative norm of the leading profile and the
positive factors and inverse coordinates constructed below is uniform in
$0\le\beta\le\beta_0$. The interval is independent of the derivative
order. Each finite derivative constant is chosen before the final
positive value of $\beta$.

The mass is
\begin{equation}
 m_\delta(z)=4\pi\int_0^z w_\delta(s)^3s^2\dd s,\qquad M(\delta)=m_\delta(1).
 \label{mat:seed-3}
\end{equation}
We prove
\begin{equation}
 \delta z+4w_\delta'(z)+\frac{m_\delta(z)}{z^2}=0,\qquad
 M(\delta)+\delta=4\varkappa_\delta,
 \label{mat:seed-4}
\end{equation}
and
\begin{equation}
 M(0)=M_{\rm Ch},\qquad M(\delta)>M_{\rm Ch}\quad(-\delta_*\le\delta<0),
 \qquad \lim_{\delta\uparrow0}M(\delta)=M_{\rm Ch}.
 \label{mat:seed-5}
\end{equation}
The quotient in \eqref{mat:seed-4} has its regular continuous value at $z=0$.

The sharp Newtonian inequality and the identification of the critical
mass are due to Cheng--Cheng--Lin
\cite[Theorem~3.1 and Lemma~5.1]{ChengChengLin}.
We prove the local construction, shooting, and coordinate estimates below.

\subsubsection{The initial-value problem at the center}\label{mat:seed-part-2}

Let $w(z)=V(\eta)$, $\eta=z^2$, with the central value $a$ and
$\delta$ varying independently in bounded neighborhoods. Then
\begin{equation}
 4\eta V''+6V'+\pi V^3=-\frac34\delta,\qquad V(0)=a.
 \label{mat:seed-6}
\end{equation}
Its integral formulation is
\begin{equation}
 \begin{split}
 V(\eta)&=a-\frac12\int_0^\eta\int_0^1
          t^2\left\{\pi V(t^2s)^3+\frac34\delta\right\}\dd t\dd s,\\
 V'(\eta)&=-\frac12\int_0^1
          t^2\left\{\pi V(t^2\eta)^3+\frac34\delta\right\}\dd t.
 \end{split}
 \label{mat:seed-7}
\end{equation}
The outer integral is oriented for $\eta<0$. Choose $B>\sup|a|$
and then $\eta_c>0$ so that
\[
 \frac{\eta_c}{6}\left(\pi B^3+\frac34\sup|\delta|\right)
       \le B-\sup|a|.
\]
The first line of \eqref{mat:seed-7} maps the ball
$\|V\|_\infty\le B$ in $C([-\eta_c,\eta_c])$ into itself,
with Lipschitz constant
\begin{equation}
 \frac{\pi B^2\eta_c}{2}.
 \label{mat:seed-8}
\end{equation}
Take this constant at most $1/2$. The contraction theorem gives
a unique solution, uniformly in the parameters; the integral
equation then gives $V\in C^1$.

\emph{Parameter and spatial derivatives.}
The linearized fixed-point map is
\begin{equation}
 (\mathcal K H)(\eta)
 =-\frac{3\pi}{2}\int_0^\eta\int_0^1
       t^2V(t^2s)^2H(t^2s)\dd t\dd s,\qquad \|\mathcal K\|\le\frac12.
 \label{mat:seed-9}
\end{equation}
Thus $\|(I-\mathcal K)^{-1}\|\le2$. For
$U_\alpha=\partial_{(a,\delta)}^\alpha V$, $\alpha\ne0$, one obtains
\begin{equation}
 \begin{split}
 (I-\mathcal K)U_\alpha
 ={}&\mathbf1_{\alpha=(1,0)}
       -\frac{\eta}{8}\mathbf1_{\alpha=(0,1)}\\
 &-\frac{\pi}{2}\int_0^\eta\int_0^1t^2
       \sum_{\substack{\alpha_1+\alpha_2+\alpha_3=\alpha\\
                       |\alpha_i|<|\alpha|\ (i=1,2,3)}}
       \frac{\alpha!}{\alpha_1!\alpha_2!\alpha_3!}
       \prod_{i=1}^3U_{\alpha_i}(t^2s)\dd t\dd s ,
 \end{split}
 \label{mat:seed-10}
\end{equation}
where $U_0=V$, and the sum is empty at order one. For
$p=(a,\delta)$ and $\Delta_hV=V(p+he_i)-V(p)$,
\[
 \|\Delta_hV\|_\infty\le C|h|,\qquad
 (V+\Delta_hV)^3-V^3-3V^2\Delta_hV
      =3V(\Delta_hV)^2+(\Delta_hV)^3.
\]
Consequently
\[
 \left\|\frac{\Delta_hV}{h}-U_{e_i}\right\|_\infty
 \le2\left\|(I-\mathcal K)\frac{\Delta_hV}{h}
                -(I-\mathcal K)U_{e_i}\right\|_\infty
 \le C|h|.
\]
At higher orders the right side of \eqref{mat:seed-10} contains
only lower parameter derivatives. Its difference quotients converge
by induction, and the same inverse bound proves existence,
continuity, and uniform bounds for every $U_\alpha$.

Differentiating the second line of \eqref{mat:seed-7} gives
\begin{equation}
 \begin{split}
 V^{(j+1)}(\eta)
 ={}&-\frac{\pi}{2}\int_0^1 t^{2j+2}
       \sum_{i_1+i_2+i_3=j}
       \frac{j!}{i_1!i_2!i_3!}
       \prod_{\nu=1}^3V^{(i_\nu)}(t^2\eta)\dd t\\
 &-\frac{\delta}{8}\mathbf1_{j=0}.
 \end{split}
 \label{mat:seed-11}
\end{equation}
This maps $C^j$ bounds to $C^{j+1}$ bounds. Applying parameter
derivatives and using \eqref{mat:seed-10} proves all mixed bounds
on a smaller closed parameter box, with the same $\eta_c$.
Thus $w(z)=V(z^2)$ is smooth and even at zero, and
\[
 6V'(0)+\pi a^3=-\frac34\delta.
\]

\emph{Continuation away from zero.}
For $z\ge z_c>0$, set $Y=(w,w')$. Then
\begin{equation}
 Y'=\left(w',-\frac2z w'-\pi w^3-\frac34\delta\right).
 \label{mat:seed-12}
\end{equation}
On a compact tube about a fixed solution, Gronwall's inequality gives
\[
 |Y(z)-\widetilde Y(z)|
 \le C e^{C(z-z_c)}
 \bigl(|Y(z_c)-\widetilde Y(z_c)|
                         +|\delta-\widetilde\delta|\bigr)
\]
up to the first exit. A smaller parameter neighborhood excludes
such an exit on the chosen compact $z$-interval.
Differentiation in the parameters gives linear equations with bounded
coefficients and sources involving only lower derivatives; Gronwall's
inequality yields their uniform bounds. Spatial derivatives satisfy
\begin{equation}
 \begin{split}
 w^{(j+2)}
 ={}&-2\sum_{\ell=0}^j \binom{j}{\ell}
                  (\partial_z^\ell z^{-1})w^{(j-\ell+1)}\\
 &-\pi\sum_{i_1+i_2+i_3=j}
          \frac{j!}{i_1!i_2!i_3!}\prod_{\nu=1}^3w^{(i_\nu)}
       -\frac34\delta\,\mathbf1_{j=0}.
 \end{split}
 \label{mat:seed-13}
\end{equation}
Uniqueness identifies this solution with the center construction on
the overlap. Equations \eqref{mat:seed-10}--\eqref{mat:seed-13}
therefore give joint smoothness and every fixed finite derivative
bound on the interval used in the shooting argument.

\subsubsection{A finite simple zero for the zero-parameter solution}\label{mat:seed-part-3}

Let $W$ be the solution constructed above for $a=1,\delta=0$.
While $W$ is positive,
\begin{equation}
 -r^2W'(r)=\pi\int_0^r t^2W(t)^3\dd t,\qquad W'(r)<0\quad(r>0).
 \label{mat:seed-14}
\end{equation}
The energy
\begin{equation}
 E(r)=\frac12W'(r)^2+\frac\pi4W(r)^4,\qquad
 E'(r)=-\frac2rW'(r)^2\le0
 \label{mat:seed-15}
\end{equation}
bounds both $W$ and $W'$ on every interval away from the center.
The continuation theorem applied to \eqref{mat:seed-12} therefore
extends $W$ to every finite radius, including beyond any zero.

Suppose $W(r)>0$ for all $r\ge0$. Monotonicity and \eqref{mat:seed-14} give
\begin{equation}
 -W'(r)\ge\frac{\pi r}{3}W(r)^3,\qquad
 (W^{-2})'(r)\ge\frac{2\pi r}{3},\qquad
 W(r)\le\left(1+\frac{\pi r^2}{3}\right)^{-1/2}.
 \label{mat:seed-16}
\end{equation}
Splitting the integral in \eqref{mat:seed-14} at $1$ now yields, for $r\ge2$,
\begin{equation}
 |W'(r)|\le C\frac{1+\log r}{r^2}.
 \label{mat:seed-17}
\end{equation}
Consider the exact Pohozaev quantity
\begin{equation}
 \mathcal P(r)
 =r^3\left(\frac12W'^2+\frac\pi4W^4\right)+\frac12r^2WW'.
 \label{mat:seed-18}
\end{equation}
Using $W''=-2W'/r-\pi W^3$, direct differentiation cancels the
$r^2W'^2$ and $r^3W^3W'$ terms and leaves
\begin{equation}
 \mathcal P'(r)=\frac\pi4r^2W(r)^4>0,\qquad \mathcal P(0)=0.
 \label{mat:seed-19}
\end{equation}
On the other hand, \eqref{mat:seed-16}--\eqref{mat:seed-17} give
\[
 r^3W'^2=O(r^{-1}(1+\log r)^2),\quad
 r^3W^4=O(r^{-1}),\quad r^2WW'=O(r^{-1}(1+\log r)),
\]
so $\mathcal P(r)\to0$. This contradicts \eqref{mat:seed-19}. The first zero $R_0$
is finite. Its derivative is explicitly
\begin{equation}
 W'(R_0)=-\frac{\pi}{R_0^2}\int_0^{R_0}t^2W(t)^3\dd t<0.
 \label{mat:seed-20}
\end{equation}
Thus $R_0$ is a simple first zero.

\subsubsection{Normalization of the first zero and conversion of parameters}\label{mat:seed-part-4}

Denote by $w_{a,\delta}$ the regular-center solution of \eqref{mat:seed-1} with central
value $a$. The zero-parameter scaling is
\begin{equation}
 w_{a,0}(z)=aW(az).
 \label{mat:seed-21}
\end{equation}
At $a_0=R_0$ it has its first zero at one. Subsubsections~\ref{mat:seed-part-2}--\ref{mat:seed-part-3} give a common
interval $[0,1+\varepsilon]$ on which $w_{a,\delta}$ is jointly smooth
for $(a,\delta)$ near $(a_0,0)$. For
$G(a,\delta)=w_{a,\delta}(1)$,
\begin{equation}
 G(a_0,0)=0,\qquad
 G_a(a_0,0)=W(R_0)+R_0W'(R_0)=R_0W'(R_0)\ne0.
 \label{mat:seed-22}
\end{equation}
The implicit-function theorem therefore gives a smooth $a(\delta)$
near zero, with $a(0)=a_0$ and $G(a(\delta),\delta)=0$.
The only divisor in its successive derivatives is $G_a$, bounded
away from zero after shrinking the neighborhood. For example,
\begin{equation}
 a'=-G_\delta/G_a,\qquad
 a''=-\{G_{aa}(a')^2+2G_{a\delta}a'+G_{\delta\delta}\}/G_a.
 \label{mat:seed-23}
\end{equation}
The higher derivatives are obtained by differentiating
$G(a(\delta),\delta)=0$; the top derivative occurs only as
$G_a a^{(j)}$, with all other terms finite products of lower
derivatives of $a$ and derivatives of $G$. Thus the branch is smooth
and has every finite derivative bound on a compact subinterval.

To show that $z=1$ is the first zero, choose $\varepsilon_e>0$ for which
$w_0'\le-\varkappa_0/2<0$ on $[1-\varepsilon_e,1+\varepsilon_e]$.
Joint $C^1$ continuity preserves a strict negative bound there for
nearby $\delta$. Then, for $1-\varepsilon_e\le z<1$,
\begin{equation}
 w_\delta(z)=-\int_z^1w_\delta'(s)\dd s>0.
 \label{mat:seed-24}
\end{equation}
On $[0,1-\varepsilon_e]$, the minimum of $w_0$ is positive and
joint continuity preserves a positive lower bound. Also $a(\delta)>0$.
Hence the zero at one is the first zero, and
$\varkappa_\delta=-w_\delta'(1)>0$.

Choose a closed interval $[-\delta_*,\delta_*]$ strictly inside this
neighborhood. All subsequent fixed positive lower bounds can be
chosen after reducing $\delta_*$ a finite number of times. None of
these reductions depends on the final positive $\beta$ or on a
derivative order. On the resulting compact interval, for some
$0<c<C<\infty$,
\begin{equation}
 c\le a_\delta,\varkappa_\delta\le C,\qquad
 c(1-z)\le w_\delta(z)\le C(1-z)\quad(0\le z\le1).
 \label{mat:seed-25}
\end{equation}
Indeed, the quotient at $z=1$ has the continuous positive value
$\varkappa_\delta$, and elsewhere positivity follows from \eqref{mat:seed-24} and the
compact-interior lower bound.

If $v_d(0)=1$, $v_d'(0)=0$ and
$v_d''+2r^{-1}v_d'+\pi v_d^3=-3d/4$, its simple first zero
$R(d)$ is smooth near zero by \eqref{mat:seed-20}. The normalized function satisfies
\begin{equation}
 w(z)=R(d)v_d(R(d)z),\qquad
 w''+\frac2zw'+\pi w^3=-\frac34R(d)^3d.
 \label{mat:seed-26}
\end{equation}
Thus the normalized forcing is
\begin{equation}
 \delta=R(d)^3d,\qquad
 \left.\frac{\dd\delta}{\dd d}\right|_{d=0}=R_0^3>0.
 \label{mat:seed-27}
\end{equation}
This is an invertible local reparameterization. The branch in
\eqref{mat:seed-22} is parametrized by the normalized forcing
$\delta$, which is the parameter set equal to $-\beta^2/2$ in
\eqref{mat:seed-2}. The same change of scale preserves the mass:
\begin{equation}
 4\pi\int_0^1[R(d)v_d(R(d)z)]^3z^2\dd z
 =4\pi\int_0^{R(d)}v_d(r)^3r^2\dd r.
 \label{mat:seed-28}
\end{equation}

\subsubsection{Positive factors and inverse mass-coordinate maps}\label{mat:seed-part-5}

All bounds below include any fixed finite number of parameter
derivatives on the compact $\delta$-interval, and hence on the
closed $\beta$-interval under $\delta=-\beta^2/2$.

\paragraph{Positive factors of the leading profile}

At the center, $w_\delta(z)=V_\delta(z^2)$ by
\eqref{mat:seed-7}. At the edge, with $q=1-z$,
\begin{equation}
 \omega_\delta(q)
 =-\int_0^1w_\delta'(1-tq)\dd t,\qquad
 w_\delta(1-q)=q\omega_\delta(q),\qquad
 \omega_\delta(0)=\varkappa_\delta.
 \label{mat:seed-29}
\end{equation}
Differentiation gives
\begin{equation}
 \partial_q^j\omega_\delta(q)
 =-\int_0^1(-t)^j w_\delta^{(j+1)}(1-tq)\dd t.
 \label{mat:seed-30}
\end{equation}
Hence \eqref{mat:seed-25} implies
\[
 0<c\le\omega_\delta\le C,\qquad
 \|\partial_q^j\omega_\delta\|_\infty
 \le\frac1{j+1}\|w_\delta^{(j+1)}\|_\infty.
\]
The flux equation similarly gives
\begin{equation}
 -\frac{w_\delta'(z)}z
 =\int_0^1t^2\left\{\pi w_\delta(tz)^3+\frac34\delta\right\}\dd t
 =:H_\delta(z).
 \label{mat:seed-31}
\end{equation}
This factor is smooth and even at zero. Since $H_0>0$ on
$[0,1]$, one reduction of the parameter interval gives
$c\le H_\delta\le C$, and therefore $w_\delta'<0$ for $z>0$.
This choice uses only the leading $C^1$ lower bound.

For a positive smooth factor $c\le U\le C$, the chain rule reads
\begin{equation}
 \partial^j(U^\alpha)
 =\sum_{\mathcal P\in\Pi_j}
   (\alpha)_{|\mathcal P|}\,U^{\alpha-|\mathcal P|}
        \prod_{B\in\mathcal P}\partial^{|B|}U,
 \label{mat:seed-32}
\end{equation}
where $\Pi_j$ consists of partitions of the $j$ labeled derivatives,
and $(\alpha)_k=\alpha(\alpha-1)\cdots(\alpha-k+1)$.
The same formula with the derivatives of $\log$ proves the logarithmic
bound. Thus reciprocals, fixed powers, and logarithms preserve all
fixed finite derivative bounds, also for mixed derivatives.

\paragraph{Center mass and the inverse radius map}

Factoring the mass integral gives
\begin{equation}
 m_\delta(z)=z^3B_\delta(z^2),\qquad
 B_\delta(\eta)=4\pi\int_0^1t^2V_\delta(t^2\eta)^3\dd t,\qquad
 B_\delta(0)=\frac{4\pi}{3}a_\delta^3.
 \label{mat:seed-33}
\end{equation}
The factor $B_\delta$ is smooth and uniformly positive near zero.
For $\xi=x^{1/3}$,
\begin{equation}
 \xi=zB_\delta(z^2)^{1/3},\qquad
 \sigma:=\xi^2=\eta B_\delta(\eta)^{2/3}=:F_\delta(\eta),\qquad
 F_\delta'(0)=B_\delta(0)^{2/3}>0.
 \label{mat:seed-34}
\end{equation}
Choose a fixed interval on which $F_\delta'\ge c>0$ and let
$G_\delta=F_\delta^{-1}$ on a common image interval. Then
\begin{equation}
 z=\xi b_\delta(\xi^2),\qquad
 b_\delta(\sigma)
 =\left(\int_0^1G_\delta'(t\sigma)\dd t\right)^{1/2},\qquad
 b_\delta(0)=B_\delta(0)^{-1/3}.
 \label{mat:seed-35}
\end{equation}
The inverse bounds follow from
\begin{equation}
 G'=1/F'(G),
 \qquad
 G^{(j)}
 =-\frac1{F'(G)}
   \sum_{\substack{\mathcal P\in\Pi_j\\|\mathcal P|\ge2}}
       F^{(|\mathcal P|)}(G)
       \prod_{B\in\mathcal P}G^{(|B|)},\quad j\ge2.
 \label{mat:seed-36}
\end{equation}
Indeed every derivative of $G$ on the right has order at most
$j-1$, and
\[
 \partial_\delta G=-\frac{(\partial_\delta F)(G)}{F'(G)}.
\]
Induction gives every mixed derivative bound. Equations
\eqref{mat:seed-32}--\eqref{mat:seed-36} then show that
$z/\xi=b_\delta(\xi^2)$ and
$w_\delta(z(\xi))=V_\delta(G_\delta(\xi^2))$ are smooth,
even, and uniformly positive near the center.

\paragraph{Endpoint mass, its fourth-root inverse, and finite remainders}

At the boundary,
\begin{equation}
 \begin{split}
 d:=M(\delta)-m_\delta(1-q)
 &=4\pi\int_0^q s^3\omega_\delta(s)^3(1-s)^2\dd s
   =q^4A_\delta(q),\\
 A_\delta(q)
 &=4\pi\int_0^1t^3\omega_\delta(tq)^3(1-tq)^2\dd t,\qquad
 A_\delta(0)=\pi\varkappa_\delta^3.
 \end{split}
 \label{mat:seed-37}
\end{equation}
Choose $q_e<1/3$ so that $A_\delta$ is uniformly positive on
$[0,3q_e]$. The coordinate $X=d^{1/4}$ satisfies
\begin{equation}
 X=qA_\delta(q)^{1/4},\qquad
 \left.\frac{\dd X}{\dd q}\right|_{q=0}=(\pi\varkappa_\delta^3)^{1/4}>0.
 \label{mat:seed-38}
\end{equation}
After decreasing $q_e$, $X_q\ge c>0$ on this interval. Its inverse is
\begin{equation}
 q=Q_\delta(X)=Xc_\delta(X),\qquad
 c_\delta(X)=\int_0^1Q_\delta'(tX)\dd t,\qquad
 c_\delta(0)=(\pi\varkappa_\delta^3)^{-1/4}.
 \label{mat:seed-39}
\end{equation}
The inverse formula \eqref{mat:seed-36} and the positive-factor
formula \eqref{mat:seed-32} give uniform bounds for
$Q_\delta,c_\delta,c_\delta^{-1}$ and all their fixed finite derivatives.
For integers $p,r\ge0$,
\begin{equation}
 \begin{split}
 c_\delta(X)&=\sum_{j=0}^p\frac{\partial_X^jc_\delta(0)}{j!}X^j
                       +X^{p+1}R_{\delta,p}(X),\\
 R_{\delta,p}(X)
 &=\frac1{p!}\int_0^1(1-t)^p
               \partial_X^{p+1}c_\delta(tX)\dd t .
 \end{split}
 \label{mat:seed-40}
\end{equation}
In particular,
\[
 \|R_{\delta,p}\|_{C^r}\le C_{p,r}
                         \|c_\delta\|_{C^{p+1+r}},\qquad
 (X\partial_X)^j
       =\sum_{\ell=0}^j c_{j\ell}X^\ell\partial_X^\ell.
\]
The product rule therefore gives
\begin{equation}
 \left|(X\partial_X)^j\{X^{p+1}R_{\delta,p}(X)\}\right|
 \le C_{p,r}X^{p+1}\quad(0\le j\le r).
 \label{mat:seed-41}
\end{equation}
The same estimates apply to the smooth term in
\begin{equation}
 \log q=\log X+\log c_\delta(X)
 \label{mat:seed-42}
\end{equation}
Moreover,
\begin{equation}
 q\partial_q
 =\left(1+\frac{qA_\delta'(q)}{4A_\delta(q)}\right)X\partial_X.
 \label{mat:seed-43}
\end{equation}
The coefficient equals one at zero. It and its reciprocal have
uniform finite derivative bounds, so either operator gives the same
remainder estimates. The normalized coefficients
\begin{equation}
 \frac{\omega_\delta(q)^4(1-q)^4}{\varkappa_\delta^4},
 \qquad
 \frac{\omega_\delta(q)^3(1-q)^4}{\varkappa_\delta^3},
 \label{mat:seed-44}
\end{equation}
are likewise positive smooth factors, equal to one at zero.
Their nonconstant terms are retained when comparing higher
coefficients.

\paragraph{Covering geometry and derivative bounds}

Choose $z_c,q_e>0$ with $3z_c+3q_e<1$. At $\delta=0$, the
images of $2z_c,2q_e$ lie strictly inside those of $3z_c,3q_e$.
After one reduction of the parameter interval, the minimum of each
larger image exceeds the maximum of the corresponding smaller images.
The center and edge inverses therefore exist on common intervals
containing all the cutoffs below.

Take fixed smooth cutoffs $\widehat\chi_c(z)$ and
$\widehat\chi_e(1-z)$, equal to one for $z\le z_c$ and
$q\le q_e$, and supported where $z<2z_c$ and $q<2q_e$.
Their supports are disjoint. With
$\widehat\chi_i=1-\widehat\chi_c-\widehat\chi_e$,
\begin{equation}
 \widehat\chi_c+\widehat\chi_i+\widehat\chi_e=1\quad\hbox{on }[0,1].
 \label{mat:seed-45}
\end{equation}
On $\operatorname{supp}\widehat\chi_i\subset[z_c,1-q_e]$,
\begin{equation}
 m_\delta'(z)=4\pi z^2w_\delta(z)^3\ge c>0.
 \label{mat:seed-46}
\end{equation}
Formula \eqref{mat:seed-36} bounds the interior inverse on its
mass image, equivalently after pullback to this fixed $z$-interval.
Enlarging the interval gives a common image for nearby parameters.
These unweighted inverse bounds are used only away from the two
singular mass endpoints.

The enlarged inverse intervals also give uniform derivative bounds
for the pulled-back hatted cutoffs. The graph and source norms use
instead the fixed-mass cutoffs of Definition~\ref{at:atlas};
no high derivative comparison between the two partitions is needed.

The quotient and inverse formulas give, for each fixed $L$,
\[
 \begin{aligned}
 \|\omega_\delta\|_{C^L}
     &\le C_L\|w_\delta\|_{C^{L+1}},\\
 \|q/X\|_{C^L}
     &\le C_L\bigl(1+\|w_\delta\|_{C^{L+2}}\bigr)^{N_L},\\
 \|b_\delta\|_{C^L}
     &\le C_L\bigl(1+\|V_\delta\|_{C^{L+1}}\bigr)^{N_L}.
 \end{aligned}
\]
Here $N_L$ is a finite integer; the constants also use the fixed
positive lower bounds. Thus leading-profile bounds through
$2S+62$ control the leading factors through $2S+60$.
For a degree-$p$ Taylor remainder with $r$ derivatives, the required
order is $p+1+r$ in \eqref{mat:seed-40}. All these constants
are uniform on the fixed compact parameter interval; higher matching
coefficients may depend on the selected $\beta>0$.

Finally, define the leading acoustic coordinate by
\begin{equation}
 \ell_\delta(z)=\frac{\sqrt3}{2}\int_0^z w_\delta(s)^{-1/2}\dd s.
 \label{mat:seed-47}
\end{equation}
Its length is finite and uniformly bounded above and below by
\eqref{mat:seed-25}. At the edge and center, respectively,
\begin{equation}
 \ell_\delta(1)-\ell_\delta(1-q)
 =\sqrt3\,q^{1/2}\int_0^1\omega_\delta(qt^2)^{-1/2}\dd t,
 \label{mat:seed-48}
\end{equation}
\begin{equation}
 \ell_\delta(z)
 =\frac{\sqrt3}{2}z\int_0^1V_\delta(z^2t^2)^{-1/2}\dd t.
 \label{mat:seed-49}
\end{equation}
The quotients by $\sqrt q$ and $z$ are smooth positive functions of
$q$ and $z^2$, with all fixed finite bounds. Parameter derivatives
of the total length pass under the integrable $q^{-1/2}$ weight
by \eqref{mat:seed-29}--\eqref{mat:seed-32}. Estimates for
the perturbed operators will be proved later.

\subsubsection{Exact integrated mass identity}\label{mat:seed-part-6}

Multiplying \eqref{mat:seed-1} by $z^2$ and integrating from zero gives
\begin{equation}
 z^2w_\delta'(z)
 =-\pi\int_0^z w_\delta(s)^3s^2\dd s-\frac{\delta}{4}z^3
 =-\frac14m_\delta(z)-\frac{\delta}{4}z^3.
 \label{mat:seed-50}
\end{equation}
This proves \eqref{mat:seed-4} for $z>0$; \eqref{mat:seed-33} extends the expression smoothly to
zero. At $z=1$, it gives
\begin{equation}
 M(\delta)+\delta=-4w_\delta'(1)=4\varkappa_\delta.
 \label{mat:seed-51}
\end{equation}
This is the mass identity used in both the leading core and layer
equations. It follows from the profile equation and requires no further
shooting condition.

\subsubsection{Critical mass and strict supercriticality}\label{mat:seed-part-7}

Let
\begin{equation}
 \rho_\delta(y)=
 \begin{cases}w_\delta(|y|)^3,&|y|<1,\\0,&|y|\ge1,\end{cases}
 \qquad
 I_{4/3}=4\pi\int_0^1w_\delta^4z^2\dd z,\qquad
 I_2=4\pi\int_0^1w_\delta^3z^4\dd z.
 \label{mat:seed-52}
\end{equation}
Then $\rho_\delta\ge0$ is bounded and compactly supported,
$\rho_\delta\in L^1\cap L^{4/3}(\mathbb R^3)$, and
$I_{4/3},I_2>0$. Define
\begin{equation}
 \mathcal W(\rho)=\frac12\iint_{\mathbb R^3\times\mathbb R^3}
                    \frac{\rho(y)\rho(\widetilde y)}{|y-\widetilde y|}
                     \dd y\dd\widetilde y.
 \label{mat:seed-53}
\end{equation}
Angular integration gives
\begin{equation}
 2\pi\int_{-1}^{1}(r^2+s^2-2rsu)^{-1/2}\dd u
 =\frac{4\pi}{\max(r,s)}.
 \label{mat:seed-54}
\end{equation}
Splitting the radial integral into $r>s$ and $s>r$ and using symmetry yields
\begin{equation}
 \mathcal W(\rho_\delta)
 =\int_0^1\frac{m_\delta(z)}z\dd m_\delta(z)
 =4\pi\int_0^1m_\delta(z)w_\delta(z)^3z\dd z.
 \label{mat:seed-55}
\end{equation}
All these integrals are finite.

We first check the normalization in \cite{ChengChengLin}. Our pressure is
\begin{equation}
 P(\rho)=\frac12\left[
  \rho^{1/3}(2\rho^{2/3}-3)\sqrt{1+\rho^{2/3}}
  +3\operatorname{arsinh}(\rho^{1/3})\right].
 \label{mat:seed-56}
\end{equation}
Their notation is $P(\rho)=A f((\rho/B)^{1/3})$, where
\[
 f(t)=t\sqrt{1+t^2}(2t^2-3)+3\log(t+\sqrt{1+t^2}).
\]
Thus
\[
 A=\frac12,\qquad B=1,\qquad K_{\rm pol}=2AB^{-4/3}=1.
\]
At $\delta=0$, \eqref{mat:seed-50} is the hydrostatic equation
for $P=\rho^{4/3}$ and $\rho=w_0^3$. Its mass is unchanged by
the central-value scaling \eqref{mat:seed-21}, as computed in
\eqref{mat:seed-28}. Hence \cite[Lemma~5.1]{ChengChengLin}
identifies this Lane--Emden mass with $M_{\rm Ch}$.

The sharp constant in \cite[Theorem~3.1]{ChengChengLin} for
the double integral without the factor $1/2$ is
$6K_{\rm pol}/M_{K_{\rm pol}}^{2/3}$. In our normalization,
for nonnegative $\rho\in L^1\cap L^{4/3}$,
\begin{equation}
 \mathcal W(\rho)\le
 3\left(\frac{\int\rho}{M_{\rm Ch}}\right)^{2/3}\int\rho^{4/3}.
 \label{mat:seed-57}
\end{equation}
To apply this to $\rho_\delta$, multiply \eqref{mat:seed-4}
by $4\pi z^3w_\delta^3$ and integrate. Since
$[z^3w_\delta^4]_0^1=0$,
\begin{equation}
 16\pi\int_0^1z^3w_\delta^3w_\delta'\dd z
 =4\pi\int_0^1z^3(w_\delta^4)'\dd z
 =-12\pi\int_0^1z^2w_\delta^4\dd z
 =-3I_{4/3}.
 \label{mat:seed-58}
\end{equation}
The other two terms are $\mathcal W(\rho_\delta)$ by
\eqref{mat:seed-55} and $\delta I_2$. Therefore
\begin{equation}
 3I_{4/3}-\mathcal W(\rho_\delta)=\delta I_2.
 \label{mat:seed-59}
\end{equation}
Combining \eqref{mat:seed-57} and \eqref{mat:seed-59} yields
\begin{equation}
 M(\delta)\ge M_{\rm Ch}
       \left(1-\frac{\delta I_2}{3I_{4/3}}\right)^{3/2}
       >M_{\rm Ch}\qquad(\delta<0).
 \label{mat:seed-60}
\end{equation}
Smooth parameter dependence gives
\[
 \lim_{\delta\uparrow0}M(\delta)=M(0)=M_{\rm Ch}.
\]
For each permitted $\beta_1>0$, continuity of
$M_\beta=M(-\beta^2/2)$ and the intermediate value theorem give
\[
 (M_{\rm Ch},M_{\beta_1})
       \subset\{M_\beta:0<\beta<\beta_1\}.
\]
No monotonicity of the mass branch is needed here.

\begin{lemma}
\label{mat:mass-tangent}
Let $I_2(0)=4\pi\int_0^1w_0^3z^4\dd z$. The smooth branch of leading profiles,
normalized to have their first zero at one, satisfies
\begin{equation}
 M'(0)=-\frac{I_2(0)}{M_{\rm Ch}}<0,
 \qquad
 M_\beta=M_{\rm Ch}
       +\frac{I_2(0)}{2M_{\rm Ch}}\beta^2+O(\beta^4)
       \quad(\beta\downarrow0),
 \label{mat:mass-tangent-expansion}
\end{equation}
where $M'(0)$ is the derivative with respect to $\delta$.
Moreover
\begin{equation}
 \frac{\dd M_\beta}{\dd\beta}
       =\frac{I_2(0)}{M_{\rm Ch}}\beta+O(\beta^3)>0
 \label{mat:mass-beta-derivative}
\end{equation}
for all sufficiently small positive $\beta$.
The monotonicity concerns the leading profiles only. No continuous
selection of the corresponding nonlinear collapse solutions is asserted.
\end{lemma}

\begin{proof}
All derivatives used below exist by the smooth parameter construction
on the fixed interval $[0,1]$. Put $w=w_0$,
$\varkappa_0=-w'(1)=M_{\rm Ch}/4$,
$v=w+zw'$, and
$\dot w=\left.\partial_\delta w_\delta\right|_{\delta=0}$.
For
\[
 Lf=f''+2z^{-1}f'+3\pi w^2f,
\]
differentiation of the equation for the leading profile and of its central-value scaling
at $\delta=0$ gives
\[
 L\dot w=-\tfrac34,\qquad Lv=0,\qquad
 \dot w(1)=0,\qquad v(1)=-\varkappa_0.
\]
The center derivatives vanish, and the Green identity with weight
$z^2$ has no center boundary term. Therefore
\[
 -\varkappa_0\dot w'(1)
   =-\frac34\int_0^1z^2(w+zw')\dd z
   =\frac32\int_0^1z^2w\dd z.
\]
Differentiating \eqref{mat:seed-51} now yields
\[
 M'(0)=-4\dot w'(1)-1
       =\frac6{\varkappa_0}\int_0^1z^2w\dd z-1.
\]
On the other hand, the leading profile equation at zero parameter and two integrations
by parts give
\[
 \begin{split}
 I_2(0)
 &=-4\int_0^1(z^4w''+2z^3w')\dd z\\
 &=4\varkappa_0-24\int_0^1z^2w\dd z.
 \end{split}
\]
Since $4\varkappa_0=M_{\rm Ch}$, these identities imply the first
formula in \eqref{mat:mass-tangent-expansion}. Taylor expansion of
the smooth function $M(\delta)$ at zero, with $\delta=-\beta^2/2$,
gives the second. Finally the chain rule and
$M'(\delta)=M'(0)+O(|\delta|)$ prove
\eqref{mat:mass-beta-derivative}; its leading coefficient is positive.
The monotonicity interval may be smaller than the full interval
permitted by the nonlinear construction.
\end{proof}

The normalized leading profile, its mass, and both endpoint coordinate inverses
now have uniform bounds at every prescribed finite order on the
fixed parameter interval.

\Needspace{6\baselineskip}
\subsection{Construction of the finite matched family}\label{mat:connection}

We connect the $4/3$ leading profile to the $5/3$ vacuum behavior
of the exact pressure by solving the core and layer coefficient
equations. Their finite expansions are compared through the recurrence
\eqref{mat:connection-coefficient-recursion} in the two overlap indices.
Subsubsection~\ref{mat:connection-part-8}
(p.~\pageref{mat:worked-connection}) computes the first logarithmic
coefficient and the first singular matching condition.

We must fix a center amplitude and an additive layer constant.
A center-regular core
solution can contain an endpoint $q^{-3}$ branch; its center amplitude
is chosen to match the singular tail of a previously constructed
layer coefficient. The zero-flux layer solution retains an additive
constant, which matches the regular core branch.

Inverting the source at either indicial root can produce logarithms.
We therefore allow finite logarithmic polynomials. Their nonconstant
logarithmic terms are determined by the source, while the two
normalizations above remain free. The different degrees of the two
branches allow these constants to be determined successively.

\subsubsection{The finite family and the parameter choices}\label{mat:connection-part-1}

Fix $0<\beta_{\rm match}\le\beta_0$ using only the positive
lower bounds for the compact family of leading profiles. This interval
is independent of the finite order $K\ge86$. Put $S=K+1$ and fix
$e\ge0$. For $0<\beta\le\beta_{\rm match}$, the leading profile
has its first zero at one and satisfies
\begin{equation}
 w''+2z^{-1}w'+\pi w^3=-3\delta/4,\quad
 w'(0)=0,\quad w(1)=0,\quad
 \delta=-\beta^2/2,\quad \varkappa=-w'(1)>0,
 \label{mat:connection-1}
\end{equation}
and
\begin{equation}
 m(z)=4\pi\int_0^z w^3s^2\dd s,\quad M=m(1),\quad
 \delta z+4w'(z)+m(z)/z^2=0,\quad M+\delta=4\varkappa.
 \label{mat:connection-2}
\end{equation}
The positive lower bounds and the required finite derivative bounds
for the leading profiles and leading layers are fixed uniformly on
$[0,\beta_{\rm match}]$ before selecting
a positive $\beta$. For that choice, the higher coefficient bounds
may depend on $(K,e,\beta)$, but not on $\lambda$ or the lower
terminal scale. The inverse connection estimates are used only at
$\beta>0$.

After completing stage $n$, let $d_n$ be the larger of the degrees
of $\Phi_n$ and $Z_{n-1}$ as polynomials in $L$, and put
\begin{equation}
 \Lambda_n=\{0,1,\ldots,d_n\},\qquad \Lambda_1=\{0\}.
 \label{mat:finite-log-index}
\end{equation}
Coefficients outside the indicated logarithmic range are set to zero.
At stage $n$, choose the core degree large enough to include both the
source and the singular coefficient already prescribed by the layer.
After solving the core equation, match its regular coefficient by the
additive layer polynomial; this fixes $d_n$.

Products add logarithmic degrees. Integration at an indicial root
increases the degree in the endpoint logarithm by at most one, while
solving the core equations in decreasing logarithmic degree does not
increase the degree in $L$. The polynomial inverses below therefore
give a finite $d_n$ at every stage.

There are finite logarithmic polynomials $\Phi_n(L,z)$, $2\le n\le S$,
and $Z_k(L,D)$, $0\le k\le K$, for which the finite radii
\begin{samepage}
\begin{equation}
 u_c=z\left(1+\sum_{n=2}^{K}\lambda^n
                         \Phi_n(\log\lambda,z)\right),\qquad
 u_v=1-\sum_{k=0}^{K}\lambda^{k+1}
                         Z_k(\log\lambda,D),\quad
 D=\frac{M-m(z)}{\lambda^4},
 \label{mat:connection-3}
\end{equation}
\equationalias{in:finite-core}{mat:connection-3}
\end{samepage}
obey, on every fixed enlarged shell
$c_1\sqrt\lambda\le q=1-z\le c_2\sqrt\lambda$,
\begin{equation}
 \sum_{a+b\le24}
 \left|(\lambda\partial_\lambda|_z)^a
                   (q\partial_q)^b(u_c-u_v)\right|
 \le C\lambda^{(K+1)/2}(1+|\log\lambda|)^{p_{\log}}.
 \label{mat:connection-4}
\end{equation}
The physical core and layer radii are $\lambda u_c$ and $\lambda u_v$,
respectively, and each
$\Phi_n(L,z)=\sum_{\ell=0}^{d_n}L^\ell\phi_{n\ell}(z)$.
Here $\Phi_S$ is an auxiliary coefficient used to fix the constant of
$Z_K$; it is not inserted into $u_c$. Its omission starts beyond
the total degrees needed in \eqref{mat:connection-4}. The core and layer equations
cancel through scale degree $K$.
The source estimates below use these radius derivatives through order 24.

The leading profile in \eqref{mat:connection-1} and its uniform bounds are supplied by
Subsection~\ref{mat:seed}. We now derive the equations for the higher
coefficients, solve them with the required endpoint conditions, and
match the two expansions.

\subsubsection{Equations for the core and layer coefficients}\label{mat:connection-part-2}

To derive the core equations, write the normalized radius as $z\psi$
and keep $L$ independent of $\lambda$ during coefficient extraction. Set
\begin{equation}
 T=\lambda\partial_\lambda+\partial_L,\quad
 J=\psi+z\psi_z,\quad \theta=w^3/(\psi^2J).
 \label{mat:connection-5}
\end{equation}
The exact reduced core residual is
\begin{equation}
 \begin{aligned}
 F_c[\psi]
 ={}&\beta^2\big((1+T)^2-\tfrac32(1+T)\big)(z\psi)\\
 &{}+e\lambda T(1+T)(z\psi)\\
 &{}+\frac{\psi^2}{w^3}
       \partial_z\{\lambda^4P(\lambda^{-3}\theta)\}
       +\frac{m(z)}{z^2\psi^2}.
 \end{aligned}                                                \label{mat:connection-6}
\end{equation}
Its Taylor coefficient at a new relative-radius degree $n\ge2$ is
\begin{equation}
 z\left\{B_n+\beta^2(2n+\tfrac12)\partial_L+
                        \beta^2\partial_L^2\right\}\Phi_n,
 \quad
 B_n=-\frac4{3w^3z^4}\partial_z(w^4z^4\partial_z)
                      +(2n^2+n-3)|\delta|.
 \label{mat:connection-7}
\end{equation}
The source is minus $z^{-1}$ times the coefficient of $\lambda^n$
in \eqref{mat:connection-6} with only lower scale degrees inserted. Therefore the
equation for $\Phi_n$ is a finite polynomial equation in $L$.
Write this known source as $\sum_\ell L^\ell\mathfrak F_{n\ell}$
and $\Phi_n=\sum_\ell L^\ell\phi_{n\ell}$. After fixing the
finite logarithmic degree, the downward equations are exactly
\begin{equation}
 B_n\phi_{n\ell}=\mathfrak F_{n\ell}
 -\beta^2(2n+\tfrac12)(\ell+1)\phi_{n,\ell+1}
 -\beta^2(\ell+1)(\ell+2)\phi_{n,\ell+2}=:F_{n\ell},
 \label{mat:core-log-recursion}
\end{equation}
with out-of-range coefficients zero. The last two terms involve only
already solved higher logarithmic coefficients at that degree.

For the layer put
\begin{equation}
 R=1-\lambda Z,\quad
 \rho=(4\pi R^2Z_D)^{-1},\quad
 T_v=\lambda\partial_\lambda+\partial_L-4D\partial_D.
 \label{mat:connection-8}
\end{equation}
The exact reduced layer residual is
\begin{equation}
 \begin{aligned}
 F_v[Z]
 ={}&\beta^2\big((1+T_v)^2-\tfrac32(1+T_v)\big)R\\
 &{}+e\lambda T_v(1+T_v)R\\
 &{}-4\pi R^2\partial_DP(\rho)
       +(M-\lambda^4D)R^{-2}.
 \end{aligned}                                                \label{mat:connection-9}
\end{equation}
At the leading degree choose
\begin{equation}
 P(\rho_0)=\varkappa D/\pi,\quad
 Z'_0=(4\pi\rho_0)^{-1},\quad
 Z_0=\frac{h(\rho_0)}{4\varkappa}+\frac1\varkappa.
 \label{mat:connection-10}
\end{equation}
The additive constant in \eqref{mat:connection-10} makes the
coefficient of degree one in the core radius vanish.
At degree $k\ge1$, let $Z_{<k}$ denote the sum of the previously
constructed layer terms, through degree $k-1$ in $\lambda$.
The new coefficient is governed by
\begin{equation}
 (A_0Z'_{k\ell})'=S_{k\ell},\quad
 A_0=16\pi^2\rho_0^2P'(\rho_0),\quad
 S_{k\ell}=-[\lambda^kL^\ell]F_v[Z_{<k}].
 \label{mat:connection-11}
\end{equation}
The extraction includes all previously determined logarithmic
coefficients. To check the complete nonlinear pressure source, set
\begin{equation}
 u=\lambda Z=\sum_{j\ge1}\lambda^ju_j,\qquad
 v=Z_D/Z'_0-1=\sum_{j\ge1}\lambda^jv_j.
 \label{mat:connection-12}
\end{equation}
With independent logarithmic-polynomial coefficients, the exact finite
convolutions are
\begin{equation}
 \begin{split}
 U_n&=\sum_{r=1}^n(r+1)
       \sum_{j_1+\cdots+j_r=n}u_{j_1}\cdots u_{j_r},\\
 V_n&=\sum_{r=1}^n(-1)^r
       \sum_{j_1+\cdots+j_r=n}v_{j_1}\cdots v_{j_r},\\
 R_n&=\sum_{a+b=n}U_aV_b,\qquad U_0=V_0=1,\\
 \Pi_n&=\sum_{r=1}^n
       \frac{\rho_0^rP^{(r)}(\rho_0)}{r!}
       \sum_{j_1+\cdots+j_r=n}R_{j_1}\cdots R_{j_r}.
 \end{split}                                                    \label{mat:connection-13}
\end{equation}
All composition indices in the inner sums are positive. In particular
\begin{equation}
 R_n=-Z'_n/Z'_0+\widehat R_n,\quad
 \Pi_n=-P'(\rho_0)\rho_0 Z'_n/Z'_0+\widehat\Pi_n,\quad
 -4\pi\partial_D(-P'(\rho_0)\rho_0 Z'_n/Z'_0)
                         =(A_0Z'_n)'.
 \label{mat:connection-14}
\end{equation}
The hatted expressions use only previously determined coefficients and
the radius coefficient $u_n=Z_{n-1}$. The full source is obtained
by multiplication by $(1-u)^2$, one $D$ derivative, the two
time-differentiation polynomials in \eqref{mat:connection-9}, and
the gravity coefficients $U_n$.

At a fixed scale degree, every operation in
\eqref{mat:connection-6} and \eqref{mat:connection-13} involves only
finitely many coefficients. The $\lambda^4\log\lambda$ term in the
high-density pressure has a spatially constant coefficient and vanishes
under $\partial_z$. All other logarithms are retained as polynomials
in the independent variables. If the coordinate conversion increases
the logarithmic degree, the coefficient list is extended by zeros;
the equations already solved are unchanged.

\subsubsection{Core solutions and their endpoint expansions}\label{mat:connection-part-3}

Put $c_n=(2n^2+n-3)|\delta|>0$. For a smooth even center source $F$,
the center solution with $f(0)=a$, $f'(0)=0$ is the fixed point of
\begin{equation}
 f(z)=a+\frac34\int_0^z\frac1{w(s)^4s^4}
                \int_0^s w(t)^3t^4\{c_nf(t)-F(t)\}\dd t\dd s.
 \label{mat:connection-15}
\end{equation}
On a sufficiently short center interval the coefficient of $f$ has
norm at most $Cz_0^2$. The differentiated fixed-integral formula
obtained by $t=s\eta$, $s=z\xi$ preserves smooth even functions.
It yields the unique center-regular solution. Ordinary ODE continuation
extends it over every compact subset of $z<1$. The source extracted
from \eqref{mat:connection-6} is even: its undivided radial force is odd in $z$, while
$w,\psi,J$ are even, $w'$ is odd and $m/z^2$ is odd. Thus division
by $z$ in the source is removable at the center.

Near $q=1-z=0$, define
\begin{equation}
 a(q)=\frac{w(1-q)^4(1-q)^4}{\varkappa^4q^4},\qquad
 b(q)=\frac{w(1-q)^3(1-q)^4}{\varkappa^3q^3}.
 \label{mat:connection-16}
\end{equation}
Both are smooth positive factors, with value one at zero. The equation
$B_nf=F$ is exactly
\begin{equation}
 (q^4af_q)_q=q^3(c f+G),\quad
 c=\frac{3c_n}{4\varkappa}b,\qquad G=-\frac3{4\varkappa}bF.
 \label{mat:connection-17}
\end{equation}
The regular homogeneous solution $U_n^r(0)=1$ is obtained by
integrating \eqref{mat:connection-17} twice from zero. Its integral operator gains a
factor $Cq_0$. Once $U_n^r>0$, a second exact solution is
\begin{equation}
 U_n^s(q)=-3U_n^r(q)
              \int_{q_0}^q\frac{\dd s}{s^4a(s)U_n^r(s)^2}.
 \label{mat:connection-18}
\end{equation}
It has leading term $q^{-3}$; adding a multiple of $U_n^r$
fixes its regular constant if desired. The Wronskian is
\begin{equation}
 q^4a\{U_n^r(U_n^s)' -U_n^s(U_n^r)'\}=-3.
 \label{mat:connection-19}
\end{equation}
Taylor expansion of the positive smooth integrand factor in \eqref{mat:connection-18} gives
the complete regular and singular expansions. In particular
\begin{equation}
 U_n^s=q^{-3}
 -\frac{3(2n^2+n+3)|\delta|}{8\varkappa}q^{-2}
 +O(q^{-1}(1+|\log q|)).
 \label{mat:connection-20}
\end{equation}
The logarithm in \eqref{mat:connection-18} comes from the $s^{-1}$ term. Thus its
coefficient multiplies the entire $U_n^r$, not just its constant.
The branch is defined by the complete formula
\eqref{mat:connection-18}.

For an inhomogeneous source, variation of parameters uses primitives
based at $q_0>0$:
\begin{equation}
 f_p(q)=\frac13U_n^r(q)
           \int_{q_0}^q U_n^s(t)t^3G(t)\dd t
       -\frac13U_n^s(q)
           \int_{q_0}^q U_n^r(t)t^3G(t)\dd t .
 \label{mat:connection-21}
\end{equation}
The positive base point makes these integrals well defined for every
$q>0$, including sources that are not integrable at zero. The general
solution is obtained by adding the two homogeneous branches. For a source term $s q^\alpha(\log q)^p$, the term of highest
logarithmic degree in the particular solution is
\begin{equation}
 -\frac{3s}{4\varkappa(\alpha+1)(\alpha+4)}
 q^{\alpha+1}(\log q)^p,\qquad \alpha\ne-1,-4.
 \label{mat:connection-22}
\end{equation}
At the two resonances the top terms are respectively
\begin{equation}
 -\frac{s}{4\varkappa(p+1)}(\log q)^{p+1},\qquad
 \frac{s}{4\varkappa(p+1)}q^{-3}(\log q)^{p+1}.
 \label{mat:connection-23}
\end{equation}
Lower logarithmic coefficients follow by downward polynomial inversion.
All other coefficients of \eqref{mat:connection-17} raise the endpoint power by at least
one in this recursion. Hence each finite endpoint block is determined
by \eqref{mat:connection-22}--\eqref{mat:connection-23} and two complete homogeneous branches.

The center-normalized homogeneous solution $u_n(0)=1$ has
\begin{equation}
 (w^4z^4u_n')'=\tfrac34c_nw^3z^4u_n.
 \label{mat:connection-24}
\end{equation}
Before a possible first zero the right side is positive; therefore
$u_n'>0$, and no first zero occurs. By \eqref{mat:connection-18}--\eqref{mat:connection-21} it has at most
the $q^{-3}$ endpoint growth. Denote by $[U_n^s]u_n$ the
coefficient of $U_n^s$ in its endpoint decomposition.
The integrated flux gives
\begin{equation}
 [U_n^s]u_n=\frac{c_n}{4\varkappa^4}
                   \int_0^1w^3z^4u_n\dd z>0.
 \label{mat:connection-25}
\end{equation}
The integral is finite. Let $v_{n\ell}$ be the center-regular
particular solution with value zero at the center. Varying
$a_{n\ell}$ in $v_{n\ell}+a_{n\ell}u_n$ changes the singular
coefficient by the positive multiple in \eqref{mat:connection-25}.
For a polynomial in $L$, solve from the highest degree downwards.
The resulting triangular matrix has this same nonzero diagonal and
is invertible for each selected $\beta>0$. Its inverse need not
remain bounded as $\beta\downarrow0$.

\subsubsection{Removal of the initial singular branches and the core index bound}\label{mat:connection-part-4}

Direct extraction from \eqref{mat:connection-6} gives
\begin{equation}
 B_2\phi_{20}=\frac2z\frac{w'}{w^2},\qquad
 B_3\phi_{30}=-6e\phi_{20}.
 \label{mat:connection-26}
\end{equation}
Choose $a_{20}$ to make the coefficient of the complete branch
$U_2^s$ zero. Equations \eqref{mat:connection-22} and \eqref{mat:connection-26} then give
\begin{equation}
 \phi_{20}=-\frac3{4\varkappa^2}q^{-1}+O(\log q).
 \label{mat:connection-27}
\end{equation}
This removes both the possible $q^{-3}$ term and its forced
$q^{-2}$ successor from the homogeneous branch. It does not
remove the forced $q^{-1}$ coefficient.

For the leading layer, if $v=\rho_0^{1/3}$ and
$\mathfrak s=(\varkappa D/\pi)^{1/4}$, the exact pressure law gives
\begin{equation}
 P(v^3)=v^4-v^2+O(\log v),\qquad
 Z_0=\varkappa^{-1}\sqrt{1+v^2}
 =\varkappa^{-1}(\mathfrak s+\tfrac34\mathfrak s^{-1})
                         +O(\mathfrak s^{-3}\log\mathfrak s).
 \label{mat:connection-28}
\end{equation}
The constant $1/\varkappa$ in \eqref{mat:connection-10} exactly cancels the constant
$-1/\varkappa$ in $h/(4\varkappa)$; no degree-one core coefficient is
introduced. Since $M-m(z)=\pi\varkappa^3q^4(1+O(q))$, the second term
of $-\lambda Z_0$ is the coefficient in \eqref{mat:connection-27}.
At $n=3$ choose the complete $U_3^s$ coefficient zero as well.
The source in \eqref{mat:connection-26} is no worse than $q^{-1}\log^{p_{\log}}q$, so its
particular solution has no pole below $q^{-2}$. This second choice is required because the corresponding branch is
absent from the layer expansion; see
Subsubsection~\ref{mat:connection-part-8}.

\noindent\emph{The first logarithmic term.}
The degree-two correction already contains the nonzero logarithmic term
$\{15\delta/(8\varkappa^3)\}\log q$, following the forced pole in
\eqref{mat:connection-27}. Its coefficient is determined by the next
leading profile and flux terms, even though the complete homogeneous singular
branch has been removed. The inhomogeneous equation therefore still
requires a logarithmic term. The full calculation
\eqref{mat:first-coupled-grades}, and the first singular-branch example
\eqref{mat:first-singular-log}, use the same normalizations as the
induction below.

We next determine the range of endpoint powers in the core expansion,
using only the coefficient equations.
Suppose, for all $2\le j<n$, that
\begin{equation}
 \begin{aligned}
 \Phi_j&=\sum_{p\ge1-j}q^pP_{jp}(L,\log q)\\
       &\quad+\text{a differentiated finite remainder}.
 \end{aligned}
 \label{mat:connection-29}
\end{equation}
Then a coefficient of degree $j$ in $J-1$ has smallest
endpoint power $-j$, while a coefficient in $\psi-1$
has smallest power $1-j$. Reciprocals of these factors at
$\lambda=0$ are finite convolution polynomials at every
positive scale degree. Their degree-$r$ coefficients
therefore have smallest endpoint power at least $-r$.

For the pressure term of scale degree $2a$,
spatial differentiation followed by division by $w^3$ gives
smallest endpoint power $-2a$. The remaining relative factors have scale degree
$n-2a$ and smallest power at least $-(n-2a)$. Thus the
smallest possible power in their product is given by
\[
                  q^{-2a}q^{-(n-2a)}=q^{-n}.
\]
At degree four, $\partial_z\log w/w^3=O(q^{-4})$, while
the spatially constant $L$ term differentiates to zero.
Negative powers in the pressure expansion satisfy the same count.

The temporal and gravitational terms are less singular, and the
terms with lower logarithmic powers in \eqref{mat:connection-7}
do not decrease the endpoint power. Consequently the source has
\begin{equation}
 \begin{aligned}
 F_{n\ell}&=\sum_{\alpha\ge-n}q^\alpha
                 P_{n\ell\alpha}(\log q)\\
          &\quad+\text{remainder}.
 \end{aligned}
 \label{mat:connection-30}
\end{equation}
The particular inverse \eqref{mat:connection-22}--\eqref{mat:connection-23} starts at $q^{1-n}$.
For $n\ge4$, the complete singular branch $q^{-3}$ is
already in this range. For $n=2,3$ it was removed above.
Thus \eqref{mat:connection-29} holds at every new degree, independently of any
agreement with the layer.

\subsubsection{Layer solutions, endpoint conditions, and the index bound}\label{mat:connection-part-5}

For a fixed additive constant, the zero-flux solution of
\eqref{mat:connection-11} is
\begin{equation}
 Z_{k\ell}(D)=\gamma_{k+1,\ell}
       +\int_0^D A_0(s)^{-1}\int_0^sS_{k\ell}(t)\dd t\dd s.
 \label{mat:connection-31}
\end{equation}
\equationalias{in:volterra}{mat:connection-31}
At zero, writing $y=D^{1/5}$, the exact pressure law and \eqref{mat:connection-10} give
$\rho_0=y^3b_0(y^2)$ and $A_0=y^8a_0(y)$, with
positive smooth factors. If all previous coefficients have
$Z_j=\gamma_{j+1}+y^2\widehat Z_j(y)$, then
\begin{equation}
 \begin{aligned}
 Z_D/Z'_0&=1+\text{a finite smooth coefficient series},\\
 \rho&=y^3\{\text{a series with a smooth positive leading factor}\},\\
 P(\rho)&=y^5\{\text{a smooth series}\}.
 \end{aligned}
 \label{mat:connection-32}
\end{equation}
Since $y=D^{1/5}$,
$\partial_D(y^5f)=f+(y/5)f_y$. Formula \eqref{mat:connection-13} therefore
gives a smooth bounded source $\widehat S(y)=S(y^5)$.
There is no logarithm of $y$ at this endpoint. Exactly,
\begin{equation}
 \begin{split}
 Z(y^5)-\gamma
 &=25\int_0^y\frac{\xi}{a_0(\xi)}
                     \int_0^1t^4\widehat S(t\xi)\dd t\dd\xi\\
 &=25y^2\int_0^1\frac{s}{a_0(sy)}
                     \int_0^1t^4\widehat S(tsy)\dd t\dd s.
 \end{split}                                                    \label{mat:connection-33}
\end{equation}
Thus $\widehat Z$ has every finite ordinary derivative
requested of the source, by differentiation under fixed integrals.
This both constructs the solution and proves zero vacuum flux.
The difference of two solutions with the same constant and source
satisfies $A_0W'=\text{constant}$; zero flux and $W(0)=0$
give $W=0$.

At infinity put $Q=D^{1/4}$. The leading factors satisfy
$\rho_0\sim Q^3$, $Z'_0\sim Q^{-3}$, and
$A_0\sim Q^7$. If $Z_j=O(Q^{j+1}\log^{p_{\log}} Q)$, then
\begin{equation}
 Z'_j/Z'_0=O(Q^j\log^{p_{\log}} Q),\quad
 [\lambda^r](\rho/\rho_0-1)=O(Q^r\log^{p_{\log}} Q).
 \label{mat:connection-34}
\end{equation}
In \eqref{mat:connection-13}, $P^{(r)}(\rho_0)\rho_0^r$ has largest power
$Q^4$. This follows directly from the exact identity
\begin{equation}
 \rho^rP^{(r)}(\rho)=
       \prod_{\nu=0}^{r-1}(\rho\partial_\rho-\nu)P(\rho).
 \label{mat:connection-35}
\end{equation}
One $D$ derivative lowers the $Q$-power by four. The
temporal and gravity terms at scale degree $k$ have
largest power $Q^k$ as well; in particular the term
$\lambda^4D$ has exactly that count. Consequently
\begin{equation}
 S_k=O(Q^k\log^{p_{\log}} Q),\qquad
 Z_k=O(Q^{k+1}\log^{p_{\log}} Q).
 \label{mat:connection-36}
\end{equation}
The last implication follows from the two integrals in \eqref{mat:connection-31}.
Their extra constant and $Q^{-3}$ terms obey the same bound.
The same statement holds after the iterates of $E_Q$
specified next. Hence every layer radius monomial has positive
$X$-degree and nonnegative $Y$-degree. These index bounds follow
from the layer equation and are available before the matching argument.

To preserve the positive radius and Jacobian factors, estimate the
derivative corrections on the regions where each expansion is used.
The core and long-layer bounds give, respectively,
\[
 \begin{aligned}
 \lambda^nq^{-n}(1+|\log q|)^{p_{\log}}
 &\le C\lambda^{n/2}(1+|\log\lambda|)^{p_{\log}},
 &&q\ge c\sqrt\lambda,\\
 \lambda^kQ^k(1+\log Q)^{p_{\log}}
 &\le C\lambda^{k/2}(1+|\log\lambda|)^{p_{\log}},
 &&1\le Q\le C\lambda^{-1/2}.
 \end{aligned}
\]
The leading layer radius change is
$O(\lambda Q)=O(\sqrt\lambda)$. The bounded $y$-chart is
controlled by \eqref{mat:connection-32}--\eqref{mat:connection-33},
and the compact core by \eqref{mat:connection-15}.

After fixing the finite family, we decrease the upper bound for
$\lambda$ until these corrections preserve positivity in every
region. This choice may depend on the nonleading coefficients; the
uniform leading bounds were fixed before selecting $\beta$.

\subsubsection{Finite asymptotic expansions and differentiated remainders}\label{mat:connection-part-6}

The scale degree and the endpoint exponent are independent indices.
We work with polynomials in an auxiliary scale variable $t$ modulo
$t^{N+1}$, whose coefficients are functions of the endpoint variable
and polynomials in $L$. At each retained scale degree, coefficient
extraction precedes truncation of the endpoint expansion.

Let $E_Q=Q\partial_Q$, $Q\ge Q_b>1$. For integers $r<m$ and $d\ge0$, a
finite expansion of type $(m,r,d)$ means
\begin{equation}
 \begin{gathered}
 f(Q,L)=\sum_{p=r+1}^{m}Q^p f_p(L,\log Q)+f_{\rm rem}(Q,L),\\
 \sup_{Q\ge Q_b}\frac{Q^{-r}|E_Q^j[L^\ell]f_{\rm rem}|}
                         {(1+\log Q)^{p_{\log}}}<\infty,
 \qquad j\le d.
 \end{gathered}
 \label{mat:finite-endpoint-class}
\end{equation}
Only finitely many coefficients in $L$ are present and each $f_p$ is a
finite polynomial; missing powers have zero coefficients. Near $q=0$
the corresponding type $(a,H,d)$ is a sum over $a\le p\le H$ plus a
remainder $O_{E_q^{\le d}}(q^{H+1}(1+|\log q|)^{p_{\log}})$.
The notation $O_{E^{\le d}}$ includes every indicated derivative, not
only a zeroth-order estimate.

\begin{lemma}
\label{mat:finite-polyhomogeneous}
The finite expansion classes defined above are closed under finite
products and coefficient extraction, after truncation at the common
order determined by the remainder estimates. They are also closed
under powers of positive factors, reciprocals, logarithms,
and composition with the exact pressure, truncated in the scale variable,
when the leading factors are independent of $L$. These are the operations in \eqref{mat:connection-6} and \eqref{mat:connection-13}.

The identities
\[
 \partial_D=(4Q^4)^{-1}E_Q,\qquad
 \partial_q=q^{-1}E_q,\qquad E_q=q\partial_q,
\]
show that these derivatives lower the endpoint exponent by four and
one, respectively. At each retained scale degree, the full force uses
at most two derivatives of a coefficient.
\end{lemma}

The proof is given in Subsubsection~\ref{app:profile-finite-calculus}.

\begin{lemma}
\label{mat:finite-volterra}
Let $d\in\mathbb N_0$, $A\in C_{\mathrm{loc}}^{d+1}(0,\infty)$,
and $S\in C_{\mathrm{loc}}^d(0,\infty)$, with $A>0$.
Suppose $S$ is bounded near $D=0$ and
$A(D)=D^{8/5}$ times a positive smooth function of $D^{1/5}$ there.
At infinity put $Q=D^{1/4}$ and suppose
$A(Q^4)^{-1}=Q^{-7}a(Q)$, with $a$ bounded above and below by
positive constants. Define
\[
 B(Q)=3+E_Q\log\!\left(\frac{A(Q^4)}{Q^7}\right).
\]
Assume the following.
\begin{itemize}[leftmargin=1.6em]
\item The expansion of the full $A^{-1}$ (not of $a$) has remainder
exponent $\sigma_A$, while $S$ has largest exponent $k$ and remainder
exponent $\sigma_S<-4$.
\item The remainder estimates hold after applying $E_Q^j$, for
$j\le d+1$ to $A^{-1}$ and $j\le d$ to $S$; the leading factors
satisfy the same bounds. The coefficient $B$ has $d+1$ bounded $E_Q$ derivatives.
\item The remainder exponents satisfy
\[
 R=\max\{\sigma_S+1,\ k+\sigma_A+8,\
                \sigma_S+\sigma_A+8,\ \sigma_A+4\}<-3.
\]
\end{itemize}
Then the zero-flux solution
\[
 Z(D)=\gamma+\int_0^D A(s)^{-1}\int_0^s S(t)\dd t\dd s
\]
has a finite expansion with remainder $O_{E_Q^{\le d+2}}(Q^R
(1+\log Q)^{p_{\log}})$. Its expansion includes a constant term and
all retained terms in the expansion of the $Q^{-3}$ homogeneous branch.
The constants include the full finite-interval and asymptotic-cutoff
contributions.
\end{lemma}

The proof is given in Subsubsection~\ref{app:profile-volterra}.
The singular coefficient contains the full contribution of the source
from the bounded vacuum interval and the cutoff transition, as shown
by \eqref{mat:volterra-actual-moment}. These contributions are retained
when imposing the connection conditions.

\paragraph{Two normalized inverse formulas.}
For the leading layer coefficient $A=a_0Q^7$ and source $S=sQ^k$,
the primitive formula gives
\[
 Z_{\rm p}=\frac{16s}{a_0(k+1)(k+4)}Q^{k+1},
 \qquad k\ne-1,-4.
\]
For $k=-1$ it gives $Z_{\rm p}=16s\log Q/(3a_0)$; for $k=-4$
it gives $Z_{\rm p}=-16sQ^{-3}\log Q/(3a_0)$.
These are identities for the normalized equation
$E_Q(E_Q+3)Z=16QS/a_0$. They are asymptotic terms above a positive
base point, not sources integrated to zero when nonintegrable there.
Actual zero-flux solutions add the regular and homogeneous branches
specified in \eqref{mat:volterra-actual-moment}. The worked
calculation at \ref{mat:worked-connection} checks a nonresonant pole
and both resonant freedoms in the full equation, with its leading profile normalization.

To compare the core and layer solutions on the moving shell, we need
endpoint expansions with differentiated integral remainders. At each
scale degree the expansion must extend beyond both the constant and
$Q^{-3}$ branches, since both coefficients enter the matching
conditions. The following lemma chooses the required endpoint orders
and proves them simultaneously with the derivative bounds.

\begin{lemma}
\label{mat:finite-remainder-schedules}
At each finite stage and for any finite connection amplitudes, the
core and layer solutions have finite expansions in the classes of
Lemma~\ref{mat:finite-polyhomogeneous}, with endpoint orders
\begin{equation}
 H_n=24+2(S-n),\quad 2\le n\le S,\qquad
 R_k=k-2K-10,\quad 0\le k\le K.
 \label{mat:connection-37}
\end{equation}
More precisely:
\begin{enumerate}[label=(\roman*),leftmargin=2em]
\item The core expansion has powers $1-n\le p\le H_n$ and
remainder $O(q^{H_n+1}\log^{p_{\log}}q)$. The complete singular
branch is removed at $n=2,3$.
\item The layer expansion has powers $R_k<p\le k+1$ and
remainder $O(Q^{R_k}\log^{p_{\log}}Q)$.
\end{enumerate}
The bounds hold after $E_q^j$ in the core and $E_Q^j$ in the layer,
for $0\le j\le24$. Since $R_k\le-K-10<-3$, the expansions include
all constant and singular branches used in the matching.
\end{lemma}

The proof is given in Subsubsection~\ref{app:profile-remainder-orders}.

\subsubsection{Common expansion variables and the residual coefficients}\label{mat:connection-part-7}

Let
\begin{equation}
 X=(M-x)^{1/4},\quad Y=\lambda/X=Q^{-1},\quad
 \mathcal{Z}=X\partial_X-Y\partial_Y,\quad \mathcal{E}=Y\partial_Y.
 \label{mat:connection-53}
\end{equation}
The smooth endpoint inverse gives
$q=c_*X(1+O(X))$, $c_*=(\pi\varkappa^3)^{-1/4}$.
Its logarithm is
$\log q=\log X+\log c_*+O(X)$.
The leading-profile regularity \eqref{mat:connection-48} permits
Taylor expansion of this inverse to the orders in
\eqref{mat:connection-37}, with the remainder estimates preserved
under $(X\partial_X)^a(Y\partial_Y)^b$, $a+b\le24$. In particular
\begin{equation}
 \lambda^nq^p=c_*^pX^{n+p}Y^n(1+O(X)),\qquad
 \lambda^{k+1}Q^p=X^{k+1}Y^{k+1-p}.
 \label{mat:connection-54}
\end{equation}
By \eqref{mat:connection-29} and \eqref{mat:connection-36}, every nonconstant radius monomial has
$X$-degree at least one and $Y$-degree nonnegative.
After the two initial removals, the only degree-one term
is $-c_*X$. Hence both the core and layer expansions belong to
\begin{equation}
 R=1-c_*X+XB,\qquad
 -\mathcal{Z}R/X=c_*-(1+\mathcal{Z})B,
 \label{mat:connection-55}
\end{equation}
where $B$ has strictly positive total degree. The second expression is the normalized Jacobian factor. Both factors
are positive on a sufficiently small region where the two overlap
variables are comparable.

In the overlap variables, the exact reduced residual is
\begin{equation}
 \begin{split}
 F_{\rm cor}[R]={}&
 \beta^2\big((1+\mathcal{E})^2-\tfrac32(1+\mathcal{E})\big)R
 +eXY\mathcal{E}(1+\mathcal{E})R\\
 &-\pi R^2(\mathcal{Z}+4)
       \{Y^4P(Y^{-3}\Theta)\}+(M-X^4)R^{-2},\\
 \Theta={}&-\frac{X}{\pi R^2\mathcal{Z}R}.
 \end{split}                                                    \label{mat:connection-56}
\end{equation}
The changes of variables from \eqref{mat:connection-6} and \eqref{mat:connection-9} give \eqref{mat:connection-56}
identically at positive $\lambda,q,D$, before any expansion.
The exact pressure expansion of $Y^4P(Y^{-3}\Theta)$ begins with
\begin{equation}
 \Theta^{4/3}-Y^2\Theta^{2/3}
 +Y^4\{\tfrac12\log\Theta-\tfrac32\log Y
                         +\tfrac32\log2-\tfrac78\}
 +\sum_{j\ge1}c_jY^{4+2j}\Theta^{-2j/3}.
 \label{mat:connection-57}
\end{equation}
\equationalias{in:pressure-series}{mat:connection-57}
For any retained finite degree its remainder is obtained
from \eqref{mat:connection-38} by the positive factor Taylor formula. Because
the remainder of \eqref{mat:connection-55} is divided by $X$ just once,
one extra radius order is retained before this operation.
The orders in \eqref{mat:connection-37} include this additional derivative.

The constant terms of the two factors in \eqref{mat:connection-55} are $1,c_*$.
Finite reciprocal and composition formulas therefore
operate in nonnegative $X,Y$ powers, with logarithmic
polynomial coefficients. The division by $X$ in the normalized Jacobian variation lowers
the $X$-degree by one and leaves the $Y$-degree unchanged. For a new degree-$d$ radius term
$H$, $d\ge2$,
\begin{equation}
 [DF_{\rm cor}(R)H]_{d-1}
 =-\frac{4\varkappa}{3c_*}X^{-1}(\mathcal{Z}+3)\mathcal{Z}H.
 \label{mat:connection-58}
\end{equation}
This follows from
\begin{equation}
 D\Theta(R)H=-\Theta(2H/R+\mathcal{Z}H/\mathcal{Z}R),
 \qquad \pi(\pi c_*)^{-4/3}=\varkappa.
 \label{mat:connection-59}
\end{equation}
The time and gravity terms in \eqref{mat:connection-56} have degree $d$.
Each additional occurrence of $H$ has degree at least
$d-1$, so a term with two occurrences cannot enter
degree $d-1$.

For this coefficient comparison let $R_c^+$ also contain the
auxiliary term $\lambda^S z\Phi_S$. The actual radius in
\eqref{mat:connection-3} is recovered by deleting that term. The two coefficient equations cancel different indices in the overlap
expansion. A core term with factor $\lambda^n$ has exactly
$Y$-degree $n$ in \eqref{mat:connection-54},
including after the endpoint inverse, which depends only
on $X$. Exact auxiliary core cancellation through $n=S$ therefore
implies
\begin{equation}
 [Y^b]F_{\rm cor}[R_c^+]=0,\qquad 0\le b\le S
 \label{mat:connection-60}
\end{equation}
at every $X$-coefficient retained below. A layer
coefficient $Z_k$ has radius $X$-degree $k+1$, and
its first force occurrence has $X$-degree $k$.
Exact cancellation of \eqref{mat:connection-11} through $k=K$ consequently gives
\begin{equation}
 [X^a]F_{\rm cor}[R_v]=0,\qquad 0\le a\le K.
 \label{mat:connection-61}
\end{equation}
The pressure and source index bounds proved in
Subsubsections~\ref{mat:connection-part-4}--\ref{mat:connection-part-5}
ensure that there are no negative powers through
which an omitted scale degree could enter \eqref{mat:connection-60} or
\eqref{mat:connection-61}. Logarithms do not alter either degree. Thus the two coefficient ODEs imply
\eqref{mat:connection-60}--\eqref{mat:connection-61} before the
overlap estimate is proved.

\subsubsection{The two connection constants}\label{mat:connection-part-8}

Write the core and layer expansions as
\[
 1+\sum_{a\ge1,\ n\ge0}X^aY^n C^\nu_{a,n}(L,\ell),
 \qquad \nu=c,v,\qquad L=\log\lambda,\quad\ell=\log X.
\]
Each coefficient is a polynomial in $L,\ell$.
The core coefficient $\Phi_n$ determines row $n$, and
$Z_{a-1}$ determines column $a$. We match the two arrays by
successively choosing the center amplitude and the additive layer
constant.

\emph{The coefficient equation.}
The common leading coefficient is $C_{1,0}=-c_*$. Put
\[
 U=R-1,\qquad
 W=-\frac{\mathcal Z(R-1+c_*X)}{c_*X},\qquad
 \Theta_*=(\pi c_*)^{-1}.
\]
Then
\[
 \Theta=\Theta_*R^{-2}(1+W)^{-1},\qquad
 [X^{a-1}Y^n]W=-c_*^{-1}(a-n+\partial_\ell)C_{a,n}
 \quad((a,n)\ne(1,0)).
\]
The index bounds proved above give nonnegative bidegrees for
$U,W$, with zero constant coefficients. Consequently every
coefficient of a product or composition is a finite sum. For example,
\[
 [X^uY^v](1+U)^\alpha(1+W)^\eta
 =\sum_{\substack{r,s\ge0\\r+s\le u+v}}
       \binom{\alpha}{r}\binom{\eta}{s}
 \sum_{\substack{\mu_1+\cdots+\mu_r+\nu_1+\cdots+\nu_s=(u,v)\\
                  \mu_i,\nu_j\in\mathbb N_0^2\setminus\{(0,0)\}}}
       \prod_{i=1}^r U_{\mu_i}\prod_{j=1}^s W_{\nu_j}.
\]
Here $U_\mu,W_\nu$ denote coefficients, and empty products equal one.
Logarithms of positive factors obey
the same finite Taylor rule, while $\log Y=L-\ell$.
Differentiation preserves the bidegree:
\[
 \mathcal Z(X^aY^nC)=X^aY^n(a-n+\partial_\ell)C,\qquad
 \mathcal E(X^aY^nC)=X^aY^n(n+\partial_L)C.
\]

Let $R_{<a,n}$ contain $1$ and the coefficients with
\[
 1\le a'\le a,\qquad 0\le n'\le n,\qquad(a',n')\ne(a,n).
\]
For $(a,n)\ne(1,0)$, the force coefficient at $X^{a-1}Y^n$ gives
\begin{equation}
 \begin{split}
 &(a-n+\partial_\ell)(a-n+3+\partial_\ell)C_{a,n}
       =G_{a,n},\\
 &G_{a,n}=\frac{3c_*}{4\varkappa}
       [X^{a-1}Y^n]F_{\rm cor}[R_{<a,n}].
 \end{split}
 \label{mat:connection-coefficient-recursion}
\end{equation}
\equationalias{in:corner-recursion}{mat:connection-coefficient-recursion}
Indeed $C_{a',n'}$ enters $U$ at $(a',n')$ and $W$ at
$(a'-1,n')$; the equation of state adds only $(0,2j)$,
$j\ge0$. Thus only $a'\le a$, $n'\le n$ can contribute.
The current coefficient occurs only linearly in $W$, since every
additional nonconstant factor increases an index. Its contribution is
\[
 \frac{4\varkappa}{3}(\mathcal Z+4)W
 =-\frac{4\varkappa}{3c_*}
       X^{-1}(\mathcal Z+3)\mathcal Z(R-1+c_*X).
\]
The time and gravity terms use $U$; the factor $eXY$ increases
both indices. They therefore give no further current coefficient.
This reasoning also applies on $a=1$ and $n=0$.

\emph{Polynomial inversion.}
Set $m=a-n$, $D_\ell=\partial_\ell$, and
$h=\deg_\ell G_{a,n}$. For $m\ne0,-3$, the unique polynomial solution is
\[
 C_{a,n}=\frac1{m(m+3)}
 \sum_{j=0}^{h}
 \left\{-\frac{(2m+3)D_\ell+D_\ell^2}{m(m+3)}\right\}^{j}G_{a,n}.
\]
The numerator lowers degree, so the sum terminates. For $m=0,-3$,
put $\sigma=3,-3$, respectively. Since
$(D_\ell+m)(D_\ell+m+3)=D_\ell(D_\ell+\sigma)$, all polynomial
solutions are
\[
 C_{a,n}(L,\ell)=\gamma(L)+
 \int_0^\ell\sum_{j=0}^{h}
       \frac{(-1)^j}{\sigma^{j+1}}
                 \partial_s^jG_{a,n}(L,s)\dd s.
\]
For zero source the sums are empty. Each resonance leaves exactly
one polynomial $\gamma(L)$ free.

When the preceding coefficients agree, the cancellations
\eqref{mat:connection-60}--\eqref{mat:connection-61} give
\begin{equation}
 \{a-n+\partial_\ell\}\{a-n+3+\partial_\ell\}
          (C^c_{a,n}-C^v_{a,n})=0.
 \label{mat:connection-62}
\end{equation}
In particular, the difference vanishes unless $a-n=0,-3$;
at either resonance it is independent of $\ell$.

\emph{Matching induction.}
At stage $N$ we establish
\begin{equation}
 C^c_{a,n}=C^v_{a,n},
 \qquad 1\le a\le N,\quad 0\le n\le N.
 \label{mat:connection-rectangle}
\end{equation}
The leading expansion of $Z_0$ gives
$C_{1,0}=-c_*$ and $C_{1,1}=0$, agreeing with the core.
Assume stage $N-1$ is complete.

First solve the degree-$N$ core equation, leaving its center
amplitude free, and compare row $n=N$ in the order
$a=1,\ldots,N-1$. The source of each comparison uses the old
rectangle and the preceding part of this row. The force equations
used are
\[
 \text{core degree }N,\qquad
 \text{layer degree }a-1\le N-2.
\]
For $N=2,3$, the complete singular branch is removed because
its column $a=N-3$ is nonpositive. All comparisons are
nonresonant. For $N\ge4$, the only resonance is
$a=N-3$, whose layer coefficient is supplied by $Z_{N-4}$.
By \eqref{mat:connection-62},
\[
 C^c_{N-3,N}-C^v_{N-3,N}=\gamma(L).
\]
The nonzero connection coefficient \eqref{mat:connection-25}
determines the center amplitude that cancels $\gamma$.
The change starts at column $N-3$, preserving the earlier
comparisons; the remaining columns agree by polynomial inversion.

Next construct $Z_{N-1}$ from \eqref{mat:connection-31}, using
only $Z_0,\ldots,Z_{N-2}$ and leaving its additive polynomial
$\gamma_N(L)$ free. Compare column $a=N$ in the order
$n=0,\ldots,N$. The predecessors lie in the old rectangle,
the new row, or the preceding part of this column. Now the force
equations used are
\[
 \text{layer degree }N-1,\qquad \text{core degree }n\le N.
\]
For $n<N$, $N-n>0$ gives equality. At $n=N$, first set
$\gamma_N=0$ and write the difference as $\gamma(L)$. Then
\[
 (u_c-u_v)|_{\gamma_N}-(u_c-u_v)|_{\gamma_N=0}
       =+\lambda^N\gamma_N(L).
\]
Choose $\gamma_N=-\gamma$; no other coefficient in the column
changes. This proves \eqref{mat:connection-rectangle} at stage $N$.

\emph{The finite linear system.}
For $N\ge4$, write $a_{N\ell}=\phi_{N\ell}(0)$ and order the
external logarithmic degrees as $\ell=d,d-1,\ldots,0$.
Equation~\eqref{mat:core-log-recursion} is triangular in this order.
By \eqref{mat:connection-25},
\[
 \frac{\partial [L^\ell]C^c_{N-3,N}}{\partial a_{N\ell}}
 =c_*^{-3}[U_N^s]u_N
 =\frac{c_*^{-3}c_N}{4\varkappa^4}
                  \int_0^1w^3z^4u_N\dd z
 =:c_N^{\rm sing}>0,
\]
where
$\lambda^Nzq^{-3}=c_*^{-3}X^{N-3}Y^N(1+O(X))$.
The center data therefore solve
\[
 c_N^{\rm sing}a_{N\ell}
       +\sum_{r=\ell+1}^{d}t_{\ell r}a_{Nr}
       =b_\ell,\qquad \ell=d,d-1,\ldots,0.
\]
The right sides depend on the particular core solution and the
previously constructed layer column. The regular difference changes
by $+\gamma_N(L)$, so the full matching matrix is
\begin{equation}
 \begin{pmatrix}T_N&0\\ *&I\end{pmatrix},
 \qquad \operatorname{diag}T_N=c_N^{\rm sing},
 \qquad c_N^{\rm sing}>0.
 \label{mat:connection-63}
\end{equation}
\equationalias{in:connection-map}{mat:connection-63}
Later core or layer coefficients add new rows or columns and preserve
all previous equalities. In particular, the singular branch in
column $N$ is used at stage $N+3$; no condition at infinity is
added to the zero-flux condition and the additive constant.

The logarithmic degrees remain finite. Writing $L_\lambda,L_X,L_Y,L_q$
for the respective logarithms,
\begin{equation}
 L_\lambda=L_X+L_Y,\qquad
 L_q=L_X+\log c_*+O(X),\qquad \log Q=-L_Y.
 \label{mat:connection-64}
\end{equation}
At each core degree choose the external $L$-degree to dominate
those of the source and prescribed singular coefficient, then solve
\eqref{mat:core-log-recursion} downwards. The two resonances add
finitely many internal logarithms. The resulting regular coefficient
determines $\gamma_N$. Increasing its degree changes neither current
equation: an additive layer constant first enters the force one scale
degree later.

Induction through $S=K+1$ constructs
$\Phi_2,\ldots,\Phi_S$ and $Z_0,\ldots,Z_K$.
The auxiliary $\Phi_S$ fixes the regular constant in $Z_K$ and is
then omitted from the physical core truncation. Every nonleading
Jacobian increment has positive degree after division by $X$;
the leading normalized coefficient in \eqref{mat:connection-55}
remains $c_*$.

\paragraph{The first coupled degrees and the logarithmic terms.}\label{mat:worked-connection}
The following calculation illustrates the induction just proved; its
constants are obtained from the exact pressure law. Write $q=1-z$ and retain
$\delta=-\beta^2/2$, $\varkappa=-w'(1)$. The equation for the leading profile gives
\[
 w(1-q)=\varkappa q
  \left\{1+\left(1-\frac{3\delta}{8\varkappa}\right)q+O(q^2)\right\}.
\]
In \eqref{mat:connection-17} the corresponding first coefficients are
$a(q)=1-3\delta q/(2\varkappa)+O(q^2)$ and
$b(q)=1+(-1-9\delta/(8\varkappa))q+O(q^2)$.
Substituting $Aq^{-1}+B\log q$ in the degree-two flux equation,
and using the zero complete singular amplitude, gives
\begin{equation}
 \begin{split}
 \phi_{20}
 &=-\frac{3}{4\varkappa^2}q^{-1}
     +\frac{15\delta}{8\varkappa^3}\log q+O(1),\\
 \phi_{30}
 &=-\frac{9e}{8\varkappa^3}\log q+O(1).
 \end{split}
 \label{mat:first-coupled-grades}
\end{equation}
Indeed the $q$ and $q^2$ coefficients of the degree-two flux
equation give
\[
 -2A=\frac{3}{2\varkappa^2},\qquad
 3\left(B+\frac{3\delta A}{2\varkappa}\right)
                              =\frac{9\delta}{4\varkappa^3}.
\]
The second line follows from $B_3\phi_{30}=-6e\phi_{20}$
and the $\alpha=-1$ inverse in \eqref{mat:connection-23}.
The integral remainders of Lemma~\ref{mat:finite-remainder-schedules}
justify these expansions after each of the retained iterates of $E_q$ and $E_Q$. In the overlap variables the
first line yields
\[
 C_{1,2}=-\frac{3}{4\varkappa^2c_*},\qquad
 C_{2,2}=\frac{15\delta}{8\varkappa^3}\ell+\gamma(L).
\]
Thus the regular resonance $a-n=0$ fixes the logarithmic coefficient
and leaves exactly the additive polynomial determined by the layer
normalization, as in the general induction.

The first singular-branch comparison occurs at $(a,n)=(1,4)$.
Put $v=\rho_0^{1/3}$ and
$\mathfrak s=(\varkappa D/\pi)^{1/4}$. Expanding the exact implicit
equation through its constant term gives
\[
 \begin{split}
 P(v^3)&=v^4-v^2+\tfrac32\log(2v)-\tfrac78+O(v^{-2}),\\
 v&=\mathfrak s+\tfrac14\mathfrak s^{-1}
       +\left\{\tfrac14-\tfrac38\log(2\mathfrak s)\right\}
          \mathfrak s^{-3}
       +O(\mathfrak s^{-5}\log\mathfrak s),\\
 Z_0&=\varkappa^{-1}\left\{
       \mathfrak s+\tfrac34\mathfrak s^{-1}
       -\tfrac38\mathfrak s^{-3}\log(2\mathfrak s)\right\}
       +O(\mathfrak s^{-5}\log\mathfrak s).
 \end{split}
\]
The constants at order $\mathfrak s^{-3}$ in the last line cancel.
Since $\mathfrak s=\varkappa c_*/Y$, the layer coefficient is
\begin{equation}
 C^v_{1,4}=A_4\{\ell-L+\log(2\varkappa c_*)\},\qquad
 A_4=\frac{3}{8\varkappa^4c_*^3}.
 \label{mat:first-singular-log}
\end{equation}
This also follows directly from the force equation in the matching coordinates. At $X^0Y^4$
every algebraic term of that bidegree is killed by
$\mathcal Z+4$. The explicit term
$-\tfrac32Y^4\log Y$ in \eqref{mat:connection-57} leaves the force
coefficient $-3\pi/2$, so
\[
 G_{1,4}=-\frac{9\pi c_*}{8\varkappa}=-3A_4,
 \qquad
 \partial_\ell(\partial_\ell-3)C_{1,4}=-3A_4.
\]
Here $\pi\varkappa^3c_*^4=1$. Hence $C_{1,4}=A_4\ell+\gamma(L)$,
and the singular matching condition fixes
$\gamma(L)=A_4\{-L+\log(2\varkappa c_*)\}$.
At core degree four the original force source has no external $L$
term. Its newly required top logarithmic coefficient therefore solves
$B_4\phi_{41}=0$. The corresponding center datum is
$\phi_{41}(0)=-A_4/c_4^{\rm sing}$, which is finite at each selected
$\beta>0$. The coefficient of the lower logarithmic power is then determined by
the additional source in \eqref{mat:core-log-recursion}.
Thus matching requires a logarithmic term even though the degree-four
core source has no external logarithm. Its center amplitude contains
the inverse of the singular connection coefficient, which may be
unbounded as the profile parameter tends to zero. Higher degrees are
constructed by the same rectangular induction.

\subsubsection{Matching of coefficients and the overlap estimate}\label{mat:connection-part-9}

By \eqref{mat:connection-rectangle}, the two radius expansions
agree at every total degree at most $K$. Indeed,
$a+n\le K$ implies $1\le a\le K$, $0\le n\le K$, and the
corresponding force degree is $a-1\le K-1$.
The cancellations \eqref{mat:connection-60}--\eqref{mat:connection-61}
and the index bounds exclude any contribution from omitted scale
degrees.

Write $E_X=X\partial_X$ and $E_Y=Y\partial_Y$.
On $c\le X/Y\le C$ with $XY=\lambda$,
\begin{equation}
 |E_X^aE_Y^b\{X^\alpha Y^\eta
                     P(\log X,\log Y)\}|
 \le C\lambda^{(\alpha+\eta)/2}
                      (1+|\log\lambda|)^{p_{\log}},\quad a+b\le24.
 \label{mat:connection-65}
\end{equation}
Every unmatched retained monomial has total degree at least $K+1$.
The core remainder in \eqref{mat:connection-37} has total degree
\begin{equation}
 2n+H_n+1=2S+25>K+1.
 \label{mat:connection-66}
\end{equation}
The layer remainder has total degree
\begin{equation}
 2(k+1)-R_k=k+2K+12>K+1.
 \label{mat:connection-67}
\end{equation}
The omitted auxiliary term $\Phi_S$ has minimal degree
$S+1=K+2$ by \eqref{mat:connection-29}; the leading-coordinate
remainder is taken beyond $K+1$ using \eqref{mat:connection-48}.
The differentiated integral remainder bounds therefore yield
\[
 \sum_{a+b\le24}
       |E_X^aE_Y^b(u_c-u_v)|
 \le C\lambda^{(K+1)/2}(1+|\log\lambda|)^{p_{\log}}.
\]

On $q\sim\sqrt\lambda$, the endpoint inverse gives
$X,Y\sim\sqrt\lambda$. At fixed $z$,
$\lambda\partial_\lambda=\mathcal E=E_Y$; at fixed $\lambda$,
\begin{equation}
 q\partial_q=
 \frac{qX_q}{X}\,(X\partial_X-Y\partial_Y),
 \qquad \frac{qX_q}{X}=1+O(q).
 \label{mat:connection-68}
\end{equation}
The coefficient $qX_q/X$ and its required derivatives through
order $24$ are uniformly bounded. Iterating
\eqref{mat:connection-68} therefore expresses every derivative
in \eqref{mat:connection-4} as a bounded linear combination of
$E_X^aE_Y^b$, $a+b\le24$. The displayed estimate proves
\eqref{mat:connection-4}.

\paragraph{The matched radius.}
Choose one fixed smooth function $\chi_{\rm m}$, equal to zero on
$(-\infty,c_1]$, equal to one on $[c_2,\infty)$, with
$0<c_1<c_2$ and $0\le\chi_{\rm m}\le1$. Set
\begin{equation}
 A(\lambda,m(z))=
 \lambda\left\{u_v(\lambda,m(z))+
   \chi_{\rm m}(q/\sqrt\lambda)
        (u_c(\lambda,z)-u_v(\lambda,m(z)))\right\},
 \qquad q=1-z.
 \label{mat:matched-radius}
\end{equation}
Near the vacuum boundary we use the layer radius; toward the center
we use the core radius. This defines the matched radius on the whole
mass interval. The matching cutoff $\chi_{\rm m}$ is distinct from
the localization cutoffs $\chi_c,\chi_i,\chi_5,\chi_8$ in
Definition~\ref{at:atlas}.
Write $u=A/\lambda$ and $\chi=\chi_{\rm m}(q/\sqrt\lambda)$.
On the matching shell, $u_{c,z},u_{v,z}\ge c>0$ and
\eqref{mat:connection-4} gives
\[
 \begin{aligned}
 u_z&=(1-\chi)u_{v,z}+\chi u_{c,z}+\chi_z(u_c-u_v),\\
 |\chi_z(u_c-u_v)|
 &\le C\lambda^{-1/2}\lambda^{(K+1)/2}
                  (1+|\log\lambda|)^{p_{\log}}=o(1).
 \end{aligned}
\]
Thus $u_z\ge c/2$ after decreasing the upper bound for $\lambda$.
Together with positivity on the uncut regions, this gives $A_x>0$
on $0<x<M$ and defines the density on the whole mass interval.

\subsection{Derivative bounds for the constructed profiles}\label{mat:native}

The following growth bounds provide the regularity of the terminal
data and control the normalized coefficients on the entire cover.
They do not yet give decay of the residual: that requires the
cancellations in the coefficient equations and is proved in
Subsections~\ref{mat:endpoint}--\ref{mat:material}.

We use $\|\cdot\|_{m,a}$ with the cutoffs of
Definition~\ref{at:atlas}. These are derivative bounds for the
prescribed profiles; no boundary condition from $X^m$ is imposed
on their forces. The factors $5$ and $8$ occur in the exact mass
measures and are omitted from this norm. Local supremum norms are
taken on the enlarged coordinate intervals before localization.

Fix $\beta>0$, $e\ge0$, the finite matching orders $K,S=K+1$,
and the connection amplitudes constructed above. Write $A(\lambda,x)$
for the matched radius, with core terms through degree $K$ and layer
coefficients through $Z_K$; the auxiliary core term of degree $S$
is omitted. Constants in this subsection may depend on these data and the fixed
cutoffs. Choose $A_*>0$ sufficiently small for this family. The
constants are then independent of $0<\lambda\le A_*$, but need not
be uniform as $\beta\downarrow0$.

The covering atlas has a fixed center chart, a fixed compact interior
chart, the bounded chart $y_5=D^{1/5}$, and the high-density chart
\[
 D=(M-x)/\lambda^4,\qquad y_8=D^{1/8},\qquad
 c<y_8<C_e\lambda^{-1/2}.
\]
The layer formula is used only through $y_8=O(\lambda^{-1/4})$.
On the remainder of the high-density chart we use the finite core radius.
The corresponding mass measures are $5\lambda^4y_5^4\dd y_5$ and
$8\lambda^4y_8^7\dd y_8$. The cutoffs and their enlarged chart intervals
are fixed as specified in the construction of the coordinate cover.

We use ordinary derivatives in the bounded local charts, and retain the
radial-vector lift at the center. In the unbounded layer coordinate put $Q=D^{1/4}$
and $E_Q=Q\partial_Q$; on the core edge put $q=1-z$ and $E_q=q\partial_q$.
These conormal derivatives are distinct from ordinary endpoint
derivatives. In the following lemma and its proof in
Subsubsection~\ref{app:profile-coefficient-growth}, $t_n$
denotes the derivative order stated in the lemma. For the physical-coordinate
estimates we shall use the larger orders in \eqref{mat:units-5}.

\begin{lemma}
\label{mat:native-lem:wave2-finite-native-growth}
For any fixed integer $31\le t_*\le60$ and $2\le n\le S$, set
\[
 t_n=t_*+2(S-n).
\]
At the center put $\Phi_{n\ell}(\eta)=\phi_{n\ell}(\sqrt\eta)$,
and at vacuum write
$Z_{n-1,\ell}=\gamma_{n\ell}+y_5^2\widehat Z_{n-1,\ell}$.
The fixed coefficient functions have the following bounds, with finite
logarithmic powers $p_{\log,n}$:
\begin{align}
 \|\Phi_{n\ell}\|_{C^{t_n+1}}
       +\|\widehat Z_{n-1,\ell}\|_{C^{t_n}}&\le C,
 \label{mat:native-eq:wave2-native-fixed-coefficient}\\
 |E_q^j\phi_{n\ell}(1-q)|&\le
        Cq^{1-n}(1+|\log q|)^{p_{\log,n}},
 &0\le j&\le t_n,
 \label{mat:native-eq:wave2-core-coefficient-growth}\\
 |E_Q^jZ_{n-1,\ell}(Q^4)|&\le
        CQ^n(1+\log Q)^{p_{\log,n}},
 &Q\ge1,\quad 0\le j&\le t_n.
 \label{mat:native-eq:wave2-layer-coefficient-growth}
\end{align}
The first norms are on the fixed enlarged center and bounded vacuum
intervals. Ordinary compact-interior coefficient norms through $t_n$ are
bounded as well. The constants may depend on the fixed connection
amplitudes chosen in the finite construction.
\end{lemma}

The proof is given in Subsubsection~\ref{app:profile-coefficient-growth}.

\begin{proposition}
\label{mat:native-prop:wave2-native-profile-source}
Let $\overline R_i=\partial_\tau^iA$ at fixed enclosed mass, evaluated
at the selected scale, and pull these functions back by the uncorrected
initial-radius map only after differentiation. Let $|\cdot|_{m,\infty,a}$
be the maximum of the enlarged local chartwise $C^m$ norms, with the
regular radial-vector lift at the center. Then
\begin{align}
 |\overline R_i|_{m,\infty,a}
      &\le C\lambda^{1-3i/2},&&i+m\le27,
 \label{mat:native-eq:wave2-native-matched-rows}\\
 |\partial_s^r\mathcal M[A]|_{m,\infty,a}
      +\|\partial_s^r\mathcal M[A]\|_{m,a}
      &\le C\lambda^{-2},&&r+m\le25.
 \label{mat:native-eq:wave2-direct-profile-force-source}
\end{align}
These estimates bound the prescribed profile and its force. The
residual decay estimates are proved separately below; only time derivatives through order $24$ are required in the solution space.
\end{proposition}

The proof is given in Subsubsection~\ref{app:profile-radius-force}.

\begin{corollary}
\label{mat:native-cor:wave2-selected-viscous-source}
For any scalar $\nu\in C^{22}$ in similarity time,
\begin{equation}
 \left\|\partial_s^{22}
       \bigl(\nu\lambda^3\partial_s\mathcal M[A]\bigr)\right\|_{0,a}
 \le C\lambda\sum_{j=0}^{22}|\partial_s^j\nu|.
 \label{mat:native-eq:wave2-selected-viscous-source}
\end{equation}
\end{corollary}

\begin{proof}
By Leibniz's rule,
\[
 \begin{aligned}
 &\left\|\partial_s^{22}
       (\nu\lambda^3\partial_s\mathcal M[A])\right\|_{0,a}\\
 &\quad\le\sum_{j+k+\ell=22}\frac{22!}{j!k!\ell!}
       |\partial_s^j\nu|\,|\partial_s^k\lambda^3|
       \|\partial_s^{\ell+1}\mathcal M[A]\|_{0,a}\\
 &\quad\le C\lambda^3\lambda^{-2}
                    \sum_{j=0}^{22}|\partial_s^j\nu|.
 \end{aligned}
\]
Here $|\partial_s^k\lambda^3|\le C_k\lambda^3$, and
$\ell+1\le23$ lies in the force-derivative range of
\eqref{mat:native-eq:wave2-direct-profile-force-source}.
\end{proof}

The construction of compatible initial data uses these profile and
force derivatives at a fixed positive scale. We next express their regularity in physical
coordinates and prove uniform bounds for the normalized operator
coefficients.

\subsection{Physical regularity and normalized coefficients}\label{mat:units}

\label{mat:units-part-1}

Fix the matched trajectory $A(\lambda,x)$, with $\beta>0$,
coefficients through stage $S=K+1$, and the prescribed cutoffs.
To prepare initial data at a fixed $\lambda_* >0$, use the radius
$a=A(\lambda_*,x)$ as coordinate, and denote the resulting enthalpy
and radial velocity by $h_*$ and $V_*$. In the physical vacuum
coordinate $y_{\rm ph}=(A(\lambda_*,M)-a)^{1/2}$, their regularity is

\begin{equation}
 h_*(y_{\rm ph})\in C^{58},\qquad
 V_*(y_{\rm ph})\in C^{57},\qquad
 h_*=y_{\rm ph}^2p_*,\quad p_*>0,\qquad V_*'(0)=0.       \label{mat:units-1}
\end{equation}

The center functions have the corresponding regular scalar/vector
Cartesian lifts; the fixed interior functions have the same finite
orders. These regularity assertions persist under the finite smooth
corrections that impose compatibility on the inviscid and
Kelvin--Voigt initial data. The constants in these physical-coordinate bounds may depend on
$\lambda_*$.

For the \emph{uncorrected reference trajectory}, on the covering edge charts

\begin{equation}
 x=M-\lambda^4y^n,\qquad n\in\{5,8\},\qquad
 r=A/\lambda,\quad
 \gamma=-A_y/(\lambda^2y),\quad T_\gamma=y\gamma_y,
                                                        \label{mat:units-2}
\end{equation}

the exact profile coefficients

\begin{equation}
 a_n=P'(\rho)/(\gamma^2y^2),\quad b_n=y\partial_ya_n,
 \quad v_n\quad\hbox{defined in \eqref{mat:units-17} below}               \label{mat:units-3}
\end{equation}

have the following full expressions:
\[
 \begin{gathered}
 \rho=\frac{ny^{n-2}}{4\pi r^2\gamma},\qquad
 a_n=\frac{P'(\rho)}{\gamma^2y^2},\\
 b_n=y(G_n)_y-\lambda y^2\gamma(G_n)_r+(G_n)_\gamma T_\gamma,
 \qquad G_n(y,r,\gamma)=a_n.
 \end{gathered}
\]

\begin{equation}
 \begin{aligned}
 v_n={}&\lambda\biggl[
 \frac{2\lambda\gamma^2y^2a_n}{r^2}\\
 &\quad+\frac{2k(\rho)\gamma a_n}{r}
     \left(n-2+\frac{2\lambda y^2\gamma}{r}
                           -\frac{T_\gamma}{\gamma}\right)
 -\frac{2x}{r^3}\biggr],\\
 k(\rho)={}&\frac{2+\rho^{2/3}}{3(1+\rho^{2/3})}.
 \end{aligned}
 \label{mat:units-17}
\end{equation}

They satisfy

\begin{equation}
 0<c\le r,\gamma,a_n\le C,
 \qquad
 \|a_n\|_{W^{21,\infty}}+
 \|a_n^{-1}\|_{W^{21,\infty}}+
 \|b_n\|_{W^{21,\infty}}+
 \|v_n\|_{W^{21,\infty}}\le C.                           \label{mat:units-4}
\end{equation}

Here the $n=5$ interval is fixed and bounded, while the $n=8$ interval
is $y_-<y<C_e\lambda^{-1/2}$. The constant in \eqref{mat:units-4} is independent of
$0<\lambda\le A_*$, after $A_*$ is decreased for the already selected
member. Past the matching shell, \eqref{mat:units-4} uses the core coefficients.
These constants may depend on the chosen profile; the leading constants
used in the energy estimate will be fixed uniformly before selecting
$\beta$.

The derivative and coordinate estimates are proved in
Subsection~\ref{app:profile-physical}. The proof treats the bounded
vacuum chart, the growing edge chart, and its core portion separately;
the positive Jacobian is preserved across the matching shell.

\subsection{The bounded endpoint residual and its physical time derivatives}
\label{mat:endpoint}

On a fixed enlarged interval $0\le y=y_5\le Y_0$, the finite layer profile is
\begin{equation}
 Z(\lambda,y)=Z_0(y)+
 \sum_{n=1}^K\sum_\ell\lambda^n(\log\lambda)^\ell Z_{n\ell}(y),
 \qquad Z_{n\ell}=\gamma_{n+1,\ell}+y^2G_{n\ell}(y).
 \label{mat:endpoint-1}
\end{equation}
All $G_{n\ell}$ are $C^{60}$ by the coefficient bounds in
Subsubsection~\ref{mat:units-part-2}.
For this finite sum, define the quotients
\[
 G(\lambda,y)=\frac{Z(\lambda,y)-Z(\lambda,0)}{y^2},
 \qquad Q(\lambda,y)=\frac{Z_y(\lambda,y)}y=2G+yG_y.
\]
The leading quotients are positive. The expansion
\eqref{mat:endpoint-1} gives convergence of $G,Q$ to these
quotients in $C^{59}$ as $\lambda\downarrow0$.
Their positive lower bounds therefore persist for a sufficiently small
upper bound on $\lambda$. No sign is required of an individual
nonleading $G_{n\ell}$.

Put $u=1-\lambda Z$, $A=\lambda u$, and
$T=\lambda\partial_\lambda-\frac45y\partial_y$.
The reduced residual and the normalized forcing residual are
\begin{equation}
 F_v=\lambda^2(A_{\tau\tau}+\mathcal{M}[A]),\qquad
 \mathcal{R}_u=\lambda^{-3}F_v,\qquad
 \partial_\tau=\lambda^{-3/2}bT.
 \label{mat:endpoint-2}
\end{equation}
We prove the endpoint derivative estimates
\begin{equation}
 \|T^q\mathcal{R}_u\|_{C^{28}([0,Y_0])}
 \le C\lambda^{K-2}(1+|\log\lambda|)^{p_{\log}},\qquad 0\le q\le22.
 \label{mat:endpoint-3}
\end{equation}
All constants here are for the selected positive $\beta$, finite family,
and fixed enlarged interval, and are independent of $\lambda$.

Formula \eqref{mat:endpoint-7} uses at most two derivatives of $G$,
and each physical-time derivative at fixed mass adds at most one.
This gives the higher regularity needed for compatible data at a
fixed positive scale. The bounds may depend on that scale; the decay
estimate \eqref{mat:endpoint-3} is uniform in $\lambda$ in its stated
derivative range. Subsection~\ref{app:profile-endpoint} proves both
assertions from the coefficient equations and Taylor's integral
remainder.

To impose compatibility, we first take physical-time derivatives,
then pull back the full residual to the uncorrected initial-radius
coordinate, and finally take its endpoint coefficients. See
\eqref{mat:endpoint-14} and \eqref{mat:endpoint-19}.
The resulting coefficients satisfy \eqref{mat:endpoint-20} for
$K\ge2J^\sharp+86$, where $J^\sharp\ge3$ is the preparation order.

\subsection{Differentiated remainders in the core and layer}\label{mat:face}

The core and layer equations cancel the coefficients of the exact
reduced residual through scale degree $K$. At each new degree, the
known part of the residual is moved to the right-hand side of the
coefficient equation. The resulting identities \eqref{mat:face-36}
and \eqref{mat:face-38} hold for the coefficient functions before any
endpoint expansion or localization.

The small parameters are
\[
 \frac{\lambda}{q}\le C\sqrt\lambda\quad\hbox{in the core},
 \qquad
 \lambda Q\le C\sqrt\lambda\quad\hbox{in the layer}.
\]
At a fixed core or layer coordinate, the geometric factors remain
positive along the corresponding Taylor segment. Material and
conormal differentiation preserve the coefficient cancellations,
so Taylor's integral formula \eqref{mat:face-10} gives the
differentiated residual bounds below.

\label{mat:face-part-1}

Fix $\beta>0$, $e\ge0$, the finite matching orders $K,S=K+1$, and the coefficient
functions and connection amplitudes constructed above. Constants may
depend on these choices and on the fixed enlarged charts. After
decreasing $A_*$, they are independent of the lower terminal scale
and of $0<\lambda\le A_*$, but need not be uniform as
$\beta\downarrow0$.

With $D=(M-x)/\lambda^4$, the two uncut radii are exactly
\[
 A_c=\lambda u_c,\qquad
 u_c=z\left(1+\sum_{n=2}^K\sum_{\ell\in\Lambda_n}
                    \lambda^n(\log\lambda)^\ell\phi_{n\ell}(z)\right),
\]
\begin{equation}
 A_v=\lambda u_v,\qquad
 u_v=1-\lambda Z_0(D)
       -\sum_{n=1}^K\sum_{\ell\in\Lambda_{n+1}}
                   \lambda^{n+1}(\log\lambda)^\ell Z_{n\ell}(D).
 \label{mat:face-1}
\end{equation}
In particular $Z_K$ is retained and $\phi_{K+1}$ is not added.
Define their exact reduced residuals by
\begin{equation}
 \mathcal F_c=\lambda^2(A_{c,\tau\tau}+\mathcal M[A_c]),
 \qquad
 \mathcal F_v=\lambda^2(A_{v,\tau\tau}+\mathcal M[A_v]).
                                                               \label{mat:face-2}
\end{equation}
The material derivative is $\partial_s$ at fixed mass, with
$\lambda_s=b\lambda$, $b^2=\beta^2+e\lambda$.
In the core let $q=1-z$, and on the long layer let $Q=D^{1/4}$.
Conormal derivatives mean $E_q=q\partial_q$ or $E_Q=Q\partial_Q$
at fixed $\lambda$.

\paragraph{Residual estimates.}
For every composition $W$ of at most 21 derivatives chosen from
$\partial_s,E_q$ in the core, or $\partial_s,E_Q$ in the layer,
\begin{equation}
 |W\mathcal F_c|
 \le C(\lambda/q)^{K+1}
              (1+|\log\lambda|+|\log q|)^{p_{\log}},
 \quad c\sqrt\lambda\le q\le q_0<1,                    \label{mat:face-3}
\end{equation}
\begin{equation}
 |W\mathcal F_v|
 \le C(\lambda Q)^{K+1}
              (1+|\log\lambda|+\log Q)^{p_{\log}},
 \quad 1\le Q\le C_0\lambda^{-1/2}.                    \label{mat:face-4}
\end{equation}
On fixed compact subsets of the core interior, every mixed material
and ordinary spatial derivative of total order at most 21 is bounded by
\begin{equation}
 C\lambda^{K+1}(1+|\log\lambda|)^{p_{\log}}.                    \label{mat:face-5}
\end{equation}
At the center, \eqref{mat:face-5} holds for the complete Cartesian radial-vector lift
of the residual, with ordinary Cartesian spatial derivatives. No
estimate for separately singular radial summands is asserted.

The coefficient equations in Subsection~\ref{mat:connection} hold
for the full coefficient functions. We combine them with the derivative
bounds from Subsections~\ref{mat:native}--\ref{mat:units}, which give,
through $r=23$,
\begin{equation}
 |E_q^r\phi_{n\ell}|\le Cq^{1-n}(1+|\log q|)^{p_{\log}},\qquad
 |E_Q^r Z_{n\ell}(Q^4)|\le CQ^{n+1}(1+\log Q)^{p_{\log}}.        \label{mat:face-6}
\end{equation}
The second estimate is used through one additional coefficient derivative
when forming $Z_{n\ell}'/Z_0'$; the total count in the proof never exceeds 23.
Fixed-center scalar coefficient functions in the squared radius and
fixed-interior coefficient functions have the same required finite
regularity. The derivative counts
\begin{equation}
 t_n=60+2(S-n),\quad t_{n-1}-2=t_n,\quad
 t_2+2=2S+58\le2S+60                                  \label{mat:face-7}
\end{equation}
give the required coefficient regularity. The differentiated remainder
estimates are proved in Subsection~\ref{app:profile-core-layer}.

\subsection{The moving cutoff and mixed source estimates}\label{mat:cutoff}

We now estimate the errors introduced by the matching cutoff.
On the matching shell, the pressure expression
\eqref{mat:cutoff-11} depends smoothly on $u,u_z,qu_{zz}$, where
$u$ ranges over normalized radii with positive Jacobian. The overlap
estimate loses one factor $q^{-1}=O(\lambda^{-1/2})$, and hence gives
\[
 \lambda^{(K+1)/2}q^{-1}(1+|\log\lambda|)^{p_{\log}}
 \le C\lambda^{K/2}(1+|\log\lambda|)^{p_{\log}}.
\]
This is the pressure error in \eqref{mat:cutoff-13}. Time
differentiation gives the commutator \eqref{mat:cutoff-14}.

We add these errors to the uncut residuals and integrate against the
coordinate measures to obtain the mixed source bound. The additional
material derivative needed for the energy estimate is treated
separately in $L^2(\dd x)$ in Subsection~\ref{mat:material}.

\label{mat:cutoff-part-1}

Fix the selected $\beta>0$, $e\ge0$, $K$, and finite cutoffs.
All constants may depend on that selection but not on
$0<\lambda\le A_*$, the lower terminal scale or $\kappa$.
Write $x=m(z)$, $q=1-z$, $m_z=4\pi w^3z^2$, and
$w(1-q)=q\,\omega(q)$, with $\omega>0$ smooth on a fixed
enlarged endpoint interval. Let
\begin{equation}
 r_\nu=\lambda u_\nu,\quad \nu=c,v,\qquad
 u_a=(1-\chi)u_v+\chi u_c,\quad
 \chi=\chi_{\rm m}(q/\sqrt\lambda).
 \label{mat:cutoff-1}
\end{equation}
The derivative support of $\chi$ lies in a shell
$c_1\sqrt\lambda\le q\le c_2\sqrt\lambda$, with fixed larger coordinate neighborhoods.

In this subsection, $D$ denotes the scale derivative and $E$
the core-edge operator $E=q\partial_q$ defined below; the normalized mass
variable is denoted $D_v$. The finite construction gives on this
shell
\begin{equation}
 \sum_{a+b\le23}|D^aE^b(u_c-u_v)|
 \le C\lambda^{(K+1)/2}(1+|\log\lambda|)^{p_{\log}}
 =:\varepsilon_\lambda,\quad
 D=\lambda\partial_\lambda\big|_z,\ E=q\partial_q.
 \label{mat:cutoff-2}
\end{equation}
This is the order-23 part of Subsection~\ref{mat:connection}.

Let $\mathcal F[u]=\lambda^2\{(\lambda u)_{\tau\tau}
+\mathcal M[\lambda u]\}$. Then
\begin{equation}
 \sum_{a+b\le21}
 \left|D^aE^b\left\{
 \mathcal F[u_a]-(1-\chi)\mathcal F[u_v]
                         -\chi\mathcal F[u_c]\right\}\right|
 \le C\lambda^{(K-1)/2}(1+|\log\lambda|)^{p_{\log}}.
 \label{mat:cutoff-3}
\end{equation}
For the cutoff error alone, the proof gives the stronger power
$K/2$; we shall use only \eqref{mat:cutoff-3}.

The preceding core and layer estimates read
\begin{equation}
 \begin{aligned}
 |W\mathcal F[u_c]|&\le
 C\lambda^{K+1}q^{-(K+1)}
       (1+|\log\lambda|+|\log q|)^{p_{\log}},\\
 &\hspace{1em}q\ge c\sqrt\lambda,\\
 |W\mathcal F[u_v]|&\le
 C\lambda^{K+1}Q^{K+1}
       (1+|\log\lambda|+\log Q)^{p_{\log}},\\
 &\hspace{1em}1\le Q=D_v^{1/4}\le C\lambda^{-1/2}.
 \end{aligned}
 \label{mat:cutoff-4}
\end{equation}
for compositions of material and conormal derivatives of total order
at most 21.
Here $D_v=(M-x)/\lambda^4$; it is distinct from the scale
derivative $D$. We also have the compact center/interior
$C^m$ residual bound $C\lambda^{K+1}(1+|\log\lambda|)^{p_{\log}}$
after $j$ material derivatives, $j+m\le21$, with the full
radial-vector center lift. The bounded endpoint itself was treated in Subsection~\ref{mat:endpoint}.

Keep the background $A=\lambda u_a$ uncorrected. Put
\begin{equation}
 \begin{gathered}
 \nu=\kappa\lambda^{-3/2},\qquad
 \mathcal R_u=\lambda^{-3}\mathcal F[u_a],\\
 F=-\mathcal R_u-\nu\lambda^{-1}\partial_s\mathcal M[A]
                       -\nu\lambda^{-1}\Lambda A_s.
 \end{gathered}
 \label{mat:cutoff-5}
\end{equation}
where $\Lambda\ge0$ is constant in time and
$\lambda_s=b\lambda$, $b^2=\beta^2+e\lambda$.
Then the source terms satisfy
\begin{equation}
  \sum_{j+m\le21}
 \mathcal S_m\!\left(\lambda^{-4}
                    \partial_s^j(\lambda^4F)\right)
 \le C\left[
 \lambda^{(K-5)/2}(1+|\log\lambda|)^{p_{\log}}
       +\kappa\lambda^{-9/2}
       +\kappa\Lambda\lambda^{-3/2}\right].
 \label{mat:cutoff-6}
\end{equation}
The source norm $\mathcal S_m$ and its cutoffs are those of
Definition~\ref{at:atlas}; their use in \eqref{mat:cutoff-6} is
detailed in Subsubsection~\ref{mat:cutoff-part-3}. The compatibility correction
changes only the initial data. Since the reference trajectory remains
$A$, the evolving source is still \eqref{mat:cutoff-5}.

The cutoff and weighted source estimates are proved in
Subsection~\ref{app:profile-cutoff}, with integration over every chart,
including the part of the growing chart where the core radius is used.

\subsection{Top-order material forcing}\label{mat:material}

\label{mat:material-part-1}

We now estimate one additional material derivative in the mass norm.
Keep the same positive leading profile, $\beta>0$, $e\ge0$,
finite family with $K\ge86$, $S=K+1$, and covering cutoffs.
We use the local bounds of Subsections~\ref{mat:units} and
\ref{mat:native}, and retain the scale derivative $D$ from
Subsection~\ref{mat:cutoff}:
\[
 D=\lambda\partial_\lambda\big|_z,\quad E=q\partial_q,
 \quad q=1-z,\quad Q=((M-m(z))/\lambda^4)^{1/4}.
\]
The overlap estimate for these radii is
\begin{equation}
 \sum_{a+b\le24}|D^aE^b(u_c-u_v)|
 \le C\lambda^{(K+1)/2}(1+|\log\lambda|)^{p_{\log}},
 \qquad q\sim\sqrt\lambda.
 \label{mat:material-1.1}
\end{equation}
On the joined trajectory $A=\lambda\{(1-\chi)u_v+\chi u_c\}$
let $\mathcal F=\lambda^2(A_{\tau\tau}+\mathcal M[A])$,
$\mathcal R_u=\lambda^{-3}\mathcal F$, $\partial_s=bD$
at fixed leading mass, and $b^2=\beta^2+e\lambda$. Then
\begin{equation}
 \sum_{j=0}^{22}\left\|\lambda^{-4}
                \partial_s^j(\lambda^4\mathcal R_u)\right\|_H
 \le C\lambda^{(K-7)/2}(1+|\log\lambda|)^{p_{\log}},
 \qquad H=L^2((0,M),\dd x).
 \label{mat:material-1.2}
\end{equation}
For the source with $\nu=\kappa\lambda^{-3/2}$ and
fixed $\Lambda\ge0$,
\[
 F=-\mathcal R_u-\nu\lambda^{-1}\partial_s\mathcal M[A]
                        -\nu\lambda^{-1}\Lambda A_s,
\]
we also have
\begin{equation}
 \sum_{j=0}^{22}\left\|\lambda^{-4}
                         \partial_s^j(\lambda^4F)\right\|_H
 \le C\{\lambda^{(K-7)/2}(1+|\log\lambda|)^{p_{\log}}
                  +\kappa\lambda^{-9/2}
                  +\kappa\Lambda\lambda^{-3/2}\}.
 \label{mat:material-1.3}
\end{equation}
The constants may depend on the selected family and on $\beta$;
they are uniform in the sufficiently small scale and subsequent
$\kappa$. These prescribed-source estimates are separate from
the energy inequality and the construction of a solution.

The overlap estimate \eqref{mat:material-1.1} is precisely
\eqref{mat:connection-4}. Its proof already supplies $D^aE^b$ bounds for $a+b\le24$, with leading derivatives available through
$\max\{H_n+30,H_n+n+25\}\le2S+50<2S+60$.
Thus the two derivatives in the force leave the 22 material
derivatives required here.

Subsection~\ref{app:profile-material} gives the additional cutoff
differentiation and the integration in the exact mass norm. The
latter has a different edge prefactor from the mixed source norm.
Thus \eqref{mat:material-1.2} is established separately, rather than
by extending the derivative range of \eqref{mat:cutoff-6}.

\subsection{Profile strains and the collected estimates}
\label{mat:exports}

The material strain estimates also distinguish constants uniform in the
leading family from those depending on the chosen profile. Put $A_i=\partial_s^iA$, impose
$e\lambda\le\beta^2$, and set
\[
 \vartheta=\frac{e\lambda}{b^2},\qquad
 \mathcal{E}_n=\lambda\partial_\lambda-\frac4n y\partial_y,\qquad
 r_i=\frac{A_i}{\lambda b^i},\qquad
 \gamma_i=-\frac{(A_i)_y}{\lambda^2yb^i}\quad(n=5,8).
\]
Exact differentiation at fixed mass gives
\begin{equation}
 \begin{split}
 r_{i+1}&=(\mathcal{E}_n+1+i\vartheta/2)r_i,\\
 \gamma_{i+1}&=(\mathcal{E}_n+2-8/n+i\vartheta/2)\gamma_i,\qquad
 \mathcal{E}_n\vartheta=\vartheta(1-\vartheta).
 \end{split}
 \label{mat:strain-recurrence}
\end{equation}
Here $[\mathcal{E}_n,\partial_y]=(4/n)\partial_y$,
$\mathcal{E}_ny=-(4/n)y$, and $0\le\vartheta\le1/2$.

\begin{lemma}
\label{mat:profile-strains}
For $0\le i\le23$,
\begin{equation}
 \left\|\frac{A_i}{A}\right\|_\infty+
 \left\|\frac{(A_i)_x}{A_x}\right\|_\infty
 \le \{C_{0,i}+\epsilon_A(\lambda)\}b^i,\qquad
 \epsilon_A(\lambda)=C_{K,e,\beta}
       \lambda^{1/2}(1+|\log\lambda|)^{p_{\log}}.
 \label{mat:strain-estimate}
\end{equation}
The constants $C_{0,i}$ depend only on the fixed compact interval
of leading profiles, their profile and layer bounds through order 24,
and the fixed cutoffs. They are chosen before $K$ and the positive
value of $\beta$. The constant $C_{K,e,\beta}$ may depend on the
selected finite family.
\end{lemma}

\begin{proof}
The leading layer satisfies
$P(\rho_0)=\varkappa D/\pi$ and
$Z_0=\varkappa^{-1}\sqrt{1+\rho_0^{2/3}}$.
On the bounded $y_5$ interval its regular quotients are smooth positive
functions of $y_5$ and $\varkappa$, with $\varkappa$ in a fixed
compact subset of $(0,\infty)$. On $Q\ge1$, implicit differentiation
gives uniform bounds for
\[
 \frac{\rho_0}{Q^3},\qquad
 \frac{Z_0}{Q},\qquad Q^3Z_0',
\]
and their required iterates under $E_Q$.
Together with the bounds of Subsection~\ref{mat:seed}, these control
the leading terms in the recurrence \eqref{mat:strain-recurrence}
through order 23, using at most 24 coefficient derivatives.
The corresponding leading denominators $r,\gamma$ are bounded below.

For the nonleading relative Jacobian terms, the coefficient bounds give
\[
 \lambda^nq^{-n}\le C\lambda^{n/2}\quad(q\ge c\sqrt\lambda),
 \qquad
 (\lambda Q)^k\le C\lambda^{k/2}\quad(Q\le C\lambda^{-1/2}),
\]
for $n\ge2$ and $k\ge1$, including the required iterates of $E_q$
and $E_Q$. The finite logarithmic factors are bounded by
$C_{K,e,\beta}(1+|\log\lambda|)^{p_{\log}}$. The bounded vacuum
interval is controlled by the finite sum of the $G_{n\ell}$.

On the matching shell, the leading asymptotics give
\[
 \lambda Z_0(D)=q+O(q^2+\lambda),\qquad
 q^{-1}O(q^2+\lambda)=O(\sqrt\lambda),
\]
with the same differentiated bounds. Since $\partial_z=-q^{-1}E_q$,
the second estimate also controls the leading Jacobian error from
the cutoff. Material derivatives of $\chi(q/\sqrt\lambda)$ are
bounded. In $(A_{23})_x$, a derivative on this cutoff leaves at most
23 material derivatives on the mismatch; an uncut term uses at most
24 coefficient derivatives. Subsections~\ref{mat:units} and
\ref{mat:material} supply these bounds. At the center we use the
full radial-vector lift, and on the interior the fixed coordinate.

The recurrence \eqref{mat:strain-recurrence} and the preceding
estimates now give the conclusion from the exact identities on the edge charts
\[
 \frac{A_i}{A}=b^i\frac{r_i}{r},\qquad
 \frac{(A_i)_x}{A_x}=b^i\frac{\gamma_i}{\gamma}.
\]
The leading quotients are uniformly bounded, and all correction
terms are bounded by $\epsilon_A(\lambda)$.
\end{proof}

\begin{proposition}
\label{mat:matched-family}
\label{exlim:upstream-profile}
\label{en:input-matched-profile}
For every $0<\beta\le\beta_{\rm match}$, $e\ge0$, and finite
$K\ge86$, $S=K+1$, the preceding construction gives a matched
radius $A(\lambda,x)$, $0<\lambda\le A_*$, with
$M=M_\beta>M_{\rm Ch}$. The upper bound for the scale may depend on this family.
All estimates below use this radius and the cutoffs of
Definition~\ref{at:atlas}.

\begin{enumerate}[label=(\roman*),leftmargin=2em]
\item \emph{Geometry.} The upper bound for the scale can also be decreased so that $e\lambda\le\beta^2$.
The finite sums are \eqref{mat:connection-3}: the auxiliary
$\Phi_S$ is not inserted in the core, whereas $Z_K$ is retained
in the layer. Their normalized radius and Jacobian factors, and the
endpoint quotients $G,Q$ of the full sum, are positive. The radius itself
vanishes at the center and is positive for $x>0$.
\item \emph{Finite regularity.} At every fixed positive scale, the
physical enthalpy and velocity are respectively $C^{58},C^{57}$
as in \eqref{mat:units-1}. On the entire covering edge, including
the portion where the core radius is used, the normalized factors
\[
 r=A/\lambda,\qquad\gamma=-A_y/(\lambda^2y),\qquad T_\gamma=y\gamma_y
\]
have bounded ordinary derivatives through order 21. The coefficients
$a_n,a_n^{-1},b_n,v_n$ obey \eqref{mat:units-4}.
At the center, the scalar and vector lifts have the corresponding
Cartesian regularity.

\item \emph{Profile derivatives.} The physical-time and material
derivatives of the profile satisfy
\[
 \begin{aligned}
 |\partial_\tau^i A|_{m,\infty,a}&\le
 C\lambda^{1-3i/2}\qquad(i+m\le27),\\
 |\partial_s^r\mathcal{M}[A]|_{m,\infty,a}
   +\|\partial_s^r\mathcal{M}[A]\|_{m,a}
 &\le C\lambda^{-2}\qquad(r+m\le25).
 \end{aligned}
\]
The norm $\|\cdot\|_{m,a}$ is the ordinary derivative norm of
Definition~\ref{at:atlas}, with the same cutoffs; it imposes no
positive-order compatible vacuum trace on $\mathcal M[A]$.
For $\lambda\le1$, comparison of the weights gives
\[
 \mathcal S_m(g)\le C\|g\|_{m,a},\qquad
 \|g\|_H\le C\|g\|_{0,a}.
\]
The first bound is used for the viscous source, and the second in
the global energy comparison.

In particular, material derivatives of $A$ through order 24 are
bounded. Equation~\eqref{mat:strain-estimate} separates the uniform
leading strain bounds from the error for the selected profile.
The orders 25 and 27 above apply to the profile and force; residual
decay is stated in Proposition~\ref{mat:matched-source}.
\end{enumerate}
\end{proposition}

\begin{proof}
The mass inequality follows from Subsection~\ref{mat:seed}. The
construction in Subsection~\ref{mat:connection}, including positivity
after the moving cutoff, proves (i). The physical-coordinate estimates
and normalized coefficient bounds of Subsection~\ref{mat:units}
give (ii). For (iii), use the profile and force derivative estimates of
Subsection~\ref{mat:native} and the additional material estimate in
Subsection~\ref{mat:material}. The strain estimate is
\eqref{mat:strain-estimate}.
\end{proof}

For the matched radius, set $H=L^2((0,M),\dd x)$, fix
$\Lambda\ge0$ independently of time, and put
\[
 \begin{gathered}
 \nu=\kappa\lambda^{-3/2},\qquad
 \mathcal{R}_u=\lambda^{-1}(A_{\tau\tau}+\mathcal{M}[A]),\\
 F=-\mathcal{R}_u-\nu\lambda^{-1}\partial_s\mathcal{M}[A]
                    -\nu\lambda^{-1}\Lambda A_s .
 \end{gathered}
\]

\begin{proposition}
\label{mat:matched-source}
\label{en:input-material-source}
Under the hypotheses of Proposition~\ref{mat:matched-family}, the
source norms with the cutoffs of Definition~\ref{at:atlas} satisfy
\begin{align}
 \sum_{j+m\le21}\mathcal{S}_m
       \bigl(\lambda^{-4}\partial_s^j(\lambda^4F)\bigr)
 &\le C\{\lambda^{(K-5)/2}(1+|\log\lambda|)^{p_{\log}}
          +\kappa\lambda^{-9/2}
          +\kappa\Lambda\lambda^{-3/2}\},
 \label{mat:strong-source-export}\\
 \sum_{j=0}^{22}\bigl\|
       \lambda^{-4}\partial_s^j(\lambda^4F)\bigr\|_H
 &\le C\{\lambda^{(K-7)/2}(1+|\log\lambda|)^{p_{\log}}
          +\kappa\lambda^{-9/2}
          +\kappa\Lambda\lambda^{-3/2}\}.
 \label{mat:material-source-export}
\end{align}
At the vacuum endpoint, \eqref{mat:endpoint-3} gives $C^{28}$ bounds
through material order 22. The source coefficients in the fixed
initial-radius coordinate satisfy \eqref{mat:endpoint-20} for
$K\ge2J^\sharp+86$.
\end{proposition}

\begin{proof}
The source estimates are \eqref{mat:cutoff-6} and
\eqref{mat:material-1.3}. Their long-chart factors are respectively
$\lambda^3$ and $\lambda^2$, so the residual powers are distinct.
For the viscous terms, finite Leibniz expansion and the profile-derivative
bounds give
\[
 \begin{aligned}
 \lambda^{-4}\partial_s^j
       (\lambda^4\nu\lambda^{-1}\partial_s\mathcal M[A])
     &=\kappa\lambda^{-5/2}
          \sum_{r=0}^j c_{jr}\partial_s^{r+1}\mathcal M[A],\\
 \lambda^{-4}\partial_s^j(\lambda^4\nu\lambda^{-1}\Lambda A_s)
     &=\kappa\Lambda\lambda^{-5/2}
          \sum_{r=0}^j c_{jr}\partial_s^{r+1}A.
 \end{aligned}
\]
Here
\[
 c_{jr}=\binom jr\lambda^{-3/2}\partial_s^{j-r}(\lambda^{3/2}),
 \qquad 0\le r\le j\le22.
\]
The identities $\lambda_s=b\lambda$ and $b_s=e\lambda/2$ bound
these coefficients uniformly. The two sums are
respectively $O(\lambda^{-2})$ and $O(\lambda)$ in the requested
norms. They therefore contribute $\kappa\lambda^{-9/2}$ and
$\kappa\Lambda\lambda^{-3/2}$. The endpoint estimates are
\eqref{mat:endpoint-3} and \eqref{mat:endpoint-20}.
\end{proof}

\begin{corollary}
\label{mat:matched-limits}
The profile in Proposition~\ref{mat:matched-family} has the following limits.
\begin{enumerate}[label=(\roman*),leftmargin=2em]
\item On every compact core interval,
\[
 A(\lambda,m(z))/\lambda\longrightarrow z,\qquad
 \lambda^3\rho[A](\lambda,m(z))\longrightarrow w_\beta(z)^3
\]
with compact spatial derivatives through order 23. More precisely,
the normalized differences tend to zero after $j$ material and
$m$ compact spatial derivatives whenever $j+m\le23$, with the
appropriate center scalar or vector lift. No assertion of convergence
after unweighted physical-time differentiation is made here.
\item On a fixed bounded interval in $y_5=D^{1/5}$, the regular
layer quantities satisfy
\[
 Z\longrightarrow Z_0,\qquad Z_D/Z_0'\longrightarrow1,
 \qquad \rho[A]/y_5^3\longrightarrow\rho_0/y_5^3
 \quad\hbox{in }C^{23}.
\]
The quotients in this display have their continuous values at zero.
At the endpoint the finite sum gives the sharper statements
\begin{equation}
 \begin{split}
 Z(\lambda,0)&=\varkappa_\delta^{-1}
                   +O(\lambda(1+|\log\lambda|)^{p_{\log}}),\\
 A(\lambda,M)&=\lambda-\lambda^2\varkappa_\delta^{-1}
                   +O(\lambda^3(1+|\log\lambda|)^{p_{\log}}).
 \end{split}
 \label{mat:boundary-profile-export}
\end{equation}
\end{enumerate}
These are asymptotics of the prescribed profile.
\end{corollary}

\begin{proof}
For the core limits, each nonleading term is a finite linear
combination of $\lambda^n(\log\lambda)^\ell\phi_{n\ell}$, $n\ge2$.
Writing $L=\log\lambda$, one has
\[
 (\lambda\partial_\lambda)^j
      \{\lambda^n(\log\lambda)^\ell\}
 =\lambda^n(n+\partial_L)^jL^\ell
 =O\!\left(\lambda^n(1+|\log\lambda|)^\ell\right).
\]
On a compact core, spatial derivatives act only on the fixed
coefficients. Since $\partial_s=b\lambda\partial_\lambda$ there
and the finitely many iterates of $\lambda\partial_\lambda$ applied to $b$ are bounded, every
stated mixed derivative of the nonleading sum tends to zero.
The same applies to the regular center lifts. The exact identity
\[
 \lambda^3\rho[A]
   =\frac{w^3}{\psi^2(\psi+z\psi_z)},\qquad
 \psi\longrightarrow1,\quad \psi+z\psi_z\longrightarrow1,
\]
then gives the density limit by differentiation of the positive
reciprocal factors.

On a fixed bounded layer, the finite expansions for the regular
factors $G,Q$, together with \eqref{mat:endpoint-5}, give the stated
$C^{23}$ limits. At $D=0$, \eqref{mat:connection-10} and
$\rho_0(0)=0$ give $Z_0(0)=\varkappa_\delta^{-1}$. Thus
\[
 \begin{aligned}
 Z(\lambda,0)-\varkappa_\delta^{-1}
     &=O\!\left(\lambda(1+|\log\lambda|)^{p_{\log}}\right),\\
 A(\lambda,M)-\lambda+\lambda^2\varkappa_\delta^{-1}
     &=-\lambda^2\{Z(\lambda,0)-\varkappa_\delta^{-1}\}
       =O\!\left(\lambda^3(1+|\log\lambda|)^{p_{\log}}\right).
 \end{aligned}
\]
This proves \eqref{mat:boundary-profile-export}.
\end{proof}

For $K\ge2J+86$, $J\ge3$, and
$0<\lambda_*\le\lambda\le A_*$, the two bounds
\[
 \kappa\le c\lambda_*^{J+5},\qquad
 \kappa\Lambda\le c\lambda_*^{J+2}
\]
give
\[
 \kappa\lambda^{-9/2}\le c\lambda^{J+1/2},\qquad
 \kappa\Lambda\lambda^{-3/2}\le c\lambda^{J+1/2}.
\]
The residual exponents satisfy
\[
 \frac{K-7}{2}-\left(J+\frac12\right)\ge39,\qquad
 \frac{K-5}{2}-\left(J+\frac12\right)\ge40.
\]
These positive margins absorb every finite logarithmic factor with
any prescribed small coefficient after decreasing $A_*$.
For $K=2J+3512$, the material margin is 1752.
The leading constants are unchanged by this final scale choice.

\Needspace{7\baselineskip}
\noindent\textbf{Where the profile estimates are used.}\phantomsection\label{mat:consumer-map}
We summarize the estimates required in the subsequent sections.
Throughout, $A$ is the matched radius constructed above.
\begin{enumerate}[label=(\roman*),leftmargin=2.4em]
\item \emph{Compatible initial data.} Section~\ref{loc:section} uses
fixed-scale physical regularity, the local radius and force derivatives, and
the ordinary range of endpoint derivatives \eqref{mat:endpoint-20}. Its finite
correction changes initial data, while $A$ and its evolving residual
remain fixed.
\item \emph{Spatial recovery.} Section~\ref{rec:chapter} uses the
positive principal coefficients and the coefficient bounds in
\eqref{mat:units-4}, together with the mixed source
\eqref{mat:strong-source-export}. The range of source derivative orders $j+m\le21$ produces
recovery estimates at orders $j+(m+2)\le23$.
\item \emph{Energy estimates} in Section~\ref{en:section} use the
leading strain constants \eqref{mat:strain-estimate} and the separate
material source \eqref{mat:material-source-export} through order $22$.
The latter supplies the forcing paired with the twenty-third time derivative in the kinetic energy; it is
not obtained by extending the range of mixed derivative orders.
\item \emph{Boundary-layer limits} in Subsection~\ref{exlim:actual-layer}
combine the profile asymptotics on a fixed rescaled layer with continuity
of the nonlinear correction in the lower-order weighted spaces. The ordinary quotient
identities pass this convergence to the vacuum traces; positive lower
bounds permit the reciprocal density factors.
\item \emph{Global physical energy} in
Proposition~\ref{exlim:clock-energy} uses the finite core coefficient
growth, the radius and density factors on the entire coordinate cover, the strain bound, and
the global $\lambda^{-2}$ force bound. These control the entire star,
including the joining shell, when comparing the solution with the
homologous reference. Fixed-core and fixed-layer convergence alone
would not justify that energy calculation.
\end{enumerate}

\section{Weighted estimates and the linearized operator}
\label{at:calculus-section}

The profile bounds give a uniformly positive principal coefficient
on each chart. We first prove an elliptic estimate for the profile
operator and then estimate the change in the operator under a radius
perturbation. Since its coefficients depend on two radius derivatives,
the product estimates must preserve the derivative count in the
spatial induction.

In \eqref{at:current-high}, derivatives on the argument have a small
coefficient controlled by a low norm of the radius perturbation.
Derivatives on that perturbation instead give a high norm times a low
norm of the argument; Section~\ref{rec:chapter} absorbs this term.
The estimates hold on the closed strong domain. The distributional
limit is placed in that domain in Subsection~\ref{exlim:realization}.

\subsection{One-sided multiplication and composition}
\label{at:weighted-calculus}

At the vacuum we use the one-sided derivative norm
\begin{equation}
 \|f\|_{\mathcal B_5^m}^2
     =\sum_{q=0}^m\int_0^L\xi^4|\partial_\xi^qf|^2\dd\xi.
 \label{in:ordinary-weight}
\end{equation}
The measure $\xi^4\dd\xi$ agrees with the mass measure up to a
constant. The norm uses one-sided derivatives, with no evenness
condition at the vacuum. For example, a cutoff of $f(\xi)=\xi^3$
satisfies $f'(0)=0$ and has finite weighted derivatives of every
finite order, whereas $f(|X|)=|X|^3$ need not have the corresponding
high Sobolev regularity on a five-dimensional ball. We therefore
prove the product estimates by extending to a cone. The radial lift
used in the local parabolic construction applies only at low order.

We first prove the estimates in one-sided spaces with neither parity
nor a trace condition at zero. The quotient estimate in
Lemma~\ref{at:natural-quotient} then gives the bounds needed on the
compatible domain.

For $d\in\{2,\ldots,8\}$, $L>0$, and an integer $m\ge0$, define
$\mathcal B_d^m(0,L)$ as the completion of $C^\infty([0,L])$ for
\[
 \|F\|_{\mathcal B_d^m(0,L)}^2
 =\sum_{q=0}^m\int_0^L|F^{(q)}(y)|^2y^{d-1}\dd y.
\]
No parity condition is imposed at zero.

\begin{lemma}
\label{at:one-sided-multiplication}
Let $d\in\{2,\ldots,8\}$, $L>0$, and let $m\ge6$ be an integer.
Then
\begin{align}
 \|FG\|_{\mathcal B_d^m}
 &\le C_{m,d,L}\bigl(
   \|F\|_{\mathcal B_d^m}\|G\|_{\mathcal B_d^6}
  +\|G\|_{\mathcal B_d^m}\|F\|_{\mathcal B_d^6}\bigr),
 \label{at:global-tame-product}\\
 \|F\|_{W^{1,\infty}(0,L)}
 &\le C_{d,L}\|F\|_{\mathcal B_d^6}.
 \label{at:global-low-embedding}
\end{align}
The same estimates hold for the
derivative-sum norm $\mathcal C_d^m(c,L)$ on $(c,L)$, with constants
depending on $m,d,c$ but independent of $L\ge2c$.
\end{lemma}

The cone extension and the uniform annular argument are given in
Subsection~\ref{app:one-sided-calculus}.

Write $\|\cdot\|_r$ for either $\|\cdot\|_{\mathcal B_d^r(0,L)}$
or $\|\cdot\|_{\mathcal C_d^r(c,L)}$.

\begin{lemma}
\label{at:low-rung-multiplier}
For each integer $0\le r\le6$,
\[
 \|GF\|_r\le C_{r,d,L}\|G\|_6\|F\|_r.
\]
In the annular case the constant depends on $r,d,c$ and is independent
of $L\ge2c$.
\end{lemma}

The estimate follows by interpolation between orders zero and six;
see Subsection~\ref{app:one-sided-calculus}.

In either of these spaces, let $G:\mathbb R^a\to\mathbb R^b$ be
smooth with bounded derivatives of every order used below.

\begin{lemma}
\label{at:coefficient-composition}
Fix $m\ge6$ and $M<\infty$.
\begin{enumerate}[label=(\roman*),leftmargin=2em]
\item \emph{Composition and differences.} The estimates are
\begin{align}
 \|G(F)-G(0)\|_m
 &\le C_{m,G,M}\|F\|_m, \qquad \|F\|_6\le M,
 \label{at:composition}\\
 \|G(F)-G(\widetilde F)\|_m
 &\le C_{m,G,M}\bigl\{\|F-\widetilde F\|_m
  +(\|F\|_m+\|\widetilde F\|_m)
                     \|F-\widetilde F\|_6\bigr\},
 \label{at:composition-difference}
\end{align}
where the second assertion assumes
$\|F\|_6+\|\widetilde F\|_6\le M$.
\Needspace{7\baselineskip}
\item \emph{Higher differentials.} For $p\ge1$ and $\|F\|_6\le M$,
\begin{align}
 \|D^pG(F)[h_1,\ldots,h_p]\|_m
 \le C_{m,p,G,M}\biggl\{
 &\sum_{i=1}^p\|h_i\|_m\prod_{j\ne i}\|h_j\|_6
 \notag\\*
 &+(1+\|F\|_m)\prod_{i=1}^p\|h_i\|_6\biggr\}.
 \label{at:coefficient-high-branch}
\end{align}
The second line accounts for derivatives falling on the coefficient
$D^pG(F)$ rather than on a direction $h_i$.
\end{enumerate}
The constants have the same
geometric dependence as in Lemma~\ref{at:one-sided-multiplication};
in particular they are independent of the growing upper endpoint
in the annular case.

If $G$ is defined only on an open set, these conclusions apply when
the ranges in question lie in a fixed compact subset of that set,
by first extending $G$ smoothly from a neighborhood of that compact
set.
\end{lemma}

A dyadic composition estimate, followed by the product bound, gives
these inequalities. The proof is in
Subsection~\ref{app:one-sided-calculus}.

\subsection{Localized product estimates}
\label{at:localization}

To estimate products in the graph norm, we retain its defining cutoff
throughout the interpolation. Since the cutoff is a high power, its
remaining powers can be distributed among the differentiated factors.
This gives a high norm of only one factor in each product term. On
the growing eighth-root chart, we also keep the volume factor explicit.

\begin{lemma}
\label{at:power-interpolation}
Let $\dd\mu$ be Lebesgue measure on a fixed Euclidean chart, or
$y^7\dd y$ on $c<y<L$, where $c>0$ is fixed. Suppose $0\le\theta\le1$
vanishes at the artificial boundary and has the prescribed bounded
derivatives. For $2\le m\le Q$, put $\chi=\theta^Q$ and
$\mu_0=\int \dd\mu$. Then
\begin{align}
 \|\theta^{Qk/m}D^ku\|_{L^{2m/k}(\dd\mu)}
 &\le C\|u\|_\infty^{1-k/m}
   \bigl(\|\chi D^mu\|_2+\mu_0^{1/2}\|u\|_\infty\bigr)^{k/m},
 \quad 1\le k\le m,
 \label{at:joint-interpolation}\\
 \sum_{k<m}\|\theta^{Q-m+k}D^ku\|_2
 &\le\epsilon\|\chi D^mu\|_2+C_\epsilon\|u\|_2,
 \quad\epsilon>0.
 \label{at:linear-interpolation}
\end{align}
The derivatives are the full ordered tensors. On the growing annulus
the constants are independent of $L$; its volume occurs only where
displayed.
\end{lemma}

Integration by parts preserves the common cutoff power and gives
a finite recurrence between derivative orders. The proof, with
constants independent of the annular length, is in
Subsection~\ref{app:localized-products}.

Since $|D^p\chi|\le C\theta^{Q-p}$, the last lemma gives
\begin{equation}
 \sum_{p=1}^q\|(D^p\chi)D^{q-p}u\|_2
 \le\epsilon\|\chi D^qu\|_2+C_\epsilon\|u\|_2,
 \qquad 2\le q\le23.
 \label{at:cutoff-commutators}
\end{equation}
In particular, for $l=0,1,2$ and $m+l\le23$,
\begin{equation}
 \|\chi D^lu\|_{H^m(\dd\mu)}
 \le C\{\|\chi u\|_{H^{m+l}(\dd\mu)}+\|u\|_{L^2(\dd\mu)}\}.
 \label{at:localized-jet}
\end{equation}
The unlocalized remainder has order zero only.

\begin{lemma}
\label{at:same-cutoff-product}
In the preceding charts, fix $1\le m\le23$.
\begin{enumerate}[label=(\roman*),leftmargin=2em]
\item \emph{Products.} One has
\begin{equation}
 \|\chi uv\|_{H^m(\dd\mu)}
 \le C_m\left\{
   \|\chi u\|_{H^m}\|v\|_\infty
  +\|\chi v\|_{H^m}\|u\|_\infty
  +\mu_0^{1/2}\|u\|_\infty\|v\|_\infty\right\}.
 \label{at:cutoff-product}
\end{equation}
The same estimate holds for a finite product, with one high localized
factor in each summand and one volume remainder.
\item \emph{Coefficient differences.} Suppose $U_a$ has
bounded ordinary derivatives through $m$, and the segments
$U_a+tu$, $0\le t\le1$, lie in a fixed compact state set.
For the smooth coefficient $G(x,U)$, assume that
\[
 \partial_x^\alpha D_U^{j+1}G(x,U),\qquad |\alpha|+j\le m,
\]
are bounded on the entire chart and that state set. On a growing
annulus these bounds and the profile bounds are required to be
independent of $L$. Then
\begin{equation}
 \begin{split}
 \|\chi\{G(x,U_a+u)-G(x,U_a)\}v\|_{H^m}
 \le C_m\{&\|\chi u\|_{H^m}\|v\|_\infty
       +\|u\|_\infty\|\chi v\|_{H^m}\\
       &+\mu_0^{1/2}\|u\|_\infty\|v\|_\infty\}.
 \end{split}
 \label{at:cutoff-coefficient}
\end{equation}
The constant depends on these mixed coefficient bounds, the finite
profile bounds, the compact state set and the cutoff bounds. It
uses no unlocalized high norm of $u$.
\end{enumerate}
\end{lemma}

The proof distributes the remaining cutoff powers among the
differentiated factors and applies
Lemma~\ref{at:power-interpolation}; see
Subsection~\ref{app:localized-products}.

A radius variation $f$ changes the normalized edge Jacobian through
$f_y/y$. We first estimate this quotient in the compatible weighted
completion. We then record its endpoint values under classical
one-sided regularity. The same integral identity gives both statements,
but only the second assigns a value at the endpoint.

\begin{lemma}
\label{at:natural-quotient}
For a compatible one-sided function and a fixed vacuum cutoff,
\begin{equation}
 \|\chi f_y/y\|_{B_5^m}
 \le C_m\{\|\chi f\|_{B_5^{m+2}}+\|f\|_{B_5^0}\},
 \qquad 0\le m\le21.
 \label{at:closed-quotient}
\end{equation}
The quotient extends as a closed bounded map on the corresponding
compatible completion. The assertion prescribes no value of
$f_y/y$ at zero.
\end{lemma}

Dilation estimates and Hardy's inequalities give the weighted
quotient bound. The proof in Subsection~\ref{app:vacuum-quotients}
also identifies its extension to the compatible completion.

\begin{lemma}
\label{at:regular-factor-calculus}
Let $k\ge0$ be an integer. On a fixed interval $[0,L]$, the maps
\begin{align}
 \mathcal Qg(y)&=\frac{g'(y)}y=\int_0^1g''(ty)\dd t,
 \label{exlim:ordinary-quotients}\\
 \mathcal Gg(y)&=\frac{g(y)-g(0)}{y^2}
                    =\int_0^1(1-t)g''(ty)\dd t.
 \label{exlim:layer-division}
\end{align}
are bounded linear maps from
$\{g\in C^{k+2}([0,L]):g'(0)=0\}$ to $C^k([0,L])$.
The endpoint values are $\mathcal Qg(0)=g''(0)$ and
$\mathcal Gg(0)=g''(0)/2$.
On any bounded set in $C^k([0,L])$ with a common strictly positive
lower bound, taking a reciprocal is locally Lipschitz in $C^k$.
These assertions use ordinary one-sided derivatives and impose no
parity condition at zero.
\end{lemma}

The integral formulas and Leibniz's rule give the stated bounds;
the proof is in Subsection~\ref{app:vacuum-quotients}.

\subsection{Coefficients of the linearized operator}
\label{at:current-coefficients}

Let $R=A+\lambda^4Z$ and $L_R=\lambda^3D\mathcal M[R]$.
The coefficients are affine in the second derivatives of the radius.
We compute them on compatible smooth functions and then use the
closed quotient of Lemma~\ref{at:natural-quotient} to pass to the
completed domains. On the edge, the normalized Jacobian $\gamma$
and the quantity $T=y\gamma_y$ give the pressure coefficients
without an extra unweighted derivative of $\gamma$. At the center,
the full vector formulation retains the angular derivatives.

On the fixed interior put $r=R/\lambda$, $p=r_x$, and
$c=(4\pi r^2p)^{-1}$. Define
\begin{equation}
 a_\lambda(r,p)=\frac{\lambda P'(\lambda^{-3}c)}{p^2},
 \qquad
 b_\lambda(r,p,x)=-\frac{2\lambda P'(\lambda^{-3}c)}r+\frac{x}{r^2}.
 \label{at:interior-atoms}
\end{equation}
The force and its derivative are exactly
\begin{align}
 \mathcal M[\lambda r]
   &=\lambda^{-2}\{-a_\lambda(r,r_x)r_{xx}+b_\lambda(r,r_x,x)\},
 \label{at:interior-force}\\
 L_Rf&=-a_\lambda f_{xx}
       +(-\partial_pa_\lambda\,r_{xx}+\partial_pb_\lambda)f_x
       +(-\partial_ra_\lambda\,r_{xx}+\partial_rb_\lambda)f.
 \label{at:interior-hessian}
\end{align}
These coefficients are smooth in $(r,r_x)$ and affine in $r_{xx}$.
Their constitutive factor is
\begin{equation}
 \lambda P'(\lambda^{-3}c)
       =\frac{4c^{1/3}}{3\sqrt{1+\lambda^2c^{-2/3}}}.
 \label{at:scaled-pressure}
\end{equation}
It extends smoothly to $\lambda=0$ for $c$ in any compact subinterval
of the positive half-line.

On an edge, with $n=5$ or $8$, write
\begin{equation}
 r=R/\lambda,\qquad
 \gamma=-R_y/(\lambda^2y),\qquad T=y\gamma_y,\qquad
 \rho=\frac{ny^{n-2}}{4\pi r^2\gamma}.
 \label{at:edge-units}
\end{equation}
Let $D_nf=f_{yy}+(n-1)f_y/y$ and
$\varkappa(\rho)=\rho P''(\rho)/P'(\rho)$. Direct differentiation
of the exact force gives
\begin{equation}
 \lambda L_Rf=-A_nD_nf-B_nf_y/y+V_nf,
 \label{at:edge-hessian}
\end{equation}
where, writing $G_n(y,r,\gamma)=P'(\rho)/(\gamma^2y^2)$,
\begin{align}
 A_n&=G_n(y,r,\gamma),\qquad
 B_n=y(G_n)_y-\lambda y^2\gamma(G_n)_r+(G_n)_\gamma T,
 \label{at:edge-principal}\\
 V_n&=\lambda\left\{
   \frac{2\lambda\gamma^2y^2A_n}{r^2}
   +\frac{2\varkappa(\rho)\gamma A_n}{r}
        \left(n-2+\frac{2\lambda y^2\gamma}{r}-\frac T\gamma\right)
   -\frac{2x}{r^3}\right\}.
 \label{at:edge-potential}
\end{align}
In particular $B_n=y\partial_yA_n$ along the current radius.
Before multiplication by the edge factor $\lambda$, the zeroth
coefficient is
\[
 \frac{2\lambda^3P'}{R^2}
  -8\pi\lambda^3R\rho P''\rho_x-\frac{2\lambda^3x}{R^3},
 \qquad
 \frac{\rho_y}{\rho}
  =\frac{n-2+2\lambda y^2\gamma/r-T/\gamma}{y}.
\]
Substitution yields \eqref{at:edge-potential}.
The principal coefficients can be written as
\begin{align}
 A_5&=\frac{4C_5^{2/3}}{3\gamma^2\sqrt{1+C_5^{2/3}y^2}},
       &C_5&=\frac5{4\pi r^2\gamma},\label{at:fifth-atom}\\
 A_8&=\frac{4C_8^{1/3}}{3\gamma^2\sqrt{1+C_8^{-2/3}y^{-4}}},
       &C_8&=\frac2{\pi r^2\gamma}.\label{at:eighth-atom}
\end{align}
Their mixed derivatives in $y$ and the positive factors $(r,\gamma)$
are bounded on the prescribed ranges. On the long chart, the bounds
$y\ge c>0$, $y\le C_e\lambda^{-1/2}$ give
\[
 |\partial_y^j y^{-4}|
 +|\partial_y^j(\lambda y^2)|
 +|\partial_y^j(M-\lambda^4y^8)|\le C_j,\qquad j\ge0.
\]
Here the positive factors range over the prescribed compact set and
$0<\lambda\le A_0\le1$. The remaining coefficients in
\eqref{at:edge-principal}--\eqref{at:edge-potential} are smooth
combinations of these bounded arguments, affine in $T$. Treating
$T$ as a separate factor gives exactly the mixed coefficient bounds
required in Lemma~\ref{at:same-cutoff-product}.

At the center use the normalized full vector
$\Phi(X)=\lambda^{-1}R(|X|^3)X/|X|=r_c(X)X$. Put
\[
 J=D\Phi,\qquad B=J^{-T},\qquad c=\frac3{4\pi\det J},\qquad
 p_\lambda=\lambda P'(\lambda^{-3}c),\qquad
 h_\lambda=\lambda h(\lambda^{-3}c).
\]
Here $h_\lambda=4(\sqrt{\lambda^2+c^{2/3}}-\lambda)$.
For scalars and vectors define
$D_i^B=B_{ik}\partial_k$,
$\operatorname{div}_B u=D_i^Bu_i$, and
$(\operatorname{curl}_Bu)_i=\varepsilon_{ijk}D_j^Bu_k$.
The scaled radial Hessian is the restriction of
\begin{equation}
 \mathcal L_\Phi u
  =-B(Du)^TB\nabla h_\lambda
    -B\nabla(p_\lambda\operatorname{div}_Bu)-2r_c^{-3}u.
 \label{at:center-variation}
\end{equation}
Indeed the scaled pressure force is $B\nabla h_\lambda$,
$\delta B=-B(Du)^TB$, and
$\delta c=-c\operatorname{div}_Bu$.
Add $p_\lambda\operatorname{curl}_B\operatorname{curl}_Bu$
to \eqref{at:center-variation}. For radial $u$, the symmetric
matrices $Du$ and $B$ commute, so $\operatorname{curl}_Bu=0$.
Thus this addition changes no radial equation.
The complete operator has the scalar principal part
\begin{equation}
 (\widehat{\mathcal L}_\Phi u)_i
 =-p_\lambda(B^TB)_{k\ell}\partial_{k\ell}u_i
                        +C_{ia\ell}\partial_\ell u_a-2r_c^{-3}u_i,
 \label{at:center-completion}
\end{equation}
with
\begin{align}
 C_{ia\ell}={}&-B_{i\ell}(B\nabla h_\lambda)_a
 -(B_{ik}\partial_kp_\lambda)B_{a\ell}
 -p_\lambda B_{ik}\partial_kB_{a\ell}\notag\\
 &+p_\lambda B_{ak}\partial_kB_{i\ell}
 -p_\lambda\delta_{ia}B_{bk}\partial_kB_{b\ell}.
 \label{at:center-lower-tensor}
\end{align}
This follows by expanding the two curls and cancelling their
mixed component principal derivatives. The identities
\[
 \partial_kc=-c\operatorname{tr}(J^{-1}\partial_kJ),\qquad
 \partial_kB=-B(\partial_kJ)^TB
\]
show that the first-order tensor is smooth in $D\Phi$ and affine
in $D^2\Phi$. The principal matrix in \eqref{at:center-completion} is positive
on the prescribed deformation range. The radial factor can be
recovered without division by $|X|$:
\begin{equation}
 r_c(X)=\int_0^1t^2\operatorname{div}\Phi(tX)\dd t.
 \label{at:center-quotient}
\end{equation}
On a center ball $B$, Minkowski's inequality and dilation give
\[
 \begin{aligned}
 \|D^kr_c\|_{L^2(B)}
 &\le\int_0^1t^{k+2-3/2}
                 \|D^k\operatorname{div}\Phi\|_{L^2(tB)}\dd t\\
 &\le\frac1{k+3/2}\|D^k\operatorname{div}\Phi\|_{L^2(B)}.
 \end{aligned}
\]
The cutoff is one near the center. On its fixed outer overlap,
$r_c$ is a smooth multiple of the radius, so the same localized
bounds follow by the product rule.

The perturbations of the radius derivatives and quotients appearing
in these coefficients are, in the interior and at an edge respectively,
\begin{equation}
 u_i=\lambda^3(Z,Z_x,Z_{xx}),\qquad
 u_e=(\lambda^3Z,-\lambda^2Z_y/y,
                         -\lambda^2(Z_{yy}-Z_y/y)).
 \label{at:increment-jets}
\end{equation}
At the center they are $\lambda^3D U_Z$, $\lambda^3D^2U_Z$
and the radial integral in \eqref{at:center-quotient}.
Thus the coefficient estimates through order $q-2$ use radius
derivatives only through order $q$, with the prescribed localization.
The identity
\begin{equation}
 y\gamma^{(q-1)}
    =-\lambda^{-2}R^{(q)}-(q-1)\gamma^{(q-2)}
 \label{at:weighted-high-atom}
\end{equation}
retains its factor $y$; it is not an estimate for an additional
unweighted quotient derivative.

\subsection{The elliptic estimate for the profile operator}
\label{at:inverse-section}

For the model normal operator the derivative gain is explicit.
If $f''+4\xi^{-1}f'=g$ and the singular flux is excluded, then
\[
 \frac{f'(\xi)}{\xi}=\int_0^1 t^4g(t\xi)\dd t,\qquad
 f''=g-4f'/\xi.
\]
Differentiating the first identity $r$ times and using
Minkowski's inequality gives
\[
 \left\|\partial_\xi^r(f'/\xi)\right\|_{L^2(\xi^4\dd\xi)}
 \le \frac{1}{r+5/2}
          \|\partial_\xi^rg\|_{L^2(\xi^4\dd\xi)}.
\]
The homogeneous solution with $|f'|\sim\xi^{-4}$ is excluded by
the pressure form domain. The integral formula therefore selects
zero form flux and controls $f'/\xi$. For the current Hessian, we
integrate with the positive coefficient kept inside the flux.
Weighted product estimates then control its differentiated
coefficients and yield the same gain of two one-sided derivatives.

Proposition~\ref{exlim:upstream-profile} supplies the positive
normalized deformation and the finite radius bounds on the compact
charts. On the edges we use the following coefficient bounds, with
$m=q-2$:
\begin{equation}
 A_n,\ A_n^{-1},\ B_n=y\partial_yA_n,\ V_n
       \in W^{m,\infty},\qquad A_n\ge a_*>0,
 \qquad 2\le q\le23.
 \label{at:profile-budget}
\end{equation}
The bounds hold for the assembled profile on the entire cover,
with constants depending on the selected finite profile. In particular,
the estimate controls $y\partial_yA_n$ through order $m$;
an unweighted derivative of $A_n$ of order $m+1$ is not required.

\begin{proposition}
\label{at:profile-inverse}
Under \eqref{at:profile-budget} and the compact-chart profile
bounds through radius order 23, compatible functions satisfy
\begin{equation}
 \|f\|_{X^q}\le C_q\{
       \mathcal S_{q-2}(L_Af)+\|f\|_{X^{q-2}}\},
 \qquad 2\le q\le23 .
 \label{at:profile-inverse-bound}
\end{equation}
The constant is uniform in the scale and the growing chart
length. It uses only the displayed finite coefficient bounds.
The assertion extends to the closure of the compatible core in
$\|f\|_{X^q}+\mathcal S_{q-2}(L_Af)$.
\end{proposition}

\eqref{at:profile-inverse-bound} is an a priori estimate with a lower
norm on the right; it does not assert invertibility or surjectivity.
Lemma~\ref{loc:finite-strong-graph} places the functions used below
in its closed strong domain.

\begin{proof}
\emph{Vacuum.}
Write $Pf=-y^{-4}(y^4Af_y)'+Vf=F$ and $g=F-Vf$.
The zero natural flux gives
\[
 f_y=-A^{-1}Tg,\qquad
 Tg=ySg,\qquad Sg=\int_0^1t^4g(ty)\dd t.
\]
Minkowski's inequality and dilation yield
\[
 \begin{gathered}
 \|\partial_y^jSg\|_{B_5^0}
 \le(j+5/2)^{-1}\|\partial_y^jg\|_{B_5^0},\\
 (Tg)'=g-4Sg,\qquad
 \|Tg\|_{B_5^{m+1}}\le C_m\|g\|_{B_5^m}.
 \end{gathered}
\]
Put $a=A^{-1}$. The highest coefficient derivative occurs as
\[
 a^{(m+1)}Tg=(ya^{(m+1)})Sg,\qquad
 ya^{(m+1)}=\partial_y^m(-B/A^2)-m\partial_y^m(A^{-1}).
\]
Thus \eqref{at:profile-budget} controls this term, including at
$m=0$, and Leibniz's rule gives
\[
 \|\partial_y^{m+1}(aTg)\|_{B_5^0}
 \le C_m\{\|Tg\|_{B_5^{m+1}}+\|Sg\|_{B_5^0}\}
 \le C_m\|g\|_{B_5^m}.
\]
Since $g=F-Vf$,
\begin{equation}
 \|f\|_{B_5^{m+2}}
       \le C\{\|F\|_{B_5^m}+\|f\|_{B_5^m}\}.
 \label{at:vacuum-inverse}
\end{equation}
For $v=\chi_5f$,
\[
 Pv=\chi_5F+[P,\chi_5]f,\qquad
 [P,\chi_5]f=-2A\chi_5'f_y
 -\left(A\chi_5''+\left(A'+\frac{4A}{y}\right)\chi_5'\right)f.
\]
The cutoff derivatives lie away from zero and retain powers at
least $Q-2-m$ after $m$ differentiations. By
\eqref{at:cutoff-commutators} and $|\dd x|=5\lambda^4y^4\dd y$,
\[
 \lambda^2\|[P,\chi_5]f\|_{B_5^m}
 \le\epsilon\lambda^2\|\chi_5f\|_{B_5^{m+2}}
                         +C_\epsilon\|f\|_H.
\]
On this transition interval the identity for $ya^{(m+1)}$ also
bounds $a^{(m+1)}$. Apply \eqref{at:vacuum-inverse} to $v$ and
absorb the first term. No value of $f(0)$ is prescribed.

\emph{Long chart and interior.}
Here
\[
 Pf=-Af_{yy}-\frac{7A+B}{y}f_y+Vf,
 \qquad c<y<C_e\lambda^{-1/2}.
\]
The differentiated equation is
\[
 \begin{aligned}
 Af^{(m+2)}
 ={}&-\partial_y^mF+\partial_y^m(Vf)
       -\sum_{k=1}^m\binom mk A^{(k)}f^{(m+2-k)}\\
 &-\partial_y^m\left(\frac{7A+B}{y}f_y\right).
 \end{aligned}
\]
The coefficient bounds \eqref{at:profile-budget} and $y\ge c>0$
give uniform bounds for $1/y$ and $\partial_y\log(y^7)=7/y$.
Apply \eqref{at:linear-interpolation} to the lower derivatives and
\eqref{at:cutoff-commutators} to $\chi_8$. We obtain
\[
 \lambda^2\|\chi_8f\|_{C_8^{m+2}}
 \le C\lambda^2\|\chi_8F\|_{C_8^m}
       +\epsilon\lambda^2\|\chi_8f\|_{C_8^{m+2}}
       +C_\epsilon\|f\|_H.
\]
Absorb the middle term. The identity
$|\dd x|=8\lambda^4y^7\dd y$ controls the order-zero remainder
uniformly in the chart length. The interior estimate follows in
the same way, with Lebesgue measure and $L_A$ in place of $\lambda L_A$.

\emph{Center.}
Use the full vector operator \eqref{at:center-completion}.
On finitely many balls, freeze its principal matrix with oscillation
at most $\epsilon$. The constant-coefficient Fourier estimate gives
\[
 \begin{gathered}
 \|D^2u\|_2
 \le C\{\|\widehat{\mathcal L}_\Phi u\|_2
       +\epsilon\|D^2u\|_2+\|\nabla u\|_2+\|u\|_2\},\\
 \|\nabla u\|_2\le\epsilon\|D^2u\|_2+C_\epsilon\|u\|_2.
 \end{gathered}
\]
Absorption gives the $H^2$ estimate. Commuting $m$ Cartesian
derivatives and repeating it gives the $H^{m+2}/H^m$ estimate.
The first-order coefficients use at most two radius derivatives,
so order $m+2=q$ suffices. The base Lipschitz principal bound
comes from the fixed low profile regularity. Localization by
$\chi_c$ adds only \eqref{at:cutoff-commutators}.

Sum the four estimates and absorb their cutoff errors.
The source norm \eqref{at:source-norm} uses precisely
$\lambda L_A$ on the edges, and the order-zero remainders are
controlled by $X^0$. This proves \eqref{at:profile-inverse-bound}.
Convergence in the stated graph norm extends the estimate to the
compatible closure; a distributional solution must first be
identified with an element of that domain.
\end{proof}

\subsection{The linearized operator at the perturbed radius}

We estimate $L_R-L_A$ by putting the high norm on either the
argument or the radius perturbation. The coefficient of the high
argument norm is small. The other terms contain a low argument norm,
including the volume contribution from the growing chart. Through
order eight we use enlarged charts; at higher orders the product
estimate retains the original cutoff.

\begin{lemma}
\label{at:low-chart-comparisons}
The enlarged center, fixed interior and local edge norms obey
\begin{equation}
 \|U_f\|_{H^8(B_c)}
   +\lambda^4\|f\|_{H^8(I_i^+)}
   +\lambda^2\bigl(\|f\|_{B_5^8(I_5^+)}
                          +\|f\|_{C_8^8(I_8^+)}\bigr)
       \le C\|f\|_{X^8}.
 \label{at:low-comparison}
\end{equation}
On the fixed outer overlap, for $0\le k\le7$,
\begin{equation}
 \|\partial_x^kf\|_\infty
           \le C\lambda^{-k/2-1/4}\|f\|_{X^8}.
 \label{at:outer-low-sup}
\end{equation}
The regions are enlarged by fixed amounts in their
own coordinates; the long chart is enlarged by a bounded
number of local intervals.
\end{lemma}

\begin{proof}
On every overlap one cutoff is bounded below. At the two
compact center--interior interfaces all coordinate changes
are fixed smooth maps. At the inner edge interface both
$y_5$ and $y_8$ stay in a fixed compact subset of $(0,\infty)$.
At the outer edge interface,
\[
 y_8\sim\lambda^{-1/2},\qquad
 \partial_x=O(\lambda^{-1/2})\partial_{y_8}.
\]
Eight derivatives therefore cost at most $\lambda^{-4}$.
The derivative bounds of the cutoff reciprocals hold on the
sets where that cutoff is bounded below. A subordinate
finite partition gives \eqref{at:low-comparison}; equivalence
of the order-zero measures supplies the factor $\lambda^2$
on the local edge. This argument is used at this fixed
low order only.

For the sharper assertion, if the interior cutoff is bounded
below, one-dimensional Sobolev embedding on a fixed physical
interval applies. Otherwise the long-chart cutoff is bounded
below. On an interval of unit length in the chart coordinate there,
$y_8^{7/2}\sim\lambda^{-7/4}$, so
\[
 \|f\|_{H^8(\text{unit interval})}
       \le C\lambda^{7/4}\lambda^{-2}\|f\|_{X^8}.
\]
The one-dimensional $H^8$ embedding controls seven
derivatives, and each $x$ derivative costs $\lambda^{-1/2}$.
These two alternatives prove \eqref{at:outer-low-sup}.
\end{proof}

\begin{proposition}
\label{at:current-perturbation}
Let $R=A+\lambda^4Z$ have the same center and vacuum type
as $A$, and suppose
\[
             \lambda^{-2}\|Z\|_{X^8}\le\eta_*.
\]
For sufficiently small $\eta_*$, the normalized radius and Jacobian
retain fixed positive lower and upper bounds. For every smooth
displacement $f$ satisfying the center and vacuum conditions of
Definition~\ref{at:atlas}, the following estimates hold.
\begin{enumerate}[label=(\roman*),leftmargin=2em]
\item For $2\le q\le8$,
\begin{equation}
 \mathcal S_{q-2}((L_R-L_A)f)
       \le C_q\lambda^{-1}\|Z\|_{X^8}\|f\|_{X^q}.
 \label{at:current-low}
\end{equation}
\item For $9\le q\le23$,
\begin{equation}
 \mathcal S_{q-2}((L_R-L_A)f)
   \le C_q\left(
       \|Z\|_{X^8}\|f\|_{X^q}
       +\lambda^{-2}\|Z\|_{X^q}\|f\|_{X^8}\right).
 \label{at:current-high}
\end{equation}
\end{enumerate}
The constants depend on the finite profile coefficient bounds and
the fixed positive bounds for the normalized radius and Jacobian.
All dependence on higher norms of the current radius is displayed
in \eqref{at:current-high}.
\end{proposition}

\begin{proof}
\emph{Orders two through eight.} Put $m=q-2$. The exact coefficients
\eqref{at:edge-hessian}--\eqref{at:center-lower-tensor}
are smooth in the first normalized radius factors and affine
in the indicated second derivative. In the interior their
increment is $\lambda^3(Z,Z_x,Z_{xx})$. By
\eqref{at:low-comparison}, the coefficient difference has
$H^6$ norm at most
\[
 C\lambda^3\|Z\|_{H^8(I_i^+)}
                  \le C\lambda^{-1}\|Z\|_{X^8}.
\]
The center bound is $C\lambda^3\|Z\|_{X^8}$. On the edge
the increment is exactly \eqref{at:increment-jets}; the
quotient estimate of Lemma~\ref{at:natural-quotient} gives
\[
 \|\operatorname{coef}(R)-\operatorname{coef}(A)\|_6
              \le C\lambda^2\|Z\|_8
              \le C\|Z\|_{X^8}.
\]
On the growing chart, apply the composition estimate to the
coefficient difference on intervals of fixed length. Every term
contains a radius increment. Multiplying by the local weight
$y^{7/2}$ and summing with finite overlap therefore gives the bound
without a norm of the constant function on the growing interval.

If $m\le6$, multiplication by such a coefficient is bounded
on the order-$m$ test space by
Lemma~\ref{at:low-rung-multiplier}. For example,
\[
 \|\chi_i(a_R-a_A)f_{xx}\|_{H^m}
 \le C\lambda^{-1}\|Z\|_{X^8}
            \bigl(\|\chi_if\|_{H^{m+2}}+\|f\|_2\bigr).
\]
The localized derivative estimate controls $\chi_if_{xx}$ with
the original cutoff. The edge terms $\chi D_nf$, $\chi f_y/y$
and $\chi f$ are treated with that cutoff and the natural quotient.
Multiplying the edge source norms by $\lambda^2$ proves
\eqref{at:current-low}, using only $X^q$ regularity of the argument
also when $q<8$.

\emph{Orders nine through twenty-three.} For $m\ge7$, apply
\eqref{at:cutoff-coefficient} with the same cutoff. On the interior,
the coefficient increment and the differentiated argument satisfy
\[
\begin{array}{ll}
 \|\chi_i u_i\|_{H^m}\le C\lambda^3\|Z\|_{X^q},
 &\|u_i\|_\infty\le C\lambda^{7/4}\|Z\|_{X^8},\\
 \|\chi_i(f,f_x,f_{xx})\|_{H^m}\le C\|f\|_{X^q},
 &\|(f,f_x,f_{xx})\|_\infty
                       \le C\lambda^{-5/4}\|f\|_{X^8}.
\end{array}
\]
The supremum bounds follow from \eqref{at:outer-low-sup}
and the compact interior embedding. The resulting source
bound is
\[
 C_q\left(
 \lambda^{7/4}\|Z\|_{X^q}\|f\|_{X^8}
 +\lambda^{7/4}\|Z\|_{X^8}\|f\|_{X^q}
 +\lambda^{1/2}\|Z\|_{X^8}\|f\|_{X^8}\right).
\]
At the center the corresponding bound is
\[
 C_q\lambda^3
       \left(\|Z\|_{X^q}\|f\|_{X^8}
                         +\|Z\|_{X^8}\|f\|_{X^q}\right).
\]
The radial integral in the zeroth coefficient is bounded
on full Cartesian Sobolev spaces as explained after
\eqref{at:center-lower-tensor}. It introduces no singular
scalar center coefficient.

On the long chart the corresponding bounds are
\[
\begin{array}{ll}
 \|\chi_8u_e\|_{C_8^m}\le C\|Z\|_{X^q},
 &\|u_e\|_\infty\le C\|Z\|_{X^8},\\
 \|\chi_8(f,f_y/y,D_8f)\|_{C_8^m}
                         \le C\lambda^{-2}\|f\|_{X^q},
 &\|(f,f_y/y,D_8f)\|_\infty
                         \le C\lambda^{-2}\|f\|_{X^8}.
\end{array}
\]
The measure of this chart has square root at most
$C\lambda^{-2}$. Applying the product estimate with the same
cutoff and then multiplying by the source factor $\lambda^2$ gives
\[
 C_q\left(
 \|Z\|_{X^q}\|f\|_{X^8}
 +\|Z\|_{X^8}\|f\|_{X^q}
 +\lambda^{-2}\|Z\|_{X^8}\|f\|_{X^8}\right).
\]
Since $q\ge8$,
\[
 \lambda^{-2}\|Z\|_{X^8}\|f\|_{X^8}
       \le\lambda^{-2}\|Z\|_{X^q}\|f\|_{X^8}.
\]
Thus the volume remainder is controlled by the high norm of the
radius perturbation and the low norm of the argument.

Finally choose $0<a<b$ strictly within the region where
$\chi_5=1$, and let $\omega$ be zero below $a$ and one
above $b$. Decompose each vacuum source as
\[
 \chi_5\,\mathrm{source}
    =(1-\omega)\,\mathrm{source}
                   +\chi_5\omega\,\mathrm{source}.
\]
The first term is treated by the one-sided calculus
on a fixed interval extending slightly beyond $b$, where
the original cutoff equals one. The second is supported
away from zero; there the logarithmic derivative of $y^4$
is bounded and the power-cutoff calculation applies.
The transition derivatives of $\omega$ are again inside
$\{\chi_5=1\}$. Both high norms therefore retain the
original cutoff. The natural quotient and the source
normalization give
\[
 C_q\left(\|Z\|_{X^q}\|f\|_{X^8}
                         +\|Z\|_{X^8}\|f\|_{X^q}\right).
\]
Combining the four charts proves \eqref{at:current-high}.
\end{proof}

\begin{proposition}
\label{at:current-inverse}
\label{exlim:upstream-calculus}
Assume the hypotheses of Proposition~\ref{at:current-perturbation},
including $\lambda^{-2}\|Z\|_{X^8}\le\eta_*$, and fix the finite
profile bounds in \eqref{at:profile-budget} on $0<\lambda\le A_0$.
There is $\varepsilon_{\mathrm{inv}}>0$ such that, if
$\lambda^{-1}\|Z\|_{X^8}\le\varepsilon_{\mathrm{inv}}$,
\begin{align}
 \|f\|_{X^q}
  &\le C_q\bigl(\mathcal S_{q-2}(L_Rf)
                              +\|f\|_{X^{q-2}}\bigr),
                    &&2\le q\le8, \label{at:current-inverse-low}\\
 \|f\|_{X^q}
  &\le C_q\bigl(\mathcal S_{q-2}(L_Rf)
                 +\|f\|_{X^{q-2}}
                 +\lambda^{-2}\|Z\|_{X^q}\|f\|_{X^8}\bigr),
                    &&9\le q\le23. \label{at:current-inverse-high}
\end{align}
These estimates hold first on the compatible core and then
on its closed strong domain.
\end{proposition}

\begin{proof}
Write $L_Af=L_Rf-(L_R-L_A)f$ in
\eqref{at:profile-inverse-bound}. By
\eqref{at:current-low}--\eqref{at:current-high},
\[
 \begin{aligned}
 \|f\|_{X^q}
 &\le C_q\{\mathcal S_{q-2}(L_Rf)+\|f\|_{X^{q-2}}\}
       +C_q\lambda^{-1}\|Z\|_{X^8}\|f\|_{X^q},
 &&2\le q\le8,\\
 \|f\|_{X^q}
 &\le C_q\{\mathcal S_{q-2}(L_Rf)+\|f\|_{X^{q-2}}
                  +\lambda^{-2}\|Z\|_{X^q}\|f\|_{X^8}\}\\
 &\quad+C_q\lambda^{-1}\|Z\|_{X^8}\|f\|_{X^q},
 &&9\le q\le23.
 \end{aligned}
\]
Choose $\varepsilon_{\mathrm{inv}}$ so that
$\varepsilon_{\mathrm{inv}}\max_{2\le q\le23}C_q\le1/4$.
Absorbing the last term proves both estimates on the compatible
core. Convergence in the closed strong graph gives the domain assertion.

The later low-norm bound gives
$\lambda^{-1}\|Z\|_{X^8}\le\lambda\eta_*$. Once the profile
parameters are fixed, decreasing the upper bound $A_0$ for $\lambda$
therefore enforces the threshold without increasing a profile
constant. The energy growth constants are selected before the small
collapse parameter; the constant in this elliptic estimate affects
only the subsequent choice of $A_0$.
\end{proof}

\section{The closed pressure form and the homologous direction}
\label{pr:section}

To separate the homologous direction in the energy estimates, we use
a projection fixed in mass coordinates. It therefore commutes with
time differentiation. The homologous mode has a fractional expansion
at vacuum, so we first close the pressure form and approximate the
mode by compatible smooth functions.

We also extend the mechanical Hessian of each admissible prepared
radius to this form domain and bound its absolute value. Coercivity
on the complementary subspace is proved later in
Lemma~\ref{en:complement-coercivity}.

Fix the selected parameter $\beta>0$, the mass $M=M_\beta$, and
$0<\lambda\le A_*\le1$. Let $A$ be the corresponding matched radius
from Proposition~\ref{exlim:upstream-profile}. We use precisely
the covering atlas of Definition~\ref{at:atlas}; in particular,
\begin{equation}
 \chi_c+\chi_i+\chi_5+\chi_8=1,\qquad
 \chi_5=\chi_e\zeta((M-x)/\lambda^4),\qquad
 \chi_8=\chi_e(1-\zeta((M-x)/\lambda^4)).
 \label{pr:partition}
\end{equation}
The center norm is the norm of the full Cartesian radial vector.
The edge variables and measures are
\begin{equation}
 x=M-\lambda^4y_n^n,\qquad
 |\dd x|=n\lambda^4y_n^{n-1}\dd y_n,\qquad n=5,8.
 \label{pr:edge-measures}
\end{equation}
The fifth-root chart is bounded, whereas
$c<y_8<C_e\lambda^{-1/2}$. All derivatives of each localized
function $\chi_\nu f$ are included in $X^q$. Outside the matching
region, the long chart uses the core profile.

In this section $C$ may depend on the selected member, its finite
coefficient bounds, and the fixed cutoffs, but not on $\lambda$.
The leading growth constants in the energy argument will be fixed
separately, before the profile parameters are chosen.

\subsection{The pressure operator on the compatible core}

Write
\begin{equation}
 H=L^2((0,M),\dd x),\qquad
 \rho_R=(4\pi R^2R_x)^{-1},\qquad
 P'(\rho)=\frac{4\rho^{2/3}}{3\sqrt{1+\rho^{2/3}}},\quad P(0)=0.
 \label{pr:pressure-law}
\end{equation}
The normalized profile factors on either edge are
\begin{equation}
 A=\lambda r,\qquad A_y=-\lambda^2y\gamma,\qquad
 T_\gamma=y\gamma_y,\qquad
 \rho_A=\frac{ny^{n-2}}{4\pi r^2\gamma},\qquad
 a_n=\frac{P'(\rho_A)}{y^2\gamma^2}.
 \label{pr:profile-units}
\end{equation}
The matched-profile bounds give
\begin{equation}
 0<c\le r,\gamma,a_n\le C,\qquad
 |T_\gamma|+\lambda y^2\le C.
 \label{pr:profile-geometry}
\end{equation}
These statements hold on the entire prescribed charts. At the center,
write $z=x^{1/3}$ and $A=\lambda z a(z^2)$. The normalized
deformation eigenvalues are positive and the required finite derivatives
are bounded. On the
fixed interior the radius and its normalized Jacobian are positive factors with the corresponding bounds.

Let $\mathscr C_\lambda$ be the common displacement core. Its elements
have smooth regular radial-vector lifts at the center and are smooth
one-sided functions at the vacuum with $f_{y_5}(0)=0$; in the
interior they are ordinary smooth functions. In particular, no even
scalar extension is imposed at the vacuum. Define on this core
\begin{equation}
 \mathcal T_Af=\lambda^{3/2}\sqrt{P'(\rho_A)}
                    \left(\frac{2f}{A}+\frac{f_x}{A_x}\right),
 \qquad
 \mathcal P_A[f,g]=(\mathcal T_Af,\mathcal T_Ag)_H,\qquad
 \mathcal P_A[f]=\mathcal P_A[f,f].
 \label{pr:pressure-operator}
\end{equation}

\begin{lemma}
\label{pr:upper-pressure}
The prescribed zero-order norm and the pressure operator satisfy
\begin{equation}
 C^{-1}\|f\|_H\le\|f\|_{X^0}\le C\|f\|_H,\qquad
 \|\mathcal T_Af\|_H\le C\lambda^{-1/2}\|f\|_{X^1},
 \qquad f\in\mathscr C_\lambda.
 \label{pr:upper-pressure-bound}
\end{equation}
\end{lemma}

\begin{proof}
The squared-cutoff sum lies between $1/4$ and $1$.
Together with \eqref{pr:edge-measures} and the fixed positive compact
chart densities, this proves the first comparison.

The ratio $f_x/A_x$ is unchanged by the edge coordinate change.
Substitution of \eqref{pr:profile-units} gives the exact identity
\begin{equation}
 \mathcal T_Af=-\lambda^{-1/2}\sqrt{a_n}\,f_y
              +\frac{2\lambda^{1/2}\sqrt{P'(\rho_A)}}{r}\,f.
 \label{pr:edge-pressure-operator}
\end{equation}
The first coefficient is bounded by $C\lambda^{-1/2}$. The
square of the second is
\[
 \frac{4\lambda P'(\rho_A)}{r^2}
       =\frac{4\lambda a_n\gamma^2y^2}{r^2}\le C.
\]
For a function supported in the corresponding edge chart, this gives
\begin{equation}
 \|\mathcal T_Af\|_H^2
 \le C\lambda^4\int
       \bigl(\lambda^{-1}|f_y|^2+|f|^2\bigr)y^{n-1}\dd y.
 \label{pr:edge-pressure-bound}
\end{equation}

At the center the strain in \eqref{pr:pressure-operator} is
\[
 \lambda^{-1}\left(\frac{f_z}{a+za_z}+\frac{2f}{za}\right).
\]
Since $\rho_A=\lambda^{-3}$ times a positive factor,
$\lambda P'(\rho_A)$ is uniformly bounded. Both $f_z$ and
$f/z$ are controlled by the complete Cartesian derivative of
$f(z)e_z$. The center estimate therefore has no loss of
$\lambda$. The fixed interior has the same calculation without
a coordinate singularity.

Linearity and \eqref{pr:partition} give exactly
\begin{equation}
                 \mathcal T_Af=\sum_{\nu=c,i,5,8}\mathcal T_A(\chi_\nu f).
 \label{pr:pressure-partition}
\end{equation}
Each cutoff derivative on the right is already included in the
corresponding $X^1$ norm. Applying the preceding four estimates
and summing proves \eqref{pr:upper-pressure-bound}.
\end{proof}

\begin{lemma}
\label{pr:closure}
The operator $\mathcal T_A:\mathscr C_\lambda\subset H\to H$ has
the following properties.
\begin{enumerate}[label=(\roman*),leftmargin=2em]
\item It is closable, and its graph closure defines the Hilbert space
\begin{equation}
 V_A=\overline{\mathscr C_\lambda}^{\,\|\cdot\|_{V_A}},
 \qquad \|f\|_{V_A}^2=\|f\|_H^2+\mathcal P_A[f],
 \label{pr:closed-domain}
\end{equation}
which embeds injectively into $H$.
\item The closure of $\mathscr C_\lambda$ in $X^1$ embeds into
$V_A$, and \eqref{pr:upper-pressure-bound} holds on this completion.
\item For every $a>0$, the two completions coincide, with equivalent
norms uniformly for $a\le\lambda\le A_*$.
\end{enumerate}
\end{lemma}

\begin{proof}
The core is dense in $H$ because it contains $C_c^\infty(0,M)$.
If $f_k\to0$ and $\mathcal T_Af_k\to g$ in $H$, write
$\mathcal T_A=c_1\partial_x+c_0$ on a compact interior interval.
Its coefficients are $C^1$, so for $\phi\in C_c^\infty(0,M)$,
\[
 (g,\phi)_H
 =\lim_k\int f_k\{-\partial_x(c_1\phi)+c_0\phi\}\dd x=0.
\]
Hence $g=0$: the graph closure is a Hilbert subspace of
$H\times H$ with injective first projection. Moreover,
\eqref{pr:upper-pressure-bound} gives
\[
 \|f\|_{V_A}^2=\|f\|_H^2+\|\mathcal T_Af\|_H^2
                    \le C_\lambda\|f\|_{X^1}^2.
\]
Thus an $X^1$-Cauchy core sequence is $V_A$-Cauchy, with the
same limit in $H$.

For the reverse inclusion, fix $a\le\lambda\le A_*$ and prove
\begin{equation}
 \|f\|_{X^1}\le C_a\bigl(\|\mathcal T_Af\|_H+\|f\|_H\bigr).
 \label{pr:fixed-positive-graph}
\end{equation}
At the center set $r=A(z)$ and $F(r)=f(z)$. Then
\[
 \dd x=4\pi\rho_A(r)r^2\dd r,\qquad
 \mathcal T_Af=\lambda^{3/2}\sqrt{P'(\rho_A)}
                         \left(F_r+\frac{2F}{r}\right).
\]
Here $\rho_A,P'(\rho_A),A/z,A_z$ have positive upper and lower
bounds depending on $a$. Since $F(0)=0$,
\[
 \int_0^{r_c}\left|F_r+\frac{2F}{r}\right|^2r^2\dd r
 =\int_0^{r_c}r^2|F_r|^2\dd r
       +2\int_0^{r_c}|F|^2\dd r+2r_c|F(r_c)|^2.
\]
This controls the full radial-vector $H^1$ norm.
On an edge, solving \eqref{pr:edge-pressure-operator} for $f_y$
gives
\[
 \|f_y\|_{L^2(\dd x)}
 \le C_a(\|\mathcal T_Af\|_H+\|f\|_H),
\]
because $a_n\ge c>0$ and the zeroth coefficient is bounded.
The interior coefficient of $f_x$ is also bounded away from zero.
For the prescribed cutoffs,
\[
 \mathcal T_A(\chi_\nu f)=\chi_\nu\mathcal T_Af+
                    c_1(\partial_x\chi_\nu)f,\qquad
 \|c_1(\partial_x\chi_\nu)f\|_H\le C_a\|f\|_H;
\]
their derivatives are supported away from the physical endpoints.
Summing proves \eqref{pr:fixed-positive-graph}. Together with
\eqref{pr:upper-pressure-bound}, it identifies the two completions
uniformly on $[a,A_*]$, without adding a classical endpoint trace.
\end{proof}

The condition $f_{y_5}(0)=0$ on the core need not persist as a
classical trace in $V_A$. For example, let $f(y)=y$ near zero,
with a smooth cutoff in the vacuum chart. Choose $\eta=0$ below
$1$ and $\eta=1$ above $2$, and set $f_\epsilon=\eta(y/\epsilon)f$.
Then $f_\epsilon$ belongs to the compatible core, and
\[
 \begin{aligned}
 \|f-f_\epsilon\|_{V_A}^2
 &\le C\lambda^4\int_0^{2\epsilon}
              (\lambda^{-1}+y^2)y^4\dd y\\
 &\le C(\lambda^3\epsilon^5+\lambda^4\epsilon^7)
                  \longrightarrow0
 \end{aligned}
\]
at every fixed positive scale, by \eqref{pr:edge-pressure-bound}.
Thus $f\in V_A$, although $f_y(0)=1$. One-sided cubic terms are
also allowed in the core.

\subsection{The fixed homologous direction}

Let $x=m_\beta(z_\beta)$ be the reference mass map, and put
\begin{equation}
 \varphi_\beta=\frac{z_\beta}{\|z_\beta\|_H},\qquad
 f^\perp=f-(f,\varphi_\beta)_H\varphi_\beta.
 \label{pr:fixed-mode}
\end{equation}
The reference law is $\dd m_\beta=4\pi w_\beta^3z^2\dd z$, with
$w_\beta\sim1-z$ near $z=1$. The mode is independent of
time, so commuting time derivatives with the projection produces no
additional terms.
It is smooth in the fourth-root mass coordinate but has a fractional
expansion in the fifth-root vacuum coordinate. The next lemma
approximates it in the form domain by compatible smooth functions.

\begin{lemma}
\label{pr:mode-membership}
The fixed mode belongs to the $X^1$ completion of
$\mathscr C_\lambda$ and satisfies
\begin{equation}
 \|\varphi_\beta\|_{X^1}\le C,\qquad
                       \mathcal P_A[\varphi_\beta]\le C.
 \label{pr:mode-bounds}
\end{equation}
\end{lemma}

\begin{proof}
Put $d=M-x$. Integration of the reference mass law gives
\begin{equation}
 1-z_\beta\sim d^{1/4},\qquad
 |\partial_x\varphi_\beta|\le Cd^{-3/4},\qquad
 |\varphi_\beta|\le C.
 \label{pr:mode-mass-tail}
\end{equation}
Consequently,
\begin{equation}
 |(\varphi_\beta)_{y_5}|\le C\lambda y_5^{1/4},\qquad
 |(\varphi_\beta)_{y_8}|\le C\lambda y_8.
 \label{pr:mode-native-tail}
\end{equation}
The reference inverse mass map has a regular radial-vector center
lift, so the compact-chart $X^1$ norms are bounded.
On the long chart,
\begin{equation}
 C\lambda^4\int_c^{C_e\lambda^{-1/2}}y^7\dd y\le C,\qquad
 C\lambda^6\int_c^{C_e\lambda^{-1/2}}y^9\dd y\le C\lambda.
 \label{pr:mode-long-graph}
\end{equation}
Cutoff derivatives contribute only the first bound. The bounded
fifth-root chart is controlled by \eqref{pr:mode-native-tail}.

For the pressure form, \eqref{pr:edge-pressure-operator} bounds each
long-chart term in squared norm by
\begin{equation}
 C\lambda^5\int_c^{C_e\lambda^{-1/2}}y^9\dd y\le C.
 \label{pr:mode-long-pressure}
\end{equation}
On the bounded chart the bound is
$C\lambda^5\int_0^C(y^{9/2}+y^6)\dd y$; on compact charts use
the regular mode. Summing proves both estimates in
\eqref{pr:mode-bounds}.

For the compatible approximation, use the cutoff $\eta$ above and set
\begin{equation}
 \varphi_{\beta,\epsilon}(y)
 =\varphi_\beta(0)+
       \eta(y/\epsilon)\{\varphi_\beta(y)-\varphi_\beta(0)\}
 \label{pr:mode-approximation}
\end{equation}
on the bounded vacuum chart, leaving the mode unchanged elsewhere.
By \eqref{pr:mode-mass-tail},
\[
 |\varphi_\beta(y)-\varphi_\beta(0)|\le C\lambda y^{5/4}.
\]
The cutoff derivative on $\epsilon<y<2\epsilon$ is
$O(\lambda\epsilon^{1/4})$, with squared pressure norm
\begin{equation}
 C\lambda^3\lambda^2\epsilon^{1/2}
             \int_\epsilon^{2\epsilon}y^4\dd y
                         \le C\lambda^5\epsilon^{11/2}.
 \label{pr:mode-approximation-error}
\end{equation}
The discarded derivative tail obeys the same estimate. The squared
$H$ norm of the error is $O(\lambda^6\epsilon^{15/2})$, and the zeroth-order
pressure coefficient is bounded. Thus, for $0<\epsilon\le1$,
\[
 \begin{aligned}
 \|\varphi_{\beta,\epsilon}-\varphi_\beta\|_{X^1}^2
      &\le C\lambda^6\epsilon^{11/2},\\
 \|\varphi_{\beta,\epsilon}-\varphi_\beta\|_{V_A}^2
      &\le C\lambda^5\epsilon^{11/2}.
 \end{aligned}
\]
The approximants are constant near vacuum and unchanged at the
regular center, hence lie in $\mathscr C_\lambda$. They converge
in both completions.
\end{proof}

\begin{proposition}
\label{pr:pressure-domain}
\label{en:input-pressure-domain}
Let $V_A$ be the pressure form domain in \eqref{pr:closed-domain}.
\begin{enumerate}[label=(\roman*),leftmargin=2em]
\item The projection $f\mapsto f^\perp$ is bounded on $V_A$ and on
the compatible $X^1$ completion, and
\begin{equation}
 \|\mathcal T_Af^\perp\|_H\le C\lambda^{-1/2}\|f\|_{X^1}.
 \label{pr:projected-pressure-bound}
\end{equation}

\item For every $b>0$, the domain
$\{f\in H:f^\perp\in V_A\}$ is exactly $V_A$, and
\begin{equation}
 \|f\|_{Z_A}^2=
 \mathcal P_A[f^\perp]+\|f^\perp\|_H^2
                         +b^2|(f,\varphi_\beta)_H|^2
 \label{pr:mode-weighted-domain}
\end{equation}
is a Hilbert norm on it. For bounded positive $b$,
\begin{equation}
 \|f\|_{Z_A}\le C\|f\|_{V_A},\qquad
 \|f\|_{V_A}\le C\max(1,b^{-1})\|f\|_{Z_A}.
 \label{pr:domain-comparison}
\end{equation}
\end{enumerate}
\end{proposition}

\begin{proof}
By Lemma~\ref{pr:mode-membership} and the mass comparison,
\[
 \|f^\perp\|_{X^1}
 \le\|f\|_{X^1}+|(f,\varphi_\beta)_H|
                               \|\varphi_\beta\|_{X^1}
 \le C\|f\|_{X^1}.
\]
Both summands lie in the compatible completion, and
Lemma~\ref{pr:upper-pressure} gives
\eqref{pr:projected-pressure-bound}. The same calculation in $V_A$
gives $\|f^\perp\|_{V_A}\le C\|f\|_{V_A}$.

Since $\varphi_\beta\in V_A$ and $\|\varphi_\beta\|_H=1$,
\[
 f\in V_A\ \Longleftrightarrow\ f^\perp\in V_A,
 \qquad
 V_A=\{g\in V_A:(g,\varphi_\beta)_H=0\}
                    \oplus\operatorname{span}\{\varphi_\beta\}.
\]
The first summand is closed, so \eqref{pr:mode-weighted-domain}
defines a Hilbert norm for $b>0$. If $a=(f,\varphi_\beta)_H$, then
\[
 \begin{aligned}
 \|f\|_{Z_A}
 &\le\|f^\perp\|_{V_A}+b|a|\le C\|f\|_{V_A},\\
 \|f\|_{V_A}
 &\le\|f^\perp\|_{V_A}+|a|\|\varphi_\beta\|_{V_A}
 \le C\max(1,b^{-1})\|f\|_{Z_A}.
 \end{aligned}
\]
Here $b$ is bounded above. In the application $b\ge\beta>0$;
no norm at $b=0$ is asserted.
\end{proof}

\subsection{The Hessian on the form domain of the pressure operator}

For a positive radius $R$, define
\begin{equation}
 \mathcal U[R]=\int_0^M\mathfrak e(\rho_R)\dd x
                         -\int_0^M\frac{x}{R}\dd x.
 \label{pr:potential}
\end{equation}
Here $\mathfrak e$ is the specific internal energy with the
vacuum normalization \eqref{eq:internal-energy-normalization}.
An admissible radius has a regular positive Cartesian deformation
at the center, a positive interior Jacobian, and the vacuum form
\[
 R=R_b-\xi^2G,\qquad
 Q=-R_\xi/\xi=2G+\xi G_\xi>0,\qquad \xi=(M-x)^{1/5}.
\]
The normalized factors have the finite one-sided regularity of the
prepared radii. In particular, $G,Q$ have two continuous derivatives,
which suffice for the first and second variations below. At each fixed
positive scale the density is $\xi^3$ times a positive regular factor.
These endpoint bounds justify differentiation of the potential under
the integral.

\begin{lemma}
\label{pr:hessian-identity}
For an admissible radius and compatible core directions $f,g$,
\begin{equation}
 D\mathcal U[R]f=(\mathcal M[R],f)_H,\qquad
 \mathcal M[R]=4\pi R^2\partial_xP(\rho_R)+\frac{x}{R^2},
 \label{pr:first-variation}
\end{equation}
and
\begin{align}
 Q_R[f,g]:=\lambda^3D^2\mathcal U[R][f,g]
   &=\mathcal P_R[f,g]+\int_0^M U_Rfg\dd x,
 \label{pr:hessian-form}\\
 U_R&=\lambda^3\left\{8\pi R\partial_xP(\rho_R)
                                      -\frac{2x}{R^3}\right\}.
 \label{pr:hessian-potential}
\end{align}
Here $\mathcal P_R$ is defined by
\eqref{pr:pressure-operator} with $A$ replaced by $R$.
Both physical endpoint fluxes vanish.
\end{lemma}

\begin{proof}
Put $G_f=2f/R+f_x/R_x$. The exact variation of density is
$D\rho_R[f]=-\rho_RG_f$. Therefore
\[
 D\!\left(\int\mathfrak e(\rho_R)\dd x\right)f
 =-\int4\pi P(\rho_R)(R^2f)_x\dd x.
\]
Integration by parts gives the pressure part of
\eqref{pr:first-variation}, and differentiation of
$-\int x/R\dd x$ gives its gravitational part.
The boundary term is $-4\pi R^2P(\rho_R)f$.

At the vacuum write $\xi=(M-x)^{1/5}$.
The positivity and boundary compatibility conditions give
\[
 R=R_b-\xi^2G,\qquad \rho_R=\xi^3\Theta,\qquad
 P(\rho_R)=O(\xi^5),\qquad D\rho_R[f]=O(\xi^3).
\]
Core directions are bounded there and satisfy $f_\xi=O(\xi)$.
Thus the first-variation flux is $O(\xi^5)$.
At the center $R=O(z)$, the density is a regular even factor,
and a regular radial direction is $O(z)$, giving flux $O(z^3)$.
These calculations also justify taking a second finite variation:
the resulting boundary factors have the same vanishing orders.

More explicitly,
\[
 D\mathcal M[R]f
 =8\pi Rf\,P_x
    +4\pi R^2\partial_x\{P'D\rho_R[f]\}
                                  -\frac{2x}{R^3}f.
\]
For the middle term, integration by parts gives
\[
 \begin{aligned}
 \int 4\pi R^2g\,\partial_x(P'D\rho_R[f])\dd x
 &=-\int 4\pi(R^2g)_xP'D\rho_R[f]\dd x\\
 &=\int P'G_fG_g\dd x.
 \end{aligned}
\]
Here $D\rho_R[f]=-\rho_RG_f$, $(R^2g)_x=R^2R_xG_g$, and
$4\pi R^2R_x\rho_R=1$. The boundary term vanishes by the endpoint
estimates above. The remaining two terms give
\eqref{pr:hessian-form}--\eqref{pr:hessian-potential}.
\end{proof}

\begin{proposition}
\label{pr:current-hessian}
Let $R$ be an admissible radius with the following properties:
\begin{enumerate}[label=(\roman*),leftmargin=2em]
\item $R=A$ in a fixed neighborhood of the center.
\item On the edges, the factors
\[
 R=\lambda r_R,\qquad R_y=-\lambda^2y\gamma_R,\qquad
 T_{\gamma_R}=y(\gamma_R)_y,\qquad
 a_{R,n}=\frac{P'(\rho_R)}{y^2\gamma_R^2}
\]
satisfy the bounds in \eqref{pr:profile-geometry}.
\item On the fixed interior, $R/\lambda$ and $\partial_x(R/\lambda)$
satisfy the same positive lower and upper bounds as the corresponding
quantities for $A$, and $|\partial_x^2(R/\lambda)|$ satisfies the
same fixed bound.
\end{enumerate}
Then the closures of $\mathcal T_R$ and $\mathcal T_A$ have the same
domain, with equivalent norms after adding the squared mass norm to
the pressure forms. On this domain the Hessian extends continuously
and satisfies
\begin{equation}
 |Q_R[f,f]|
 \le C\{\mathcal P_A[f]+\|f\|_H^2\}
 \le C\{\mathcal P_A[f^\perp]+\|f\|_H^2\}.
 \label{pr:current-upper-bound}
\end{equation}
The constant is uniform on any family satisfying these bounds. This is an absolute upper bound, not a coercivity statement.
\end{proposition}

\begin{proof}
\emph{Comparison of the pressure forms.} For $\alpha=A_x/R_x$,
direct algebra gives
\begin{equation}
 \frac{2f}{R}+\frac{f_x}{R_x}
 =\alpha\left(\frac{2f}{A}+\frac{f_x}{A_x}\right)
                         +2\left(\frac1R-\frac{\alpha}{A}\right)f.
 \label{pr:strain-comparison}
\end{equation}
The ratios $\alpha$, $\rho_R/\rho_A$, and
$P'(\rho_R)/P'(\rho_A)$ are bounded above and below.
For the last ratio use \eqref{pr:pressure-law} directly.
On the edges the coefficient of $f$ in
\eqref{pr:strain-comparison}, after multiplication by
$\lambda^{3/2}\sqrt{P'(\rho_R)}$, has size at most
\[
 C\lambda^{1/2}\sqrt{P'(\rho_R)}
    =C\lambda^{1/2}y\gamma_R\sqrt{a_{R,n}}\le C.
\]
It is also bounded on the fixed interior by the high-density
scaling, and it is zero near the center because $R=A$.
Consequently
\begin{equation}
 \mathcal P_R[f]\le C\{\mathcal P_A[f]+\|f\|_H^2\}.
 \label{pr:current-pressure-comparison}
\end{equation}
Interchanging $R,A$ gives the reverse comparison after the
squared mass norm is added to both sides. The closability proof of Lemma~\ref{pr:closure}
applies equally to $\mathcal T_R$. Equivalent graph norms on the same
core therefore have the same completion in $H$.

\emph{The potential term.} Use the decreasing edge coordinate:
\[
 \partial_x=-\frac1{n\lambda^4y^{n-1}}\partial_y,\qquad
 \frac{y(\rho_R)_y}{\rho_R}
 =n-2+\frac{2\lambda y^2\gamma_R}{r_R}
                                      -\frac{T_{\gamma_R}}{\gamma_R}.
\]
Substitution into \eqref{pr:hessian-potential} gives exactly
\begin{equation}
 U_R=-\frac{2a_{R,n}\gamma_R}{r_R}
       \left(n-2+\frac{2\lambda y^2\gamma_R}{r_R}
                         -\frac{T_{\gamma_R}}{\gamma_R}\right)
                         -\frac{2x}{r_R^3}.
 \label{pr:edge-hessian-potential}
\end{equation}
All factors are bounded throughout the chart. This estimate uses
$P'$ itself and requires no derivative of $P'$.
The fixed interior follows from the normalized second-derivative
bounds. At the center, writing
$x=z^3$, $R=\lambda z a(z^2)$, and
$\rho_R=\lambda^{-3}b_0(z^2)$, one has
\[
 \partial_x\rho_R=O(\lambda^{-3}/z),\qquad
 P'(\rho_R)=O(\lambda^{-1}),\qquad
 \lambda^3R\partial_xP(\rho_R)=O(1),\qquad
 \lambda^3x/R^3=O(1).
\]
Hence $\|U_R\|_\infty\le C$.
Equations \eqref{pr:hessian-form} and
\eqref{pr:current-pressure-comparison} prove the first bound in
\eqref{pr:current-upper-bound}. Decompose $f$ as in
\eqref{pr:fixed-mode} and use \eqref{pr:mode-bounds} for the second.
Polarization, or Cauchy--Schwarz in the pressure term and in $H$,
gives the continuous bilinear extension to the completed domain.
\end{proof}

In particular, for a sufficiently large fixed $c_R$,
\[
 Q_R[f,f]+c_R\|f\|_H^2\sim\|f\|_{V_A}^2.
\]
Thus the shifted form is closed. Its natural $H$-realization,
when needed, is specified without extra ordinary endpoint traces by
\begin{equation}
 \begin{split}
 \operatorname{Dom}_H(L_R)
 =\{f\in V_A:\ &\text{there is }h\in H\text{ such that}\\
              &Q_R[f,g]=(h,g)_H\text{ for every }g\in V_A\},
 \qquad L_Rf=h.
 \end{split}
 \label{pr:weak-operator-domain}
\end{equation}
The representative $h$ is unique because $V_A$ is dense in $H$.
On the common core it is exactly $\lambda^3D\mathcal M[R]f$,
by Lemma~\ref{pr:hessian-identity}. On the full domain its
interior distributional expression is the same, by testing against
compactly supported smooth functions. No identification with powers
of a reference spectral operator is involved.

\subsection{Application to finite prepared data}

The initial-energy estimate applies the Hessian bound along the
whole curve of compatible radii. The correction constructed in
Lemma~\ref{loc:finite-map} and
Subsection~\ref{loc:common-mass-correction} has fixed total mass;
we now check the hypotheses of Proposition~\ref{pr:current-hessian}
along its curve.

Write this common-mass radius curve as $R^\theta$,
$0\le\theta\le1$, with finite parameter curve $\Xi(\theta)$,
$R^0=A$, and $R^1=R^\sharp$. Put
$D_{\rm par}=\sup_\theta|\Xi'(\theta)|$.
Its exact common-mass derivative on the bounded collar is
\begin{equation}
 \partial_\theta R^\theta=\lambda^6y_5^3B_\theta(y_5),
 \qquad \|B_\theta\|_{C^2}\le CD_{\rm par}.
 \label{pr:prepared-radius-path}
\end{equation}
On the long edge outside this collar it is
$\lambda^6y_8^{-6}\widetilde B_\theta$, with the analogous
bounded derivatives in the chart coordinates; the transition obeys the weaker
uniform bounds. Integrating and differentiating
\eqref{pr:prepared-radius-path} twice yields
\begin{equation}
 \delta r=O(\lambda^5y^3D_{\rm par}),\qquad
 \delta\gamma=O(\lambda^4yD_{\rm par}),\qquad
 \delta T_\gamma=O(\lambda^4yD_{\rm par})
 \label{pr:prepared-unit-increments}
\end{equation}
on the bounded chart. On the long chart these bounds may be
weakened to $C\lambda^5D_{\rm par}$,
$C\lambda^4D_{\rm par}$, and
$C\lambda^4D_{\rm par}$, respectively. The preparation estimate
makes $\lambda^4D_{\rm par}$ a sufficiently small multiple of the
positive lower bounds for the profile after the selected $A_*$ is decreased.
Thus all segments from $A$ to $R^\sharp$ have the required
positive radius and Jacobian factors and bounded $T_\gamma$.
The positive bounds for $a_{R,n}$ follow from
$\rho_R/\rho_A\sim1$, the pressure law, and the profile
bound for $a_n$.

Both the correction and its mass compensator are supported away
from the center. Below their support the exact mass integrals
coincide, so their inverse mass maps coincide as well:
\[
                        R^\sharp=A\quad\text{near the center}.
\]
On the fixed interior, the normalized radius and Jacobian retain
their positive bounds, and the second derivative remains bounded.
Together with the edge estimates, these bounds verify the hypotheses
of Proposition~\ref{pr:current-hessian} along the prepared segments.
The Hessian bound therefore applies along the whole preparation
curve, with constants determined by the profile and correction bounds.

\begin{corollary}
\label{pr:finite-array}
Let $D_\lambda\ge0$, and let $\psi_0,\ldots,\psi_{23}$ belong to
the corresponding compatible completions. Suppose
\[
 \|\psi_j\|_{X^q}\le C\lambda^{-L_j}D_\lambda,\qquad
 j+q\le23,\qquad L_0=L_1=3,\quad L_j=2j-1\quad(j\ge2).
\]
Then
\begin{equation}
 \sum_{j+q\le23}\|\psi_j\|_{X^q}^2
       +\sum_{j=0}^{22}\mathcal P_A[\psi_j^\perp]
                         \le C\lambda^{-90}D_\lambda^2.
 \label{pr:finite-array-bound}
\end{equation}
\end{corollary}

\begin{proof}
The graph and projected-pressure estimates give
\[
 \begin{aligned}
 \sum_{j+q\le23}\|\psi_j\|_{X^q}^2
 &\le CD_\lambda^2\sum_{j=0}^{23}(24-j)\lambda^{-2L_j}
 \le C\lambda^{-90}D_\lambda^2,\\
 \sum_{j=0}^{22}\mathcal P_A[\psi_j^\perp]
 &\le C\lambda^{-1}\sum_{j=0}^{22}\|\psi_j\|_{X^1}^2
 \le CD_\lambda^2\sum_{j=0}^{22}\lambda^{-2L_j-1}
 \le C\lambda^{-87}D_\lambda^2.
 \end{aligned}
\]
Here $L_{23}=45$, $L_{22}=43$, and $0<\lambda\le1$.
This proves \eqref{pr:finite-array-bound}. The pressure sum stops
at time order $22$; no twenty-fourth derivative or trace of the
twenty-fifth source derivative is used.
\end{proof}

The next section constructs compatible initial time derivatives satisfying
these estimates. The absolute Hessian bound
\eqref{pr:current-upper-bound} controls the quadratic part of their
initial energy; the relative potential and commutator corrections are
estimated there from their defining formulas. The differential energy
inequality is proved in Section~\ref{en:section}.

\section{Compatible data and local solutions of the regularized equation}
\label{loc:section}

Fix the matched trajectory of Proposition~\ref{exlim:upstream-profile}
and its mass interval. We construct terminal data with the prescribed
mass whose time derivatives, generated by the equation, satisfy the
vacuum compatibility conditions. We then solve the regularized equation
with these data. The correction is estimated in the localized weighted
norms of Definition~\ref{at:atlas} and the pressure form domain of
Proposition~\ref{pr:pressure-domain}.

For polytropic Euler flow, Jang--Masmoudi proved local well-posedness
in one dimension \cite{JangMasmoudi2009} and in three dimensions
\cite{JangMasmoudi}, using spaces adapted to the physical vacuum.
Coutand--Lindblad--Shkoller established three-dimensional a priori
estimates combining weighted energy bounds, vorticity estimates, and
elliptic control of the normal derivatives \cite{CoutandLindbladShkoller}.
Coutand--Shkoller \cite{CoutandShkoller} obtained three-dimensional
solutions by a degenerate parabolic regularization that retains the
degeneracy of the pressure at vacuum, with estimates uniform in the regularization
parameter. Gu--Lei \cite{GuLei2016} treated three-dimensional
polytropic Euler--Poisson flow for $1<\gamma<3$, while
Luo--Xin--Zeng \cite{LuoXinZeng} constructed spherical solutions with
or without self-gravity using separate estimates near the center and
the vacuum boundary. Ifrim--Tataru \cite{IfrimTataru} established an
Eulerian well-posedness theory for polytropic Euler flow in weighted
Sobolev spaces, including continuous dependence and a continuation
criterion.

The exact pressure law requires compatibility estimates that remain
uniform as the terminal scale tends to zero. We first impose the Euler
compatibility conditions. Their Jacobian is lower triangular, but its
off-diagonal entries may grow as the scale decreases.
Lemma~\ref{loc:finite-map} bounds the full inverse; the endpoint residual
decays fast enough to solve the constraints. We then compare the corrected
data with the fixed reference trajectory in the same mass coordinate
(Subsection~\ref{loc:common-mass-correction}). This comparison retains the
distinct core and surface scales. Finally,
Proposition~\ref{loc:selector} imposes viscous compatibility by a further
small correction at each fixed positive scale. The allowed viscosity
interval may depend on that scale.

For fixed positive viscosity, the local construction starts with the
radius and velocity in H\"older spaces. Higher time derivatives are
recovered on the same interval and justify the spatial recovery and
energy identities. The initial energy bound in
Proposition~\ref{loc:prepared-energy} is uniform in the terminal scale;
the local lifespan may depend on that scale and on the viscosity.
Section~\ref{ct:section} uses the a priori estimates and continuation
to obtain a common interval of existence.

\subsection{Fixed-mass domains and differentiation of the force}
\label{loc:domains}

The enclosed mass is $x\in[0,M]$. At the center use $z=x^{1/3}$
and the equivariant vector lift of a displacement; at the vacuum use

\[
 \xi=(M-x)^{1/5}.
\]
These charts and their enlarged overlaps are fixed. On each fixed
positive-scale interval, bounded smooth coordinate maps relate them
to the scale-dependent charts used in the uniform weighted estimates.
Fix $0<\alpha<1$, and put $\theta=\alpha/2$. For each integer
$m\ge0$, define the little H\"older domain $\mathcal D^m$ by the
full equivariant space $h^{m+2,\alpha}(B^3;\mathbb R^3)$ at the
center, the ordinary $h^{m+2,\alpha}$ spaces in the interior,
and the following space at the vacuum:
\begin{equation}
 \left\{f\in h^{m+2,\alpha}[0,\ell]:
 f_\xi(0)=0,\quad
 \frac{f_\xi}{\xi},\ \frac{f-f(0)}{\xi^2}\in h^{m,\alpha}[0,\ell]
 \right\}.
 \label{loc:domain}
\end{equation}
The representatives agree on overlaps. The source space $\mathcal Y^m$
has the corresponding $h^{m,\alpha}$ components, with no prescribed
first-derivative trace at the vacuum. Write $\mathcal D=\mathcal D^0$,
$\mathcal Y=\mathcal Y^0$. The quotient identities
\begin{equation}
 \frac{f_\xi}{\xi}=\int_0^1f_{\xi\xi}(t\xi)\dd t,
 \qquad
 \frac{f-f(0)}{\xi^2}=\int_0^1(1-t)f_{\xi\xi}(t\xi)\dd t
 \label{loc:quotients}
\end{equation}
prove boundedness of both closed quotient maps at every stated finite
order. They impose no evenness; a one-sided cubic term is allowed.

The pressure law, enthalpy, and force are
\begin{equation}
 \begin{gathered}
 P(0)=0,\qquad P'(\rho)=\frac{4\rho^{2/3}}{3\sqrt{1+\rho^{2/3}}},
 \qquad h(\rho)=4\bigl(\sqrt{1+\rho^{2/3}}-1\bigr),\\
 \mathcal M[R]=4\pi R^2\partial_xP((4\pi R^2R_x)^{-1})+\frac{x}{R^2}.
 \end{gathered}
 \label{loc:force}
\end{equation}
We write $H_R=D\mathcal M[R]$. The admissible radii have
positive lower bounds and finite upper bounds for
\[
 R/z,\ R_z\quad\text{at the center},\qquad
 R,\ R_x\quad\text{in the interior},\qquad
 R,\ -R_\xi/\xi\quad\text{at the vacuum}.
\]

For later use we record the exact nonsingular coefficient fields.
At the vacuum put
\begin{equation}
 q=R_\xi/\xi<0,\quad a=-\frac5{4\pi R^2q},\quad
 A=\frac{4a^{2/3}}{3q^2\sqrt{1+\xi^2a^{2/3}}},\qquad
 B_0=-\frac{2\xi^2Aq^2}{R}+\frac{M-\xi^5}{R^2}.
 \label{loc:vacuum-coefficients}
\end{equation}
Within each local vacuum coefficient calculation, $A$ denotes
only the scalar coefficient in \eqref{loc:vacuum-coefficients};
elsewhere $A$, and in particular the time derivatives $A_j$, denote
the matched trajectory.
Direct substitution of $\rho=\xi^3a$ in \eqref{loc:force} gives
\begin{equation}
 \mathcal M[R]=-A\xi q_\xi+3Aq+B_0,
 \qquad \xi q_\xi=R_{\xi\xi}-q.
 \label{loc:vacuum-force}
\end{equation}
Thus the force is a smooth function of $\xi,R,q,R_{\xi\xi}-q$
on bounded sets where $R$ and $-q$ have positive lower bounds.
At the center let
\[
 \Phi(X)=R(|X|^3)\frac{X}{|X|}=r(X)X,
 \qquad J=D\Phi,\qquad B=J^{-T},\qquad
 \rho=\frac3{4\pi\det J}.
\]
Here and below a profile written as $R(z)$ in a center formula
means the pullback $R(z^3)$. The identities
\begin{equation}
 r(X)=\int_0^1t^2\operatorname{div}\Phi(tX)\dd t,
 \qquad
 \boldsymbol{\mathcal M}[\Phi]=B\nabla h(\rho)+r^{-2}X
 \label{loc:center-force}
\end{equation}
follow by integrating the derivative of $t^3r(tX)$ and by the
radial eigenvalue $R_z$ of $J$. The first formula avoids a
separate singular scalar quotient. The second uses at most two
Cartesian derivatives of $\Phi$. On the interior, with $p=R_x$,
\begin{equation}
 \mathcal M[R]=-\frac{P'((4\pi R^2p)^{-1})}{p^2}R_{xx}
              -\frac{2P'((4\pi R^2p)^{-1})}{R}+\frac{x}{R^2}.
 \label{loc:interior-force}
\end{equation}
Consequently $\mathcal M:\mathcal D^m\to\mathcal Y^m$ is
$C^\infty$ on the open set of radii satisfying the stated positivity
conditions. At each fixed finite spatial order, its differentials are
finite sums of products of the displayed derivative and quotient maps
with smooth functions of the positive factors.
Leibniz' rule and
$[fg]_\alpha\le\|f\|_\infty[g]_\alpha+[f]_\alpha\|g\|_\infty$
prove the norm bounds, including continuity of all finite
differentials. Applying the same product rule to spatial and temporal
differences
gives the parabolic source estimates used below.

\subsection{Compatible initial data at fixed mass}
\label{loc:euler-preparation}

Fix $J\ge3$ and the matching order
\begin{equation}
 K=2(J+1713)+86=2J+3512.
 \label{loc:matching-order}
\end{equation}
All constants in this subsection may depend on the selected
$\beta>0$, the finite matched family, and the prescribed cutoffs,
but not on $0<\lambda\le A_0$. We use the following profile bounds:
\begin{itemize}[leftmargin=2em]
\item physical $C^{58}$ enthalpy and $C^{57}$ velocity, and local
derivatives through order sixty;
\item endpoint source estimates through time order 22 and endpoint
order 23, to impose the compatibility conditions;
\item strong source estimates of total order 21, to bound the correction
in the localized weighted norms.
\end{itemize}

On the fixed initial-radius collar write
$a=R_\lambda-\lambda^2z^2$, $h=z^2p$, $p>0$. Choose twelve
enthalpy correction directions and eleven velocity correction directions
whose collar germs are
\begin{equation}
 \delta_{\eta_n}h=\frac{\lambda^4z^{2n+3}}{(2n+1)!}\quad(0\le n\le11),
 \qquad
 \delta_{\vartheta_n}V=
 \frac{\lambda^{5/2}z^{2n+3}}{(2n+1)!(2n+3)}\quad(0\le n\le10).
 \label{loc:columns}
\end{equation}
The cutoffs equal one on the enlarged collar. As the remaining
correction direction, choose a nonnegative enthalpy bump
$\lambda^{-1}\psi(a/\lambda)$ supported on a fixed positive
normalized interior interval, away from both endpoints. On its support,
$h\sim\lambda^{-1}$, $\rho'(h)\sim\lambda^{-2}$, and
$a^2\dd a\sim\lambda^3\dd r$. Its mass derivative is therefore
bounded above and below by positive constants.

Subtract a multiple of this bump from each collar enthalpy direction
so that its mass derivative vanishes at the base datum. These are
column operations on the parameter Jacobian and leave the endpoint
germs unchanged. Let $\Xi\in\mathbb R^{24}$ denote the parameters
and $x_\Xi(a)$ the enclosed-mass integral. A collar enthalpy
correction changes the mass tail by $O(\lambda^8z^{2n+6})$;
we include this contribution in the force.

We first compute the compatibility conditions as coefficients of
finite Laurent expansions. These coefficients are defined even when
the provisional time derivatives have no classical endpoint trace.
Freeze $\lambda$, put
$t=\tau/\lambda^2$, $R=a+\lambda^2u(t,z)$, and set
\[
 U=\frac{\lambda^2}{a}u,\qquad W=\frac{u_z}{2z},\qquad
 Q=(1+U)^{-4/3}(1-W)^{-2/3}.
\]
Mass conservation gives the exact enthalpy relation
\begin{equation}
 H+4=4\sqrt{1+\frac{h_\Xi(h_\Xi+8)}{16}Q}.
 \label{loc:finite-enthalpy}
\end{equation}
For $0\le\sigma\le1$, write $u_j=\partial_t^ju|_{t=0}$ and
$u=\sum_{j=0}^{24}t^ju_j/j!$, interpreted as a polynomial truncated
at time order 24. Determine these coefficients from
\begin{equation}
 u_{tt}=\frac{H_z}{2z(1-W)}
 -\frac{\lambda^2x_\Xi}{a^2}(1+U)^{-2}
 +(1-\sigma)\sum_{m=0}^{22}\frac{t^m}{m!}\lambda^{2m+2}S_m.
 \label{loc:finite-recursion}
\end{equation}
Here $S_m=\partial_\tau^m(A_{\tau\tau}+\mathcal M[A])$ is
pulled back once to the uncorrected initial-radius variable and
is independent of $\Xi$. At $\sigma=0,\Xi=0$, the recursion
gives the time derivatives of the matched forced trajectory; at
$\sigma=1$, it gives the derivatives determined by the Euler equation.
For a finite expansion in $z$, let $[f]_r$ denote the coefficient
of $z^r$. Set $\gamma_j=4$ for even $j$ and
$\gamma_j=5/2$ for odd $j$, and define the normalized endpoint source bound
\[
 \epsilon_{\rm ep}
 =\max_{\substack{0\le m\le22\\0\le r\le23}}
       \left|[\lambda^{\,2(m+2)-\gamma_{m+2}}S_m]_r\right|.
\]
To estimate this quantity, we use the bound for the full endpoint
source before restricting its derivative range. In the physical-to-material
identity \eqref{mat:endpoint-14}, the normalized scale exponent is
\[
 \chi_m=2(m+2)-\gamma_{m+2}+1-\tfrac32m,\qquad
 \chi_{2n}=n+1,\quad\chi_{2n+1}=n+3.
\]
Since $\chi_m\ge1$, the bounded coefficients in the change of
time variable, the $C^{28}$ bound \eqref{mat:endpoint-3}, and
the fixed-radius pullback \eqref{mat:endpoint-19} give
\[
 \epsilon_{\rm ep}
 \le C\lambda^{K-1}(1+|\log\lambda|)^{p_{\log}}.
\]
Here the exponent $K-1$ uses the full endpoint source bound,
before restriction to the weaker rectangle of derivative orders.

The quotient $f_z/z$ is regular only after the first Taylor coefficient
of $f$ vanishes. We therefore impose this condition on the 23 generated
time derivatives, together with the mass constraint, and order the
correction directions by the first time derivative they affect.

The mass constraint and the coefficients of $z^1$ in the finite
Laurent expansions of the time derivatives of orders 2 through 24
define a smooth map $F(\sigma,\Xi)\in\mathbb R^{24}$. Off its zero
set, $[u_j]_1$ denotes coefficient extraction even when negative
powers are present; it need not be a classical endpoint derivative.

\begin{lemma}
\label{loc:finite-map}
The finite map $F$ has the following properties.
\begin{enumerate}[label=(\roman*),leftmargin=2em]
\item \emph{Differential bounds.} Fix a sufficiently small parameter
ball before choosing the terminal scale. With the row normalization
specified below, the Jacobian $T=D_\Xi F(0,0)$ is invertible, and
the following bounds hold uniformly for $0\le\sigma\le1$:
\begin{equation}
 \begin{gathered}
 \|T^{-1}\|\le C\lambda^{-851},\qquad
 \|D_\Xi F\|+\|D_\Xi^2F\|\le C\lambda^{-4},\\
 \|\partial_\sigma F\|+\|D_\Xi\partial_\sigma F\|
       \le C\lambda^{-7/2}\epsilon_{\rm ep},\qquad
 \epsilon_{\rm ep}\le C\lambda^{K-1}(1+|\log\lambda|)^{p_{\log}}.
 \end{gathered}
 \label{loc:finite-map-bounds}
\end{equation}
\item \emph{Compatible curve.} There is a $C^1$ curve $\Xi(\sigma)$, $F(\sigma,\Xi(\sigma))=0$,
with $\Xi(0)=0$ and
\begin{equation}
 D_{\rm par}:=\sup_{0\le\sigma\le1}|\Xi'(\sigma)|
 \le C\lambda^{J+1775/2}(1+|\log\lambda|)^{p_{\log}}.
 \label{loc:parameter-size}
\end{equation}
Every member has exact mass, satisfies the stated positivity conditions
on the radius and density, and has compatible one-sided time derivatives
through order 24.
\end{enumerate}
\end{lemma}

\begin{proof}
\emph{Finite coefficient recurrence.}
Write $\widehat u_m=u_m/m!$, $[f]_{m,r}=[t^mz^r]f$, and
\[
 \mathscr P_bf=\sum_{0\le m\le24,\,m+r\le b}[f]_{m,r}t^mz^r,
 \qquad
 [fg]_{m,r}=\sum_{i=0}^m\sum_{a+b=r}[f]_{i,a}[g]_{m-i,b}.
\]
The lower supports proved below make these sums finite. Normalize
each scalar argument by its positive constant coefficient. Its remaining
factor is $1+v$, where every monomial of $v$ has positive weight $m+r$;
hence
\begin{equation}
 \mathscr P_b(1+v)^c
  =\mathscr P_b\sum_{k=0}^{b}\binom ck v^k,
 \qquad c\in\{-2,-4/3,-2/3,-1,1/2,3/2,1/5\}.
 \label{loc:finite-power-rule}
\end{equation}
For a smooth scalar function the finite Taylor polynomial gives the
same rule, with remainder of weight greater than $b$. No convergence
of an infinite power series is used.

Retain $h_\Xi$ through degree 25 and $V_\Xi$ through degree 24,
and set $\widehat u_0=0$, $\widehat u_1=V_\Xi$. Compute
$M_\Xi$ from the entire datum and the local mass tail from
\begin{equation}
 \begin{gathered}
 p_\Xi=h_\Xi/z^2,\qquad
 \mu_\Xi=\{p_\Xi(z)(z^2p_\Xi(z)+8)/16\}^{3/2},\\
 x_\Xi(z)=M_\Xi-8\pi\lambda^2\int_0^z
            (R_\lambda-\lambda^2s^2)^2s^4\mu_\Xi(s)\dd s .
 \end{gathered}
 \label{loc:finite-tail-rule}
\end{equation}
A requested coefficient is obtained by integrating its Taylor monomial,
$s^r\mapsto z^{r+1}/(r+1)$. For $0\le m\le22$, the known
$\widehat u_0,\ldots,\widehat u_{m+1}$ then determine
\begin{equation}
 \begin{aligned}
 U&=\mathscr P_{23}(\lambda^2u/a),\qquad
 W=\mathscr P_{23}(u_z/(2z)),\\
 Q&=\mathscr P_{23}\{(1+U)^{-4/3}(1-W)^{-2/3}\},\\
 H&=\mathscr P_{25}\{4(1+h_\Xi(h_\Xi+8)Q/16)^{1/2}-4\},\\
 \Phi&=\mathscr P_{23}\left\{
       \frac{H_z}{2z}(1-W)^{-1}
       -\frac{\lambda^2x_\Xi}{a^2}(1+U)^{-2}\right\},\\
 [\widehat u_{m+2}]_r
  &=\frac{[\Phi]_{m,r}
       +(1-\sigma)\lambda^{2m+2}[S_m]_r/m!}{(m+1)(m+2)}.
 \end{aligned}
 \label{loc:complete-coefficient-recursion}
\end{equation}
Use the last line for $m+2+r\le25$. The source remains fixed in
the original radius chart, independently of $\Xi$. Thus parameter
differentiation of these finite convolutions gives the full Jacobian,
including gravity and the mass tail.

\emph{Truncation.}
The lower support induction is
\begin{equation}
 \operatorname{val}_z\widehat u_m\ge\min(0,3-m),\qquad
 \operatorname{val}_z W_m\ge1-m,\qquad
 \operatorname{val}_z H_m\ge3-m\quad(m\ge1).
 \label{loc:finite-support-induction}
\end{equation}
Indeed,
\[
 U_0=W_0=0,\qquad H_0=h_\Xi,\qquad
 V_\Xi=v_*(\lambda)+z^2v_\Xi(z).
\]
The constant velocity is killed by $\partial_z/z$. Every positive-time
term of $U,W$ has weight at least one. Rationalizing the square root
writes $H$ as $h_\Xi$ times such a factor; consequently its
positive-time pressure terms retain weight at least one after
$\partial_z/z$. Gravity has no smaller weight. The constant allowed
in $\widehat u_2$ is killed in $W_2$, and an ordinary source term
at time $m$ has weight at least $m$ before two time integrations.
These observations close the support induction.

For the upper truncation, an omitted radius monomial of weight $w>25$
can enter pressure through $U$ or $W$. The largest loss follows
\[
 w\xrightarrow{\ \partial_z/z\ }w-2
  \xrightarrow{\ \times h_\Xi\ }w
  \xrightarrow{\ \partial_z/z\ }w-2>23.
\]
Its contribution through $(1-W)^{-1}$ also has weight greater than
23, because the other pressure factor has nonnegative weight. An
omitted enthalpy term has weight greater than 23 after
$\partial_z/z$. All other products add nonnegative weights and
the time integrations add two. In the mass tail the weight change is
\[
 -2+4+1=3
\]
from division by $z^2$ and integration against $z^4\dd z$.
Thus the discarded enthalpy terms cannot affect retained gravity
coefficients. Together with $0\le r\le23-m$ for the source, this
proves that the retained enthalpy degree 25 and velocity degree 24
suffice. The finite power rule introduces no negative-weight factors.

\emph{Parameter bounds.}
Parameter differentiation leaves the time and spatial indices unchanged.
For a scalar composition and a reciprocal, use
\[
 \begin{aligned}
 D_a\phi(f)&=\phi'(f)D_af,\\
 D_aD_b\phi(f)&=\phi''(f)D_afD_bf+\phi'(f)D_aD_bf,\\
 D(X^{-1})&=-X^{-2}DX,\\
 D_1D_2(X^{-1})&=2X^{-3}D_1X D_2X-X^{-2}D_1D_2X.
 \end{aligned}
\]
The differentiated recurrences therefore have the same truncation and
support bounds, through three parameter derivatives, with no additional
spatial derivative. A source derivative starts with
$-\lambda^{2m+2}[S_m]_r/m!$ in the force coefficient; the last
line of \eqref{loc:complete-coefficient-recursion} supplies its two
time-integration factors.

The apparently unbounded constant velocity is harmless, since
\[
 \begin{gathered}
 |v_*|\le C\lambda^{-1/2},\qquad
 \|v_\Xi\|_{C^{27}}\le C\lambda^{1/2},\\
 U_1=(\lambda^2/a)V_\Xi=O(\lambda^{1/2}),\\
 W_1=V_{\Xi,z}/(2z)=O(\lambda^{1/2}).
 \end{gathered}
\]
All remaining input coefficients, their first three parameter
derivatives, and the reciprocals of their leading factors are uniformly
bounded on the parameter ball. For the retained finite sums,
\[
 \sum_{m,r}|[fg]_{m,r}|
 \le\left(\sum_{m,r}|[f]_{m,r}|\right)
      \left(\sum_{m,r}|[g]_{m,r}|\right),\qquad
 \sum_{m,r}|[\partial_zf]_{m,r}|
 \le\sum_{m,r}|r|\,|[f]_{m,r}|.
\]
The finite power rule is bounded by
$\sum_{k=0}^{25}|\binom ck|C^k$ when the coefficient sum of its
argument is at most $C$. Applying these bounds successively to
\eqref{loc:complete-coefficient-recursion}, with the positive divisors
$(m+1)(m+2)$, gives uniform coefficient and parameter bounds in 23 steps.
After one source differentiation, each product contains exactly one
source coefficient; the resulting bound is therefore linear in the
source norm.

\emph{Mass coordinates and row normalization.}
On the collar,
\[
 y=\bigl((M_\Xi-x_\Xi)/\lambda^4\bigr)^{1/5}=zc_\Xi(z),
 \qquad
 c_\Xi(0)^5=\frac{8\pi}{5}(R_\lambda/\lambda)^2(p_\Xi(0)/2)^{3/2}.
\]
The correction directions preserve $p_\Xi(0)=p_*$, so
$c_\Xi(0)=c_*>0$ is independent of $\Xi$; both have uniform
positive bounds. Write
$c_\Xi(z)=c_*\{1+\sum_{r=1}^{23}c_{\Xi,r}z^r\}$,
$\zeta=y/c_*$ and
$z=B_\Xi(\zeta)=\sum_{r=1}^{24}b_r\zeta^r$. Then $b_1=1$ and
\begin{equation}
 b_n=-\sum_{k=1}^{n-1}c_{\Xi,k}
              [\zeta^n]B_\Xi(\zeta)^{k+1},\qquad 2\le n\le24.
 \label{loc:finite-pullback-rule}
\end{equation}
The right side uses only $b_1,\ldots,b_{n-1}$. Substitution into
negative Laurent powers uses the same finite power rule. The smallest
power in the twenty-fourth time derivative is $-21$, so its first
coefficient needs at most $b_{23}$.

Before imposing $M_\Xi=M$, regard $M_\Xi$ as a parameter in the
force. On the mass constraint, $y=((M-x)/\lambda^4)^{1/5}$ is common
to the family; differentiating the change of coordinates gives
\eqref{loc:source-transport}.

The components of $F$ are defined using the original $z$ coordinate:
\[
 M_\Xi-M,\qquad
 F_j=c_*^{-1}\lambda^{-\gamma_j}[u_j]_1,
 \quad 2\le j\le24.
\]
If the preceding constraints vanish, the next time derivative has an
ordinary one-sided expansion. Its first coefficient in $y$ is then
$c_*^{-1}[u_j]_1$. Induction therefore identifies the compatibility
conditions in the two coordinates. Away from their common zero set,
negative Laurent powers can contribute to the first coefficient after
substitution, so the two coefficient maps need not agree.

Use maximum norms on the rows and parameters. The largest row factor
is $\lambda^{-4}$, while
\[
 \lambda^{2m+2}S_m
 =\lambda^{\gamma_{m+2}-2}
       \{\lambda^{2(m+2)-\gamma_{m+2}}S_m\},\qquad
 \gamma_{m+2}-2\ge\tfrac12.
\]
The finite coefficient bounds thus give $C\lambda^{-4}$ for
parameter derivatives and $C\lambda^{-7/2}\epsilon_{\rm ep}$
for derivatives containing one source differentiation, including
$D_\Xi\partial_\sigma F$.

For mass derivatives of orders $1\le r\le3$, the interior density,
bump and volume factors give
\[
 \lambda^{r-3}\lambda^{-r}\lambda^3=1,
\]
and the collar integrand is bounded by
$C\lambda^{4r+4}z^{r+4}\dd z$. These integrable bounds and the
same finite recurrence also yield
\[
 \|D_\Xi^3F\|\le C\lambda^{-4},\qquad
 \|D_\Xi^2\partial_\sigma F\|
       \le C\lambda^{-7/2}\epsilon_{\rm ep}.
\]
No additional endpoint derivative is used.

\emph{Triangular Jacobian.}
Linearize the complete recurrence at the matched forced trajectory.
For a displacement $f$, \eqref{loc:finite-enthalpy} gives at time zero
\[
 D_uH[f]=\frac{h(h+8)}{2(h+4)}
        \left(-\frac{4\lambda^2f}{3a}+\frac{f_z}{3z}\right).
\]
Since $h=p_*z^2+O(z^3)$, the only terms lowering the endpoint
degree by two in the pressure differential are
\[
 \frac1{2z}\partial_z\left(\frac{p_*z}{3}f_z\right)
       +p_*\frac{f_z}{2z}
   =L_0f,\qquad
 L_0=\frac{p_*}{6}(\partial_z^2+4z^{-1}\partial_z).
\]
The $U$ term and gravity lower degree by at most zero.  A
nonconstant spatial coefficient in the displayed principal part
improves weight by at least one.  The time derivatives of positive order at the base are ordinary
compatible functions, so their $z$-quotients are regular. A coefficient
with time index $k\ge1$ improves weight by at least $k$.  Thus all terms in the linearized recursion preserve the weight of the parameter variation, and equality is possible only for the constant,
time-zero $L_0$ followed by two time integrations.

An enthalpy variation $z^p$ first contributes $(p/2)z^{p-2}$ to the
second time derivative through the pressure term. A velocity variation
$z^p$ is already present in the first time derivative.

Each collar enthalpy direction has zero mass derivative, so
$D x_\Xi=-D(M_\Xi-x_\Xi)$. In
\eqref{loc:finite-tail-rule} the enthalpy variation gives
$D\mu_\Xi=O(z^{p-2})$, hence a mass tail $O(z^{p+3})$.
Its first contribution to gravity has weight $p+5$, five above the
pressure contribution. It therefore cannot contribute to a
compatibility condition of lower weight.
The interior mass column has zero enthalpy and velocity germs;
its gravity contribution is retained with no support restriction
beyond the mass row being first.

Consequently an $\eta_n$ variation has no radius monomial of weight
less than $2n+3$, and a $\vartheta_n$ variation has none below
$2n+4$. Since the obstruction in row $j$ has weight $j+1$, the
possible nonzero entries of the $24$ by $24$ Jacobian are as follows:
\begin{enumerate}[label=(\roman*),leftmargin=2em]
\item The mass row $M_\Xi-M$ has only the interior mass column,
whose entry is $D_{\rm mass}M_\Xi>0$.
\item In row $F_{2n+2}$, $0\le n\le11$, only the mass column and
the columns
\[
 \eta_0,\ldots,\eta_n;\qquad \vartheta_0,\ldots,\vartheta_{n-1}
\]
can be nonzero. The diagonal entry is $c_*^{-1}d_n$.
\item In row $F_{2n+3}$, $0\le n\le10$, only the mass column and
the columns
\[
 \eta_0,\ldots,\eta_n;\qquad \vartheta_0,\ldots,\vartheta_n
\]
can be nonzero. The diagonal entry is $c_*^{-1}e_n$.
\end{enumerate}
An empty velocity range means no velocity column. All unlisted
entries vanish; the remaining entries, including their lower-order
terms, are determined by the differentiated recurrence
\eqref{loc:complete-coefficient-recursion}.
Order the columns and rows as
\[
 (\mathrm{mass},\eta_0,\vartheta_0,\ldots,
       \eta_{10},\vartheta_{10},\eta_{11}),\qquad
 (\mathrm{mass},F_2,\ldots,F_{24}).
\]
All 276 entries above the diagonal vanish, so the Jacobian is lower triangular.
The identity
$(\partial_z^2+4z^{-1}\partial_z)z^p=p(p+3)z^{p-2}$
gives its non-mass diagonal, apart from $c_*^{-1}$, as
\[
 \begin{aligned}
 d_n&=p_*^n\frac{(n+1)(n+2)(2n+3)}{4\,6^n},
                                  &&0\le n\le11,\\
 e_n&=p_*^{n+1}\frac{(n+1)(n+2)(n+3)}{6^{n+1}},
                                  &&0\le n\le10.
 \end{aligned}
\]
For example $d_0=3/2$, $d_1=5p_*/4$,
$e_0=p_*$, and $e_1=2p_*^2/3$. The factorial normalization is the one in \eqref{loc:columns}.

The first lower entries also follow from the exact force, and show
where the nonconstant background is retained. Write its initial germs
as $h=p_*z^2+h_3z^3+\cdots$ and
$V=v_0+v_2z^2+v_3z^3+\cdots$, and set
\[
 a_0=R_\lambda,\qquad
 \ell_1=-\frac{4\lambda^2v_0}{3a_0}+\frac83v_2,
 \qquad g_1=-\frac{4\lambda^2v_0}{3a_0}+\frac53v_2.
\]
Differentiating \eqref{loc:finite-enthalpy} once at $t=0$ gives
\[
 H_1=\frac{h(h+8)}{2(h+4)}
           \left(-\frac43U_1+\frac23W_1\right).
\]
It follows that the coefficient of $h_3$ in $[u_3]_1$ is
$-2\lambda^2v_0/a_0+5v_2/2$, and the coefficient of $v_3$
is $3p_*$. More generally, at $(\sigma,\Xi)=(0,0)$ the
first subdiagonal entries of the Jacobian are
\begin{equation}
 \begin{aligned}
 \partial_{\eta_n}F_{2n+3}
   &=c_*^{-1}\lambda^{3/2}d_n
                 \{g_1+n(n+2)\ell_1\},&&0\le n\le10,\\
 \partial_{\vartheta_n}F_{2n+4}
   &=c_*^{-1}\lambda^{-3/2}e_n
                    (n+1)(n+2)\ell_1,&&0\le n\le10.
 \end{aligned}
 \label{loc:first-lower-entries}
\end{equation}
Indeed, the principal coefficient and its logarithmic derivative are
\[
 \ell(t)=\frac{p_*}{6}(1+U(t,0))^{-4/3}(1-W(t,0))^{-8/3},
 \qquad \frac{\ell'(0)}{\ell(0)}=\ell_1.
\]
The coefficient of the enthalpy variation has logarithmic derivative
$g_1$. Distributing one additional time derivative over this coefficient
and the normal operators gives the factors
\[
 \sum_{k=1}^n(2k+1)=n(n+2),\qquad
 \sum_{k=1}^{n+1}2k=(n+1)(n+2).
\]
The support induction excludes higher endpoint coefficients from these
entries and gives the stated zero pattern.
The remaining lower entries, including the contributions of gravity
and the mass column, are given by
\eqref{loc:complete-coefficient-recursion} and satisfy its polynomial
majorant. In particular, the second line of
\eqref{loc:first-lower-entries} is covered by the common
$O(\lambda^{-4})$ bound despite its negative scale power.

\emph{Bound for the inverse Jacobian.}
Write $T=D+N$, with $D$ diagonal
and $N$ strictly lower triangular. The positive mass derivative
and the displayed diagonal entries give $\|D^{-1}\|\le C$.
Since $(D^{-1}N)^{24}=0$,
\begin{equation}
 T^{-1}=\sum_{k=0}^{23}(-D^{-1}N)^kD^{-1},\qquad
 \|T^{-1}\|\le C\sum_{k=0}^{23}\lambda^{-4k}
                      \le C\lambda^{-92}.
 \label{loc:triangular-inverse}
\end{equation}
For the preparation and initial energy estimates we use the common
loss exponent $851$. Since $0<\lambda\le1$, the bound
$C\lambda^{-92}$ for the full inverse implies
$C\lambda^{-851}$. The matching order \eqref{loc:matching-order}
dominates this loss, the two parameter derivatives, and the subsequent
powers lost in the weighted norm and energy estimates. These numerical exponents are sufficient
for the argument and are not optimized.

\emph{Solution of the constraints.}
Consider $\Xi\mapsto\Xi-T^{-1}F(\sigma,\Xi)$ on the ball
$|\Xi|\le r_\lambda$, where
\[
 r_\lambda=C\lambda^{-851-7/2}\epsilon_{\rm ep}
\]
with a sufficiently large fixed $C$. By
\eqref{loc:finite-map-bounds}, its derivative is bounded by
\[
 C\lambda^{-851}
       (\lambda^{-4}r_\lambda+\lambda^{-7/2}\epsilon_{\rm ep})
 \le C\{\lambda^{K-3421/2}+\lambda^{K-1711/2}\}
                     (1+|\log\lambda|)^{p_{\log}}.
\]
Both powers are positive for \eqref{loc:matching-order}. Decrease
$A_0$ so that this bound is at most $1/2$, uniformly in $\sigma$.
Since $F(0,0)=0$, the source-derivative bound also gives
\[
 |T^{-1}F(\sigma,0)|
 \le\|T^{-1}\|\int_0^\sigma|\partial_\sigma F(v,0)|\dd v
 \le r_\lambda/2.
\]
Thus the map preserves the ball and is a contraction there.
Differentiating the equation for its fixed point yields
\[
 \Xi'=-(D_\Xi F)^{-1}\partial_\sigma F,\qquad
 |\Xi'|\le C\lambda^{K-1711/2}(1+|\log\lambda|)^{p_{\log}}.
\]
Finally,
\[
 K-\frac{1711}{2}
   =J+\frac{1775}{2}+(J+1769),
\]
so $0<\lambda\le1$ gives \eqref{loc:parameter-size}.

\emph{Compatibility and regularity.}
Induct on $j$. If all earlier radius derivatives have zero first
coefficient, their quotients in \eqref{loc:quotients} are ordinary
one-sided functions, and the enthalpy derivatives are $z^2$ times
ordinary factors. The recurrence then generates a regular $j$-th
derivative. Its first Laurent coefficient is its first Taylor coefficient;
the $j$-th constraint makes it vanish and permits the next step.

The force loses at most two one-sided derivatives. The physical data
have $C^{58}$ enthalpy and $C^{57}$ velocity by
\eqref{mat:units-1}, and Subsection~\ref{app:profile-endpoint} gives,
at fixed positive scale,
\[
 \partial_\tau^m(A_{\tau\tau}+\mathcal M[A])\in C^{58-m},
 \quad 0\le m\le22,\qquad
 R_j\in C^{54-2j}\ (0\le j\le24),\quad R_{24}\in C^6.
\]
Here the source retains its regularity after the initial-radius
pullback, and the last two conclusions follow by the finite mass
inverse and the recurrence. The constants may depend on the positive
scale. Differentiating the vanishing constraints
gives the same first trace for the total parameter derivatives.

The collar enthalpy change is $\lambda^4z^3O(r_\lambda)$,
so $h/z^2$ remains positive. The interior correction is a relative
$O(r_\lambda)$ change and the center is unchanged. Since the
physical initial radius remains $a$, the positive mass derivative
gives a positive inverse radius map with exact total mass. Thus all
positivity conditions persist. Only on this curve are the formal
coefficient constraints used as classical trace conditions.

\end{proof}

\subsection{Estimates for the correction in mass coordinates}
\label{loc:common-mass-correction}

We now estimate the correction at fixed mass. Differentiating in
$\sigma$ changes both the initial data and the inverse mass map.
We first estimate this change of coordinates and the transported source,
then differentiate the force recurrence to bound the corrected data.

\emph{Variation of the mass coordinate.}
Write $a_\sigma=x_\sigma^{-1}$ for the inverse initial mass map.
Since total mass is constant, its exact variation is
\begin{equation}
 \dot a_\sigma=-\frac{\dot x_\sigma}{(x_\sigma)_a}\circ a_\sigma,
 \qquad
 \dot x_\sigma(a)=-4\pi\int_a^{R_\lambda}
                r^2\rho'(h_\sigma(r))\dot h_\sigma(r)\dd r.
 \label{loc:mass-variation}
\end{equation}
On the bounded collar, $\dot h=\lambda^4z^3O(|\Xi'|)$, whence
\[
 \dot x=\lambda^8z^6O(|\Xi'|),\qquad
 x_a\sim\lambda^2z^3,\qquad
 \dot a=\lambda^6z^3O(|\Xi'|).
\]
The factors obey these bounds through 28 ordinary derivatives by
finite differentiation of the positive mass inverse. Write
$\sigma_{\rm mass}$ for the coefficient of the interior mass correction
and $\eta=(\eta_0,\ldots,\eta_{11})$ for the enthalpy parameters.
The mass equation gives
$ |\sigma_{\rm mass}'|\le C\lambda^8|\eta'|$.
Thus the variation is zero near the unchanged center, has size
$O(\lambda^8|\Xi'|)$ in fixed mass on the compact interior,
and has long-chart shift $O(\lambda^4y_8^{-7}|\Xi'|)$.
In particular
\begin{equation}
 \|\dot a_\sigma\|_{X^r}+|\dot a_\sigma|_{\mathcal C^r}
       \le C\lambda^6|\Xi'|,\qquad r\le26.
 \label{loc:mass-shift-bound}
\end{equation}
Here $\mathcal C^r$ denotes the enlarged ordinary norms in the chart coordinates,
with the full vector center component and no volume prefactor.

\emph{The transported source.}
Write the physical time derivatives at common mass as
$Y_j(\sigma,x)=R_j(\sigma,a_\sigma(x))$, and set
$V_j=\partial_\sigma Y_j$. Since the initial radius is
$Y_0=a_\sigma$, its variation is $V_0=\dot a_\sigma$.
For the transported source, put
$\widetilde S_m(\sigma,x)=S_m(a_\sigma(x))$.
The strong source estimate and the triangular change of time
derivatives give
\begin{equation}
 E_{\rm src}:=\max_{m+q\le21}
       \|\lambda^{-1+3m/2}S_m\|_{q,a}
 \le C\lambda^{(K-7)/2}(1+|\log\lambda|)^{p_{\log}}.
 \label{loc:correction-source}
\end{equation}
The loss of one power relative to the strong-source norm is
explicit: on the same cutoffs,
$\mathcal S_q(f)\le\|f\|_{q,a}\le\lambda^{-1}\mathcal S_q(f)$.
These are ordinary source norms; no vacuum first trace is required.
The source transport requires only the separate local growth bound
$|S_m|_{\mathcal C^{q+1}}\le C\lambda^{-2-3m/2}$, $m+q\le21$,
whose total radius derivative count is $m+q+3\le24$. Formula
\eqref{loc:mass-variation} and the fundamental theorem of calculus give
\begin{equation}
 \|\partial_\sigma\widetilde S_m\|_{q,a}
 \le C\lambda^{2-3m/2}|\Xi'|,\qquad m+q\le21.
 \label{loc:source-transport}
\end{equation}
The long-chart volume is included, since
$\lambda^2(\int_c^{C\lambda^{-1/2}}y^7\dd y)^{1/2}\le C$.

\emph{Differentials of the force.}
We estimate one direction in a weighted norm and the others in
ordinary chart norms. For each $a\in\{1,\ldots,p\}$,
\begin{equation}
 \|D^p\mathcal M[Y_0][f_1,\ldots,f_p]\|_{q,a}
 \le C_{p,q}\lambda^{-2p-2}\|f_a\|_{X^{q+2}}
                    \prod_{i\ne a}|f_i|_{\mathcal C^{q+2}}.
 \label{loc:family-force-estimate}
\end{equation}
Here $1\le p\le22$, $0\le q\le21$, and all directions
have the displayed compatible ordinary regularity. Its counterpart in the chart norms is used at the finite orders through 26 below.
Indeed on an edge, $Y_0=\lambda r$,
$(Y_0)_y=-\lambda^2y\gamma$, and the exact Hessian is
\[
 \lambda^4H_{Y_0}f=-A_n(f_{yy}+(n-1)f_y/y)-B_nf_y/y+V_nf.
\]
Its coefficients depend smoothly on the positive factors $r,\gamma$,
affine in $y\gamma_y$; every explicit $\lambda y^2$ is bounded
on the entire long chart. A further variation inserts
\[
 \delta r=f/\lambda,\quad
 \delta\gamma=-\lambda^{-2}f_y/y,\quad
 \delta(y\gamma_y)=-\lambda^{-2}(f_{yy}-f_y/y).
\]
Thus $p-1$ variations cost at most $\lambda^{-2(p-1)}$
in addition to $\lambda^{-4}$, with only two grouped derivatives
per direction. Apply the closed one-sided quotient of
Lemma~\ref{at:natural-quotient} and the localized derivative estimate
\eqref{at:localized-jet} to the designated argument; estimate the
other arguments in their displayed ordinary chart norms.

At the center, use the full vector formula \eqref{loc:center-force}.
The scalar radial factor is its divergence integral, whose
$k$-th derivative has the integrable dilation factor
$t^{k+2-3/2}$. On the interior use \eqref{loc:interior-force}.
These calculations prove \eqref{loc:family-force-estimate} and
the corresponding chart estimate, with the high derivatives of the
background included in the constants.

For chart orders beyond 21, differentiate the integral identity
\eqref{loc:quotients} in $C^r$ to bound $f_y/y$ by the
$C^{r+2}$ norm. The weighted quotient lemma is used only at
source orders $q\le21$.

\emph{The corrected time derivatives.}
For $2\le j\le24$, the compatible recurrence at common mass is
\[
 Y_j=-\sum_{\pi\in\Pi_{j-2}}D^{|\pi|}\mathcal M[Y_0]
            [Y_{|B|}:B\in\pi]+(1-\sigma)\widetilde S_{j-2}.
\]
Here $\Pi_m$ is the set of partitions of a labeled $m$-element
set; the empty partition gives the force itself. The local bound
\begin{equation}
 |Y_i|_{\mathcal C^q}\le C_i\lambda^{2-7i/2},\qquad
 1\le i\le24,\quad i+q\le26
 \label{loc:native-corrected-rows}
\end{equation}
follows by induction on $i$. For a partition with $p$ blocks,
$\sum_B|B|=i-2$, and hence
\[
 \lambda^{-2p-2}\prod_B\lambda^{2-7|B|/2}
                         =\lambda^{5-7i/2}.
\]
The source has the better local growth $\lambda^{1-3i/2}$.

Indeed, for $i\ge2$, put $m=i-2$. With time derivatives taken before the
uncorrected initial-radius pullback,
\eqref{mat:native-eq:wave2-native-matched-rows} and
\eqref{mat:native-eq:wave2-direct-profile-force-source} give
\[
 \begin{gathered}
 |\overline R_i|_{q,\infty,a}
 +|\partial_\tau^m\mathcal M[A]|_{q,\infty,a}
 \le C\lambda^{1-3i/2},\\
 i+q\le26<27,\quad m+q\le24<25.
 \end{gathered}
\]
For the force, iterate $\partial_\tau=\lambda^{-3/2}\partial_s$;
after factoring out $\lambda^{-3m/2}$, the scalar coefficients
are bounded. The positive mass inverses in \eqref{loc:mass-variation}
transfer this bound to $\widetilde S_m$.
These are growth bounds; the mixed decay estimate
\eqref{loc:correction-source} is used only for $m+q\le21$.
Every earlier derivative required by the recurrence satisfies
$|B|+q+2\le i+q\le26$, so induction gives the estimate through
time order 24.

Differentiating the recurrence at fixed mass gives exactly
\begin{align}
 V_j={}&-\sum_{\pi\in\Pi_{j-2}}D^{p+1}\mathcal M[Y_0]
                  [V_0,Y_{|B|}:B\in\pi]\notag\\
 &-\sum_{\pi\in\Pi_{j-2}}\sum_{B_*\in\pi}
       D^p\mathcal M[Y_0][V_{|B_*|},Y_{|B|}:B\ne B_*]
       -\widetilde S_{j-2}+(1-\sigma)\partial_\sigma\widetilde S_{j-2},
 \quad p=|\pi|.
 \label{loc:homotopy-row}
\end{align}
Differentiating the vanishing constraints gives the compatible
first-derivative trace for each total variation $V_i$. This trace
need not hold for the variations of the datum and of the coordinate
separately, so the quotient is taken only after those two variations
have been combined. We estimate the individual force and source
terms in \eqref{loc:homotopy-row} in the ordinary norms; their sum
has the compatible trace and is therefore measured in $X^q$.

The initial bounds from \eqref{loc:mass-shift-bound} and the velocity
correction directions are
\[
 \|V_0\|_{X^{23}}\le C\lambda^6D_{\rm par},\qquad
 \|V_1\|_{X^{22}}\le C\lambda^{5/2}D_{\rm par}.
\]
Induction in \eqref{loc:homotopy-row} then gives
\begin{equation}
 \|V_j\|_{X^q}\le C_j\lambda^{5-7j/2}(D_{\rm par}+E_{\rm src}),
 \qquad 1\le j\le23,\quad j+q\le23.
 \label{loc:homotopy-graph}
\end{equation}
For a partition with $p$ blocks and a varied argument of time order
$k$, the two force terms have scale factors
\[
 \begin{aligned}
 \lambda^{-2p-4}\lambda^6\lambda^{2p-7(j-2)/2}
     &=\lambda^{9-7j/2},\\
 \lambda^{-2p-2}\lambda^{5-7k/2}
             \lambda^{2(p-1)-7(j-2-k)/2}
     &=\lambda^{8-7j/2}.
 \end{aligned}
\]
The source bounds \eqref{loc:correction-source} and
\eqref{loc:source-transport} give, respectively,
\[
 \lambda^{4-3j/2}E_{\rm src}
   =\lambda^{5-7j/2}\lambda^{2j-1}E_{\rm src},\qquad
 \lambda^{5-3j/2}D_{\rm par}
   =\lambda^{5-7j/2}\lambda^{2j}D_{\rm par}.
\]
All exponents improve the asserted bound. The required indices satisfy
\[
 k+q+2\le j+q\le23\quad(k\le j-2),\qquad
 (j-2)+q\le21,
\]
so every force and source estimate is used within its stated range.

\emph{The mixed Sobolev and pressure bounds.}
Integrating in $\sigma$ and converting the physical time derivatives
to derivatives in rescaled time gives, for the Euler-compatible datum,
\begin{equation}
 \begin{gathered}
 \psi_j=\lambda^{-4}(R_j^0-A_j),\qquad
 \|\psi_j\|_{X^q}\le C\lambda^{-L_j}D_\lambda,\quad j+q\le23,\\
 D_\lambda=D_{\rm par}+E_{\rm src},\qquad
 L_0=L_1=3, L_j=2j-1 (j\ge2).
 \end{gathered}
 \label{loc:prepared-row-bound}
\end{equation}
The coefficients in the change of time variable are bounded for
the selected $b\ge\beta>0$.
Write $\|f\|_{\dot V_A}:=\mathcal P_A[f]^{1/2}$ for the
closed pressure seminorm of Proposition~\ref{pr:pressure-domain},
without the Hilbert norm term. The maximum loss is
$L_{23}=45$. The pressure upper bound
$\|f^\perp\|_{\dot V_A}\le C\lambda^{-1/2}\|f\|_{X^1}$
adds only $L_{22}+1/2=87/2<45$. Therefore
\begin{equation}
 \begin{aligned}
 \sum_{j+q\le23}\|\psi_j\|_{X^q}^2
   +\sum_{j\le22}\|\psi_j^\perp\|_{\dot V_A}^2
 &\le C\lambda^{-90}D_\lambda^2,\\
 \left(\sum_{j+q\le23}\|\psi_j\|_{X^q}^2
          +\sum_{j\le22}\|\psi_j^\perp\|_{\dot V_A}^2\right)^{1/2}
 &\le C\lambda^{J+1685/2}(1+|\log\lambda|)^{p_{\log}}.
 \end{aligned}
 \label{loc:prepared-graph}
\end{equation}
Thus the correction is small both in the mixed norms and in the
pressure seminorm. Coercivity is proved separately in the energy
argument.

The corrected data and reference profile now use the same mass
coordinate, and \eqref{loc:prepared-graph} is uniform in the terminal
scale. We next impose viscous compatibility by the implicit function
theorem at each fixed positive scale; the inverse bound may depend on
that scale.

\subsection{Compatibility conditions for the regularized equation}
\label{loc:viscosity-selector}

Fix a positive terminal scale and a finite $\Lambda\ge0$.
For an initial pair $(r_0,r_1)$, generate the regularized time derivatives by
\begin{equation}
 r_j=-\left.\partial_\tau^{j-2}\mathcal M[R]\right|_0
       -\kappa\left.\partial_\tau^{j-1}\mathcal M[R]\right|_0
       -\kappa\Lambda r_{j-1},\qquad 2\le j\le25.
 \label{loc:initial-recursion}
\end{equation}
The derivatives on the right are finite Bell polynomials in the
previously generated time derivatives; their highest positive index
is $j-1$. Thus the recursion is triangular also for positive $\kappa$.

\begin{proposition}
\label{loc:selector}
Fix the Euler-compatible datum of Lemma~\ref{loc:finite-map}.
There are $\kappa_1(\lambda,\Lambda)>0$ and a $C^1$ family of
initial data with the same total mass such that their generated time
derivatives satisfy
\begin{equation}
 \kappa\longmapsto
 (r_0^\kappa,\ldots,r_{24}^\kappa,r_{25}^\kappa)
 \ \in C^1\left([0,\kappa_1];
       \prod_{j=0}^{24}\mathcal D^{48-2j}\times\mathcal Y\right).
 \label{loc:initial-cone}
\end{equation}
The vacuum trace satisfies
\[
 \partial_\xi r_j^\kappa(0)=0,\qquad 0\le j\le24.
\]
No derivative trace is imposed on $r_{25}^\kappa$. At the fixed
scale, the family differs from the Euler-compatible time derivatives
by $O(\kappa)$ in the displayed product space.
\end{proposition}

\begin{proof}
All bounds below are at a fixed positive scale.

\emph{Mass-coordinate regularity.}
Near the vacuum write $a=a_b-y^2$. For the $C^{58}$ enthalpy,
$h_\Xi(0)=h_{\Xi,y}(0)=0$ gives
\[
 p_\Xi=h_\Xi/y^2=\int_0^1(1-t)h_{\Xi,yy}(ty)\dd t,
 \qquad
 \mu_\Xi=\left\{\frac{p_\Xi(y)(y^2p_\Xi(y)+8)}{16}\right\}^{3/2}.
\]
Both are $C^1_\Xi C^{56}_y$, and $\mu_\Xi>0$. The exact tail is
\begin{equation}
 M_\Xi-x_\Xi(y)=y^5c_\Xi(y)^5,\qquad
 c_\Xi(y)^5=8\pi\int_0^1t^4(a_b-y^2t^2)^2\mu_\Xi(ty)\dd t.
 \label{loc:finite-mass-inverse}
\end{equation}
On a common shorter collar $\partial_y(yc_\Xi)>0$. Differentiating
$y_\Xi(\xi)c_\Xi(y_\Xi(\xi))=\xi$ gives a $C^{56}$ inverse and
\[
 \partial_{\Xi_a}y_\Xi
 =-\left(\frac{\partial_{\Xi_a}(yc_\Xi)}
                    {\partial_y(yc_\Xi)}\right)\circ y_\Xi\in C^{55}.
\]
At each higher order the unknown inverse derivative has the same
positive coefficient; all other terms contain lower derivatives.
The difference recurrence gives continuity. Thus radius and velocity
are $C^1_\Xi C^{55}$, in particular
$C^1_\Xi(C^{51}\times C^{49})$. Interior inversion is ordinary,
and the supported corrections preserve the Cartesian center class.
Before imposing the mass equation use $(M_\Xi-x)^{1/5}$ and keep
$M_\Xi$ as a force parameter; on its zero set the coordinate is the
fixed $\xi$.

\emph{Finite map and its Jacobian.}
Set
\[
 \mathcal F(\kappa,\Xi)=(M_\Xi-M,[r_2]_1,\ldots,[r_{24}]_1).
\]
For formal initial germs agreeing in radius through degree $L-1$
and velocity through degree $L-3$, the operators
$I,\partial_\xi/\xi,\partial_\xi^2$ in
\eqref{loc:vacuum-force} and the recursion
\eqref{loc:initial-recursion} give
\[
 \operatorname{val}r_k\ge2-2k,\qquad
 \operatorname{val}\Delta r_k\ge L-2k,\qquad1\le k\le25.
\]
Indeed a partition term of total time index $k-1$ sums the losses
$-2i$; changing one direction or a base coefficient improves the
bound by $L-2$. Expanding product differences proves the second
inequality. With $L=50$, the rows through $k=24$ therefore depend
only on radius coefficients through degree 49, velocity coefficients
through degree 47, and $M_\Xi$. Expanding reciprocals about positive
constant factors makes $\mathcal F$ a $C^1$ finite-dimensional map.

We compare its Jacobian at the Euler-compatible zero with that of
Lemma~\ref{loc:finite-map}. In $y=\lambda z$, the directions
\eqref{loc:columns} are nonzero fixed-scale multiples of
$y^{2n+3}$ and $y^{2n+3}/(2n+3)$. Their cutoff differences have
zero retained Taylor coefficients. The interior mass column removes
the enthalpy mass derivative; velocity differences affect neither
mass nor these coefficients. Thus the column change is triangular
and invertible on the retained coefficient and mass space.

For the rows, put $f_j=[R_j(z)]_1$, $g_j=[R_j(\xi)]_1$ in
physical time. Since $R_j=\lambda^{2-2j}u_j$, the old row is a
nonzero multiple of $f_j$. On the common kernel of the mass row and
the preceding $Df_i$, the earlier linearized derivatives are ordinary
and compatible. The substitution $\xi=c_{\rm ph}z+O(z^2)$,
$c_{\rm ph}>0$, gives $Dg_j=c_{\rm ph}^{-1}Df_j$ there.
Its coordinate-variation term has no first coefficient: that variation
vanishes at zero and the background derivative has zero first trace.
The preceding rows are independent, so
\[
 Dg_j=c_{\rm ph}^{-1}Df_j+\sum_{2\le i<j}a_{ji}Df_i
                                     +b_jD(M_\Xi-M).
\]
This is an invertible lower triangular row change, proving the required
Jacobian invertibility.

\emph{Viscosity curve.}
Recenter the Euler zero at $\Xi=0$ and put
$J_0=\partial_\Xi\mathcal F(0,0)$,
$\mathcal T_\kappa(\Xi)=\Xi-J_0^{-1}\mathcal F(\kappa,\Xi)$.
On a fixed small ball, for sufficiently small $\kappa$,
\[
 \|D_\Xi\mathcal T_\kappa\|\le\tfrac12,\qquad
 |\mathcal T_\kappa(0)|\le C_0\kappa,\qquad
 |\mathcal T_\kappa(\Xi)|\le C_0\kappa+\tfrac12|\Xi|.
\]
Hence the ball $|\Xi|\le2C_0\kappa$ is invariant and contains
a unique fixed point. Differentiating its equation gives
\[
 \Xi'(\kappa)=-
 \{\partial_\Xi\mathcal F(\kappa,\Xi(\kappa))\}^{-1}
                      \partial_\kappa\mathcal F(\kappa,\Xi(\kappa)),
\]
which is continuous.

Along this zero set, generate each time derivative before taking the
next quotient. The two-derivative force loss and vanishing constraints
give
\[
 r_j^\kappa\in C^1_\kappa C^{51-2j}\quad(0\le j\le25),
 \qquad \partial_\xi r_j^\kappa(0)=0\quad(0\le j\le24).
\]
Parameter differentiation preserves these traces; \eqref{loc:quotients}
then gives the domain quotients. In particular $r_{24}$ is $C^3$
and $r_{25}$ is $C^1$. The embedding
$C^{k+1}\hookrightarrow h^{k,\alpha}$, with little modulus
$\delta^{1-\alpha}\|D^{k+1}f\|_\infty\to0$, proves
\eqref{loc:initial-cone}. Integrating the continuous parameter
derivative gives the fixed-scale $O(\kappa)$ bound.
\end{proof}

\subsection{The modified energy of the initial data}
\label{loc:initial-energy}

We estimate the initial energy in absolute value; its coercivity will
be proved later. We first compute the derivatives of the potential
that enter the modified energy. Define
\begin{equation}
 \mathcal{U}[R]=\int_0^M\left\{\mathfrak e(\rho_R)-\frac{x}{R}\right\}\dd x,
 \qquad \mathfrak e(0)=0,\quad \mathfrak e'(\rho)=P(\rho)/\rho^2.
 \label{loc:potential}
\end{equation}
For a direction $f$, write $a_f=f/R$, $c_f=f_x/R_x$.
The exact identity
\begin{equation}
 \rho_{R+\sum t_if_i}
   =\rho_R(1+\sum t_i a_{f_i})^{-2}(1+\sum t_i c_{f_i})^{-1}
 \label{loc:strain-identity}
\end{equation}
contains only first spatial strains. If
$u=\rho^{2/3}/(1+\rho^{2/3})$, then
\[
 (\rho\partial_\rho)\log P'=\tfrac23-\tfrac13u,\qquad
 (\rho\partial_\rho)u=\tfrac23u(1-u).
\]
Induction, followed by
$P(\rho)=\rho\int_0^1P'(\rho t)\dd t$, gives for every fixed
finite $r\ge1$
\begin{equation}
 |(\rho\partial_\rho)^r\mathfrak e|\le C_rP/\rho,\qquad
 |D^r\mathfrak e(\rho_R)[f_1,\ldots,f_r]|
   \le C_r(P/\rho)\prod_i(|a_{f_i}|+|c_{f_i}|).
 \label{loc:potential-majorant}
\end{equation}
At the vacuum $P/\rho=O(\xi^2)$ and the domain strains are
bounded; the measure is $5\xi^4\dd\xi$. At the center
$f_i=O(z)$, $R\sim z$, and the gravitational integrand
is bounded by $C_rz^2\prod_i\|f_i\|_{\mathcal D}$.
These bounds also control the next-order remainders on a smaller
neighborhood satisfying the positivity conditions. Integrating the
Taylor remainder proves that $\mathcal{U}$ is $C^\infty$ there.

Integration on $[\epsilon_c^3,M-\epsilon_v^5]$ gives
\[
 D\mathcal{U}[R]f
 =-\int4\pi P(\rho_R)\partial_x(R^2f)\dd x
       +\int xR^{-2}f\dd x.
\]
The boundary term $-[4\pi R^2P f]$ is
$O(\epsilon_v^5)+O(\epsilon_c^3)$. Therefore
\begin{equation}
 D\mathcal{U}[R]f=(\mathcal M[R],f)_H,\qquad
 D^2\mathcal{U}[R][f,g]=\mathfrak q_R[f,g],
 \quad H=L^2((0,M),\dd x).
 \label{loc:potential-force}
\end{equation}
The exact Green form is
\begin{equation}
 \mathfrak q_R[f,g]
 =\int P'(\rho_R)(2f/R+f_x/R_x)(2g/R+g_x/R_x)\dd x
  +\int(8\pi RP_x-2x/R^3)fg\dd x.
 \label{loc:green-form}
\end{equation}
Indeed the truncated boundary flux is
$-[4\pi R^2\rho P'(2f/R+f_x/R_x)g]$, with the same endpoint
powers. Even at base center regularity, $\rho\in C^{1,\alpha}$
is invariant and $\rho_z=O(z^\alpha)$; hence the apparently
singular $RP_x$ term times $fg\dd x$ is
$O(z^{3+\alpha})\dd z$, integrable. A pointwise bound for
$RP_x$ at the center is not required.

On a segment of prepared radii satisfying the normalized positivity
bounds, we estimate two directions in $X^1$ and the others in
$\mathcal C^2$. On the edge charts $n=5,8$,
\[
 \rho=\frac{ny^{n-2}}{4\pi r^2\gamma},\quad
 P/\rho\le Cy^2,\quad
 |a_f|+|c_f|
 \le C\lambda^{-2}y^{-1}(|f_y|+\lambda y|f|).
\]
The density bound uses $P/\rho\le C\rho^{2/3}$ on the bounded
chart and $P/\rho\le C\rho^{1/3}$ on the entire long chart.
For the two directions measured in $X^1$, the two factors
$y^{-1}$ are canceled by the factor $y^2$ in $P/\rho$.
Each remaining direction is bounded by
$C\lambda^{-2}|f|_{\mathcal C^2}$, using its zero first
derivative at $y=0$. The center calculation uses the full
Cartesian first derivative, which includes $f/z$. Gravity
has no worse bound, since its two center denominators cancel
the $z^2$ in $\lambda^3x/R$.
To localize, write $S=\sum\chi_\mu^2\ge1/4$ and partition the
integral by $\chi_\mu^2/S$. Both high factors receive their
own cutoff; $\chi f_y=(\chi f)_y-\chi_yf$, with the error
controlled by $H$. Cauchy--Schwarz with the exact measure
$n\lambda^4y^{n-1}\dd y$ proves
\begin{equation}
 |D^r(\lambda^3\mathcal{U})[R][f_1,\ldots,f_r]|
 \le C_r\lambda^{3-2r}\|f_a\|_{X^1}\|f_b\|_{X^1}
              \prod_{i\ne a,b}|f_i|_{\mathcal C^2},
 \qquad r\ge2\ \hbox{fixed}.
 \label{loc:two-high-potential}
\end{equation}
In particular orders through 26 are allowed; this changes neither
the number of time derivatives nor the endpoint traces.

Use ordinary derivatives in rescaled time $R_i,A_i$, and put
$Y_i=R_i-A_i=\lambda^4\psi_i$. Define
\[
 Q_R=\lambda^3D^2\mathcal{U}[R],\qquad
 Q_R^{\langle q\rangle}=\lambda^3\partial_s^qD^2\mathcal{U}[R].
\]
The outer $\lambda^3$ is held fixed in the second definition.
The finite partition identity is
\begin{equation}
 Q_R^{\langle q\rangle}[f,g]
 =\sum_{\pi\in\Pi_q}D^{|\pi|+2}(\lambda^3\mathcal{U})[R]
                         [R_{|B|}:B\in\pi,f,g].
 \label{loc:form-bell}
\end{equation}
The commutator corrections in the energy are
\begin{equation}
 \begin{aligned}
 \mathfrak A_j=\lambda^{-8}\Biggl\{&
 (Q_R-Q_A)[A_j,Y_j]\\
 &+\sum_{q=1}^{j-1}\binom{j-1}{q}
  \bigl(Q_R^{\langle q\rangle}[R_{j-q},Y_j]
                  -Q_A^{\langle q\rangle}[A_{j-q},Y_j]\bigr)
                         \Biggr\},\quad 1\le j\le22.
 \end{aligned}
 \label{loc:primitive}
\end{equation}
For these initial time derivatives, set
\[
 G_8^2=\sum_{q=0}^{8}\sum_{j=0}^{23-q}\|\psi_j\|_{X^q}^2,
 \qquad
 \mathcal X_8=\sum_{q=0}^{8}\sum_{j=0}^{23-q}\|\psi_j\|_{X^q}.
\]
There are 180 summands, so
$G_8\le\mathcal X_8\le\sqrt{180}\,G_8$.
The potential at order zero is the nonlinear relative potential
\[
 \mathfrak V_0
 =\lambda^{-5}\{\mathcal{U}[R]-\mathcal{U}[A]-D\mathcal{U}[A](R-A)\}
 =\lambda^3\int_0^1(1-t)\mathfrak q_{A+tY_0}[\psi_0,\psi_0]\dd t.
\]
For fixed bounded weights $\omega_j$ and
$c_{\rm en}=13+\Lambda_{\rm en}$, set
\begin{align}
 \widetilde H_0&=\tfrac12\|\psi_1\|_H^2+\mathfrak V_0
       +2b(\psi_0,\psi_1)_H+\tfrac12c_{\rm en}b^2\|\psi_0\|_H^2,
       \notag\\
 H_j&=\tfrac12\|\psi_{j+1}\|_H^2+\tfrac12Q_R[\psi_j,\psi_j]
       +2b(\psi_j,\psi_{j+1})_H+\tfrac12c_{\rm en}b^2\|\psi_j\|_H^2,
       \notag\\
 \widetilde E&=\omega_0\widetilde H_0
                +\sum_{j=1}^{22}\omega_j(H_j+\mathfrak A_j).
 \label{loc:modified-energy}
\end{align}

\begin{proposition}
\label{loc:prepared-energy}\label{en:input-prepared-energy}
There is a constant $C$, independent of the terminal scale, with the
following properties.
\begin{enumerate}[label=(\roman*),leftmargin=2em]
\item The Euler-compatible data satisfy
\begin{equation}
 |\widetilde E^0|\le C\lambda^{2J+1644}(1+|\log\lambda|)^{2p_{\log}},
 \qquad
 \mathcal X_8^0\le C\lambda^{J+1685/2}(1+|\log\lambda|)^{p_{\log}}.
 \label{loc:strong-preparation-margin}
\end{equation}
\item There is
$0<\kappa_0(\lambda,\Lambda)\le\kappa_1(\lambda,\Lambda)$
such that, for $0\le\kappa\le\kappa_0(\lambda,\Lambda)$,
\begin{equation}
 \begin{gathered}
 |\widetilde E^\kappa|
       \le C\lambda^{2J+1644}(1+|\log\lambda|)^{2p_{\log}},\\
 \mathcal X_8^\kappa\le C\lambda^{J+1685/2}(1+|\log\lambda|)^{p_{\log}},\\
 |\widetilde E^\kappa|+(G_8^\kappa)^2\le C\lambda^{2J+1},\\
 \kappa_0\le c\lambda^{J+5},\qquad
 \kappa_0\le\frac{c\lambda^{J+2}}{1+\Lambda}.
 \end{gathered}
 \label{loc:full-energy-margin}
\end{equation}
The same constant $C$ applies throughout this viscosity interval.
\end{enumerate}
The endpoint $\kappa_0(\lambda,\Lambda)$ may depend on the terminal scale.
\end{proposition}

\begin{proof}
\emph{Quadratic terms.}
The segment of prepared radii is positive by \eqref{loc:mass-variation}
and equals $A$ near the center. On the bounded collar,
\[
 \delta r=O(\lambda^5y^3D_{\rm par}),\quad
 \delta\gamma=O(\lambda^4yD_{\rm par}),\quad
 \delta(y\gamma_y)=O(\lambda^4yD_{\rm par}),
\]
with smaller long-edge bounds. In the fixed interior, the normalized radius and first derivative have
uniform positive upper and lower bounds, and the second derivative is
uniformly bounded. Proposition~\ref{pr:current-hessian} applies on the
whole segment. Its absolute Hessian bound and
\eqref{loc:prepared-graph} give
\[
 |\widetilde H_0|+\sum_{j=1}^{22}|H_j|
       \le C\lambda^{-90}D_\lambda^2.
\]
Here the bound for $\mathfrak V_0$ follows by integration along the
segment. Converting \eqref{loc:native-corrected-rows} to rescaled time
also gives
\[
 |R_i|_{\mathcal C^2}+|A_i|_{\mathcal C^2}
       \le C_i\lambda^{2-2i},\qquad1\le i\le22.
\]

\emph{Commutator corrections.}
For a partition in \eqref{loc:form-bell} with $p$ blocks, interpolate
the base point and directions from $A$ to $R$. In a base-point
variation, place $Y_0,Y_j$ in the form norms; the other $p+1$
directions have total time index $j$. Including $\lambda^{-8}$,
\[
 \lambda^{-8}\lambda^{3-2(p+3)}\lambda^8
                  \lambda^{2(p+1)-2j}=\lambda^{-1-2j}.
\]
An argument variation $Y_i$ costs $\lambda^{-1-2j+2i}$.
The first term in \eqref{loc:primitive} has the base-point bound
with $p=0$. Consequently
\[
 \sum_{j=1}^{22}|\mathfrak A_j|
 \le C\lambda^{-45}
       \left(\sum_{i=0}^{22}\|\psi_i\|_{X^1}\right)^2
 \le C\lambda^{-131}D_\lambda^2.
\]
This uses potential order at most 24 and coefficient time order at most
21; the separate first term uses the twenty-second time derivative.
Since $D_\lambda\le C\lambda^{J+1775/2}(1+|\log\lambda|)^{p_{\log}}$,
\[
 |\widetilde E^0|
 \le C\lambda^{2J+1644}(1+|\log\lambda|)^{2p_{\log}}
 \le C\lambda^{2J+1}.
\]
The squared norms in \eqref{loc:prepared-graph} have a larger
positive scale power, proving \eqref{loc:strong-preparation-margin}.

\emph{Viscous family.}
At fixed scale, Proposition~\ref{loc:selector} gives
\[
 |\widetilde E^\kappa-\widetilde E^0|\le C_E(\lambda)\kappa,
 \qquad
 |\mathcal X_8^\kappa-\mathcal X_8^0|\le C_G(\lambda)\kappa.
\]
The parameter derivative uses potential order 25, allowed by
\eqref{loc:two-high-potential}. Time derivatives through order 22
are $C^7$ in space, and the order-23 kinetic derivative is $C^5$.
These regularities give the displayed Lipschitz bounds without taking
a square root of the energy. Restrict the family by
\[
 \kappa\le
 \frac{\lambda^{2J+1644}(1+|\log\lambda|)^{2p_{\log}}}{1+C_E(\lambda)},
 \qquad
 \kappa\le
 \frac{\lambda^{J+1685/2}(1+|\log\lambda|)^{p_{\log}}}{1+C_G(\lambda)},
\]
and by the two viscosity bounds in \eqref{loc:full-energy-margin}.
Their minimum with $\kappa_1(\lambda,\Lambda)$ is positive at every
fixed scale. The resulting estimates have one enlarged
scale-independent constant. Decreasing $A_0$ absorbs the logarithms
in the excess scale powers and proves \eqref{loc:full-energy-margin}.
\end{proof}

We now solve the regularized equation with these data. The local
lifespan depends on the radius and velocity norms in $\mathcal D$;
the higher initial time derivatives do not enter its choice.

\enlargethispage{3pt}
\subsection{Local solutions}
\label{loc:local-solutions}

In the local evolution below we use $t=\tau-\tau_*\ge0$, where
$\tau_*$ is the prepared time and $\tau$ is reversed physical time.
This $t$ is distinct from the normalized time in the finite compatibility
calculation.

Let $\mathcal P_T,\mathcal F_T$ denote the little parabolic
domain and source spaces with exponents $2+\alpha,1+\theta$
and $\alpha,\theta$ on the fixed charts. Their vacuum component
uses the radial lift $f(\xi)\mapsto f(|X|)$ on a five-dimensional
ball, as in Lemma~\ref{loc:base-lift}. In particular,
\[
 \mathcal P_T\subset C([0,T];\mathcal D),\qquad
 \partial_t\mathcal P_T\subset\mathcal F_T
                            \subset C([0,T];\mathcal Y).
\]

On a Cartesian chart, the source norm is
\[
 \|F\|_\infty+[F]_{X;\alpha}+[F]_{t;\theta}.
\]
The domain norm is the sum of the spatial $C^{2,\alpha}$ norm,
the source norm of $u_t$, $[D^2u]_{t;\theta}$, and
$[Du]_{t;(1+\alpha)/2}$. The temporal Hölder seminorm of the
spatial Hessian is taken with values in $C^0$.

The little spaces are the closures of smooth functions in these
norms, with the indicated symmetry and, for a local domain, the
homogeneous artificial normal trace. The initial value lies in
$\mathcal D$, while $u_t(0)$ lies
in $\mathcal Y$. At the vacuum, the five-dimensional lift is
used at this base regularity; it imposes no higher even trace
conditions on the one-sided radius.

Fix $(r_0,r_1)\in\mathcal D\times\mathcal D$, with
$\|r_0\|_{\mathcal D}+\|r_1\|_{\mathcal D}\le L$ and all the
stated positive factors of $r_0$ bounded below by $\eta$. Write $V=R_t$ and $R_j=\partial_t^jR$. With initial data
$R(0)=r_0$ and $V(0)=r_1$, the regularized equation and the
required base regularity are
\begin{equation}
 R_t=V,\qquad V_t+\kappa(H_R+\Lambda)V=-\mathcal M[R],
 \quad V\in\mathcal P_T,\quad R=r_0+\int_0^tV\in C_t^1\mathcal D.
 \label{loc:actual-velocity}
\end{equation}
\equationalias{in:viscous-equation}{loc:actual-velocity}

\begin{proposition}
\label{loc:high-local}
Fix $\kappa>0$, and let $(r_0,r_1)$ satisfy the norm and positivity
bounds specified above. Generate $r_2,\ldots,r_{25}$ by
\eqref{loc:initial-recursion}, and assume
\[
 r_j\in\mathcal D^{48-2j}\quad(0\le j\le24),\qquad
 r_{25}\in\mathcal Y.
\]
There is a time
\[
 T=T(L,\eta,\kappa,\Lambda,M,\alpha,\mathrm{atlas})>0
\]
with the following properties.
\begin{enumerate}[label=(\roman*),leftmargin=2em]
\item Equation \eqref{loc:actual-velocity} has a unique solution with
$V\in\mathcal P_T$ and $R\in C_t^1\mathcal D$.

\item Its time derivatives satisfy $R_j(0)=r_j$ for $0\le j\le25$ and
\begin{equation}
 R_j\in\mathcal P_T\cap C_t\mathcal D^{48-2j}
       \quad(0\le j\le24),\qquad R_{25}\in\mathcal F_T.
 \label{loc:actual-finite-cone}
\end{equation}
No first-derivative trace at the vacuum is imposed on $R_{25}$.

\item For fixed parameters, convergence of the compatible initial time
derivatives in
\[
 \prod_{j=0}^{24}\mathcal D^{48-2j}\times\mathcal Y,
\]
with common bounds $L,\eta$, implies convergence of the solutions in
all the displayed spaces on a common interval.

\item A finite maximal time can occur only if
$\|R(t)\|_{\mathcal D}+\|R_t(t)\|_{\mathcal D}$ becomes unbounded
or one of the stated positive factors has no uniform positive lower bound.
\end{enumerate}
In particular, every fixed positive-viscosity member of
\eqref{loc:initial-cone} satisfies these hypotheses.
\end{proposition}

Appendix~\ref{app:local-theory} proves the proposition in three steps.
First, the frozen inverse and the factor $T^{1-\theta}$ in
\eqref{loc:integrated-gain} give a contraction for the velocity.
The restrictions \eqref{loc:base-ball-time} depend only on the base
norms and fixed parameters.

Second, $H_R=D\mathcal M[R]$ gives the coupled system
\eqref{loc:radius-history}. Its differentiated equations are linear in
the new time derivative, with sources involving only earlier derivatives.
Third, \eqref{loc:elliptic-implications} recovers the spatial regularity
in the order \eqref{loc:triangle-order}. Both inductions use the same
interval. Higher initial norms enter the estimates, but not the lifespan;
this permits continuation from bounds on the base norms and the positive
factors. The viscosity remains fixed and positive throughout the local
argument.

\subsection{Justification of the identities at finite regularity}
\label{loc:finite-legality-section}

Fix a positive scale, the covering cutoffs, and $2\le q\le23$.
Write $\mathscr C$ for the compatible smooth core of
Definition~\ref{at:atlas}, and $L_R=\lambda^3H_R$. At order $q$,
we take the closure of
\[
 \{(f,L_Rf):f\in\mathscr C\}
\]
with respect to the product norm
\[
 (f,g)\longmapsto\|f\|_{X^q}+\mathcal S_{q-2}(g).
\]
The second component lies in the source completion defined by
\eqref{at:source-norm}; its edge term is $\lambda L_Rf$, with
no first-derivative trace condition on the source.

\begin{lemma}
\label{loc:finite-strong-graph}
At a fixed positive scale, the following approximation statements hold
for $2\le q\le23$.

\begin{enumerate}[label=(\roman*),leftmargin=2em]
\item \emph{Finite local class.} If $R\in\mathcal D^{q-2}$ is positive
and $f\in\mathcal D^{q-2}$, there are $f_n\in\mathscr C$ such
that
\begin{equation}
 \|f_n-f\|_{X^q}
       +\mathcal S_{q-2}(L_Rf_n-L_Rf)\longrightarrow0.
 \label{loc:finite-graph-approximation}
\end{equation}
Simultaneous approximation of a finite number of radius and direction
fields in $\mathcal D^{q-2}$ also gives convergence of every fixed
finite differential of $\mathcal M$ in $\mathcal S_{q-2}$.

\item \emph{Compatible weighted class.} Suppose instead that
$R=A+\lambda^4Z$, the low-norm and positivity hypotheses of
Proposition~\ref{at:current-perturbation} hold, and
$Z\in X^{\max(q,8)}$. Then every $f$ in the compatible
completion $X^q$ satisfies \eqref{loc:finite-graph-approximation}.
For $q\ge9$, the only high-order norms needed are
$\|Z\|_{X^q}$ and $\|f\|_{X^q}$.
\end{enumerate}
The constants in these approximation bounds may depend on the
positive scale.
\end{lemma}

The proof is given in Subsection~\ref{app:current-domain-approximation}.

\begin{proposition}
\label{loc:finite-calculus}\label{en:input-finite-calculus}
Fix the positive viscosity and a compact time interval of the solution
in Proposition~\ref{loc:high-local} on which the scale remains positive.
Then:
\begin{enumerate}[label=(\roman*),leftmargin=2em]
\item The mixed derivative equations hold in the stated source space
for every pair $q\ge2$, $j+q\le23$.
\item The force identities, Green's formula, and the identities for the
rescaled Hessian and modified energy hold through time order 22.
The kinetic energy differentiation uses time order 24. The highest
form argument has time order 23; the derivative of order 25 belongs
only to the source space and has no prescribed derivative trace.
\item These identities persist under convergence in the full range of
spaces in \eqref{loc:actual-finite-cone}.
\end{enumerate}
\end{proposition}

The proof is given in Subsection~\ref{app:finite-regularity-identities}.

\begin{proposition}
\label{loc:prepared-local-family}\label{exlim:upstream-preparation}
Fix the matched profile above and choose the upper bound for the scale sufficiently
small. For every $0<\lambda_*\le A_0$ and finite $\Lambda\ge0$,
the following assertions hold.
\begin{enumerate}[label=(\roman*),leftmargin=2em]
\item The family of compatible data with exact mass \eqref{loc:initial-cone} exists. Its
zero-viscosity correction satisfies \eqref{loc:prepared-graph}.
There is $\kappa_0(\lambda_*,\Lambda)>0$ such that the whole modified
initial energy and $G_8$ satisfy \eqref{loc:full-energy-margin} for
$0\le\kappa\le\kappa_0(\lambda_*,\Lambda)$.
\item For each positive $\kappa$ in this interval, the finite local
solution satisfies \eqref{loc:actual-finite-cone}, the identities of
Proposition~\ref{loc:finite-calculus}, and the forward continuation
criterion of Proposition~\ref{loc:high-local}.
\end{enumerate}
The correction, source estimates, and evolution have the same mass,
profile, and covering cutoffs. The lifespan may depend on $\lambda_*$
and $\kappa$.
\end{proposition}

\begin{proof}
Lemma~\ref{loc:finite-map} and
\eqref{loc:mass-variation}--\eqref{loc:prepared-graph} construct the
Euler-compatible data and bound their correction at fixed mass.
Proposition~\ref{loc:selector} gives the viscous family, and
Proposition~\ref{loc:prepared-energy} proves its initial energy and
mixed Sobolev bounds. For each positive viscosity,
Propositions~\ref{loc:high-local} and \ref{loc:finite-calculus} give the
local solution, its continuation criterion, and the asserted identities.
The lifespan depends on the base-space bounds; the initial energy and
mixed Sobolev estimates have the scale-independent constants stated above.
\end{proof}

\section{Recovery of the mixed spatial derivatives}
\label{rec:chapter}

The energy controls time derivatives and one spatial derivative.
We recover the remaining spatial derivatives from the equation, using
Proposition~\ref{at:current-inverse} for the linearization at $R$.
The order of the argument is dictated by two terms: at spatial order
$q$, the inertial term uses spatial order $q-2$ and two more time
derivatives, while a profile commutator uses the same spatial order
and fewer time derivatives.

We therefore induct first on the spatial order and then on the time
order. Terms containing the highest radius norm carry a small
low-order strain and are absorbed after the induction. The
first-order quantity $B$ remains on the right-hand side until
Corollary~\ref{en:base-bound} bounds it by the energy.

We use the weighted spaces of Definition~\ref{at:atlas}, the matched profile
of Proposition~\ref{exlim:upstream-profile}, and the compatible
solution class of Proposition~\ref{en:input-finite-calculus}.
Proposition~\ref{rec:recovery} gives the full inviscid estimate and
the lower spatial estimate uniform in viscosity used in
Section~\ref{ct:section}. Both apply to compatible strong
representatives and extend to the closed strong operator domains.
For the inviscid limit, Subsection~\ref{exlim:realization} first
uses the limiting equation to show that the derivatives lie in these
domains. We then apply the uniform recovery estimate.

\subsection{The equation and the quantities to be recovered}

In the rescaled time $s$, the scale satisfies
\begin{equation}
 \lambda_s=b\lambda,\qquad b^2=\beta^2+e\lambda,\qquad
 0\le e\lambda\le\beta^2,\qquad
 0<\beta\le b\le\sqrt2\,\beta\le b_{\max},
 \label{rec:clock}
\end{equation}
We restrict the upper bound for the scale to at most one. Here $b_{\max}$ is a
fixed upper bound; it provides no positive lower bound independent
of the selected $\beta$. Write
\begin{equation}
 Y=R-A,\qquad \Psi_j=\lambda^{-4}\partial_s^jY,\qquad
 L_R=\lambda^3D\mathcal M[R],\qquad
 \nu=\kappa\lambda^{-3/2}.
 \label{rec:rows}
\end{equation}
The fixed nonnegative shift is denoted by $\Lambda$. The displacement
equation is exactly
\begin{equation}
 \begin{split}
 Y_{ss}-\frac32bY_s
 &+\lambda^3\{\mathcal M[R]-\mathcal M[A]\}\\
 &+\nu\lambda^3\partial_s\{\mathcal M[R]-\mathcal M[A]\}
       +\nu\lambda^3\Lambda Y_s=\lambda^4F .
 \end{split}
 \label{rec:equation}
\end{equation}
The prescribed forcing $F$ is given in
Proposition~\ref{en:input-material-source}. It contains the matched
residual and the profile terms arising from viscosity and the shift;
the perturbative force remains on the left. The finite correction
changes only the initial data, leaving $A$ and the evolving source
unchanged. We may also include an additional prescribed forcing
whose strong norms below are finite.

For $Q_*\in\{8,23\}$ set
\begin{align}
 \mathcal X_{Q_*}
   &=\sum_{q=0}^{Q_*}\sum_{j=0}^{23-q}\|\Psi_j\|_{X^q},
 &G_8^2&=\sum_{q=0}^{8}\sum_{j=0}^{23-q}\|\Psi_j\|_{X^q}^2,
 \label{rec:cones}\\
 B&=\sum_{j=0}^{23}\|\Psi_j\|_{X^0}
                  +\sum_{j=0}^{22}\|\Psi_j\|_{X^1},
 &\eta&=\lambda^{-2}G_8.
 \label{rec:base}
\end{align}
Thus $B$ measures one spatial derivative and is defined independently
of the energy estimate. The elementary comparison
$G_8\le\mathcal X_8\le\sqrt{180}\,G_8$ will allow either low
norm to be used in the bootstrap. The strong forcing quantities are
\begin{equation}
 F_j=\lambda^{-4}\partial_s^j(\lambda^4F),\qquad
 \mathcal R_{Q_*}
   =\sum_{q=2}^{Q_*}\sum_{j=0}^{23-q}
                         \mathcal S_{q-2}(F_j).
 \label{rec:source}
\end{equation}
Here the source is measured in the strong norm
$\mathcal S_{q-2}$, with the edge factor $\lambda$ from
\eqref{at:source-norm}, rather than in a dual form norm.

Throughout the proof the current radius and the segments
$A+tY$, $0\le t\le1$, have the positive lower bounds of the
current inverse. We impose $\eta\le\eta_*$, with $0<\eta_*\le1$ to be
chosen after the finite profile and inverse constants. Since
\[
 \lambda^{-1}\|\Psi_0\|_{X^8}\le\lambda\eta_*,
\]
decreasing the upper bound for the scale, if necessary, enforces the smallness
required by Proposition~\ref{at:current-inverse}.

\subsection{Material derivatives of the profile coefficients}

We first bound the time derivatives of the profile operator in terms
of the prescribed radius and Jacobian factors. These bounds are
independent of the high norm of the unknown solution. Time derivatives
in this subsection are ordinary derivatives at fixed mass; the
rescaled Hessian derivatives in Section~\ref{en:section} hold the
exterior factor fixed.

\begin{lemma}
\label{rec:profile-columns}
For $2\le q\le23$, $\ell\ge0$, and $\ell+q\le24$,
\begin{equation}
 \mathcal S_{q-2}\bigl(L_A^{(\ell)}u\bigr)
             \le C_{q,\ell}b^\ell\|u\|_{X^q},
 \qquad L_A^{(\ell)}=\partial_s^\ell L_A.
 \label{rec:profile-column-bound}
\end{equation}
The constants depend on the selected finite profile bounds, the
fixed cover, and the prescribed positive lower bounds. They are
independent of the terminal scale and the length of the long chart.
\end{lemma}

\begin{proof}
In the edge chart with $n=5$ or $8$, put
\begin{equation}
 c_n=\frac4n,\qquad
 \mathscr E=\lambda\partial_\lambda-c_ny\partial_y,\qquad
 r_\ell=\frac{\partial_s^\ell A}{\lambda b^\ell},\qquad
 \gamma_\ell
   =-\frac{\partial_y\partial_s^\ell A}{\lambda^2yb^\ell}.
 \label{rec:profile-units}
\end{equation}
Let $\vartheta=e\lambda/b^2$. Direct differentiation, using
$b_s/b^2=\vartheta/2$, gives
\begin{align}
 r_{\ell+1}
    &=\left(\mathscr E+1+\frac{\ell\vartheta}{2}\right)r_\ell,
 \label{rec:radius-recurrence}\\
 \gamma_{\ell+1}
    &=\left(\mathscr E+2-2c_n+
                          \frac{\ell\vartheta}{2}\right)\gamma_\ell,
 &\mathscr E\vartheta&=\vartheta(1-\vartheta).
 \label{rec:jacobian-recurrence}
\end{align}
Indeed $\partial_y\mathscr E=(\mathscr E-c_n)\partial_y$.
Consequently the ordinary derivatives through order $m$ of
$r_\ell,\gamma_\ell,y\partial_y\gamma_\ell$ involve mixed derivatives of the radius of total order at most
$\ell+m+2$.

The profile theorem supplies these bounds on the prescribed charts.
On the bounded vacuum chart, the material derivative $b\mathscr E$
maps each coefficient $\lambda^a(\log\lambda)^h$ to a sum with
the same power of $\lambda$. On the long boundary chart, a correction
of matching degree $r\ge2$ remains bounded by
\[
 C(\lambda y_8^2)^{r-1}
              (1+|\log\lambda|+\log y_8)^{p_{\log}}.
\]
Here $r$ is the matching degree, not the chart dimension $n$.
On the outer core, the bounds for iterates of $q\partial_q$
transfer to the required iterates of $y\partial_y$ on the
normalized Jacobian factors. Ordinary derivatives follow from
\[
 \partial_y^m=y^{-m}\prod_{r=0}^{m-1}(y\partial_y-r).
\]
On the matching shell $q\sim\sqrt\lambda$, or equivalently
$y\sim\lambda^{-1/4}$, the cutoff is a fixed smooth function
of $q/\sqrt\lambda$. Its required derivatives under
$\mathscr E$ and $y\partial_y$ remain bounded. Applying these
operators therefore preserves the small power in the connection
remainder and uses only the finite coefficients already constructed.

For an index pair in the recovery argument $m=q-2$, the largest profile time index
needed by the viscous equation is $j+1$. Hence
\begin{equation}
 \ell+m+2\le j+q+1\le24.
 \label{rec:profile-order-count}
\end{equation}
The same count applies to the Cartesian derivatives of the radial vector
at the center and the ordinary derivatives of the radius in the fixed interior.

Now use the exact coefficient identity
\[
 \lambda L_A=-A_nD_n-B_n\frac{\partial_y}{y}+V_n,
 \qquad D_n=\partial_y^2+\frac{n-1}{y}\partial_y,
\]
from \eqref{at:edge-hessian}. For a function held fixed in mass,
\begin{equation}
 \partial_s(\partial_y^ku)
      =c_nkb\,\partial_y^ku .
 \label{rec:physical-derivative}
\end{equation}
Thus an ordinary time derivative of the operator differentiates
its coefficients, its exterior $\lambda^{-1}$, and these
coordinate factors. Equations
\eqref{rec:radius-recurrence}--\eqref{rec:jacobian-recurrence}
bound every resulting coefficient through order $q-2$ by
$C_{q,\ell}b^\ell$. Known profile coefficients are used as
$W^{q-2,\infty}$ multipliers. Lemma~\ref{at:same-cutoff-product}
retains the prescribed cutoff, and
Lemma~\ref{at:natural-quotient} extends the quotient
$u_y/y$ to the weighted completion at the vacuum. At the center use the full Cartesian operator, including its angular
part; in the interior use \eqref{at:interior-hessian}. This proves
\eqref{rec:profile-column-bound}. Each time derivative of the profile
operator still has spatial order two, so its application requires
the $X^q$ norm of $u$.
\end{proof}

\subsection{Nonlinear products at the current radius}

In the differentiated force, the highest spatial derivative may
fall on a direction or on the current coefficient. The following
count ensures that the remaining perturbation factors have low norms.
It will also be used in the form estimates, with their different
spatial hypotheses. Suppose the perturbation factors have time
indices $i_1,\ldots,i_p\ge0$, $p\ge2$, with sum at most $T\le23$,
and arrange them in decreasing order. Then
\begin{equation}
 i_{(2)}\le\lfloor T/2\rfloor\le11,
 \qquad
 i_{(1)}\ge16\ \Longrightarrow\
                 \sum_{a\ge2}i_{(a)}\le7.
 \label{rec:single-high-rule}
\end{equation}
Profile indices also count toward $T$, so omitting them only weakens
these inequalities. Taylor expansion of a base-point difference adds
a perturbation factor of time index zero. Thus, if one perturbation factor has a time index too high for its
$X^8$ norm to occur in the controlled range, every other factor has
its $X^8$ norm. The spatial product estimate must also include the
term in which the highest derivative falls on the current coefficient.

\Needspace{13\baselineskip}
\begin{lemma}
\label{rec:current-bell}
Consider a term in one of the finite differentiated equations
below, with perturbation directions $u_1,\ldots,u_p$ and profile
directions $A_{\ell_1},\ldots,A_{\ell_k}$, where
$A_\ell=\partial_s^\ell A$ and $\ell_i\ge1$. Assume $q\ge9$.
\begin{enumerate}[label=(\roman*),leftmargin=2em]
\item \emph{At least two perturbation directions.} For $p\ge2$, its strong coefficient estimate is
\begin{equation}
 \begin{split}
 &\mathcal S_{q-2}\!\left(
  \lambda^{4p-1}D^{p+k}\mathcal M[R]
   [u_1,\ldots,u_p,A_{\ell_1},\ldots,A_{\ell_k}]\right)\\
 &\quad\le
 Cb^{\ell_1+\cdots+\ell_k}\lambda^{-2}
 \left\{
  \sum_{a=1}^p\|u_a\|_{X^q}
                         \prod_{h\ne a}\|u_h\|_{X^8}
  +\|\Psi_0\|_{X^q}\prod_{h=1}^p\|u_h\|_{X^8}
 \right\}.
 \end{split}
 \label{rec:many-direction-bound}
\end{equation}
\item \emph{One perturbation direction.} For one perturbation direction the corresponding estimate is
\begin{equation}
 \begin{split}
 &\mathcal S_{q-2}\!\left(
   \lambda^3D^{1+k}\mathcal M[R]
                   [u,A_{\ell_1},\ldots,A_{\ell_k}]\right)\\
 &\quad\le Cb^{\ell_1+\cdots+\ell_k}
 \left\{\|u\|_{X^q}
    +\lambda^{-2}\bigl(
       G_8\|u\|_{X^q}+\|\Psi_0\|_{X^q}\|u\|_{X^8}\bigr)\right\}.
 \end{split}
 \label{rec:one-direction-bound}
\end{equation}
After subtracting the profile value, only the last two terms
in braces are needed.
\end{enumerate}
The same estimates hold uniformly on
$A+tY$, with the current increment replaced by $t\Psi_0$.
All products here have the finite material and spatial orders
of \eqref{rec:profile-order-count}.
\end{lemma}

\begin{proof}
\emph{Coefficient formulas.}
On a fixed interior interval,
\begin{equation}
 \begin{split}
 L_Rf={}&-a_\lambda(r,r_x)f_{xx}\\
 &+\{-\partial_pa_\lambda(r,r_x)r_{xx}
                          +\partial_pb_\lambda(r,r_x,x)\}f_x\\
 &+\{-\partial_ra_\lambda(r,r_x)r_{xx}
                          +\partial_rb_\lambda(r,r_x,x)\}f,
 \qquad r=R/\lambda .
 \end{split}
 \label{rec:interior-coefficients}
\end{equation}
A displacement $\lambda^4u$ changes its normalized arguments by
\begin{equation}
 \lambda^3(u,u_x,u_{xx}).
 \label{rec:interior-increment}
\end{equation}
On either edge,
\begin{equation}
 \begin{gathered}
 \lambda L_Rf=-A_nD_nf-B_nf_y/y+V_nf,\\
 \delta(r,\gamma,y\gamma_y)
   =\left(\lambda^3u,-\lambda^2u_y/y,
                         -\lambda^2(u_{yy}-u_y/y)\right).
 \end{gathered}
 \label{rec:edge-increment}
\end{equation}
By \eqref{at:edge-principal}--\eqref{at:edge-potential}, these
coefficients are smooth in the positive first-order factors and
affine in $y\gamma_y$. The center formulas
\eqref{at:center-completion}--\eqref{at:center-lower-tensor}
have the same structure in the full Cartesian deformation.
Thus each variation contains at most one second spatial derivative
before the $q-2$ source derivatives are taken.

\emph{The product estimate.}
Apply the Leibniz rule and Lemma~\ref{at:same-cutoff-product}
to the coefficients and directions together. On an edge, the
increments in \eqref{rec:edge-increment} have low norm
$C\|u\|_{X^8}$ and localized high norm $C\|u\|_{X^q}$.
For the direction on which the second-order operator acts,
\[
 \|\chi_5u\|_{B_5^q}+\|\chi_8u\|_{C_8^q}\le C\lambda^{-2}\|u\|_{X^q}.
\]
After multiplication of the edge equation by $\lambda$, the
source norm contributes $\lambda^2$. These two factors cancel.
The volume term in the localized product estimate contributes
\begin{equation}
 \mu_8^{1/2}
   =\left(\int_{c}^{C_e\lambda^{-1/2}}y^7\dd y\right)^{1/2}
       \le C\lambda^{-2}.
 \label{rec:long-volume}
\end{equation}
Assign that term to one direction's $X^q$ norm, using $q\ge8$.
If the highest derivative falls on a current coefficient, retain
$\|\Psi_0\|_{X^q}$ and bound every direction in $X^8$.
This gives precisely the two terms in braces in
\eqref{rec:many-direction-bound}.

On the bounded vacuum chart the volume is fixed, and
Lemma~\ref{at:natural-quotient} controls $u_y/y$ with the
prescribed cutoff. In the interior and at the center, the
additional physical increments give respectively
\[
 \lambda^3\lambda^{-5/4}\|u\|_{X^8}
       =\lambda^{7/4}\|u\|_{X^8},\qquad
 \lambda^3\|u\|_{X^8}.
\]
Both powers are nonnegative, as in \eqref{at:current-high}.
Lemma~\ref{rec:profile-columns} bounds the profile directions by
coefficient norms with factor $b^{\ell_1+\cdots+\ell_k}$;
no global low norm over the growing interval is used.
Combining the charts proves \eqref{rec:many-direction-bound}.

\emph{One direction.}
At $R=A$, Lemma~\ref{rec:profile-columns} gives
$Cb^{\ell_1+\cdots+\ell_k}\|u\|_{X^q}$.
For the coefficient difference, integrate its variation along
$A+tY$. The same product estimate gives
\[
 Cb^{\ell_1+\cdots+\ell_k}\lambda^{-2}
 \bigl(G_8\|u\|_{X^q}
       +\|\Psi_0\|_{X^q}\|u\|_{X^8}\bigr).
\]
This proves \eqref{rec:one-direction-bound} and its subtracted
version. Replacing $\Psi_0$ by $t\Psi_0$ gives the same estimates
uniformly along the segment.
\end{proof}

\begin{lemma}
\label{rec:bell-remainder}
For $q\ge2$, $j+q\le23$, each remainder in the iterated chain rule
with at least two perturbation factors is bounded by
\begin{equation}
 C\eta\,\mathcal X_{Q_*},
 \qquad
 Q_*=8\ \text{if }q\le8,\qquad Q_*=23\ \text{if }q\ge9.
 \label{rec:bell-remainder-bound}
\end{equation}
An additional profile material index contributes its stated
power of $b$. When the high norm falls on the current coefficient,
at least one further perturbation factor supplies the low norm.
\end{lemma}

\begin{proof}
\emph{The range $2\le q\le8$.}
The $X^8$ norm of a perturbation factor is available when its time
index is at most $15$. The force and viscous force have total time
indices at most $j$ and $j+1\le22$, respectively. Thus
\[
 i_{(1)}\ge16
 \quad\Longrightarrow\quad
 \sum_{h\ge2}i_{(h)}\le22-16=6.
\]
Taylor differences add only index-zero directions. There is therefore
at most one factor without an $X^8$ bound. Estimate this factor in
$X^q$ and the others in $X^8$; if every factor has an $X^8$ bound,
choose any one as the high-order factor. For each remaining factor,
\[
 0\le r\le2
 \quad\Longrightarrow\quad
 \partial_y^r u\text{ has }8-r\ge6
 \text{ controlled spatial derivatives}.
\]
Lemma~\ref{at:low-rung-multiplier} then bounds the product at source
order $q-2\le6$. Profile directions are bounded coefficient
multipliers by Lemma~\ref{rec:profile-columns}.

In the interior an additional perturbation factor costs
at most $C\lambda^{-1}G_8$; on an edge it costs $CG_8$.
Since
\[
 G_8=\lambda^2\eta,\qquad
 \lambda^{-1}G_8=\lambda\eta\le\eta_*,
\]
the remaining finitely many factors are bounded by a fixed
polynomial of the low-norm bound. The first additional perturbation
can be bounded by $\lambda^{-2}G_8=\eta$. Summing the possible
exceptional directions gives \eqref{rec:bell-remainder-bound}.
No $X^8$ norm of that direction has been used.

\emph{The range $9\le q\le23$.}
Here $j\le14$. After removing the principal viscous derivative,
all positive perturbation time indices are at most $j$, so every
such factor has an $X^8$ bound. Apply
\eqref{rec:many-direction-bound}, or the subtracted version
of \eqref{rec:one-direction-bound}. For example the
high current-coefficient contribution with two directions
is bounded by
\[
 C\lambda^{-2}\|\Psi_0\|_{X^q}
                    \|u_1\|_{X^8}\|u_2\|_{X^8}
 \le C\eta G_8\,\|\Psi_0\|_{X^q}
 \le C\eta\,\mathcal X_{23}.
\]
All higher products are treated in the same way. A coefficient
difference paired with one direction has instead the factor
$\lambda^{-2}\|\Psi_0\|_{X^q}\|u\|_{X^8}$, also bounded
by $C\eta\mathcal X_{23}$. Thus no product contains two
independently unbounded high spatial norms.
\end{proof}

\subsection{The differentiated equations}

In the linearized equation, a time derivative of the profile
operator acts on a lower time derivative of the perturbation and still
has spatial order two. The nonlinear remainder contains at least two
perturbation factors, one of which is small in the low norm. We separate
these terms in the commutators $\mathfrak C_j$ and $\mathfrak C_j^{\rm vis}$ below,
retaining the factor $\nu$ in every viscous term.

Define the scalar coefficients by
\begin{equation}
 \alpha_r=\lambda^{-3}\partial_s^r\lambda^3,\qquad
 \mathfrak d_r=(\nu\lambda^3)^{-1}\partial_s^r(\nu\lambda^3).
 \label{rec:scalar-arrays}
\end{equation}
For $\kappa=0$, define $\mathfrak d_r$ by the polynomial
recurrence below, without using the quotient in
\eqref{rec:scalar-arrays}. Direct differentiation for $\kappa>0$
yields
\begin{equation}
 \begin{gathered}
 \alpha_0=\mathfrak d_0=1,\qquad
 \alpha_{r+1}=\partial_s\alpha_r+3b\alpha_r,\qquad
 \mathfrak d_{r+1}=\partial_s\mathfrak d_r+\frac32b\mathfrak d_r,\\
 |\alpha_r|+|\mathfrak d_r|\le C_rb^r,\qquad
 |\partial_s^\ell b|\le C_\ell b^{\ell+1}.
 \end{gathered}
 \label{rec:scalar-bounds}
\end{equation}
The constants use only the finite orders and $b_{\max}$.

Define
\begin{align}
 \mathfrak C_j
   ={}&\lambda^{-4}\sum_{\ell=0}^j\binom j\ell
       (\lambda^3)^{(j-\ell)}
       \{\partial_s^\ell\mathcal M[R]
                       -\partial_s^\ell\mathcal M[A]\}
                         -L_R\Psi_j,
 \label{rec:force-commutator}\\
 \mathfrak C_j^{\rm vis}
   ={}&\lambda^{-4}\sum_{\ell=0}^j\binom j\ell
       (\nu\lambda^3)^{(j-\ell)}
       \{\partial_s^{\ell+1}\mathcal M[R]
                       -\partial_s^{\ell+1}\mathcal M[A]\}
                         -\nu L_R\Psi_{j+1}.
 \label{rec:viscous-commutator}
\end{align}
The following bounds apply to these complete commutators.

\begin{lemma}
\label{rec:commutator-bounds}
At $q\ge2$, $j+q\le23$,
\begin{align}
 \mathcal S_{q-2}(\mathfrak C_j)
   &\le C\sum_{i<j}b^{j-i}\|\Psi_i\|_{X^q}
                                  +C\eta\mathcal X_{Q_*},
 \label{rec:force-bound}\\
 \mathcal S_{q-2}(\mathfrak C_j^{\rm vis})
   &\le C\nu(b+\eta)\mathcal X_{Q_*}.
 \label{rec:viscous-bound}
\end{align}
Here $Q_*$ is as in Lemma~\ref{rec:bell-remainder}; in the
full range of derivative orders either occurrence may be replaced by $23$.
The first sum is empty at $j=0$. In
\eqref{rec:force-bound}, the operators retain spatial order two and act on time derivatives
of order less than $j$.
All terms in \eqref{rec:viscous-bound} retain $\nu$,
including the second-order spatial terms acting on $\Psi_j$.
\end{lemma}

\begin{proof}
\emph{The force commutator.}
For $\ell\ge1$, the chain rule gives
\begin{equation}
 \partial_s^\ell\mathcal M[R]
 =\sum_{\pi\in\Pi_\ell}
       D^{|\pi|}\mathcal M[R]
              [\,\partial_s^{|B|}R:B\in\pi\,],
 \qquad
 \partial_s^iR=A_i+\lambda^4\Psi_i.
 \label{rec:partition-formula}
\end{equation}
Here $\Pi_\ell$ consists of the partitions of $\{1,\ldots,\ell\}$.
Differentiation either adds a singleton block or adds the new index
to an existing block. Thus partitions with $m_r$ blocks of size $r$
have multiplicity
\begin{equation}
 \frac{\ell!}{\prod_{r=1}^{\ell}m_r!(r!)^{m_r}},
 \qquad \sum_{r=1}^{\ell}r m_r=\ell,
 \label{rec:partition-multiplicity}
\end{equation}
Repeated directions are kept as separate arguments when expanding
$R_i=A_i+\lambda^4\Psi_i$.

For $j\ge1$, the one-block term cancels $L_R\Psi_j$.
The linear part at $Y=0$ is
\begin{equation}
 \mathfrak C_j^{\rm lin}
   =\lambda^{-4}\partial_s^j(L_AY)-L_A\Psi_j
   =\sum_{\ell=1}^j\binom j\ell
                        L_A^{(\ell)}\Psi_{j-\ell}.
 \label{rec:linear-force}
\end{equation}
At $j=0$, Taylor's formula gives directly
\begin{equation}
 \mathfrak C_0
   =-\lambda^7\int_0^1t\,
        D^2\mathcal M[A+tY][\Psi_0,\Psi_0]\dd t .
 \label{rec:zero-force}
\end{equation}
At higher orders, apply Taylor's integral remainder along
$t(Y,Y_s,\ldots,\partial_s^jY)$ to the finite partition sum.
After its linear part is removed, each term contains at least two
perturbation factors. Coefficient differences are expanded along
$A+tY$, so a high coefficient norm is accompanied by a perturbation
factor. Lemmas~\ref{rec:profile-columns} and~\ref{rec:bell-remainder}
give
\[
 \mathcal S_{q-2}(\mathfrak C_j^{\rm lin})
       \le C\sum_{i<j}b^{j-i}\|\Psi_i\|_{X^q},\qquad
 \mathcal S_{q-2}(\mathfrak C_j-\mathfrak C_j^{\rm lin})
       \le C\eta\mathcal X_{Q_*}.
\]
This proves \eqref{rec:force-bound}, with the first sum empty at
$j=0$.

\emph{The viscous commutator.}
The one-block cancellation removes $\Psi_{j+1}$; every remaining
variation has time index at most $j$. Set
\[
 K_A=\partial_sL_A-3bL_A,
 \qquad
 \lambda^3\partial_s(D\mathcal M[A]Y)=L_AY_s+K_AY.
\]
The linear terms are therefore
\begin{equation}
 \mathfrak C_j^{\rm vis,lin}
   =\sum_{r=1}^j\binom jr
              \partial_s^r(\nu L_A)\Psi_{j-r+1}
     +\sum_{r=0}^j\binom jr
              \partial_s^r(\nu K_A)\Psi_{j-r}.
 \label{rec:linear-viscosity}
\end{equation}
Using $\nu_s=-3b\nu/2$,
\eqref{rec:scalar-bounds} and Lemma~\ref{rec:profile-columns},
\[
 \begin{split}
 \mathcal S_{q-2}(\mathfrak C_j^{\rm vis,lin})
 &\le C\nu\left\{
       \sum_{r=1}^jb^r\|\Psi_{j-r+1}\|_{X^q}
       +\sum_{r=0}^jb^{r+1}\|\Psi_{j-r}\|_{X^q}\right\}\\
 &\le C\nu b\mathcal X_{Q_*}.
 \end{split}
\]
Every nonlinear term retains $\nu$ and contains two perturbation
factors. Lemma~\ref{rec:bell-remainder} bounds their sum by
$C\nu\eta\mathcal X_{Q_*}$, proving \eqref{rec:viscous-bound}.

For completeness, at order zero the exact cancellation is
\begin{align}
 \mathfrak C_0^{\rm vis}
  &=\nu\lambda^{-4}(L_R-L_A)A_s
 \nonumber\\*
  &=\nu(\partial_sL_A-3bL_A)\Psi_0
    +\nu\lambda^7\int_0^1(1-t)
       D^3\mathcal M[A+tY][\Psi_0,\Psi_0,A_s]\dd t .
 \label{rec:zero-viscosity}
\end{align}
Indeed,
\[
 \partial_sL_A-3bL_A
       =\lambda^3D^2\mathcal M[A][A_s,\cdot].
\]
Symmetry gives the displayed linear term; Taylor's formula for
$L_R$ gives the integral. Its two perturbation directions are
estimated by \eqref{rec:many-direction-bound}, and $A_s$ supplies
$b$. All operator derivatives above are taken at fixed mass,
including derivatives of the coordinate factors and exterior scale.
\end{proof}

Differentiating \eqref{rec:equation} now gives, exactly,
\begin{equation}
 \begin{split}
 L_R(\Psi_j+\nu\Psi_{j+1})
   ={}&F_j-\Psi_{j+2}+\frac32b\Psi_{j+1}
       +\frac32\sum_{\ell=1}^j\binom j\ell
                      b^{(\ell)}\Psi_{j-\ell+1}\\
     &-\mathfrak C_j-\mathfrak C_j^{\rm vis}\\
     &-\lambda^{-4}\sum_{\ell=0}^j\binom j\ell
         (\nu\lambda^3\Lambda)^{(j-\ell)}
                                  \lambda^4\Psi_{\ell+1}.
 \end{split}
 \label{rec:row-equation}
\end{equation}
In the last line the inner $\lambda^4\Psi_{\ell+1}$ is
$\partial_s^{\ell+1}Y$, not a factor outside the
differentiation defining the coefficient. Since $\Lambda$
is fixed, this last line equals
\[
 -\nu\lambda^3\Lambda
          \sum_{\ell=0}^j\binom j\ell \mathfrak d_{j-\ell}\Psi_{\ell+1}.
\]
Its strong source norm is at most
$C\nu\lambda^3\Lambda\mathcal X_{Q_*}$.
The largest derivative is $j+1$, but it occurs at spatial
source order $q-2$, so its total order is at most $22$.
The term $\ell=1$ in the first line of
\eqref{rec:row-equation} is an inertial term of the same time order. It is retained at lower spatial order.

\subsection{The triangular spatial estimate}

Each use of the second-order inverse trades two spatial derivatives
for the kinetic term two time derivatives higher. Thus the kinetic term at $(q,j)$ uses $(q-2,j+2)$, already
controlled at a lower spatial order. Profile commutators instead lead to $(q,i)$ with $i<j$.
We therefore proceed by increasing $q$ and, for each $q$, increasing $j$,
with $2\le q\le Q_*$ and $0\le j\le23-q$.
We write $(r,i)\prec(q,j)$ for this strict order, including the base index pairs with $r=0,1$. Terms involving indices not earlier in this order retain their small
coefficients until the finite substitution is complete.

\begin{lemma}
\label{rec:node}
Put $Z_j=\Psi_j+\nu\Psi_{j+1}$ and assume
$0\le\nu\le1$. For each index pair in the stated range,
\begin{equation}
 \begin{split}
 \|Z_j\|_{X^q}
  \le{}&C\left\{\mathcal S_{q-2}(F_j)+B
       +\sum_{(r,i)\prec(q,j)}\|\Psi_i\|_{X^r}\right\}\\
 &+C\{\eta+\nu\eta+\nu b+\nu\lambda^3\Lambda\}
                                      \mathcal X_{Q_*}.
 \end{split}
 \label{rec:node-bound}
\end{equation}
In the inviscid case use $Q_*=23$ and $Z_j=\Psi_j$.
For the viscous estimate below use only $Q_*=8$.
\end{lemma}

\begin{proof}
Lemma~\ref{loc:finite-strong-graph} provides compatible approximants
to $Z_j$ converging in $X^q$, whose images under $L_R$ converge in
$\mathcal S_{q-2}$. Thus $Z_j$ lies in the strong domain of $L_R$. For an inviscid limit, Subsection~\ref{exlim:realization}
establishes this domain condition. We may therefore apply
\eqref{at:current-inverse-low} to \eqref{rec:row-equation} when
$q\le8$, and \eqref{at:current-inverse-high} when $q\ge9$.

The kinetic term and the lower-order terms in the inverse estimate
use, respectively,
\[
 (q-2,j+2),\qquad (q-2,j),\qquad (q-2,j+1)\prec(q,j).
\]
The last term is multiplied by $\nu\le1$. For $q=2,3$ these
indices lie in the base norms; otherwise they have already been
estimated at a smaller spatial order. The scalar terms also have
smaller spatial order. The profile commutators in
\eqref{rec:force-bound} use
\[
 (q,i)\prec(q,j),\qquad i<j,
\]
and are therefore controlled by the induction at fixed $q$.

For $q\ge9$, the inverse estimate contains the additional term
\begin{equation}
 \lambda^{-2}\|\Psi_0\|_{X^q}\|Z_j\|_{X^8}
   \le(1+\nu)\eta\|\Psi_0\|_{X^q},
 \label{rec:inverse-high-branch}
\end{equation}
because $j+1\le15$ and both companion $X^8$ norms
occur in $G_8$. At $q\le8$ the low inverse has no
such term. Combining these facts with
\eqref{rec:force-bound}--\eqref{rec:viscous-bound}
and the exact shift estimate proves
\eqref{rec:node-bound}. The bound $\nu\le1$ makes the constants in this induction
independent of viscosity.
\end{proof}

\subsection{Time dependence of the localized norms}

The elliptic estimate controls $Z_j=\Psi_j+\nu\Psi_{j+1}$.
We recover $\Psi_j$ from the scalar relaxation equation below.
Its kernel decays as the scale increases, whereas the moving vacuum
charts can increase the $X^q$ norm. The transport estimate below
controls this increase. After division by $\lambda^{J+1/2}$, the
normalized relaxation kernel is bounded.

\begin{lemma}
\label{rec:transport}
There are constants $\Gamma_q$, depending only on the fixed cover
and $q$, such that the localized weighted norms of any fixed function of the
mass variable satisfy
\begin{equation}
 \|f\|_{X^q(s)}
    \le\exp\!\left(\Gamma_q\int_\sigma^s b(u)\dd u\right)
                                     \|f\|_{X^q(\sigma)}
 \quad(s\ge\sigma).
 \label{rec:norm-transport}
\end{equation}
These constants are fixed before $\beta,J,K$.
\end{lemma}

\begin{proof}
With $d=M-x$, the edge coordinate gives the exact
identities
\begin{equation}
 \partial_y=n\lambda^{4/n}d^{1-1/n}\partial_d,\qquad
 \lambda^4y^{n-1}\dd y=\frac{\dd d}{n}.
 \label{rec:transport-measure}
\end{equation}
The measure is independent of $s$ in mass coordinates. For a function
$f$ fixed in mass, \eqref{rec:physical-derivative} gives
\[
 \partial_s\partial_y^k(\chi_n f)
   =\frac{4k}{n}b\,\partial_y^k(\chi_n f)
      +\partial_y^k((\partial_s\chi_n)f),\qquad n=5,8.
\]
The outer cutoff is fixed in $d$. For the inner split,
\[
 \partial_s\zeta(D)=-4bD\zeta'(D),\qquad
 \operatorname{supp}\zeta'\subset[1/2,2].
\]
Restrict the upper scale so that $2\lambda^4\le d_0$. On this
support, $\chi_e=1$ and $\chi_5+\chi_8=1$. The coordinates
$y_5=D^{1/5}$ and $y_8=D^{1/8}$, their derivatives, and their
weights are uniformly comparable. Apply the product rule after
writing $f=\chi_5f+\chi_8f$ there. It gives
\[
 \begin{split}
 \lambda^2\bigl(
   \|(\partial_s\chi_5)f\|_{B_5^q}
  +\|(\partial_s\chi_8)f\|_{C_8^q}\bigr)
 &\le C_qb\lambda^2
       \bigl(\|\chi_5f\|_{B_5^q}+\|\chi_8f\|_{C_8^q}\bigr)\\
 &\le C_qb\|f\|_{X^q(s)}.
 \end{split}
\]
This estimate uses only the overlap of the two boundary charts;
it involves no division by a cutoff. The center and interior norms
of $f$ are independent of $s$.

Differentiating the full squared norm therefore gives
\begin{equation}
 \frac{\dd}{\dd s}\|f\|_{X^q(s)}^2
       \le2\Gamma_qb\|f\|_{X^q(s)}^2.
 \label{rec:transport-differential}
\end{equation}
Its integration proves \eqref{rec:norm-transport}.
The constants depend only on the fixed cutoffs and the finite
derivative orders; both norms use the prescribed cover.
\end{proof}

With $p=2J+1$, for a nonnegative quantity $H$ write
\[
 H^\sharp(s)=\sup_{s_*\le\sigma\le s}
                     \lambda(\sigma)^{-p/2}H(\sigma).
\]
\begin{corollary}
\label{rec:relaxation}
Let $p=2J+1$, and suppose
\begin{equation}
 4+\frac p2\ge\max_{0\le q\le8}\Gamma_q .
 \label{rec:transport-threshold}
\end{equation}
If $\kappa>0$, then for $q\le8$,
\begin{equation}
 \|\Psi_j\|_{X^q}^{\sharp}(s)
 \le\lambda_*^{-p/2}\|\Psi_j(s_*)\|_{X^q(s_*)}
                         +\|Z_j\|_{X^q}^{\sharp}(s),
 \qquad \lambda_*=\lambda(s_*).
 \label{rec:relaxation-bound}
\end{equation}
\end{corollary}

\begin{proof}
From \eqref{rec:rows},
$\partial_s\Psi_j=\Psi_{j+1}-4b\Psi_j$.
Consequently
\begin{equation}
 \partial_s\Psi_j+(\nu^{-1}+4b)\Psi_j
                              =\nu^{-1}Z_j.
 \label{rec:relaxation-equation}
\end{equation}
Applying variation of constants, Minkowski's inequality, and
\eqref{rec:norm-transport}, we obtain
\[
 \begin{split}
 \lambda(s)^{-p/2}\|\Psi_j(s)\|_{X^q(s)}
 &\le K_q(s,s_*)\lambda_*^{-p/2}\|\Psi_j(s_*)\|_{X^q(s_*)}\\
 &\quad+\int_{s_*}^s K_q(s,\sigma)\nu(\sigma)^{-1}
        \lambda(\sigma)^{-p/2}\|Z_j(\sigma)\|_{X^q(\sigma)}\dd\sigma,
 \end{split}
\]
where
\begin{equation}
 K_q(s,\sigma)
  =\exp\!\left[-\int_\sigma^s
       \left\{\nu(u)^{-1}
                +(4+p/2-\Gamma_q)b(u)\right\}\dd u\right].
 \label{rec:relaxation-kernel}
\end{equation}
Under \eqref{rec:transport-threshold},
\[
 0\le K_q(s,\sigma)
       \le\exp\!\left(-\int_\sigma^s\nu(u)^{-1}\dd u\right),
 \qquad
 \int_{s_*}^sK_q(s,\sigma)\nu(\sigma)^{-1}\dd\sigma\le1 .
\]
The initial coefficient is also at most one. Taking the normalized
supremum proves \eqref{rec:relaxation-bound}, with no loss in $\nu$.

\end{proof}

We use the relaxation estimate only for $q\le8$. Extending it through
order $23$ would require the stronger condition
$4+p/2\ge\max_{q\le23}\Gamma_q$.

\subsection{Recovery of spatial derivatives}

The following estimates bound the mixed spatial norms in terms of
the first-order quantity $B$ and the source. The estimate
$B\le C_0b^{-1}\widetilde E^{1/2}$ will follow from the pressure
form in Corollary~\ref{en:base-bound}.

\begin{proposition}
\label{rec:recovery}
\label{exlim:upstream-recovery}
Fix the finite profile, its prescribed positive lower bounds,
and the weighted spaces on the coordinate cover. There are positive constants
$\eta_*,\varepsilon_{\mathrm{rec}}$ and $C_{\mathrm{rec}}$
with the following properties.

\begin{enumerate}[label=(\roman*),leftmargin=2em]
\item \emph{Inviscid estimate.} If $\kappa=0$, the compatible strong solutions of
\eqref{rec:equation} satisfying $\eta\le\eta_*$ obey
\begin{equation}
 \mathcal X_{23}(s)
       \le C_{\mathrm{rec}}\{B(s)+\mathcal R_{23}(s)\}.
 \label{rec:inviscid-recovery}
\end{equation}
\item \emph{Viscous estimate.} If $\kappa>0$, assume $0<\nu\le1$ and the transport
threshold \eqref{rec:transport-threshold}. On a forward
interval $[s_*,s]$, suppose
\begin{equation}
 \sup_{s_*\le\sigma\le s}
    \{\eta+\nu\eta+\nu b+\nu\lambda^3\Lambda\}(\sigma)
                                \le\varepsilon_{\mathrm{rec}} .
 \label{rec:small-coefficients}
\end{equation}
Then
\begin{equation}
 \mathcal X_8^\sharp(s)
  \le C_{\mathrm{rec}}\left\{
     \lambda_*^{-p/2}\mathcal X_8(s_*)
                       +B^\sharp(s)+\mathcal R_8^\sharp(s)\right\}.
 \label{rec:viscous-recovery}
\end{equation}
\end{enumerate}
The constants may depend on $J,K,\beta$, the selected
finite coefficient bounds, the fixed cutoffs, and the
positive lower bounds. They do not depend on
$\lambda_*$, $\kappa$, the length of the long chart,
or an uncontrolled high norm of the current radius.
The base estimates do not require an unweighted norm of $\Psi_{24}$.
\end{proposition}

There are $47$ pairs with $q=0,1$, $133$ with $2\le q\le8$,
and $120$ with $9\le q\le23$, in each case with $j+q\le23$.
Thus the viscous estimate involves $180$ norms and the inviscid estimate
involves $300$. The constant in the absorption argument is determined
by the estimates below.

\begin{proof}
\emph{The inviscid estimate.}
Apply \eqref{rec:node-bound} in increasing $q$ and then increasing
$j$, with $Q_*=23$ and $Z_j=\Psi_j$. Apart from the source, each
quantity in the first line has already been estimated or belongs to
the base norm. Arrange the remaining norms in a vector $U$ in this order. The inequalities
have the componentwise form
\[
 U\le W+TU+C\eta\,\mathcal X_{23}\mathbf1 ,
\]
where $T$ is strictly lower triangular with
nonnegative bounded entries, and the sum of the
components of $W$ is at most $C(B+\mathcal R_{23})$.
If there are $N$ index pairs, then
\[
 T^N=0,\qquad (I-T)^{-1}=\sum_{k=0}^{N-1}T^k,\qquad
 U\le\left(\sum_{k=0}^{N-1}T^k\right)
                  (W+C\eta\mathcal X_{23}\mathbf1).
\]
The inverse has nonnegative entries bounded in terms of the component estimates. Adding the base norms gives
\begin{equation}
 \mathcal X_{23}\le C_*(B+\mathcal R_{23})
                                  +C_*\eta\mathcal X_{23}.
 \label{rec:inviscid-substitution}
\end{equation}
The constant $C_*$ depends on the entries of the finite triangular
system through the displayed inverse. Choose
$\eta_*\le(2C_*)^{-1}$. Absorption proves
\eqref{rec:inviscid-recovery}.

\emph{The viscous estimate.}
For $\kappa>0$, use only the index pairs with $q\le8$.
Take the normalized supremum of each component estimate and apply
\eqref{rec:relaxation-bound}.
For nonnegative component norms $a_1,\ldots,a_N$,
\[
 \sup_{[s_*,s]}\lambda^{-p/2}\sum_{i=1}^Na_i
 \le\sum_{i=1}^N\sup_{[s_*,s]}\lambda^{-p/2}a_i
 \le N\sup_{[s_*,s]}\lambda^{-p/2}\sum_{i=1}^Na_i.
\]
Thus the same finite inverse applies to the individual time-dependent norms.
Their initial terms sum to
$\lambda_*^{-p/2}\mathcal X_8(s_*)$.
Consequently
\begin{equation}
 \begin{split}
 \mathcal X_8^\sharp
 \le{}&C_*\{
    \lambda_*^{-p/2}\mathcal X_8(s_*)
                   +B^\sharp+\mathcal R_8^\sharp\}\\
 &+C_*\sup_{[s_*,s]}
       \{\eta+\nu\eta+\nu b+\nu\lambda^3\Lambda\}
                                      \mathcal X_8^\sharp .
 \end{split}
 \label{rec:viscous-substitution}
\end{equation}
Choose $\varepsilon_{\mathrm{rec}}\le(2C_*)^{-1}$,
also small enough to preserve the positivity bounds and the inverse estimates.
Absorption proves \eqref{rec:viscous-recovery}.
All choices follow the selected profile and recovery
constants, whereas the transport constants were
already fixed before $J,\beta,K$.

For the finite local solutions, Proposition~\ref{en:input-finite-calculus}
and the approximation in Lemma~\ref{loc:finite-strong-graph} justify
the preceding applications of the strong inverse estimates to $Z_j$
and its differentiated force. The term $\nu\Psi_{j+1}$ can require
spatial derivatives beyond the uniformly controlled range. These
derivatives are supplied by the finite local class, and their higher
norms enter neither the constants nor the conclusions of the estimates.

Approximation in the closed strong domains gives the stated
extension. For an inviscid limit, the limiting equation must first
show that the derivatives lie in these domains.
\end{proof}

\begin{corollary}
\label{rec:actual-source}
For the matched forcing, with no additional prescribed
source, one has
\begin{equation}
 \mathcal R_{23}
 \le C_{J,K,\beta}
 \left\{\lambda^{(K-5)/2}(1+|\log\lambda|)^{p_{\log}}
       +\kappa\lambda^{-9/2}
       +\kappa\Lambda\lambda^{-3/2}\right\}.
 \label{rec:strong-source-budget}
\end{equation}
If $K>2J+6$, the upper bound for the scale may be decreased so that
the first term is at most $C\lambda^{J+1/2}$.
The bounds on viscosity at the terminal scale
\begin{equation}
 \kappa\le c\lambda_*^{J+5},\qquad
 \kappa\Lambda\le c\lambda_*^{J+2}
 \label{rec:terminal-caps}
\end{equation}
make the other two terms at most
$Cc\lambda^{J+1/2}$ for $\lambda\ge\lambda_*$.
Thus
\begin{equation}
 \mathcal R_{23}\le C\lambda^{J+1/2},
 \qquad
 \mathcal R_8^\sharp\le C .
 \label{rec:normalized-source-budget}
\end{equation}
\end{corollary}

\begin{proof}
The substitution $m=q-2$ is a bijection from
$2\le q\le23$, $0\le j\le23-q$, onto
$m,j\ge0$, $j+m\le21$. Consequently
\eqref{mat:strong-source-export} gives
\eqref{rec:strong-source-budget} in the same topology.
The exponent gap
$(K-2J-6)/2>0$ absorbs its finite logarithm.
For the other terms,
\[
 \kappa\lambda^{-9/2}
 \le c\lambda^{J+1/2}
             \left(\frac{\lambda_*}{\lambda}\right)^{J+5},
 \qquad
 \kappa\Lambda\lambda^{-3/2}
 \le c\lambda^{J+1/2}
             \left(\frac{\lambda_*}{\lambda}\right)^{J+2}.
\]
This proves the first conclusion. Since $p/2=J+1/2$
and $\mathcal R_8\le\mathcal R_{23}$, the normalized supremum bound follows.
\end{proof}

After the constants have been fixed, the same bounds make $\nu$,
$\nu b$ and $\nu\lambda^3\Lambda$ small. Control of $\eta$ and
$\nu\eta$ also requires the energy bound for $B$, proved in
Section~\ref{en:section}. The energy estimate measures the forcing
in the mass Hilbert norm through time order $22$; spatial recovery
uses the strong source norms above. Section~\ref{ct:section} closes
these coupled estimates on a common interval.
\section{The nonlinear modified energy}
\label{en:section}

At order zero we use the nonlinear relative potential, whose first
variation gives the exact force difference. For the differentiated
equations we use the Hessian at the current radius. The highest
commutator is bounded in the dual of the pressure form domain, but
the highest time derivative is controlled only in $L^2(\dd x)$.
We therefore add the commutator paired with the displacement to the
energy. Its time derivative cancels the term that cannot be estimated
directly.

The spatial coercivity estimate controls the form norm on the orthogonal
complement of the homologous direction. It does not imply stability of collapse:
Lemma~\ref{en:homogeneous-acoustics} gives infinitely many growing
relative acoustic modes in the forward homogeneous reference problem.
In the reversed time used for construction,
Theorem~\ref{en:differential-energy} allows the energy to grow. We fix
its growth constant before choosing the preparation order in
Section~\ref{ct:section}, so that the initial error and forcing remain
small after integration.
Corollary~\ref{en:base-bound} bounds the first-order quantity $B$
from Section~\ref{rec:chapter} by the energy. Together with spatial
recovery, this estimate closes the low-order bootstrap on a common
interval.

We work with the matched radius $A$ of
Proposition~\ref{en:input-matched-profile} and the localized weighted norms
of Definition~\ref{at:atlas}. The selected
Goldreich--Weber parameter is $\delta=-\beta^2/2$. Its conserved
mass is denoted interchangeably by $M=M_\beta$, so the mass interval
is $(0,M_\beta)$, and
\begin{equation}
 H=L^2((0,M_\beta),\dd x),\qquad
 \rho_R=(4\pi R^2R_x)^{-1},\qquad
 \mathcal M[R]=4\pi R^2\partial_xP(\rho_R)+\frac{x}{R^2}.
 \label{en:force}
\end{equation}
The pressure is normalized by $P(0)=0$ and
\begin{equation}
 P'(\rho)=\frac{4\rho^{2/3}}{3\sqrt{1+\rho^{2/3}}}.
 \label{en:eos}
\end{equation}
Its mechanical potential is
\begin{equation}
 \mathcal U[R]=\int_0^{M_\beta}\varepsilon(\rho_R)\dd x
                    -\int_0^{M_\beta}\frac{x}{R}\dd x,
 \qquad \varepsilon'(\rho)=\frac{P(\rho)}{\rho^2}.
 \label{en:potential}
\end{equation}
For the potential differences and variations in this section the
additive constant in $\varepsilon$ cancels. When discussing absolute
conserved energy we use the normalization in \eqref{loc:potential}.

Fix $e\ge0$ and $0<\beta\le\beta_{\rm env}$, and take the positive
branch of $b$. The scale equations and normalized perturbations are
\begin{equation}
 \lambda_s=b\lambda,\qquad b^2=\beta^2+e\lambda,
 \qquad b_s=\frac12e\lambda,\qquad e\lambda\le\beta^2,
 \qquad
 Y=R-A,\quad Y_j=\partial_s^jY=\lambda^4\Psi_j.
 \label{en:clock}
\end{equation}
In particular $\partial_s\Psi_j=\Psi_{j+1}-4b\Psi_j$.
Here $s$ increases with reversed physical time $\tau=-t$, with
$\dd\tau/\dd s=\lambda^{3/2}$. The construction starts at a small
positive value of $\lambda$ and proceeds toward larger values.

We first establish the identities for the finite-regularity solutions
of Proposition~\ref{en:input-finite-calculus}. Their consecutive
physical time derivatives belong to $C_t\mathcal D^{48-2j}$ for
$0\le j\le24$; the derivative of order 25 belongs only to the source
space. The smooth change of time variable preserves these orders on every
compact interval with $\lambda>0$. The profile also has material
derivatives through order 24.

The same proposition identifies $D\mathcal U$ with the force and
$D^2\mathcal U$ with its Hessian form. Its Green fluxes vanish at the
physical endpoints, with orders $O(\xi^5)$ at the vacuum and $O(z^3)$
at the center. The constants in these finite differentiation arguments
may depend on the interval of positive scales. The uniform dependence
of the constants in the energy estimates is specified below.

\subsection{The pressure form and the homologous direction}

Write
\begin{equation}
 Q_R=\lambda^3D^2\mathcal U[R],\qquad
 \mathcal P_A[f]=\lambda^3\int P'(\rho_A)
                    \left(2\frac fA+\frac{f_x}{A_x}\right)^2\dd x.
 \label{en:forms}
\end{equation}
The form domain is the closed pressure form domain of
Proposition~\ref{en:input-pressure-domain}. At each positive scale it is
the completion, in the norm $\{\mathcal P_A[f]+\|f\|_H^2\}^{1/2}$,
of regular radial vectors at the center and one-sided smooth vacuum
functions with $f_{y_5}(0)=0$. This condition specifies the approximating
core; it is not an extra trace imposed on its completion. Put
\begin{equation}
 \varphi_\beta=\frac{z_\beta}{\|z_\beta\|_H},\qquad
 x=m_\beta(z_\beta),\qquad
 f=h+a\varphi_\beta,\qquad h\perp_H\varphi_\beta,
 \label{en:mode}
\end{equation}
and define, once and for all,
\begin{equation}
 V(f)^2=\mathcal P_A[f]+\|f\|_H^2,
 \qquad
 Z(f)^2=\mathcal P_A[h]+\|h\|_H^2+b^2|a|^2.
 \label{en:norms}
\end{equation}
The vector $\varphi_\beta$ is fixed in time. It is not the normalized
current radius.

Fix positive weights $\omega_0,\ldots,\omega_{22}$. We distinguish
constants uniform over the leading family from those depending on
the selected matched profile.

The constants $c_0,C_0$, with a finite index when needed, depend
only on the leading profile bounds, the covering cutoffs, and the
weights. They use profile time derivatives through order 23,
first spatial strains through spatial order 24, and constitutive
derivatives through order 26. They are fixed before $J$ and before
selecting $\beta>0$. The larger Taylor order, chosen after $J$,
does not enter these constants.

The nonleading terms contribute the error
\begin{equation}
 \epsilon_A(\lambda)
       =C_{J,K,\beta}\lambda^{1/2}(1+|\log\lambda|)^{p_{\log}}.
 \label{en:profile-error}
\end{equation}
After choosing the profile and $\beta$, we decrease the upper
scale bound until $\epsilon_A\le c_0\beta^2$.

We also choose $b_{\max}\le1$ using only the leading bounds, fixed
cutoffs, and weights, and restrict $\beta_{\rm env}\le b_{\max}/\sqrt2$.
The scale restriction in \eqref{en:clock} then gives
\[
 \beta\le b=(\beta^2+e\lambda)^{1/2}
       \le\sqrt2\,\beta\le b_{\max}.
\]
The bounds $b_{\max}$ and $\beta_{\rm env}$ are fixed before $J$
and the selected matched profile.

For a positive radius set $a_f=f/R$, $d_f=f_x/R_x$, and $G_f=2a_f+d_f$.

\begin{lemma}
\label{en:strain-control}
For each admissible radius $R$,
\begin{align}
 Q_R[f,g]&=I_R[f,g]-\int\frac{2\lambda^3x}{R^3}fg\dd x,
 \label{en:component-hessian}\\
 I_R[f,g]&=\lambda^3\int
 \left[\left(P'-\frac P\rho\right)G_fG_g
                  +\frac P\rho(2a_fa_g+d_fd_g)\right]\dd x,
 \label{en:internal-form}\\
 I_R[f,f]&\sim\lambda^3\int P'(\rho_R)(a_f^2+d_f^2)\dd x.
 \label{en:component-comparison}
\end{align}
For the matched profile,
\begin{equation}
 I_A[f,f]=\mathcal P_A[f]
                 +\int8\pi\lambda^3A P_x(\rho_A)f^2\dd x
       \le \mathcal P_A[f]+C_0\|f\|_H^2.
 \label{en:profile-component-control}
\end{equation}
\end{lemma}

\begin{proof}
The exact density variations are
\[
 D\rho[f]=-\rho G_f,\qquad
 D^2\rho[f,g]=\rho(G_fG_g+2a_fa_g+d_fd_g).
\]
Substitution in \eqref{en:potential} gives
\eqref{en:component-hessian}--\eqref{en:internal-form}.
Moreover,
\[
 \frac13\le\frac{uP''(u)}{P'(u)}\le\frac23,
 \qquad
 (u/\rho)^{2/3}\le\frac{P'(u)}{P'(\rho)}
                         \le(u/\rho)^{1/3},\quad 0<u\le\rho.
\]
Integrating the second inequality yields
$3/5\le P/(\rho P')\le3/4$, which proves the positive quadratic-form
comparison \eqref{en:component-comparison}.

The first variation and its endpoint flux are given in
\eqref{loc:potential-force}. For the second variation, the difference
between the internal form and the pressure square is
$-8\pi\lambda^3\int P(\rho_A)\partial_x(Af^2)\dd x$.
Integration by parts gives the equality in
\eqref{en:profile-component-control}; its boundary term vanishes by
the same compatible endpoint orders.

It remains to bound the coefficient of the zeroth-order term.
In either vacuum chart it is a smooth expression in
\[
 A/\lambda,\qquad -A_y/(\lambda^2y),\qquad
 y\partial_y\{-A_y/(\lambda^2y)\}.
\]
Choose the upper bound for $\lambda$ so that the errors in the radius,
Jacobian, and first spatial strains are at most one. The matched-profile
bounds then keep these factors in a fixed range determined by the
leading family. At the center, $A=\lambda z a(z^2)$ gives
$\lambda^3AP_x=O(1)$; the fixed interior satisfies the same bound
without a coordinate singularity. These estimates prove the last
inequality with $C_0$ uniform in the leading family. The selected
higher coefficient constants enter only the outer-scale restriction.
\end{proof}

For $0\le\beta\le\beta_{\rm env}$, let
$T_\beta f(z)=f(m_\beta(z))$ and
\[
 \dd\omega_\beta=4\pi w_\beta(z)^3z^2\dd z,\qquad
 \dd\omega_*=4\pi(1-z)^3z^2\dd z,\qquad
 H_*=L^2((0,1),\dd\omega_*).
\]
\begin{lemma}
\label{en:common-seed-pivot}
The following statements hold uniformly for $0\le\beta\le\beta_{\rm env}$.
\begin{enumerate}[label=(\roman*),leftmargin=2em]
\item The map $T_\beta:L^2(0,M_\beta)\to L^2(\dd\omega_\beta)$ is
an isometry. The norms $L^2(\dd\omega_\beta)$ are uniformly equivalent
to the norm of $H_*$, and the pulled-back reference form is
\[
 \mathfrak p_\beta[f]=\frac{16\pi}{3}
        \int_0^1w_\beta^4z^4\left|\partial_z(T_\beta f/z)\right|^2\dd z.
\]
\item If $0\le\beta_n\le\beta_{\rm env}$ and
\[
 \sup_n\{\mathfrak p_{\beta_n}[f_n]+\|f_n\|_{L^2(0,M_{\beta_n})}^2\}
       <\infty,
\]
then $T_{\beta_n}f_n$ has a strongly convergent subsequence in $H_*$.
\item Suppose $\beta_n\to\bar\beta$, the form norms of $f_n$ are bounded,
and $T_{\beta_n}f_n\to f$ strongly in $H_*$. Then
\[
 f_n\perp\varphi_{\beta_n}\quad\Longrightarrow\quad
 f\perp T_{\bar\beta}\varphi_{\bar\beta}
       \quad\text{in }L^2(\dd\omega_{\bar\beta}).
\]
\end{enumerate}
\end{lemma}

\begin{proof}
The leading profiles satisfy $w_\beta(z)=(1-z)v_\beta(z)$, with
\[
 0<c\le v_\beta\le C,\qquad
 \dd\omega_\beta=v_\beta^3\dd\omega_*,\qquad
 T_\beta\varphi_\beta=z/N_\beta,\quad
 N_\beta^2=\int_0^1z^2\dd\omega_\beta.
\]
The factors $v_\beta$ and their required derivatives depend
continuously on $\beta$. The change of variables $x=m_\beta(z)$
proves the isometry, the norm comparison, and the formula for the
reference form in (i). It also gives $0<c\le N_\beta\le C$ and
continuity of $N_\beta$.

To prove (ii), put $g=T_\beta f/z$, $\mu_\beta=w_\beta^3z^4$,
and $a_\beta=w_\beta^4z^4$. On a fixed interior interval, the
form norm controls the $H^1$ norm and an interior value $g(z_0)$.
Cauchy--Schwarz gives
\[
 \begin{aligned}
 |g(z)|^2
 &\le2|g(z_0)|^2+
     2\left(\int_0^1a_\beta|g'|^2\dd r\right)
      \left|\int_{z_0}^z a_\beta^{-1}\dd r\right|\\
 &\le C\{\mathfrak p_\beta[f]+\|f\|_H^2\}
                 \{1+z^{-3}+(1-z)^{-3}\}.
 \end{aligned}
\]
Since $\mu_\beta\sim z^4(1-z)^3$, integration yields
\begin{align}
 \int_0^\epsilon\mu_\beta|g|^2\dd z
  &\le C\epsilon^2\{\mathfrak p_\beta[f]+\|f\|_H^2\},
 \label{en:center-tail}\\
 \int_{1-\epsilon}^1\mu_\beta|g|^2\dd z
  &\le C\epsilon\{\mathfrak p_\beta[f]+\|f\|_H^2\}.
 \label{en:vacuum-tail}
\end{align}
The constants are uniform for $0\le\beta\le\beta_{\rm env}$.
Interior Rellich compactness and these tail bounds give a Cauchy
subsequence in $H_*$. Approximation extends the estimates to the
closed form domain and gives
\begin{equation}
 \mathcal V_\beta:=
 \overline{\{zg:g\in C^\infty([0,1])\}}
 ^{\,\{\mathfrak p_\beta+\|\cdot\|_H^2\}^{1/2}}
       \Subset L^2(\dd\omega_\beta).
 \label{en:reference-compactness}
\end{equation}

Finally, under the hypotheses of (iii),
\[
 \begin{aligned}
 0=(f_n,\varphi_{\beta_n})_{L^2(0,M_{\beta_n})}
 &=\int_0^1(T_{\beta_n}f_n)\frac{zv_{\beta_n}^3}{N_{\beta_n}}
                                      \dd\omega_*\\
 &\longrightarrow
   \int_0^1f\frac{z}{N_{\bar\beta}}\dd\omega_{\bar\beta}.
 \end{aligned}
\]
Indeed $T_{\beta_n}f_n\to f$ in $H_*$ and
$zv_{\beta_n}^3/N_{\beta_n}\to zv_{\bar\beta}^3/N_{\bar\beta}$
in $H_*$. This proves (iii). The reference estimates are uniform
in $\beta$; the comparison with the exact pressure form below is
only one-sided.
\end{proof}

\begin{lemma}
\label{en:complement-coercivity}
There are $\beta_{\rm env}>0$ and $c_0,C_0>0$, determined by
the leading profile bounds and fixed cutoffs, with the following
property. For each selected $0<\beta\le\beta_{\rm env}$ and its
finite matched profile, the upper scale bound can be chosen so that
\begin{equation}
 c_0V(h)^2\le Q_A[h,h]\le C_0V(h)^2,
 \qquad h\perp_H\varphi_\beta.
 \label{en:complement-bound}
\end{equation}
The scale bound may depend on the selected matching coefficients;
the constants $c_0,C_0$ do not.
\end{lemma}

\begin{proof}
\emph{The reference form.}
Had\v{z}i\'c--Jang \cite[Lemma~3.1 and Appendix~A]{HadzicJangExpansion}
proved a spectral gap uniform for sufficiently small negative $\delta$.
We record a weighted variance proof in the present normalization
before comparing the reference and matched forms. Let
$Q_{\rm ref}$ be the normalized Hessian at $\lambda z_\beta$ for
$P(\rho)=\rho^{4/3}$. With $f=zg$ and
$\dd x=4\pi w_\beta^3z^2\dd z$, the Goldreich--Weber equation gives
\begin{equation}
 Q_{\rm ref}[f,f]=\mathfrak p_\beta[f,f]-\beta^2\|f\|_H^2,
 \qquad
 \mathfrak p_\beta[f,f]
       =\frac{16\pi}{3}\int_0^1w_\beta^4z^4|g'|^2\dd z.
 \label{en:reference-form}
\end{equation}
Indeed,
\[
 \frac{16\pi}{3}\int w_\beta^4z^2(zg'+3g)^2\dd z
 =\mathfrak p_\beta[f,f]
                  -64\pi\int w_\beta^3w_\beta'z^3g^2\dd z;
\]
the remaining pressure and gravity terms cancel the last integral
and leave $-\beta^2\|f\|_H^2$.

Put $\mu=w_\beta^3z^4$, $a=w_\beta^4z^4$, and
$\bar g=(\int\mu g)/(\int\mu)$. The variance identity gives
\[
 \begin{aligned}
 \int\mu|g-\bar g|^2
 &=\frac1{\int\mu}\int_{0<r<t<1}
       \mu(r)\mu(t)|g(t)-g(r)|^2\dd r\dd t\\
 &\le\frac{\int a|g'|^2}{\int\mu}
      \int_0^1\frac{(\int_0^u\mu)(\int_u^1\mu)}{a(u)}\dd u.
 \end{aligned}
\]
Thus
\begin{equation}
 \int\mu|g-\bar g|^2\le C_0\int a|g'|^2,
 \qquad
 C_0\le C\frac{\displaystyle\int_0^1
       \frac{(\int_0^t\mu)(\int_t^1\mu)}{a(t)}\dd t}
                         {\displaystyle\int_0^1\mu}.
 \label{en:variance}
\end{equation}
Uniformly in the leading family,
$\mu\sim z^4(1-z)^3$ and $a\sim z^4(1-z)^4$; the last
integrand is $O(u)$ at zero and $O(1)$ at one. Since
$h\perp_H\varphi_\beta$ is equivalent to $\bar g=0$, decreasing
$\beta_{\rm env}$ gives
\begin{equation}
 Q_{\rm ref}[h,h]\ge c_0\{
                  \mathfrak p_\beta[h,h]+\|h\|_H^2\},
 \qquad h\perp_H\varphi_\beta.
 \label{en:reference-gap}
\end{equation}

\emph{Comparison of the form domains.}
We claim
\begin{equation}
 \mathfrak p_\beta[f,f]+\|f\|_H^2
                       \le C_0\{\mathcal P_A[f]+\|f\|_H^2\}.
 \label{en:one-sided-domain}
\end{equation}
Write $d=M_\beta-x$, $y_n=(d/\lambda^4)^{1/n}$ and
$\gamma=-A_{y_n}/(\lambda^2y_n)$. The coefficient of $f_x^2$ is
\begin{equation}
 \frac{\lambda^3P'(\rho_A)}{A_x^2}
     =n^2\lambda^7y_n^{2n-2}a_n,
 \qquad a_n=\frac{P'(\rho_A)}{y_n^2\gamma^2}.
 \label{en:actual-principal-coefficient}
\end{equation}
For $n=8$ it is comparable to $d^{7/4}$ throughout the long chart.
For $n=5$ its factor $y_5^8$ dominates the reference factor
$y_5^{35/4}$ on the bounded chart. Away from the center, the other
reference coefficient is bounded. At the center, the positive
normalized deformations and \eqref{en:component-comparison}
control both reference strains. Now
\eqref{en:profile-component-control} proves
\eqref{en:one-sided-domain}. Hence the compactness of
Lemma~\ref{en:common-seed-pivot} applies to every bounded sequence
in the matched pressure form norm.

\emph{The matched form.}
The upper bound in \eqref{en:complement-bound} follows from
Lemma~\ref{en:strain-control}. If no uniform lower bound held for
small scales and small selected profile errors, there would be
$\lambda_n\to0$, errors tending to zero, $\beta_n\to\bar\beta$,
and $h_n\perp\varphi_{\beta_n}$ with
\[
 \mathcal P_{A_n}[h_n]+\|h_n\|_H^2=1,
 \qquad Q_{A_n}[h_n,h_n]\le1/n.
\]
Set $\widehat h_n=T_{\beta_n}h_n$ and
$\widehat A_n=A_n\circ m_{\beta_n}$. After passage to a subsequence,
\[
 \begin{gathered}
 \widehat h_n\to h\text{ in }H_*,
 \qquad \widehat h_n\rightharpoonup h
       \text{ in }H^1_{\rm loc}(0,1),\\
 h\perp T_{\bar\beta}\varphi_{\bar\beta}
       \text{ in }L^2(\dd\omega_{\bar\beta}).
 \end{gathered}
\]
On compact interior intervals, the positive matrices of $I_{A_n}$
and their positive square roots converge uniformly to those of
$I_{\rm ref}$. Weak lower semicontinuity and exhaustion give
\begin{equation}
 I_{\rm ref}[h,h]\le\liminf_n I_{A_n}[h_n,h_n].
 \label{en:internal-liminf}
\end{equation}
The gravity coefficients
\[
 g_n=\frac{2\lambda_n^3m_{\beta_n}}{\widehat A_n^3}
\]
are uniformly bounded, including at the center, since
$m_\beta(z)\sim z^3$ and $\widehat A_n\sim\lambda_n z$.
They converge on compact interior intervals to
$2m_{\bar\beta}(z)/z^3$. Strong convergence in $H_*$, uniform
convergence of the measure ratios, and the tails
\eqref{en:center-tail}--\eqref{en:vacuum-tail} imply
\begin{equation}
 \int g_n|\widehat h_n|^2\dd\omega_{\beta_n}
 \longrightarrow
 \int\frac{2m_{\bar\beta}(z)}{z^3}|h|^2\dd\omega_{\bar\beta}.
 \label{en:gravity-convergence}
\end{equation}
By \eqref{en:internal-liminf} and \eqref{en:gravity-convergence},
$Q_{\rm ref}[h,h]\le0$, and
\eqref{en:reference-gap} gives $h=0$. It follows that
\[
 \begin{gathered}
 \|h_n\|_H\longrightarrow0,\\
 0\le I_{A_n}[h_n,h_n]
 =Q_{A_n}[h_n,h_n]+
     \int g_n|\widehat h_n|^2\dd\omega_{\beta_n}
 \le n^{-1}+o(1)\longrightarrow0,\\
 \mathcal P_{A_n}[h_n]
 \le C I_{A_n}[h_n,h_n]\longrightarrow0.
 \end{gathered}
\]
The last inequality uses $|2a+d|^2\le5(a^2+d^2)$ and
\eqref{en:component-comparison}. This contradicts the normalization.
All constants are uniform over the leading family; the selected
matching coefficients determine only the required upper bound for
$\lambda$.
\end{proof}

The preceding coercivity is a spatial property of the pressure form.
To distinguish it from forward dynamics, consider the homogeneous
law $P_{\rm h}(\rho)=\rho^{4/3}$ and its zero-energy homologous
collapse. Introduce the forward rescaled time $\theta$ by
$\dd t/\dd\theta=\lambda^{3/2}$ and $\lambda_\theta=-\beta\lambda$.
For this reference solution, $\theta$ has the opposite orientation
to the construction time $s$.
\begin{lemma}
\label{en:homogeneous-acoustics}
For this reference problem, the following assertions hold.
\begin{enumerate}[label=(\roman*),leftmargin=2em]
\item In the relative displacement $R=\lambda z(1+g)$ the linearized equation is
\begin{equation}
 g_{\theta\theta}-\frac\beta2 g_\theta
       +(\mathcal L_\beta-3\beta^2/2)g=0,
 \label{en:homogeneous-linearization}
\end{equation}
\equationalias{in:linearized-core}{en:homogeneous-linearization}
where $\mathcal L_\beta$ is the nonnegative self-adjoint realization
of $\mathfrak p_\beta$ in $L^2(4\pi w_\beta^3z^4\dd z)$, expressed in
$g$. Its kernel consists of the constants and its positive eigenvalues
$\mu_j$ tend to infinity.
\item The high-mode characteristic roots are
\begin{equation}
 \sigma_j^\pm=\frac\beta4
             \pm\sqrt{\frac{25\beta^2}{16}-\mu_j}.
 \label{en:homogeneous-roots}
\end{equation}
Thus infinitely many homogeneous relative-displacement modes have
oscillatory envelope $\lambda^{-1/4}$.
\end{enumerate}
These conclusions concern the linear reference problem above.
\end{lemma}

\begin{proof}
In the $g$ variable, use the inner product and norm of
$L^2(4\pi w_\beta^3z^4\dd z)$. The closed nonnegative form defines
$\mathcal L_\beta$ by minimization of
\[
 \frac{16\pi}{3}\int_0^1w_\beta^4z^4|g'|^2\dd z
                         +\|g\|^2-2(F,g).
\]
Its resolvent $(\mathcal L_\beta+1)^{-1}$ maps boundedly into the
form domain and is compact by \eqref{en:reference-compactness}.
It is also self-adjoint and injective. The compact self-adjoint
spectral theorem gives a complete eigenbasis and $\mu_j\to\infty$.
Moreover,
\[
 \int_0^1w_\beta^4z^4|g'|^2\dd z=0
 \quad\Longleftrightarrow\quad g'=0\text{ on }(0,1),
\]
so the kernel consists of the constants.

For the evolution, write $R=\lambda(z+f)$, $f=zg$. Since
$\partial_t=\lambda^{-3/2}\partial_\theta$ and
$\lambda_\theta=-\beta\lambda$,
\[
 R_{tt}=\lambda^{-2}
   \left\{f_{\theta\theta}-\frac\beta2f_\theta
                         -\frac{\beta^2}{2}(z+f)\right\}.
\]
The homogeneous force satisfies
\[
 \mathcal M_{\rm h}[\lambda(z+f)]
       =\lambda^{-2}\mathcal M_{\rm h}[z+f],\qquad
 \mathcal M_{\rm h}[z]=\frac{\beta^2}{2}z.
\]
Thus the exact equation is
\[
 f_{\theta\theta}-\frac\beta2f_\theta
 -\frac{\beta^2}{2}(z+f)+\mathcal M_{\rm h}[z+f]=0.
\]
By \eqref{en:reference-form},
$D\mathcal M_{\rm h}[z](zg)=z(\mathcal L_\beta-\beta^2)g$.
Linearization therefore gives \eqref{en:homogeneous-linearization}.

For $\mathcal L_\beta e_j=\mu_je_j$, the ansatz
$g=e^{\sigma\theta}e_j$ yields
\[
 \sigma^2-\frac\beta2\sigma+\mu_j-\frac{3\beta^2}{2}=0.
\]
Its roots are \eqref{en:homogeneous-roots}; at $\mu_j=0$ they are
$3\beta/2$ and $-\beta$. For every $\mu_j>25\beta^2/16$,
\[
 \Re\sigma_j^\pm=\frac\beta4,\qquad
 |e^{\sigma_j^\pm\theta}|=e^{\beta\theta/4}
              =\left(\frac{\lambda(\theta)}{\lambda(0)}\right)^{-1/4}.
\]
There are infinitely many such eigenvalues, proving (ii).
\end{proof}

This growth concerns the relative displacement $g$ in the linear
reference problem. Eliminating finitely many growing modes leaves
infinitely many growing acoustic modes, which motivates prescribing
the asymptotic solution at small terminal scales. This calculation
does not establish nonlinear instability, stability, or a codimension
for the solutions of the exact pressure law constructed here.

\subsection{Time derivatives of the profile and the Hessian}

For $A_i=\partial_s^iA$ define the rescaled Hessian derivatives by
\begin{equation}
 Q_A^{\langle k\rangle}=\lambda^3\partial_s^kD^2\mathcal U[A].
 \label{en:covariant-form-definition}
\end{equation}
Thus the exterior factor $\lambda^3$ is not differentiated in this
notation; in particular $Q_{A,s}=3bQ_A+Q_A^{\langle1\rangle}$.

\begin{lemma}
\label{en:profile-columns}
On the scale range fixed above, the following estimates hold for
all $f,g$ in the pressure form domain.
\begin{enumerate}[label=(\roman*),leftmargin=2em]
\item For $0\le i,k\le23$,
\begin{align}
 \left\|\frac{A_i}{A}\right\|_\infty+
 \left\|\frac{(A_i)_x}{A_x}\right\|_\infty
       &\le C_{0,i}b^i,
 \label{en:profile-strains}\\
 |Q_A^{\langle k\rangle}[f,g]|
       &\le C_{0,k}b^kZ(f)Z(g),
 \qquad |Q_{A,s}[f,g]|\le C_0bZ(f)Z(g).
 \label{en:profile-form-columns}
\end{align}
\item There is $\Lambda_{\rm en}\ge4$, depending only on the leading
profile bounds and fixed cutoffs, such that
\begin{equation}
 c_0Z(f)^2\le Q_A[f,f]+\Lambda_{\rm en}b^2\|f\|_H^2
                                                   \le C_0Z(f)^2.
 \label{en:profile-shifted-coercivity}
\end{equation}
\end{enumerate}
\end{lemma}

\begin{proof}
The coefficient bounds give estimates in the pressure form norm.
To obtain the weaker weight on the homologous component, we compare
the Hessian with the homogeneous reference and estimate the contribution
of the boundary layer to that component. The finite expansion below
also keeps track of the material derivatives.

\emph{Material differentiation.} In either vacuum chart set
$E_n=\lambda\partial_\lambda-(4/n)y\partial_y$,
$\theta=e\lambda/b^2$, and
\[
 r_i=\frac{A_i}{\lambda b^i},\qquad
 \gamma_i=-\frac{(A_i)_y}{\lambda^2yb^i}.
\]
Direct differentiation yields
\begin{equation}
 r_{i+1}=(E_n+1+i\theta/2)r_i,\qquad
 \gamma_{i+1}=(E_n+2-8/n+i\theta/2)\gamma_i,
 \qquad E_n\theta=\theta(1-\theta).
 \label{en:profile-rate-recurrence}
\end{equation}
Here $0\le\theta\le1/2$. The matched-profile estimates control these
finite recurrences. In the core region $q\ge c\sqrt\lambda$, where
$q=1-z$, the nonleading relative Jacobian terms are bounded by powers
of $\lambda/q$. In the layer they are bounded by
$(\lambda D^{1/4})^{n-1}$ for $n\ge2$. Both bounds are included in
\eqref{en:profile-error}. They persist under the required iterates of
$E_n$ in the vacuum charts and $\lambda\partial_\lambda|_z$ in the core.

On the matching shell the leading relation
$\lambda Z_0(D)=q+O(q^2+\lambda)$ gives relative cutoff errors
$O(q+\lambda/q)=O(\sqrt\lambda)$. In $(A_{23})_x$, a derivative
on the cutoff leaves only 23 material derivatives on the mismatch;
an uncut term uses spatial order 24. The center and fixed interior use
their ordinary finite profile bounds. This proves
\eqref{en:profile-strains} with the asserted leading constants.

\emph{Hessian derivatives.}
With the exterior scale factor held fixed, form at each $(s,x)$ the
truncated Taylor polynomials
\[
 a(t)=1+\sum_{i=1}^{23}\frac{A_i}{i!A}t^i,\qquad
 d(t)=1+\sum_{i=1}^{23}\frac{(A_i)_x}{i!A_x}t^i,
 \quad u(t)=a(t)^{-1},\quad v(t)=d(t)^{-1}\pmod {t^{24}}.
\]
Their coefficients are determined without any infinite expansion by
\begin{equation}
 u_0=v_0=1,\qquad
 u_m=-\sum_{i=1}^m\frac{A_i}{i!A}u_{m-i},\qquad
 v_m=-\sum_{i=1}^m\frac{(A_i)_x}{i!A_x}v_{m-i}.
 \label{en:column-reciprocal}
\end{equation}
Set $n(t)=u(t)^2v(t)=1+\sum_{m=1}^{23}n_mt^m$. For either
$F=P'$ or $F=P/\rho$ the coefficient of its exact composition is
\begin{equation}
 [t^m]F(\rho_A n(t))
 =\sum_{\substack{\alpha_1,\ldots,\alpha_m\ge0\\
                         \sum i\alpha_i=m}}
   \rho_A^{|\alpha|}F^{(|\alpha|)}(\rho_A)
              \prod_{i=1}^m\frac{n_i^{\alpha_i}}{\alpha_i!},
 \qquad m\ge1,
 \label{en:column-bell}
\end{equation}
with constant coefficient $F(\rho_A)$. Write
$a_f(t)=(f/A)u(t)$, $d_f(t)=(f_x/A_x)v(t)$ and
$G_f(t)=2a_f(t)+d_f(t)$. For $0\le k\le23$, the rescaled Hessian derivative is
\begin{equation}
\begin{split}
 Q_A^{\langle k\rangle}[f,g]=\lambda^3 k![t^k]\int
 \bigg[&\left(P'(\rho_A n)-\frac{P(\rho_A n)}{\rho_A n}\right)
                                      G_f(t)G_g(t)\\
 &+\frac{P(\rho_A n)}{\rho_A n}
                  \{2a_f(t)a_g(t)+d_f(t)d_g(t)\}\\
 &-\frac{2x}{A^3}u(t)^3fg\bigg]\dd x.
\end{split}
 \label{en:all-column-formula}
\end{equation}
Only degrees at most $k$ enter \eqref{en:all-column-formula}.
The finite recursions \eqref{en:column-reciprocal} and
\eqref{en:column-bell}, together with \eqref{en:profile-strains}, give
\[
 |u_m|+|v_m|+|n_m|\le C_mb^m,\qquad 1\le m\le23.
\]
Every degree-$m$ coefficient in the Hessian thus contains the
factor $b^m$; truncation discards no coefficient used in the proof.

For the finite orders used here, the exact pressure law satisfies
\[
 |\rho^r(P')^{(r)}(\rho)|+
 |\rho^r(P/\rho)^{(r)}(\rho)|\le C_rP'(\rho).
\]
Indeed,
\[
 (\rho\partial_\rho)P'(\rho)
 =\left(\frac23-\frac{\rho^{2/3}}{3(1+\rho^{2/3})}\right)P'(\rho).
\]
Repeated application of $\rho\partial_\rho$ gives the same bound at
each fixed order, since the coefficient and its iterates are bounded.
The identity $P(\rho)/\rho=\int_0^1P'(t\rho)\dd t$ gives the second estimate. The remaining factors in
\eqref{en:all-column-formula} are therefore bounded by the integrand
in \eqref{en:component-comparison}, which controls the individual
strains, and by a bounded coefficient in the gravity term.

One spatial derivative of $A_i$, $i\le23$, has total derivative
order at most 24. No time derivative of order 24 or 25 enters this formula.

In particular the exterior $\lambda^3$ in
\eqref{en:all-column-formula} has been held fixed. Differentiating also the exterior factor gives the identity
\begin{equation}
 \partial_s Q_A^{\langle k\rangle}
       =3bQ_A^{\langle k\rangle}+Q_A^{\langle k+1\rangle},
 \qquad 0\le k\le22.
 \label{en:column-transport}
\end{equation}
For the homologous reference the unnormalized Hessian has degree
$-3$ in $\lambda$; hence $p_0=1$,
$p_{k+1}=\partial_sp_k-3bp_k$ gives
$Q_{\rm ref}^{\langle k\rangle}=p_kQ_{\rm ref}$ and
$|p_k|\le C_kb^k$.

\emph{Coefficient comparison.}
Put $q=1-z_\beta(x)$, $D=(M_\beta-x)/\lambda^4$, and
\[
 \mathcal J_A(f)=\lambda^3P'(\rho_A)
        \{(f/A)^2+(f_x/A_x)^2\}+f^2.
\]
We subtract the matched and reference formulas one factor at a time
in the finite products \eqref{en:column-reciprocal}--\eqref{en:all-column-formula}.
The factor $b^k$ is common to all terms. The remaining coefficient
estimates depend on the region as follows.
\begin{itemize}[leftmargin=2em]
\item \emph{Core.}
On $q\ge c_2\sqrt\lambda$, write
\[
 A=\lambda z\psi,\qquad J_c=\psi^2(\psi+z\psi_z),\qquad
 \rho_A=\lambda^{-3}w_\beta^3/J_c.
\]
The normalized radius and Jacobian errors are bounded by $\epsilon_A$,
with the required material derivatives. The additional pressure factor is
\[
 \left(1+(\lambda/w_\beta)^2J_c^{2/3}\right)^{-1/2}
       =1+O((\lambda/q)^2).
\]
The same error bound holds after the required iterates of
$\lambda\partial_\lambda|_z$ and $q\partial_q$.
Thus the difference of the two integrands is bounded by
\[
 C_kb^k\epsilon_A\sqrt{\mathcal J_A(f)\mathcal J_A(g)}.
\]
This region includes the outer part of the long chart.
\item \emph{Matching region.}
Here $q\sim\sqrt\lambda$. For a cutoff $\chi(q/\sqrt\lambda)$,
\[
 q\partial_q\chi(q/\sqrt\lambda)
   =\frac q{\sqrt\lambda}\chi'(q/\sqrt\lambda),\qquad
 \lambda\partial_\lambda|_q\chi(q/\sqrt\lambda)
   =-\frac q{2\sqrt\lambda}\chi'(q/\sqrt\lambda).
\]
These derivatives and their required iterates are bounded.
When a spatial derivative falls on the cutoff, it multiplies the
matched difference. The resulting relative error is
$O(q+\lambda/q)=O(\sqrt\lambda)$, through total derivative order $24$.
After enlarging $\epsilon_A$, the core majorant therefore also
bounds the difference on the shell.
\item \emph{Exact-pressure region.}
On $q\le c_1\sqrt\lambda$, equivalently $0\le D\le C\lambda^{-2}$
up to fixed comparison constants, we use an absolute estimate.
The positive normalized factors and the pressure derivative bounds
above give the majorant
\[
 C_kb^k\sqrt{\mathcal J_A(f)\mathcal J_A(g)}.
\]
The reference integrand satisfies the same bound: its principal
coefficient is dominated by \eqref{en:actual-principal-coefficient},
and its remaining coefficient contributes at most $Cf^2$ to the
quadratic form. This region includes the inner part of the long chart;
no small relative pressure error is used there.
\end{itemize}
The radius derivatives have total order at most $k+1\le24$ in all
three regions. Cauchy--Schwarz and \eqref{en:component-comparison}
make the displayed majorants integrable on the pressure form domain.
The constants use only the positive leading factors, fixed cutoffs,
and $k\le23$. The selected higher matching coefficients occur in
$\epsilon_A$ and in the upper restriction on $\lambda$, not in $C_k$.

\emph{The contribution near the vacuum.} The coefficient difference
need not be small on
$\mathcal L_\lambda=\{M_\beta-x\le C\lambda^2\}$. Define
\begin{equation}
 W_\lambda(f)^2=\lambda^3\int_{\mathcal L_\lambda}P'(\rho_A)
       \left\{\left|\frac fA\right|^2+
                    \left|\frac{f_x}{A_x}\right|^2\right\}\dd x
                         +\int_{\mathcal L_\lambda}|f|^2\dd x.
 \label{en:layer-form}
\end{equation}
Then $W_\lambda(f)\le C_0V(f)$. For the fixed mode and
$D=(M_\beta-x)/\lambda^4$,
\[
 \left|\frac{\varphi_\beta}{A}\right|\le C\lambda^{-1},
 \qquad
 \left|\frac{(\varphi_\beta)_x}{A_x}\right|
 \le C\lambda^{-1}
 \begin{cases}D^{-3/20},&D\le1,\\1,&D\ge1.
 \end{cases}
\]
The exact pressure law gives $P'(\rho_A)\le CD^{2/5}$ for $D\le1$ and
$P'(\rho_A)\le CD^{1/4}$ for $D\ge1$. Consequently
\begin{equation}
 \begin{split}
 W_\lambda(\varphi_\beta)^2
 &\le C\lambda^5\int_0^1(1+D^{1/10})\dd D
       +C\lambda^5\int_1^{C\lambda^{-2}}D^{1/4}\dd D
       +C\lambda^4\int_0^{C\lambda^{-2}}\dd D\\
 &\le C_0\lambda^2.
 \end{split}
 \label{en:mode-layer-tail}
\end{equation}
The reference component form has no larger contribution on
$\mathcal L_\lambda$. The outer covering annulus lies outside
$\mathcal L_\lambda$ and satisfies the core estimate.

For $B_k=Q_A^{\langle k\rangle}-p_kQ_{\rm ref}$, the preceding
coefficient subtraction gives
\begin{equation}
 |B_k[f,g]|\le C_{0,k}b^k
       \{W_\lambda(f)W_\lambda(g)+\epsilon_A V(f)V(g)\}.
 \label{en:profile-form-difference}
\end{equation}
For $h,k\perp\varphi_\beta$, the three types of pairings, including
the reference term, therefore satisfy
\begin{align}
 |Q_A^{\langle r\rangle}[h,k]|
  &\le C_rb^rV(h)V(k),\notag\\
 |Q_A^{\langle r\rangle}[h,\varphi_\beta]|
  &\le C_rb^r(\lambda+\epsilon_A)V(h),\notag\\
 |Q_A^{\langle r\rangle}[\varphi_\beta,\varphi_\beta]|
  &\le C_rb^r(\beta^2+\lambda^2+\epsilon_A),
       \qquad 0\le r\le23.
 \label{en:three-form-blocks}
\end{align}
The identities $Q_{\rm ref}[h,\varphi_\beta]=0$ and
$Q_{\rm ref}[\varphi_\beta,\varphi_\beta]=-\beta^2$ give the
smaller coefficients in the last two estimates.
Enlarge $\epsilon_A$ to include $\lambda+\lambda^2$, and select
the upper bound for the scale so that $\epsilon_A\le c_0\beta^2$. Since
$b\ge\beta$ and $b\le1$,
\[
 \frac{\lambda+\epsilon_A}{b}\le\frac{C\beta^2}{b}\le C,
 \qquad
 \frac{\beta^2+\lambda^2+\epsilon_A}{b^2}\le C.
\]
Thus the mixed and quadratic homologous terms have respectively
the one and two factors of $b$ required by the $Z$ norm. The three
estimates prove the first bound in \eqref{en:profile-form-columns},
with constants uniform in the family of leading profiles. The identity
$Q_{A,s}=3bQ_A+Q_A^{\langle1\rangle}$ proves the second.

At $k=0$ the same calculation gives
\[
 |Q_A[\varphi_\beta,\varphi_\beta]+\beta^2|\le C_0\epsilon_A,
 \qquad |Q_A[h,\varphi_\beta]|\le C_0\epsilon_A V(h).
\]
For $f=h+a\varphi_\beta$ with $h\perp_H\varphi_\beta$, Young's
inequality and \eqref{en:complement-bound} give
\[
 \begin{split}
 Q_A[f,f]+\Lambda_{\rm en}b^2\|f\|_H^2
 &\ge \frac{c_0}{2}V(h)^2\\
 &\quad+\bigl(\Lambda_{\rm en}b^2-\beta^2
                    -C_0\epsilon_A-C_0\epsilon_A^2\bigr)|a|^2.
 \end{split}
\]
Choose $\Lambda_{\rm en}\ge4$ as a fixed leading constant and then
make $\epsilon_A/\beta^2$ sufficiently small. Since $b\ge\beta$,
the right side controls $V(h)^2+b^2|a|^2=Z(f)^2$.
The upper bound follows from the same decomposition and the already
proved form estimate. This proves \eqref{en:profile-shifted-coercivity}.
\end{proof}

\subsection{Estimates with two directions in the form domain}

We estimate two arguments in the pressure form norm and the strains
of the remaining arguments in $L^\infty$. In the differentiated
equations, the latter arguments have lower time orders and are
controlled by the available $X^8$ bounds.

Assume throughout the rest of the section that
\begin{equation}
 G_8^2=\sum_{q=0}^8\sum_{j=0}^{23-q}\|\Psi_j\|_{X^q}^2,
 \qquad \eta=\lambda^{-2}G_8\le\varepsilon_0b^3,
 \label{en:low-cone}
\end{equation}
and that the radii $A+tY$, $0\le t\le1$, have the positive normalized
radius and Jacobian bounds of the covering atlas. The latter bounds
also follow from sufficiently small low-order derivative bounds and the positive profile.
For a direction with bounded strains put
\[
 \sigma_R(l)=\|l/R\|_\infty+\|l_x/R_x\|_\infty.
\]

\begin{lemma}
\label{en:two-form-slots}
Let $2\le r\le26$, let $f,g$ lie in the pressure form domain, and
let $l_1,\ldots,l_{r-2}$ be compatible directions with bounded strains.
Then
\begin{equation}
 |\lambda^3D^r\mathcal U[A+tY][f,g,l_1,\ldots,l_{r-2}]|
 \le C_{0,r}V(f)V(g)\prod_{\alpha=1}^{r-2}
                                \sigma_{A+tY}(l_\alpha).
 \label{en:two-form-bound}
\end{equation}
Moreover, for $0\le i\le15$,
\begin{equation}
 \sup_{0\le t\le1}\sigma_{A+tY}(Y_i)\le C_0\eta.
 \label{en:low-strains}
\end{equation}
\end{lemma}

\begin{proof}
For a perturbed radius, the density satisfies
\begin{equation}
 \rho_{R+\sum t_il_i}
 =\rho_R\left(1+\sum t_i\frac{l_i}{R}\right)^{-2}
          \left(1+\sum t_i\frac{(l_i)_x}{R_x}\right)^{-1}.
 \label{en:multilinear-density}
\end{equation}
Every mixed variation has one strain from each direction. The explicit
formula \eqref{en:eos} gives
$|(\rho\partial_\rho)^kP'(\rho)|\le C_kP'(\rho)$ at every fixed
order used here. From
$\rho\varepsilon'=P/\rho$ one obtains
$| (\rho\partial_\rho)^r\varepsilon|\le C_rP/\rho$ for $r\ge1$.
Put $a_f=f/R$ and $d_f=f_x/R_x$, and likewise for $g$.
The internal-energy contribution is bounded by
\[
 \begin{aligned}
 &C_r\lambda^3\int_0^M P'(\rho_R)
       (|a_f|+|d_f|)(|a_g|+|d_g|)\dd x
                       \prod_\alpha\sigma_R(l_\alpha)\\
 &\quad\le C_r
  \left(\lambda^3\int_0^MP'(\rho_R)(a_f^2+d_f^2)\dd x\right)^{1/2}
  \left(\lambda^3\int_0^MP'(\rho_R)(a_g^2+d_g^2)\dd x\right)^{1/2}
                       \prod_\alpha\sigma_R(l_\alpha)\\
 &\quad\le C_{0,r}V(f)V(g)\prod_\alpha\sigma_R(l_\alpha).
 \end{aligned}
\]
The last step uses Lemma~\ref{en:strain-control} and the fixed
positive bounds for the density, radius and Jacobian ratios along
$R=A+tY$.

For gravity the corresponding integrand is bounded by
\[
 C_r\frac{\lambda^3x}{R^3}|fg|
                              \prod_\alpha|l_\alpha/R|.
\]
Its coefficient is bounded also at the center, where
$x\sim z^3$ and $R\sim\lambda z$. This proves
\eqref{en:two-form-bound}.

For \eqref{en:low-strains}, $Y_i=\lambda^4\Psi_i$ and $i\le15$
give the required $X^8$ bound. On the bounded edge,
\[
 \frac{(Y_i)_y}{A_y}
       =-\frac{\lambda^2(\Psi_i)_y}{y\gamma}.
\]
The local $C^2$ estimate controls the compatible quotient, with the
measure normalization included in the weighted embedding. The same
estimate applies on fixed-size intervals of the long chart.

In the fixed interior, the powers of the scale in $Y_i$ compensate
for the loss in the embedding. At the center, use the full vector
$C^2$ estimate and the quotient $Y_i/z$. The prescribed cutoffs cover
the mass interval. Combining these local bounds gives
\[
 \sup_{0\le t\le1}\sigma_{A+tY}(Y_i)
       \le C_0\lambda^{-2}G_8=C_0\eta.
\]
\end{proof}

\begin{lemma}
\label{en:current-forms}
For sufficiently small $\varepsilon_0$, uniformly for $0\le t\le1$
and for $f,g$ in the pressure form domain,
\begin{align}
 c_0Z(f)^2&\le Q_{A+tY}[f,f]+\Lambda_{\rm en}b^2\|f\|_H^2
                                                    \le C_0Z(f)^2,
 \label{en:current-coercivity}\\
 |Q_{A+tY}[f,g]|&\le C_0Z(f)Z(g),
 \qquad |Q_{R,s}[f,g]|\le C_0bZ(f)Z(g).
 \label{en:current-transport}
\end{align}
\end{lemma}

\begin{proof}
The decomposition $f=h+a\varphi_\beta$ gives
\begin{equation}
 V(f)\le C_0b^{-1}Z(f),\qquad b\|f\|_H\le C_0Z(f).
 \label{en:mode-loss}
\end{equation}
Integrating the third variation along the segment from the profile to the solution and using \eqref{en:two-form-bound} gives
\begin{equation}
 |(Q_R-Q_A)[f,g]|\le C_0\eta V(f)V(g).
 \label{en:current-form-difference}
\end{equation}
Consequently, uniformly for $0\le t\le1$,
\[
 |(Q_{A+tY}-Q_A)[f,g]|
 \le C_0\eta b^{-2}Z(f)Z(g)
 \le C_0\varepsilon_0bZ(f)Z(g).
\]
Absorb this error in \eqref{en:profile-shifted-coercivity} to obtain
\eqref{en:current-coercivity}. The profile form bound gives the first
estimate in \eqref{en:current-transport}.

At fixed arguments, differentiation of the same integral gives exactly
\begin{align}
 (Q_R-Q_A)_s[f,g]
 ={}&3b(Q_R-Q_A)[f,g]
 \nonumber\\
 &+\lambda^3\int_0^1D^4\mathcal U[A+tY]
                               [Y,A_s+tY_s,f,g]\dd t
 \nonumber\\
 &+\lambda^3\int_0^1D^3\mathcal U[A+tY][Y_s,f,g]\dd t.
 \label{en:transport-identity}
\end{align}
Since $Y_s=Y_1$, the strain bounds are
\[
 \sigma_{A+tY}(A_s)\le C_0b,\qquad
 \sigma_{A+tY}(Y)+\sigma_{A+tY}(Y_s)\le C_0\eta.
\]
The three terms in \eqref{en:transport-identity} are therefore bounded
in sum by
\[
 C_0\{b\eta+\eta(b+\eta)+\eta\}V(f)V(g)
 \le C_0\eta b^{-2}Z(f)Z(g)
 \le C_0\varepsilon_0bZ(f)Z(g).
\]
Adding the profile derivative estimate in
\eqref{en:profile-form-columns} proves the last bound in
\eqref{en:current-transport}.
\end{proof}

\subsection{The commutator correction}

The highest commutator is bounded in the dual of the form domain,
whereas the highest time derivative is controlled only in $H$.
Pairing the commutator with the displacement gives the energy correction.
Its derivative cancels the kinetic pairing in
\eqref{en:primitive-cancellation}. The same identity retains the
contribution from differentiating the exterior factor $\lambda^3$.

Set $R_j=\partial_s^jR$, $A_j=\partial_s^jA$, and extend
\eqref{en:covariant-form-definition} to the current radius:
\[
 Q_R^{\langle q\rangle}=\lambda^3\partial_s^qD^2\mathcal U[R].
\]
Ordinary differentiation of the force gives, for $j\ge1$,
\begin{equation}
\begin{aligned}
 \lambda^3\partial_s^j\mathcal M[R]
     &=Q_RR_j+\mathcal B_j[R],\\
 \mathcal B_j[R]
     &=\sum_{q=1}^{j-1}\binom{j-1}{q}
                                  Q_R^{\langle q\rangle}R_{j-q}.
\end{aligned}
 \label{en:force-partition}
\end{equation}
Consequently
\begin{equation}
\begin{aligned}
 C_j&=(Q_R-Q_A)A_j+\mathcal B_j[R]-\mathcal B_j[A],\\
 \lambda^3\partial_s^j(\mathcal M[R]-\mathcal M[A])
     &=Q_RY_j+C_j.
\end{aligned}
 \label{en:commutator-definition}
\end{equation}
The operator expressions in these formulas are interpreted in the dual
of the form domain. For $1\le j\le22$ define
\begin{equation}
 \mathfrak A_j=\lambda^{-8}C_j[Y_j],
 \qquad
 \mathfrak R_j=\lambda^{-8}(C_{j,s}-3bC_j)[Y_j].
 \label{en:primitive-definition}
\end{equation}
Using $Y_j=\lambda^4\Psi_j$ and $\partial_sY_j=Y_{j+1}$, we have
\[
 \begin{aligned}
 \partial_s\mathfrak A_j
  &=-8b\mathfrak A_j+\lambda^{-8}C_{j,s}[Y_j]
                           +\lambda^{-8}C_j[Y_{j+1}],\\
 -(\lambda^{-4}C_j)[\Psi_{j+1}+2b\Psi_j]
  &=-\lambda^{-8}C_j[Y_{j+1}]-2b\mathfrak A_j.
 \end{aligned}
\]
Since $\lambda^{-8}C_{j,s}[Y_j]=\mathfrak R_j+3b\mathfrak A_j$,
the terms containing $Y_{j+1}$ cancel and $-8-2+3=-7$ gives
\begin{equation}
 -(\lambda^{-4}C_j)[\Psi_{j+1}+2b\Psi_j]
                  +\partial_s\mathfrak A_j
       =\mathfrak R_j-7b\mathfrak A_j.
 \label{en:primitive-cancellation}
\end{equation}
Thus the highest derivative is eliminated before any estimate in the
form domain is used.

For the fixed positive weights $\omega_j$, set
\begin{equation}
 E=\sum_{j=0}^{22}\omega_j
               \{\|\Psi_{j+1}\|_H^2+Z(\Psi_j)^2\}.
 \label{en:positive-energy}
\end{equation}
\begin{lemma}
\label{en:primitive-bounds}
The commutator corrections satisfy
\begin{align}
 \left|\sum_{j=1}^{22}\omega_j\mathfrak A_j\right|
       &\le C_0(b_{\max}+\eta b^{-2})E,
 \label{en:primitive-smallness}\\
 \sum_{j=1}^{22}\omega_j|\mathfrak R_j|
       &\le C_0(b+\eta b^{-2})E\le C_0bE,
 \label{en:primitive-derivative-bound}
\end{align}
In addition, for $1\le j\le23$ and every form test $f$,
\begin{equation}
 |(\lambda^{-4}C_j)[f]|
           \le C_0b\sum_{i=0}^{j-1}Z(\Psi_i)Z(f).
 \label{en:commutator-form-bound}
\end{equation}
\end{lemma}

\begin{proof}
Differentiate \eqref{en:commutator-definition}, holding the final test
$Y_j$ fixed. The expression $(C_{j,s}-3bC_j)[Y_j]$ is the sum of
\begin{equation}
 (Q_R^{\langle1\rangle}-Q_A^{\langle1\rangle})[A_j,Y_j]
                     +(Q_R-Q_A)[A_{j+1},Y_j]
 \label{en:primitive-first-terms}
\end{equation}
and, for $1\le q\le j-1$, $\binom{j-1}{q}$ times
\begin{align}
 &Q_R^{\langle q+1\rangle}[R_{j-q},Y_j]
                   -Q_A^{\langle q+1\rangle}[A_{j-q},Y_j]
 \nonumber\\
 &\qquad
 +Q_R^{\langle q\rangle}[R_{j-q+1},Y_j]
                   -Q_A^{\langle q\rangle}[A_{j-q+1},Y_j].
 \label{en:primitive-remaining-terms}
\end{align}

\emph{Linear terms.}
Differentiating $D\mathcal M[A]Y$ gives
\begin{align}
 C_j^{\rm lin}
    &=\sum_{k=1}^j\binom jk Q_A^{\langle k\rangle}Y_{j-k},
 \label{en:commutator-linearization}\\
 \mathfrak A_j^{\rm lin}
    &=\sum_{k=1}^j\binom jk
                  Q_A^{\langle k\rangle}[\Psi_{j-k},\Psi_j],
 \label{en:primitive-linearization}\\
 \mathfrak R_j^{\rm lin}
    &=\sum_{k=1}^j\binom jk
       \bigl\{Q_A^{\langle k+1\rangle}[\Psi_{j-k},\Psi_j]
                    +Q_A^{\langle k\rangle}[\Psi_{j-k+1},\Psi_j]\bigr\}.
 \label{en:remainder-linearization}
\end{align}
The last identity differentiates the first before pairing. Its
same-order term is
\begin{equation}
 jQ_A^{\langle1\rangle}[\Psi_j,\Psi_j]
            =j(Q_{A,s}-3bQ_A)[\Psi_j,\Psi_j].
 \label{en:same-row-second-order}
\end{equation}
It is $O(b)$ in the form norm. For example, the leading $5/3$
vacuum coefficient satisfies
\[
 \lambda^7y_5^8=\lambda^{3/5}(M_\beta-x)^{8/5},\qquad
 (\partial_s-3b)(\lambda^7y_5^8)
       =-\frac{12}{5}b\lambda^7y_5^8.
\]
Thus the coefficient in \eqref{en:same-row-second-order} is
$-12jb/5$; this second-order term is retained in the form estimate.

\emph{Nonlinear terms.}
Write $U_T^{(r)}=\lambda^3D^r\mathcal U[T]$. For a partition
$\pi$ of $\{1,\ldots,j\}$ and $S\subset\pi$, put
\[
 \mathbf A_S=(A_{|B|})_{B\in S},\qquad
 \mathbf Y_S=(Y_{|B|})_{B\in S}.
\]
Blocks with equal sizes remain distinct. With $\Pi_j$ denoting all
such partitions, \eqref{rec:partition-formula}--%
\eqref{rec:partition-multiplicity} give
\begin{equation}
 \lambda^3\partial_s^j\mathcal M[R][f]
   =\sum_{\pi\in\Pi_j}U_R^{(|\pi|+1)}
                      [(R_{|B|})_{B\in\pi},f].
 \label{en:form-bell}
\end{equation}
The single-block term leaves $(Q_R-Q_A)[A_j,f]$ after subtraction
of the profile and $Q_RY_j$. Hence every remaining perturbation
index is at most $j-1$. Taylor expansion in the base radius and
expansion of the arguments give, for $C_j^{\rm nl}=C_j-C_j^{\rm lin}$
and $k=|\pi|$,
\begin{align}
 C_j^{\rm nl}[f]={}&
 \sum_{\pi\in\Pi_j}\int_0^1(1-t)
 U_{A+tY}^{(k+3)}[Y_0,Y_0,\mathbf A_\pi,f]\dd t\notag\\*
 &+\sum_{\substack{\pi\in\Pi_j\\k\ge2}}
   \sum_{\substack{S\subset\pi\\|S|=1}}\int_0^1
   U_{A+tY}^{(k+2)}
       [Y_0,\mathbf Y_S,\mathbf A_{\pi\setminus S},f]\dd t\notag\\*
 &+\sum_{\substack{\pi\in\Pi_j\\k\ge2}}
   \sum_{\substack{S\subset\pi\\|S|\ge2}}
   U_R^{(k+1)}[\mathbf Y_S,\mathbf A_{\pi\setminus S},f].
 \label{en:nonlinear-commutator-expansion}
\end{align}
The three sums contain zero, one, and at least two perturbations
among the original block arguments; after the base expansion, each
summand has at least two perturbation factors. At the first two orders,
\begin{align*}
 C_1[f]&=(Q_R-Q_A)[A_1,f],\\
 C_2[f]&=(Q_R-Q_A)[A_2,f]
     +U_R^{(3)}[A_1+Y_1,A_1+Y_1,f]
     -U_A^{(3)}[A_1,A_1,f].
\end{align*}
Their linear parts, as in \eqref{en:commutator-linearization}, are
$Q_A^{\langle1\rangle}[Y_0,f]$ and
$Q_A^{\langle2\rangle}[Y_0,f]+2Q_A^{\langle1\rangle}[Y_1,f]$;
the factor two counts the two singleton subsets.

For the differentiated remainder in \eqref{en:primitive-definition}, use
\[
 (\partial_s-3b)U_{A+tY}^{(r)}
   =U_{A+tY}^{(r+1)}[A_1+tY_1,\,\cdot\,],\qquad
 \partial_sA_i=A_{i+1},\quad\partial_sY_i=Y_{i+1}.
\]
The test remains fixed. The maximal perturbation index and total
time index change from $j-1,j$ to $j,j+1$; at least two
perturbations remain. The potential orders are at most $j+3$ for
$C_j$ and $j+4$ for its differentiated remainder.

\emph{Estimates.}
Choose the perturbation of largest index and the final test as the
two form directions. The sum of all time indices is at most $23$;
hence the second-largest perturbation index is at most $11$, as in
\eqref{rec:single-high-rule}. Since $11+8\le23$, each remaining
perturbation has the $X^8$ bound in \eqref{en:low-strains}, and at
least one supplies a factor $\eta$. Profile directions satisfy
\eqref{en:profile-strains}. For $C_{23}^{\rm nl}[f]$ and
$(\partial_s-3b)C_{22}^{\rm nl}[f]$ the counts are
\begin{equation}
 \begin{aligned}
 \text{largest perturbation time index}&\le22,\\
 \text{every remaining perturbation time index}&\le11,\\
 \text{order of the potential variation}&\le26.
 \end{aligned}
 \label{en:top-slot-budget}
\end{equation}
In particular $A_{23}$ requires spatial order $24$; neither a
form estimate for $Y_{23}$ nor a potential derivative of order $27$
is used.

Apply \eqref{en:two-form-bound}. In each primitive the two form
directions contribute $\lambda^4$ each, canceling $\lambda^{-8}$:
\begin{align}
 |\mathfrak A_j-\mathfrak A_j^{\rm lin}|
     &\le C_0\eta\sum_{i=0}^{j-1}V(\Psi_i)V(\Psi_j),
 \label{en:nonlinear-primitive}\\
 |\mathfrak R_j-\mathfrak R_j^{\rm lin}|
     &\le C_0\eta\sum_{i=0}^{j}V(\Psi_i)V(\Psi_j).
 \label{en:nonlinear-remainder}
\end{align}
For the fixed weights,
\[
 \sum_{j=1}^{22}\omega_j\sum_{i=0}^{j}Z(\Psi_i)Z(\Psi_j)
 \le C_0\sum_{i=0}^{22}\omega_iZ(\Psi_i)^2\le C_0E.
\]
Together with \eqref{en:profile-form-columns} and
\eqref{en:mode-loss}, this gives
\[
 \begin{aligned}
 \left|\sum_{j=1}^{22}\omega_j\mathfrak A_j\right|
       &\le C_0(b_{\max}+\eta b^{-2})E,\\
 \sum_{j=1}^{22}\omega_j|\mathfrak R_j|
       &\le C_0(b+\eta b^{-2})E\le C_0bE.
 \end{aligned}
\]
For an arbitrary form test $f$, the same estimates give
\[
 \begin{aligned}
 |(\lambda^{-4}C_j)[f]|
 &\le C_0\sum_{i=0}^{j-1}
       \{b^{j-i}Z(\Psi_i)Z(f)+\eta V(\Psi_i)V(f)\}\\
 &\le C_0b\sum_{i=0}^{j-1}Z(\Psi_i)Z(f),
 \end{aligned}
\]
since $b\le1$ and $\eta b^{-2}\le\varepsilon_0b$.
This proves \eqref{en:commutator-form-bound} through $j=23$;
the removed term $Q_RY_{23}$ does not enter this estimate.
\end{proof}

\subsection{The relative potential and energy coercivity}

At order zero, the relative potential along $A+tY$ differentiates
to the exact force difference paired with $\Psi_1$. The quadratic
form $Q_R[\Psi_0,\Psi_0]/2$ does not have this property.
At positive time orders, the principal term is linear in the highest
direction, and the current Hessian gives its energy. Set
\begin{align}
 V_0&=\lambda^{-5}
       \{\mathcal U[R]-\mathcal U[A]-D\mathcal U[A]Y\}
 \nonumber\\
     &=\int_0^1(1-t)Q_{A+tY}[\Psi_0,\Psi_0]\dd t,
 \label{en:relative-potential}\\
 F_0^{\rm nl}&=\lambda^{-1}(\mathcal M[R]-\mathcal M[A]),
 \nonumber\\*
 (F_0^{\rm nl},\Psi_0)_H
     &=\int_0^1Q_{A+tY}[\Psi_0,\Psi_0]\dd t.
 \label{en:zero-force}
\end{align}
Direct differentiation gives
\begin{align}
 \partial_sV_0&=-5bV_0+(F_0^{\rm nl},\Psi_1)_H+N_0,
 \label{en:relative-potential-derivative}\\
 N_0&=\lambda^{-5}
       (\mathcal M[R]-\mathcal M[A]-D\mathcal M[A]Y)[A_1]
 \nonumber\\
    &=\int_0^1(1-t)\lambda^3D^3\mathcal U[A+tY]
                                      [\Psi_0,\Psi_0,A_1]\dd t.
 \label{en:zero-remainder}
\end{align}
At the profile the last integrand is
$Q_A^{\langle1\rangle}[\Psi_0,\Psi_0]$. Subtracting this base value
and using the bound with two directions in the form norm gives
\begin{equation}
 |N_0|\le C_0bZ(\Psi_0)^2+C_0\eta bV(\Psi_0)^2
                                      \le C_0bZ(\Psi_0)^2.
 \label{en:zero-remainder-bound}
\end{equation}
Thus no homogeneity identity for the exact pressure law is required.

Define
\begin{align}
 \widetilde H_0
   &=\frac12\|\Psi_1\|_H^2+V_0+2b(\Psi_0,\Psi_1)_H
                     +\frac{13+\Lambda_{\rm en}}2b^2\|\Psi_0\|_H^2,
 \label{en:zero-energy}\\
 H_j
   &=\frac12\|\Psi_{j+1}\|_H^2+\frac12Q_R[\Psi_j,\Psi_j]
       +2b(\Psi_j,\Psi_{j+1})_H
                     +\frac{13+\Lambda_{\rm en}}2b^2\|\Psi_j\|_H^2,
 \quad 1\le j\le22,
 \label{en:positive-row-energy}\\
 \widetilde E
   &=\omega_0\widetilde H_0+
                  \sum_{j=1}^{22}\omega_j(H_j+\mathfrak A_j).
 \label{en:modified-energy}
\end{align}
The relative potential and commutator corrections at the prepared
initial time were estimated in
Proposition~\ref{en:input-prepared-energy}.

\begin{proposition}
\label{en:energy-coercivity}
Choose $b_{\max}$ and $\varepsilon_0$ sufficiently small in terms of
the leading profile bounds, fixed cutoffs, and weights, and restrict
$\beta_{\rm env}\le b_{\max}/\sqrt2$. Then
\begin{equation}
                       c_0E\le\widetilde E\le C_0E.
 \label{en:energy-equivalence}
\end{equation}
These choices and the constants $c_0,C_0$ are made before $J$
and before selecting $\beta>0$.
\end{proposition}

\begin{proof}
Completing the square gives, for $j\ge1$,
\begin{equation}
 H_j=\frac12\|\Psi_{j+1}+2b\Psi_j\|_H^2
       +\frac12\{Q_R[\Psi_j,\Psi_j]
                           +\Lambda_{\rm en}b^2\|\Psi_j\|_H^2\}
       +\frac92b^2\|\Psi_j\|_H^2.
 \label{en:energy-square}
\end{equation}
At order zero the same identity holds with $Q_R[\Psi_0,\Psi_0]$
replaced by $2V_0$. In particular, for $u=\Psi_0$ the exact
Taylor formula with integral remainder is
\begin{equation}
 2V_0+\Lambda_{\rm en}b^2\|u\|_H^2
 =2\int_0^1(1-t)
   \{Q_{A+tY}[u,u]+\Lambda_{\rm en}b^2\|u\|_H^2\}\dd t.
 \label{en:row-zero-segment-coercivity}
\end{equation}
Since $2\int_0^1(1-t)\dd t=1$, \eqref{en:current-coercivity} gives
\[
 c_0Z(u)^2\le 2V_0+\Lambda_{\rm en}b^2\|u\|_H^2
                         \le C_0Z(u)^2.
\]
For the kinetic term, use
\[
 \|v\|_H^2\le2\|v+2bu\|_H^2+8b^2\|u\|_H^2,
 \qquad b\|u\|_H\le C_0Z(u).
\]
Applying these bounds to \eqref{en:energy-square} and its order-zero
version yields
\[
 c_*E\le\omega_0\widetilde H_0+
                  \sum_{j=1}^{22}\omega_jH_j\le C_*E.
\]
The correction estimate \eqref{en:primitive-smallness} gives
\[
 \left|\sum_{j=1}^{22}\omega_j\mathfrak A_j\right|
 \le C_0(b_{\max}+\eta b^{-2})E
 \le C_0(1+\varepsilon_0)b_{\max}E\le\frac{c_*}{2}E,
\]
provided $C_0(1+\varepsilon_0)b_{\max}\le c_*/2$.
Combining the last two inequalities proves \eqref{en:energy-equivalence},
with constants uniform in the leading family.
\end{proof}

\subsection{The differentiated equation}

Suppose that $Y$ satisfies the exact forced equation
\begin{equation}
 \begin{split}
 Y_{ss}-\frac32bY_s+\lambda^3(\mathcal M[R]-\mathcal M[A])
 &+\nu\lambda^3\partial_s(\mathcal M[R]-\mathcal M[A])\\
 &+\nu\lambda^3\Lambda Y_s=\lambda^4F,
 \qquad \nu=\kappa\lambda^{-3/2},
 \end{split}
 \label{en:forced-equation}
\end{equation}
where $\Lambda\ge0$ is fixed on the positive-scale interval. Let
\begin{equation}
 K_R[f,g]=Q_R[f,g]+\lambda^3\Lambda(f,g)_H,\qquad
 F_\ell=\lambda^{-4}\partial_s^\ell(\lambda^4F).
 \label{en:shift-and-source}
\end{equation}
The operator associated with $Q_R$ is denoted by $L_R$:
$L_R=\lambda^3D\mathcal M[R]$ and $L_Rf[g]=Q_R[f,g]$.

The force has the exterior factor $\lambda^3$. Applying
$\mathbf D_3=\partial_s-3b$ to the equation removes derivatives of
this factor from the force, giving the commutators used above.
The differentiated inertial, viscous, shift, and source terms also
contain lower derivatives with scalar coefficients bounded by
powers of $b$. To compute these terms, define
\begin{equation}
 d_{0,0}=1,\qquad
 d_{j+1,\ell}=d_{j,\ell-1}+\partial_sd_{j,\ell}-3bd_{j,\ell},
 \qquad d_{j,\ell}=0\quad\hbox{outside }0\le\ell\le j.
 \label{en:scalar-recurrence}
\end{equation}
Then
\begin{equation}
 \mathbf D_3^jW=\sum_{\ell=0}^jd_{j,\ell}\partial_s^\ell W,
 \quad d_{j,j}=1,\quad d_{j,j-1}=-3jb,
 \quad |d_{j,\ell}|\le C_0b^{j-\ell}.
 \label{en:scalar-expansion}
\end{equation}
The bounds follow by induction from \eqref{en:clock}. In particular
the lower inertial terms are exactly
\begin{align}
 I_j={}&-\sum_{\ell=0}^{j-1}d_{j,\ell}\Psi_{\ell+2}
 \nonumber\\*
 &+\frac32\left\{\sum_{\ell=0}^jd_{j,\ell}
       \sum_{a=0}^{\ell}\binom\ell a
            b^{(\ell-a)}\Psi_{a+1}-b\Psi_{j+1}\right\},
 \qquad I_0=0,
 \label{en:inertial-terms}\\
 \widehat F_j={}&\sum_{\ell=0}^jd_{j,\ell}F_\ell,
 \qquad b^{(r)}=\partial_s^rb.
 \label{en:transformed-source}
\end{align}
The only term of the highest kinetic order in $I_j$ is $3jb\Psi_{j+1}$.
Every other derivative of order $i\le j$ has coefficient $O(b^{j+2-i})$. Hence
\begin{equation}
 \|I_j\|_H\le C_0b
       \left(\|\Psi_{j+1}\|_H+\sum_{i=0}^jZ(\Psi_i)\right),
 \qquad
 \sum_{j=0}^{22}\|\widehat F_j\|_H^2
                         \le C_0\sum_{\ell=0}^{22}\|F_\ell\|_H^2.
 \label{en:scalar-bounds}
\end{equation}

The viscosity has a different scalar recurrence:
\begin{equation}
 v_r=\nu^{-1}\partial_s^r\nu,\qquad v_0=1,\qquad
 v_{r+1}=v_r'-\frac32bv_r,\qquad |v_r|\le C_0b^r.
 \label{en:viscosity-recurrence}
\end{equation}
For $\kappa=0$ the same scalar polynomials define $v_r$ and every
term multiplied by $\nu$ is zero. The identity
$\mathbf D_3^j(\lambda^3W)=\lambda^3\partial_s^jW$ gives
\begin{equation}
 \begin{split}
 &\lambda^{-4}\mathbf D_3^j
       \{\nu\lambda^3\partial_s(\mathcal M[R]-\mathcal M[A])\}\\
 &\qquad
 =\nu\sum_{\ell=1}^{j+1}\binom j{\ell-1}v_{j+1-\ell}
                 \{L_R\Psi_\ell+\lambda^{-4}C_\ell\}.
 \end{split}
 \label{en:full-viscous-expansion}
\end{equation}
After extracting only $\nu L_R\Psi_{j+1}$, the remaining form is
$\nu T_j$, where
\begin{equation}
 T_j=\sum_{\ell=1}^{j}\binom j{\ell-1}v_{j+1-\ell}L_R\Psi_\ell
       +\sum_{\ell=1}^{j+1}\binom j{\ell-1}v_{j+1-\ell}
                                            \lambda^{-4}C_\ell.
 \label{en:viscous-remainder}
\end{equation}
In particular
\begin{equation}
 T_0=\lambda^{-4}C_1
           =\lambda^{-4}(Q_R-Q_A)A_1,\qquad
 |T_j[f]|\le C_0b\sum_{i=0}^jZ(\Psi_i)Z(f).
 \label{en:viscous-form-bound}
\end{equation}
The bound follows from \eqref{en:commutator-form-bound},
\eqref{en:current-transport}, and \eqref{en:viscosity-recurrence}.
It is a form bound, not a strong spatial source estimate. At the
highest energy order, it remains to estimate $C_{23}$ in $T_{22}$.
The one-block term $Q_RY_{23}$ has already been removed, so
\eqref{en:top-slot-budget} bounds $T_{22}[f]$ using perturbation
derivatives in the form domain only through order $22$.

The remaining occurrence of $\Psi_{23}$ is the extracted dissipative
term $\nu K_R\Psi_{23}$. Its estimate retains the factor $\nu$ and
does not require a viscosity-independent pressure norm for $\Psi_{23}$.

The constant shift has the separate exact expansion
\begin{align}
 \lambda^{-4}\mathbf D_3^j(\nu\lambda^3\Lambda Y_s)
       &=\nu\lambda^3\Lambda(\Psi_{j+1}+B_j),
 \label{en:full-shift-expansion}\\
 B_j&=\sum_{\ell=1}^{j}\binom j{\ell-1}v_{j+1-\ell}\Psi_\ell,
 \qquad B_0=0,\qquad
 \|B_j\|_H\le C_0\sum_{i=0}^jZ(\Psi_i).
 \label{en:shift-remainder}
\end{align}
Every lower shift coefficient contains a factor of $b$, which compensates for the loss in the Hilbert norm estimate in \eqref{en:mode-loss}.

Thus, for $1\le j\le22$, the exact differentiated equations are
\begin{equation}
 \begin{split}
 u_s&=v-4bu,\qquad u=\Psi_j,\quad v=\Psi_{j+1},\\
 v_s+\frac52bv+L_Ru+\nu K_Rv
     &=-\lambda^{-4}C_j+I_j+\widehat F_j
                          -\nu T_j-\nu\lambda^3\Lambda B_j.
 \end{split}
 \label{en:row-equation}
\end{equation}
At $j=0$, replace $L_Ru+\lambda^{-4}C_j$ by $F_0^{\rm nl}$
and use $I_0=B_0=0$. The same invertible triangular transformation
has been applied to the inertial, force, viscous, and source terms.
In particular, all scalar terms arising from the undifferentiated
force are retained.

\subsection{The differential estimate}

\begin{theorem}
\label{en:differential-energy}
Let $R=A+Y$ be a solution of \eqref{en:forced-equation} in the
finite-regularity class specified at the beginning of this section.
Assume the following conditions.
\begin{enumerate}[label=(\roman*),leftmargin=2em]
\item The bound \eqref{en:low-cone} holds, and the normalized radius
and Jacobian along $A+tY$, $0\le t\le1$, have the positive bounds
used above.
\item $F_\ell\in H$ for $0\le\ell\le22$.
\item The fixed shift and viscosity satisfy
\begin{equation}
 \lambda^3\Lambda\ge\Lambda_{\rm en}b^2,\qquad
              \nu b+\frac{\nu\lambda^3\Lambda}{b}\le1.
 \label{en:viscous-cap}
\end{equation}
\end{enumerate}
Then the modified energy \eqref{en:modified-energy} satisfies
\begin{equation}
 \partial_s\widetilde E+
       c_0\nu\sum_{j=0}^{22}\omega_j
                            K_R[\Psi_{j+1},\Psi_{j+1}]
 \le \Gamma_0b\widetilde E+
                     C_0b^{-1}\sum_{\ell=0}^{22}\|F_\ell\|_H^2.
 \label{en:energy-inequality}
\end{equation}
The constants $c_0,C_0,\Gamma_0$ depend only on the fixed leading
bounds and weights. They are independent of the terminal scale,
the subsequent viscosity, and the selected nonleading coefficients.
The latter enter only the outer-scale restriction
\eqref{en:profile-error}.
\end{theorem}

\begin{proof}
We derive the energy identities, estimate the inviscid terms, and
then absorb the viscous cross terms into the dissipation.

\emph{Step 1: the energy identities.}
Put $u=\Psi_j$ and $v=\Psi_{j+1}$. For $j\ge1$, write
\[
 G_j=-\lambda^{-4}C_j+I_j+\widehat F_j
                                -\nu T_j-\nu\lambda^3\Lambda B_j.
\]
Differentiate \eqref{en:positive-row-energy}, using
\eqref{en:row-equation}. The cross form $Q_R[u,v]$ cancels, and
the remaining identity is
\begin{align}
 \partial_sH_j={}&-\frac12b\|v\|_H^2-6bQ_R[u,u]
                             +\frac12Q_{R,s}[u,u]
 \nonumber\\
 &+(2b_s+\Lambda_{\rm en}b^2)(u,v)_H
       +(13+\Lambda_{\rm en})(bb_s-4b^3)\|u\|_H^2
 \nonumber\\
 &+G_j[v+2bu]-\nu\{K_R[v,v]+2bK_R[u,v]\}.
 \label{en:positive-energy-identity}
\end{align}
Add $\partial_s\mathfrak A_j$ and use
\eqref{en:primitive-cancellation}. The pairing with $\Psi_{j+1}$
cancels, leaving $\mathfrak R_j-7b\mathfrak A_j$.

At order zero, \eqref{en:relative-potential-derivative} instead gives
\begin{align}
 \partial_s\widetilde H_0={}&-\frac12b\|v\|_H^2-5bV_0
                              -2b(F_0^{\rm nl},u)_H+N_0
 \nonumber\\
 &+(2b_s+\Lambda_{\rm en}b^2)(u,v)_H
       +(13+\Lambda_{\rm en})(bb_s-4b^3)\|u\|_H^2
 \nonumber\\
 &+(\widehat F_0-\nu T_0)[v+2bu]
                       -\nu\{K_R[v,v]+2bK_R[u,v]\}.
 \label{en:zero-energy-identity}
\end{align}
\emph{Step 2: the inviscid terms and the source.}
The current form bounds, \eqref{en:primitive-derivative-bound},
\eqref{en:zero-remainder-bound}, and \eqref{en:scalar-bounds} bound
all inviscid terms after the cancellation by $C_0bE$. For example,
\[
 b^2|(u,v)_H|\le C_0bZ(u)\|v\|_H,\qquad
 |bb_s-4b^3|\|u\|_H^2\le C_0bZ(u)^2.
\]
The second-order term \eqref{en:same-row-second-order} has the
same bound. For the forcing, Cauchy--Schwarz gives
\[
 \begin{aligned}
 \sum_{j=0}^{22}\omega_j|\widehat F_j[\Psi_{j+1}+2b\Psi_j]|
 &\le C_0\left(\sum_{j=0}^{22}\|\widehat F_j\|_H^2\right)^{1/2}E^{1/2}\\
 &\le C_0bE+C_0b^{-1}\sum_{j=0}^{22}\|\widehat F_j\|_H^2.
 \end{aligned}
\]
The finite scalar coefficients in $\widehat F_j$ are bounded,
so
\begin{equation}
 \sum_{j=0}^{22}\omega_j
       |\widehat F_j[\Psi_{j+1}+2b\Psi_j]|
       \le C_0bE+C_0b^{-1}\sum_{\ell=0}^{22}\|F_\ell\|_H^2.
 \label{en:strong-forcing-bound}
\end{equation}

\emph{Step 3: the viscous terms.}
The shift assumption and \eqref{en:current-coercivity} imply
\begin{equation}
 \begin{split}
 K_R[f,f]&\ge c_0Z(f)^2,\\
 \lambda^3\Lambda\|f\|_H^2&\le C_0K_R[f,f],\\
 K_R[f,f]&\le C_0Z(f)^2+\lambda^3\Lambda\|f\|_H^2.
 \end{split}
 \label{en:shifted-form-bounds}
\end{equation}
For the middle inequality, use
$\lambda^3\Lambda\|f\|_H^2=K_R[f,f]-Q_R[f,f]$ and the first one.
Its constant therefore contains no factor of the terminal shift.
The Cauchy--Schwarz inequality for the positive form $K_R$ gives
for every prescribed $\delta>0$,
\begin{align}
 2b\nu|K_R[u,v]|
   &\le\delta\nu K_R[v,v]
         +C_\delta\nu(b^2+\lambda^3\Lambda)Z(u)^2,
 \label{en:viscous-cross-young}\\
 \nu|T_j[v+2bu]|
   &\le\delta\nu K_R[v,v]
         +C_\delta\nu b^2\left(\sum_{i=0}^jZ(\Psi_i)\right)^2,
 \label{en:viscous-remainder-young}\\
 \nu\lambda^3\Lambda|(B_j,v+2bu)_H|
   &\le\delta\nu K_R[v,v]
         +C_\delta\nu\lambda^3\Lambda
                         \left(\sum_{i=0}^jZ(\Psi_i)\right)^2.
 \label{en:shift-young}
\end{align}
The second line uses the full bound for $T_j$, including $T_0$.
The third uses the mass part of \eqref{en:shifted-form-bounds};
it is identically zero when $\Lambda=0$.

Choose $0<\delta\le1/6$. At each derivative order the three terms
to be absorbed have sum at most
$3\delta\nu K_R[v,v]\le\nu K_R[v,v]/2$.
Moreover,
\[
 \nu(b^2+\lambda^3\Lambda)\le b,\qquad
 \sum_{j=0}^{22}\omega_j
          \left(\sum_{i=0}^jZ(\Psi_i)\right)^2\le C_0E.
\]
Thus the remaining terms in
\eqref{en:viscous-cross-young}--\eqref{en:shift-young} have sum
at most $C_0bE$. Summing
\eqref{en:positive-energy-identity} and
\eqref{en:zero-energy-identity}, together with the derivatives of
the commutator corrections, proves
\[
 \partial_s\widetilde E+
      c_0\nu\sum_{j=0}^{22}\omega_jK_R[\Psi_{j+1},\Psi_{j+1}]
 \le C_0bE+C_0b^{-1}\sum_{\ell=0}^{22}\|F_\ell\|_H^2.
\]
By \eqref{en:energy-equivalence}, the energy on the right is
bounded by a constant times the modified energy. This proves
\eqref{en:energy-inequality}.

All differentiations are justified by the finite class stated at
the beginning of the section. The differentiated energy pairs orders 22 and 23 in the form domain,
and orders 23 and 24 in $H$. The dissipation uses order 23 in
the pressure form domain. The constitutive derivatives through order 26
come from the potential variations. They require neither additional time
derivatives nor a first trace at order 25; the finite compatibility
conditions suffice.

\end{proof}

On a fixed positive-scale interval, the common bounds of
Section~\ref{ct:section} give $\nu\le C_a\kappa$. The dissipation
estimate and \eqref{en:shifted-form-bounds} imply
\[
 \int\nu K_R[\Psi_{23},\Psi_{23}]\dd s\le C_a,\qquad
 \kappa\int Z(\Psi_{23})^2\dd s\le C_a.
\]
For a fixed square-integrable form test $f$, Cauchy--Schwarz yields
\[
 \begin{aligned}
 \left|\int\nu K_R[\Psi_{23},f]\dd s\right|
 &\le\left(\int\nu K_R[\Psi_{23},\Psi_{23}]\dd s\right)^{1/2}
       \left(\int\nu K_R[f,f]\dd s\right)^{1/2}\\
 &\le C_{a,f}\sqrt\kappa.
 \end{aligned}
\]
The lower-order terms have size $O_a(\kappa)$ by
\eqref{en:viscous-form-bound}. Section~\ref{exlim:section}
identifies the force in the inviscid limit and proves its spatial
regularity.

\begin{corollary}
\label{en:base-bound}
For the same covering atlas,
\begin{equation}
 B:=\sum_{j=0}^{23}\|\Psi_j\|_{X^0}
                    +\sum_{j=0}^{22}\|\Psi_j\|_{X^1}
          \le C_0b^{-1}E^{1/2}
          \le C_0b^{-1}\widetilde E^{1/2}.
 \label{en:energy-to-base}
\end{equation}
\end{corollary}

\begin{proof}
The mass measures give $\|f\|_{X^0}\le C_0\|f\|_H$.
On either edge, Lemma~\ref{en:strain-control} and $a_n\ge c_0$ give
\[
 \int|f_{y_n}|^2\dd x
 \le C_0\lambda I_A[f,f]
 \le C_0\lambda V(f)^2\le C_0V(f)^2,
\]
since $\lambda\le1$. The cutoff derivatives are bounded in their
chart coordinates and contribute at most $C_0\|f\|_H^2$.
At the center the strain form controls the complete Cartesian
derivative of the radial vector; in the fixed interior use the
one-dimensional derivative estimate. Summing gives
\[
 \|f\|_{X^1}\le C_0V(f)\le C_0b^{-1}Z(f).
\]
Consequently, by Cauchy--Schwarz in the finite sums,
\[
 \begin{aligned}
 B
 &\le C_0b^{-1}
    \left(\sum_{j=0}^{22}
       \{\|\Psi_{j+1}\|_H^2+Z(\Psi_j)^2\}\right)^{1/2}\\
 &\le C_0b^{-1}E^{1/2}
 \le C_0b^{-1}\widetilde E^{1/2}.
 \end{aligned}
\]
This proves \eqref{en:energy-to-base}. The factor
$b^{-1}\le\beta^{-1}$ is used after parameter selection and does
not enter the energy growth constant $\Gamma_0$.
\end{proof}

\begin{corollary}
\label{en:energy-recovery}
Under the respective hypotheses of Proposition~\ref{rec:recovery}
and Corollary~\ref{en:base-bound},
\begin{equation}
 B(s)\le C_0b(s)^{-1}E(s)^{1/2}
             \le C_0b(s)^{-1}\widetilde E(s)^{1/2}.
 \label{rec:energy-base}
\end{equation}
For the inviscid solution this gives
\begin{equation}
 \mathcal X_{23}(s)\le C_{\mathrm{rec},\beta}
      \{\widetilde E(s)^{1/2}+\mathcal R_{23}(s)\}.
 \label{rec:inviscid-energy-recovery}
\end{equation}
\end{corollary}

\begin{proof}
The first estimate is Corollary~\ref{en:base-bound}. Substitution
in \eqref{rec:inviscid-recovery}, with $b^{-1}\le\beta^{-1}$, proves
the second. The resulting constant depends on the chosen $\beta$;
the energy growth constant $\Gamma_0$ was fixed before this step.
Section~\ref{ct:section} combines these bounds with the integral
estimate for the viscous equation.
\end{proof}

\begin{corollary}
\label{en:prescribed-source}
For the normalized profile residual and viscosity source,
\begin{equation}
 F=-\mathcal R_u-\nu\lambda^{-1}\partial_s\mathcal M[A]
                           -\nu\lambda^{-1}\Lambda A_s,
 \qquad
 \mathcal R_u=\lambda^{-1}(A_{\tau\tau}+\mathcal M[A]),
 \label{en:actual-prescribed-source}
\end{equation}
Proposition~\ref{en:input-material-source} gives
\begin{equation}
 \sum_{\ell=0}^{22}\|F_\ell\|_H
 \le C_{J,K,\beta}
       \{\lambda^{(K-7)/2}(1+|\log\lambda|)^{p_{\log}}
                 +\kappa\lambda^{-9/2}
                 +\kappa\Lambda\lambda^{-3/2}\}.
 \label{en:source-twenty-two}
\end{equation}
Consequently the last term of \eqref{en:energy-inequality} is at most
\begin{equation}
 C_{J,K,\beta}b^{-1}
       \{\lambda^{K-7}(1+|\log\lambda|)^{2p_{\log}}
                 +\kappa^2\lambda^{-9}
                 +\kappa^2\Lambda^2\lambda^{-3}\}.
 \label{en:explicit-energy-source}
\end{equation}
\end{corollary}

\begin{proof}
Insert \eqref{en:source-twenty-two} in
\eqref{en:energy-inequality} and apply Cauchy--Schwarz to the finite sum.
The factor $\lambda^{-1}$ in the definition of $\mathcal R_u$ is
part of the normalization. The source constant is allowed to depend
on the selected profile and is not included in $\Gamma_0$.
\end{proof}

The energy growth constant $\Gamma_0$ is now fixed. The parameter selection
in Section~\ref{ct:section} chooses $J$ above its growth threshold,
then the leading profile and a common endpoint for the scale interval
on which \eqref{en:profile-error} and the source errors are small. These later
choices do not change $\Gamma_0$. Substituting the bound of
Corollary~\ref{en:base-bound} into the recovery estimate of
Section~\ref{rec:chapter} gives the low spatial bound used to
continue the solution throughout that interval.

The two source requirements remain distinct. The energy estimate uses the
pure material estimate \eqref{mat:material-source-export} through time order 22;
the mixed recovery source \eqref{mat:strong-source-export} covers only
$j+m\le21$. The form bound \eqref{en:viscous-form-bound} controls the
energy commutator, while strong spatial recovery uses its own strong
coefficient estimates. Section~\ref{exlim:section} uses the retained
dissipation to remove the highest viscous derivative in the form domain,
and then establishes the strong spatial domains for the inviscid limit.
\section{A common interval for the regularized solutions}
\label{ct:section}

The local solution starts at $\lambda=a>0$. We extend it to
$\lambda=A_0$, where $A_0$ is independent of the terminal scale
and the permitted viscosity; the mass and matched profile stay fixed.

On a bootstrap interval, the low spatial bound supplies the small
strains required by the energy estimate. Since
\[
 p=2J+1>\Gamma_0,
\]
the initial error and source decay faster than the growth allowed by
the energy inequality. Integration improves the energy bound, and
spatial recovery then improves the low spatial bound. These estimates
also control the local continuation norm and preserve positivity of
the normalized radius and Jacobian. The high local norms and the
local lifespan may depend on $(a,\kappa)$.

Section~\ref{exlim:section} removes viscosity using the low spatial
bounds and retained dissipation. The inviscid equation then gives
the full spatial estimates before $a\downarrow0$.

In the finite profile bounds below, $p_{\log}$ is a nonnegative
integer, enlarged when finitely many logarithmic degrees are combined.
Here $a>0$ denotes the terminal scale, $s_a$ its similarity time, and
$\Lambda_a$ the constant regularizing shift. We keep the energy
$\widetilde E$, the weights $\omega_j$, and the forms $K_R$ of
Theorem~\ref{en:differential-energy}. Write
\begin{align}
 \mathcal D&=\sum_{j=0}^{22}\omega_j
                         K_R[\Psi_{j+1},\Psi_{j+1}],\qquad
 \nu=\kappa\lambda^{-3/2},
 \label{ct:dissipation}\\
 \mathcal X_8&=\sum_{q=0}^{8}\sum_{j=0}^{23-q}
                                      \|\Psi_j\|_{X^q},\qquad
 G_8^2=\sum_{q=0}^{8}\sum_{j=0}^{23-q}
                                      \|\Psi_j\|_{X^q}^2.
 \label{ct:low-graphs}
\end{align}
There are 180 terms in these sums; hence
$G_8\le\mathcal X_8\le\sqrt{180}\,G_8$.
The source norms $\mathcal S_m$ below are the strong source norms
in the same covering atlas, not form-dual norms.

\subsection{Recovery and initial bounds}

We use the local construction of Proposition~\ref{loc:prepared-local-family}
and the viscous low-order estimate \eqref{rec:viscous-recovery}.
Set $p=2J+1$. For any nonnegative scalar function $H(s)$, write
\[
 H^\sharp(s)=\sup_{s_a\le\sigma\le s}
                         \lambda(\sigma)^{-p/2}H(\sigma).
\]
Define the mixed source norm by
\[
 \mathcal R_8=\sum_{q=2}^{8}\sum_{j=0}^{23-q}
                                      \mathcal S_{q-2}(F_j).
\]
On the interval where the strains are small, the recovery estimate reads
\begin{equation}
 \mathcal X_8^\sharp(s)
 \le C_{\rm rec}\left\{
      a^{-p/2}\mathcal X_8(s_a)+B^\sharp(s)
                                  +\mathcal R_8^\sharp(s)\right\}.
 \label{ct:viscous-recovery}
\end{equation}
The first-order quantity $B$ is defined in Corollary~\ref{en:base-bound}.
Here $\Gamma_q^{\rm cov}$ denotes the covering transport constant
from Lemma~\ref{rec:transport}, including when $q=0$; it is distinct
from the energy growth constant $\Gamma_0$.
The transport constants entering this estimate depend on
the fixed cutoffs and $q$, and are chosen before $\beta$.
The scalar integral kernel has exponent
$-\int[\nu^{-1}+(4+p/2-\Gamma_q^{\rm cov})b]$; thus its use requires
$4+p/2\ge\Gamma_q^{\rm cov}$. Only order-eight spatial recovery is used here.

The initial bounds follow from Proposition~\ref{loc:prepared-energy}.
The Euler-compatible initial data satisfy
\begin{equation}
 |\widetilde E(s_a)|
       \le C_\beta a^{2J+1644}(1+|\log a|)^{p_{\log}},
 \qquad
 \mathcal X_8(s_a)
       \le C_\beta a^{J+1685/2}(1+|\log a|)^{p_{\log}}.
 \label{ct:preparation-margin}
\end{equation}
The family of viscous compatible data is continuous at that datum in all
these graph norms and in the value of the full modified energy.
The energy bound in \eqref{ct:preparation-margin} includes the
relative potential and all commutator corrections. The same bounds
hold for sufficiently small positive viscosity at each fixed terminal
scale.

\begin{proposition}
\label{exlim:upstream-tube}
\label{ct:proposition}
Fix the leading parameter $\beta>0$, its finite matched profile,
and integers $J\ge3,K$ satisfying \eqref{ct:orders}. Choose
$d_0,\epsilon,c,A_0$ as in
Subsection~\ref{ct:parameter-admissibility-section}.

Let $0<a\le A_0$ and $0<\kappa\le\kappa_{\rm tube}(a)$, with
$\kappa_{\rm tube}(a)$ and the prepared data chosen as in that
subsection, so that \eqref{ct:initial-small} holds. The prepared
solution of
\begin{equation}
 R_{\tau\tau}+\mathcal M[R]
                +\kappa\{D\mathcal M[R]+\Lambda_a\}R_\tau=0,
 \qquad
 \Lambda_a=\Lambda_{\rm en}b(a)^2a^{-3},
 \label{ct:regularized-equation}
\end{equation}
has the following properties.
\begin{enumerate}[label=(\roman*),leftmargin=2em]
\item It exists throughout $a\le\lambda\le A_0$, with the stated
positivity bounds and finite local regularity.
\item For the selected $d_0>0$, its energy, spatial norm, and
dissipation satisfy
\begin{equation}
 0\le\widetilde E(s)\le\frac{d_0^2}{4}\lambda(s)^p,
 \qquad G_8(s)\le\frac12\lambda(s)^{p/2},
 \qquad
 \int_{s_a}^{s}\nu\mathcal D\dd\sigma
                              \le C d_0^2\lambda(s)^p.
 \label{ct:conclusion}
\end{equation}
\end{enumerate}
The constants may depend on the selected parameters, but are independent
of $a$ and the permitted viscosity. The high local norms are finite for
each fixed $a,\kappa$; no uniform bound for them is asserted.
\end{proposition}

\begin{proof}
Fix constants satisfying the stated restrictions; their choice is
proved in Lemma~\ref{ct:parameter-admissibility}.

\emph{Energy and spatial bounds.}
Consider the connected interval on which $\lambda\le A_0$ and
\[
 \widetilde E\le d_0^2\lambda^p,\qquad G_8\le\lambda^{p/2}.
\]
The initial bounds are strict. Moreover,
$\eta=\lambda^{-2}G_8\le A_0^{J-3/2}$, so the strain restriction
and \eqref{ct:coefficient-small} give the hypotheses of the energy
theorem. Multiply \eqref{en:energy-inequality} by
$\lambda^{-\Gamma_0}$ and use
$\dd s=\dd\lambda/(b\lambda)$ and \eqref{ct:source-small}:
\begin{align}
 \widetilde E(\lambda)
 &\le\widetilde E(a)(\lambda/a)^{\Gamma_0}
    +C_E\epsilon^2\lambda^{\Gamma_0}
       \int_a^\lambda\frac{u^{p-\Gamma_0-1}}{b(u)^2}\dd u\notag\\*
 &\le\left\{\frac{d_0^2}{8}
       +\frac{C_E\epsilon^2}{\beta^2(p-\Gamma_0)}\right\}\lambda^p
 \le\frac{d_0^2}{4}\lambda^p.
 \label{ct:energy-improvement}
\end{align}
Proposition~\ref{en:energy-coercivity} also gives
$\widetilde E\ge0$. The strong $H$ source bound is used in the
integral above. From the weaker bootstrap energy bound,
\begin{align}
 \mathcal X_8^\sharp(s)
 &\le C_{\rm rec}\{
    a^{-p/2}\mathcal X_8(s_a)+B^\sharp(s)+\mathcal R_8^\sharp(s)\}
 \le\frac18+\frac18+\frac18=\frac38.
 \label{ct:graph-improvement}
\end{align}
Since $G_8\le\mathcal X_8$, this improves the spatial bound.
Keeping the dissipation in the same integration gives
\[
 c_0\int_{s_a}^{s}
   \left(\frac{\lambda(s)}{\lambda(\sigma)}\right)^{\Gamma_0}
                         \nu(\sigma)\mathcal D(\sigma)\dd\sigma
                              \le\frac{d_0^2}{4}\lambda(s)^p.
\]
Its weight is at least one, giving the last estimate in
\eqref{ct:conclusion}.

\emph{The continuation norm.}
Fix $a>0$ and $0<\kappa\le\kappa_{\rm tube}(a)$. On $[a,A_0]$,
the transitions from the moving atlas to an enlarged fixed
center--interior--vacuum atlas have bounded derivatives, with
constants depending on $a$. On a region where a prescribed cutoff
is bounded below, division by that cutoff bounds the fixed local
norms by $X^8$.

At the vacuum, set $\xi=(M-x)^{1/5}$ on a fixed collar $[0,c]$.
The $B_5^8$ norm bounds
$\|f^{(r)}\|_{L^2(\xi^4\dd\xi)}$, $r\le8$, and the interior traces
through order seven. Taylor's formula for $f^{(5)}$ at an interior
point $c_0$ has remainder
\begin{align}
 \frac12\left|\int_\xi^{c_0}(t-\xi)^2f^{(8)}(t)\dd t\right|
 &\le C\|f^{(8)}\|_{L^2(t^4\dd t)}
       \left(\int_\xi^{c_0}(t-\xi)^4t^{-4}\dd t\right)^{1/2}
 \le C\|f\|_{B_5^8}.
 \label{ct:endpoint-embedding}
\end{align}
Thus
\[
 \|f\|_{C^5([0,c])}\le C\|f\|_{B_5^8},
 \qquad
 \|f\|_{C^{2,\alpha}([0,c])}\le C_\alpha\|f\|_{B_5^8},
 \quad 0<\alpha<1.
\]
For compatible $f$, the condition $f_\xi(0)=0$ gives
\begin{equation}
 \frac{f_\xi(\xi)}{\xi}=\int_0^1f_{\xi\xi}(t\xi)\dd t,
 \qquad
 \frac{f(\xi)-f(0)}{\xi^2}
                     =\int_0^1(1-t)f_{\xi\xi}(t\xi)\dd t.
 \label{ct:quotients}
\end{equation}
These identities bound both source quotients in $C^{0,\alpha}$.
The finite local solution already lies in the required little
spaces; this estimate bounds its norm without imposing a first
trace on a rough limit.

At the center we use the full equivariant Cartesian vector.
After a cutoff in an enlarged ball, $H^8(\mathbb R^3)$ controls
three uniform derivatives, since
\[
 \int_{\mathbb R^3}|\zeta|^6(1+|\zeta|^2)^{-8}\dd\zeta<\infty.
\]
It therefore controls $C^{2,\alpha}$ for $0<\alpha<1$.
The fixed interior uses the one-dimensional embedding. Apply these
bounds to $\Psi_0,\Psi_1$ in
\[
 R=A+\lambda^4\Psi_0,\qquad
 R_\tau=\lambda^{-3/2}(A_s+\lambda^4\Psi_1)
\]
to obtain
\begin{equation}
 \sup_{a\le\lambda\le A_0}
       \{\|R\|_{\mathcal D^0}+\|R_\tau\|_{\mathcal D^0}\}
                           \le L_{a,A_0,\beta,J}<\infty.
 \label{ct:base-ball}
\end{equation}
The strain bounds preserve positive lower bounds for the normalized
radius and Jacobian. These and \eqref{ct:base-ball} are the
fixed-viscosity continuation hypotheses; a bound for all 48 initial
derivatives is not required.

\emph{Continuation.}
Let $\tau_+$ be the right endpoint of the maximal connected
bootstrap interval, and suppose $\lambda(\tau_+)<A_0$.
The preceding bounds give a lifespan $T_{\rm loc}>0$ at every
earlier slice, depending on $a,\kappa,\Lambda_a$ but not on the
slice. Choose
\[
 \max\{\tau_a,\tau_+-T_{\rm loc}/2\}<\tau_0<\tau_+.
\]
The compatible data at $\tau_0$ extend the solution through
$\tau_0+T_{\rm loc}>\tau_+$; fixed-viscosity uniqueness identifies
the extension on the overlap. There,
\[
 \lambda^{-p}\widetilde E\le d_0^2/4<d_0^2,\qquad
 \lambda^{-p/2}G_8\le3/8<1
\]
by \eqref{ct:energy-improvement} and \eqref{ct:graph-improvement}.
Continuity preserves the bootstrap bounds beyond $\tau_+$, a
contradiction. The same argument supplies the endpoint
$\lambda=A_0$ and proves \eqref{ct:conclusion}.
\end{proof}

\subsection{Choice of parameters}
\label{ct:parameter-admissibility-section}

We now choose the parameters required by
Proposition~\ref{ct:proposition}. Leading constants determine $J,K$
on a common parameter interval for the leading profiles. After selecting $\beta>0$, a decrease
of the upper bound for the scale makes the remaining profile, preparation, and
recovery errors small.

\begin{lemma}
\label{ct:parameter-admissibility}
Fix $e\ge0$, the covering cutoffs, and the finite energy weights.
There are $\beta_*>0$ and integers $J\ge3$, $K=2J+3512$, chosen
before the value of $\beta$ is fixed, with the following properties
for every $0<\beta\le\beta_*$.
\begin{enumerate}[label=(\roman*),leftmargin=2em]
\item There is $A_0=A_0(\beta,e)>0$ on which the matched profile,
exact-mass preparation, recovery estimates, energy estimates, and
Proposition~\ref{ct:proposition} apply. In particular, the
estimates uniform in the terminal scale hold for $0<\lambda\le A_0$.
\item For each $0<a\le A_0$, the permitted viscosities contain
$(0,\kappa_{\rm tube}(a)]$. The same $A_0$ works for every such $a$ and viscosity. Thus both the
viscosity limit at fixed $a$ and the subsequent subsequence limit as
$a\downarrow0$ are taken at the same mass $M_\beta$.
\end{enumerate}
No positive lower bound for $A_0$ or $\kappa_{\rm tube}(a)$ is
asserted uniformly as $\beta\downarrow0$.
\end{lemma}

\begin{proof}
\emph{Leading parameters and orders.}
The smooth leading profiles in
Subsubsections~\ref{mat:seed-part-1}--\ref{mat:seed-part-5}
exist on a compact parameter interval independent of the subsequent
derivative order. The pullbacks after Definition~\ref{at:atlas}
use the measures $4\pi w_\beta^3z^2\dd z$, with common endpoint
powers, and require no unweighted endpoint inverse.

Let $C_{\rm var}$ be a common constant in \eqref{en:variance}.
For $h=zg\perp_H\varphi_\beta$,
\[
 \frac{16\pi}{3}\int_0^1w_\beta^4z^4|g'|^2\dd z
 \ge\sigma_*\|h\|_H^2,\qquad
 \sigma_*=\frac4{3C_{\rm var}}>0.
\]
The variance integrand is $O(z)$ at the center and $O(1)$ at
vacuum. Choose
\[
 0<\beta_*\le
 \min\{\beta_{\rm match},\beta_{\rm env},b_{\max}/\sqrt2,
                          (\sigma_*/2)^{1/2}\}.
\]
Lemma~\ref{en:complement-coercivity} then requires only a later
reduction of the scale bound. The energy constants $\Gamma_0$ and
$\Lambda_{\rm en}$ use profile time derivatives through order $23$,
spatial derivatives through $24$, and constitutive derivatives
through $26$. The error
$C_{J,K,\beta}\lambda^{1/2}(1+|\log\lambda|)^{p_{\log}}$
is kept separate. The constants $\Gamma_q^{\rm cov}$, $q\le8$,
depend only on the fixed cover.

Choose $J\ge3$ and $K$ so that
\begin{equation}
 \begin{gathered}
 p=2J+1>\Gamma_0,\qquad
 4+p/2\ge\max_{q\le8}\Gamma_q^{\rm cov},\\
 J^\sharp:=J+1713,\qquad K=2J^\sharp+86=2J+3512.
 \end{gathered}
 \label{ct:orders}
\end{equation}
The endpoint order in \eqref{mat:endpoint-20} is this $J^\sharp$:
\[
 J^\sharp+\frac{55}{2}=J+\frac{3481}{2},\qquad
 K-1-\left(J^\sharp+\frac{55}{2}\right)=J^\sharp+\frac{115}{2}>0.
\]
The material source, mixed source, initial energy and initial
spatial norm have excess powers
\begin{equation}
 \begin{aligned}
 (K-7)/2-p/2&=1752,& (K-5)/2-p/2&=1753,\\
 (2J+1644)-p&=1643,& (J+1685/2)-p/2&=842.
 \end{aligned}
 \label{ct:exact-margins}
\end{equation}
Smooth dependence bounds all higher leading-profile derivatives
needed at this $K$ on the same interval. Their constants, and the
recovery and compatible-data inverses, do not change $\Gamma_0$
or \eqref{ct:orders}.

\emph{Constants for a fixed positive parameter.}
Fix $0<\beta\le\beta_*$. For $2\le n\le K+1$, put
$c_n=(2n^2+n-3)\beta^2/2>0$, and let $u_n$ be the center-regular
homogeneous solution with $u_n(0)=1$. By
\eqref{mat:connection-24},
\[
 w^4z^4u_n'(z)=\frac34c_n\int_0^z w^3t^4u_n(t)\dd t.
\]
While $u_n>0$ the right side is positive, excluding a first zero.
The fundamental pair
\eqref{mat:connection-18}--\eqref{mat:connection-19} gives
$u_n=O(q^{-3})$, $q=1-z$. Thus $w^3z^4u_n$ is $O(z^4)$ at
the center and $O(1)$ at vacuum. For $s_n=[U_n^s]u_n$, the
limiting flux identity is
\[
 3\varkappa^4s_n=\frac34c_n\int_0^1w^3z^4u_n\dd z>0.
\]
For $n=2,3$, the coefficient $s_n>0$ determines the amplitude
which removes the singular branch. For $n\ge4$, the two matching
diagonals are
\[
 c_*^{-3}s_n>0\quad\text{and}\quad1.
\]
The recursions \eqref{mat:core-log-recursion} and
\eqref{mat:connection-63} use these diagonals at each logarithmic
degree, through the auxiliary coefficient of degree $K+1$.
The fixed index resonances require only finite logarithmic
primitives and exclude no $\beta>0$. We use no inverse at
$\beta=0$, where $s_n$ vanishes.

The preparation map has positive mass derivative by
Lemma~\ref{loc:finite-map}; its remaining diagonal coefficients
$d_n,e_n$ are positive multiples of powers of $p_*>0$.
Thus \eqref{loc:finite-map-bounds} controls its triangular inverse.
All matching, preparation, elliptic, and recovery constants are
finite for this fixed $\beta$. Fix them on $0<\lambda\le A_*$;
decreasing the upper endpoint preserves their bounds.

\emph{Smallness conditions.}
For $\alpha>0$, set
\[
 H_{\alpha,p_{\log}}(A)=\sup_{0<\lambda\le A}
               \lambda^\alpha(1+|\log\lambda|)^{p_{\log}};
 \qquad H_{\alpha,p_{\log}}(A)\longrightarrow0
 \quad(A\downarrow0).
\]
Choose $d_0>0$ and then $\epsilon>0$ so that
\begin{equation}
 C_{\rm rec}C_B\beta^{-1}d_0\le\frac18,
 \qquad C_{\rm rec}\epsilon\le\frac18,
 \qquad
 \frac{C_E\epsilon^2}{\beta^2(p-\Gamma_0)}\le\frac{d_0^2}{8}.
 \label{ct:small-constants}
\end{equation}
Here $B\le C_B\beta^{-1}\widetilde E^{1/2}$ is
Corollary~\ref{en:base-bound}; $C_E$ contains only the finite
strong-source comparison factors. Fix $C_F$ covering
\eqref{en:source-twenty-two} and \eqref{rec:strong-source-budget},
and choose $c>0$ with $2C_Fc\le\epsilon/2$.
After these choices, decrease $A_0$ to satisfy:
\begin{enumerate}[label=(\roman*),leftmargin=2em]
\item \emph{Scale and covering.}
For \eqref{rec:clock} and Lemma~\ref{rec:transport},
\[
 A_0\le1,\qquad eA_0\le\beta^2,\qquad 2A_0^4\le d_0^{\rm edge}.
\]
\item \emph{Profile and form comparison.}
For \eqref{en:profile-error} and the positive lower bounds in
Proposition~\ref{mat:matched-family},
\[
 C_{J,K,\beta}H_{1/2,p_{\log}}(A_0)\le c_0\beta^2.
\]
\item \emph{Prepared-data Newton iteration.}
The radius $CH_{K-1711/2,p_{\log}}(A_0)$ in
\eqref{loc:finite-map-bounds} lies inside the fixed positive
parameter ball, and
\[
 C\{H_{K-3421/2,p_{\log}}(A_0)+H_{K-1711/2,p_{\log}}(A_0)\}\le1/2.
\]
\item \emph{Prescribed source.}
For pure material derivatives through order $22$ and mixed
derivatives with $j+m\le21$,
\[
 C_F\{H_{1752,p_{\log}}(A_0)+H_{1753,p_{\log}}(A_0)\}\le\epsilon/2.
\]
\item \emph{Initial energy and spatial norm.}
For \eqref{ct:preparation-margin},
\[
 C_PH_{1643,p_{\log}}(A_0)\le d_0^2/16,\qquad
 C_PH_{842,p_{\log}}(A_0)\le(16C_{\rm rec})^{-1}.
\]
These retain a factor of one half for continuity from $\kappa=0$.
\item \emph{Strains and inviscid absorption.}
For \eqref{rec:inviscid-substitution}, $A_0^{J-3/2}$ is below
half each selected positivity and energy threshold, including
$\varepsilon_0\beta^3$, and below both $(2C_*)^{-1}$ and
$\varepsilon_{\rm rec}/4$.
\item \emph{The strong operator estimate.}
For Proposition~\ref{at:current-inverse},
\[
 A_0^{J-1/2}\le\varepsilon_{\rm inv}.
\]
\item \emph{Viscous recovery and energy.}
For \eqref{rec:small-coefficients} and \eqref{en:viscous-cap},
\[
 \begin{gathered}
 cA_0^{J+7/2}\le1,\qquad
 c(b_{\max}+1)A_0^{J+7/2}\le\varepsilon_{\rm rec}/2,\\
 c(b_{\max}+\beta^{-1})A_0^{J+7/2}\le1.
 \end{gathered}
\]
\end{enumerate}
Enlarge the finitely many constants and logarithmic degrees before
choosing $A_0$. Every power above is positive by \eqref{ct:orders}
and \eqref{ct:exact-margins}. In item (iii), the Newton radius
has order $\lambda^{K-1711/2}(1+|\log\lambda|)^{p_{\log}}$;
the inverse and nonlinear derivative in
\eqref{loc:finite-map-bounds} give the other exponent
$K-3421/2$. Also,
\[
 G_8\le\lambda^{p/2}
 \quad\Longrightarrow\quad
 \eta\le A_0^{J-3/2},\qquad
 \lambda^{-1}\|\Psi_0\|_{X^8}\le A_0^{J-1/2}.
\]
Thus a single positive $A_0$ meets all restrictions. The order of
choice is $d_0,\epsilon,c,A_0$, followed by $a,\kappa$.

\emph{Terminal scale and viscosity.}
For any $0<a\le A_0$, set
$\Lambda_a=\Lambda_{\rm en}b(a)^2a^{-3}$. Since
\[
 \frac{\mathrm d}{\mathrm d\lambda}
       \{b(\lambda)^2\lambda^{-3}\}
 =-3\beta^2\lambda^{-4}-2e\lambda^{-3}<0,
\]
we have $\lambda^3\Lambda_a\ge\Lambda_{\rm en}b(\lambda)^2$
for $a\le\lambda\le A_0$. Impose
\begin{equation}
 \kappa\le c a^{J+5},\qquad
 \kappa\Lambda_a\le c a^{J+2}.
 \label{ct:viscosity-caps}
\end{equation}
or equivalently
\[
 0<\kappa\le\kappa_{\rm cap}(a)
 :=ca^{J+5}\min\{1,(\Lambda_{\rm en}b(a)^2)^{-1}\}.
\]
Intersect this bound with the positive compatible-data interval
at $(a,\Lambda_a)$ and the continuity interval on which
\begin{equation}
 |\widetilde E(s_a)|\le\frac{d_0^2}{8}a^p,
 \qquad
 \mathcal X_8(s_a)\le\frac{1}{8C_{\rm rec}}a^{p/2}.
 \label{ct:initial-small}
\end{equation}
Item (v) and Proposition~\ref{loc:prepared-energy} give strict
margins at zero viscosity. The intersection therefore contains
$(0,\kappa_{\rm tube}(a)]$ for some $\kappa_{\rm tube}(a)>0$.
Its size need not be uniform in $a$.

On $[a,A_0]$ the caps give
\[
 \begin{aligned}
 \nu&=\kappa\lambda^{-3/2}
       \le ca^{J+7/2}\le cA_0^{J+7/2},\\
 \nu\lambda^3\Lambda_a
       &=\kappa\Lambda_a\lambda^{3/2}
       \le ca^{J+2}\lambda^{3/2}\le cA_0^{J+7/2}.
 \end{aligned}
\]
Using $\beta\le b\le b_{\max}$, we obtain
\begin{equation}
 \nu\le cA_0^{J+7/2},\quad
 \nu b\le cb_{\max}A_0^{J+7/2},\quad
 \frac{\nu\lambda^3\Lambda_a}{b}
                   \le c\beta^{-1}A_0^{J+7/2}.
 \label{ct:coefficient-small}
\end{equation}
Since $\nu\le1$, also $\nu\eta\le\eta$.
Items (vi)--(viii) therefore control the full sum in
\eqref{rec:small-coefficients} and imply \eqref{en:viscous-cap}.
Only the product of the shift and viscosity enters these bounds.

The prescribed source is
\begin{equation}
 F=-\lambda^{-1}(A_{\tau\tau}+\mathcal M[A])
       -\nu\lambda^{-1}\partial_s\mathcal M[A]
       -\nu\lambda^{-1}\Lambda_a A_s,
 \qquad F_j=\lambda^{-4}\partial_s^j(\lambda^4F).
 \label{ct:source}
\end{equation}
Its viscous terms satisfy
\[
 \kappa\lambda^{-9/2}\le c\lambda^{J+1/2},\qquad
 \kappa\Lambda_a\lambda^{-3/2}\le c\lambda^{J+1/2}.
\]
Since $p/2=J+1/2$, item (iv) and $2C_Fc\le\epsilon/2$ yield
\begin{equation}
 \sum_{j=0}^{22}\|F_j\|_H\le\epsilon\lambda^{p/2},
 \qquad \mathcal R_8\le\epsilon\lambda^{p/2}.
 \label{ct:source-small}
\end{equation}
The first estimate uses the pure material bound; the second uses
only the stated mixed-derivative range. This verifies all
hypotheses of Proposition~\ref{ct:proposition}.

\emph{Uniformity.}
The bootstrap constants depend on the selected profile but not on
$a$ or the permitted viscosity. The coordinate comparisons in
\eqref{ct:base-ball} and the local continuation time may depend
on the fixed $a,\kappa,\Lambda_a$.

Compactness and strong spatial regularity at fixed $a$ may also
use $a$-dependent constants. Once the force is identified and the
limiting derivatives lie in the strong operator domains,
\eqref{exlim:recovery-input} uses the original $C_{\rm rec}$ and
the threshold in item (vi); the fixed-$a$ regularity constants
do not enter \eqref{exlim:full-decay}. Remove viscosity first and
then take a nested subsequence as $a\downarrow0$, keeping $A_0$
and $M_\beta$ fixed.
\end{proof}
\section{Removal of viscosity and the terminal limit}
\label{exlim:section}

The mass, profile, and upper endpoint of the scale interval remain fixed
during both limits.
At a fixed terminal scale $a>0$, the low-order bounds and dissipation
of Section~\ref{ct:section} give compactness as the artificial
viscosity tends to zero and identify an Euler--Poisson solution.
The limiting equation gives higher spatial regularity and shows that
the limiting derivatives lie in the closed strong operator domains.
These arguments may use constants depending on $a$. We then apply
the original recovery estimate to obtain high-order bounds uniform
in $a$.

A diagonal subsequence as $a\downarrow0$ gives one solution on all
positive scales. Its norms are essentially bounded at total order $23$
and continuous through total order $22$; no convergence of the full
terminal sequence or uniqueness is asserted. The continuous lower
derivatives give classical endpoint regularity. Energy conservation
and a global comparison with the matched profile then determine the
conserved energy. The same lower-order bounds yield the core and
boundary-layer limits at every time before collapse.

\subsection{Uniform bounds and weighted spaces at fixed mass}
\label{exlim:inputs}

Fix the selected leading profile, mass $M$, parameter $\beta>0$, and $e\geq0$.
The profile $A$, the integers $J\geq3$ and $K=2J+3512$,
and the upper bound for the scale $A_0$ are fixed throughout. Write
\begin{equation}
 \lambda_s=b\lambda,\qquad b^2=\beta^2+e\lambda,\qquad
 \partial_\tau=\lambda^{-3/2}\partial_s,\qquad
 \rho[R]=(4\pi R^2R_x)^{-1}.
 \label{exlim:clock}
\end{equation}
The exact pressure and force are
\begin{equation}
 P(0)=0,\qquad
 P'(\rho)=\frac{4\rho^{2/3}}{3\sqrt{1+\rho^{2/3}}},\qquad
 \mathcal M[R]=4\pi R^2\partial_xP(\rho[R])+\frac{x}{R^2}.
 \label{exlim:force}
\end{equation}
We write $H_R=D\mathcal M[R]$, $L_R=\lambda^3H_R$, and
\[
 Q_R[f,g]:=\lambda^3\langle H_Rf,g\rangle
\]
on compatible form directions.
The potential underlying its closed form is
\[
 \mathcal U[R]=\int_0^M\left\{\mathfrak e(\rho[R])-\frac{x}{R}\right\}\dd x,
\]
with the specific internal energy of \eqref{eq:internal-energy-normalization}.
On compatible test directions its first variation is
$\mathcal M[R]$, with the vanishing boundary flux identified below.

We use the weighted spaces of the fixed-mass cover. At the center
$z_c=x^{1/3}$, let $U_f$ denote the equivariant lift
$f(|X|^3)X/|X|$ of a displacement. Its norm is the Cartesian
$H^q$ norm.
The interior norm is ordinary $H^q$. The two coordinates near the
vacuum are
\[
 y_5=((M-x)\lambda^{-4})^{1/5},\qquad
 y_8=((M-x)\lambda^{-4})^{1/8}.
\]
Their norms are the ordinary derivative sums with measures
$y_5^4\dd y_5$ and $y_8^7\dd y_8$, denoted by $B_5^q$ and
$C_8^q$. If $\chi_c,\chi_i,\chi_5,\chi_8$ are the prescribed
nonnegative power cutoffs, their sum is one and
\begin{equation}
 \|f\|_{X^q(s)}^2
 =\|\chi_cU_f\|_{H^q(B;\mathbb R^3)}^2+\|\chi_if\|_{H^q}^2
  +\lambda^4\|\chi_5f\|_{B_5^q}^2
  +\lambda^4\|\chi_8f\|_{C_8^q}^2.
 \label{exlim:graph}
\end{equation}
All derivatives of the cutoffs in these localized norms are included. The long chart ends
at $C_e\lambda^{-1/2}$; the actual core profile is used beyond the
matching region. At the center we complete smooth equivariant vectors in these norms.
At the vacuum we complete smooth one-sided functions satisfying
$f_{y_5}(0)=0$. These completions define $X^q(s)$;
no even scalar extension at the vacuum is imposed.
The source norm $\mathcal S_m(g)$ is given by
\eqref{exlim:graph} with $g$ in the center and interior terms and
$\lambda g$ in the two edge terms.

For a solution set
\begin{equation}
 R_j=\partial_s^jR,\quad A_j=\partial_s^jA,\quad
 \Psi_j=\lambda^{-4}(R_j-A_j),\quad
 \mathcal X_{23}=\sum_{q=0}^{23}\sum_{j=0}^{23-q}
                              \|\Psi_j\|_{X^q}.
 \label{exlim:rows}
\end{equation}
The lower quantities needed below are
\begin{equation}
 B=\sum_{j=0}^{23}\|\Psi_j\|_{X^0}
       +\sum_{j=0}^{22}\|\Psi_j\|_{X^1},\qquad
 G_8^2=\sum_{q=0}^{8}\sum_{j=0}^{23-q}\|\Psi_j\|_{X^q}^2,
 \qquad \eta=\lambda^{-2}G_8.
 \label{exlim:base}
\end{equation}

We recall the conclusions used below.
For every fixed $0<a<A_0$, Proposition~\ref{exlim:upstream-tube}
supplies fixed-mass solutions on the whole interval
$I_a=\{s:a\leq\lambda(s)\leq A_0\}$, for sufficiently small
positive $\kappa$, of
\begin{equation}
 R^\kappa_{\tau\tau}
 +\kappa(H_{R^\kappa}+\Lambda_a)R^\kappa_\tau
 +\mathcal M[R^\kappa]=0.
 \label{exlim:viscous}
\end{equation}
The shift $\Lambda_a$ is constant in time. The allowed $\kappa$
and the local existence times may depend on $a$. The family has
uniformly positive normalized deformation and vacuum factors.
With $H=L^2((0,M),\dd x)$, let $V(s)$ be the closed pressure form
norm of the profile, including its $H$ norm, and put
\[
 \mathcal K_R[f,g]=Q_R[f,g]
                         +\lambda^3\Lambda_a(f,g)_H,\qquad
 \nu=\kappa\lambda^{-3/2}.
 \]
By Lemma~\ref{pr:closure} and \eqref{pr:fixed-positive-graph},
the form norms on $I_a$ are equivalent to the compatible
$X^1(s)$ norms. The norm comparison on $I_a$ in
Lemma~\ref{exlim:spatial-compactness} identifies them with one
fixed Hilbert space $V_a$. We use $V_a^*$ for its dual with
pivot $H$. The bounds supplied by Proposition~\ref{ct:proposition} are
\begin{align}
 &\sup_{I_a}\left\{\sum_{j=0}^{23}\|\Psi_j^\kappa\|_H^2
          +\sum_{j=0}^{22}\|\Psi_j^\kappa\|_{V_a}^2
          +(G_8^\kappa)^2\right\}\leq C_a,
 \label{exlim:low-uniform}\\
 &\int_{I_a}\nu\sum_{j=0}^{22}
       \mathcal K_{R^\kappa}[\Psi_{j+1}^\kappa,\Psi_{j+1}^\kappa]
                                                        \dd s\leq C_a,
 \label{exlim:dissipation}\\
 & B^\kappa\leq C_\beta\lambda^{J+1/2},\qquad
   G_8^\kappa\leq C_\beta\lambda^{J+1/2},\qquad
   C_{\rm rec}\eta^\kappa\leq\tfrac12.
 \label{exlim:decay-input}
\end{align}
Here $C_a$ may depend on the fixed positive interval, whereas
$C_\beta$, $C_{\rm rec}$, and the chosen $A_0$ do not depend
on $a$ or on $\kappa$. As throughout the argument for the chosen profile,
the notation $C_\beta$ also permits dependence on the fixed $e$,
finite orders, coefficient bounds, and prescribed cutoffs. These
parameters remain fixed during both limits.

The pressure form domain of Proposition~\ref{pr:pressure-domain},
the strain estimates of Lemma~\ref{en:strain-control}, and
\eqref{en:current-coercivity} give
\begin{equation}
 |Q_{R^\kappa}[f,g]|+\lambda^3\Lambda_a|(f,g)_H|
       \leq C_a\|f\|_{V_a}\|g\|_{V_a},\qquad
 \mathcal K_{R^\kappa}[f,f]\geq c_a\|f\|_{V_a}^2.
 \label{exlim:form-bounds}
\end{equation}
In particular, positivity here concerns the shifted form; no
unshifted coercivity assertion is being used.

Proposition~\ref{exlim:upstream-preparation} gives convergence of
the initial data in every graph norm used below. The prepared
derivatives $R_0,\ldots,R_{24}$ belong to the compatible domain,
whereas $R_{25}$ belongs only to the source space. No first vacuum
trace of $R_{25}$ is prescribed. These are regularity assertions
for the prepared data at fixed positive viscosity; their high norms
are not asserted to be uniform in $\kappa$.

The two source estimates in Proposition~\ref{mat:matched-source}
have different uses. Put
\[
 \mathcal R_u=\lambda^{-1}(A_{\tau\tau}+\mathcal M[A]),\qquad
 F^\kappa=-\mathcal R_u-\nu\lambda^{-1}\partial_s\mathcal M[A]
                         -\nu\lambda^{-1}\Lambda_a A_s,\qquad
 F_j^\kappa=\lambda^{-4}\partial_s^j(\lambda^4F^\kappa).
\]
The energy construction uses the estimate for material time derivatives
\begin{equation}
 \sum_{j=0}^{22}\|F_j^\kappa\|_H
 \leq C_\beta\left\{
  \lambda^{(K-7)/2}(1+|\log\lambda|)^{p_{\log}}
       +\kappa\lambda^{-9/2}+\kappa\Lambda_a\lambda^{-3/2}\right\}.
 \label{exlim:source22}
\end{equation}
The inviscid spatial recurrence uses instead
\begin{equation}
 \mathcal R_{23}^0:=
 \sum_{q=2}^{23}\sum_{j=0}^{23-q}\mathcal S_{q-2}(F_j^0)
 \leq C_\beta\lambda^{(K-5)/2}(1+|\log\lambda|)^{p_{\log}},
 \qquad F^0=-\mathcal R_u.
 \label{exlim:source21}
\end{equation}
Every derivative on the right satisfies $j+(q-2)\leq21$.
Estimate~\eqref{exlim:source21} does not control material derivatives of order $22$. The chosen $K$ makes its right-hand side smaller than
a fixed multiple of $\lambda^{J+1/2}$.

Finally, Proposition~\ref{exlim:upstream-recovery} is the uniform
estimate for compatible derivatives with strong sources. Once all
the norms involved are finite, its finite component inequalities give
\begin{equation}
 \mathcal X_{23}
 \leq C_{\rm rec}(B+\mathcal R_{23}^0)
                      +C_{\rm rec}\eta\mathcal X_{23}.
 \label{exlim:recovery-input}
\end{equation}
The inverse of the current linearized operator and the multiplier
estimates act on these completed spaces. Before applying
\eqref{exlim:recovery-input}, we show that the limiting derivatives
lie in the required strong operator domains.

\subsection{Completion, endpoint regularity, and compactness}
\label{exlim:compactness}

Both limiting operations use compactness on a fixed interval $I_a$,
$a>0$. The following lemmas give endpoint regularity, equivalence
with fixed Hilbert norms, and compactness in time. The equivalence
constants may depend on $a$. In a fixed coordinate neighborhood of the
vacuum boundary, write
\[
 \|f\|_{B_5^q(0,c)}^2=\sum_{r=0}^q
                              \int_0^c y^4|f^{(r)}(y)|^2\dd y.
 \]

\begin{lemma}
\label{exlim:endpoint}
On the bounded vacuum chart, the following statements hold.
\begin{enumerate}[label=(\roman*),leftmargin=2em]
\item For every integer $q\geq3$,
\begin{equation}
 \|f\|_{C^{q-3}([0,c])}\leq C_q\|f\|_{B_5^q(0,2c)}.
 \label{exlim:endpoint-embedding}
\end{equation}
Interior terms on the enlarged collar are included in the norm.
\item For $q\geq4$, the compatible completion is the maximal
distributional $B_5^q$ space with $f'(0)=0$.
For $q\leq3$, that completion imposes no continuous first-derivative trace
condition. 
\item In the high-order completion, the quotient identities
\eqref{exlim:ordinary-quotients} hold, and the quotient maps are bounded
$B_5^{k+5}\to C^k$ by Lemma~\ref{at:regular-factor-calculus}.
\end{enumerate}
\end{lemma}

The proof in Subsection~\ref{app:endpoint-completion} combines a
weighted Taylor estimate with approximation in the compatible
completion.

\begin{lemma}
\label{exlim:spatial-compactness}
Fix $0<a<A_0$. There are time-independent Hilbert norms
$F_a^q$, equivalent to $X^q(s)$ uniformly on $I_a$ for
$0\leq q\leq23$, such that
\begin{equation}
 F_a^{q+1}\Subset F_a^q\quad(0\leq q\leq22),\qquad
 V_a\Subset H\Subset V_a^*.
 \label{exlim:compact-inclusions}
\end{equation}
The equivalence constants may depend on $a$.
\end{lemma}

At fixed $a$, the moving charts are uniformly equivalent to a fixed
physical cover. An estimate for the weighted norm near the endpoint
then gives compactness; see
Subsection~\ref{app:fixed-scale-compactness}.

\begin{lemma}
\label{exlim:time-compactness}
On $I_a$, the following compactness criteria hold.
\begin{enumerate}[label=(\roman*),leftmargin=2em]
\item If $f^n$ is bounded in $L^\infty F_a^{q+1}$ and
$\partial_sf^n$ in $L^\infty F_a^q$, its continuous representatives
are precompact in $C(I_a;F_a^q)$.
\item If $f^n$ is bounded in $L^\infty H$ and $\partial_sf^n$ in
$L^2V_a^*$, its continuous representatives are precompact in
$C(I_a;V_a^*)$.
\end{enumerate}
\end{lemma}

The preceding spatial compactness and equicontinuity in time
give this result. The complete argument is in
Subsection~\ref{app:fixed-scale-compactness}.

Whenever $C(I_a;X^q(s))$ is used below, continuity is understood
after identification with a fixed space $F_a^q$ from
Lemma~\ref{exlim:spatial-compactness}. Uniform norm equivalence on $I_a$
makes this convention independent of that choice.

\Needspace{6\baselineskip}
We distinguish the derivative orders used in the viscosity limit:
\begin{itemize}[leftmargin=2em]
\item The material index $j$ counts derivatives $\partial_s^j$, through
$22$ in the commuted equation; the principal force term is $H_RR_j$.
The index $r+1$ includes the extra material derivative already present
in $(H_R+\Lambda_a)R_\tau$. It reaches $23$ only in the one-block term.
\item The potential-variation index $\ell$ denotes $D^\ell\mathcal U$.
A term with $k$ force directions, paired with a test, has
$\ell=k+1\le24$.
\item The spatial order $q$ is distinct from these time and variation
orders. The derivatives are continuous for $j+q\le22$ and essentially
bounded for $j+q\le23$.
\item The reversed physical-time derivative $\partial_\tau$ is
$\lambda^{-3/2}\partial_s$, and collapsing time is $t=-\tau$.
Its relation to $\partial_s$ therefore includes the displayed scale factor.
\end{itemize}

\subsection{Identification of the inviscid limit at fixed terminal scale}
\label{exlim:viscosity-limit}

\begin{lemma}
\label{exlim:top-viscous-slots}
Let the fixed-$a$ solutions of \eqref{exlim:viscous} satisfy
\eqref{exlim:low-uniform}--\eqref{exlim:form-bounds}.
\begin{enumerate}[label=(\roman*),leftmargin=2em]
\item The differentiated artificial viscosity through material order
$22$ tends to zero in $L^2(I_a;V_a^*)$ at rate $O_a(\sqrt\kappa)$.
\item If the radius derivatives through order $15$ converge strongly
in $C(I_a;X^7)$, the differentiated inviscid force through order $22$
has the weak form limit \eqref{exlim:bell-dual-limit}.
\end{enumerate}
Neither conclusion requires a viscosity-independent $V_a$ bound at
time order $23$ or a boundary trace at time order $25$. The latter
derivative lies only in the source space.
\end{lemma}

\begin{proof}
\emph{Variations of the force.}
On the prepared solutions, Proposition~\ref{loc:finite-calculus} gives
$\partial_s^j\mathcal M[R]=H_RR_j+\mathcal N_j[R]$ for
$1\le j\le23$, where
$\Pi_j$ denotes the partitions of $\{1,\ldots,j\}$ and
\begin{equation}
 \mathcal N_j[R]=
 \sum_{\substack{\pi\in\Pi_j\\|\pi|\geq2}}
 D^{|\pi|}\mathcal M[R]\,[R_{|B|}:B\in\pi].
 \label{exlim:bell}
\end{equation}
For regular compatible directions,
\begin{equation}
 \rho[R+\textstyle\sum t_if_i]
 =\rho[R]\left(1+\sum t_if_i/R\right)^{-2}
          \left(1+\sum t_i(f_i)_x/R_x\right)^{-1}.
 \label{exlim:strain-identity}
\end{equation}
Set $Z_R(f)=\|f/R\|_\infty+\|f_x/R_x\|_\infty$, using the regular
endpoint values. Since $0<\rho\le C_a$ on a fixed positive-scale
interval, differentiation of
\[
 \frac{P(\rho)}\rho
 =\frac43\rho^{2/3}\int_0^1
       \frac{t^{2/3}}{\sqrt{1+t^{2/3}\rho^{2/3}}}\dd t
\]
gives
\[
 |\rho^l\partial_\rho^l(P(\rho)/\rho)|\le C_{a,l}P'(\rho).
\]
Cauchy--Schwarz therefore controls two strains in the form norm.
The gravitational variation
\[
 D^p\mathcal U_{\rm grav}[R][f_1,\ldots,f_p]
 =(-1)^{p+1}p!\int_0^M xR^{-p-1}\prod_{i=1}^pf_i\dd x
\]
has the same bound. Thus, for $2\le\ell\le24$,
\begin{equation}
 |D^\ell\mathcal U[R][f_1,\ldots,f_\ell]|
 \le C_{a,\ell}\|f_1\|_{V_a}\|f_2\|_{V_a}
                     \prod_{i=3}^{\ell}Z_R(f_i).
 \label{exlim:two-slot-bound}
\end{equation}
The two form arguments may occupy any positions, by symmetry;
approximation in the closed form domain removes their strain
boundedness assumption. This is the estimate of
Section~\ref{en:section}.

For a proper partition of $j\le23$, order the block sizes as
$h_1\ge h_2\ge\cdots$. Then
\[
 h_1\le j-1\le22,\qquad 2h_2\le j\le23,\qquad h_2\le11.
\]
Choose $R_{h_1}$ and the test function as the two form arguments.
All other strains are bounded by \eqref{exlim:low-uniform};
hence $\|\mathcal N_j[R]\|_{V_a^*}\le C_a$ for $j\le23$.

\emph{Viscosity.}
Since $(H_R+\Lambda_a)R_\tau
=\lambda^{-3/2}(\partial_s\mathcal M[R]+\Lambda_aR_s)$,
the differentiated viscous term is
\begin{equation}
 \kappa\sum_{r=0}^j\binom jr
 \partial_s^{j-r}(\lambda^{-3/2})
 \left\{H_RR_{r+1}+\mathcal N_{r+1}[R]
                            +\Lambda_aR_{r+1}\right\}.
 \label{exlim:viscous-bell}
\end{equation}
Every summand is $O_a(\kappa)$ in $L^\infty V_a^*$ except the
term with $j=r=22$. For this term,
\[
 R_{23}^\kappa=A_{23}+\lambda^4\Psi_{23}^\kappa,\qquad
 \kappa\int_{I_a}\|R_{23}^\kappa\|_{V_a}^2\dd s\le C_a
\]
by \eqref{exlim:dissipation}--\eqref{exlim:form-bounds}. For
$0<\kappa\le1$, it follows that
\begin{equation}
 \begin{split}
 &\sum_{j=0}^{22}
 \left\|\partial_s^j
  \{\kappa(H_{R^\kappa}+\Lambda_a)R_\tau^\kappa\}
                    \right\|_{L^2(I_a;V_a^*)}\\
 &\qquad\leq C_a\kappa+
 C_a\sqrt\kappa
 \left(\kappa\int_{I_a}\|R_{23}^\kappa\|_{V_a}^2\dd s\right)^{1/2}
 \leq C_a\sqrt\kappa\longrightarrow0.
 \end{split}
 \label{exlim:all-viscosity-zero}
\end{equation}
Only the proper-partition terms use the potential variation of
order $24$; their largest material index is at most $22$.

\emph{Force convergence.}
For $j\le22$, write each chain-rule term as $T^\kappa R_h^\kappa$,
with $h$ the largest block size. The other indices are at most
$11$, so their strains converge uniformly. By
\eqref{exlim:two-slot-bound},
\[
 T^\kappa\to T\quad\hbox{in }L^\infty\mathcal L(V_a,V_a^*),
 \qquad R_h^\kappa\rightharpoonup^*R_h
                 \quad\hbox{in }L^\infty V_a.
\]
For $g\in L^1(I_a;V_a)$,
\[
 \begin{split}
 \left|\int_{I_a}\langle T^\kappa R_h^\kappa-TR_h,g\rangle\dd s\right|
 &\le \|T^\kappa-T\|_{L^\infty\mathcal L(V_a,V_a^*)}
       \|R_h^\kappa\|_{L^\infty V_a}\|g\|_{L^1V_a}\\
 &\quad+\left|\int_{I_a}
          \langle T(R_h^\kappa-R_h),g\rangle\dd s\right|
 \longrightarrow0.
 \end{split}
\]
Summing the finite partition expansion gives
\begin{equation}
 \partial_s^j\mathcal M[R^\kappa]
 \rightharpoonup^*H_RR_j+\mathcal N_j[R]
 \quad\hbox{in }L^\infty(I_a;V_a^*),\qquad 1\leq j\leq22.
 \label{exlim:bell-dual-limit}
\end{equation}
\end{proof}

\begin{lemma}
\label{exlim:identified-limit}
Fix $a>0$. A subsequence of the family
\eqref{exlim:viscous} converges to an inviscid solution $R$
on $I_a$, with the limiting prepared data, such that
\begin{equation}
 R_{\tau\tau}+\mathcal M[R]=0,\qquad
 R_j^\kappa\to R_j\ \hbox{in }C(I_a;H)\quad(0\leq j\leq22).
 \label{exlim:identified-euler}
\end{equation}
Moreover,
\begin{equation}
 R_j^\kappa\to R_j\ \hbox{in }C(I_a;X^7(s))\quad(j\leq15),
 \qquad
 R_{23}^\kappa\to R_{23}\ \hbox{in }C(I_a;V_a^*),
 \label{exlim:low-strong}
\end{equation}
and $R_{23}$ is weakly continuous in $H$. The inherited
low-order weighted bounds hold for these distributional derivatives.
\end{lemma}

\begin{proof}
\emph{Compactness of the radius derivatives.}
The profile and scale are fixed smooth functions on $I_a$.
Hence \eqref{exlim:low-uniform} gives equivalent fixed-scale
bounds for the unnormalized derivatives $R_j^\kappa$.
For $0\le j\le22$,
\[
 \sup_\kappa\bigl(
 \|R_j^\kappa\|_{L^\infty V_a}
   +\|\partial_sR_j^\kappa\|_{L^\infty H}\bigr)\le C_a,
 \qquad \partial_sR_j^\kappa=R_{j+1}^\kappa.
\]
The compact inclusion $V_a\Subset H$ and equicontinuity in $H$
give the first convergence in \eqref{exlim:identified-euler}.
For $j\le15$, the stronger bounds are
\[
 \sup_\kappa\bigl(
 \|R_j^\kappa\|_{L^\infty X^8}
   +\|R_{j+1}^\kappa\|_{L^\infty X^7}\bigr)\le C_a,
 \qquad j+1+7\le23.
\]
Lemma~\ref{exlim:time-compactness} then gives strong convergence in
$C(I_a;X^7(s))$.

The remaining weighted derivatives in \eqref{exlim:low-uniform}
have weak-star subsequential limits, identified by testing in the
interior.

\emph{Passage to the limit in the equation.}
Lemma~\ref{exlim:endpoint} shows that the first convergence in
\eqref{exlim:low-strong} is ordinary $C^4$ convergence of the
vacuum radius and sufficient Cartesian $C^3$ convergence at the
center. It therefore gives strong $C^2$ convergence of the
normalized vacuum quotient and the deformation factors. Their
positive lower bounds persist. To see directly that the force
converges, put
\[
 \xi=(M-x)^{1/5},\quad Q=-R_\xi/\xi,\quad
 a_v=\rho/\xi^3=\frac5{4\pi R^2Q},\quad
 p(\xi^2,a_v)=\frac{4a_v^{2/3}}
                         {3\sqrt{1+\xi^2a_v^{2/3}}}.
 \]
On the vacuum collar its nonsingular expression is
\begin{equation}
 \mathcal M[R]=-\frac{4\pi}{5}R^2
            p(\xi^2,a_v)(3a_v+\xi(a_v)_\xi)+\frac{x}{R^2}.
 \label{exlim:vacuum-force}
\end{equation}
The quantity $\xi(a_v)_\xi$ is a smooth function of $R,Q$,
$\xi R_\xi$, and $R_{\xi\xi}$, affine in the last argument.
At the center, with
$\Phi(X)=R(|X|^3)X/|X|=r_c(X)X$, the expression is
\begin{equation}
 U_{\mathcal M[R]}(X)
  =\frac{4\pi}{3}r_c^2\nabla P(\rho)+\frac{X}{r_c^2},
 \qquad
 \rho=\frac3{4\pi\det D\Phi},\qquad
 r_c=\int_0^1t^2\operatorname{div}\Phi(tX)\dd t.
 \label{exlim:center-force}
\end{equation}
These formulas, and the interior formula \eqref{exlim:force},
imply $\mathcal M[R^\kappa]\to\mathcal M[R]$ in $C(I_a;H)$.

By \eqref{exlim:dissipation}--\eqref{exlim:form-bounds},
\begin{equation}
 \left\|\kappa(H_{R^\kappa}+\Lambda_a)R_\tau^\kappa
                         \right\|_{L^2(I_a;V_a^*)}
 \leq C_a\sqrt{\kappa}.
 \label{exlim:viscous-zero}
\end{equation}
The profile contribution in
$R_\tau^\kappa=\lambda^{-3/2}(A_s+\lambda^4\Psi_1^\kappa)$
is bounded on $I_a$ and, after multiplication by $\kappa$, tends to zero.
Thus \eqref{exlim:viscous} passes to the Euler equation
\eqref{exlim:identified-euler}, first in distributions and then in
$C(I_a;H)$.

\emph{The highest time derivative.}
Lemma~\ref{exlim:top-viscous-slots} now applies to the convergence at
low orders just established. In particular all
commuted viscosity terms vanish by \eqref{exlim:all-viscosity-zero}.
Using $\partial_\tau=\lambda^{-3/2}\partial_s$ gives
\[
 \partial_s^{22}R_{\tau\tau}^\kappa
  =\lambda^{-3}R_{24}^\kappa
       +\sum_{\ell=1}^{23}c_\ell(s)R_\ell^\kappa,
\]
where the coefficients $c_\ell$ are smooth and bounded on $I_a$.
The inviscid force derivative uses only derivatives in the form domain through $22$;
hence the differentiated equation gives
$\|\partial_sR_{23}^\kappa\|_{L^2(I_a;V_a^*)}\leq C_a$.
Lemma~\ref{exlim:time-compactness} proves the last convergence in
\eqref{exlim:low-strong}. Testing first against the dense space
$V_a$, and then using the uniform $H$ bound, gives weak
$H$ continuity. No order-$23$ strong source is used.

Passing the identities for consecutive derivatives through these
convergences identifies every $R_j$ as the corresponding
distributional derivative of $R$. Convergence of the prepared
data identifies their initial values. The argument does not require
$R_{24}$ to lie in a strong operator domain or a boundary trace
for the source derivative $R_{25}$.
\end{proof}

\begin{lemma}
\label{exlim:strong-rows}
Let $R$ be the limit in Lemma~\ref{exlim:identified-limit}.
\begin{enumerate}[label=(\roman*),leftmargin=2em]
\item For $0\le j\le21$, the derivative
$\partial_s^j\mathcal M[R]$ is given in $H$ by the iterated chain
rule. The corresponding form identities hold in $V_a^*$ through
$j=22$.
\item For $0\le j\le14$, the same identities follow by
differentiating the force as a map $F_a^7\to H$ on the open set
where the normalized radius and Jacobian are positive.
\end{enumerate}
\end{lemma}

\begin{proof}
Combining the weak force limit \eqref{exlim:bell-dual-limit} of
Lemma~\ref{exlim:top-viscous-slots} with the vanishing viscosity
estimate \eqref{exlim:all-viscosity-zero} gives the differentiated
Euler equation. At the highest order, its integral form is the
following identity in $V_a^*$, for every $s,t\in I_a$:
\begin{equation}
 R_{23}(s)-R_{23}(t)
 =-\int_t^s\lambda^3
 \left\{\sum_{\ell=1}^{23}c_\ell R_\ell
                 +H_RR_{22}+\mathcal N_{22}[R]\right\}\dd u.
 \label{exlim:actual-top-row}
\end{equation}
Here the coefficients are evaluated at $u$; the strong
$C(I_a;V_a^*)$ convergence of the derivative of order $23$ identifies both traces.

\emph{Identification in $H$.}
For $j\leq21$, the radius derivative with the largest time index
is bounded in $X^2$, since $j+2\leq23$. The nonsingular force
coefficients are smooth functions of the radius, its first derivative,
and its regular quotients, and are affine in its second derivative.
Each differentiated monomial therefore contains at most one second
spatial derivative of a direction. The factor with the largest time
index converges weakly in $B_5^2$, or in Cartesian $H^2$, while
the remaining factors converge uniformly. The quotient extends continuously
to the weighted completion by the estimate
\[
 \|\chi f_y/y\|_{B_5^m}
 \leq C_m\bigl(\|\chi f\|_{B_5^{m+2}}+\|f\|_{B_5^0}\bigr)
 \]
from Lemma~\ref{at:natural-quotient}. Thus the limit is the
chain-rule expression in $H$. Distributional differentiation
identifies its higher spatial derivatives once the corresponding
weighted Sobolev bounds have been established.

\emph{Direct differentiation at lower orders.}
For $j\leq14$, identify the moving spaces
with $F_a^7$. By \eqref{exlim:low-strong}, the consecutive derivatives
$R_0,\ldots,R_{15}$ are continuous in this space and
$\partial_sR_i=R_{i+1}$. In particular, $R\in C^{15}(I_a;F_a^7)$.
On the positive-radius domain, the nonsingular formulas
\eqref{exlim:vacuum-force} and \eqref{exlim:center-force}, together
with the quotient estimate above, define a smooth map into $H$.
Its first derivative is justified by
\[
 \mathcal M[R+h]-\mathcal M[R]-H_Rh
       =\int_0^1(H_{R+th}-H_R)h\dd t,
\]
whose $H$ norm is $o(\|h\|_{F_a^7})$ as $h\to0$ in $F_a^7$.
The coefficient formulas give the same continuity for each of the
finitely many higher variations. Applying the ordinary chain rule
to $s\mapsto R(s)$ therefore gives the stated expressions through
$j=14$.
\end{proof}

\subsection{Spatial regularity and strong operator domains}
\label{exlim:realization}

The low-order compactness argument does not yet give the full
spatial regularity required for recovery. For the radius, the radial
force equation determines an integral of the pressure. The strictly
increasing pressure law then determines the density and the radius
Jacobian. Once the radius has gained one derivative, the linearized
flux identity gives the same gain for each identified time derivative.
The next two lemmas prove these gains on a fixed interval of positive
scales, using the norms already controlled there. They also place the
derivatives in the strong operator domains needed for the uniform
inverse estimate.

In reversed physical time, the corresponding chain-rule identity
for a smooth inviscid radius is
\begin{equation}
 \partial_\tau^{j+2}R+\mathcal H_R\partial_\tau^jR
 =-\mathcal N_j[R;\partial_\tau R,\ldots,\partial_\tau^{j-1}R],
 \qquad j\ge1,
 \label{in:time-rows}
\end{equation}
where $\mathcal H_R=H_R=D\mathcal M[R]$ and the remainder is the
partition sum in \eqref{exlim:bell}, with each $R_i$ replaced by
$\partial_\tau^iR$. It is zero for $j=1$. For the finite-regularity limit,
Lemma~\ref{exlim:strong-rows} justifies the differentiated identities
used in the spatial induction below, in rescaled time.

\begin{lemma}
\label{exlim:nonlinear-realization}
Fix $a>0$, $9\leq q\leq23$, and put $m=q-2$.
Suppose that $R$ satisfies the positivity and lower-order bounds
above. Then
\[
 R\in X^{q-1},\qquad \mathcal S_m(\mathcal M[R])<\infty
       \quad\Longrightarrow\quad R\in X^q.
\]
The $X^q$ norm is bounded on bounded sets of the two norms on the
left and the positive lower-order factors. The bound may depend on
$a$, but not on the norm being established.
\end{lemma}

\begin{proof}
All calculations are on enlarged fixed physical charts. We treat
the vacuum, center, and interior separately, then verify convergence
of compatible approximations.

\emph{Vacuum.} Let $F=\mathcal M[R]$ and use $Q,a_v$ as in
\eqref{exlim:vacuum-force}. Define
\begin{equation}
 g=-\frac5{4\pi R^2}\left(F-\frac{x}{R^2}\right),\qquad
 Sf(\xi)=\int_0^1t^4f(t\xi)\dd t,\qquad
 \Pi=\frac{P(\rho)}{\xi^5}.
 \label{exlim:pressure-primitive}
\end{equation}
The force equation says $P_\xi=\xi^4g$. Since the actual
pressure is zero at $\xi=0$,
\begin{equation}
 \Pi=Sg,\qquad \xi\Pi_\xi=g-5\Pi,\qquad
 \|(Sf)^{(r)}\|_{B_5^0}
              \leq\frac1{r+5/2}\|f^{(r)}\|_{B_5^0}.
 \label{exlim:S-bound}
\end{equation}
The last bound follows by Minkowski's inequality: the $r$-th
derivative gives $t^r$, and dilation in $L^2(\xi^4\dd\xi)$
costs $t^{-5/2}$, leaving
$\int_0^1t^{r+3/2}\dd t=(r+5/2)^{-1}$.

The exact pressure law gives the nonsingular scalar relation
\begin{equation}
 \Pi=\mathscr P(\xi^2,a_v),\qquad
 \mathscr P(\sigma,a)=\frac43a^{5/3}
       \int_0^1\frac{v^{2/3}}{\sqrt{1+\sigma a^{2/3}v^{2/3}}}\dd v,
 \qquad
 \mathscr P_a=\frac{4a^{2/3}}{3\sqrt{1+\sigma a^{2/3}}}>0.
 \label{exlim:EOS-inverse}
\end{equation}
On the known positive range its inverse $a=\mathscr J(\sigma,\Pi)$
is smooth. This follows recursively by differentiating
$\mathscr P(\sigma,\mathscr J)=\Pi$ and dividing only by the
displayed positive derivative. The weighted composition and product
estimates at $m\geq7$ control the derivatives of these coefficient
functions through order $m$. Hence $a_v,Q\in B_5^m$, and
\begin{align}
 \xi(a_v)_\xi
   &=2\xi^2\mathscr J_\sigma
                     +\mathscr J_\Pi(g-5\Pi),\nonumber\\
 \xi Q_\xi
   &=-2Q\,\xi R_\xi/R-Q\,\xi(a_v)_\xi/a_v,\qquad
 R_{\xi\xi}=-Q-\xi Q_\xi\in B_5^m.
 \label{exlim:nonlinear-edge-gain}
\end{align}
The assumed lower-order bound supplies $R_\xi\in B_5^m$; every other
factor has just been obtained from the source. This proves
$R\in B_5^{m+2}$, with no appeal to an unknown high norm.

\emph{Center.} At the center use \eqref{exlim:center-force}.
Initially $\Phi\in H^{m+1}$ gives only $r_c\in H^m$ through
the divergence integral. Since $m\geq7$, this is sufficient
for multiplication of the $H^m$ source. Therefore
$\nabla P(\rho)\in H^m$, and hence $P(\rho)\in H^{m+1}$;
its $L^2$ norm is controlled by the known lower-order bounds for the radius.
At fixed $a>0$, the lower-order estimates give
\[
 0<c_a\le\rho\le C_a,\qquad
 \inf_{[c_a,C_a]}P'>0.
\]
Thus $P^{-1}$ is smooth on a neighborhood of the pressure range, and
Sobolev composition yields
\[
 \rho=P^{-1}(P(\rho))\in H^{m+1},\qquad
 J_0=\frac3{4\pi\rho}=\det D\Phi\in H^{m+1}.
\] The radial volume identity is
\begin{equation}
 r_c^3(X)=3\int_0^1t^2J_0(tX)\dd t.
 \label{exlim:center-volume}
\end{equation}
Boundedness of $r_c$ at the center forces the integration constant to vanish.
Minkowski's inequality gives, for $0\le\ell\le m+1$,
\[
 \|D^\ell(r_c^3)\|_{L^2(B)}
 \le3\int_0^1t^{\ell+2-3/2}\dd t\,
                         \|D^\ell J_0\|_{L^2(B)}
 =\frac3{\ell+3/2}\|D^\ell J_0\|_{L^2(B)}.
\]
Taking the positive cubic root yields $r_c\in H^{m+1}$, and hence
\[
 \operatorname{div}\Phi=2r_c+J_0/r_c^2\in H^{m+1},
 \qquad \operatorname{curl}\Phi=0.
\]
Choose $\chi=1$ on $B'\Subset B$, supported in $B$. The Fourier identity
\[
 |\zeta|^2|\widehat{\chi\Phi}|^2
   =|\zeta\cdot\widehat{\chi\Phi}|^2
                         +|\zeta\times\widehat{\chi\Phi}|^2
\]
gives
\[
 \|\Phi\|_{H^{m+2}(B')}
 \le C\bigl(\|\operatorname{div}\Phi\|_{H^{m+1}(B)}
                    +\|\Phi\|_{H^{m+1}(B)}\bigr),\qquad B'\Subset B.
\]
The last term controls the cutoff commutators.

\emph{Interior and completion.} On a fixed interior interval integrate
$P_x=(F-x/R^2)/(4\pi R^2)$, invert the positive pressure law,
and use $R_x=(4\pi R^2\rho)^{-1}$.
These steps give successively $P,\rho\in H^{m+1}$ and
$R\in H^{m+2}$. On the part of the long chart away from the
vacuum, the same calculation applies after the fixed-scale coordinate
change.

To approximate $R$ by compatible functions, use its first vacuum
trace, already supplied by $X^8$, and the one-sided approximation
of Lemma~\ref{exlim:endpoint}. At the center, mollify the equivariant
vector and average over rotations; smooth the remaining pieces away
from the endpoints. The fixed-$a$ coordinate changes and all
derivatives of the localization cutoffs are bounded. An enlarged
mass partition therefore gives $R_n\in\mathscr C$ such that, with
$v_n=R_n-R$,
\[
 \begin{aligned}
 \|v_n\|_{X^q}^2
 ={}&\|\chi_cU_{v_n}\|_{H^q}^2+\|\chi_iv_n\|_{H^q}^2\\
 &+\lambda^4\|\chi_5v_n\|_{B_5^q}^2
       +\lambda^4\|\chi_8v_n\|_{C_8^q}^2\longrightarrow0.
 \end{aligned}
\]
The same approximation
applies to the linear gain below.
\end{proof}

\begin{lemma}
\label{exlim:linear-realization}
Let $9\leq q\leq23$, and suppose $R\in X^q$ satisfies the
positivity and lower-order bounds above. Suppose that
\[
 u\in X^{q-1},\qquad \mathcal S_{q-2}(H_Ru)<\infty,
\]
where $H_Ru$ is interpreted in the sense of distributions. Assume
the compatible first-derivative trace at the vacuum, which is
defined by the lower-order regularity of $u$. Then $u\in X^q$. Its norm is bounded
in terms of the source, the preceding norm of $u$, the radius
$X^q$ norm, and the fixed lower-order and positivity bounds.
\end{lemma}

\begin{proof}
\emph{Vacuum.} At the vacuum the exact self-adjoint expression, with fixed-scale factors absorbed into its source, is
\begin{equation}
 -\xi^{-4}(\xi^4\mathsf a\,u_\xi)_\xi+Wu=F,\qquad
 \mathsf a=\frac{P'(\rho)}{\xi^2Q^2}>0.
 \label{exlim:linear-edge}
\end{equation}
The coefficient identities \eqref{at:edge-principal}--\eqref{at:edge-potential}
and the product estimates of Lemma~\ref{at:same-cutoff-product}
give, for $m=q-2$,
\begin{equation}
 \alpha=\mathsf a^{-1},\quad \xi\alpha_\xi,\quad W\in B_5^m,
 \qquad \mathsf a\in C^1,\qquad
 \mathsf a,\mathsf a^{-1}\in L^\infty.
 \label{exlim:edge-coefficients}
\end{equation}
Only $R$ through order $q$ and its natural quotient through
order $q-2$ enter these coefficients. In particular the term
containing the highest quotient derivative is retained as
\[
 \xi Q^{(q-1)}
       =-R^{(q)}-(q-1)Q^{(q-2)};
 \]
an extra unweighted derivative of $Q$ is not required.

Let $g=F-Wu$. Since $u\in B_5^1$ and $g\in B_5^0$, the
distributional flux in \eqref{exlim:linear-edge} has the form
\[
 \xi^4\mathsf a\,u_\xi=C-\int_0^\xi t^4g(t)\dd t.
 \]
By Cauchy--Schwarz,
\[
 \left|\int_0^\xi t^4g(t)\dd t\right|
 \le\left(\int_0^\xi t^4\dd t\right)^{1/2}
       \|g\|_{B_5^0(0,\xi)}
 =\frac{\xi^{5/2}}{\sqrt5}\|g\|_{B_5^0(0,\xi)}=o(1).
\]
If $C\ne0$, boundedness of $\mathsf a$ would give
$|u_\xi|\ge c\xi^{-4}$ for small $\xi$, and hence
\[
 \int_0^\varepsilon\xi^4|u_\xi|^2\dd\xi
   \ge c^2\int_0^\varepsilon\xi^{-4}\dd\xi=\infty,
\]
contrary to $u\in B_5^1$.
Thus $C=0$, and
\begin{equation}
 u_\xi=-\alpha\,\xi Sg,\qquad
 u_{\xi\xi}=-\alpha g+(4\alpha-\xi\alpha_\xi)Sg.
 \label{exlim:linear-flux-gain}
\end{equation}
At $m\geq7$, the products on the right belong to $B_5^m$
by \eqref{exlim:S-bound}--\eqref{exlim:edge-coefficients}.
The term $Wu$ uses only $u\in B_5^{q-1}$.
Hence $u\in B_5^q$. Its first-derivative trace vanishes by the
lower-order regularity. Lemma~\ref{exlim:endpoint} therefore places
$u$ in the compatible completion at order $q$. The flux argument
alone does not give a classical first trace at order two.

\emph{Center.} At the center put $G=D\Phi$, $B_c=G^{-T}$, and define
$\nabla_{B_c}=B_c\nabla$, with its corresponding divergence
and curl. The pressure principal part is
$-P'(\rho)\nabla_{B_c}\operatorname{div}_{B_c}u$.
Add $P'(\rho)\operatorname{curl}_{B_c}
\operatorname{curl}_{B_c}u$. For radial $u$, the matrices
$Du$ and $B_c$ have the same radial and tangential eigenspaces,
so $\operatorname{curl}_{B_c}u=0$. The added term vanishes,
also distributionally by radial $H^1$ approximation. After this
addition, the principal part acts diagonally on the vector components:
\[
 -a^{kl}\partial_k\partial_lu,\qquad
 a^{kl}=P'(\rho)(B_c^TB_c)_{kl}.
 \]
The matrix is uniformly positive on the fixed center ball.
Its coefficients are in $H^{m+1}$; the lower coefficients
are smooth functions of $D\Phi$, affine in $D^2\Phi$, and
belong to $H^m$. The gravitational coefficient is a smooth
function of the nonsingular quotient $r_c$.
The hypothesis $u\in X^{q-1}$ gives $u\in H^{m+1}$ on the
center chart. Together with $a\in H^{m+1}$ and $m=q-2\geq7$,
this supplies the Lipschitz control of the leading coefficient
used in the difference-quotient argument below.

Commute $m$ Cartesian derivatives and put $w=D^mu\in H^1$.
For the principal commutators, $m\ge7$ gives
\[
 \begin{aligned}
 \sum_{k=1}^m\|D^ka\,D^{m+2-k}u\|_2
 &\le
 \sum_{k\le(m+1)/2}\|D^ka\|_\infty\|D^{m+2-k}u\|_2\\
 &\quad+
 \sum_{k>(m+1)/2}\|D^ka\|_2\|D^{m+2-k}u\|_\infty\\
 &\le C_m\|a\|_{H^{m+1}}\|u\|_{H^{m+1}}.
 \end{aligned}
\]
The Sobolev orders used for the $L^\infty$ factors are respectively
$m+1-k>3/2$ and $k-1>3/2$. Lower-order commutators have one fewer
derivative. The complete commuted source is therefore in $L^2$.

To prove the local gain from $H^1$ to $H^2$ for $w$,
write the principal part in divergence form and move
$(\partial_ka^{kl})\partial_lw$ to the $L^2$ source. Test the weak equation with
$-\delta_{-h}(\chi^2\delta_hw)$, where $\chi$ has compact
support in a larger ball and
$\delta_hw=(w(\,\cdot+h e_l)-w)/h$ is a coordinate difference
quotient. Ellipticity, the Lipschitz coefficient
bound, Cauchy--Schwarz, and absorption give
\[
 \|\chi\nabla\delta_hw\|_2^2
       \leq C\bigl(\|F_{\rm comm}\|_2^2+\|w\|_{H^1}^2\bigr).
 \]
The bound is uniform for small $h$. Weak compactness identifies
the second distributional derivatives of $w$ in $L^2$,
which proves $u\in H^{m+2}$ on the smaller ball.
The same one-dimensional argument applies on fixed interior
intervals and on the part of the long chart away from the vacuum.
Combining the estimates on the enlarged charts proves the lemma.
\end{proof}

\begin{lemma}
\label{exlim:realization-induction}
For the identified inviscid limit on $I_a$,
\[
 R_j\in L^\infty(I_a;X^q),\qquad j+q\le23.
\]
When $q\ge2$, these derivatives also lie in the strong domains of
the linearized operator. The bounds may depend on $a$; their constants
do not enter the uniform recovery inequality.
\end{lemma}

\begin{proof}
The bounds for $q\le8$, $j+q\le23$, are already known.
For the remaining pairs, induct in the order
\begin{equation}
 q=9,\ldots,23,\qquad j=0,\ldots,23-q.
 \label{exlim:induction-order}
\end{equation}
At $j=0$, the equation gives
\[
 \mathcal M[R]=-\lambda^{-3}(R_2-\tfrac32bR_1),\qquad
 \underset{s\in I_a}{\operatorname{ess\,sup}}\,
       \mathcal S_{q-2}(\mathcal M[R])<\infty.
\]
Indeed, $R_1,R_2$ are controlled at spatial order $q-2$.
Lemma~\ref{exlim:nonlinear-realization} yields
$R\in L^\infty X^q$.

For $j\ge1$, Lemma~\ref{exlim:strong-rows} gives
\begin{equation}
 H_RR_j=-\partial_s^jR_{\tau\tau}-\mathcal N_j[R].
 \label{exlim:realized-row}
\end{equation}
The source has the following bounds:
\[
 \begin{aligned}
 \partial_s^jR_{\tau\tau}
     &=\sum_{\ell=1}^{j+2}c_{j\ell}(s)R_\ell,
 &\ell+(q-2)&\le j+q\le23,\\
 \mathcal N_j[R]
     &=\sum_{\substack{\pi\in\Pi_j\\|\pi|\ge2}}
       D^{|\pi|}\mathcal M[R]
                  [R_{|B|}:B\in\pi],
 &\max_{B\in\pi}|B|&<j.
 \end{aligned}
\]
The first line uses lower spatial orders. In the second, each
direction has already been recovered in $X^q$ at an earlier time
index. The product estimates give
\[
 \underset{s\in I_a}{\operatorname{ess\,sup}}\,
 \bigl\{\mathcal S_{q-2}(\partial_s^jR_{\tau\tau})
        +\mathcal S_{q-2}(\mathcal N_j[R])\bigr\}<\infty.
\]
Since $j\le23-q\le14$, the direct chain rule of
Lemma~\ref{exlim:strong-rows} applies. Lemma~\ref{exlim:linear-realization}
now gives $R_j\in L^\infty X^q$.

At each step the new bound depends only on $a$, the fixed scale
function, the positive lower-order factors and previously bounded
norms. The induction is finite. Distributional identification of
the derivatives in separable Hilbert spaces gives strong
measurability.

Finally Lemma~\ref{loc:finite-strong-graph} gives, for every
recovered direction $f$ and $q\ge2$, compatible approximants with
\[
 \|f_n-f\|_{X^q}
       +\mathcal S_{q-2}(L_Rf_n-L_Rf)\longrightarrow0.
\]
For $q\le8$ use the order-$23$ radius regularity; for $q\ge9$
the radius and direction have the same spatial regularity.
The asymmetric product bounds give convergence of the differentiated
sources as well. Thus all required derivatives lie in the closed
strong domains. None of these fixed-$a$ constants enters the
uniform recovery estimate.
\end{proof}

\begin{proposition}
\label{exlim:fixed-terminal}
For every $0<a<A_0$, the viscosity limit is an inviscid
solution on $I_a$ with the limiting prepared data. Its normalized
derivatives have the following properties.
\begin{enumerate}[label=(\roman*),leftmargin=2em]
\item \emph{Regularity.}
\begin{equation}
 \Psi_j\in L^\infty(I_a;X^{23-j})\quad(0\leq j\leq23),\qquad
 \Psi_j\in C(I_a;X^{22-j})\quad(0\leq j\leq22),
 \label{exlim:full-cone}
\end{equation}
They are weakly continuous in the spaces of total order $23$.
The normalized derivative of order $23$ is also strongly continuous
in $V_a^*$.
\item \emph{Uniform bound.}
\begin{equation}
 \mathcal X_{23}(s)\leq C_\beta\lambda(s)^{J+1/2}.
 \label{exlim:full-decay}
\end{equation}
The constant and $A_0$ are independent of $a$.
\end{enumerate}
\end{proposition}

\begin{proof}
\emph{Uniform spatial bounds.}
Lemma~\ref{exlim:realization-induction} gives essentially bounded
norms in time and the required strong domains. For
$2\leq q\leq23$ and $j+q\leq23$, the differentiated forces
have been identified in their completed source spaces. At $q\leq8$,
the one possible high time factor has its $X^q$ bound and the
others have $X^8$ bounds, so the closed asymmetric multiplier
estimate applies. At $q\geq9$, the preceding induction gives the
full spatial regularity.

The finite recovery inequalities now apply to the identified
derivatives and give \eqref{exlim:recovery-input}. Every norm in
this inequality is finite, and the chosen smallness threshold allows
us to absorb its last term. Thus
\begin{equation}
 \mathcal X_{23}\leq2C_{\rm rec}(B+\mathcal R_{23}^0).
 \label{exlim:uniform-recovery-after-realization}
\end{equation}
Along the retained viscosity subsequence $\kappa_n\to0$, the
bound for $B$ passes to the limit as a finite sum of norms:
\begin{equation}
 B(s)\leq\liminf_{n\to\infty}B^{\kappa_n}(s)
                         \leq C_\beta\lambda(s)^{J+1/2}.
 \label{exlim:base-lsc}
\end{equation}
For derivatives through order $22$, identify the weak $X^1$ limits by
their strong limits in the Hilbert space; for the derivative of order $23$, use its weak $H$
limit. Lower semicontinuity applies to the finite sum defining $B$; no
passage to the limit in the modified nonlinear energy is needed.
The same reasoning preserves the small $G_8$ bound.

Estimate~\eqref{exlim:source21} completes \eqref{exlim:full-decay}.
Only the already selected $C_{\rm rec}$, the original bound for
$B$, and the profile source constant enter this last step.
The fixed-scale regularity and approximation estimates may depend
on $a$, but their constants do not enter the choice of $A_0$.

\emph{Time continuity.}
For $j\le22$, Bochner integration in $F_a^{22-j}$ gives
\[
 R_j(t)-R_j(s)=\int_s^tR_{j+1}(\sigma)\dd\sigma,\qquad
 \|R_j(t)-R_j(s)\|_{F_a^{22-j}}
 \le |t-s|\|R_{j+1}\|_{L^\infty F_a^{22-j}}.
\]
The weaker prepared-data convergence identifies the initial value.
If $s_n\to s$, the highest-order bound gives weakly convergent
subsequences. Strong convergence at lower orders identifies every such
limit with $R_j(s)$. Hence the representative at the highest order is
weakly continuous and
\[
 \|\Psi_j(s)\|_{X^{23-j}(s)}
 \le\liminf_{n\to\infty}\|\Psi_j(s_n)\|_{X^{23-j}(s_n)},
\]
using the continuous fixed-scale coefficients. This extends the
almost-everywhere bounds to all times. For the derivative of order $23$,
use the established $C_sV_a^*$ and weak $C_sH$ representatives.
\end{proof}

\subsection{One solution on the full interval of positive scales}
\label{exlim:terminal}

The energy and weighted spatial bounds are uniform in $a$.
We apply compactness on nested intervals bounded away from
$\lambda=0$ and choose a subsequence whose limits agree on overlaps.
Weak lower semicontinuity preserves the high-order bounds, while
strong convergence preserves the lower-order bounds and their
continuity in time.

\begin{proposition}
\label{exlim:terminal-diagonal}
There is an inviscid solution $R$ on $0<\lambda\leq A_0$,
on the same mass interval, with the stated positivity bounds in
normalized coordinates and consecutive distributional derivatives.
It has the following regularity and decay.
\begin{enumerate}[label=(\roman*),leftmargin=2em]
\item On every compact interval of positive scales,
\begin{equation}
 \Psi_j\in C_sX^q(s)\quad(j+q\leq22),\qquad
 \Psi_j\in L^\infty_sX^{23-j}\quad(j\leq23).
 \label{exlim:terminal-regularity}
\end{equation}
\item The total-order-$23$ bound \eqref{exlim:full-decay} holds
almost everywhere. The continuous lower-order derivatives satisfy
\begin{equation}
 \left(\sum_{q=0}^8\sum_{j=0}^{22-q}
                         \|\Psi_j(s)\|_{X^q(s)}^2\right)^{1/2}
                   \leq C_\beta\lambda(s)^{J+1/2}
 \label{exlim:terminal-low-decay}
\end{equation}
at every positive scale.
\end{enumerate}
The construction uses subsequential convergence and makes no
uniqueness assertion.
\end{proposition}

\begin{proof}
\emph{Compactness.}
Set $a_n=A_02^{-n}$ and take the solutions of
Proposition~\ref{exlim:fixed-terminal}, all with the same profile
and mass interval. Fix $0<a<A_0$ and discard the finitely many $n$
for which $a_n>a$. On $I_a$, for $j+q\le22$,
\[
 \sup_n\bigl(\|\Psi_j^n\|_{L^\infty F_a^{q+1}}
             +\|\Psi_{j+1}^n\|_{L^\infty F_a^q}\bigr)\le C_a,
\]
and
\begin{equation}
 \partial_s\Psi_j^n=\Psi_{j+1}^n-4b\Psi_j^n.
 \label{exlim:normalized-row-identity}
\end{equation}
Lemma~\ref{exlim:time-compactness} and the fixed-scale norm
equivalence give a common subsequence such that
\begin{equation}
 \Psi_j^n\longrightarrow\Psi_j
            \quad\hbox{in }C(I_a;X^q(s)),\qquad j+q\leq22.
 \label{exlim:terminal-strong}
\end{equation}
The total-order-$23$ bound also gives
\[
 \Psi_j^n\rightharpoonup^*\Psi_j
       \quad\hbox{in }L^\infty(I_a;F_a^{23-j}),\qquad 0\le j\le23.
\]
These limits remain in the closed compatible spaces of
Lemma~\ref{exlim:spatial-compactness}. For
$\phi\in C_c^\infty(I_a;H)$, passing to the limit in
\eqref{exlim:normalized-row-identity} yields
\[
 \int_{I_a}(\Psi_j,\phi_s)_H\dd s
 =-\int_{I_a}(\Psi_{j+1}-4b\Psi_j,\phi)_H\dd s,\qquad j\le22.
\]
At $j=22$ the next derivative converges weak-star; all other
terms converge strongly. Thus $R_j=A_j+\lambda^4\Psi_j$
are consecutive distributional derivatives.

\emph{Equation and boundary conditions.}
The strong convergence controls more than two ordinary spatial
derivatives on compact interior intervals and controls $R_2^n$.
In \eqref{exlim:force} all denominators stay positive, so
\begin{equation}
 R_{ss}-\tfrac32bR_s+\lambda^3\mathcal M[R]=0.
 \label{exlim:terminal-equation}
\end{equation}
At vacuum, Lemma~\ref{exlim:endpoint} preserves the radius and
its regular quotients; at the center,
\eqref{exlim:center-force} preserves the deformation and determinant.
Their positive lower bounds persist. For a smooth compatible test,
the pressure flux satisfies
\[
 R^2P(\rho)f=
 \begin{cases}
 O_a(\xi^5),&\xi\to0 \quad\hbox{at vacuum},\\
 O_a(z_c^3),&z_c\to0 \quad\hbox{at the center}.
 \end{cases}
\]
Both boundary terms vanish. Density and the two-form bound extend
the equation to the closed form domain.

\emph{Diagonal limit and bounds.}
On $I_\ell=\{s:a_\ell\le\lambda(s)\le A_0\}$ choose nested
subsequences $N_{\ell+1}\subset N_\ell\subset\{n:n\ge\ell\}$.
Uniqueness of the strong limit on overlaps gives
\[
 R^{[\ell+1]}|_{I_\ell}=R^{[\ell]},\qquad
 R|_{I_\ell}:=R^{[\ell]}.
\]
An increasing diagonal $n_\ell\in N_\ell$ converges on every $I_k$.
For $0\le\chi\in C_c^\infty(I_a)$, lower semicontinuity in the
original, scale-dependent norms gives
\[
 \int_{I_a}\chi\mathcal X_{23}\dd s
 \le\liminf_n\int_{I_a}\chi\mathcal X_{23}^n\dd s
 \le C_\beta\int_{I_a}\chi\lambda^{J+1/2}\dd s.
\]
This proves the total-order-$23$ bound almost everywhere.
For the continuous lower derivatives, \eqref{exlim:terminal-strong}
gives at every $s\in I_a$
\[
 \left(\sum_{q=0}^8\sum_{j=0}^{22-q}
             \|\Psi_j(s)\|_{X^q(s)}^2\right)^{1/2}
 \le\liminf_n\mathcal X_{23}^n(s)
 \le C_\beta\lambda(s)^{J+1/2}.
\]
The constants for fixed-scale compactness do not enter these
inequalities. Strong convergence also preserves the trace at $A_0$.
\end{proof}

At total order $23$, the terminal limit is only known to have
essentially bounded norms in time. The physical conclusions below
use lower derivatives, which are continuous in time, and therefore
hold without exceptional times. The required orders $(j,q)$ are:
\begin{itemize}[leftmargin=2em]
\item The physical energy comparison uses $(0,1)$ and $(1,0)$,
and $(0,8)$ for positivity along the radius segment.
\item Core density and relative velocity use $(0,8)$ and $(1,8)$,
together with \eqref{exlim:compact-core-embedding}.
\item The boundary radius and concentration use $(0,8)$ and the fixed mass.
\item The boundary-layer limit in $C^3([0,Y])$ uses $(0,6)$ and the
fixed-$Y$ estimate \eqref{exlim:layer-native-loss}.
\item The regular layer factors and vacuum slope use the bounded
quotient maps and the spatial endpoint trace of the same $(0,6)$
representative.
\end{itemize}

\subsection{Classical physical-vacuum regularity}
\label{exlim:physical-class}

Define physical time and the collapsing radius by
\begin{equation}
 t=-\tau(\lambda),\qquad
 \tau(\lambda)=\int_0^\lambda
             \frac{u^{1/2}}{\sqrt{\beta^2+eu}}\dd u,\qquad
 r(t,x)=R(s(\lambda(t)),x).
 \label{exlim:physical-time}
\end{equation}
The coefficients relating physical and rescaled time derivatives are bounded
on compact negative-time intervals.

\begin{proposition}
\label{exlim:classical}
The radius $r$ is an order-$22$ classical radial
physical-vacuum solution. 
\begin{enumerate}[label=(\roman*),leftmargin=2em]
\item At the center,
$r/z_c,\rho_L,u_L/z_c$ have even one-sided extensions and
bounded mixed ordinary derivatives through total order $16$.
\item At the vacuum,
\begin{equation}
 r(t,x)=R_b(t)-\xi^2G(t,\xi),\qquad
 \rho_L(t,x)=\xi^3\Theta(t,\xi),\qquad
 \xi=(M-x)^{1/5},
 \label{exlim:pv-factorization}
\end{equation}
where $G,\Theta$ are positive, and they and their reciprocals
have the same mixed derivative bounds.
\item The Eulerian enthalpy
is $C^1$ in physical radius up to the boundary and has strictly
negative outward derivative there.
\end{enumerate}
\end{proposition}

\begin{proof}
\emph{Vacuum.}
On a fixed collar, Lemma~\ref{exlim:endpoint} gives
\[
 \begin{gathered}
 B_5^{k+5}\hookrightarrow C^{k+2},\\
 \left\|(f-f(0))/\xi^2\right\|_{C^k}
       +\|f_\xi/\xi\|_{C^k}
 \le C_k\|f\|_{B_5^{k+5}},\qquad f_\xi(0)=0.
 \end{gathered}
\]
For $j+k\le16$,
\begin{equation}
 j+(k+5)\leq21,\qquad (j+1)+(k+5)\leq22
 \label{exlim:mixed-count}
\end{equation}
Thus both consecutive time derivatives are continuous in the
required spatial space. Integration in time permits differentiation
of the bounded quotient maps. Since $r_\xi(t,0)=0$,
\[
 G=\frac{R_b-r}{\xi^2},\qquad
 Q=-\frac{r_\xi}{\xi}=2G+\xi G_\xi,\qquad
 \Theta=\frac5{4\pi r^2Q},\qquad G(t,0)=\tfrac12Q(t,0)>0.
\]
On a sufficiently small collar, these factors and their reciprocals
are uniformly bounded on compact time intervals, with the asserted
mixed derivatives.

\emph{Center and interior.}
For $\Phi(t,X)=r(t,|X|^3)X/|X|$, localize on an enlarged center ball and use
\[
 \begin{gathered}
 H^{k+4}(\mathbb R^3)\hookrightarrow C^{k+2},\qquad
 \int_{\mathbb R^3}|\zeta|^{2(k+2)}
        (1+|\zeta|^2)^{-k-4}\dd\zeta<\infty,\\
 \frac r{z_c}=\int_0^1\theta^2\operatorname{div}\Phi(t,\theta X)\dd\theta,\qquad
 \rho_L=\frac3{4\pi\det D\Phi}.
 \end{gathered}
\]
The quotient costs one spatial derivative; velocity costs one time
derivative. The available orders are
\[
 (j+1)+(k+4)\le21,\qquad (j+2)+(k+4)\le22
                   \quad(j+k\le16).
\]
Equivariance gives even scalar quotients along each diameter.
In the interior, $H^{k+1}\hookrightarrow C^k$ and the formulas
for density and velocity use total order at most $18$.
Positivity of $r,r_x$ permits reciprocal differentiation and the
change to Eulerian coordinates.

Writing $u_L=r_t$, the mass identity and
\eqref{exlim:terminal-equation} give
\[
 \partial_t\rho_L=-\rho_L(2r_t/r+r_{tx}/r_x),\qquad
 r_{tt}=-4\pi r^2\partial_xP(\rho_L)-x/r^2.
\]
If $x_E$ is the inverse mass coordinate, then
\[
 \dd x=4\pi\rho_Er^2\dd r,\qquad
 \Phi_{E,r}=\frac{x_E}{r^2},\qquad
 (r^2\Phi_{E,r})_r=4\pi r^2\rho_E.
\]
These are the Eulerian equations, with total mass $M$ and exterior
field $M/r^2$. All identities are pointwise by the established
ordinary regularity.

\emph{Physical vacuum.}
The exact enthalpy is
\begin{equation}
 h_L=4\bigl(\sqrt{1+\xi^2\Theta^{2/3}}-1\bigr)
   =\xi^2\,\frac{4\Theta^{2/3}}
                  {\sqrt{1+\xi^2\Theta^{2/3}}+1}.
 \label{exlim:enthalpy-factor}
\end{equation}
The regular coefficient of $\xi^2$ has mixed derivatives through
order $16$. Dividing $\partial_\xi h_L$ by $r_\xi=-\xi Q$ gives
\begin{equation}
 -\partial_rh_E(t,R_b(t))
       =\frac{4\Theta(t,0)^{2/3}}{Q(t,0)}>0.
 \label{exlim:enthalpy-slope}
\end{equation}
The quotient is continuous at $\xi=0$, which proves the $C^1$
Eulerian extension. No even extension at vacuum is used.
\end{proof}

\subsection{Conserved energy and the scale parameter}
\label{exlim:physical-energy}

The leading kinetic, internal, and gravitational energies are of
order $\lambda^{-1}$, so core convergence alone does not determine
the finite part of their sum. We first prove conservation in the
classical endpoint class. Lemma~\ref{exlim:global-profile-energy}
then compares the full matched profile with the homologous solution,
including the joining region; the small mass of this region alone
does not control its energy.

With the vacuum normalization \eqref{eq:internal-energy-normalization},
the linear term $-4\rho$ contributes $-4M$. The global comparison
bounds the remaining equation-of-state error and the profile
corrections. The continuous lower-order bounds show that the
solution and profile energies have the same limit, which yields the
identity in Proposition~\ref{exlim:clock-energy} at the same
$(\delta,e)$.

\begin{proposition}
\label{exlim:energy-conservation}
Every solution in Definition~\ref{def:global-solution}
has finite conserved energy on its classical interval, with
\begin{equation}
 \mathscr E(t)=\int_0^M
 \left\{\frac12u_L(t,x)^2+\mathfrak e(\rho_L(t,x))
                                  -\frac{x}{r(t,x)}\right\}\dd x.
 \label{exlim:conserved-energy}
\end{equation}
Here $\mathfrak e$ is the specific internal energy defined in
\eqref{eq:internal-energy-normalization}, with its continuous vacuum
value zero.
\end{proposition}

\begin{proof}
Fix a compact time interval $I$ before collapse. At the vacuum,
\[
 r=R_b-\xi^2G,\qquad \rho_L=\xi^3\Theta,\qquad
 u_L=R_b'-\xi^2G_t.
\]
The boundary radius has two bounded time derivatives: for any fixed
$\xi_*>0$ in the collar,
$R_b(t)=r(t,M-\xi_*^5)+\xi_*^2G(t,\xi_*)$.
The endpoint class and the pressure law give, uniformly on $I$,
\[
 \begin{gathered}
 |u_L|+|\partial_tu_L|\le C_I,\qquad
 P(\rho_L)=O_I(\xi^5),\qquad
 \mathfrak e(\rho_L)=O_I(\xi^2),\\
 \partial_t\mathfrak e(\rho_L)
   =\mathfrak e'(\rho_L)\partial_t\rho_L
   =O_I(\xi^{-1})O_I(\xi^3)=O_I(\xi^2).
 \end{gathered}
\]
At the center the positive deformation gives
\[
 \begin{gathered}
 c_Iz_c\le r\le C_Iz_c,\qquad |u_L|+|\partial_tu_L|\le C_Iz_c,
 \qquad c_I\le\rho_L\le C_I,\\
 \frac{x}{r}+\frac{x|u_L|}{r^2}=O_I(z_c^2).
 \end{gathered}
\]
These bounds make the energy density and its time derivative integrable.

At fixed mass, $\rho_L=(4\pi r^2r_x)^{-1}$ and $r_t=u_L$ imply
\[
 \frac{\partial_t\rho_L}{\rho_L}
       =-\frac{2u_L}{r}-\frac{(u_L)_x}{r_x},\qquad
 \partial_t\mathfrak e(\rho_L)
       =\frac{P(\rho_L)}{\rho_L^2}\partial_t\rho_L
       =-4\pi P(\rho_L)\partial_x(r^2u_L).
\]
Multiplying the momentum equation by $u_L$ now gives
\[
 \begin{aligned}
 \partial_t\{\tfrac12u_L^2+\mathfrak e(\rho_L)-x/r\}
 &=-4\pi r^2u_L\partial_xP(\rho_L)-\frac{xu_L}{r^2}\\
 &\quad-4\pi P(\rho_L)\partial_x(r^2u_L)+\frac{xu_L}{r^2}\\
 &=-4\pi\partial_x\{r^2u_LP(\rho_L)\}.
 \end{aligned}
\]
For $0<x_-<x_+<M$ and $t_1,t_2\in I$, integration gives
\[
 \left[\int_{x_-}^{x_+}
   \{\tfrac12u_L^2+\mathfrak e(\rho_L)-x/r\}\dd x\right]_{t_1}^{t_2}
 =-4\pi\int_{t_1}^{t_2}[r^2u_LP(\rho_L)]_{x_-}^{x_+}\dd t.
\]
The right side is bounded in absolute value by
$C_I|t_2-t_1|\{x_-+(M-x_+)\}$, since the endpoint fluxes are
$O_I(z_c^3)$ and $O_I(\xi^5)$, respectively. Letting
$x_-\downarrow0$ and $x_+\uparrow M$ proves conservation.
\end{proof}

\begin{lemma}
\label{exlim:global-profile-energy}
Fix the selected profile. Put
$I_2=\int_0^M z(x)^2\dd x$ and
$\mathscr E_A=\frac12\|A_\tau\|_H^2+\mathcal U[A]$.
For a finite
$p_{\log}$,
\begin{equation}
       \left|\mathscr E_A(\lambda)+4M-\frac e2I_2\right|
 \le C\lambda(1+|\log\lambda|)^{p_{\log}+1}.
 \label{exlim:profile-energy-limit}
\end{equation}
The comparison with the homologous profile uses the following
estimates.
\begin{enumerate}[label=(\roman*),leftmargin=2em]
\item Replacing the exact internal-energy density by
$3\rho^{4/3}-4\rho$ gives an integrated error $O(\lambda)$.
\item The integrated core error is
$O(\lambda(1+|\log\lambda|)^{p_{\log}+1})$.
\item The complementary tail contains the entire joining shell.
Its mass is $O(\lambda^2)$, its kinetic and gravitational energies
are $O(\lambda)$, and its polytropic internal energy is
$O(\lambda^{3/2})$.
\end{enumerate}
All constants are independent of the terminal scale and viscosity.
\end{lemma}

\begin{proof}
\emph{The equation of state.}
Define
\[
 \mathscr E_A^{\rm pol}
 =\frac12\|A_\tau\|_H^2
       +\int_0^M\{3\rho_A^{1/3}-x/A\}\dd x.
\]
By \eqref{eq:threshold-exact-law-defect}, H\"older's inequality,
and \eqref{mat:boundary-profile-export},
\begin{equation}
 0\le\mathscr E_A+4M-\mathscr E_A^{\rm pol}
 \le3\int_{B_{A(\lambda,M)}}\rho_A^{2/3}\dd y
 \le3M^{2/3}|B_{A(\lambda,M)}|^{1/3}\le C\lambda.
 \label{exlim:energy-eos-defect}
\end{equation}
For $A^{\rm hom}=\lambda z(x)$, the density is
$\lambda^{-3}w(z)^3$. The virial identity \eqref{mat:seed-59} gives
\begin{equation}
 \mathscr E_{\rm hom}^{\rm pol}
 =\frac{\beta^2 I_2}{2\lambda}+\frac e2 I_2
          +\frac{3I_{4/3}-\mathcal W(\rho_\delta)}\lambda
 =\frac e2 I_2,
 \qquad \delta=-\frac{\beta^2}{2}.
 \label{exlim:homologous-energy}
\end{equation}
We estimate the difference between the two polytropic energies.
All constants below may depend on the fixed $(\beta,e)$, matching
order and cutoffs, but not on $\lambda$, terminal scale or viscosity.

\emph{The core.}
Fix $c>c_2$, with $c_2$ from \eqref{mat:matched-radius}.
On $q=1-z\ge c\sqrt\lambda$, the profile is the finite core sum.
Write
\[
 A=\lambda z(1+v),\qquad d=v+zv_z,\qquad
 A_\tau=b\lambda^{-1/2}z(1+v+\lambda\partial_\lambda|_zv).
\]
The coefficient bound
\eqref{mat:native-eq:wave2-core-coefficient-growth} gives
\begin{equation}
 \begin{split}
 |v|+|\lambda\partial_\lambda|_zv|
 &\le C\sum_{n=2}^K\lambda^n q^{1-n}\mathcal L(\lambda,q)^{p_{\log}},\\
 |d|&\le C\sum_{n=2}^K\lambda^n q^{-n}\mathcal L(\lambda,q)^{p_{\log}},
 \qquad
 \mathcal L(\lambda,q)=1+|\log\lambda|+|\log q|.
 \end{split}
 \label{exlim:energy-core-budget}
\end{equation}
On the compact core the same quantities are
$O(\lambda^2(1+|\log\lambda|)^{p_{\log}})$, including the
regular center quotients. Hence $v,d=o(1)$ uniformly on the region
under consideration. Using
\[
 \rho_A^{1/3}=\lambda^{-1}w(1+v)^{-2/3}(1+d)^{-1/3},
 \qquad A_\tau^{\rm hom}=b\lambda^{-1/2}z,
\]
we obtain
\[
 \begin{aligned}
 |A_\tau^2-(A_\tau^{\rm hom})^2|
   &\le C\lambda^{-1}z^2(|v|+|\lambda\partial_\lambda|_zv|),\\
 |x/A-x/(\lambda z)|&\le C\lambda^{-1}(x/z)|v|,\\
 |\rho_A^{1/3}-\lambda^{-1}w|
   &\le C\lambda^{-1}w(|v|+|d|).
 \end{aligned}
\]
Since $|\dd x|\sim q^3\dd q$ and $w\sim q$ on the boundary
collar, each integrated difference is bounded by the compact-core
error $C\lambda(1+|\log\lambda|)^{p_{\log}}$ plus
\begin{equation}
 C\sum_{n=2}^K\lambda^{n-1}
       \int_{c\sqrt\lambda}^{q_0}
                      q^{4-n}\mathcal L(\lambda,q)^{p_{\log}}\dd q
       \le C\lambda(1+|\log\lambda|)^{p_{\log}+1}.
 \label{exlim:energy-integrated-core}
\end{equation}
Indeed,
\[
 \int_{c\sqrt\lambda}^{q_0}q^{4-n}\dd q
 \le C_n\begin{cases}
 1,&2\le n\le4,\\
 1+|\log\lambda|,&n=5,\\
 \lambda^{(5-n)/2},&n\ge6,
 \end{cases}
\]
and $\lambda^{n-1}\lambda^{(5-n)/2}=\lambda^{(n+3)/2}$.

\emph{The joining region and vacuum.}
The complementary region $q\le c\sqrt\lambda$ contains the
joining shell and has mass
\[
 \int_{q\le c\sqrt\lambda}\dd x
 =4\pi\int_{1-c\sqrt\lambda}^1w(z)^3z^2\dd z
 \le C\int_0^{c\sqrt\lambda}q^3\dd q\le C\lambda^2.
\]
By \eqref{mat:strain-estimate} and \eqref{pr:profile-geometry},
\[
 A\sim\lambda,\qquad |A_\tau|\le C\lambda^{-1/2},\qquad
 \int_{q\le c\sqrt\lambda}(A_\tau^2+x/A)\dd x\le C\lambda.
\]
On the bounded fifth-root chart, $\rho_A\le C$; on the rest of
this region,
\[
 y_8\le C\lambda^{-1/4},\qquad
 \rho_A\le Cy_8^6\le C\lambda^{-3/2},\qquad
 \int_{q\le c\sqrt\lambda}\rho_A^{1/3}\dd x
       \le C\lambda^{-1/2}\lambda^2=C\lambda^{3/2}.
\]
The homologous radius satisfies the same tail bounds, since
$\rho_{\rm hom}^{1/3}=\lambda^{-1}w\le C\lambda^{-1/2}$ there.
Adding the core and tail estimates to
\eqref{exlim:energy-eos-defect} proves
\eqref{exlim:profile-energy-limit}.
\end{proof}

\Needspace{11\baselineskip}
\begin{proposition}
\label{exlim:clock-energy}
For the collapsing solution constructed above, put
\[
 I_2(\delta)=\int_0^M z(x)^2\dd x
               =4\pi\int_0^1w_\delta(z)^3z^4\dd z.
\]
With the normalization in \eqref{exlim:conserved-energy}, its energy is
\begin{equation}
                  \mathscr E+4M=\frac e2 I_2(\delta).
 \label{exlim:clock-energy-identity}
\end{equation}
\end{proposition}

\begin{proof}
Put $Y=R-A$, and let $G(\lambda)$ be the continuous lower-order
bound in \eqref{exlim:terminal-low-decay}, so that
$G\le C_\beta\lambda^{J+1/2}$.
We use only $(j,q)=(0,1),(1,0),(0,8)$ from that bound.
The last pair and \eqref{en:low-strains} preserve the positive
normalized radius and Jacobian along $A+\theta Y$, $0\le\theta\le1$.
The first two pairs, together with
$Y_\tau=\lambda^{-3/2}Y_1=\lambda^{5/2}\Psi_1$, give
\begin{equation}
 \|Y\|_H\le C\lambda^4G,\qquad
 \|Y_\tau\|_H\le C\lambda^{5/2}G,\qquad
 \|A_\tau\|_H\le C\lambda^{-1/2}.
\label{exlim:energy-nonlinear-norms}
\end{equation}
Thus the kinetic energies satisfy
\[
 \begin{split}
 \frac12\left|\|R_\tau\|_H^2-\|A_\tau\|_H^2\right|
 &\le \|A_\tau\|_H\|Y_\tau\|_H+\frac12\|Y_\tau\|_H^2\\
 &\le C(\lambda^2G+\lambda^5G^2).
 \end{split}
\]
For the potential, Proposition~\ref{mat:matched-family} at material
and spatial order zero and \eqref{pr:upper-pressure-bound} give
\[
 \|\mathcal M[A]\|_H\le C\lambda^{-2},\qquad
 V(Y)\le C\lambda^{7/2}G.
\]
Here we use only the bound
$\|g\|_H\le C\|g\|_{0,a}$ from Proposition~\ref{mat:matched-family}.
The ordinary derivative norm $\|\cdot\|_{0,a}$ imposes no compatible
vacuum trace.
Taylor's formula and \eqref{en:two-form-bound} now yield
\[
 \begin{split}
 |\mathcal U[R]-\mathcal U[A]|
 &\le |(\mathcal M[A],Y)_H|
       +\int_0^1(1-\theta)
              |D^2\mathcal U[A+\theta Y][Y,Y]|\dd\theta\\
 &\le C\lambda^2G+C\lambda^{-3}V(Y)^2
 \le C(\lambda^2G+\lambda^4G^2).
 \end{split}
\]
Combining the two bounds gives
\[
 |\mathscr E-\mathscr E_A|
 \le C(\lambda^2G+\lambda^4G^2)
 \le C_\beta\left(\lambda^{J+5/2}+\lambda^{2J+5}\right)
 \longrightarrow0.
\]
Every norm used here is in the continuous range, so this limit holds
at every positive scale. The constants come from the global covering
and pressure estimates and are independent of the fixed-scale
regularity constants. Lemma~\ref{exlim:global-profile-energy}
identifies the limit of $\mathscr E+4M$ as $eI_2(\delta)/2$;
energy conservation from Proposition~\ref{exlim:energy-conservation}
then proves \eqref{exlim:clock-energy-identity}.
\end{proof}

\subsection{The collapse asymptotics}
\label{exlim:asymptotics}

We now apply Corollary~\ref{mat:matched-limits}. Let $x=m(z)$, $0\leq z\leq1$,
be the leading mass coordinate with leading density profile $w(z)^3$.
The leading profile is positive for $z<1$, has the regular even center
structure, and its mass satisfies $M>M_{\rm Ch}$.
On each compact core $0\leq z\leq z_1<1$, the matched profile is
eventually the uncut finite expansion
\begin{equation}
 A(s,m(z))=\lambda z\left\{1+
                 \sum_{n=2}^{K}\lambda^n\Phi_n(\log\lambda,z)\right\}.
 \label{exlim:profile-core}
\end{equation}
Its needed finite coefficient derivatives are bounded by finite
logarithmic polynomials. At the boundary,
\begin{equation}
 A(s,M)=\lambda\left\{1-\frac{\lambda}{\varkappa_\delta}
                 +O\bigl(\lambda^2(1+|\log\lambda|)^{p_{\log}}\bigr)\right\},
 \qquad \varkappa_\delta>0.
 \label{exlim:profile-boundary}
\end{equation}
The coefficient in this expression is the leading layer value
at its vacuum endpoint; it is not inferred from a core expansion
outside its domain.

\begin{lemma}
\label{exlim:compact-core-ordinary}
For $z_1<1$ and every compatible displacement $f\in X^8(s)$,
\begin{equation}
 \left\|\frac{f(m(z))}{z}\right\|_{C^1([0,z_1])}
        \le C_{z_1,\beta}\lambda^{-3}\|f\|_{X^8(s)},
 \qquad 0<\lambda\le A_0\le1.
 \label{exlim:compact-core-embedding}
\end{equation}
The quotient has its removable center value. The constant depends
on the selected leading profile, compact core, and fixed covering cutoffs, but
not on the terminal scale, viscosity, or sufficiently small $\lambda$.
\end{lemma}

\begin{proof}
Use the exact partition $f=\sum_{\nu=c,i,5,8}\chi_\nu f$.
Near the center $\chi_c=1$. The equivariant lift obeys
\[
 \frac{f(x)}{x^{1/3}}
   =\int_0^1 t^2\operatorname{div}U_f(tX)\dd t,
 \qquad |X|=x^{1/3}.
\]
The localized $H^8$ norm controls three Cartesian derivatives.
Differentiating the integral through order two gives, on a smaller
center interval $0\le z_c\le c$,
\[
 \left\|\frac{f(z_c^3)}{z_c}\right\|_{C^2([0,c])}
 \le C\|\chi_cU_f\|_{H^8}
 \le C\|f\|_{X^8}.
\]
The leading mass map is a smooth odd change of variables near zero,
and $m(z)^{1/3}/z$ is smooth, positive, and even. In
\[
 \frac{f(m(z))}{z}
   =\frac{f(m(z))}{m(z)^{1/3}}\,
      \frac{m(z)^{1/3}}{z},
\]
both factors therefore have the required $C^1$ bounds. On the
remaining fixed center and interior charts, $z$ is bounded away
from zero and ordinary Sobolev embedding applies.

For small $\lambda$, the bounded fifth-root chart does not meet
$[0,z_1]$, since $M-m(z_1)>0$. The only moving contribution is the
long chart. On its intersection with this compact core,
\[
 y_8=\lambda^{-1/2}(M-m(z))^{1/8},\qquad
 |\partial_z y_8|+|\partial_z^2y_8|
                 \le C_{z_1,\beta}\lambda^{-1/2}.
\]
On each interval of unit length in the chart coordinate the weight $y_8^7$ has a fixed positive
lower bound, so the ordinary $H^4\to C^2$ estimate and the edge
prefactor give
\[
 \|\chi_8 f\|_{C^2(\text{unit interval})}
            \le C\lambda^{-2}\|f\|_{X^8}.
\]
The constant is uniform over these unit intervals; taking their
supremum does not introduce their number. With $h=\chi_8f$, the
chain rule gives
\[
 \partial_z(h(y_8(z)))=h_y\partial_zy_8,\qquad
 \partial_z^2(h(y_8(z)))
    =h_{yy}(\partial_zy_8)^2+h_y\partial_z^2y_8.
\]
Consequently the bound above yields
\[
 \|h(y_8(\cdot))\|_{C^2}
       \le C_{z_1,\beta}\lambda^{-3}\|f\|_{X^8}.
\]
On this part of the compact core, $z$ has a fixed positive lower
bound, so division by $z$ preserves the estimate.

All derivatives fall on the localized function $\chi_8 f$, extended
by zero across artificial chart ends where it vanishes. The other
localized functions are controlled on the fixed charts. Summing
the four localized terms proves the bound for smooth compatible
functions; convergence in $X^8$ extends it to the completion.
Enlarging the constant covers any remaining compact interval of
positive $\lambda$ in the selected outer range.
\end{proof}

\begin{proposition}
\label{exlim:physical-asymptotics}
As $t\uparrow0$, the solution above has the following properties.
\begin{enumerate}[label=(\roman*),leftmargin=2em]
\item \emph{Scale.}
\begin{equation}
 \lambda(t)=\left(\frac{3\beta}{2}(-t)\right)^{2/3}
                                \{1+O((-t)^{2/3})\}.
 \label{exlim:scale-asymptotic}
\end{equation}
\item \emph{Interior density and velocity.} For every $0<z_0<1$,
\begin{align}
 \rho_E(t,r)&=\lambda^{-3}w(r/\lambda)^3\{1+o(1)\}
       \quad\hbox{uniformly for }0\leq r/\lambda\leq z_0,
 \label{exlim:density-asymptotic}\\
 \sup_{0<r\leq z_0\lambda}
 \left|\frac{u_E(t,r)}{(\dot\lambda/\lambda)r}-1\right|
                                      &\longrightarrow0.
 \label{exlim:velocity-asymptotic}
\end{align}
At the center,
\begin{equation}
 \rho_E(t,0)=\frac{4w(0)^3}{9\beta^2}(-t)^{-2}\{1+o(1)\}.
 \label{exlim:center-asymptotic}
\end{equation}
\item \emph{Support and mass concentration.}
\begin{equation}
 R_b(t)=\lambda\left\{1-\frac{\lambda}{\varkappa_\delta}
                  +O\bigl(\lambda^2(1+|\log\lambda|)^{p_{\log}}\bigr)\right\}.
 \label{exlim:boundary-asymptotic}
\end{equation}
The density, extended by zero outside its support, satisfies
$\rho_E(t,|\cdot|)\dd y\rightharpoonup^* M\delta_0$ as finite
Radon measures.
\end{enumerate}
\end{proposition}

\begin{proof}
\emph{The scale.}
Taylor's formula in the exact integral \eqref{exlim:physical-time}
gives
\[
 \tau(\lambda)=\frac2{3\beta}\lambda^{3/2}
       \left\{1-\frac{3e}{10\beta^2}\lambda+O(\lambda^2)\right\}.
 \]
Monotone inversion proves \eqref{exlim:scale-asymptotic}, and
differentiation of \eqref{eq:main-collapse-clock} gives
$\dot\lambda=-b\lambda^{-1/2}$.

\emph{The interior density and velocity.}
Write $r(t,m(z))=\lambda zU(t,z)$.
Apply Lemma~\ref{exlim:compact-core-ordinary} to the continuous
norms at $(0,8)$ and $(1,8)$. It gives
\begin{equation}
 \sum_{i=0}^{1}\|\Psi_i(s,m(\cdot))/z\|_{C^1([0,z_1])}
       \le C_{z_1,\beta}\lambda^{-3}\lambda^{J+1/2}.
 \label{exlim:core-perturbation}
\end{equation}
Equations~\eqref{exlim:profile-core} and
\eqref{exlim:core-perturbation} imply
\begin{equation}
 U\longrightarrow1\quad\hbox{in }C^1([0,z_1]),\qquad
 \frac{R_s(s,m(z))}{b\lambda z}\longrightarrow1
       \quad\hbox{uniformly for }0<z\leq z_1.
 \label{exlim:normalized-core}
\end{equation}
For the velocity perturbation, the bound is
$C_{z_1,\beta}\lambda^3\lambda^{-3}\lambda^{J+1/2}/b$,
which tends to zero because $b\geq\beta>0$.
The quotients extend continuously to the center.

The density formula in the leading mass coordinate is exact:
\begin{equation}
 \rho_L(t,m(z))
  =\lambda^{-3}\frac{w(z)^3}{U(t,z)^2(U(t,z)+zU_z(t,z))}.
 \label{exlim:core-density-exact}
\end{equation}
Choose $z_0<z_1<1$. The map $z\mapsto r/\lambda=zU$
is a monotone $C^1$ perturbation of the identity on
$[0,z_1]$, and its image contains $[0,z_0]$ for small scale.
Its inverse converges uniformly to the identity. Positivity of
$w$ on the larger compact interval and
\eqref{exlim:core-density-exact} give
\eqref{exlim:density-asymptotic}.
Moreover
\[
 u_E=-\lambda^{-3/2}R_s,\qquad
 (\dot\lambda/\lambda)r=-b\lambda^{-3/2}r.
 \]
Their ratio is $R_s/(br)$, which tends uniformly to one by
\eqref{exlim:normalized-core}. This proves
\eqref{exlim:velocity-asymptotic}; its denominator is not used
at $r=0$. At the center, combine the removable value in
\eqref{exlim:core-density-exact} with the scale law to obtain
\eqref{exlim:center-asymptotic}.

\emph{The boundary radius and concentration.}
On the bounded vacuum chart, Lemma~\ref{exlim:endpoint} at
order eight gives an endpoint trace bound with factor
$\lambda^{-2}$, the inverse of the norm prefactor. Therefore
\begin{equation}
 |R(s,M)-A(s,M)|
 \leq C\lambda^4\lambda^{-2}\|\Psi_0\|_{X^8}
 \leq C_\beta\lambda^{J+5/2}=o(\lambda^3).
 \label{exlim:boundary-perturbation}
\end{equation}
Together with \eqref{exlim:profile-boundary}, this proves
\eqref{exlim:boundary-asymptotic}. In particular the entire
support shrinks to the origin, since every mass shell lies
between zero and $R_b(t)$.

For a continuous compactly supported test $\varphi$, exact
mass conservation and the support bound give
\[
 \left|\int_{\mathbb R^3}\varphi(y)\rho_E(t,|y|)\dd y
                                     -M\varphi(0)\right|
 \leq M\sup_{|y|\leq R_b(t)}|\varphi(y)-\varphi(0)|
 \longrightarrow0.
 \]
This proves the measure limit. The strict inequality
$M>M_{\rm Ch}$ is the leading profile conclusion of Proposition~\ref{mat:matched-family};
it does not follow merely from concentration.
\end{proof}

Restriction to $0<\lambda\le A$, $A\le A_0$, gives the
initial trace at $t=-\tau(A)$. All the preceding asymptotics
hold for this restriction of the same solution.

\subsection{The boundary layer of the solution}
\label{exlim:actual-layer}

We prove Corollary~\ref{cor:actual-layer} by comparing the solution
with the matched layer in its continuous low-order weighted norms.
The convergence holds on each fixed bounded rescaled interval
$[0,Y]$, including the vacuum endpoint; uniform convergence out to
the moving joining shell is not asserted. The endpoint values of the limiting layer
determine the physical-vacuum factors and the boundary-radius
correction.

In this subsection
$\Theta_\lambda$ denotes the density factor in the chart coordinate
$y=(M-x)^{1/5}\lambda^{-4/5}$; it is distinct from the physical
$\xi$-coordinate factor $\Theta(t,\xi)$ in
\eqref{exlim:pv-factorization}. All constants may depend on the fixed
selected family and on $Y$. They are independent of the terminal
subsequence and of sufficiently small $\lambda$.

\begin{lemma}
\label{exlim:fixed-layer-embedding}
Fix $0<Y<\infty$. For all sufficiently small $\lambda$ and every
$f\in X^6(s)$, let $f^{[\lambda]}(y)=f(M-\lambda^4y^5)$.
Then
\begin{equation}
 \|f^{[\lambda]}\|_{C^3([0,Y])}
       \le C_Y\lambda^{-2}\|f\|_{X^6(s)}.
 \label{exlim:layer-native-loss}
\end{equation}
The constant depends on $Y$, the fixed enlarged coordinate intervals, and the
prescribed cutoffs. It is independent of the terminal scale,
viscosity, and sufficiently small $\lambda$. The compatible
completion gives $(f^{[\lambda]})'(0)=0$.
\end{lemma}

\begin{proof}
Choose $0<c<Y$ with $(2c)^5<1/2$. For small $\lambda$,
$\chi_e=1$ on the layer and its fixed enlarged coordinate intervals,
and $\chi_5=1$ on $[0,2c]$. Lemma~\ref{exlim:endpoint} and
\eqref{at:graph-norm} give
\[
 \|f^{[\lambda]}\|_{C^3([0,c])}
 \le C\|\chi_5f\|_{B_5^6}
 \le C\lambda^{-2}\|f\|_{X^6(s)}.
\]
On $[c,Y]$, use $\chi_5+\chi_8=1$. Each weight is bounded
above and away from zero on the relevant fixed annuli. For the
second term the coordinate change is
\[
 y_8=y^{5/8},\qquad
 0<c_Y\le y,y_8\le C_Y
       \quad\hbox{on its support in }[c,Y].
\]
This map and its inverse have bounded derivatives there.
Applying $H^4\hookrightarrow C^3$ on finitely many enlarged
intervals yields
\[
 \begin{split}
 \|f^{[\lambda]}\|_{C^3([c,Y])}
 &\le C_Y\bigl(\|\chi_5f\|_{B_5^6}
                         +\|\chi_8f\|_{C_8^6}\bigr)\\
 &\le C_Y\lambda^{-2}\|f\|_{X^6(s)}.
 \end{split}
\]
For small $\lambda$, these intervals lie below the outer end
$C_e\lambda^{-1/2}$ of the long chart. All cutoff derivatives
are included in the norms. Only localized terms vanishing near
artificial chart ends are extended by zero.

The estimate proves convergence in $C^3([0,Y])$ for compatible
approximants. Their zero first trace passes to the limit, giving
$(f^{[\lambda]})'(0)=0$ and \eqref{exlim:layer-native-loss}.
\end{proof}

\begin{proof}[Proof of Corollary~\ref{cor:actual-layer}]
\emph{The leading profile.}
At degree zero in the exact layer force \eqref{mat:connection-9},
the normalized radius is $R=1$ and the inertial term is
$-\beta^2/2=\delta$. Thus
\[
 M+\delta-4\pi\partial_DP(\rho_0)=0.
\]
The leading profile identity \eqref{mat:connection-2} and the vacuum value
$P(\rho_0(0))=0$ give
$P(\rho_0(D))=\varkappa_\delta D/\pi$.
Since $h'(\rho)=P'(\rho)/\rho$,
\[
 \frac{\dd}{\dd D}\left(\frac{h(\rho_0(D))}{4\varkappa_\delta}\right)
       =\frac1{4\pi\rho_0(D)}=Z_0'(D),\qquad D>0.
\]
The integration constant is $1/\varkappa_\delta$, as fixed by
the normalization in \eqref{mat:connection-10}.

For the normalized pressure, direct integration gives
\begin{equation}
 P(q^3)=4\int_0^q\frac{u^4}{\sqrt{1+u^2}}\dd u
       =q^5p(q^2),\qquad
 p(v)=4\int_0^1\frac{a^4}{\sqrt{1+va^2}}\dd a,
 \quad p(0)=\frac45.
 \label{exlim:layer-pressure-unit}
\end{equation}
Consequently $\rho_0(y^5)=y^3\Theta_0(y)$, where near zero the
positive factor solves
\[
 \Theta_0^{5/3}p(y^2\Theta_0^{2/3})
                         =\frac{\varkappa_\delta}{\pi}.
\]
The derivative of the left side with respect to $\Theta_0$ is
positive at $y=0$. The ordinary implicit function theorem gives
a smooth positive function of $y^2$ there, with the value stated
in \eqref{eq:actual-layer-endpoint-values}. Away from zero,
$P'(\rho)>0$ gives smoothness and positivity by ordinary inversion.
In particular these properties hold on every fixed $[0,Y]$.
Furthermore,
\begin{equation}
 \mathcal Q_0(y)=\frac{5}{4\pi\Theta_0(y)},\qquad
 \mathcal G_0(y)=\int_0^1 a\mathcal Q_0(ay)\dd a.
 \label{exlim:leading-layer-quotients}
\end{equation}
They are positive and regular, including at zero, and
$\mathcal Z_0'(0)=0$. This proves all leading endpoint values.

\emph{Convergence on a fixed layer.}
To pass from the weighted norm to an ordinary bound on the rescaled
layer, put
\[
 f_\lambda(y)=\Psi_0(s,M-\lambda^4y^5),\qquad
 \mathcal Z_A(\lambda,y)
     =\frac{\lambda-A(s,M-\lambda^4y^5)}{\lambda^2}.
\]
Lemma~\ref{exlim:fixed-layer-embedding} gives
$\|f_\lambda\|_{C^3([0,Y])}\le
C_Y\lambda^{-2}\|\Psi_0(s)\|_{X^6(s)}$.

The displacement identity $R-A=\lambda^4\Psi_0$ now yields the
exact difference
\begin{equation}
 \mathcal Z_\lambda-\mathcal Z_A=-\lambda^2f_\lambda,
 \qquad
 \|\mathcal Z_\lambda-\mathcal Z_A\|_{C^3([0,Y])}
       \le C_{\beta,Y}\lambda^{J+1/2}.
 \label{exlim:layer-nonlinear-error}
\end{equation}
Here we used \eqref{exlim:terminal-low-decay} at $(j,q)=(0,6)$.
The factor $\lambda^2$ in the difference exactly cancels the
$\lambda^{-2}$ loss in the weighted embedding.

For fixed $Y$, the layer eventually lies in the region where the cutoff leaves the layer expansion unchanged in
\eqref{mat:connection-3}: its distance from the vacuum boundary in the leading coordinate is $1-z=O_Y(\lambda)$,
while the joining shell has $1-z\sim\sqrt\lambda$.
The finite sum there gives
\begin{equation}
 \|\mathcal Z_A(\lambda,\cdot)-\mathcal Z_0\|_{C^3([0,Y])}
       \le C_{\beta,Y}\lambda(1+|\log\lambda|)^{p_{\log}}
 \label{exlim:layer-profile-error}
\end{equation}
for a finite $p_{\log}$. Every nonleading term contains at least one
factor $\lambda$; its coefficient has the needed ordinary derivatives in the chart coordinates by Subsection~\ref{mat:native} and
Proposition~\ref{mat:matched-family}. The finite logarithmic degree
in this estimate may be enlarged. It concerns only spatial
derivatives at fixed $\lambda$.
Combining \eqref{exlim:layer-nonlinear-error} and
\eqref{exlim:layer-profile-error} proves the $C^3$ convergence,
with error bounded by
\begin{equation}
 C_{\beta,Y}\left\{
       \lambda(1+|\log\lambda|)^{p_{\log}}+\lambda^{J+1/2}\right\}.
 \label{exlim:layer-error-budget}
\end{equation}
The profile truncation and the nonlinear correction are distinct
terms in this bound.

\emph{Convergence of the regular factors.}
The compatible $X^6$ completion and the reference layer give
$\mathcal Z_\lambda'(0)=\mathcal Z_0'(0)=0$.
For $g=\mathcal Z_\lambda-\mathcal Z_0$, the quotient identities
of Lemma~\ref{at:regular-factor-calculus} are
\[
 \frac{g'(y)}y=\int_0^1g''(ty)\dd t,\qquad
 \frac{g(y)-g(0)}{y^2}=\int_0^1(1-t)g''(ty)\dd t.
\]
They extend to $y=0$, and differentiation under the integral gives
\[
 \|\mathcal Q_\lambda-\mathcal Q_0\|_{C^1([0,Y])}
 +\|\mathcal G_\lambda-\mathcal G_0\|_{C^1([0,Y])}
 \le C_Y\|g\|_{C^3([0,Y])}.
\]
Hence both factors converge with the error in
\eqref{exlim:layer-error-budget}. Their positive limiting minima
give uniform positive lower bounds for small $\lambda$.
Only the vanishing first derivative at zero is used; no evenness
condition is imposed at the vacuum.

For $y>0$, differentiate the exact mass label and radius:
\[
 \frac{\dd x}{\dd y}=-5\lambda^4y^4,
 \qquad
 \frac{\dd r}{\dd y}=-\lambda^2y\mathcal Q_\lambda(y),
 \qquad
 r_x=\frac{\mathcal Q_\lambda(y)}{5\lambda^2y^3}.
\]
The exact density identity consequently becomes
\begin{equation}
 \Theta_\lambda(y)
   =\frac{5}{4\pi(1-\lambda\mathcal Z_\lambda(y))^2
                           \mathcal Q_\lambda(y)}.
 \label{exlim:layer-density-exact}
\end{equation}
The right side extends positively in $C^1$ to zero and agrees with
the density factor in the ordinary endpoint class. Since
$1-\lambda\mathcal Z_\lambda\to1$ and
$\mathcal Q_\lambda\to\mathcal Q_0>0$ in $C^1$, we obtain
\[
 \begin{split}
 &\left\|\Theta_\lambda-\frac5{4\pi\mathcal Q_0}\right\|_{C^1([0,Y])}\\
 &\qquad\le C_{\beta,Y}\bigl(
   \lambda\|\mathcal Z_\lambda\|_{C^1([0,Y])}
    +\|\mathcal Q_\lambda-\mathcal Q_0\|_{C^1([0,Y])}\bigr).
 \end{split}
\]
By \eqref{exlim:leading-layer-quotients} the limit is $\Theta_0$.
This proves its uniform positive lower bound and the error estimate
\eqref{exlim:layer-error-budget} for all three factors.

On each compact interval of positive scales,
$\Psi_0\in C_sX^6(s)$ by \eqref{exlim:terminal-regularity}.
The bounded coordinate changes and the preceding embeddings give its
ordinary continuous representative in the layer. The estimates
therefore hold at every positive scale and give the stated limits as
$t\uparrow0$ for the solution already constructed, without further
subsequence extraction. The constants in \eqref{exlim:layer-native-loss}
and the decay bound are independent of the fixed-scale comparisons
used for compactness.

\emph{Layer mass and vacuum slope.}
The exact mass coordinate gives the mass
$\lambda^4Y^5$ and the thickness identity
\eqref{eq:actual-layer-thickness}. Finally, at each fixed time the
exact enthalpy in this coordinate is
\[
 h_L(t,M-\lambda^4y^5)
   =4\left(\sqrt{1+y^2\Theta_\lambda(y)^{2/3}}-1\right).
\]
Its derivative divided by $y$ tends to
$4\Theta_\lambda(0)^{2/3}$ as $y\downarrow0$.
Since $r_y=-\lambda^2y\mathcal Q_\lambda$, the ordinary
physical-vacuum trace therefore gives the exact relation
\begin{equation}
 \lambda^2[-\partial_rh_E(t,R_b(t))]
       =\frac{4\Theta_\lambda(0)^{2/3}}{\mathcal Q_\lambda(0)}.
 \label{exlim:layer-slope-exact}
\end{equation}
Taking the already proved endpoint limits and using
$\Theta_0(0)^{5/3}=5\varkappa_\delta/(4\pi)$ gives
\[
 \frac{4\Theta_0(0)^{2/3}}{\mathcal Q_0(0)}
       =\frac{16\pi}{5}\Theta_0(0)^{5/3}
       =4\varkappa_\delta=M+\delta.
\]
This proves \eqref{eq:actual-layer-slope} from spatial traces,
without differentiating the remainder in the boundary position in time.
\end{proof}

We now combine the preceding construction and asymptotics to prove the
main theorems.

\subsection{Proof of the main theorems}
\label{exlim:main-theorem-proof}

\begin{proof}[Proof of Theorem~\ref{thm:main-collapse}]
Fix $e\ge0$. Choose $\beta_*,J,K$ by
Lemma~\ref{ct:parameter-admissibility}, so that
\[
 p:=2J+1>\Gamma_0,\qquad
 4+\frac p2\ge\max_{q\le8}\Gamma_q^{\rm cov},\qquad
 J^\sharp=J+1713,\qquad K=2J+3512.
\]
For $0<\beta\le\beta_*$, set
\[
 \delta=-\beta^2/2,\qquad
 M=4\pi\int_0^1w_\delta^3z^2\dd z,\qquad
 \varkappa_\delta=-w_\delta'(1).
\]
Equations \eqref{mat:seed-51} and \eqref{mat:seed-60} give
$M+\delta=4\varkappa_\delta$ and $M>M_{\rm Ch}$.
All subsequent radii are functions on the same interval $[0,M]$.

\emph{Regularized solutions and the two limits.}
Proposition~\ref{mat:matched-family} gives the radius $A$.
Lemma~\ref{ct:parameter-admissibility} fixes the common upper bound for the scale
$A_0=A_0(\beta,e)>0$ and the permitted viscosity bound
$\kappa_{\rm tube}(a)>0$. At each terminal scale $a>0$ in this
range, Lemma~\ref{loc:finite-map} and Proposition~\ref{loc:selector}
give exact-mass compatible data. Proposition~\ref{ct:proposition}
then gives regularized solutions on $a\le\lambda\le A_0$ for
$0<\kappa\le\kappa_{\rm tube}(a)$. These solutions satisfy
\[
 0\le\widetilde E^{a,\kappa}\le\frac{d_0^2}{4}\lambda^p,
 \qquad G_8^{a,\kappa}\le\frac12\lambda^{p/2},\qquad
 \int_{s_a}^{s}\nu\mathcal D\dd\sigma\le C d_0^2\lambda(s)^p.
\]
The positivity bounds and the constants in these inequalities
are independent of $a$ and $\kappa$.

At fixed $a$, Proposition~\ref{exlim:fixed-terminal} removes
viscosity and gives
\[
 \mathcal X_{23}^{a}(s)\le C_\beta\lambda(s)^{J+1/2},\qquad
 \Psi_j^{a}\in C(I_a;X^{22-j})\quad(0\le j\le22).
\]
Apply Proposition~\ref{exlim:terminal-diagonal} to
$a_n=A_02^{-n}$. The diagonal subsequence defines one radius $R$ on
$0<\lambda\le A_0$ satisfying
\begin{equation}
 R_{ss}-\tfrac32bR_s+\lambda^3\mathcal M[R]=0,
 \qquad
 \Psi_j\in C_{\rm loc}(X^q),\quad j+q\le22,
 \label{exlim:main-selected-orbit}
\end{equation}
and
\[
 \left(\sum_{q=0}^8\sum_{j=0}^{22-q}
             \|\Psi_j(s)\|_{X^q(s)}^2\right)^{1/2}
 \le C_\beta\lambda(s)^{J+1/2}.
\]
The last estimate holds at every positive scale; the total-order-$23$
estimate holds almost everywhere. The limits defining $R$ on
overlapping scale intervals agree because the subsequences are
nested and the mass interval is fixed.

\emph{Physical time and regularity.}
Take $A_{\rm out}=A_0$, decreasing it if necessary, and define
\[
 t(\lambda)=-\int_0^\lambda
       \frac{\mu^{1/2}}{\sqrt{\beta^2+e\mu}}\dd\mu,
 \qquad
 s(\lambda)=\int_{A_{\rm out}}^\lambda
       \frac{\dd\mu}{\mu\sqrt{\beta^2+e\mu}}.
\]
Then $\dd t/\dd\lambda<0$ and
\[
 \frac{\dd s}{\dd t}=-\lambda^{-3/2},\qquad
 \dot\lambda=-b\lambda^{-1/2},\qquad
 r(t,x):=R(s(\lambda(t)),x),\qquad
 r_{tt}=\lambda^{-3}(R_{ss}-\tfrac32bR_s).
\]
Equation \eqref{exlim:main-selected-orbit} therefore gives
$r_{tt}+\mathcal M[r]=0$. Define
\[
 \rho_L=(4\pi r^2r_x)^{-1},\qquad
 u_L=-\lambda^{-3/2}R_s.
\]
Proposition~\ref{exlim:classical} gives the center and vacuum
regularity of Definition~\ref{def:global-solution}, with ordinary
mixed order $22-6=16$, positive interior density, and
$0<-\partial_rh_E(t,R_b(t))<\infty$. In particular,
\[
 4\pi\int_0^{R_b(t)}\rho_E(t,r)r^2\dd r
   =\int_0^M4\pi\rho_Lr^2r_x\dd x=M.
\]
The continuous trace at $\lambda=A_{\rm out}$ includes the left
endpoint $t=-T_{A_{\rm out}}$.

\emph{The collapse limits.}
The scale relation \eqref{eq:main-collapse-clock} gives
\[
 -t=\frac{2}{3\beta}\lambda^{3/2}
       \left(1-\frac{3e}{10\beta^2}\lambda+O(\lambda^2)\right),
 \qquad
 \lambda(t)=\left(\frac{3\beta}{2}(-t)\right)^{2/3}
                     (1+O((-t)^{2/3})).
\]
Write $r(t,m(z))=\lambda zU(t,z)$. By
\eqref{exlim:normalized-core}, on every $[0,z_1]$, $z_1<1$,
\[
 U\to1\ \hbox{in }C^1,
 \qquad \frac{R_s(s,m(z))}{b\lambda z}\to1
       \ \hbox{uniformly, with the removable center value}.
\]
The mass identity and the time change give
\[
 \lambda^3\rho_L(t,m(z))
   =\frac{w_\delta(z)^3}{U^2(U+zU_z)},\qquad
 \frac{u_L(t,m(z))}{(\dot\lambda/\lambda)r(t,m(z))}
   =\frac{R_s(s,m(z))}{b\lambda zU(t,z)}.
\]
For $z_0<z_1$, the maps $z\mapsto zU(t,z)$ have positive
derivative, cover $[0,z_0]$ at small scale, and converge with their
inverses to the identity. The two formulas prove
\eqref{eq:main-collapse-density}--\eqref{eq:main-collapse-velocity}.
Center parity gives $u_E(t,0)=0$, and the first formula yields
\[
 \rho_E(t,0)=\lambda^{-3}w_\delta(0)^3(1+o(1))
    =\frac{4w_\delta(0)^3}{9\beta^2}(-t)^{-2}(1+o(1)).
\]

At the vacuum, \eqref{exlim:boundary-perturbation} gives
$|R(s,M)-A(s,M)|\le C_\beta\lambda^{J+5/2}$.
Since $J\ge3$, the profile expansion implies
\[
 R_b(t)=\lambda\left\{1-\frac{\lambda}{\varkappa_\delta}
       +O\bigl(\lambda^2(1+|\log\lambda|)^{p_{\log}}\bigr)\right\}.
\]
Choose $C_{\log}=p_{\log}$ large enough for the finitely many
retained polynomials. For $\varphi\in C_c(\mathbb R^3)$,
\[
 \left|\int\varphi(X)\rho_E(t,|X|)\dd X-M\varphi(0)\right|
 \le M\sup_{|X|\le R_b(t)}|\varphi(X)-\varphi(0)|\longrightarrow0.
\]
This proves the measure convergence. The estimates above use only
the continuous lower-order bounds in \eqref{exlim:main-selected-orbit}.
Restricting the same solution to $[-T_A,0)$ gives the assertion for
every $0<A\le A_{\rm out}$.
\end{proof}

\begin{proof}[Proof of Theorem~\ref{thm:physical-main}]
Fix $e\ge0$ and $0<\beta_1\le\beta_*(e)$. The mass curve of the leading profiles satisfies
\[
 M_0=M_{\rm Ch},\qquad M_{\beta_1}>M_{\rm Ch},\qquad
 M_\beta=M_{\rm Ch}+\frac{I_2(0)}{2M_{\rm Ch}}\beta^2+O(\beta^4).
\]
Set $\varepsilon_M(e)=M_{\beta_1}-M_{\rm Ch}$. For
$M\in(M_{\rm Ch},M_{\rm Ch}+\varepsilon_M(e))$, continuity gives
$\beta\in(0,\beta_1)$ with $M_\beta=M$.
Theorem~\ref{thm:main-collapse} supplies a solution of this mass
with the core asymptotics and mass concentration in (i)--(ii).
Propositions~\ref{exlim:energy-conservation} and~\ref{exlim:clock-energy}
give, for this solution,
\[
 \mathscr E(t)=\mathscr E(t_0),\qquad
 \mathscr E+4M=\frac e2 I_2(\delta).
\]
This proves the remaining assertion in (ii).

For the same solution, fix $Y<\infty$.
The estimates proved in Subsection~\ref{exlim:actual-layer} give
\[
 \|\mathcal Z_\lambda-\mathcal Z_0\|_{C^3([0,Y])}
 +\sum_{\mathcal F\in\{\mathcal Q,\mathcal G,\Theta\}}
       \|\mathcal F_\lambda-\mathcal F_0\|_{C^1([0,Y])}
 \le C_{\beta,e,Y}\{
        \lambda(1+|\log\lambda|)^{p_{\log}}+\lambda^{J+1/2}\}.
\]
The limiting factors are positive, and the exact mass and radius
identities are
\[
 \int_{M-\lambda^4Y^5}^{M}\dd x=\lambda^4Y^5,
 \qquad
 R_b(t)-r(t,M-\lambda^4y^5)=\lambda^2y^2\mathcal G_\lambda(y).
\]
These prove (iii). The limiting layer is specified by
$P(\rho_0)=\varkappa_\delta D/\pi$ and
$Z_0=h(\rho_0)/(4\varkappa_\delta)+1/\varkappa_\delta$.
Its spatial vacuum trace also gives
\[
 \lambda^2[-\partial_rh_E(t,R_b(t))]
   =\frac{4\Theta_\lambda(0)^{2/3}}{\mathcal Q_\lambda(0)}
   \longrightarrow4\varkappa_\delta.
\]
The regularity follows from the same order-$22$ solution; no further
selection of a subsequence is needed for these limits.
\end{proof}

\appendix

\section{Weighted estimates and compactness}
\label{app:weighted-estimates}

We give the proofs of the weighted estimates stated in
Subsections~\ref{at:weighted-calculus}--\ref{at:localization}
and of the endpoint and compactness lemmas in
Subsection~\ref{exlim:compactness}. The multiplication and quotient
arguments use one-sided derivatives at the vacuum. The compactness
arguments use the equivariant center condition and the endpoint
conditions of the compatible completions.

\subsection{One-sided multiplication and composition}
\label{app:one-sided-calculus}

We use the spaces and abbreviated norms of
Subsection~\ref{at:weighted-calculus}. No parity condition is
imposed at the vacuum. The proofs distinguish the bounded
one-sided interval from the annulus with a growing upper endpoint.

\begin{proof}[Proof of Lemma~\ref{at:one-sided-multiplication}]
On the bounded Lipschitz cone
\[
 \Omega=\{(y,z)\in\mathbb R\times\mathbb R^{d-1}:
                   0<y<L,\ |z|<y\},
 \qquad U_F(y,z)=F(y),
\]
the full derivative-sum norm satisfies
\begin{equation}
 \|U_F\|_{H^r(\Omega)}^2
   =|B^{d-1}_1|\|F\|_{\mathcal B_d^r(0,L)}^2
   \qquad(r\in\mathbb N_0).
 \label{at:cone-norm-identity}
\end{equation}
This lift allows arbitrary one-sided derivatives of $F$ at zero;
it is not the radial lift $F(|X|)$.

The scalar case of \cite[Theorem~3.6]{HiptmairLiZou2012} gives
a single extension operator $E$ bounded at all integer orders $r\ge1$:
\[
 (EU)|_\Omega=U,\qquad
 \|EU\|_{H^r(\mathbb R^d)}\le C_{r,\Omega}\|U\|_{H^r(\Omega)}.
\]
Indeed, for zero-forms the graph norm in that theorem is equivalent to
$\|U\|_{H^{r-1}}+\|\nabla U\|_{H^{r-1}}\sim\|U\|_{H^r}$.
For $d\le8$, Cauchy--Schwarz gives
\[
 \begin{gathered}
 \|\widehat U\|_1\le C_d\|U\|_{H^6},\qquad
 \||\xi|\widehat U\|_1\le C_d\|U\|_{H^6},\\
 \int_{\mathbb R^d}(\langle\xi\rangle^{-12}
                  +\langle\xi\rangle^{-10})\dd\xi<\infty.
 \end{gathered}
\]
Using
$\langle\xi\rangle^m\le C_m(\langle\eta\rangle^m+
\langle\xi-\eta\rangle^m)$, Young's inequality and Fourier inversion yield
\[
 \begin{aligned}
 \|UV\|_{H^m}
 &\le C_m\bigl(\|U\|_{H^m}\|\widehat V\|_1
                +\|\widehat U\|_1\|V\|_{H^m}\bigr),\\
 \|U\|_{W^{1,\infty}}&\le C_d\|U\|_{H^6}.
 \end{aligned}
\]
Apply these bounds to $EU_F,EU_G$ and restrict to $\Omega$.
Equation~\eqref{at:cone-norm-identity} proves
\eqref{at:global-tame-product}--\eqref{at:global-low-embedding}.
Density extends the estimates to the completions, and uniform
convergence identifies their product with the pointwise product.

For the annulus, set $\ell=\min(c/4,1)$ and partition $[c,L]$
into intervals $I_i$ of lengths in $[\ell,2\ell]$. Since $L-c\ge4\ell$,
such a partition exists. With $a_i=\inf I_i$,
\[
 y^{d-1}\sim_c a_i^{d-1},\qquad
 \|F\|_{H^6(I_i)}
 \le c^{-(d-1)/2}\|F\|_{\mathcal C_d^6(c,L)}.
\]
The one-dimensional product and embedding constants are uniform on
these intervals. Hence
\[
 \begin{aligned}
 \|FG\|_{\mathcal C_d^m(c,L)}^2
 &\le C\sum_i a_i^{d-1}
 \bigl(\|F\|_{H^m(I_i)}^2\|G\|_{H^6(I_i)}^2
       +\|G\|_{H^m(I_i)}^2\|F\|_{H^6(I_i)}^2\bigr)\\
 &\le C_c\bigl(\|F\|_{\mathcal C_d^m}^2\|G\|_{\mathcal C_d^6}^2
             +\|G\|_{\mathcal C_d^m}^2\|F\|_{\mathcal C_d^6}^2\bigr),\\
 \|F\|_{W^{1,\infty}(c,L)}
 &\le C_c\sup_i\|F\|_{H^6(I_i)}
 \le C_c\|F\|_{\mathcal C_d^6(c,L)}.
 \end{aligned}
\]
The constants are independent of $L$.
\end{proof}

\begin{proof}[Proof of Lemma~\ref{at:low-rung-multiplier}]
For $r=0$ this is the $L^\infty$ embedding of
Lemma~\ref{at:one-sided-multiplication}.
On $\mathbb R^d$, let $T$ be multiplication by a fixed $V\in H^6$.
The preceding estimates give
\[
 \|T\|_{\mathcal L(L^2)}+\|T\|_{\mathcal L(H^6)}
       \le C_d\|V\|_{H^6}.
\]
Interpolation gives the same bound on $H^r$, $0\le r\le6$.
Indeed, $[L^2,H^6]_\theta=H^{6\theta}$ follows from the
weighted Fourier $L^2$ norms. Apply the three-lines inequality to
$\langle D\rangle^{6z}T\langle D\rangle^{-6z}$ on smooth
Fourier-truncated inputs, then remove the truncations by density.

For $r\ge1$, extend $U_G$ at order six and $U_F$ at order $r$,
multiply, restrict to the cone, and use
\eqref{at:cone-norm-identity}. This requires only the order-$r$
norm of $F$ and imposes no parity condition.

On the fixed-size annular intervals, the one-dimensional estimate gives
\[
 \begin{aligned}
 \|FG\|_{\mathcal C_d^r(c,L)}^2
 &\le C_c\sum_i a_i^{d-1}
             \|F\|_{H^r(I_i)}^2\|G\|_{H^6(I_i)}^2\\
 &\le C_c\|G\|_{\mathcal C_d^6(c,L)}^2
             \sum_i a_i^{d-1}\|F\|_{H^r(I_i)}^2\\
 &\le C_c\|G\|_{\mathcal C_d^6(c,L)}^2
             \|F\|_{\mathcal C_d^r(c,L)}^2.
 \end{aligned}
\]
The interval lengths and weight comparisons are uniform in $L$.
\end{proof}

\begin{proof}[Proof of Lemma~\ref{at:coefficient-composition}]
For $U\in H^m(\mathbb R^d)$ with $\|U\|_\infty\le M$, choose
$S_j=\chi(2^{-j}D)$, where $\chi$ is smooth, compactly supported,
and equal to one near zero. Set $V_j=(S_{j+1}-S_j)U$. Then
\[
 \begin{gathered}
 G(U)-G(0)=G(S_0U)-G(0)+\sum_{j\ge0}A_jV_j,\\
 A_j=\int_0^1DG(S_jU+tV_j)\dd t,\qquad
 \|\partial^\alpha A_j\|_\infty\le C_{\alpha,G,M}2^{j|\alpha|}.
 \end{gathered}
\]
The last bound follows from the uniformly bounded convolution kernels
of $S_j$ and the chain rule. Bernstein's inequality gives
\[
 \begin{gathered}
 \|A_jV_j\|_2\le C\|V_j\|_2,\qquad
 \|A_jV_j\|_{H^N}\le C_N2^{jN}\|V_j\|_2,\\
 \|G(S_0U)-G(0)\|_{H^N}\le C_N\|U\|_2.
 \end{gathered}
\]
For the last estimate, place one differentiated factor in $L^2$
and the others in $L^\infty$; at order zero use the mean-value formula.
With an inhomogeneous dyadic decomposition $\Delta_k$ and $N>m$,
\[
 2^{km}\|\Delta_k(A_jV_j)\|_2
 \le C_{N,G,M}
 \begin{cases}
  2^{-m(j-k)}2^{jm}\|V_j\|_2,&k\le j,\\
  2^{-(N-m)(k-j)}2^{jm}\|V_j\|_2,&k>j.
 \end{cases}
\]
Both sequence kernels are summable. Young's inequality and Plancherel,
using the finite overlap of the Fourier annuli, therefore give
\[
 \begin{aligned}
 \|G(U)-G(0)\|_{H^m}
 &\le C_{m,G,M}\left(\|U\|_2+
              \Bigl[\sum_j2^{2jm}\|V_j\|_2^2\Bigr]^{1/2}\right)\\
 &\le C_{m,G,M}\|U\|_{H^m}.
 \end{aligned}
\]
These bounds also justify the series by approximation from smooth $U$.

Apply this estimate to the cone extension $EU_F$. Since
\[
 \|EU_F\|_\infty\le C\|F\|_6,\qquad
 \|EU_F\|_{H^m}\le C\|F\|_m,
\]
restriction and \eqref{at:cone-norm-identity} prove
\eqref{at:composition}. If $G$ is defined only near the prescribed
compact range, use a compactly supported smooth extension of $G$;
neither $E$ nor the segments below need preserve that range.

For $H=F-\widetilde F$, write
\[
 G(F)-G(\widetilde F)
 =DG(0)H+
   \int_0^1\{DG(\widetilde F+tH)-DG(0)\}H\dd t.
\]
The product estimate and the composition bound give
\[
 \begin{aligned}
 \|G(F)-G(\widetilde F)\|_m
 &\le C\|H\|_m
 +C\int_0^1\bigl(
      \|\widetilde F+tH\|_6\|H\|_m
     +\|\widetilde F+tH\|_m\|H\|_6\bigr)\dd t\\
 &\le C_M\bigl\{\|H\|_m+
       (\|F\|_m+\|\widetilde F\|_m)\|H\|_6\bigr\}.
 \end{aligned}
\]
This is \eqref{at:composition-difference}.
On each annular interval from the preceding proof,
$\|F\|_{H^6(I_i)}\le C_cM$. Apply the same estimates there,
multiply by $a_i^{(d-1)/2}$, and square and sum. The low factors
are bounded globally; subtracting $G(0)$ eliminates the volume term.

Finally decompose
\[
 D^pG(F)[h_1,\ldots,h_p]
 =D^pG(0)[h_1,\ldots,h_p]
   +(D^pG(F)-D^pG(0))[h_1,\ldots,h_p].
\]
The first term is estimated by the finite-product inequality.
For the second, \eqref{at:composition} bounds its coefficient at
orders $6$ and $m$; the same product inequality gives
\eqref{at:coefficient-high-branch}. The difference bound passes
all these estimates to the completions.
\end{proof}

\subsection{Localized interpolation and products}
\label{app:localized-products}

The cutoffs and measures are those of
Subsection~\ref{at:localization}. We retain the prescribed cutoff
in each high norm and keep the annular volume explicit.

\begin{proof}[Proof of Lemma~\ref{at:power-interpolation}]
Set $p_k=2m/k$,
$N_k=\|\theta^{Qk/m}D^ku\|_{p_k}$, and $N_0=\|u\|_\infty$.
All integrands in $N_k^{p_k}$ have the same cutoff power $2Q$.
For $1\le k<m$, transfer one derivative from $D^ku$ in this
integral. The resulting terms are bounded by
\[
 C\int\theta^{2Q}|D^{k-1}u||D^ku|^{p_k-2}|D^{k+1}u|\dd\mu
 +C\int\theta^{2Q-1}|D^{k-1}u||D^ku|^{p_k-1}\dd\mu.
\]
On the weighted annulus the derivative of the measure is included
in the second term, since $7/y$ is bounded. A regularization of
$|D^ku|^{p_k-2}$ justifies this calculation when the exponent is
not an integer. The positive cutoff power eliminates the artificial
boundary terms. The principal term has exactly the required powers:
\[
 \frac{Q(k-1)}m+\frac{Q(k+1)}m+(p_k-2)\frac{Qk}m=2Q,
 \qquad
 \frac1{p_{k-1}}+\frac1{p_{k+1}}+\frac{p_k-2}{p_k}=1.
\]
Here $1/p_0=0$. In the second term the remaining cutoff exponent
is $Q/m-1\ge0$. The remaining factor in H\"older's inequality is the
constant function in $L^{2m}$, whose norm is $\mu_0^{1/(2m)}$.
Thus
\begin{equation}
 N_k^2\le C N_{k-1}N_{k+1}
                   +C\mu_0^{1/(2m)}N_{k-1}N_k.
 \label{at:joint-recurrence}
\end{equation}
Divide each norm by $\mu_0^{1/p_k}$, so the recurrence has no volume
factor. Unless $N_0=0$, set
\[
 a=\left(1+\frac{\mu_0^{-1/2}N_m}{N_0}\right)^{1/m},
 \qquad M_k=\frac{\mu_0^{-1/p_k}N_k}{N_0a^k}.
\]
Then $M_0=1$, $M_m\le1$, and
\[
 M_k^2\le C M_{k-1}M_{k+1}+C M_{k-1}M_k.
\]
Young's inequality gives
\[
 M_k\le\epsilon M_{k+1}+C_\epsilon M_{k-1}.
\]
Inductively substitute $M_{k-1}\le tM_k+C_tM_0$, choose
$t=(2C_\epsilon)^{-1}$, and absorb the $M_k$ term. Thus, at each
derivative order,
\[
 M_k\le\epsilon M_{k+1}+C_\epsilon M_0.
\]
Iteration to $m$, followed by $M_m\le1$, bounds every $M_k$.
Undoing the normalization proves \eqref{at:joint-interpolation}.

For \eqref{at:linear-interpolation}, use instead the $L^2$ norms
$I_k=\|\theta^{Q-m+k}D^ku\|_2$. Integration by parts gives
$I_k^2\le C I_{k-1}I_{k+1}+C I_{k-1}I_k$; the two principal
weights add to $2(Q-m+k)$, and a differentiated weight has exponent
$2(Q-m+k)-1$. The same finite absorption proves
$\sum_{k<m}I_k\le\epsilon I_m+C_\epsilon I_0$.
Finally $I_0\le\|u\|_2$. The argument applies first to smooth
fields and then to the indicated derivative norms by approximation.
\end{proof}

\begin{proof}[Proof of Lemma~\ref{at:same-cutoff-product}]
A derivative of $\chi uv$ has terms
$(D^\alpha\chi)(D^\beta u)(D^\gamma v)$ with
$|\alpha|+|\beta|+|\gamma|\le m$.
Because
\[
 Q-|\alpha|-Q(|\beta|+|\gamma|)/m\ge0,
\]
\eqref{at:joint-interpolation} applies to both factors with the
same cutoff. Write $U=\|u\|_\infty$, $V=\|v\|_\infty$ and
$H_u=\|\chi D^mu\|_2+\mu_0^{1/2}U$, and similarly for $v$.
H\"older's inequality bounds the term by
\[
 C\mu_0^{(m-|\beta|-|\gamma|)/(2m)}
 U^{1-|\beta|/m}V^{1-|\gamma|/m}
 H_u^{|\beta|/m}H_v^{|\gamma|/m}.
\]
If $U,V>0$, divide by $\mu_0^{1/2}UV$ and use
\[
 \left(\frac{H_u}{\mu_0^{1/2}U}\right)^{|\beta|/m}
 \left(\frac{H_v}{\mu_0^{1/2}V}\right)^{|\gamma|/m}
 \le \frac{H_u}{\mu_0^{1/2}U}
       +\frac{H_v}{\mu_0^{1/2}V}.
\]
Both bases are at least one, and $|\beta|+|\gamma|\le m$, so this
is the weighted arithmetic-geometric mean inequality. Multiplying
back and using \eqref{at:localized-jet} proves
\eqref{at:cutoff-product}. If $U=0$ or $V=0$, the product vanishes.
For $m=1$ the endpoint norms and the product rule give the same
conclusion directly. Applying this reasoning to each factor proves
the finite-product version.

Finally write the coefficient difference as
\[
 G(x,U_a+u)-G(x,U_a)
   =\int_0^1D_UG(x,U_a+tu)u\dd t.
\]
The chain rule leaves at least one increment factor in every term.
Profile derivatives and explicit $x$ derivatives are bounded coefficients.
Apply the product estimate to the increment factors and $v$ with the
same cutoff. This proves \eqref{at:cutoff-coefficient}, including the
term in which the highest derivative falls on the coefficient.
\end{proof}

\subsection{Quotients at the vacuum boundary}
\label{app:vacuum-quotients}

We first prove the quotient estimate in the compatible weighted
completion. The ordinary integral formulas
\eqref{exlim:ordinary-quotients}--\eqref{exlim:layer-division}
then give the endpoint values when the function has the stated
classical regularity.

\begin{proof}[Proof of Lemma~\ref{at:natural-quotient}]
First take a compactly localized smooth function with $f_y(0)=0$.
For $j\ge2$,
\[
 \partial_y^j(f_y/y)=\int_0^1t^jf^{(j+2)}(ty)\dd t.
\]
Dilation in $L^2(y^4\dd y)$ costs $t^{-5/2}$, so Minkowski gives
the bound with constant $(j-3/2)^{-1}$.

For orders zero and one, integration by parts gives the backward
Hardy estimates
\[
 \int y^2|h|^2\dd y\le\frac49\int y^4|h'|^2\dd y,
 \qquad
 \int |h|^2\dd y\le4\int y^2|h'|^2\dd y,
\]
when the artificial endpoint terms vanish. They yield
\[
 \begin{aligned}
 \|f_y/y\|_{B_5^0}
 &=\|f_y\|_{L^2(y^2\dd y)}
       \le\frac23\|f_{yy}\|_{B_5^0},\\
 \|(f_y/y)'\|_{B_5^0}
 &\le\|f_{yy}\|_{L^2(y^2\dd y)}+\|f_y\|_{L^2(\dd y)}\\
 &\le3\|f_{yy}\|_{L^2(y^2\dd y)}
       \le2\|f_{yyy}\|_{B_5^0}.
 \end{aligned}
\]

Boundary terms away from zero are controlled on the fixed overlap
of the localization. Since
$(\chi f)_y/y=\chi f_y/y+\chi'f/y$, the additional term is
supported away from zero and is estimated by
\eqref{at:localized-jet}. This proves \eqref{at:closed-quotient}.
Passing through the compatible approximating core defines the
quotient in the completion. Interior distributional tests identify
the result, so its value is independent of the approximation.
\end{proof}

\begin{proof}[Proof of Lemma~\ref{at:regular-factor-calculus}]
The fundamental theorem of calculus gives the two integral formulas.
For $0\le j\le k$, differentiation under the integral yields
\[
 \begin{aligned}
 (\mathcal Qg)^{(j)}(y)&=\int_0^1t^jg^{(j+2)}(ty)\dd t,
 &\| (\mathcal Qg)^{(j)}\|_\infty
     &\le\frac{\|g^{(j+2)}\|_\infty}{j+1},\\
 (\mathcal Gg)^{(j)}(y)&=\int_0^1(1-t)t^jg^{(j+2)}(ty)\dd t,
 &\| (\mathcal Gg)^{(j)}\|_\infty
     &\le\frac{\|g^{(j+2)}\|_\infty}{(j+1)(j+2)}.
 \end{aligned}
\]
Evaluation at $y=0$ gives the stated traces.

Suppose $v,w\ge c>0$ and $\|v\|_{C^k}+\|w\|_{C^k}\le M$.
Differentiating $vv^{-1}=1$ gives, for $j\ge1$,
\[
 (v^{-1})^{(j)}=-v^{-1}\sum_{i=1}^j\binom ji
                         v^{(i)}(v^{-1})^{(j-i)}.
\]
Induction, followed by Leibniz's rule, gives
\[
 \begin{aligned}
 \|v^{-1}\|_{C^k}+\|w^{-1}\|_{C^k}&\le C_k(c,M),\\
 \|v^{-1}-w^{-1}\|_{C^k}
 &=\|(w-v)v^{-1}w^{-1}\|_{C^k}
 \le C_k(c,M)\|v-w\|_{C^k}.
 \end{aligned}
\]
\end{proof}

\subsection{Endpoint regularity and the compatible completion}
\label{app:endpoint-completion}

We use the vacuum norm of Subsection~\ref{exlim:compactness}
and the compatible core of Definition~\ref{at:atlas}. The density argument below
determines exactly when the completion retains the first-derivative
trace at the vacuum.

\begin{proof}[Proof of Lemma~\ref{exlim:endpoint}]
On a fixed interior annulus the weight is nondegenerate. The
fundamental theorem of calculus and an average over the annulus
control the derivative traces needed below. Taylor's formula for
$f^{(q-3)}$, starting from an interior point $c_0$, has remainder
\[
 \frac12\int_y^{c_0}(t-y)^2f^{(q)}(t)\dd t.
 \]
The squared dual norm of this kernel in $L^2(t^4\dd t)$ is bounded
by $\frac14\int_y^{c_0}(t-y)^4t^{-4}\dd t\leq c_0/4$.
The kernels are continuous in $y$ in that dual norm. The remainder
is therefore bounded and continuous up to zero; integration gives
the lower derivatives and proves \eqref{exlim:endpoint-embedding}.
In particular the first-derivative trace is continuous when $q\geq4$.

It remains to identify the compatible completion. We first
approximate by smooth one-sided functions, then impose the
first-derivative condition at zero.

For $q\geq1$, subtract at $y=\epsilon$ the Taylor polynomial
$P_\epsilon$ of degree $q-1$.
Integration from $\epsilon$ satisfies
\begin{equation}
 \left\|\int_\epsilon^y v(t)\dd t\right\|_{L^2(y^4\dd y;(0,2\epsilon))}
       \leq C\epsilon\|v\|_{L^2(y^4\dd y;(0,2\epsilon))}.
 \label{exlim:anchored}
\end{equation}
Indeed, after scaling to $\epsilon=1$, Cauchy--Schwarz bounds
the integral by $C(1+y^{-3/2})\|v\|_{L^2(y^4\dd y)}$;
its square is integrable against $y^4\dd y$.
Iteration of \eqref{exlim:anchored} yields
\[
 \|(f-P_\epsilon)^{(r)}\|_{B_5^0(0,2\epsilon)}
 \leq C\epsilon^{q-r}\|f^{(q)}\|_{B_5^0(0,2\epsilon)}
 \quad(0\leq r\leq q).
 \]
Choose a smooth cutoff $\chi$ that vanishes on $[0,1]$ and equals
one on $[2,\infty)$, and set
\[
 f_\epsilon=P_\epsilon+\chi(y/\epsilon)(f-P_\epsilon).
\]
For $0\le r\le q$ and $0<\epsilon\le1$, Leibniz's rule gives
\[
 \begin{aligned}
 \|(f-f_\epsilon)^{(r)}\|_{B_5^0}
 &\le C\sum_{l=0}^r\epsilon^{-l}
       \|(f-P_\epsilon)^{(r-l)}\|_{B_5^0(0,2\epsilon)}\\
 &\le C\epsilon^{q-r}\|f^{(q)}\|_{B_5^0(0,2\epsilon)}
       \longrightarrow0.
 \end{aligned}
\]
Convolution away from zero produces smooth one-sided approximants
$g_n\to f$ in $B_5^q$. At order zero, cutoff and weighted $L^2$
approximation give the same density conclusion.

Let $\eta$ be a fixed smooth cutoff equal to one near zero.
If $q\ge4$ and $f'(0)=0$, trace continuity gives $g_n'(0)\to0$.
Thus
\[
 \widetilde g_n=g_n-g_n'(0)y\eta(y),\qquad
 \widetilde g_n'(0)=0,\qquad
 \|\widetilde g_n-g_n\|_{B_5^q}
       \le C_q|g_n'(0)|\longrightarrow0.
\]
For $q\le3$, the correction can instead be localized at scale
$0<\delta\le1$:
\[
 \|c_1y\eta(y/\delta)\|_{B_5^q}^2
       \le C_qc_1^2\sum_{r=0}^q\delta^{7-2r}
       \le C_qc_1^2\delta^{7-2q}\longrightarrow0.
\]
Choosing $\delta$ after $c_1$ imposes the derivative condition
without changing the limit. Thus the completion imposes no
continuous first-derivative trace at these orders; in particular,
$f(y)=y$ belongs to the order-two completion.
Finally, \eqref{exlim:endpoint-embedding} gives
$B_5^{k+5}\hookrightarrow C^{k+2}$, and
Lemma~\ref{at:regular-factor-calculus} gives the two quotient bounds
on the subspace with the first-derivative trace just identified. No evenness was imposed.
\end{proof}

\subsection{Compactness at fixed positive scale}
\label{app:fixed-scale-compactness}

Fix $0<a<A_0$ and use the interval $I_a$ and pressure-form space
$V_a$ of Subsection~\ref{exlim:inputs}. The norm equivalences in
the following proof may depend on $a$. They are used to take each
limit on a fixed interval; the uniform estimates are obtained
separately in the main argument.

\begin{proof}[Proof of Lemma~\ref{exlim:spatial-compactness}]
We reduce to fixed coordinates, prove compactness in the weighted
Sobolev spaces, and then obtain the compact inclusion into $V_a^*$
by duality.

\emph{Fixed coordinates.} Choose a boundary neighborhood $M-x<c a^4$, so the fifth-root chart
is available throughout $I_a$. There the coordinate change is
the bounded dilation $y_5=\lambda^{-4/5}(M-x)^{1/5}$.
The complement is covered by a fixed center ball and finitely
many intervals away from the endpoints.

To compare the norms, use that the four cutoffs sum to one. At
each point some cutoff is at least $1/4$; on a sufficiently small
neighborhood in space and scale, its reciprocal has bounded
derivatives through the required finite order. A finite cover then
compares \eqref{exlim:graph} in both directions with the fixed
localized weighted norm in the fixed physical coordinates. These comparison constants may depend on
$a$.

\emph{Spatial compactness.} For $v\in B_5^1$, an interior trace and
Cauchy--Schwarz give
\[
 |v(y)|^2\leq C(1+y^{-3})\|v\|_{B_5^1}^2,\qquad
 \int_0^\epsilon y^4|v(y)|^2\dd y
                       \leq C\epsilon^2\|v\|_{B_5^1}^2.
 \]
Apply this to $v=f^{(r)}$, $0\le r\le q$. On
$[\epsilon,c]$, let $P_{\epsilon,h}$ replace each component of the
derivative tuple by its average on intervals of length at most $h$,
and set it to zero on $(0,\epsilon)$. Then
\[
 \sum_{r=0}^q\|f^{(r)}-P_{\epsilon,h}f^{(r)}\|_{B_5^0}^2
 \le C\bigl(\epsilon^2+C_\epsilon h^2\bigr)
                                      \|f\|_{B_5^{q+1}}^2.
\]
Its range is finite dimensional. Choose $\epsilon$, then $h$, to
make this error small. Interior distributional identities identify
any limit of the derivative tuples with one function.

At the center, localize the full vector inside a larger ball.
Fourier frequencies above $N$ contribute at most $N^{-2}$
times its squared $H^{q+1}$ bound to the $H^q$ norm. The
low-frequency derivatives are uniformly bounded and uniformly
continuous on the smaller ball, so finite spatial grids give
compactness there. On the fixed intervals away from the endpoints,
the same approximation by local averages applies without an
endpoint tail. Summing the
fixed partition proves the first inclusion in
\eqref{exlim:compact-inclusions}. The equivalence of $V_a$ with
$F_a^1$ gives $V_a\Subset H$.

\emph{The dual inclusion.} Choose a finite-rank orthogonal projection
$P_\epsilon$ in $H$ such that
\[
 \sup_{\|v\|_{V_a}\le1}\|(I-P_\epsilon)v\|_H\le\epsilon.
\]
For $f\in H$,
\[
 \|(I-P_\epsilon)f\|_{V_a^*}
 =\sup_{\|v\|_{V_a}\le1}|(f,(I-P_\epsilon)v)_H|
 \le\epsilon\|f\|_H.
\]
The finite-rank maps $P_\epsilon:H\to V_a^*$ approximate the
inclusion in operator norm, proving $H\Subset V_a^*$.
\end{proof}

\begin{proof}[Proof of Lemma~\ref{exlim:time-compactness}]
Write $(B_1,B_0)=(F_a^{q+1},F_a^q)$ in the first case and
$(B_1,B_0)=(H,V_a^*)$ in the second. Bochner integration gives
\[
 \|f^n(t)-f^n(s)\|_{B_0}\le C|t-s|^\gamma,
 \qquad
 \gamma=1\ \hbox{or}\ \tfrac12,
\]
respectively; the second estimate uses Cauchy--Schwarz in time.
To extend the $B_1$ bound to every time, choose times $t_k\to t$
at which the bound holds. Extract a weakly convergent subsequence
in the bounded $B_1$ ball and identify its limit with $f^n(t)$ in $B_0$.
Thus $\{f^n(t)\}_n$ is relatively compact in $B_0$ by
Lemma~\ref{exlim:spatial-compactness}.

For a grid $t_0<\cdots<t_N$ of mesh at most $h$,
\[
 \|f^n-f^m\|_{C(I_a;B_0)}
 \le\max_{0\le k\le N}\|f^n(t_k)-f^m(t_k)\|_{B_0}
                         +2Ch^\gamma.
\]
Extract at the finitely many grid points and then diagonalize over
$h\downarrow0$. The displayed bound makes the resulting subsequence
Cauchy in $C(I_a;B_0)$. A finite collection of derivatives is treated by
the same extraction.
\end{proof}

\section{Derivative and remainder estimates for the matched profile}
\label{app:profile-estimates}

We prove the auxiliary estimates for the finite profile constructed in
Section~\ref{mat:chapter}. The coefficient functions and connection
amplitudes are those of Subsection~\ref{mat:connection}. At fixed
$\beta>0$, $e\ge0$, and finite orders $K,S=K+1$, the constants
may depend on that family and the prescribed cutoffs. We then reduce
the upper bound for the scale as needed.

\subsection{Finite asymptotic expansions and integral remainders}
\label{app:profile-finite}

Use the finite expansion classes of
Subsubsection~\ref{mat:connection-part-6}. The scale variable and the
logarithm are independent during coefficient extraction. In the layer,
$Q=D^{1/4}$ and $E_Q=Q\partial_Q$; in the core,
$q=1-z$ and $E_q=q\partial_q$. We retain the constants of integration
and estimate the derivatives of the remainders.

\subsubsection{Finite products and compositions}\label{app:profile-finite-calculus}

We prove Lemma~\ref{mat:finite-polyhomogeneous}.

\begin{proof}
For products, let $f$ have largest exponent $m$ and remainder
exponent $r$, and let $g$ have largest exponent $n$ and remainder
exponent $s$. Expand
\[
 fg=f_{\rm fin}g_{\rm fin}
       +f_{\rm rem}g_{\rm fin}+f_{\rm fin}g_{\rm rem}
       +f_{\rm rem}g_{\rm rem}.
\]
For $0\le j\le\min(d_f,d_g)$, Leibniz's rule gives
\[
 \left|E_Q^j\{fg-f_{\rm fin}g_{\rm fin}\}\right|
 \le C\bigl(Q^{r+n}+Q^{s+m}+Q^{r+s}\bigr)
                 (1+\log Q)^{p_{\log}}.
\]
Retain from $f_{\rm fin}g_{\rm fin}$ the powers strictly above
$\max\{r+n,s+m,r+s\}$ and include the remaining terms in the
remainder. Near $q=0$, the same calculation uses the minimum of
the three error exponents and retains the powers below it.

For an $L$-independent normalized positive factor $c(1+v)$ with $c>0$ and
$v=O_{E^{\le d}}(Q^{-\eta}(1+\log Q)^p)$, $\eta>0$, use Taylor's
formula with integral remainder for $(1+v)^\alpha$ or $\log(1+v)$.
Choose $N$ so that $(N+1)\eta$ is beyond the requested endpoint
order. Every derivative of the remainder through $d$ still has at
least $N+1$ factors of the form $E_Q^jv$, because
\[
 (1+v)^\alpha-\sum_{j=0}^{N}\binom\alpha jv^j
 =\frac{v^{N+1}}{N!}\int_0^1(1-t)^N
       (\alpha)_{N+1}(1+tv)^{\alpha-N-1}\dd t.
\]
Reduce the asymptotic interval so that $|v|\le1/2$. The denominators in
this formula are then positive. On the complementary finite interval, the function is retained
without asymptotic expansion.
The reciprocal can equivalently use the exact finite geometric identity.
The same proof at $q=0$ applies with $v=O(q)$.

For $L$-polynomial perturbations, truncate the composition at the
prescribed scale degree. Its coefficients remain polynomial in $L$,
although the untruncated reciprocal need not be polynomial in $L$. A factor
$U_0+\sum_{j=1}^Nt^jU_j$, with $U_0$ independent of $L$, has inverse
coefficients
\[
 V_0=U_0^{-1},\qquad
 V_n=-U_0^{-1}\sum_{j=1}^nU_jV_{n-j}.
\]
For a smooth outer function and $1\le n\le N$,
\[
 [t^n]G\!\left(U_0+\sum_{j=1}^Nt^jU_j\right)
 =\sum_{r=1}^n\frac{G^{(r)}(U_0)}{r!}
      \sum_{\substack{j_1+\cdots+j_r=n\\j_1,\ldots,j_r\ge1}}
      U_{j_1}\cdots U_{j_r}.
\]
Every index on the right is at most $n$, so discarded scale
degrees do not affect a retained coefficient.

Apply this rule to the exact pressure at the positive leading density
$\rho_0$, the core Jacobian factors, and the positive leading layer
coefficient $A_0$. For the endpoint expansion, factor the pressure law into
its leading power, powers of positive factors, and the explicit
logarithm in \eqref{mat:connection-38}. The preceding Taylor estimate
controls the remainder. At low density,
\eqref{mat:endpoint-6} is smooth after its stated power is removed.

The required density derivative of the outer function may have
finite order depending on $N$. It acts on the explicit pressure law and
does not require another derivative of a profile coefficient.

Coefficient extraction and differentiation in the independent $L$
are linear maps on a finite polynomial space and cost no endpoint
derivative. Finally $\partial_D=(4Q^4)^{-1}E_Q$ and
$\partial_q=q^{-1}E_q$ prove the rules for spatial differentiation. Inspecting
\eqref{mat:connection-6} and \eqref{mat:connection-13}, a Jacobian
uses one derivative, differentiating the pressure adds one more, and
the inertia uses two; there is no third coefficient derivative.
\end{proof}

\subsubsection{Asymptotic integration of the layer equation}\label{app:profile-volterra}

We prove Lemma~\ref{mat:finite-volterra}.

\begin{proof}
Near $D=0$, boundedness of $S$ and the positive leading factor of $A$ give
\[
 \begin{gathered}
 A(s)^{-1}\int_0^s S(t)\dd t=O(s^{-3/5}),\qquad
 Z(D)-\gamma=O(D^{2/5}),\\
 A(D)Z'(D)=\int_0^D S(t)\dd t\longrightarrow0.
 \end{gathered}
\]
On compact subintervals of $(0,\infty)$, positivity of $A$
and the stated local regularity make both integrals classical.

We first retain the two constants of integration, then estimate
the remainder and its derivatives.

\emph{The first primitive.} Let $I_{\rm fin}$ be the sum of the
primitives of $4Q^3S_{\rm fin}$ based at $Q_b$, computed by
\eqref{mat:connection-46}. Choose $\chi=0$ below $Q_b$ and
$\chi=1$ above $2Q_b$, and set
$\widehat S_{\rm rem}(Q)=S(Q^4)-\chi(Q)S_{\rm fin}(Q)$.
Since $\sigma_S<-4$, the limit
$\mu=\lim_{Q\to\infty}(I-I_{\rm fin})$ exists and
\begin{equation}
 \begin{split}
 I(Q)&:=\int_0^{Q^4}S(t)\dd t
       =I_{\rm fin}(Q)+\mu-J(Q),\\
 J(Q)&=\int_Q^\infty4t^3 \widehat S_{\rm rem}(t)\dd t
       =O_{E_Q^{\le d+1}}(Q^{\sigma_S+4}(1+\log Q)^{p_{\log}}),
 \qquad Q\ge2Q_b.
 \end{split}
 \label{mat:volterra-actual-moment}
\end{equation}
Indeed, integration bounds $J$, and
$E_QJ=-4Q^4\widehat S_{\rm rem}$ bounds its next $d+1$
derivatives. The constant $\mu$ includes both the bounded vacuum
interval and the cutoff transition. Monomials that are not integrable
at zero have only been integrated from $Q_b>0$.

\emph{The second primitive.} Multiply \eqref{mat:volterra-actual-moment}
by $4Q^3/A(Q^4)$. The retained polynomial terms have explicit
primitives, including logarithms at exponent $-1$. The constant
$\mu$ contributes $-\mu\mathcal H_A(Q)$, where
\[
 \begin{aligned}
 \mathcal H_A(Q)&=\int_Q^\infty\frac{4t^3}{A(t^4)}\dd t
                  =\frac4{3a_0}Q^{-3}(1+o(1)),\\
 A(Q^4)&=a_0Q^7(1+o(1))\qquad(Q\to\infty).
 \end{aligned}
\]
Its expansion has remainder exponent $\sigma_A+4$. The other
outer error integrals have exponents
\[
 \sigma_S+1,\qquad k+\sigma_A+8,\qquad
 \sigma_S+\sigma_A+8.
\]
All four exponents are at most $R<-3$. Write each error integral
as a constant minus its tail. The constants, including the original
$\gamma$, give the regular branch; the tails are
$O(Q^R(1+\log Q)^{p_{\log}})$. This retains every required term
of the $Q^{-3}$ branch. Changing $\chi$ changes only the decomposition,
not either exact primitive.

\emph{Differentiated remainders.} The integral formulas already
control the remainder and its first $E_Q$ derivative. In the $Q$
coordinate the equation is
\[
 E_Q^2Z+\left(3+E_Q\log\frac{A(Q^4)}{Q^7}\right)E_QZ
                   =\frac{16Q^8}{A(Q^4)}S(Q^4).
\]
Subtract the finite expansion. If $r$ denotes the remainder and
$B=3+E_Q\log(A(Q^4)/Q^7)$, the resulting equation has the form
\[
 E_Q^2r+B E_Qr=F,\qquad
 E_Q^jF=O(Q^R(1+\log Q)^{p_{\log}}),\quad j\le d.
\]
Lemma~\ref{mat:finite-polyhomogeneous} gives these bounds on $F$.
For $0\le j\le d$,
\[
 E_Q^{j+2}r
 =E_Q^jF-\sum_{i=0}^j\binom ji(E_Q^iB)E_Q^{j-i+1}r.
\]
Induction gives the asserted $d+2$ derivatives with the same
remainder exponent. Near $D=0$, \eqref{mat:connection-33} is the
zero-flux solution; for the constructed sources, differentiation
of \eqref{mat:connection-32} gives the stated ordinary regularity
in $D^{1/5}$.
\end{proof}

\subsubsection{Endpoint orders and differentiated remainders}\label{app:profile-remainder-orders}

We prove Lemma~\ref{mat:finite-remainder-schedules}.

\begin{proof}
We induct on scale degree, and at each core degree on decreasing
logarithmic degree.

\emph{Leading expansions.} The exact pressure is
\begin{equation}
 P(\rho)=\tfrac12\{\rho^{1/3}(2\rho^{2/3}-3)
       \sqrt{1+\rho^{2/3}}+3\operatorname{arsinh}(\rho^{1/3})\}.
 \label{mat:connection-38}
\end{equation}
Set $t=Q^{-1}$ and $\rho_0^{1/3}=(\varkappa/\pi)^{1/4}Qv$.
The implicit equation becomes
\[
 \mathscr P(t,v):=\frac\pi\varkappa t^4
 P\bigl((\varkappa/\pi)^{3/4}t^{-3}v^3\bigr)-1=0.
\]
Its leading term is $v^4-1$, followed by a multiple of $t^2v^2$
and terms of the forms
\begin{equation}
 t^4\log t,\quad t^4\log v,\quad t^{4+2j}v^{-2j}.
 \label{mat:connection-39}
\end{equation}
At degree $D_*=2K+12$, coefficient comparison determines
$v_{\rm fin}(t,\log t)$: each new coefficient is multiplied by four.
The pressure expansion through $j=K+6$ has remainder starting
at normalized degree $2K+18>D_*+1$. Taylor's integral formula in
\eqref{mat:connection-38}, together with
$\operatorname{arsinh}(q)=\log q+\log(1+\sqrt{1+q^{-2}})$, gives
\[
 E_t^r\mathscr P(t,v_{\rm fin})
 =O\bigl(t^{D_*+1}(1+|\log t|)^{p_{\log}}\bigr),
 \qquad 0\le r\le25.
\]
For small $t$, the mean-value coefficient
$\int_0^1\mathscr P_v(t,v_{\rm fin}+\theta(v-v_{\rm fin}))\dd\theta$
is at least two. Implicit differentiation and induction on $r$ yield
\begin{equation}
 \begin{gathered}
 |E_t^r(v-v_{\rm fin})|
       \le C t^{D_*+1}(1+|\log t|)^{p_{\log}},\\
 0\le r\le25.
 \end{gathered}
 \label{mat:connection-40}
\end{equation}
Only ordinary derivatives of the explicit pressure through $3K+70$
are used here; the profile coefficients are not differentiated further.
Substitution gives
\begin{equation}
 \begin{aligned}
 Z_0&=\text{a finite Laurent-log polynomial}\\
    &\quad+O_{E_Q^{\le25}}(Q^{-2K-12}\log^{p_{\log}} Q),\\
 A_0^{-1}&=\text{a finite Laurent-log polynomial}\\
    &\quad+O_{E_Q^{\le25}}(Q^{-2K-20}\log^{p_{\log}} Q).
 \end{aligned}
 \label{mat:connection-41}
\end{equation}
Thus $A_0^{-1}$ has largest power $-7$. The pressure factors in
\eqref{mat:connection-35}, including their density derivatives
through $r\le S$ and their first 25 iterates under $E_Q$, have
relative remainder $O(Q^{-2K-13}\log^{p_{\log}}Q)$.

\emph{Layer remainders.} A remainder $O(Q^{R_j})$ in $Z_j$
contributes $O(Q^{R_j-1})$ to $Z'_j/Z'_0$. At degree $k>j$,
the remaining factors have total power at most $k-j$, so
\begin{equation}
 Q^{R_j-1+k-j}=Q^{k-2K-11}=Q^{R_k-1}.
 \label{mat:connection-42}
\end{equation}
The radius factor gives the same exponent,
$R_j+(k-j-1)=R_k-1$. Multiple remainders and the pressure-factor
remainders decay faster. Hence
\begin{equation}
 \sigma_S=R_k-1=k-2K-11.
 \label{mat:connection-43}
\end{equation}
With $\sigma_A=-2K-20$, Lemma~\ref{mat:finite-volterra} gives
the four integrated powers
\begin{equation}
 \sigma_S-7+8=R_k,\quad
 k+\sigma_A+8=R_k-2,\quad
 \sigma_S+\sigma_A+8<R_k,\quad
 \sigma_A+4\le R_k.
 \label{mat:connection-44}
\end{equation}
All satisfy the required bound $O(Q^{R_k}\log^{p_{\log}}Q)$.

Since $\sigma_S<-4$, the source remainder is integrable in $D$:
\begin{equation}
 \int_0^D S_{\rm rem}
   =\int_0^\infty S_{\rm rem}
                          -\int_D^\infty S_{\rm rem}.
 \label{mat:connection-45}
\end{equation}
Cut off each asymptotic monomial near $D=1$ and include the finite
transition in the remainder. The actual moment in
\eqref{mat:connection-45}, multiplied by the retained terms of
$A_0^{-1}$, produces the complete homogeneous tail. For retained
monomials use \eqref{mat:connection-45} only when integrable at
infinity; otherwise integrate from a positive base point, using
\begin{equation}
 \int x^{\gamma-1}(\log x)^p\dd x
 =x^\gamma\sum_{j=0}^p
       \frac{(-1)^jp!}{(p-j)!\gamma^{j+1}}
                       (\log x)^{p-j}\quad(\gamma\ne0);
 \label{mat:connection-46}
\end{equation}
and $(\log x)^{p+1}/(p+1)$ when $\gamma=0$. These formulas retain
both resonances and every required term of the $Q^{-3}$ branch.

\emph{Core remainders.} A radius remainder of order $q^{H_j+1}$
enters the degree-$n$ source with power at least $H_j-(n-j)$, and
\begin{equation}
 H_j-(n-j)-H_n=n-j\ge1,\qquad j<n.
 \label{mat:connection-47}
\end{equation}
The smooth source factors can carry $q^{-n}$. Taylor expansion of
$q^nF_n$, including the required 22 iterates of $E_q$ and two
source derivatives, is therefore justified by
\begin{equation}
 \max\{H_n+30,\ H_n+n+25\}\le 2S+50<2S+60.
 \label{mat:connection-48}
\end{equation}
The leading-profile ODE \eqref{mat:connection-1} and the finite
Leibniz recurrence for positive reciprocals supply the stated
$2S+60$ derivatives on the fixed center and endpoint charts.

After subtracting the finite particular polynomial, solve
\begin{equation}
 f(q)=\int_0^q\frac1{s^4a(s)}
                   \int_0^s t^3\{G(t)+c(t)f(t)\}\dd t\dd s,
 \quad G=O_{E_q^{\le22}}(q^{H_n}\log^{p_{\log}}q).
 \label{mat:connection-49}
\end{equation}
In the norm
\[
 \|f\|_*=
 \sup_{0<q\le q_0}
 \frac{|f(q)|}{q^{H_n+1}(1+|\log q|)^{p_{\log}}},
\]
write the right side as $f_G+\mathcal Tf$. Since $a\ge c_0>0$
and $c$ is bounded, direct integration gives
\[
 \begin{aligned}
 |f_G(q)|&\le Cq^{H_n+1}(1+|\log q|)^{p_{\log}},\\
 |\mathcal Tf(q)|&\le Cq^{H_n+2}(1+|\log q|)^{p_{\log}}\|f\|_*.
 \end{aligned}
\]
Choose $Cq_0<1$. The fixed point has both homogeneous coefficients
zero; its difference from the center-regular solution is a linear
combination of the complete branches \eqref{mat:connection-18},
which are retained.

\emph{Differentiated remainders.} The core equation is
\begin{equation}
 E_q^2f+(3+E_q\log a)E_qf
                         =q\,a^{-1}(G+cf).
 \label{mat:connection-50}
\end{equation}
The two integrals bound $f,E_qf$ with power $q^{H_n+1}$.
Applying $E_q^j$, $0\le j\le22$, gives
\[
 \begin{split}
 E_q^{j+2}f={}&E_q^j\{qa^{-1}(G+cf)\}\\
 &-\sum_{i=0}^j\binom ji
 E_q^i(3+E_q\log a)E_q^{j-i+1}f.
 \end{split}
\]
Each right side uses only previously bounded derivatives of $f$;
the factors $qa^{-1}$ preserve the required power. This proves
the bounds through order 24. In the layer,
$A_0(Q^4)=Q^7a_\infty(Q)$ gives
\begin{equation}
 E_Q^2Z+(3+E_Q\log a_\infty)E_QZ
                         =16Q a_\infty^{-1}S(Q^4).
 \label{mat:connection-51}
\end{equation}
After subtracting the finite expansion and actual moments,
\eqref{mat:connection-44}--\eqref{mat:connection-46} bound the remainder
and its first $E_Q$ derivative. The identical Leibniz induction in
\eqref{mat:connection-51} gives the next 23 derivatives. Thus
\[
 24\ \text{coefficient derivatives}
 \ \Longrightarrow\ 22\ \text{source derivatives}
 \ \Longrightarrow\ 24\ \text{remainder derivatives}.
\]
The first implication uses exactly the two coefficient derivatives
in \eqref{mat:connection-6} and \eqref{mat:connection-13}.

For ordinary derivatives at the center and vacuum, take
\begin{equation}
 t_n=60+2(S-n),\qquad t_{n-1}-2=t_n.
 \label{mat:connection-52}
\end{equation}
The source costs at most two ordinary derivatives, and
\eqref{mat:connection-15}, \eqref{mat:connection-33} preserve the
requested order. Lower scale degrees therefore supply every source
derivative; higher logarithmic coefficients at the same scale degree
enter without spatial differentiation. This completes the finite
induction from the leading bounds
\eqref{mat:connection-40}, \eqref{mat:connection-48}.
\end{proof}

\subsection{Derivative bounds for the profile and force}
\label{app:profile-derivatives}

The local norms, enlarged chart intervals, and fixed parameter choices
are those of Subsection~\ref{mat:native}. Ordinary derivatives are used
at the center and on bounded charts; the endpoint growth estimates
use iterates of $E_q$ and $E_Q$. The radius is the finite matched sum,
with core terms through degree $K$ and the layer coefficient
$Z_K$ retained.

\subsubsection{Derivatives of the coefficient functions}\label{app:profile-coefficient-growth}

We prove Lemma~\ref{mat:native-lem:wave2-finite-native-growth}.

\begin{proof}
For $31\le t_*\le60$, the derivative orders satisfy
\[
 t_{n-1}-2=t_n,\qquad
 t_2+2=2S+t_*-2\le2S+58<2S+60.
\]
We induct on scale degree, solving each core equation in decreasing
logarithmic degree.

\emph{Ordinary derivatives.} Set
$W(\eta)=w(\sqrt\eta)$,
$\Phi(\eta)=\phi_{n\ell}(\sqrt\eta)$, and
$\mathcal F(\eta)=F_{n\ell}(\sqrt\eta)$, where $F_{n\ell}$
is the full source in \eqref{mat:core-log-recursion}.
At vacuum write $A_0(y^5)=y^8\widehat a(y^2)$ and
$\widehat S(y)=S(y^5)$. The exact integral formulas are
\begin{align}
 \Phi'(\eta)&=\frac{3}{8W(\eta)^4}
 \int_0^1t^4W(t^2\eta)^3
       \{c_n\Phi(t^2\eta)-\mathcal F(t^2\eta)\}\dd t,
 \qquad c_n=(2n^2+n-3)|\delta|,
 \label{mat:native-eq:wave2-center-native-equation}\\
 \widehat Z(y)&=25\int_0^1\frac{s}{\widehat a(y^2s^2)}
                         \int_0^1t^4\widehat S(tys)\dd t\dd s.
 \label{mat:native-eq:wave2-vacuum-native-equation}
\end{align}
\equationalias{mat:units-8}{mat:native-eq:wave2-center-native-equation}
\equationalias{mat:units-7}{mat:native-eq:wave2-vacuum-native-equation}
Differentiate under the integrals. Since
$\partial_y^j\widehat S(tys)=(ts)^j\widehat S^{(j)}(tys)$,
Leibniz's rule and the positive coefficient bounds give
\[
 \|\widehat Z\|_{C^j}\le C_j\|\widehat S\|_{C^j},\qquad
 \|\Phi\|_{C^{j+1}}
 \le C_j\bigl(|\Phi(0)|+\|\mathcal F\|_{C^j}\bigr),
 \qquad 0\le j\le t_n.
\]
For the second inequality, first apply the integral inequality to
the zeroth-order equation, then induct on $j$.
No parity of $\widehat Z$ is used. The regular interior ODE and
positive coordinate Jacobian give the compact-interior bounds.

\emph{Endpoint differential equations.} In the core,
\[
 (q^4af')'=q^3(G+cf),\qquad a\ge c_0>0,
\]
where $\varkappa=\varkappa_\delta$ and
\[
 a=\frac{w(1-q)^4(1-q)^4}{\varkappa^4q^4},\qquad
 b_0=\frac{w(1-q)^3(1-q)^4}{\varkappa^3q^3},\qquad
 G=-\frac{3b_0}{4\varkappa}F_{n\ell},\qquad
 c=\frac{3c_nb_0}{4\varkappa}.
\]
Put $B=3+(E_qa)/a$ and $C=q/a$. For $r\ge0$,
\begin{align}
 E_q^{r+2}f={}&-\sum_{j=0}^r\binom rj(E_q^jB)E_q^{r-j+1}f
  +\sum_{j=0}^r\binom rj(E_q^jC)E_q^{r-j}G\\
 &+\sum_{j=0}^r\binom rj E_q^j(Cc)E_q^{r-j}f.
 \label{mat:native-eq:wave2-Fuchs-derivative-recurrence}
\end{align}
\equationalias{mat:units-9}{mat:native-eq:wave2-Fuchs-derivative-recurrence}
Here $E_q^jC,E_q^j(Cc)=O(q)$. In the layer, set
$B_V=E_D\log A_0-1$, $C_V=D^2/A_0$, with $E_D=D\partial_D$.
Then
\begin{equation}
 E_D^{r+2}Z=-\sum_{j=0}^r\binom rj(E_D^jB_V)E_D^{r-j+1}Z
       +\sum_{j=0}^r\binom rj(E_D^jC_V)E_D^{r-j}S.
 \label{mat:native-eq:wave2-Volterra-derivative-recurrence}
\end{equation}
The bounds
\[
 E_D^jB_V=O(1),\qquad E_D^jC_V=O(Q)
\]
follow from $P(\rho_0)=\varkappa_\delta D/\pi$ and
$A_0=16\pi^2\rho_0^2P'(\rho_0)$. Indeed, put
$v=\log D$, $\theta(v)=\log\rho_0(e^v)$, and
$F_0(\theta)=\log P(e^\theta)-\log(\varkappa_\delta/\pi)$.
Then $F_0(\theta(v))=v$, and
\[
 \theta'=\frac1{F_0'(\theta)},\qquad
 F_0'(\theta)\theta^{(r)}
 =-\sum_{\substack{\pi\in\Pi_r\\|\pi|\ge2}}
 F_0^{(|\pi|)}(\theta)\prod_{B\in\pi}\theta^{(|B|)},
 \quad r\ge2,
\]
where $\Pi_r$ denotes the partitions of $r$ labelled indices.
The factor $F_0'=\rho_0P'(\rho_0)/P(\rho_0)$ is bounded
above and below for $D\ge1$. The explicit pressure remainder and
\[
 \rho^jP^{(j)}(\rho)
 =\prod_{a=0}^{j-1}(\rho\partial_\rho-a)P(\rho)
\]
therefore bound every derivative used here, $r\le t_2+2$;
also $\rho_0^{1/3}\sim Q$.

\emph{Source estimates.} A core radius coefficient of degree $k$
has power $q^{1-k}$, and its Jacobian variation has power $q^{-k}$.
The pressure term of degree $2j$, after one spatial derivative,
division by $w^3$, and a deformation product of degree $n-2j$,
has power
\[
 (4-2j)-1-3-(n-2j)=-n.
\]
The $\lambda^4\log(w/\lambda)$ term obeys the same count;
gravity and inertia have power at least $1-n$.
In the layer the reciprocal-density convolution gives
\[
 \begin{gathered}
 Z_{k-1}=O(Q^k\log^{p_{\log}}Q),\qquad
 Z'_{k-1}/Z'_0=O(Q^{k-1}\log^{p_{\log}}Q),\\
 [\lambda^m](\rho/\rho_0-1)=O(Q^m\log^{p_{\log}}Q).
 \end{gathered}
\]
At index $n-1$, the pressure has power at most $Q^{n+3}$;
$\partial_D=(4Q^4)^{-1}E_Q$ lowers it to $Q^{n-1}$.
The time and gravity terms have the same bound. The finite
compositions of $(1-u)^{-2}$, $(1+v)^{-1}$, and the pressure
Taylor coefficients exhaust the source terms. Thus
\[
 \begin{aligned}
 |E_q^jF_{n\ell}|&\le Cq^{-n}(1+|\log q|)^{p_{\log}},\\
 |E_D^jS_{n-1,\ell}|&\le CQ^{n-1}(1+\log Q)^{p_{\log}},
 \qquad 0\le j\le t_n.
 \end{aligned}
\]
Each source differentiates a lower-degree coefficient at most twice,
and $t_n+2\le t_k$ for $k<n$. Previously solved logarithmic
coefficients at the same degree enter without spatial differentiation.

\emph{Integration.} The core homogeneous powers are $0,-3$.
A source of power $-n$ gives a particular solution of power $1-n$,
with a logarithm at $n=4$. The connection conditions remove the
complete singular branch at $n=2,3$, while $q^{-3}=O(q^{1-n})$
for $n\ge4$. Hence
\[
 |f(q)|\le Cq^{1-n}(1+|\log q|)^{p_{\log,n}},\qquad
 E_qf=\frac{C_*+\int_{q_*}^q t^3(G+cf)\dd t}{q^3a(q)}.
\]
At $n=2,3$, $C_*=q_*=0$. At $n\ge4$ take $q_*>0$;
the numerator has power $q^{4-n}$, with a logarithm at $n=4$.
Thus $E_qf$ has the same bound as $f$, and
\eqref{mat:native-eq:wave2-Fuchs-derivative-recurrence} gives all
iterates through $t_n$.
For the layer,
\[
 E_DZ=\frac D{A_0(D)}\int_0^DS(t)\dd t
 =O\bigl(Q^n(1+\log Q)^{p_{\log,n}}\bigr).
\]
The second primitive, including its fixed additive constant, bounds
$Z$ by the same expression. Equation
\eqref{mat:native-eq:wave2-Volterra-derivative-recurrence} gives
the remaining derivatives. Products add logarithmic degrees and
the two primitives add at most two; their finite maximum defines
$p_{\log,n}$ and completes the induction.

These are bounds for the complete coefficients. For the particular
core remainder with $H=H_n$ in \eqref{mat:connection-37},
the same equations give only the implication
\[
 E_q^jG=O(q^H\log^{p_{\log}}q)\ (j\le25)
 \quad\Longrightarrow\quad
 E_q^jf=O(q^{H+1}\log^{p_{\log}}q)\ (j\le27),
\]
after removing both homogeneous branches. The layer requires the
corresponding differentiated source bound and subtraction of all
moments. No extra derivative of an asymptotic remainder follows
from the complete-coefficient bounds alone.
\end{proof}

\subsubsection{Derivatives of the radius and force}\label{app:profile-radius-force}

We prove Proposition~\ref{mat:native-prop:wave2-native-profile-source}.

\begin{proof}
Write $A=\lambda u$ and let $T$ be the fixed-mass scale derivative.
We first bound $u$, then $V=\lambda^2\mathcal M[A]$.

\emph{Radius.} On the compact core and interior charts the finite
coefficient sums have bounded derivatives through order 27.
At the center $u=z_ca_c(\lambda,z_c^2)$, so the vector lift is
regular. In the bounded vacuum chart,
\[
 \begin{gathered}
 u=1-\lambda Z_0
   -\sum_{n=2}^S\lambda^n\sum_\ell(\log\lambda)^\ell Z_{n-1,\ell},\\
 T=\lambda\partial_\lambda-\tfrac45y_5\partial_{y_5}.
 \end{gathered}
\]
For $L=\log\lambda$, the exact layer identity is
\[
 T^r\{\lambda^nL^\ell f(Q)\}
 =\lambda^n\sum_{a=0}^{\min(r,\ell)}
 \binom ra\frac{\ell!}{(\ell-a)!}L^{\ell-a}(n-E_Q)^{r-a}f(Q).
\]
On the bounded chart replace $E_Q$ by
$(4/5)y_5\partial_{y_5}$; in the core $q$ is fixed.
Thus $r+m\le27$ uses at most 27 of the available 31 ordinary
coefficient derivatives.

On the layer part of the long chart, $Q=y_8^2$,
$1\le Q\le C\lambda^{-1/2}$, and $T=\lambda\partial_\lambda-E_Q$.
The coefficient bounds give
\[
 \begin{gathered}
 |E_Q^jZ_0|\le C_jQ,\qquad
 \lambda^n|E_Q^jZ_{n-1,\ell}|
 \le C(\lambda Q)^n(1+\log Q)^{p_{\log,n}},\\
 \sup_{1\le Q\le C\lambda^{-1/2}}
 (\lambda Q)^{n-1}(1+\log Q)^{p_{\log,n}}=o(1),\quad n\ge2,\\
 \partial_{y_8}^m=y_8^{-m}\prod_{a=0}^{m-1}(2E_Q-a).
 \end{gathered}
\]
Since $y_8\ge c$, these inequalities bound the ordinary mixed
derivatives as well. On the core part,
\[
 \begin{gathered}
 q=c_*\lambda y_8^2U(\lambda y_8^2)\ge c\sqrt\lambda,
 \qquad U>0,\qquad y_8\partial_{y_8}q=O(q),\\
 \lambda^n|E_q^j\phi_{n\ell}|
 \le Cq(\lambda/q)^n(1+|\log q|)^{p_{\log,n}}
 \le Cq\lambda^{n/2}(1+|\log\lambda|)^{p_{\log,n}}.
 \end{gathered}
\]
The coordinate factors have bounded iterates under $E_q$.
On the matching shell, $q\sim\sqrt\lambda$ and Leibniz's rule gives
\[
 u_{\rm c}-1, u_{\rm v}-1=O(q),\qquad
 T^rE_q^j\{\chi(q/\sqrt\lambda)(u_{\rm c}-u_{\rm v})\}=O(q)
\]
through the required orders. This uses the complete coefficients,
not higher derivatives of the matching remainder.

Define $a_{00}=1$, set coefficients outside the triangle to zero, and put
\[
 a_{i+1,r}=\partial_sa_{ir}+(1-3i/2)ba_{ir}+ba_{i,r-1}.
\]
Since $\partial_\tau=\lambda^{-3/2}\partial_s$,
\[
 \partial_\tau^i(\lambda u)
 =\lambda^{1-3i/2}\sum_{r=0}^ia_{ir}T^ru.
\]
The scalar coefficients and their required derivatives are bounded
for fixed $\beta>0$. The preceding chart estimates prove
\eqref{mat:native-eq:wave2-native-matched-rows}.

\emph{Force.} On the bounded vacuum chart, the exact low-density
law $P(\rho)=\rho^{5/3}p_0(\rho^{2/3})$ gives
\[
 Z_D=y_5^{-3}b_Z,\qquad
 \rho=(4\pi u^2Z_D)^{-1}=y_5^3b_\rho,\qquad
 P(\rho)=y_5^5b_P,
\]
with positive factors and bounded ordinary derivatives. Hence
\begin{equation}
 \lambda^2\mathcal M[A]
 =-4\pi u^2\left(b_P+\tfrac15y_5\partial_{y_5}b_P\right)
           +\frac{M-\lambda^4y_5^5}{u^2}.
 \label{mat:native-eq:wave2-direct-bounded-force}
\end{equation}
The factor $b_Z$ uses one coefficient derivative and $\partial_{y_5}b_P$
uses one more: $r+m\le25$ requires at most 27 derivatives.

On the long layer write $\rho=\rho_0\theta$. The relative Jacobian
formula gives bounded mixed derivatives of $\theta,\theta^{-1}$
under $T,E_Q$. Since $P(\rho)/\rho^{4/3}$ and its iterates under
$\rho\partial_\rho$ are bounded for $\rho\ge c>0$, and
$\rho_0^{1/3}/Q$ has bounded $E_Q$ derivatives,
\[
 P(\rho_0\theta)=Q^4p,\qquad
 |T^rE_Q^jp|\le C_{rj},\qquad
 \partial_DP(\rho)=p+\tfrac14E_Qp.
\]
Thus $V=-4\pi u^2\partial_DP(\rho)+x/u^2$ has the required
mixed bounds. Conversion to $y_8$ again costs only bounded powers
of $y_8^{-1}$.

On the core and shell put
\[
 J=\frac{u^2u_z}{z^2},\qquad
 \rho=\lambda^{-3}w^3J^{-1},\qquad
 H=\lambda h(\rho)=4\sqrt{\lambda^2+w^2J^{-2/3}}-4\lambda.
\]
Then
\begin{equation}
 \lambda^2\mathcal M[A]=\frac{H_z}{u_z}+\frac{m(z)}{u^2}.
 \label{mat:native-eq:wave2-direct-core-force}
\end{equation}
Here $u_z\ge c>0$. For the cutoff term
$h=\chi(u_{\rm c}-u_{\rm v})$, the radius bounds imply
\[
 h_z=-q^{-1}E_qh=O(1),\qquad
 qh_{zz}=q^{-1}(E_q^2-E_q)h=O(1).
\]
They also hold after the required material and $E_q$ derivatives.
Consequently $J^{\pm1}$ and $qJ_z$ are bounded in these norms, and
\[
 H_z=4\frac{wJ^{-2/3}}{\sqrt{\lambda^2+w^2J^{-2/3}}}
       \left(w_z-\frac{w}{3qJ}\,qJ_z\right)
\]
is bounded: $w\sim q$, $w_z=O(1)$, and
$q\ge c\sqrt\lambda$ on both the core edge and shell.
No extension of the layer beyond the matching interval is used.

At the center, $u=z_ca_c(\lambda,z_c^2)$ and
$H=H_c(\lambda,z_c^2)$, so
\[
 \frac{H_{z_c}}{u_{z_c}}
 =z_c\frac{2\partial_2H_c}
             {a_c+2z_c^2\partial_2a_c},\qquad
 \frac{x}{u^2}=\frac{z_c}{a_c^2}.
\]
Their full vector lifts have bounded Cartesian derivatives.
On the compact interior, \eqref{mat:native-eq:wave2-direct-core-force}
is a smooth composition of the coefficient functions and their
first two derivatives. This proves the pointwise force bound.

\emph{Time derivatives and the weighted norm.} Define
\[
 d_{00}=1,\qquad
 d_{r+1,j}=\partial_sd_{rj}-2bd_{rj}+bd_{r,j-1}.
\]
Then
\[
 \partial_s^r(\lambda^{-2}V)
 =\lambda^{-2}\sum_jd_{rj}T^jV,\qquad
 \lambda^4\int_c^{C_e\lambda^{-1/2}}y_8^7\dd y_8\le C.
\]
The coefficients $d_{rj}$ are bounded and spatially constant.
Integrating the local derivative bounds against the chart measures,
with the prescribed cutoffs, proves
\eqref{mat:native-eq:wave2-direct-profile-force-source}.
\end{proof}

\subsection{Physical regularity and normalized coefficients}
\label{app:profile-physical}

We prove the physical-coordinate regularity and the coefficient bounds
\eqref{mat:units-1} and \eqref{mat:units-4}. Use the matched trajectory
and the normalized factors $r,\gamma,T_\gamma$ defined in
Subsection~\ref{mat:units}. The initial-radius coordinate is fixed at
$\lambda_*>0$; its unscaled derivative bounds may depend on
$\lambda_*$. The normalized operator estimates are uniform for
$0<\lambda\le A_*$, where the upper bound for the scale is chosen
for the fixed family.

\subsubsection{Derivative bounds for the local data}\label{mat:units-part-2}

We apply the derivative induction of
Lemma~\ref{mat:native-lem:wave2-finite-native-growth} using the
same coefficient equations and recurrences at its larger finite order.
In this subsection,
$t_n$ denotes the larger derivative order

\begin{equation}
 t_n=60+2(S-n),\qquad 2\le n\le S.                       \label{mat:units-5}
\end{equation}

The derivative counts needed for this induction are\nopagebreak

\begin{equation}
 t_{n-1}-2=t_n,\qquad
 t_2+2=2S+58\le 2S+60.                             \label{mat:units-6}
\end{equation}

Lemma~\ref{mat:native-lem:wave2-finite-native-growth} with
$t_*=60$ supplies these bounds: the fixed-integral formulas
\eqref{mat:native-eq:wave2-center-native-equation}--\eqref{mat:native-eq:wave2-vacuum-native-equation}
and the recurrences for $E_q$ and $E_D$ apply at the orders \eqref{mat:units-5}
by \eqref{mat:units-6}. The coefficients and connection
amplitudes are unchanged. These are bounds for the complete
coefficient functions; the derivative range of the overlap remainder
estimate remains unchanged.

\subsubsection{Initial data in physical coordinates}\label{mat:units-part-3}

At a fixed scale, the bounded vacuum profile has the exact form

\begin{equation}
 A(y)=A(0)-\lambda_*^2 y^2G(y),\quad G\in C^{60},\qquad
 \gamma=2G+yG_y>0.                                      \label{mat:units-10}
\end{equation}

These positive lower bounds follow from the lower-order estimates for
the constructed profile. In the normalized mass coordinate $y=y_5$, its density is

\[
 \rho(y)=y^3b_\rho(y),\qquad
 b_\rho(y)=\frac5{4\pi(A(y)/\lambda_*)^2\gamma(y)}
                                                  \in C^{59}.
\]

There is no fractional-power differentiation at zero. The exact enthalpy is\nopagebreak

\begin{equation}
 h(y)=y^2\frac{4b_\rho(y)^{2/3}}
                         {\sqrt{1+y^2b_\rho(y)^{2/3}}+1}
                                                   \in C^{59}. \label{mat:units-11}
\end{equation}

The normalized and physical boundary coordinates are related by\nopagebreak

\begin{equation}
 z=y\sqrt{G(y)},\qquad y_{\rm ph}=\lambda_*z,
 \qquad z'(y)=\frac{2G+yG_y}{2\sqrt G}>0.                 \label{mat:units-12}
\end{equation}

The positive derivative in \eqref{mat:units-12} gives a $C^{60}$
inverse. Composing \eqref{mat:units-11} with this inverse yields a
$C^{59}$ function of $z$. The rescaling $y_{\rm ph}=\lambda_*z$
preserves this order, with derivative bounds depending on
$\lambda_*^{-1}$. This proves the $C^{58}$ enthalpy assertion in
\eqref{mat:units-1}.

The exact material derivative is

\begin{equation}
 \partial_\tau A=\lambda^{-3/2}b(\lambda)
    \left(\lambda\partial_\lambda-\frac45y\partial_y\right)A.
                                                        \label{mat:units-13}
\end{equation}

For $V=\partial_\tau A$, the only spatial differentiation in
\eqref{mat:units-13} is $y\partial_y$. In particular,
\[
 y\partial_y(y^2G)=y^2(2G+yG_y),\qquad
 V(y)=V(0)+y^2G_V(y),\quad G_V\in C^{59}.
\]
Differentiation of the finite scale factors costs no spatial derivative.
Thus $V\in C^{59}$ and $V_y(0)=0$. Composition with
\eqref{mat:units-12} gives the same regularity in $y_{\rm ph}$,
which implies the $C^{57}$ assertion in \eqref{mat:units-1}.

At the center write $A=\lambda_*z_c a(z_c^2)$. The density identity is\nopagebreak

\begin{equation}
 \rho=\frac3{4\pi\lambda_*^3
                    a(z_c^2)^2\{a(z_c^2)+2z_c^2a'(z_c^2)\}}.
                                                        \label{mat:units-14}
\end{equation}

The denominator is positive and even. Hence the enthalpy is an
even scalar and the velocity is $z_c$ times an even coefficient.
Formula \eqref{mat:units-14} uses one ordinary coefficient
derivative, within the range \eqref{mat:units-5}. The center
radius map has a regular odd inverse, so composition preserves the
Cartesian scalar and vector regularity in \eqref{mat:units-1}.

On the fixed compact interior, both the radius and mass Jacobians
are positive. The same inverse derivative recurrence preserves the
local regularity before localization.

The correction directions are fixed smooth functions in these physical
coordinates. The enthalpy corrections vanish to at least order three
at vacuum, the velocity corrections have zero first derivative there,
and the mass correction is supported strictly in the interior. Their
finite sums preserve the stated regularity, scalar or vector type, and
strict positivity bounds on a sufficiently small parameter ball.
Scalar coefficients used to correct the mass do not lower spatial
regularity. These bounds therefore hold before the compatibility
parameters are determined.

\subsubsection{Normalized coefficients on the boundary charts}\label{mat:units-part-4}

The operator coefficients in \eqref{mat:units-4} are functions of
the normalized radius and Jacobian factors. We estimate these factors
on the two edge regions and the joining shell before using the exact
coefficient formulas.

On $y_5\in[0,L_5]$, \eqref{mat:units-10} immediately gives $\gamma\in C^{59}$ and
$T_\gamma\in C^{58}$, with positive lower bounds for $r,\gamma$.
All these bounds are uniform in $\lambda$ in normalized coordinates.

\Needspace{6\baselineskip}
On the long layer $Q=y_8^2\ge c>0$, the exact identity is

\begin{equation}
 \gamma=2Z_Q
 =2(Z_0)_Q+2\sum_{n=2}^S\lambda^{n-1}
       \sum_\ell(\log\lambda)^\ell(Z_{n-1,\ell})_Q.
                                                        \label{mat:units-15}
\end{equation}

The equation $P(\rho_0)=\varkappa_\delta Q^4/\pi$ gives
$c\le (Z_0)_Q\le C$ and bounded iterates of $(Z_0)_Q$ and its
reciprocal under $E_Q$. For each correction in \eqref{mat:units-15},
the coefficient bounds give, for $0\le j\le22$,
\[
 \left|E_Q^j\!\left\{
  \lambda^{n-1}(\log\lambda)^\ell
                 (Z_{n-1,\ell})_Q\right\}\right|
 \le C(\lambda Q)^{n-1}
       (1+|\log\lambda|+|\log Q|)^{p_{\log}}.
\]
Since $n\ge2$ and $c\le Q\le C\lambda^{-1/2}$, the right side
tends uniformly to zero as $\lambda\downarrow0$. Thus
$\gamma\ge c>0$ and $|E_Q^j\gamma|\le C_j$ for $j\le22$.

On the part of the edge chart where the core radius is used, put

\begin{equation}
 q=c_*\lambda y^2U(\lambda y^2),\quad
 \omega=\frac{yq_y}{2q},\qquad
 \gamma=-\frac{2q}{\lambda y^2}\omega u_q,\quad u=A/\lambda.
                                                        \label{mat:units-16}
\end{equation}

Both coordinate factors in \eqref{mat:units-16} are positive,
with bounded iterates under $y\partial_y$ on the fixed physical
endpoint interval. The complete coefficient bounds give
\[
 u-1=O(q),\qquad E_q^ju=O(q)\quad(1\le j\le23)
\]
through every order used below. Hence
$u_q=q^{-1}E_qu$ and its required $E_q$ derivatives are bounded.

On the matching shell, put $h=u_{\rm c}-u_{\rm v}$ and
$\chi=\chi_{\rm m}(q/\sqrt\lambda)$. Then
\[
 \begin{aligned}
 u_q&=(1-\chi)(u_{\rm v})_q+\chi(u_{\rm c})_q+\chi_qh,\\
 E_q^j\partial_q(\chi h)&=q^{-1}(E_q-1)^jE_q(\chi h),
 \qquad 0\le j\le22.
 \end{aligned}
\]
The complete coefficient bounds give $E_q^jh=O(q)$ at the
required orders, and every $E_q$ derivative of $\chi$ is bounded.
These identities therefore give all upper bounds for $u_q$.
For positivity, use the lower-order matching estimate
\[
 |h|+|E_qh|
 \le C\lambda^{(K+1)/2}(1+|\log\lambda|)^{p_{\log}}.
\]
Since $q\sim\sqrt\lambda$,
\[
 \begin{aligned}
 (u_{\rm c})_z&=1+o(1),\\
 (u_{\rm v})_z-(u_{\rm c})_z,\ \chi_qh
   &=O\!\left(\lambda^{K/2}(1+|\log\lambda|)^{p_{\log}}\right).
 \end{aligned}
\]
Thus $-u_q\ge c>0$ for the matched radius when $A_*$ is small.
Off the shell the same bound follows from the small core Jacobian
correction. Equation~\eqref{mat:units-16} now gives the positive
lower bound for $\gamma$.

For ordinary derivatives use

\[
 \partial_y^jf=y^{-j}\prod_{a=0}^{j-1}(y\partial_y-a)f.
\]

Since $y\ge y_-$, the bounds for $(y\partial_y)^j$ give $C^{21}$ bounds on $\gamma$
and $T_\gamma$. The last of these uses at most twenty-three conormal
derivatives of $u$. Also $r_y=-\lambda y\gamma$, and
$\lambda y^2\le C$ on the entire covering chart. All ordinary derivatives
of $r$, its reciprocal, $x=M-\lambda^4y^8$, and $\lambda y^2$
that occur below are therefore uniformly bounded.

\emph{Operator coefficients.} We use the exact formulas in
Subsection~\ref{mat:units}, in particular \eqref{mat:units-17}.

In terms of the positive factors $c_n=n/(4\pi r^2\gamma)$, the principal coefficients
have the nonsingular expressions\nopagebreak

\begin{equation}
 G_5=\frac{4c_5^{2/3}}
              {3\gamma^2\sqrt{1+c_5^{2/3}y^2}},\qquad
 G_8=\frac{4c_8^{1/3}}
              {3\gamma^2\sqrt{1+c_8^{-2/3}y^{-4}}}.        \label{mat:units-18}
\end{equation}

All derivatives of \eqref{mat:units-18} needed in \eqref{mat:units-17} are bounded on the specified
positive factor ranges and local $y$ intervals. The first expression is used on a
bounded interval; the second on $y\ge y_->0$. Their denominators also
give the positive lower bound for $a_n$. The factors
$k(\rho)$, $\lambda y^2$, $r^{\pm1}$, $\gamma^{\pm1}$ and their
ordinary derivatives are bounded. Formula \eqref{mat:units-17} is affine in $T_\gamma$.
The finite ordinary chain and Leibniz rules through order twenty-one
therefore prove \eqref{mat:units-4}. In fact $\|v_n\|_{W^{21,\infty}}\le C\lambda$,
which is stronger than its uniform boundedness.

On the fixed center set $\Psi=\lambda^{-1}A e_r$,
$J=D\Psi$, $B=J^{-T}$, and $c_0=3/(4\pi\det J)$. The normalized
principal matrix of the full vector completion is exactly
\[
 \lambda P'(\lambda^{-3}c_0)B^TB
 =\frac{4c_0^{1/3}}{3\sqrt{1+\lambda^2c_0^{-2/3}}}B^TB.
\]
The positive Jacobian and the bounds through order sixty in the
chart coordinates give uniform ellipticity and Cartesian coefficient
derivatives through order twenty-one. The first-order tensor depends
smoothly on $J$ and linearly on $D J$; its derivatives through
twenty-one therefore use at most twenty-three local radius
derivatives. The center representation \eqref{mat:units-14} gives
regular gravity coefficients as well.

These center estimates concern the full Cartesian vector operator.
They do not estimate its separately singular radial coefficients
or $(M-x)\partial_x$ at the center. On the fixed compact
interior, the corresponding one-dimensional identities have the
same derivative count, and the density factor stays bounded away
from zero.

These estimates control the normalized coefficients on the entire
cover, including the part where the core radius is used. For
compatibility at the vacuum boundary, we also need estimates for
the residual in the regular endpoint coordinate.

\subsection{The endpoint residual and physical time derivatives}
\label{app:profile-endpoint}

Use the finite bounded-vacuum profile and positive quotients $G,Q$
from Subsection~\ref{mat:endpoint}. On a fixed enlarged interval
$0\le y=y_5\le Y_0$, put $u=1-\lambda Z$, $A=\lambda u$,
and $T=\lambda\partial_\lambda-(4/5)y\partial_y$.
The reduced residual $F_v$ and normalized residual
$\mathcal R_u$ are defined in \eqref{mat:endpoint-2}.
We first prove \eqref{mat:endpoint-3}, then take physical-time
derivatives and pull them back to the fixed initial-radius coordinate.

\emph{The residual in vacuum coordinates.} The exact identities $x=M-\lambda^4y^5$, $A=\lambda u$ give
\begin{equation}
 Z_D=\frac{Q}{5y^3},\qquad \rho=y^3B,\qquad
 B=\frac5{4\pi u^2Q}.
 \label{mat:endpoint-5}
\end{equation}
Integrating the exact derivative of $P$ in \eqref{mat:eos} yields
\begin{equation}
 P(y^3B)=y^5p(y,B),\qquad
 p(y,B)=4B^{5/3}\int_0^1
          \frac{t^4}{\sqrt{1+y^2B^{2/3}t^2}}\dd t.
 \label{mat:endpoint-6}
\end{equation}
This is smooth at $y=0$ for $B>0$. In particular
\begin{equation}
 \begin{split}
 F_v={}&\beta^2\{(1+T)^2-\tfrac32(1+T)\}u
       +e\lambda T(1+T)u\\
 &-4\pi u^2\left(p(y,B)+\frac y5\partial_yp(y,B)\right)
       +\frac{M-\lambda^4y^5}{u^2}.
 \end{split}
 \label{mat:endpoint-7}
\end{equation}
Differentiation in the pressure term includes $B_y$.
Since $Q=2G+yG_y$, the reduced residual uses at most two
derivatives of $G$.
No division by $y$ remains in this formula.

\emph{Taylor cancellation.} Introduce independent $t,L$ and replace each
$\lambda^n(\log\lambda)^\ell$ in \eqref{mat:endpoint-1} by
$t^nL^\ell$. Set $u=1-tZ(t,L,y)$ and
\begin{equation}
 {\mathbb T}=t\partial_t+\partial_L-\frac45y\partial_y.
 \label{mat:endpoint-8}
\end{equation}
Let $F(t,L,y)$ be \eqref{mat:endpoint-7} with these replacements.
The term $\partial_L$ differentiates the explicit powers of $\log\lambda$.
For $L=\log\lambda$ and $0\le t\le\lambda$,
\begin{equation}
 \|tZ\|_{C^{59}}\le C\lambda(1+|L|)^{p_{\log}},\qquad
 \|Q-Q_0\|_{C^{59}}\le C\sum_{n=1}^K\lambda^n(1+|L|)^{p_{\log}}.
 \label{mat:endpoint-9}
\end{equation}
Thus $u\ge1/2$ and $Q\ge c_0/2$ on the entire Taylor segment.
All finite $t,L$ derivatives of the polynomial $Z$, its positive
reciprocal compositions, and the integral \eqref{mat:endpoint-6}
have polynomial logarithmic majorants there.
For $q\le22$ and $a\le28$, the largest coefficient spatial order is
\begin{equation}
 q+a+2\le52<60.
 \label{mat:endpoint-10}
\end{equation}
Consequently the $C^{28}_y$ norms of the $t$ derivatives through
$K+1$ of $({\mathbb T}-3)^qF$ are bounded by $C(1+|L|)^{p_{\log}}$.
The many auxiliary $t$ derivatives act on a finite scale polynomial
and smooth functions of positive factors, not on additional spatial
derivatives of a coefficient.

At a fixed positive scale, the force uses two derivatives of
$G\in C^{60}$ and each material derivative adds at most one. Hence,
before pullback to the initial-radius coordinate,
\[
 \partial_\tau^m(A_{\tau\tau}+\mathcal M[A])\in C^{58-m},
 \qquad 0\le m\le22.
\]
The positive $C^{60}$ coordinate map in \eqref{mat:endpoint-18}
preserves this regularity. These higher derivative bounds may depend
on the initial scale; the uniform decay estimate is
\eqref{mat:endpoint-3}, in its stated smaller derivative range.

The leading equation cancels degree zero. At degree $n\ge1$, the new
coefficient appears only as $\partial_D(A_0\partial_DZ_n)$;
all other terms are exactly its already constructed source.
The coefficient equations therefore give the polynomial identities
\begin{equation}
 \partial_t^nF(0,L,\cdot)=0,\qquad 0\le n\le K.
 \label{mat:endpoint-11}
\end{equation}
Since ${\mathbb T}$ preserves scale degree, the derivatives
$\partial_t^nH_q(0,L,\cdot)$, $0\le n\le K$, also vanish for
$H_q=({\mathbb T}-3)^qF$. Taylor's formula in $C^{28}$ gives
\begin{equation}
 H_q(\lambda,\log\lambda,y)
 =\frac{\lambda^{K+1}}{K!}
   \int_0^1(1-\theta)^K
       \partial_t^{K+1}H_q(\theta\lambda,\log\lambda,y)\dd\theta .
 \label{mat:endpoint-12}
\end{equation}
Since $T^q(\lambda^{-3}F_v)=\lambda^{-3}H_q(\lambda,\log\lambda,y)$,
this proves \eqref{mat:endpoint-3}.
The matching shell is at $y\sim\lambda^{-2/5}$, outside the fixed
interval, so no cutoff derivative is omitted.

\emph{Time derivatives in physical coordinates.} We now estimate the
source coefficients used to impose compatibility on the initial data. Set
$\mathfrak R=A_{\tau\tau}+\mathcal{M}[A]=\lambda\mathcal{R}_u$
and $a_j=1-3j/2$. Define
\begin{equation}
 p_{00}=1,\quad p_{jq}=0\quad(q<0\text{ or }q>j),\qquad
 p_{j+1,q}=\partial_sp_{jq}+a_jbp_{jq}+bp_{j,q-1}.
 \label{mat:endpoint-13}
\end{equation}
Using $\partial_s\lambda=b\lambda$ and
$\partial_s(T^q\mathcal{R}_u)=bT^{q+1}\mathcal{R}_u$, induction gives
\begin{equation}
 \partial_\tau^j\mathfrak R
 =\lambda^{1-3j/2}\sum_{q=0}^jp_{jq}T^q\mathcal{R}_u .
 \label{mat:endpoint-14}
\end{equation}
The coefficients are bounded through $j=22$ for the selected
positive $\beta$. For the source coefficients, use the weight
\begin{equation}
 \begin{gathered}
 \omega_\ell=2\ell-\sigma_{\ell-1},\qquad
 \sigma_a=\begin{cases}4,&a\text{ odd},\\5/2,&a\text{ even},\end{cases}\\
 \chi_j=\omega_{j+2}+1-\frac32j,\qquad
 \chi_{2n}=n+1,\quad\chi_{2n+1}=n+3.
 \end{gathered}
 \label{mat:endpoint-15}
\end{equation}
For $[f]_a=f^{(a)}(0)/a!$, put
\begin{equation}
 c_{(j+2,q,a)}=[\lambda^{\chi_j}p_{jq}T^q\mathcal{R}_u]_a,\qquad
 [\lambda^{\omega_{j+2}}\partial_\tau^j\mathfrak R]_a
       =\sum_{q=0}^jc_{(j+2,q,a)}.
 \label{mat:endpoint-16}
\end{equation}
Equations \eqref{mat:endpoint-3}, \eqref{mat:endpoint-13} imply
\begin{equation}
 |c_{(j+2,q,a)}|
 \le C\lambda^{K-2+\chi_j}(1+|\log\lambda|)^{p_{\log}},\qquad
 0\le q\le j\le22,\quad 0\le a\le28.
 \label{mat:endpoint-4}
\end{equation}
These are the Taylor coefficients of the prescribed profile residual.
They are distinct from the Laurent coefficients obtained by
differentiating the Euler equation for data that do not yet satisfy
the compatibility conditions. If $K\ge2J^\sharp+86$, $J^\sharp\ge3$, then
\[
 K-1-(J^\sharp+55/2)\ge J^\sharp+115/2>0.
\]

At a fixed initial slice define the normalized radius coordinate $z$ by
$a=A(\lambda,y)$, $A(\lambda,0)-a=\lambda^2z^2$. Exactly,
\begin{equation}
 z=y\sqrt{G(\lambda,y)},\qquad
 z_y=\frac{2G+yG_y}{2\sqrt G}=\frac Q{2\sqrt G}>0 .
 \label{mat:endpoint-18}
\end{equation}
The positive factors $G,Q$ have uniform $C^{59}$ bounds on a
smaller fixed enlarged interval. Successive differentiation of
$z(y(z))=z$ isolates the highest inverse derivative with divisor
$z_y$; its other factors have lower inverse order. Hence
\begin{equation}
 \|f(y(z))\|_{C^{28}_z}\le C\|f\|_{C^{28}_y}.
 \label{mat:endpoint-19}
\end{equation}
Apply this only after the physical time derivatives have been taken.
The same bound for the full function underlying \eqref{mat:endpoint-4} gives
\begin{equation}
 \max_{\substack{0\le m\le22\\0\le r\le23}}
 \left|[\lambda^{\omega_{m+2}}\partial_\tau^m\mathfrak R]^z_r\right|
 \le C\lambda^{J^\sharp+55/2}(1+|\log\lambda|)^{p_{\log}} .
 \label{mat:endpoint-20}
\end{equation}
The map \eqref{mat:endpoint-18} is fixed at the uncorrected initial slice;
it is not differentiated in a correction parameter. This estimate
does not assert uniformity in the unscaled physical coordinate
$\lambda z$.

The endpoint derivative estimates hold on every fixed bounded interval
in the local vacuum coordinate.
The remaining force estimates require uniform control toward the moving
overlap, which we now derive from the residuals of the uncut profiles.

\subsection{Differentiated core and layer residuals}
\label{app:profile-core-layer}

We prove \eqref{mat:face-3}--\eqref{mat:face-5} for the uncut
radii and residuals of Subsection~\ref{mat:face}, with its fixed
parameters and coefficient bounds. The auxiliary expansion parameter
is $\lambda/q$ in the core and $\lambda Q$ in the layer.
On their respective intervals it is bounded by $C\sqrt\lambda$.
The logarithm is held independent during Taylor expansion.

\subsubsection{Finite composition and Taylor remainders}\label{mat:face-part-2}

For the auxiliary families, consider finitely many terms of the form
\begin{equation}
 f(\epsilon,L,\zeta)=f_0(\zeta)
          +\sum_{n=1}^{d}\sum_\ell
                    \epsilon^n L^\ell f_{n\ell}(\zeta),          \label{mat:face-8}
\end{equation}
where $E=\zeta\partial_\zeta$, and all required $E^r f_{n\ell}$
are bounded by $C(1+|\log\zeta|)^{p_{\log}}$.
Every finite mixed derivative in $\epsilon,L,E$ of \eqref{mat:face-8} is therefore
bounded by a polynomial in $1+|L|+|\log\zeta|$, uniformly for
$0\le\epsilon\le1$. An ordinary $\epsilon$-derivative differentiates
only a finite polynomial, not the coefficient function.

This assertion is stable under products by Leibniz's rule. If a factor
is bounded below by $c>0$, the same assertion holds for any fixed
real power, reciprocal, square root, or logarithm used below. Indeed,
applying $a$ successive derivatives to $G(f_1,\ldots,f_d)$
gives a finite sum of terms
\begin{equation}
 (\partial_{i_1}\cdots\partial_{i_p}G)(f)
                   \prod_{\nu=1}^p W_\nu f_{i_\nu},
 \qquad \sum_{\nu=1}^p|W_\nu|=a,                     \label{mat:face-9}
\end{equation}
where $W_\nu$ is a composition of $|W_\nu|$ derivatives.
The sum also includes terms in which a derivative acts on an operator
coefficient.
The derivatives of $G$ are bounded on the positive compact factor
range. Thus \eqref{mat:face-9} has a polynomial logarithmic majorant. The largest
number of spatial coefficient derivatives is the largest number
landing on a single factor; auxiliary derivatives do not increase it.
All orders used are finite, even though they can depend on $K$.

We will use ordinary Taylor's formula in $\epsilon$, holding $L$
and the respective core or layer variable fixed:
\begin{equation}
 G(\epsilon)=\frac{\epsilon^{K+1}}{K!}
       \int_0^1(1-t)^K G^{(K+1)}(t\epsilon)\dd t
 \quad\hbox{if }G^{(n)}(0)=0\ (0\le n\le K).            \label{mat:face-10}
\end{equation}
We apply the material and spatial derivatives to the residual before
\eqref{mat:face-10}, so the differentiated remainder is given by the
same integral formula.

\subsubsection{Smooth extension of the core residual to zero auxiliary scale}\label{mat:face-part-3}

On the fixed edge collar write $w=q\omega(q)$, $z=1-q$.
Here $z$, $\omega$, and their reciprocals are bounded positive factors. Introduce independent $(\epsilon,L,q)$, and set
\begin{equation}
 V=\sum_{n=2}^K\sum_\ell
       \epsilon^nL^\ell q^{n-1}z\phi_{n\ell}(z),
 \quad u=z+qV,\quad
 D_c=\epsilon\partial_\epsilon+\partial_L,\quad
 Z_c=q\partial_q-\epsilon\partial_\epsilon.             \label{mat:face-11}
\end{equation}
Restriction to $\epsilon=\lambda/q,L=\log\lambda$ gives $u_c$.
The chain rule gives $D_c=\lambda\partial_\lambda|_q$ and
$Z_c=E_q|_\lambda$ on restricted functions; $[D_c,Z_c]=0$.
Consequently the physical Jacobian is
\begin{equation}
 v=u_z=-q^{-1}Z_cu=1-(1+Z_c)V,\qquad qv_z=-Z_cv.        \label{mat:face-12}
\end{equation}
The positive geometric factors used in the pressure force are
\begin{equation}
 a=(z^2/(u^2v))^{1/3},\qquad
 H=\sqrt{1+\epsilon^2/(\omega^2a^2)}.                  \label{mat:face-13}
\end{equation}

Let $\mathcal M_p$ denote the pressure contribution to
$\mathcal M$. The density is $\rho=\lambda^{-3}w^3a^3$, so the
enthalpy formula gives
\begin{equation}
 \lambda h=4\sqrt{\lambda^2+w^2a^2}-4\lambda,\qquad
 \lambda^2\mathcal M_p[A]=(\lambda h)_z/v.
 \qquad
 \frac{a_z}{a}=\frac23\left(\frac1z-\frac vu\right)
                         -\frac{v_z}{3v}.
                                                               \label{mat:face-14}
\end{equation}
Substitution, using \eqref{mat:face-12}, gives
\begin{equation}
 \mathcal P_c
 =\frac{4a}{vH}
    \left[-\omega-q\omega_q+
        \frac{2q\omega}{3}\left(\frac1z-\frac vu\right)
                         +\frac{\omega}{3v}Z_cv\right].          \label{mat:face-15}
\end{equation}
In particular there is no negative power of $q$ in \eqref{mat:face-15}.
The reduced residual in these auxiliary variables is
\begin{equation}
 F_c=\beta^2\big((1+D_c)^2-\tfrac32(1+D_c)\big)u
       +e\epsilon qD_c(1+D_c)u+\mathcal P_c+m(z)/u^2.              \label{mat:face-16}
\end{equation}
Its restriction is exactly \eqref{mat:face-2}, with the full
Chandrasekhar pressure law.

By \eqref{mat:face-6}, every $q^{n-1}z\phi_{n\ell}$ in \eqref{mat:face-11} has polynomial
logarithmic conormal bounds through 23. For
$\epsilon_0=\lambda/q$, $q\ge c\sqrt\lambda$, one has
\begin{equation}
 \epsilon_0\le c^{-1}\sqrt\lambda,\qquad
 |\log q|\le C(1+|\log\lambda|).
                                                               \label{mat:face-17}
\end{equation}
Uniformly for $0\le\epsilon\le\epsilon_0$, \eqref{mat:face-11}--\eqref{mat:face-12} give
\begin{equation}
 |V|+|(1+Z_c)V|
       \le C\epsilon_0^2(1+|\log\lambda|)^{p_{\log}}=o(1).
                                                               \label{mat:face-18}
\end{equation}
Thus $u\ge z_0/2$, $v\ge1/2$, and $a,a^{-1},H,H^{-1}$
stay in fixed bounded positive intervals throughout the Taylor segment.

The calculus of Subsubsection~\ref{mat:face-part-2} applied to \eqref{mat:face-11}--\eqref{mat:face-16} now proves
\begin{equation}
 \sup_{0\le t\le1}
 \left|\partial_\epsilon^{K+1}
             \mathcal W F_c(t\epsilon_0,L,q)\right|
 \le C(1+|L|+|\log q|)^{p_{\log}},\qquad |\mathcal W|\le21,       \label{mat:face-19}
\end{equation}
where $\mathcal W$ is a composition of $|\mathcal W|$
derivatives from $D_c,Z_c$, and $L=\log\lambda$.
For the derivative count, $v$ uses one conormal derivative of
the coefficient in $V$, $Z_cv$ uses two, and \eqref{mat:face-16} uses two material
derivatives. Applying 21 further mixed derivatives therefore uses at most
$21+2=23$ coefficient spatial derivatives. Terms obtained when
$\partial_\epsilon^{K+1}$ differentiates the coefficients of $D_c$ or $Z_c$
are still ordinary derivatives of the same finite auxiliary family;
they cost no additional coefficient spatial order.

\subsubsection{Pressure derivatives in the layer at large density}\label{mat:face-part-4}

Write $Z'_{n\ell}=\partial_DZ_{n\ell}$, and put
\begin{equation}
 h_0=Z_0(Q^4)/Q,\quad
 h_{n\ell}=Z_{n\ell}(Q^4)/Q^{n+1},\quad
 k_{n\ell}=Q^{-n}Z'_{n\ell}(Q^4)/Z'_0(Q^4).
                                                               \label{mat:face-20}
\end{equation}
The leading density is the exact positive solution of
\begin{equation}
 P(\rho_0)=\varkappa_\delta Q^4/\pi,\qquad
 Z'_0=(4\pi\rho_0)^{-1}.                               \label{mat:face-21}
\end{equation}
It follows that $h_0$ is bounded with bounded conormal derivatives.
Also $Q^3 Z'_0$ is a positive conormal factor, so \eqref{mat:face-6} and
\begin{equation}
 Z'_{n\ell}=(4Q^4)^{-1}E_Q Z_{n\ell}                  \label{mat:face-22}
\end{equation}
give polynomial logarithmic bounds for $h_{n\ell}$ through 23
and $k_{n\ell}$ through 22.

Introduce independent $(\epsilon,L,Q)$, and define
\begin{equation}
 u=1-\epsilon h_0
       -\sum_{n=1}^K\sum_\ell\epsilon^{n+1}L^\ell h_{n\ell},
 \quad
 R=1+\sum_{n=1}^K\sum_\ell\epsilon^nL^\ell k_{n\ell},
 \quad \theta=u^{-2}R^{-1},                            \label{mat:face-23}
\end{equation}
\begin{equation}
 T_v=\partial_L-Q\partial_Q,\qquad
 Z_v=Q\partial_Q+\epsilon\partial_\epsilon.
                                                               \label{mat:face-24}
\end{equation}
Restriction to $\epsilon=\lambda Q,L=\log\lambda$ gives the
actual $u_v,Z_D/Z'_0,\rho/\rho_0$. Moreover $T_v$ is exactly
$\lambda\partial_\lambda|_x$ and $Z_v=E_Q|_\lambda$.
In particular $T_v\epsilon=0$, and $[T_v,Z_v]=0$.
For $\epsilon_0=\lambda Q$, $1\le Q\le C_0\lambda^{-1/2}$,
\begin{equation}
 \epsilon_0\le C_0\sqrt\lambda,\qquad
 \log Q\le C(1+|\log\lambda|).
                                                               \label{mat:face-25}
\end{equation}
Equations \eqref{mat:face-20}--\eqref{mat:face-23} imply uniformly on $0\le\epsilon\le\epsilon_0$
that $u=1+o(1)$, $R=1+o(1)$. Hence $\theta,\theta^{-1}$
are in a fixed compact positive interval on the whole Taylor segment.

Uniform bounds for the auxiliary derivatives follow from the
density-derivative estimates for the exact pressure. Let
$\mathscr E_\rho=\rho\partial_\rho$.
The explicit law
\begin{equation}
 P(\rho)=\tfrac12\left[
 \rho^{1/3}(2\rho^{2/3}-3)\sqrt{1+\rho^{2/3}}
                       +3\operatorname{arsinh}(\rho^{1/3})\right]
                                                               \label{mat:face-26}
\end{equation}
implies, for every fixed integer $a\ge0$ and $\rho\ge c>0$,
\begin{equation}
 |\mathscr E_\rho^aP(\rho)|\le C_{a,c}\rho^{4/3}.       \label{mat:face-27}
\end{equation}
At $\rho\ge1$ factor the algebraic term by
$\rho^{4/3}$; its remaining factors are smooth functions of
$\rho^{-2/3}\in[0,1]$. The inverse-hyperbolic-sine term is
\begin{equation}
 \tfrac13\log\rho+
       \log\!\left(1+\sqrt{1+\rho^{-2/3}}\right).
 \label{mat:face-28}
\end{equation}
The logarithm is bounded by $C\rho^{4/3}$, and each of its
positive-order iterates under $\mathscr E_\rho$ is bounded.
On $[c,1]$, smoothness gives \eqref{mat:face-27}.
The identity
\begin{equation}
 \rho^aP^{(a)}(\rho)
       =\prod_{\ell=0}^{a-1}(\mathscr E_\rho-\ell)P(\rho)
                                                               \label{mat:face-29}
\end{equation}
gives the corresponding ordinary density derivatives.

The leading equation \eqref{mat:face-21} gives
$\rho_0\sim Q^3$. In its logarithmic implicit equation,
the coefficient of the highest derivative is
\[
 \frac{\rho_0P'(\rho_0)}{P(\rho_0)},
\]
which is bounded above and below for $Q\ge1$. Successive
differentiation therefore gives
\begin{equation}
 |E_Q^r\log\rho_0|\le C_r\quad(r\ge1),\qquad
 |E_Q^r(\rho_0/Q^3)|\le C_r.                           \label{mat:face-30}
\end{equation}
At each stage, the remaining terms are finite products of lower
logarithmic derivatives and bounded logarithmic derivatives of
$P$. They are controlled by \eqref{mat:face-27} and
$P(\rho)\sim\rho^{4/3}$ on $\rho\ge c$, closing
the finite induction.

Define
\begin{equation}
 p(\epsilon,L,Q)=Q^{-4}P(\rho_0(Q)\theta(\epsilon,L,Q)).
                                                               \label{mat:face-31}
\end{equation}
Before composition with $\theta$, the density derivatives satisfy
\[
 \partial_\theta^a[Q^{-4}P(\rho_0\theta)]
 =Q^{-4}\theta^{-a}
   \left.\prod_{\ell=0}^{a-1}(\mathscr E_\rho-\ell)P(\rho)
   \right|_{\rho=\rho_0\theta}.
\]
For $\theta$ in the fixed positive interval, Leibniz's rule and
\eqref{mat:face-30} therefore give
\[
 \begin{aligned}
 \left|E_Q^r\partial_\theta^a
          [Q^{-4}P(\rho_0\theta)]\right|
 &\le C_{r,a}Q^{-4}
             \sum_{j=0}^{r+a}|\mathscr E_\rho^jP(\rho_0\theta)|\\
 &\le C_{r,a}Q^{-4}(\rho_0\theta)^{4/3}\le C_{r,a}.
 \end{aligned}
\]
Here the empty product for $a=0$ is the identity. Composing with
\eqref{mat:face-23} and applying \eqref{mat:face-9} yields the
polynomial logarithmic bounds for all mixed derivatives used below.

The reduced residual of the auxiliary family is
\begin{equation}
 F_v=\beta^2\big((1+T_v)^2-\tfrac32(1+T_v)\big)u
       +e(\epsilon/Q)T_v(1+T_v)u
       -4\pi u^2\big(p+\tfrac14Z_vp\big)
       +(M-\epsilon^4)/u^2.                            \label{mat:face-32}
\end{equation}
This follows because
$\partial_DP(\rho)=p+\frac14E_Qp$ at fixed $\lambda$,
and $\lambda^4D=\epsilon^4$.
At $\epsilon=0$, $u=R=1$ and
$p=\varkappa_\delta/\pi$; thus
\begin{equation}
 F_v(0)= -\beta^2/2-4\varkappa_\delta+M=0              \label{mat:face-33}
\end{equation}
by $\delta=-\beta^2/2$ and $M+\delta=4\varkappa_\delta$.

The preceding bounds prove the exact differentiated remainder bound
\begin{equation}
 \sup_{0\le t\le1}
 \left|\partial_\epsilon^{K+1}
               \mathcal W F_v(t\epsilon_0,L,Q)\right|
 \le C(1+|L|+\log Q)^{p_{\log}},\qquad |\mathcal W|\le21,         \label{mat:face-34}
\end{equation}
for compositions of derivatives from $T_v,Z_v$. In the pressure term,
$R$ uses one
coefficient derivative through \eqref{mat:face-22}; applying $Z_v$
uses one more. Hence at most 23 coefficient derivatives are used.
The two material derivatives in the inertial term give the same
bound. Density derivatives of $P$ may have higher finite order
depending on $K$; they are controlled by
\eqref{mat:face-26}--\eqref{mat:face-29} and require no additional
derivatives of the profile coefficients.

\subsubsection{Cancellation of the auxiliary Taylor coefficients}\label{mat:face-part-5}

At a fixed core point $z<1$, the enthalpy formula gives a smooth
function of the auxiliary scale and the finite profile variables.
By \eqref{mat:face-14}, its Taylor coefficients agree with those
obtained from the scaled pressure in the coefficient equations.
The explicit $\lambda^4\log\lambda$ term has a spatially constant
coefficient and vanishes under spatial differentiation.
Once $L$ is treated as independent, the force therefore has no
additional logarithmic dependence on the scale.

More concretely, in the polynomial ring in $L$, the scale operator
satisfies
\begin{equation}
 (t\partial_t+\partial_L)(t^nL^\ell)
       =t^n(nL^\ell+\ell L^{\ell-1}).                  \label{mat:face-35}
\end{equation}
This is the rule used to extract the source coefficients in
\eqref{mat:core-log-recursion}. If all earlier
coefficients have been inserted, the core coefficient of
$t^nL^\ell$ is
\begin{equation}
 -zF_{n\ell}+zB_n\phi_{n\ell}=0              \label{mat:face-36}
\end{equation}
by \eqref{mat:core-log-recursion}, with $F_{n\ell}$ its full right-hand side after moving the higher-logarithmic terms. Higher logarithmic coefficients at the same degree
are already inserted because the block is solved in decreasing
logarithmic degree. Lower scale degrees are not altered.
The leading profile cancels degree zero. The degree-one radius
coefficient and its source vanish by the chosen normalization, and
\eqref{mat:face-36} applies to $2\le n\le K$.
Thus every coefficient through degree $K$ vanishes identically
as a function of $z$ and polynomially in $L$.

For the layer, the density variation produced by a new coefficient
$Z_n$ at its first residual degree is
\begin{equation}
 \delta\rho=-4\pi\rho_0^2 Z_n',\qquad
 -4\pi\partial_D(P'(\rho_0)\delta\rho)
                  =\partial_D(A_0Z_n'),
 \quad A_0=16\pi^2\rho_0^2P'(\rho_0).
                                                               \label{mat:face-37}
\end{equation}
Every other occurrence of that radius coefficient has a further
factor $\lambda$. The coefficient at degree $n,\ell$ is therefore
\begin{equation}
 -S_{n\ell}+\partial_D(A_0Z'_{n\ell})=0               \label{mat:face-38}
\end{equation}
by the layer equations in Subsection~\ref{mat:connection}.
The source in \eqref{mat:face-38} is the coefficient of the full
reduced residual with the preceding coefficients inserted, including
the $\partial_L$ terms in \eqref{mat:face-35}.
Equation \eqref{mat:face-33} supplies degree zero, and \eqref{mat:face-38} supplies degrees
$1,\ldots,K$.

In particular the retained $Z_K'$ is precisely
what cancels degree $K$.

At fixed $q$, the change $t=\epsilon q$ multiplies the $n$-th
core scale coefficient by $q^n$, without changing $L$.
At fixed $Q$, $t=\epsilon/Q$ multiplies the $n$-th layer
coefficient by $Q^{-n}$. Hence \eqref{mat:face-36}--\eqref{mat:face-38} imply the exact identities
\begin{equation}
 \partial_\epsilon^nF_c(0,L,q)=
 \partial_\epsilon^nF_v(0,L,Q)=0,\qquad 0\le n\le K.     \label{mat:face-39}
\end{equation}
These identities follow from the exact core and layer ODEs
\eqref{mat:face-36} and \eqref{mat:face-38}, before any endpoint
expansion or localization.

The auxiliary fields preserve the cancellations
\eqref{mat:face-39}. Indeed,
\[
 \partial_\epsilon^n(\epsilon\partial_\epsilon G)|_{\epsilon=0}
       =n\,\partial_\epsilon^nG|_{\epsilon=0},
\]
and $\partial_L$, $q\partial_q$, and $Q\partial_Q$ act on
identically zero coefficient functions. Thus, in the respective
variables, every composition $\mathcal W$ in \eqref{mat:face-19}
or \eqref{mat:face-34} satisfies
\[
 \partial_\epsilon^n\mathcal W F_c(0,L,q)
 =\partial_\epsilon^n\mathcal W F_v(0,L,Q)=0,
 \qquad 0\le n\le K.
\]
Taylor's formula \eqref{mat:face-10}, with the bounds
\eqref{mat:face-19} and \eqref{mat:face-34}, proves
\eqref{mat:face-3}--\eqref{mat:face-4} for the auxiliary derivatives.

To recover $\partial_s=b\lambda\partial_\lambda|_x$, use
\begin{equation}
 \partial_s^j f=\sum_{r=0}^j c_{jr}(\lambda)T^rf,\quad
 c_{00}=1,\quad
 c_{j+1,r}=\partial_sc_{jr}+bc_{j,r-1},                 \label{mat:face-40}
\end{equation}
where $T=D_c$ in the core and $T=T_v$ in the layer, with coefficients
outside $0\le r\le j$ set to zero. For fixed $\beta>0$,
$b=\sqrt{\beta^2+e\lambda}$ and its required derivatives are bounded,
so the recurrence bounds every $c_{jr}$ and its required material
derivatives. They depend only on $\lambda$; spatial derivatives
therefore commute with them. Every mixed derivative of total order
at most 21 is consequently a bounded linear combination of the
auxiliary derivatives just estimated. This proves
\eqref{mat:face-3}--\eqref{mat:face-4}.

\subsubsection{The compact core and the residual at the center}\label{mat:face-part-6}

On a fixed compact interior interval, $w,u,u_z$ have positive
lower bounds. Replace $\lambda$ by an auxiliary $t$, keep
$L$ independent, and use \eqref{mat:face-14} with
$\lambda h=4\sqrt{t^2+w^2a^2}-4t$.
This is a smooth function of the finite coefficient
family and at most two spatial derivatives, on the stated positive range. Mixed material/spatial
order 21 therefore uses coefficient order at most 23.
Equations \eqref{mat:face-35}--\eqref{mat:face-36}, followed by \eqref{mat:face-10}, prove \eqref{mat:face-5} there.

At the center, use the fixed mass coordinate
$r_c=x^{1/3}$, $X\in\mathbb R^3$, $|X|=r_c$.
The normalized radial map is
\begin{equation}
 U(t,L,X)=a(t,L,|X|^2)X,\qquad
 J=D_XU,\quad B=J^{-T},\quad
 c_0=\frac3{4\pi\det J}.                              \label{mat:face-41}
\end{equation}
Here $a=a_0+\sum_{n=2}^K\sum_\ell t^nL^\ell a_{n\ell}$,
with the regular scalar coefficient functions furnished by the fixed
center construction. The leading tangential and radial eigenvalues
of $J$ are positive. The finite perturbation is
$O(t^2(1+|L|)^{p_{\log}})$ in the required low coefficient norm.
Thus for $0\le t\le\lambda$, $L=\log\lambda$, the entire
segment has a uniformly positive determinant and inverse matrix.

The full normalized vector residual, including its gravity term, is
\begin{equation}
 \begin{aligned}
 \mathbf F_c={}&
 \beta^2\big((1+D_t)^2-\tfrac32(1+D_t)\big)U
   +etD_t(1+D_t)U\\
 &+B\nabla_X\!\left(4\sqrt{t^2+c_0^{2/3}}-4t\right)
   +\frac{X}{a^2},\qquad D_t=t\partial_t+\partial_L.
 \end{aligned}    \label{mat:face-42}
\end{equation}
Indeed the density is $t^{-3}c_0$, the gradient transforms by
$B$, and $x/u^2$ lifted radially is
$r_c/a^2\,(X/r_c)=X/a^2$.
Formula \eqref{mat:face-42} is defined at $X=0$ as written; no term is divided by
$|X|$. It agrees with the exact scalar radial equation away from
zero, hence everywhere by continuity.

The matrix inverse and determinant depend smoothly on the
nonsingular deformation matrix $D_XU$. The gradient in \eqref{mat:face-42} requires
one further derivative. Thus every ordinary Cartesian derivative
of order $m$, together with $j$ auxiliary material derivatives,
$j+m\le21$, is a finite smooth composition using at most 23
spatial coefficient derivatives. The finite regularity in the
squared-radius coefficients in \eqref{mat:face-41} supplies these Cartesian
derivatives directly by the chain rule. No individual singular
angular or radial coefficient is estimated.

The Taylor coefficients of \eqref{mat:face-42} through degree $K$ vanish on
$X\ne0$ by \eqref{mat:face-36}, and at $X=0$ by their continuous extension.
Their Cartesian derivatives also vanish, because they are zero
functions on the entire center ball. The same composition estimates
therefore bound the $t^{K+1}$ Taylor integral in the full
Cartesian $C^m$ norm. Applying \eqref{mat:face-40} proves \eqref{mat:face-5} for the
center-vector lift and completes the residual estimates
\eqref{mat:face-3}--\eqref{mat:face-5}.

The differentiated residuals of the core and layer radii are now
controlled. To estimate the joined profile, we bound the error introduced
by the matching cutoff and then pass to the weighted source norm,
including its scale factors and localization derivatives.

\subsection{The moving cutoff and the mixed source norm}
\label{app:profile-cutoff}

Use the joined radius, reduced residual, and source in
\eqref{mat:cutoff-1}--\eqref{mat:cutoff-5}. Here
$D=\lambda\partial_\lambda|_z$, $E=q\partial_q$, and
$D_v=(M-x)/\lambda^4$, as in Subsection~\ref{mat:cutoff}.
The overlap estimate \eqref{mat:cutoff-2} controls the mismatch.
We prove the nonlinear cutoff estimate and then the mixed source bound
\eqref{mat:cutoff-6}, with the localization of Definition~\ref{at:atlas}.

\subsubsection{Proof of the nonlinear cutoff estimate}\label{mat:cutoff-part-2}

The pressure estimate uses the mean-value formula along the
interpolating radii, so positivity of their Jacobians is needed on
the whole segment. The inertial error is the commutator of the
cutoff with the time operator.

For a shell function define its order-$r$ conormal bound to be
the supremum over the enlarged shell of the sum of
$|D^aE^bf|$, $a+b\le r$. Denote it temporarily by
$|f|_r$. Both fields are derivations and commute.
The finite Leibniz rule gives
\begin{equation}
 |fg|_r\le C_r|f|_r|g|_r.                                  \label{mat:cutoff-7}
\end{equation}
If $G$ is smooth on a neighborhood of a fixed compact set, its
derivatives through $r+1$ are bounded there. Repeated chain
rules and
\begin{equation}
 G(U)-G(V)=\int_0^1DG(V+t(U-V))(U-V)\dd t
 \label{mat:cutoff-8}
\end{equation}
give $|G(U)-G(V)|_r\le C|U-V|_r$ whenever the whole segment
lies in that set and the components of both $U$ and $V$
have bounded order-$r$ norms. The constant depends only on these finite derivative bounds and on
the compact range of the interpolating segment.

The local estimates for the finite family in Subsection~\ref{mat:units} and
Subsection~\ref{mat:native} imply, on this shell,
\begin{equation}
 u_\nu,\quad v_\nu=(u_\nu)_z,\quad q(v_\nu)_z
 \quad\hbox{have bounded order-21 conormal norms},\qquad
 u_\nu,v_\nu\ge c>0.                                      \label{mat:cutoff-9}
\end{equation}
For a derivative of $v_\nu$ or $q(v_\nu)_z$, write it first
using $\partial_z=-q^{-1}E$. At most 23 conormal derivatives of the radius are used.
The coefficient bounds provide 60 derivatives. Material differentiation acts on the finite coefficient sums before
their norms are taken.

Put $h=u_c-u_v$. Since
\[
 \partial_z=-q^{-1}E,\qquad
 \partial_z^2=q^{-2}(E^2-E),
\]
and $Dq=0$, $Eq=q$, \eqref{mat:cutoff-2} gives
\begin{equation}
 |h|_{21}+|h_z|_{21}+|q h_{zz}|_{21}
 \le C\varepsilon_\lambda/q_{\min},\qquad
 q_{\min}=c_1\sqrt\lambda.                                \label{mat:cutoff-10}
\end{equation}
Every derivative of $\chi(q/\sqrt\lambda)$ under $D,E$
is bounded; for example
\[
 D\chi=-\tfrac12\zeta\chi_{\rm m}'(\zeta),\quad
 E\chi=\zeta\chi_{\rm m}'(\zeta),\quad \zeta=q/\sqrt\lambda.
\]
The same calculation proves \eqref{mat:cutoff-10} with $h$ replaced by
$\chi h$. The derivative
\[
 (u_a)_z=(1-\chi)v_v+\chi v_c+\chi_z h
\]
is bounded below, since the first two terms are at least $c$
and the last is $O(\varepsilon_\lambda/q_{\min})=o(1)$.

Every interpolant $u_v+t h$ and $u_v+t\chi h$, $0\le t\le1$,
has the same bounds after decreasing $A_*$. This also proves
\eqref{mat:cutoff-9} for those interpolants.

\emph{The pressure term.} We use the exact enthalpy. For an arbitrary
interpolant let
\[
 a=(z^2/(u^2v))^{1/3},\qquad
 H=\sqrt{1+(\lambda/w)^2a^{-2}}.
\]
The density and pressure force obey the exact identities
\[
 \rho=\lambda^{-3}w^3a^3,\qquad
 \lambda h(\rho)=4\sqrt{\lambda^2+w^2a^2}-4\lambda,
\]
\begin{equation}
 \mathcal P[u]:=\lambda^2\mathcal M_{\rm p}[\lambda u]
 =\frac{4a}{vH}
 \left[w'+\frac{2w}{3}\left(\frac1z-\frac vu\right)
                         -\frac{w}{3v}v_z\right].
 \label{mat:cutoff-11}
\end{equation}
\equationalias{mat:material-4.2}{mat:cutoff-11}
Indeed $\mathcal M_{\rm p}=h_z/(\lambda v)$ and
$a_z/a=(2/z-2v/u-v_z/v)/3$, proving \eqref{mat:cutoff-11} directly.

Each of the functions $z,z^{-1},w,w',w/q,\lambda/w$ has bounded
conormal derivatives through order 21. For the last function,
$D(\lambda/w)=\lambda/w$ and
$E(\lambda/w)=-(\lambda/w)E\log w$, while
$\lambda/w=O(\sqrt\lambda)$. In particular $H\ge1$;
there is no small divisor from the constitutive law.
Expression \eqref{mat:cutoff-11} is a smooth function of the quantities just listed and
\begin{equation}
                       u,\ v,\ qv_z.                     \label{mat:cutoff-12}
\end{equation}
It is affine in the third argument. Applying \eqref{mat:cutoff-7}--\eqref{mat:cutoff-10} to \eqref{mat:cutoff-11}
therefore gives
\begin{equation}
 |\mathcal P[u_v+\chi h]-\mathcal P[u_v]|_{21}
 +|\mathcal P[u_c]-\mathcal P[u_v]|_{21}
 \le C\varepsilon_\lambda/q_{\min}.                       \label{mat:cutoff-13}
\end{equation}
Subtracting $\chi$ times the second difference proves the same
bound for the exact pressure cutoff defect. The mean-value formula includes the full nonlinear pressure term.
The zeroth-order map $m(z)/u^2$ satisfies \eqref{mat:cutoff-13} without the
factor $q_{\min}^{-1}$.

\emph{The inertial term.} Exact scale differentiation gives the
normalized inertial operator
\[
 I_\lambda=\beta^2(D^2+\tfrac12D-\tfrac12)
                              +e\lambda(D^2+D).
\]
Consequently
\begin{equation}
 \begin{split}
 I_\lambda(u_v+\chi h)-(1-\chi)I_\lambda u_v-\chi I_\lambda u_c
 ={}&2(\beta^2+e\lambda)(D\chi)Dh\\
 &+\left[(\beta^2+e\lambda)D^2\chi
             +(\tfrac12\beta^2+e\lambda)D\chi\right]h .
 \end{split}                                             \label{mat:cutoff-14}
\end{equation}
\equationalias{mat:material-4.4}{mat:cutoff-14}
This has order 21 conormal bound $C\varepsilon_\lambda$;
it uses at most 22 mismatch derivatives.
Equations \eqref{mat:cutoff-13}--\eqref{mat:cutoff-14} and gravity give the sharper defect
$C\varepsilon_\lambda/\sqrt\lambda=C\lambda^{K/2}$
times logarithmic factors. Since $\lambda\le1$, this implies
\eqref{mat:cutoff-3}. This estimate concerns the cutoff error on the
joining shell; the uncut residuals are controlled by \eqref{mat:cutoff-4}.

\subsubsection{Estimates in the covering source norm}\label{mat:cutoff-part-3}

Use a fixed physical mass partition
$\chi_c+\chi_i+\chi_e=1$, with $\chi_e=1$ on
$M-x\le d_0$ and supported on $M-x<2d_0$.
Each cutoff is a fixed power at least 24 of a smooth function.
Split $\chi_e$ by the fixed power-normalized transition
$\zeta(D_v)$, equal to 1 near $D_v=0$ and zero for $D_v\ge2$:
\[
 \chi_5=\chi_e\zeta(D_v),\qquad
 \chi_8=\chi_e(1-\zeta(D_v)),\qquad
 y_n=D_v^{1/n},\quad n=5,8.
\]
The center uses $x^{1/3}$ with its full radial-vector lift;
the interior chart is fixed and compact. Define
\begin{equation}
 \begin{split}
 \mathcal S_m(g)^2={}&
 \|\chi_cg\,e_r\|_{H^m(B^3;\mathbb R^3)}^2+
 \|\chi_i g\|_{H^m(I)}^2\\
 &+\lambda^6\sum_{a=0}^m\int
       |\partial_{y_5}^a(\chi_5g)|^2y_5^4\dd y_5\\
 &+\lambda^6\sum_{a=0}^m\int
       |\partial_{y_8}^a(\chi_8g)|^2y_8^7\dd y_8 .
 \end{split}                                             \label{mat:cutoff-15}
\end{equation}
Every displayed derivative acts on the cutoff as well.
This is the source norm of Definition~\ref{at:atlas}, with the same
localization cutoffs.

The long chart has $c<y_8<C_e\lambda^{-1/2}$. On its layer
part $Q=y_8^2$, so
\begin{equation}
 \partial_y^a
 =y^{-a}\prod_{\ell=0}^{a-1}(2E_Q-\ell),\qquad
 \partial_s=b(\lambda\partial_\lambda-E_Q).
 \label{mat:cutoff-16}
\end{equation}
On the core portion of this chart, the smooth inverse of the mass
deficit gives
$q=c_*\lambda y_8^2U(\lambda y_8^2)$, where $U$ is a
positive smooth fixed factor on a bounded interval. Writing
$X=(M-x)^{1/4}=\lambda y_8^2$, we have
$y_8\partial_{y_8}=2\sigma(q)E$, where
$\sigma=\dd\log q/\dd\log X$ is bounded, positive, and smooth.
Repeated differentiation expresses every $\partial_{y_8}^a$
as $y_8^{-a}$ times a finite combination of at most $a$
iterates of $E$, with bounded coefficients. The material
derivative there is $bD$. Because $y_8\ge c>0$, no
negative power of $\lambda$ is lost in either conversion.
All mixed terms use no more than $j+m\le21$ operations.

On the inner cutoff transition $D_v\sim1$ the local
cutoff derivatives are bounded. At the outer physical transition
$d=M-x\sim d_0$, the formula $d=\lambda^4y_8^8$ gives
\begin{equation}
 |\partial_{y_8}^a\chi_e(\lambda^4y_8^8)|
 \le C_a\lambda^{a/2}.
 \label{mat:cutoff-17}
\end{equation}
The finite Leibniz sums in \eqref{mat:cutoff-15} therefore preserve every local
supremum bound used here. The material derivative is applied to $g$ before the cutoffs in
\eqref{mat:cutoff-15}; the localized norm is not differentiated in time.

On the long chart, combine \eqref{mat:cutoff-3}--\eqref{mat:cutoff-4}:
all required chart derivatives of the reduced residual are bounded by
\[
 C\lambda^{(K-1)/2}(1+|\log\lambda|)^{p_{\log}}.
\]
Use the layer estimate up to the joining shell and the core estimate
in \eqref{mat:cutoff-4} on the rest of the chart. Since
$\mathcal R_u=\lambda^{-3}\mathcal F[u_a]$, the corresponding
local bound for $\mathcal R_u$ is
$C\lambda^{(K-7)/2}$ times the same logarithmic factor.

To pass from a local supremum bound $H_\lambda$ to the strong
source norm, use
\begin{equation}
 C\lambda^3
 \left(\int_c^{C_e\lambda^{-1/2}}y^7\dd y\right)^{1/2}
 H_\lambda
 \le C\lambda H_\lambda.                                \label{mat:cutoff-18}
\end{equation}
The long-chart contribution therefore has power
$\lambda^{(K-5)/2}$.

On the bounded $y_5$ chart, Subsection~\ref{mat:endpoint} gives
$C\lambda^{K-2}$ before the factor $\lambda^3$ in the source norm,
hence a contribution $C\lambda^{K+1}$. The compact center
and interior estimates of Subsection~\ref{mat:face} give
$C\lambda^{K-2}$. For $K\ge1$,
\[
 K+1\ge\frac{K-5}{2},\qquad K-2\ge\frac{K-5}{2}.
\]
Thus both contributions satisfy the long-chart bound obtained from
\eqref{mat:cutoff-18}. On every
region the logarithmic factor is bounded by a fixed power of
$1+|\log\lambda|$.

The scalar terms require no additional
derivatives. Let $T$ be the fixed-mass scale derivative
in its respective chart, and define $p_{00}=1$, with indices
outside the triangle zero, by
\[
 p_{j+1,k}=\partial_sp_{jk}+b p_{jk}+b p_{j,k-1}.
\]
Then, exactly,
\begin{equation}
 \lambda^{-4}\partial_s^j(\lambda^4\mathcal R_u)
  =\lambda^{-3}\sum_{k=0}^jp_{jk}T^k\mathcal F[u_a].
 \label{mat:cutoff-19}
\end{equation}
\equationalias{mat:material-5.2}{mat:cutoff-19}
Every coefficient is bounded for the selected $\beta>0$ and
$j\le21$. Spatial derivatives do not act on it.
Equations \eqref{mat:cutoff-15}--\eqref{mat:cutoff-19} prove the residual part of \eqref{mat:cutoff-6}.

\subsubsection{The viscosity and shift terms}\label{mat:cutoff-part-4}

The differentiated viscosity term in \eqref{mat:cutoff-5} is
\begin{equation}
 -\lambda^{-4}\partial_s^j
       \{\nu\lambda^3\partial_s\mathcal M[A]\}.
 \label{mat:cutoff-20}
\end{equation}
Each Leibniz summand has a factor bounded by $C\nu\lambda^{-1}$
and a force derivative of order $r\le j+1$.
For $j+m\le21$, precisely $r+m\le22<25$.
The local force bound and the ordinary derivative bound with the same
cutoffs, proved in Subsection~\ref{mat:native}, are
$C\lambda^{-2}$. Since $\mathcal S_m(g)\le C\|g\|_{m,a}$
with these same cutoffs and $\lambda\le1$, \eqref{mat:cutoff-20} is bounded by
\begin{equation}
                       C\kappa\lambda^{-9/2}.            \label{mat:cutoff-21}
\end{equation}
Here $\partial_s^a\nu=O(\nu)$ and
$\partial_s^a\lambda^3=O(\lambda^3)$; these follow directly
from $\nu=\kappa\lambda^{-3/2}$ and the bounded derivatives
of $b$. They contain no negative extra scale power.

The shift is fixed along a terminal trajectory. Its source is
\[
 -\Lambda\lambda^{-4}\partial_s^j(\nu\lambda^3 A_s).
\]
It uses radius derivatives $r+m\le22$, bounded by
$C\lambda$ in the same local estimates. Its source norm is therefore
bounded by $C\kappa\Lambda\lambda^{-3/2}$.
Summing over the indicated indices proves \eqref{mat:cutoff-6}.
All differentiated quantities belong to the prescribed profile.

If $K\ge2J+86$, $J\ge3$, $0<\lambda_*\le\lambda\le A_*$,
$\Lambda\le B_*$, and
\[
 \kappa\le c\lambda_*^{J+5},\qquad
 \kappa B_*\le c\lambda_*^{J+2},
\]
the last two terms of \eqref{mat:cutoff-6} are at most
$Cc\lambda^{J+1/2}$. Moreover
\[
 (K-5)/2-(J+1/2)\ge40.
\]
Consequently, for every $\delta>0$, choosing $A_*$ sufficiently
small for the selected finite family bounds the residual term by
$\delta\lambda^{J+1/2}$. The leading constants fixed before
$\beta$ are unchanged.

This proves the strong source estimate through mixed total order $21$.
The energy additionally uses the pure material derivative of order $22$ in the mass
norm. Its proof retains one more overlap derivative and uses a different
edge prefactor, so we supply it separately.

\Needspace{6\baselineskip}

\subsection{The highest material derivatives of the source}
\label{app:profile-material}

We prove \eqref{mat:material-1.2}--\eqref{mat:material-1.3} under
the hypotheses of Subsection~\ref{mat:material}. In particular,
$K\ge86$, the overlap has 24 differentiated bounds, and the
same finite profile and cutoffs are used. We retain the scale fields
$D,E$ and the normalized residual of that subsection. The final
integration uses $H=L^2((0,M),\dd x)$, rather than the mixed source norm.

\subsubsection{The reduced residual through material order 22}\label{mat:material-part-4}

Use the auxiliary families of
Subsection~\ref{app:profile-core-layer}, now through material order 22.
The expansion parameters and material derivatives are
\[
 \begin{aligned}
 \epsilon&=\lambda/q,&
 T&=t\partial_t+\partial_L&&\text{in the core},\\
 \epsilon&=\lambda Q,&
 T&=\partial_L-Q\partial_Q&&\text{in the layer}.
 \end{aligned}
\]
On both regions, $\epsilon\le C\sqrt\lambda$, and the positive
factors remain uniformly bounded away from zero for
$0\le t\le\epsilon$. The force uses two coefficient derivatives;
the local coefficient bounds and
Subsubsection~\ref{mat:connection-part-6} supply all
$22+2=24$ derivatives required here.

Each field $T$ preserves the auxiliary scale degree. Hence the
coefficient cancellations through degree $K$ persist for the
differentiated residual $H_j$, $0\le j\le22$:
\[
 \partial_t^nH_j(0,L,\cdot)=0,\qquad 0\le n\le K.
\]
Taylor's formula gives
\begin{equation}
 \frac{\epsilon^{K+1}}{K!}
   \int_0^1(1-\theta)^K
           \partial_t^{K+1}H_j(\theta\epsilon,L,\cdot)\dd\theta,
 \label{mat:material-4.1}
\end{equation}
where $H_j$ is the residual after the stated material differentiations.
The integrand is bounded by a finite logarithmic polynomial on
the full positive segment. Its many $t$ derivatives differentiate
finite polynomial amplitudes and smooth reciprocals of the positive factors; they do
not add $K$ spatial derivatives to a coefficient function.

Thus the core and layer bounds in that proof hold for these 22 material
operations, with the same powers $(\lambda/q)^{K+1}$ and
$(\lambda Q)^{K+1}$. On the compact center/interior the ordinary
parameter Taylor calculation has power $\lambda^{K+1}$, and
uses the full Cartesian vector at the center. At the bounded
vacuum, Subsection~\ref{mat:endpoint} supplies the derivatives through
material order 22 in $C^{28}$.

We next estimate the additional derivative of the cutoff.
Use the exact enthalpy pressure expression \eqref{mat:cutoff-11}.
On the shell this is a smooth function of positive bounded factors and the three arguments $u,v,qv_z$, affine in the third.
The full segment between the two profiles and the joined profile
is positive by the low-order part of \eqref{mat:material-1.1}. These arguments have 22
bounded material derivatives. For $h=u_c-u_v$, \eqref{mat:material-1.1} gives
\begin{equation}
 \sum_{j\le22}(|D^jh|+|D^jh_z|+|D^j(qh_{zz})|)
 \le C\lambda^{-1/2}\lambda^{(K+1)/2}(1+|\log\lambda|)^{p_{\log}}.
 \label{mat:material-4.3}
\end{equation}
Here $h_z=-q^{-1}Eh$, $qh_{zz}=q^{-1}(E^2-E)h$.
At $j=22$, the last term uses the mixed derivative of total order
24 in \eqref{mat:material-1.1}. Replacing $h$ by $\chi h$ has the same estimate.
The finite chain rule and the mean-value formula in these arguments
give a pressure cutoff defect of size $C\lambda^{K/2}$ up
to logarithms, after every $D^j$, $j\le22$. Gravity has
no derivative loss.

The exact inertial cutoff defect is \eqref{mat:cutoff-14}.
Its 22 material derivatives use at most 23 derivatives of $h$,
already included in \eqref{mat:material-1.1}. Consequently the reduced
residual of the joined profile and its material derivatives through order 22 are bounded on the
long chart by
\begin{equation}
 C\lambda^{(K-1)/2}(1+|\log\lambda|)^{p_{\log}}.
 \label{mat:material-4.5}
\end{equation}
We use this weaker exponent to combine the cutoff and uncut residual
bounds. The layer estimate applies up to the joining shell, and the
core estimate controls the outer portion of the long chart.

\subsubsection{The exact mass norm and scalar terms}\label{mat:material-part-5}

To obtain the $H$ bound in \eqref{mat:material-1.2}, we integrate
against the exact mass measures. Since $x=M-\lambda^4y_n^n$,
the positive measure on an edge chart is
$|\dd x|=n\lambda^4y_n^{n-1}\dd y_n$. For a local supremum bound
$B_\lambda$ on the long chart, its mass norm is at most
\begin{equation}
 C\lambda^2
       \left(\int_c^{C_e\lambda^{-1/2}}y^7\dd y\right)^{1/2}
       B_\lambda\le C B_\lambda.
 \label{mat:material-5.1}
\end{equation}
Since $\mathcal R_u=\lambda^{-3}\mathcal F$, \eqref{mat:material-4.5} gives
$B_\lambda=C\lambda^{(K-7)/2}(1+|\log\lambda|)^{p_{\log}}$.
The factor in \eqref{mat:material-5.1} is $\lambda^2$, whereas
the mixed source norm has factor $\lambda^3$.

The bounded vacuum contribution is at most
$C\lambda^2\lambda^{K-2}=C\lambda^K$, by \eqref{mat:endpoint-3}.
The compact contributions are at most $C\lambda^{K-2}$,
and the regular full-vector center norm controls its mass norm
with a fixed geometric factor. Each is smaller than \eqref{mat:material-5.1}.
The partition sums to one, so the triangle inequality for the
four localized mass norms gives the global $H$ bound. No
spatial derivative of a localized norm or moving measure is taken
in this argument.

The scalar terms are exactly \eqref{mat:cutoff-19}, now with
$j\le22$. Its finite coefficient recurrence is unchanged and all
coefficients remain bounded for the selected $\beta>0$.
Together with \eqref{mat:material-4.5} and
\eqref{mat:material-5.1}, this proves \eqref{mat:material-1.2}.

For viscosity, the differentiated force term is
\begin{equation}
 -\lambda^{-4}\partial_s^j
                  (\nu\lambda^3\partial_s\mathcal M[A]).
 \label{mat:material-5.3}
\end{equation}
Since $|\partial_s^l(\nu\lambda^3)|\le C_l\nu\lambda^3$,
Leibniz's rule and $\|f\|_H\le C\|f\|_{0,a}$ give, for $j\le22$,
\[
 \begin{aligned}
 \left\|\lambda^{-4}\partial_s^j
          (\nu\lambda^3\partial_s\mathcal M[A])\right\|_H
 &\le C\nu\lambda^{-1}
                \sum_{r=1}^{j+1}\|\partial_s^r\mathcal M[A]\|_{0,a}\\
 &\le C\kappa\lambda^{-5/2}\lambda^{-2}
       =C\kappa\lambda^{-9/2}.
 \end{aligned}
\]
The force derivatives have order $j+1\le23$, within
\eqref{mat:native-eq:wave2-direct-profile-force-source}.
Only the order-zero norm is used, so no positive-order compatibility
trace of the prescribed source is required. For the fixed shift,
the radius bounds similarly yield
\[
 \begin{aligned}
 \left\|\Lambda\lambda^{-4}\partial_s^j
                      (\nu\lambda^3 A_s)\right\|_H
 &\le C\Lambda\nu\lambda^{-1}
                  \sum_{r=1}^{j+1}\|\partial_s^rA\|_H\\
 &\le C\kappa\Lambda\lambda^{-5/2}\lambda
       =C\kappa\Lambda\lambda^{-3/2}.
 \end{aligned}
\]
This proves \eqref{mat:material-1.3}.

The order selected in Lemma~\ref{ct:parameter-admissibility} leaves
a positive power between the residual in \eqref{mat:material-1.3}
and the target $\lambda^{J+1/2}$. Choosing the same two bounds on viscosity
$\kappa\le c\lambda_*^{J+5}$,
$\kappa\Lambda\le c\lambda_*^{J+2}$ therefore makes the
right side $(\varepsilon+Cc)\lambda^{J+1/2}$ on
$[\lambda_*,A_0]$, after decreasing $A_0$ for any fixed
$\varepsilon>0$. This is the pure material source estimate at energy order 22.

This completes the forcing estimates for the finite family used in
Proposition~\ref{mat:matched-source}. The leading strain constants
and the final profile conclusions are stated in
Subsection~\ref{mat:exports}.

\section{Local solutions of the regularized equation}
\label{app:local-theory}

We prove the local existence and continuation result of
Proposition~\ref{loc:high-local}, then the approximation and
finite-regularity statements of Lemma~\ref{loc:finite-strong-graph}
and Proposition~\ref{loc:finite-calculus}. The mass interval and local
charts are those of Subsection~\ref{loc:domains}. The viscosity
$\kappa>0$ and the finite shift $\Lambda$ are fixed throughout
the local construction. Its constants may depend on these parameters
and on the positive terminal scale; the uniform estimates used in the
limiting argument are proved in Sections~\ref{rec:chapter}--\ref{ct:section}.

\subsection{Linear estimates for the regularized equation}
\label{loc:base-inverse-section}

Use the local time and parabolic spaces of
Subsection~\ref{loc:local-solutions} on the fixed charts of
Subsection~\ref{loc:domains}. The local inverses are based on Lunardi's
scalar oblique-boundary Schauder theorem
\cite[Theorem~3.1]{Lunardi2004}, applied in the autonomous case on a
smooth ball or interval with the outward normal derivative as boundary
operator. Its source class
is $C^{\alpha,\theta}$, its initial class is $C^{2,\alpha}$,
and its boundary-data class is
$C^{1+\alpha,(1+\alpha)/2}$. The compatibility condition is
$\partial_n u_0=g(0)$. Every local inverse below has $u_0=g=0$,
so this condition is automatic and the source has no corner restriction.
On subintervals of a fixed time interval, its Schauder constant has a
common upper bound. A contraction on functions with zero initial data
then gives the inverse for time-dependent coefficients.

The homogeneous Neumann condition is imposed only at the artificial
outer boundary of an enlarged chart. We multiply the local solutions by
cutoffs before extending them and summing over the charts. Near the
center we first solve an equation for general vector fields and then
restrict to equivariant fields. At the vacuum we lift a one-sided radial
function to a five-dimensional ball; the physical boundary becomes the
origin of that ball.

We construct the frozen and nonautonomous inverses in the ordinary
anisotropic Hölder spaces, with uniform bounds at exponent
$\alpha$. Approximation at the end of the subsection then restricts
them to the closed little spaces $\mathcal P_T$ and $\mathcal F_T$.

\begin{lemma}
\label{loc:base-lift}
The radial lift $f(\xi)\mapsto f(|X|)$ has the following properties.
\begin{enumerate}[label=(\roman*),leftmargin=2em]
\item It identifies $\{f\in h^{2,\alpha}[0,\ell]:f'(0)=0\}$ with
the radial subspace of $h^{2,\alpha}(B^5_\ell)$. The norms are equivalent
after adjoining the two quotients in \eqref{loc:quotients}.
\item The same correspondence holds in the stated base parabolic norms.
\item It maps $h^\alpha[0,\ell]$ onto the radial source space
$h^\alpha(B^5_\ell)$, without a derivative trace condition.
\end{enumerate}
\end{lemma}

\begin{proof}
Write $U(X)=f(|X|)$. Put $q_f=f'/\xi$, $e_f=f''-q_f$ and
$e=X/|X|$. Since $e_f(0)=0$,
\[
 \partial_iU=f'e_i,\qquad
 \partial_{ij}U=q_f\delta_{ij}+e_fe_ie_j.
\]
For $g\in C^\alpha[0,\ell]$ with $g(0)=0$, the angular factor satisfies
\[
 [g(|X|)e_ie_j]_{C^\alpha(B^5_\ell)}
       \le C[g]_{C^\alpha[0,\ell]}.
\]
Indeed, set $d=|X-Y|$. If $d\ge\tfrac12\min(|X|,|Y|)$, both
radii are at most $3d$ and $|g(r)|\le[g]_\alpha r^\alpha$.
Otherwise the radii are comparable, and the angular contribution is
bounded by
\[
 C[g]_\alpha |Y|^{\alpha-1}d\le C[g]_\alpha d^\alpha;
\]
the radial contribution is bounded by $[g]_\alpha d^\alpha$.
Together with \eqref{loc:quotients}, this proves the upper domain
norm bound. Restriction to a ray proves the reverse bound. The source
estimate follows from $\big||X|-|Y|\big|\le|X-Y|$.

For the little spaces, approximate $f''$ by smooth one-sided functions
and $e_f$ by smooth functions vanishing at zero. The corresponding
angular products are Lipschitz at the origin, hence little-$\alpha$;
the preceding estimate proves convergence. Conversely, average smooth
Cartesian approximants over rotations and restrict to a ray.

In the parabolic norm, the temporal $\theta$-seminorm of $D^2U$ is
controlled by those of $f''$ and $q_f$. The gradient has temporal
exponent $(1+\alpha)/2$, and $U_t$ is the radial lift of $f_t$.
The integral formulas \eqref{loc:quotients} bound these components,
including the source norm of $f_t$. Smooth approximation in the little
anisotropic norms proves the parabolic correspondence. In particular,
$\xi^3$ is allowed: its lift $|X|^3$ belongs to $C^{2,1}$. No
higher even trace conditions are imposed.
\end{proof}

The coefficients depend only on the following derivatives of $R$:
Cartesian derivatives through order two at the center, ordinary
derivatives through order two in the interior, and
\[
 R,\quad R_\xi/\xi,\quad R_{\xi\xi}-R_\xi/\xi
\]
at the vacuum. Denote this collection by $\mathcal J R$.
Throughout the local construction, assume
$\|\mathcal J R\|_{\mathcal F_T}\le K_0$, the positive lower
bounds $\eta$, and fixed upper bounds on these coefficients.
No third derivative of the radius is used. We write
$B_v=\xi\partial_\xi A$ for the vacuum coefficient in
\eqref{loc:vacuum-hessian}, with $A$ as in
\eqref{loc:vacuum-coefficients}.

\begin{lemma}
\label{loc:coefficient-approximation}
On each fixed enlarged chart, the coefficient and source fields admit
smooth approximants converging in the little anisotropic source norm.
They can be chosen so that:
\begin{enumerate}[label=(\roman*),leftmargin=2em]
\item the principal coefficients remain uniformly elliptic and all
coefficient bounds in the source norm are uniform;
\item covariance at the center, radial symmetry at the vacuum, and
$B_v(t,0)=0$ are preserved;
\item the corresponding operators converge in
$\mathcal L(\mathcal P_T,\mathcal F_T)$.
\end{enumerate}
The approximating operators have scalar principal part and the stated
covariance; they need not be Hessians of radius fields.
\end{lemma}

\begin{proof}
The coefficients are smooth functions of $\mathcal J R$ on a fixed
positive range. Extend them past each time endpoint by constants and
across the artificial spatial boundary by reflection; on a sphere use
\[
 F((\ell+s)e)=F((\ell-s)e).
\]
After multiplication by a fixed outer cutoff, convolution gives smooth
approximants converging in the little anisotropic source norm. Indeed,
extension and convolution are bounded in that norm, preserve the
vanishing H\"older moduli, and converge on its dense smooth subspace.

At the center, symmetrize the principal matrix and average each tensor
over rotations, acting on all its indices. For a matrix and a vector the
averages are
\[
 \int_{SO(3)}O^Ta(OX)O\dd O,\qquad
 \int_{SO(3)}O^TF(OX)\dd O.
\]
These bounded projections fix the original covariant fields and preserve
smoothness. Uniform convergence preserves at least half the original
ellipticity; rotation averaging preserves this lower bound. Apply the
same projection to all three indices of the lower-order tensor, retaining
its coupling between vector components.

At the vacuum, average the scalar approximants radially on the lifted
five-ball. For the approximation $\widetilde B_n$ of $B_v$, set
\[
 B_n(t,X)=\widetilde B_n(t,X)-\widetilde B_n(t,0),\qquad
 B_nQ=\frac{B_n}{|X|^2}X\cdot\nabla,\qquad Qu=u_\xi/\xi.
\]
Evaluation at zero is bounded into $h^\theta([0,T])$. Since
$B_v(t,0)=0$,
\[
 \|B_n-B_v\|_{\mathcal F_T}
 \le C\|\widetilde B_n-B_v\|_{\mathcal F_T}\longrightarrow0.
\]
Smooth radiality and $B_n(t,0)=0$ give
$B_n=|X|^2\widetilde b_n$ with smooth $\widetilde b_n$. Thus the
approximating quotient term is a smooth drift; the limiting operator
uses $Q:\mathcal P_T\to\mathcal F_T$ from
Lemma~\ref{loc:base-lift}. Uniform convergence preserves ellipticity.
On an ordinary interval, convolution alone suffices.

Keep the zero initial and artificial normal data unchanged. To preserve
interior support of a localized source, smooth on the enlarged chart and
multiply by a fixed cutoff equal to one on its original support. No
first-derivative trace is imposed on a source at the physical vacuum.
The same operations apply to positive averaged secants: covariance and
$B_v(t,0)=0$ are linear conditions, and their principal lower bound is
preserved by averaging.

Let $\varepsilon_n$ be the sum of the coefficient errors in the full
anisotropic source norm, including $B_n-B_v$. The product estimate and
$\|Qu\|_{\mathcal F_T}\le C_Q\|u\|_{\mathcal P_T}$ give
\begin{equation}
 \|(L_n-L)u\|_{\mathcal F_T}
       \le C\varepsilon_n\|u\|_{\mathcal P_T},
       \qquad \varepsilon_n\longrightarrow0.
 \label{loc:operator-approximation}
\end{equation}
In particular, the quotient term is bounded by
\[
 \|(B_n-B_v)Qu\|_{\mathcal F_T}
 \le C\|B_n-B_v\|_{\mathcal F_T}\|u\|_{\mathcal P_T}.
\]
Thus both spatial and temporal coefficient seminorms enter the estimate.
The constants may depend on the fixed chart and positive viscosity.
\end{proof}

\begin{lemma}
\label{loc:local-chart-inverses}
Fix $\kappa>0$, finite $\Lambda$, and the coefficient bounds above.
On each enlarged chart, for every source $F\in\mathcal F_T$, the problem
\[
 u_t+\kappa(H_R+\Lambda)u=F,\qquad
 u(0)=0,\qquad \partial_nu=0\ \text{on the artificial boundary}
\]
has a unique solution $u\in\mathcal P_T$ for a common sufficiently
small $T>0$. The inverse bound and this time are uniform over the stated
coefficient class, including positive averaged secants of the Hessian.

No compatibility condition on $F(0)$ is required. At the center, the
inverse acts on general vector fields before restriction to equivariant
fields.
\end{lemma}

\begin{proof}
We construct the inverse in ordinary H\"older spaces, then pass to the
little spaces by approximation.

\emph{Scalar inverse.}
For $u(0)=0$ and $T\le1$,
\begin{equation}
 \|u\|_{\mathcal F_T}\le CT^{1-\theta}\|u\|_{\mathcal P_T},
 \qquad
 \|Du\|_{\mathcal F_T}\le CT^{1/2}\|u\|_{\mathcal P_T}.
 \label{loc:zero-initial-gain}
\end{equation}
The first estimate follows by integration in time. For the second use
\[
 \|Du\|_{C^\alpha}
 \le C\|u\|_{C^\alpha}^{1/2}\|u\|_{C^{2,\alpha}}^{1/2}
       +C\|u\|_{C^\alpha},\qquad
 [Du]_{t;\theta}\le T^{1/2}[Du]_{t;(1+\alpha)/2}.
\]
The spatial estimate holds up to the boundary: the extension
\[
 Eu((1+s)e)=6u((1-s)e)-8u((1-2s)e)+3u((1-3s)e)
\]
matches derivatives through order two. With a fixed cutoff it is bounded
at orders $0,\alpha,2+\alpha$ and preserves zero initial value.

Write
\[
 L(t)=\partial_t-a^{ij}(t,X)\partial_{ij}
            +d^i(t,X)\partial_i+c(t,X),
\]
with fixed ellipticity and coefficient norms at most $K$ in
$C^{\alpha,\theta}$. If
$\|h\|_\infty\le\varepsilon$ and
$[h]_{x;\alpha}+[h]_{t;\theta}\le2K$, then $D^2u(0)=0$ gives
\[
 \|D^2u\|_\infty\le T^\theta[D^2u]_{t;\theta},\qquad
 \|hD^2u\|_{C^{\alpha,\theta}}
 \le C(\varepsilon+KT^\theta)
                         \|u\|_{C^{2+\alpha,1+\theta}}.
\]
Here each seminorm uses
$[hD^2u]\le\|h\|_\infty[D^2u]+[h]\|D^2u\|_\infty$.

Let $S^{\rm fr}$ be the zero-initial, homogeneous Neumann inverse for
the operator frozen at time zero. The scalar Schauder theorem gives a
bound on a fixed upper time interval; bounded extension of the source
and restriction of the solution give the same bound on shorter intervals.
For frozen operators whose principal coefficients differ by at most
$\varepsilon$ in $C^0$, the preceding product bound and
\eqref{loc:zero-initial-gain} yield
\[
 \|(L_1-L_2)u\|_{\mathcal F_T}
 \le C\{\varepsilon+KT^\theta+KT^{1/2}+KT^{1-\theta}\}
                                      \|u\|_{\mathcal P_T}.
\]
The coefficient class has compact closure in $C^0$, with the same
ellipticity and ordinary $C^\alpha$ bounds. A Neumann series therefore
bounds each inverse in a $C^0$ neighborhood after decreasing $T$.
A finite subcover gives one inverse bound and one time interval for
the class. The closure need not consist of little-H\"older coefficients;
this step takes place in ordinary H\"older spaces.

Since $\|a(t)-a(0)\|_\infty\le KT^\theta$,
\begin{equation}
 \|(a-a(0))D^2u\|_{C^{\alpha,\theta}}
       \le CK T^\theta\|u\|_{C^{2+\alpha,1+\theta}}.
 \label{loc:scalar-time-freezing}
\end{equation}
The lower-order differences cost
$CK(T^{1/2}+T^{1-\theta})$. Hence, for
$B=S^{\rm fr}(L-L^{\rm fr})$, choose the common time so that
\[
 \|B\|_{\mathcal L(C^{2+\alpha,1+\theta})}
 \le C_{S^{\rm fr}}CK(T^\theta+T^{1/2}+T^{1-\theta})\le\tfrac12.
\]
Then
\[
 u=\sum_{n=0}^\infty(-B)^nS^{\rm fr}F,\qquad
 \|u\|_{C^{2+\alpha,1+\theta}}
       \le2C_{S^{\rm fr}}\|F\|_{C^{\alpha,\theta}}.
\]
Every term has zero initial value and zero normal derivative. The
constants depend only on the stated coefficient bounds, ellipticity
and smooth chart domain.

\emph{Vacuum chart.}
Differentiating \eqref{loc:vacuum-force} gives
\begin{equation}
 H_Rf=-A(f_{\xi\xi}+4f_\xi/\xi)-B_vf_\xi/\xi+bf,\qquad
 B_v=\xi\partial_\xi A,\quad
 b=A_R(3q-\xi q_\xi)+(B_0)_R.
 \label{loc:vacuum-hessian}
\end{equation}
Partial subscripts on $A,B_0$ hold $\xi,q$ fixed. The coefficients
$A,A^{-1},B_v,b$ are smooth functions of the base fields, with
$B_v(t,0)=0$. One cannot treat $B_v/\xi$ as a bounded drift:
for $R=R_b-a_*\xi^2+c_*\xi^{2+\gamma}$, $\alpha<\gamma<1$,
the identity $A_q/A=-8/(3q)$ at zero gives
$B_v=C_*\xi^\gamma+o(\xi^\gamma)$.

Set $\xi=\ell y$, $t=t_0+\ell^2\tau/\kappa$. On the unit
five-ball the extra radial operator is $B_{v,\ell}Q$, $Qf=f_y/y$,
and
\[
 \|Qf\|_{\mathcal F}\le C_Q\|f\|_{\mathcal P},\qquad
 \|B_{v,\ell}\|_{\alpha,\theta}
 \le C\ell^\alpha\bigl([B_v]_{x;\alpha}
                       +\kappa^{-\theta}[B_v]_{t;\theta}\bigr).
\]
The scalar operator
$\partial_\tau-A_\ell\Delta_5+\ell^2(b_\ell+\Lambda)$
has common ellipticity and source bounds for $0<\ell\le1$.
Its zero-initial Neumann inverse $S_0$ thus has a common bound and
time interval before $\ell$ is decreased. Choose the collar so that
\[
 C_{S_0}C_Q\|B_{v,\ell}\|_{\alpha,\theta}\le\tfrac12.
\]
Then $f=S_0(F+B_{v,\ell}Qf)$ is a contraction. Rotational
uniqueness gives a radial solution. The physical vacuum is the
interior origin of this ball; only its artificial outer sphere carries
a Neumann condition, and no boundary condition is imposed on $F$.

\emph{Center chart.}
Use \eqref{loc:center-force}, with $c^2=P'(\rho)$, $H=h(\rho)$,
and
\[
 \operatorname{div}_BF=B_{ij}\partial_jF_i,\qquad
 (\operatorname{curl}_BF)_i=\varepsilon_{ijk}B_{j\ell}\partial_\ell F_k.
\]
The identities $\delta B=-B(DF)^TB$ and
$\delta\rho=-\rho\operatorname{div}_BF$ give the completion
\begin{equation}
 \widehat H_\Phi F
 =-B(DF)^TB\nabla H-B\nabla(c^2\operatorname{div}_BF)
       -2r^{-3}F+c^2\operatorname{curl}_B\operatorname{curl}_BF.
 \label{loc:center-completion}
\end{equation}
For radial $F$, the matrices $DF,B$ are commuting combinations of
$I,e\otimes e$, so $\operatorname{curl}_BF=0$ and
$\widehat H_\Phi(fe)=H_Rf\,e$. Moreover,
\[
 (\operatorname{curl}_B\operatorname{curl}_BF)_i
 =B_{ak}\partial_k(B_{i\ell}\partial_\ell F_a)
       -B_{ak}\partial_k(B_{a\ell}\partial_\ell F_i).
\]
Substitution cancels the mixed-component second derivatives and yields
\begin{equation}
 (\widehat H_\Phi F)_i
 =-c^2(B^TB)_{k\ell}\partial_{k\ell}F_i
                +C_{ia\ell}\partial_\ell F_a-2r^{-3}F_i,
 \label{loc:center-scalar-principal}
\end{equation}
where
\begin{align*}
 C_{ia\ell}={}&-B_{i\ell}(B\nabla H)_a
 -(B_{ik}\partial_kc^2)B_{a\ell}
 -c^2B_{ik}\partial_kB_{a\ell}\\*
 &+c^2B_{ak}\partial_kB_{i\ell}
 -c^2\delta_{ia}B_{bk}\partial_kB_{b\ell}.
\end{align*}
These are base source coefficients involving only $D\Phi,D^2\Phi$:
\[
 \partial_kB=-B(\partial_kJ)^TB,\quad
 \partial_k\rho=-\rho\,\operatorname{tr}(J^{-1}\partial_kJ),\quad
 \partial_kH=(c^2/\rho)\partial_k\rho.
\]
The bounds $\eta\le R_z,R/z\le L$ imply
\[
 c^2\zeta^TB^TB\zeta
 \ge L^{-2}\min_{3/(4\pi L^3)\le v\le3/(4\pi\eta^3)}P'(v)
          |\zeta|^2>0.
\]
Apply the scalar inverse componentwise to the principal part, with
homogeneous normal data on the artificial sphere. By
\eqref{loc:zero-initial-gain}, the lower matrix terms perturb it by
an operator of norm $C(T^{1/2}+T^{1-\theta})<1/2$ after a common
decrease of $T$. The inverse is unique and commutes with rotations.
For equivariant sources its solution is equivariant, hence radial:
the stabilizer of a nonzero point fixes only its radial line.

On interior intervals the scalar inverse applies directly. Averaging
the completed coefficients along a positive segment preserves their
principal lower bound, base bounds and, at the vacuum, $B_v(0)=0$.
Thus the same construction applies to secants.

\emph{Little H\"older spaces.}
Apply Lemma~\ref{loc:coefficient-approximation} simultaneously to the
coefficients and source. Their ordinary exponent-$\alpha$ inverses
have a common bound $C_0$, interval $[0,T]$ and vacuum collar.
For each fixed approximant choose $\alpha<\alpha'<1$. Its vacuum
coefficient satisfies $B_n=|X|^2\widetilde b_n$ with smooth radial
$\widetilde b_n$, so the full quotient term is a smooth drift.
The scalar construction at exponent $\alpha'$ therefore applies
on the fixed collar, as does the center matrix perturbation.

These higher-exponent solutions extend to $[0,T]$. Indeed the smooth
coefficient norms on this interval give a positive step length. At a
restart time $t_0$, subtract the stationary value $u(t_0)$:
\[
 L\{u-u(t_0)\}=F-Lu(t_0),\qquad
 (u-u(t_0))|_{t=t_0}=0,\qquad
 \partial_n(u-u(t_0))=0.
\]
The source is $C^{\alpha',\alpha'/2}$ and has no corner restriction;
the step length is independent of its norm. Finitely many steps
cover $[0,T]$. The equation identifies $u_t$ at each join, while the
common $C^{2,\alpha'}$ trace identifies the spatial derivatives.
Thus $u\in C^{2+\alpha',1+\alpha'/2}$ throughout the interval.

Such a solution belongs to the little exponent-$\alpha$ domain.
Use even reflection in normal boundary coordinates and tangential
partitions constant in the normal direction at the boundary. Product
smoothing, even in that direction, preserves the zero normal trace;
bounded time extension and smoothing handle the time endpoints.
The gap $\alpha'-\alpha>0$ gives convergence in the base norm.
Rotation averaging preserves covariance. Ordinary uniqueness
identifies this solution with the exponent-$\alpha$ inverse. Its
higher norms and number of restarts may depend on the approximant.

For two approximating solutions, the common inverse bound gives
\begin{equation}
 \|u_n-u_m\|_{\mathcal P_T}
 \le C_0\bigl(\|F_n-F_m\|_{\mathcal F_T}
       +\|L_n-L_m\|_{\mathcal L(\mathcal P_T,\mathcal F_T)}
                                     \|u_m\|_{\mathcal P_T}\bigr).
 \label{loc:little-inverse-limit}
\end{equation}
Since $\|u_m\|_{\mathcal P_T}\le C_0\|F_m\|_{\mathcal F_T}$,
the sequence is Cauchy. The little space, equivariance, radial
symmetry, and zero initial and normal traces are closed in this norm.
The limit solves the equation, including its bounded vacuum quotient.
The same bounds hold on shorter intervals because every perturbation
factor $T^\theta,T^{1/2},T^{1-\theta}$ decreases. The collar and
rescaling constants may depend on fixed $\kappa>0$; no
zero-viscosity estimate is used here.
\end{proof}

\begin{proposition}
\label{loc:global-base-inverse}
Under the coefficient and positivity bounds above, there are
$T_*,C_*>0$, depending only on those bounds, fixed
$\kappa,\Lambda,M,\alpha$ and the atlas, such that for $0<T\le T_*$,
\begin{equation}
 \begin{gathered}
 u\longmapsto
 (u_t+\kappa(H_R+\Lambda)u,u(0)):
 \mathcal P_T\longrightarrow\mathcal F_T\times\mathcal D
 \quad\hbox{is an isomorphism},\\
 \|u\|_{\mathcal P_T}\le C_*(\|F\|_{\mathcal F_T}+\|u_0\|_{\mathcal D}).
 \end{gathered}
 \label{loc:global-inverse}
\end{equation}
Here $F=u_t+\kappa(H_R+\Lambda)u$ and $u_0=u(0)$.
The same assertion holds for positive averaged secants and bounded
scalar shifts in the source class.
\end{proposition}

\begin{proof}
Choose a finite enlarged chart cover and fixed mass cutoffs
$\sum_j\rho_j=1$, $\psi_j=1$ near $\operatorname{supp}\rho_j$.
Both cutoffs vanish near the artificial edges, are constant near the
physical endpoints, and have derivative supports in regular interior
overlaps. Let $\mathcal T_j$ be the chart/lift map and $S_j$ its
zero-initial Neumann inverse. Define
\[
 \mathcal PF=\sum_j\mathcal T_j^{-1}
          \{\psi_jS_j\mathcal T_j(\rho_jF)\}.
\]
Each local source is supported away from its artificial boundary.
Extend the solution by zero only after multiplication by $\psi_j$;
the center symmetry and vacuum lift put it in the global domain.

For $L=\partial_t+\kappa(H_R+\Lambda)$, covariance gives
\[
 L\mathcal P=I+E,\qquad
 EF=\sum_j\mathcal T_j^{-1}[L_j,\psi_j]
                              S_j\mathcal T_j(\rho_jF).
\]
The commutators have order one:
\[
 [L_j,\psi]u=-2a^{k\ell}\psi_k u_\ell
                   -a^{k\ell}\psi_{k\ell}u+d^k\psi_ku.
\]
At the vacuum this is
\[
 -\kappa\{2A\psi'f'+[A\psi''+(4A+B_v)\psi'/\xi]f\},
\]
whose quotient coefficient is supported away from zero. Hence
\eqref{loc:zero-initial-gain} implies
\[
 \|\mathcal P\|\le C_P,\qquad \|E\|\le C_ET^{1/2}.
\]
Also $EF(0)=0$, since the localized solutions and their spatial
derivatives vanish initially. Every further localization again has
support away from the artificial boundaries. For $C_ET^{1/2}\le1/2$,
\[
 S=\mathcal P\sum_{n=0}^\infty(-E)^n,\qquad
 LS=I,\qquad \|S\|\le2C_P.
\]

If $Lu=0$, $u(0)=0$, localization gives
\[
 L_j(\rho_j\mathcal T_ju)=[L_j,\rho_j]\mathcal T_ju,
\]
with zero artificial normal data. The local estimates therefore yield
\[
 \|u\|_{\mathcal P_T}
 \le C\sum_j\|\rho_j\mathcal T_ju\|_{C^{2+\alpha,1+\theta}}
 \le C'T^{1/2}\|u\|_{\mathcal P_T}.
\]
Decrease $T$ so that $C'T^{1/2}<1$; then $u=0$.

For $u_0\in\mathcal D$, add its stationary lift globally:
\[
 u(t)=u_0+S\{F-\kappa(H_R(t)+\Lambda)u_0\},\qquad
 u_t(0)=F(0)-\kappa(H_{R(0)}+\Lambda)u_0\in\mathcal Y.
\]
The source belongs to $\mathcal F_T$. All local solves still use
zero initial data and supported sources, so their artificial Neumann
conditions restrict neither $u_0$ nor the vacuum first derivative of
$u_t(0)$. The cutoffs, localized extensions, chart maps and convergent
series preserve the little spaces. The same estimates apply to positive
averaged secants and bounded scalar shifts in the source class.
\end{proof}

\subsection{The contraction argument for the velocity}
\label{loc:velocity-construction}

Fix an initial pair $(r_0,r_1)$ satisfying the bounds $L,\eta$ in
Subsection~\ref{loc:local-solutions}.

\emph{Coefficient estimates.}
On a fixed ball satisfying the positivity conditions, integration of
the coefficient derivatives along the segment from $S$ to $R$ gives
\begin{align}
 \|\mathcal M[R]-\mathcal M[S]\|_{\mathcal F_T}
     &\le C\|\mathcal J(R-S)\|_{\mathcal F_T},\notag\\
 \|(H_R-H_S)f\|_{\mathcal F_T}
     &\le C\|\mathcal J(R-S)\|_{\mathcal F_T}\|f\|_{\mathcal P_T}.
 \label{loc:base-lipschitz}
\end{align}
The product estimate applies to both spatial and temporal H\"older
differences. For $R^V=r_0+\int_0^tV$, the linear maps in
$\mathcal J$, including the closed quotients, commute with time
integration. Thus, with $p=1-\theta>0$,
\begin{equation}
 \sup_t\|R^V-r_0\|_{\mathcal D}\le T\|V\|_{\mathcal P_T},\qquad
 \|\mathcal J(R^V-R^W)\|_{\mathcal F_T}
                    \le CT^p\|V-W\|_{\mathcal P_T}.
 \label{loc:integrated-gain}
\end{equation}
The spatial $C^\alpha$ supremum gains $T$, and the temporal
$\theta$-seminorm in $C^0$ gains $T^{1-\theta}$. This uses no
spatial derivative of $V_t$.

\emph{The contraction.}
Let $S_*$ be the zero-initial inverse of
$L_*=\partial_t+\kappa(H_{r_0}+\Lambda)$ from
Proposition~\ref{loc:global-base-inverse}, and solve
\[
 L_*V_*=-\mathcal M[r_0],\qquad V_*(0)=r_1.
\]
That proposition gives, with the same bound on shorter intervals,
\[
 \|V_*\|_{\mathcal P_T}\le C_0(M_L+L),\qquad
 M_L=\sup\{\|\mathcal M[r]\|_{\mathcal Y}:
                   \|r\|_{\mathcal D}\le L,\ \text{positive factors}\ge\eta\}.
\]
Set $B=C_0(M_L+L)+1$. On the closed radius-one ball about $V_*$
with initial value $r_1$, consider
\begin{equation}
 V\longmapsto V_*+S_*\{-\mathcal M[R^V]+\mathcal M[r_0]
                     -\kappa(H_{R^V}-H_{r_0})V\}
 \label{loc:velocity-map}
\end{equation}
Its source has norm at most $CT^p(B+\kappa B^2)$ and Lipschitz
constant at most $CT^p(1+2\kappa B)$. For the quadratic term, use
\[
 \begin{aligned}
 &(H_{R^V}-H_{r_0})V-(H_{R^W}-H_{r_0})W\\
 &\qquad=(H_{R^V}-H_{r_0})(V-W)
                   +(H_{R^V}-H_{R^W})W.
 \end{aligned}
\]
Both coefficient differences have the time-integral gain in
\eqref{loc:integrated-gain}.

Choose the coefficient ball $\|\mathcal JR\|_{\mathcal F_T}\le C_DL+1$
with positive factors at least $\eta/2$, and fix its constants
$C_0,C_N$ before decreasing $T$. Require
\begin{equation}
 \begin{gathered}
 T\le\min(1,T_{\rm lin}),\quad TB\le1,\quad
 C_{\rm geo}TB\le\eta/2,\quad C_{\rm jet}T^pB\le1,\\
 C_0C_NT^p(B+\kappa B^2)\le1,\qquad
 C_0C_NT^p(1+2\kappa B)\le\tfrac12 .
 \end{gathered}
 \label{loc:base-ball-time}
\end{equation}
The first four restrictions keep $R^V$ in that coefficient ball and
preserve positivity. The last two make \eqref{loc:velocity-map} a
contraction of the closed velocity ball. Its fixed point solves
\eqref{loc:actual-velocity}. Subtracting two solutions and applying
the same estimates on successive short intervals proves uniqueness
in the base spaces.

\emph{Continuous dependence.}
For convergent initial pairs, keep the frozen operator and enlarge the
common ball slightly. The initial-radius error contributes
$C(1+\kappa B)\|r_{0,n}-r_0\|_{\mathcal D}$ to the source and
$C\kappa\|r_{0,n}-r_0\|_{\mathcal D}$ to its Lipschitz constant.
Choose $T$ so that the contraction constant is at most $1/4$, and then
$n$ large enough that the coefficient error is at most $1/4$. The
inverse estimate and subtraction of the fixed-point equations give
\[
 \begin{aligned}
 \|V_n-V\|_{\mathcal P_T}
 &\le C\bigl((1+\kappa B)\|r_{0,n}-r_0\|_{\mathcal D}
                    +\|r_{1,n}-r_1\|_{\mathcal D}\bigr)
       +\tfrac12\|V_n-V\|_{\mathcal P_T},\\
 R_n-R&=r_{0,n}-r_0+\int_0^t(V_n-V)\dd s.
 \end{aligned}
\]
Absorption proves convergence in $\mathcal P_T\times C_t^1\mathcal D$.
The common interval depends only on $L,\eta$, the atlas and the fixed
positive parameters.

\subsection{Finite spatial regularity}
\label{loc:elliptic-realization}

\begin{lemma}
\label{loc:elliptic-lemma}
For each finite $m\ge0$, under the fixed positivity bounds,
\begin{align}
 R\in\mathcal D,\quad \mathcal M[R]\in\mathcal Y^m
       &\ \Longrightarrow\ R\in\mathcal D^m,\notag\\
 R\in\mathcal D^m,\quad f\in\mathcal D,\quad H_Rf\in\mathcal Y^m
       &\ \Longrightarrow\ f\in\mathcal D^m.
 \label{loc:elliptic-implications}
\end{align}
In each case, the higher-regularity norm is controlled by the norms
on the left-hand side. Convergence in the latter norms implies
convergence in the former.
\end{lemma}

\begin{proof}
\emph{Nonlinear equation at the vacuum.}
For $F=\mathcal M[R]$, equation \eqref{loc:vacuum-force} becomes
\begin{equation}
 R_\xi=\xi q,\qquad
 \xi q_\xi=G(\xi,R,q,F),\qquad
 G=3q+(B_0-F)/A.
 \label{loc:fuchs-system}
\end{equation}
Since $\xi q_\xi=R''-q\to0$, we have $G(0)=0$. With
$s=\xi^2a^{2/3}$, the coefficient formulas give
\[
 A_q/A=-\frac{8+7s}{3q(1+s)},\qquad
 A_R/A=-\frac{4+2s}{3R(1+s)}.
\]
At zero, $B_{0,q}=0$ and $B_0-F=-3Aq$, hence
\begin{equation}
 G_q(0)=3+3qA_q/A=-5.
 \label{loc:stable-index}
\end{equation}
Put $e=G_q+5$. Its controlled H\"older norm and $e(0)=0$
allow a collar $[0,\ell]$ on which
$\|e\|_\infty+\ell^\alpha[e]_\alpha\le3$.

For $1\le k\le m$, assume
$q\in h^{k-1,\alpha}$, $R\in h^{k+1,\alpha}$ and
$F\in h^{k,\alpha}$. Differentiation on the punctured collar gives
\[
 \xi\partial_\xi q^{(k)}+(k+5-e)q^{(k)}=f_k,\qquad
 f_k=\partial_\xi^kG-G_qq^{(k)}\in h^\alpha.
\]
The last expression omits the highest derivative of $q$ from the
chain rule. For $p=(\xi,R,F)$, for example,
\[
 f_1=G_p\,p',\qquad
 f_2=G_{pp}[p',p']+2G_{pq}p'q'+G_{qq}(q')^2+G_p\,p''.
\]
Here $G_p$ is the derivative in the vector $p$. The operator
\[
 K_\lambda g(\xi)=\int_0^1t^{\lambda-1}g(t\xi)\dd t
\]
has norm at most $\lambda^{-1}$ in the scaled H\"older norm
(the seminorm gains $(\lambda+\alpha)^{-1}$). Thus
$W_k=K_{k+5}(f_k+eW_k)$ has a unique solution in $h^\alpha$,
with contraction factor $3/(k+5)\le1/2$.

To identify this solution, first note that the numerator in the
first-order equation for $q^{(k-1)}$ must vanish at zero: a nonzero
limit would give a $\xi^{-1}$ term in $q^{(k)}$, contradicting
continuity of $q^{(k-1)}$. The numerator is $C^\alpha$, so
$q^{(k)}=O(\xi^{\alpha-1})$. For $k=1$ this follows from
\eqref{loc:fuchs-system}. Set
\[
 \mu_k(\xi)=\exp\left(\int_\ell^\xi
                         \frac{k+5-e(r)}r\dd r\right).
\]
Since $k+5-e\ge k+2$,
\[
 0<\mu_k(\xi)\le(\xi/\ell)^{k+2},\qquad
 \mu_k(\xi)q^{(k)}(\xi)=O(\xi^{k+1+\alpha})\longrightarrow0.
\]
Integration from zero therefore yields
\[
 q^{(k)}(\xi)=\mu_k(\xi)^{-1}
      \int_0^\xi\mu_k(s)f_k(s)\frac{\dd s}{s},\qquad
 \|q^{(k)}\|_\infty\le\frac{\|f_k\|_\infty}{k+2}.
\]
The bounded solutions $q^{(k)}$ and $W_k$ differ by a multiple of
the unbounded function $\mu_k^{-1}$, so they agree. Integration
from a positive endpoint and passage to zero identifies this continuous
extension with the distributional derivative. Finally,
\[
 R''=q+G\in h^{k,\alpha},\qquad
 \frac{R-R(0)}{\xi^2}=\int_0^1tq(t\xi)\dd t\in h^{k,\alpha}.
\]
Induction gives order $m$. The system is nonsingular away from zero;
subtracting the contraction equations at each stage gives continuity
in the stated inputs and preserves the little H\"older classes.

\emph{Nonlinear equation at the center.}
For radial $g\in h^{k,\alpha}(B^3)$, put
\[
 Tg(X)=X\int_0^1t^2g(tX)\dd t.
\]
This operator gains one Cartesian derivative on smaller balls. To see
this, cut off $g$ radially and convolve with $X/(4\pi|X|^3)$.
The derivative kernel, in principal value, is
\[
 \frac1{4\pi}\left(\frac{\delta_{ij}}{|X|^3}
                   -\frac{3X_iX_j}{|X|^5}\right),
\]
with the additional term $\delta_{ij}\delta_0/3$. Its spherical
mean vanishes. Subtracting the source value and splitting the difference
of two such integrals at twice the distance $h$ gives
\[
 \int_0^{2h}r^{\alpha-1}\dd r
       +h\int_{2h}^{C}r^{\alpha-2}\dd r\le C_\alpha h^\alpha.
\]
The corresponding supremum bound and commutation of derivatives yield
$\|Tg\|_{C^{k+1,\alpha}(B')}
\le C\|g\|_{C^{k,\alpha}(B)}$ for $B'\Subset B$.
The convolution equals $Tg$ on $B'$: both are bounded radial fields
of divergence $g$, and the other radial solution is a multiple of
$|X|^{-2}e$. Cartesian mollification gives the same bound and
continuity in the little spaces.

For $\Phi=rX$, let $b_c=R_z=\operatorname{div}\Phi-2r$.
If $\Phi\in h^{k,\alpha}$, $k\ge2$, then
$r,b_c\in h^{k-1,\alpha}$, and the force and mass equations give
\begin{equation}
 \nabla H=b_c(\boldsymbol F-r^{-2}X),\qquad
 r^3X=\frac9{4\pi}T(\rho(H)^{-1}),\qquad H=h(\rho).
 \label{loc:center-realization-system}
\end{equation}
The second identity is the integral of
$\partial_z(R^3)=9z^2/(4\pi\rho)$. For
$\boldsymbol F\in h^{m,\alpha}$ and $s=\min(k,m+1)$ these
identities imply
\[
 H\in h^{s,\alpha},\qquad
 r^3\in h^{s,\alpha},\qquad Xr^3\in h^{s+1,\alpha}.
\]
Here the density stays positive, so $\rho(H)^{-1}$ is smooth on
its range. To pass from $r^3X$ to $rX$, use, for $|\beta|=s+1$,
\[
 X D^\beta u=D^\beta(Xu)
       -\sum_i\beta_i e_iD^{\beta-e_i}u.
\]
If $u\in h^{s,\alpha}$ and $Xu\in h^{s+1,\alpha}$, this
identity, with the usual coordinate vectors and unit multiindices,
controls the highest chain-rule term
$X\phi'(u)D^\beta u$; all other terms in $D^\beta(X\phi(u))$
use at most $s$ derivatives of $u$. Thus
$X\phi(u)\in h^{s+1,\alpha}$. The extended derivatives are
distributional derivatives across zero, by integration on segments
away from zero and passage to the limit. Apply this to
$u=r^3$, $\phi(u)=u^{1/3}$ on the fixed positive range.
Finitely many nested balls give $\Phi\in h^{m+2,\alpha}$,
with continuous bounds. The construction preserves the invariant
scalars and equivariant vectors. In the interior, solving
\eqref{loc:interior-force} for $R_{xx}$ gives the same induction.

\emph{Linear equation at the vacuum.}
The coefficient formulas satisfy
\[
 8A+3qA_q+(B_0)_q+\xi A_\xi|_{R,q}+\xi^2qA_R=0,
 \qquad (B_0)_q=-2\xi^2(A_qq^2+2Aq)/R.
\]
The displayed logarithmic derivatives cancel the first, second and
fourth terms; the remaining two cancel. Consequently
\[
 H_Rf=-\xi^{-4}\partial_\xi(\xi^4Af_\xi)+bf.
\]
For $F=H_Rf$, the base domain excludes the homogeneous flux
$|f_\xi|\sim\xi^{-4}$. Integration gives
\begin{equation}
 \frac{f_\xi}{\xi}
  =A^{-1}\int_0^1t^4(bf-F)(t\xi)\dd t,\qquad
 f''=A^{-1}(bf-F)-(4+\xi A_\xi/A)f_\xi/\xi.
 \label{loc:linear-vacuum-realization}
\end{equation}
For $R\in\mathcal D^m$ we have
$A^{\pm1},b,\xi A_\xi\in h^{m,\alpha}$, with bounds depending
on $\|R\|_{\mathcal D^m}$. These formulas improve
$f\in h^{k,\alpha}$ to $h^{\min(k,m)+2,\alpha}$, starting
at $k=2$. Iteration gives both quotient bounds and continuity.

\emph{Linear equation at the center.}
Write \eqref{loc:center-scalar-principal} as
\[
 (\mathcal Lf)_i=-a^{kl}\partial_{kl}f_i
                  +C_{ial}\partial_lf_a+D_{ia}f_a,
\]
where $a\in h^{m+1,\alpha}$ is uniformly positive and
$C,D\in h^{m,\alpha}$. For $m\ge1$, fix
$B_0\Subset B_1\Subset B_2$ inside the chart and a coordinate
direction $e$. The difference quotients
$w_h(X)=(f(X+he)-f(X))/h$ belong to $h^{2,\alpha}$ for each
$h\ne0$ and converge to $\partial_ef$ in $h^{1,\alpha}(B_1)$
by translation continuity. Their equation is
\begin{equation}
 \mathcal Lw_h=K_h:=\delta_h F
       +(\delta_h a^{kl})(\partial_{kl}f)(X+he)
       -(\delta_h C)(Df)(X+he)-(\delta_h D)f(X+he).
 \label{loc:center-quotient-source}
\end{equation}
Its source converges in $h^\alpha(B_1)$: the coefficients and
$F$ have the required first derivatives, and $D^2f$ is translation
continuous in $h^\alpha$.

Choose $\chi\in C_c^\infty(B_1)$ equal to one on $B_0$ and
$\zeta(\tau)$ vanishing near zero and equal to one on the last
third of a short inverse interval $[0,T_1]$. At fixed physical time,
$v_h=\zeta\chi w_h$ satisfies on $B_2$
\begin{equation}
 \begin{cases}
 (\partial_\tau+\mathcal L)v_h
  =\zeta\chi K_h+\zeta[\mathcal L,\chi]w_h
                                      +\zeta'\chi w_h,\\
 v_h(0)=0,\qquad \partial_n v_h|_{\partial B_2}=0,
 \end{cases}
 \label{loc:center-auxiliary-problem}
\end{equation}
All sources are extended by zero after localization. Since
\[
 [\mathcal L,\chi]w=-2a^{kl}\chi_k\partial_lw-a^{kl}\chi_{kl}w
                                      +C_{\cdot\cdot l}\chi_lw,
\]
Lemma~\ref{loc:local-chart-inverses}, with viscosity one, gives
\[
 \|v_h-v_k\|_{\mathcal P_{T_1}}
 \le C\bigl(\|K_h-K_k\|_{h^\alpha(B_1)}
             +\|w_h-w_k\|_{h^{1,\alpha}(B_1)}\bigr)
 \longrightarrow0.
\]
The fixed $\zeta$ controls the temporal seminorms. Taking a trace
where $\zeta=1$ proves convergence of $w_h$ in $h^{2,\alpha}(B_0)$,
and hence $\partial_ef\in h^{2,\alpha}(B_0)$.

For a multiindex $\beta$ with $|\beta|=r<m$, the differentiated
equation is
\[
 \begin{aligned}
 \mathcal L D^\beta f=D^\beta F
 +\sum_{0<\gamma\le\beta}\binom\beta\gamma\bigl\{
 &(D^\gamma a^{kl})D^{\beta-\gamma}\partial_{kl}f\\
 &-(D^\gamma C)D^{\beta-\gamma}Df
   -(D^\gamma D)D^{\beta-\gamma}f\bigr\}.
 \end{aligned}
\]
If $f\in h^{r+2,\alpha}$, this source belongs to $h^{1,\alpha}$:
one additional derivative uses at most $r+2$ derivatives of $f$ and
$r+1\le m$ derivatives of the coefficients. Applying the preceding
difference-quotient argument on nested balls gives successively
$f\in h^{r+3,\alpha}$, and hence $f\in h^{m+2,\alpha}$.
The difference quotients are full Cartesian vectors; their limits are
derivatives of the equivariant field.

Subtracting the auxiliary equations gives continuity in the stated
inputs. The scalar argument applies in the interior, and a finite
enlarged cover proves both assertions of the lemma.
\end{proof}

\subsection{Time derivatives and mixed regularity}
\label{loc:finite-time-cone}

We first recover the time derivatives in the base domain on the
interval of local existence. We then use the equation to obtain their
higher spatial regularity. The following auxiliary variable separates
these two steps. Put
\[
 c=\kappa^{-1}-\kappa\Lambda,\qquad
 d=\kappa^{-2}-\Lambda,\qquad
 Z=V+\kappa\mathcal M[R]-cR .
\]
The force is smooth from the positive subset of $\mathcal D$ to
$\mathcal Y$, and $R_t=V$ holds in $\mathcal D$. Therefore
$\partial_t\mathcal M[R]=H_RV$ in $\mathcal Y$. Using the
velocity equation and $c/\kappa=d$, we obtain
\[
 Z_t=V_t+\kappa H_RV-cV
     =-\kappa^{-1}V-\mathcal M[R]
     =-\kappa^{-1}Z-dR.
\]
Together with the definition of $Z$, this gives
\begin{equation}
 R_t+\kappa\mathcal M[R]-cR=Z,\qquad
 Z_t+\kappa^{-1}Z=-dR.
 \label{loc:radius-history}
\end{equation}
In particular $Z,Z_t\in\mathcal F_T$. The base inverse applies
also to $A_R=\kappa H_R-c$, since its scalar shift is the
finite number $\Lambda-\kappa^{-2}$.

We use the following time-integration estimate, with a constant $C_I$
independent of $0<T\le1$:
\[
 \left\|\int_0^tf(s)\dd s\right\|_{\mathcal F_T}
       \le C_I T^{1-\theta}\|f\|_{C_t\mathcal Y}.
\]
The spatial seminorm gains $T$ and the temporal seminorm gains
$T^{1-\theta}$. For $f\in C_t\mathcal Y$, approximation by finite
sums of smooth spatial sources times smooth time functions puts the
integral in the little space and preserves this bound.

\begin{lemma}
\label{loc:time-identification}
Let $0<T\le1$ and $u\in\mathcal P_T$ solve $u_t+A(t)u=G$.
Assume that the base inverse has bound $C_0$ and
\[
 A\in C_t^1\mathcal L(\mathcal D,\mathcal Y),\qquad
 A_tu,G_t\in\mathcal F_T,\qquad
 w_0=G(0)-A(0)u(0)\in\mathcal D.
\]
If
\[
 C_0C_IT^{1-\theta}\sup_t\|A_t\|_{\mathcal D\to\mathcal Y}<1,
\]
then $u_t\in\mathcal P_T$ as a $\mathcal D$-valued derivative.
Writing $w=u_t$, we have
\[
 w_t+Aw=G_t-A_tu,\qquad w(0)=w_0.
\]
\end{lemma}

\begin{proof}
Construct $w$ by Proposition~\ref{loc:global-base-inverse}, and put
$v=u(0)+\int_0^tw$, $e=v-u$, $E=v_t+Av-G$. In $\mathcal Y$,
\[
 E(0)=0,\qquad E_t=A_te,\qquad e_t+Ae=E,\qquad e(0)=0.
\]
The inverse and time-integration bounds give
\[
 \|e\|_{\mathcal P_T}
 \le C_0C_IT^{1-\theta}
       \sup_t\|A_t\|_{\mathcal D\to\mathcal Y}
       \|e\|_{\mathcal P_T}.
\]
Hence $e=0$, so $u=v$ and $u_t=w$ in $\mathcal D$.

The difference quotients satisfy
\[
 (\delta_hu)_t+A(t+h)\delta_hu
       =\delta_hG-(\delta_hA)u(t),\qquad
 \delta_hu(t)=\int_0^1w(t+sh)\dd s.
\]
Translation continuity, obtained by smooth approximation in the little
parabolic norm, gives $\delta_hu\to w$ on every shorter interval
and $\delta_hu(0)\to w_0$ in $\mathcal D$. No source trace
condition is used.
\end{proof}

For the first velocity derivative,
\[
 \begin{gathered}
 (\delta_hV)_t+\kappa(H_{R(t+h)}+\Lambda)\delta_hV
       =-\delta_h\mathcal M[R]-\kappa(\delta_hH_R)V(t),\\
 \delta_hR=\int_0^1V(t+sh)\dd s.
 \end{gathered}
\]
The force differences satisfy
\[
 \begin{aligned}
 \delta_h\mathcal M[R]
 &=\int_0^1H_{R(t)+s(R(t+h)-R(t))}\delta_hR\dd s,\\
 (\delta_hH_R)f
 &=\int_0^1D^2\mathcal M[R(t)+s(R(t+h)-R(t))]
                                      [\delta_hR,f]\dd s.
 \end{aligned}
\]
Thus the source converges to
$-H_RV-\kappa D^2\mathcal M[R][V,V]$. The initial derivative is
$r_2$ from \eqref{loc:initial-recursion}.

\emph{Higher time derivatives.}
For the induction, separate the term linear in the highest time
derivative in the chain rule:
\[
 \partial_t^j\mathcal M[R]=H_RR_j+N_j,\qquad N_1=0,
\]
where
\begin{equation}
 N_j=\sum_{\substack{\sum_{\ell=1}^{j-1}\ell m_\ell=j\\
                        \sum m_\ell\ge2}}
 \frac{j!}{\prod_{\ell=1}^{j-1}m_\ell!(\ell!)^{m_\ell}}\,
 D^{\sum m_\ell}\mathcal M[R]
              [R_1^{\,m_1},\ldots,R_{j-1}^{\,m_{j-1}}].
 \label{loc:lower-bell}
\end{equation}
The superscripts in the brackets denote repeated arguments.
Differentiating this decomposition at the preceding order gives
\[
 N_j=(N_{j-1})_t+D^2\mathcal M[R][R_1,R_{j-1}].
\]
In particular $N_2=D^2\mathcal M[R][R_1,R_1]$.
The auxiliary time derivatives are determined by the recurrence
\[
 Z_j=-\kappa^{-1}Z_{j-1}-dR_{j-1},
\]
so $Z_j$ requires no radius derivative beyond $j-1$.

Suppose the time derivatives through order $j-1$ have been
identified with values in the domain, where $2\le j\le24$.
The last one satisfies
\[
 (R_{j-1})_t+A_RR_{j-1}=G_{j-1}:=Z_{j-1}-\kappa N_{j-1}.
\]
Its differentiated source is
$Z_j-\kappa(N_{j-1})_t\in\mathcal F_T$, because every index
in the latter derivative is at most $j-1$.
Moreover
\[
 (A_R)_tf=\kappa D^2\mathcal M[R][R_1,f],\qquad
 G_{j-1}(0)-A_R(0)r_{j-1}=r_j\in\mathcal D.
\]
The initial equality follows by eliminating $Z$ in
\eqref{loc:radius-history}; its recursion is exactly
\eqref{loc:initial-recursion}. Lemma~\ref{loc:time-identification}
constructs the next actual derivative and gives
\begin{equation}
 (R_j)_t+A_RR_j=Z_j-\kappa N_j,\qquad
 R_j(0)=r_j,\qquad R_j\in\mathcal P_T.
 \label{loc:actual-time-row}
\end{equation}
The time interval can be chosen uniformly for all these stages,
since the base bound gives
$\sup\|(A_R)_t\|_{\mathcal D\to\mathcal Y}
\le\kappa C_{2,L}B$. It therefore suffices to decrease the time
in \eqref{loc:base-ball-time} once so that
\[
 C_0C_IT^{1-\theta}\kappa C_{2,L}B\le\tfrac12.
\]
This choice depends only on the base bounds. At each subsequent
stage the fixed linear inverse controls the new source on the same
interval, regardless of its norm.

The last differentiated equation gives
\begin{equation}
 R_{25}=(R_{24})_t=Z_{24}-\kappa N_{24}-A_RR_{24}
                       \in\mathcal F_T .
 \label{loc:source-row-25}
\end{equation}
Its initial value is $r_{25}$ in $\mathcal Y$.
Only source regularity is obtained for $R_{25}$, with no
first-derivative trace at the vacuum. The induction stops here,
without a derivative of order 26. The estimates
\[
 \|R_j\|_{\mathcal P_T}\le
 C_0\{\|r_j\|_{\mathcal D}+\|Z_j\|_{\mathcal F_T}
                                  +\kappa\|N_j\|_{\mathcal F_T}\}
\]
are triangular in the earlier derivatives. Finite induction therefore
bounds each derivative in terms of the prescribed initial time
derivatives and the base constants.

\emph{Higher spatial derivatives.}
Integrate the auxiliary equation. The force derivative in its initial
value is evaluated by the finite chain rule on the prescribed initial
time derivatives:
\begin{equation}
 Z_j(t)=e^{-t/\kappa}z_j
             -d\int_0^te^{-(t-s)/\kappa}R_j(s)\dd s,\qquad
 z_j=r_{j+1}+\kappa\partial_t^j\mathcal M[r_0]-cr_j
                       \in\mathcal Y^{48-2j}.
 \label{loc:high-history}
\end{equation}
The force differential costs two ordinary derivatives, and
at $j=24$ the stated initial order is only $\mathcal Y$.
We induct over the spatial and temporal orders in the sequence
\begin{equation}
 m=1,\ldots,48,\qquad
 j=0,\ldots,\left\lfloor(48-m)/2\right\rfloor .
 \label{loc:triangle-order}
\end{equation}
At each step we already have $R_j,R_{j+1}\in C_t\mathcal Y^m$.
For $m=1$, this follows from the base domains. For $m\ge2$,
use the preceding spatial level and the index identity
\[
 \mathcal D^{m-2}\subset\mathcal Y^m,\qquad
 2(j+1)+(m-2)=2j+m\le48.
\]
Equation \eqref{loc:high-history} then gives
$Z_j\in C_t\mathcal Y^m$.

At $j=0$, the source
\[
 \mathcal M[R]=\kappa^{-1}(Z-R_1+cR)
\]
belongs to $C_t\mathcal Y^m$. The nonlinear implication of
Lemma~\ref{loc:elliptic-lemma} therefore gives
$R\in C_t\mathcal D^m$.

For $j\ge1$, each time derivative in $N_j$ has smaller
index and has already been recovered at the same spatial order. Thus
\[
 H_RR_j=\kappa^{-1}(Z_j-R_{j+1}+cR_j)-N_j\in C_t\mathcal Y^m,
\]
and the linear implication gives $R_j\in C_t\mathcal D^m$.
This induction never requires the twenty-fifth time derivative
in the domain at a positive spatial order.

\begin{proof}[Proof of Proposition~\ref{loc:high-local}]
Existence and regularity were proved above.

\emph{Continuous dependence.}
The base contraction gives $R_n,V_n\to R,V$. Write
$\Delta R_j=R_{n,j}-R_j$, with the same convention for
$r_j,Z_j,N_j$. Subtracting \eqref{loc:actual-time-row} gives
\[
 \begin{aligned}
 (\Delta R_j)_t+A_{R_n}\Delta R_j
   &=\Delta Z_j-\kappa\Delta N_j-(A_{R_n}-A_R)R_j,\\
 \Delta R_j(0)&=\Delta r_j,
 \end{aligned}
\]
and therefore
\[
 \begin{aligned}
 \|\Delta R_j\|_{\mathcal P_T}\le C_0\bigl(&
 \|\Delta r_j\|_{\mathcal D}+\|\Delta Z_j\|_{\mathcal F_T}
 +\kappa\|\Delta N_j\|_{\mathcal F_T}\\
 &+\|(A_{R_n}-A_R)R_j\|_{\mathcal F_T}\bigr)\longrightarrow0.
 \end{aligned}
\]
Induction and \eqref{loc:base-lipschitz} justify the limit through
order 24; \eqref{loc:source-row-25} gives order 25 in the source
space. For spatial regularity, follow \eqref{loc:triangle-order}:
\eqref{loc:high-history} gives convergence of the auxiliary sources,
and the continuity in Lemma~\ref{loc:elliptic-lemma} gives the next
domain norm. This finite induction is uniform on the compact time
interval.

\emph{Continuation.}
Suppose $T_{\max}<\infty$ but
\[
 \sup_{t<T_{\max}}\{\|R(t)\|_{\mathcal D}+\|R_t(t)\|_{\mathcal D}\}
       \le L,\qquad\text{all positive factors}\ge\eta>0.
\]
The base contraction and time-identification estimates give a common
lifespan $T_0=T(L,\eta,\kappa,\Lambda,M,\alpha,\mathrm{atlas})>0$,
independent of higher norms. Choose
$T_{\max}-a<T_0/2$, $a<T_{\max}$, and restart with
\[
 r_j=R_j(a)\in\mathcal D^{48-2j}\quad(0\le j\le24),
 \qquad r_{25}=R_{25}(a)\in\mathcal Y.
\]
These regularities follow from \eqref{loc:actual-finite-cone}.
Differentiating \eqref{loc:actual-velocity} at $a$ gives the
recursion \eqref{loc:initial-recursion}; the domain traces impose
vacuum compatibility through order 24 only. At fixed $\kappa>0$,
these data suffice without a viscosity-dependent preparation.

The restarted solution exists on $[a,a+T_0]$, with
$a<T_{\max}<a+T_0$. Lemma~\ref{loc:time-identification} and
\eqref{loc:triangle-order} give \eqref{loc:actual-finite-cone} on this
interval. Higher initial norms enter these estimates but do not shorten
the interval. Base uniqueness identifies the two solutions on
$[a,T_{\max})$, yielding an extension past $T_{\max}$ with the stated
strong regularity. The restart uses an interior time at fixed positive
viscosity; no trace at $T_{\max}$ or weak inviscid time slice is assumed.
\end{proof}

\subsection{Approximation in the domain of the linearized operator}
\label{app:current-domain-approximation}

Use the compatible core and the closed graph of $L_R$ defined in
Subsection~\ref{loc:finite-legality-section}. Fix $q,R,f$ as in the
corresponding part of Lemma~\ref{loc:finite-strong-graph}. The
approximation constants may depend on the fixed positive scale.

\begin{proof}[Proof of Lemma~\ref{loc:finite-strong-graph}]
\emph{Finite local class.}
On the vacuum collar choose $b_l$, $1\le l\le q+1$, so that
\[
 \sum_{l=1}^{q+1}b_l(-l)^r=1\quad(0\le r\le q),\qquad
 Ef(y)=\sum_{l=1}^{q+1}b_lf(-ly)\quad(y<0\text{ near }0).
\]
The Vandermonde system is invertible. With a fixed outer cutoff,
$E$ is a bounded $h^{q,\alpha}$ extension. Let
$g_n=\varrho_{1/n}*Ef$, choose $\eta=1$ near zero, and put
\[
 c_n=g_n'(0)\longrightarrow f'(0)=0,\qquad
 f_n=g_n-c_ny\eta(y).
\]
Then $f_n'(0)=0$ and $f_n\to f$ in $h^{q,\alpha}$. The identities
\eqref{loc:quotients} give
\[
 \left\|\frac{(f_n-f)'}y\right\|_{h^{q-2,\alpha}}
 +\left\|\frac{(f_n-f)-(f_n-f)(0)}{y^2}\right\|_{h^{q-2,\alpha}}
 \le C\|f_n-f\|_{h^{q,\alpha}}\longrightarrow0.
\]
At the center, extend and mollify the full vector, then apply the
bounded equivariant projection
\[
 (\mathscr RF)(X)=\int_{SO(3)}O^{-1}F(OX)\dd O.
\]
Here $\dd O$ is normalized Haar measure. Ordinary smoothing on the interior
and a fixed mass partition, constant near the endpoints, give a single
sequence $f_n\in\mathscr C$ converging in $\mathcal D^{q-2}$.
All partition derivatives lie in regular enlarged overlaps.

At the fixed scale, the chart changes and the prescribed cutoffs give
\[
 \|v\|_{X^q}\le C_{\lambda,q}\|v\|_{\mathcal D^{q-2}},\qquad
 \mathcal S_{q-2}(g)\le C_{\lambda,q}\|g\|_{\mathcal Y^{q-2}}.
\]
The second bound includes the edge factor $\lambda$.
The formulas \eqref{loc:vacuum-hessian},
\eqref{loc:center-scalar-principal} and \eqref{loc:interior-force} imply
\[
 \|H_R(f_n-f)\|_{\mathcal Y^{q-2}}
 \le C_R\|f_n-f\|_{\mathcal D^{q-2}}\longrightarrow0,
\]
which proves \eqref{loc:finite-graph-approximation}. Simultaneous
approximation of $R,f_1,\ldots,f_p$ preserves the positive factors and,
by the differentiation formulas for the force, gives
\[
 D^p\mathcal M[R_n][f_{1,n},\ldots,f_{p,n}]
 \longrightarrow D^p\mathcal M[R][f_1,\ldots,f_p]
 \quad\hbox{in }\mathcal Y^{q-2}.
\]
The fixed-scale comparison gives convergence in the source norm.

\emph{Compatible weighted class.}
Take a defining sequence $f_n\in\mathscr C$, $f_n\to f$ in $X^q$.
On the edge, apply the product estimates with the same cutoff to
\[
 -A_n\bigl(v_{yy}+(n-1)v_y/y\bigr)-B_nv_y/y+V_nv;
\]
the quotient is controlled by Lemma~\ref{at:natural-quotient}.
On the long chart $y$ is bounded away from zero. The profile
coefficients and the compact-chart formulas therefore give
\[
 \mathcal S_{q-2}(L_Av)\le C_A\|v\|_{X^q}.
\]
Combining this with \eqref{at:current-low}--\eqref{at:current-high}
yields
\begin{equation}
 \mathcal S_{q-2}(L_Rv)\le C_{R,\lambda,q}\|v\|_{X^q}.
 \label{loc:forward-current-bound}
\end{equation}
For $q\le8$ the radius has the required $X^8$ norm. For $q\ge9$,
the term involving the high norm of the radius perturbation uses $Z\in X^q$ and
$\|v\|_{X^8}\le\|v\|_{X^q}$. Approximate $Z$ in this compatible
space. Low convergence preserves the positive factors; the coefficient
difference estimates and Lemma~\ref{at:natural-quotient} pass the
operator and its closed vacuum quotient to the limit.

Thus $L_Rf_n$ is Cauchy in $\mathcal S_{q-2}$. If its limit is $g$,
then for every $\varphi\in C_c^\infty(0,M)$,
\[
 \langle g,\varphi\rangle
 =\lim_n\langle L_Rf_n,\varphi\rangle
 =\lim_n\langle f_n,L_R^*\varphi\rangle
 =\langle L_Rf,\varphi\rangle.
\]
The compatible completions fix the endpoint representatives, and
the endpoints have zero mass measure. Hence $g=L_Rf$ in the
source space. Thus $(f_n,L_Rf_n)$ converges to $(f,L_Rf)$ in the
closed graph of $L_R$. Applying \eqref{loc:forward-current-bound}
to differences proves continuity.

On a compact positive-scale interval, the chart and cutoff bounds
can be chosen in common, with constants depending on the lower
scale. These constants establish that $f$ belongs to the strong
operator domain; they are not used in the uniform inverse estimate.
\end{proof}

\subsection{Identities at finite regularity}
\label{app:finite-regularity-identities}

We work on a compact positive-scale interval of the solution in
Proposition~\ref{loc:high-local}, with positive viscosity fixed.

\begin{proof}[Proof of Proposition~\ref{loc:finite-calculus}]
\emph{Force derivatives.}
If $u_t=v$ in $\mathcal Y$, $u(a)\in E$ and $v\in C_tE$,
where $E\hookrightarrow\mathcal Y$ continuously, then
\[
 u(t)=u(a)+\int_a^tv(\tau)\dd\tau,\qquad
 \left\|\frac{u(t+h)-u(t)}h-v(t)\right\|_E
 \le\sup_{|\tau-t|\le|h|}\|v(\tau)-v(t)\|_E\longrightarrow0.
\]
Injectivity of the embedding identifies the $E$-valued integral
with $u$. Applied successively, this gives
\[
 R_j\in C_t^1\mathcal D^{46-2j}\quad(j\le23),\qquad
 R\in C_t^{24}\mathcal D.
\]
Writing $R_j^t=\partial_t^jR$, iterate
$\partial_s(cR_j^t)=c_sR_j^t+c\lambda^{3/2}R_{j+1}^t$.
On a positive-scale interval this changes time variables without
changing the derivative counts. Below, integer subscripts denote
$s$-derivatives, and $Y,\Psi_j,F$ are as in
\eqref{rec:rows}--\eqref{rec:source}.

The coefficient and quotient bounds in \eqref{loc:vacuum-force} and
\eqref{loc:center-force} give, within the positive subset of
$\mathcal D$,
\[
 \begin{aligned}
 \mathcal M[R+h]-\mathcal M[R]-H_Rh
   &=\int_0^1(H_{R+th}-H_R)h\dd t,\\
 \|\mathcal M[R+h]-\mathcal M[R]-H_Rh\|_{\mathcal Y}
   &\le\sup_{0\le t\le1}
       \|H_{R+th}-H_R\|_{\mathcal L(\mathcal D,\mathcal Y)}
       \|h\|_{\mathcal D}=o(\|h\|_{\mathcal D}).
 \end{aligned}
\]
The higher coefficient variations are continuous on the same domain.
Thus $R\in C_s^{24}\mathcal D$ permits the partition formulas
through order 24 in $\mathcal Y$.

For a recovery pair $m=q-2$, $j+q\le23$, the force and kinetic
terms require, respectively,
\[
 R_{j+1}\in\mathcal D^{q-2},\qquad
 R_{j+2}\in\mathcal Y^{q-2}.
\]
Both follow from \eqref{loc:actual-finite-cone}, since
\[
 \begin{aligned}
 (48-2(j+1))-(q-2)&=48-2j-q\ge25-j\ge4,\\
 (50-2(j+2))-(q-2)&=48-2j-q\ge4,
 \qquad j+2\le23.
 \end{aligned}
\]
Hence
\begin{equation}
 \begin{split}
 Y_{ss}-\tfrac32bY_s
 +\lambda^3(\mathcal M[R]-\mathcal M[A])
 &+\nu\lambda^3\partial_s(\mathcal M[R]-\mathcal M[A])\\
 &+\nu\lambda^3\Lambda Y_s=\lambda^4F,
 \qquad \nu=\kappa\lambda^{-3/2}.
 \end{split}
 \label{loc:strong-difference-equation}
\end{equation}
can be differentiated $j$ times in $\mathcal Y^{q-2}$.
At energy order $j=22$ the required quantities are
$\partial_s^{23}\mathcal M[R]\in\mathcal Y$ and $R_{24}$,
which are also available. These observations justify differentiation;
the quantitative source estimates are proved separately.

\emph{Strong operator domains and transport.}
Lemma~\ref{loc:finite-strong-graph} supplies, for
$f=\Psi_j+\nu\Psi_{j+1}$, compatible core functions satisfying
\[
 \|f_n-f\|_{X^q}
       +\mathcal S_{q-2}(L_Rf_n-L_Rf)\longrightarrow0.
\]
Thus the uniform inverse estimate is applied to an element of the
strong domain. The same approximation handles differentiated nonlinear
forces. The extra regularity of $\nu\Psi_{j+1}$ comes from the
finite local class; its higher norm does not enter the recovery constant.

At fixed mass $d=M-x=\lambda^4y_n^n$, direct differentiation gives
\[
 \partial_y=n\lambda^{4/n}d^{1-1/n}\partial_d,\qquad
 \partial_s(\partial_y^kf)
   =\partial_y^kf_s+\frac{4k}{n}b\,\partial_y^kf,\qquad
 \lambda^4y_n^{n-1}\dd y_n=\frac{\dd d}{n}.
\]
Also
\[
 \partial_s\zeta(d/\lambda^4)
       =-4b(d/\lambda^4)\zeta'(d/\lambda^4).
\]
This cutoff derivative is supported where the two edge charts are
uniformly comparable, and both graph contributions are retained.
The next time derivative is continuous in the required topology since
\[
 48-2j-q\ge25-j\ge2\qquad(j+q\le23).
\]
The norm transport and integral estimates therefore hold at the stated
pairs, with their positivity, forcing and smallness hypotheses.

\emph{Potential and energy identities.}
The majorant \eqref{loc:potential-majorant} permits finite
differentiation on the positive subset of $\mathcal D$. At the vacuum,
domain variations differentiate bounded factors in
$\rho=\xi^3a$, $P'=\xi^2p$; at the center they remain $O(z)$.
Thus the differentiated Green fluxes remain $O(\xi^5)$ and
$O(z^3)$ and vanish at the endpoints. For $C^1$ domain-valued
curves, \eqref{loc:potential-force} gives
\[
 \partial_s\{D\mathcal U[R]f\}
 =D^2\mathcal U[R][R_s,f]+D\mathcal U[R]f_s
 =\mathfrak q_R[R_s,f]+(\mathcal M[R],f_s)_H.
\]
Only domain arguments are differentiated; no source-only derivative
is used in a strain.

Let $C_j^{\rm form}$ be the expression in \eqref{loc:primitive}
before its last argument $Y_j$, so
$\mathfrak A_j=\lambda^{-8}C_j^{\rm form}[Y_j]$. The rule
\begin{align*}
 \partial_s\{Q_R^{\langle q\rangle}[f,g]\}
 ={}&3bQ_R^{\langle q\rangle}[f,g]
       +Q_R^{\langle q+1\rangle}[f,g]\\
 &+Q_R^{\langle q\rangle}[f_s,g]
       +Q_R^{\langle q\rangle}[f,g_s]
\end{align*}
uses coefficient order at most 22, since $q\le21$ in
\eqref{loc:primitive}. With
\[
 \mathfrak R_j=\lambda^{-8}
           (\partial_sC_j^{\rm form}-3bC_j^{\rm form})[Y_j],
\]
we obtain
\begin{equation}
 (\mathfrak A_j)_s=-5b\mathfrak A_j+\mathfrak R_j
       +(\lambda^{-4}C_j^{\rm form})[\Psi_{j+1}],
 \qquad
 -(\lambda^{-4}C_j^{\rm form})[\Psi_{j+1}+2b\Psi_j]
       +(\mathfrak A_j)_s
       =\mathfrak R_j-7b\mathfrak A_j.
 \label{loc:finite-cancellation}
\end{equation}
Here $-5=-8+3$ comes from the exterior scale factors. At $j=22$,
the $Y_{23}$ pairing cancels before estimation, while
\[
 Y_{23}\in\mathcal D^2,\qquad
 R_{22}\in\mathcal D^4,\qquad R_{24}\in\mathcal D^0.
\]
These high derivatives enter as form arguments, not as bounded
coefficients. The first telescoping uses potential order 25 and the
further Taylor subtraction order 26; both are covered by
\eqref{loc:two-high-potential}.

At order zero, set
$\mathcal F_0=\lambda^{-1}(\mathcal M[R]-\mathcal M[A])$.
The same differentiation gives
\begin{equation}
 (\mathfrak V_0)_s=-5b\mathfrak V_0
       +(\mathcal F_0,\Psi_1)_H
       +\lambda^{-5}(\mathcal M[R]-\mathcal M[A]-H_AY_0,A_1)_H .
 \label{loc:finite-zero-potential}
\end{equation}
The last pairing is the second-order force remainder against $A_1$.
The top kinetic derivative is
\[
 \tfrac12\partial_s\|\Psi_{23}\|_H^2
       =(\Psi_{23},\Psi_{24}-4b\Psi_{23})_H.
\]
Differentiate the four terms of $H_j$ and of $\widetilde H_0$
and substitute \eqref{loc:strong-difference-equation}. The
$Q_R[\Psi_j,\Psi_{j+1}]$ terms cancel; equations
\eqref{loc:finite-cancellation} and \eqref{loc:finite-zero-potential},
together with the product rule for the scalar shift and exterior
normalization, give the modified energy identities.

Under convergence in \eqref{loc:actual-finite-cone}, the finite force
polynomials converge in the required $C_t\mathcal Y^{q-2}$ norms,
hence in the corresponding source and graph norms at fixed scale.
The multilinear potential bounds give $C_s^1$ convergence of the
corrections and relative potential; both kinetic factors converge in
$C_sH$. The energies therefore converge in $C_s^1$, without a
first-derivative trace on the source.
\end{proof}

\clearpage
\section*{Index of notation}
\phantomsection\label{notation:index}

The entries are grouped by use. The last column gives the defining
formula or statement and its page. Symbols with more than one meaning
are listed separately in their respective contexts.

\begingroup
\footnotesize
\setlength{\tabcolsep}{5pt}
\renewcommand{\arraystretch}{1.00}
\begin{longtable}{@{}>{\raggedright\arraybackslash}p{.23\textwidth}>{\raggedright\arraybackslash}p{.55\textwidth}>{\raggedright\arraybackslash}p{.15\textwidth}@{}}
\toprule
Symbol & Meaning & Definition\\
\midrule
\endfirsthead
\toprule
Symbol & Meaning & Definition\\
\midrule
\endhead
\bottomrule
\endfoot
\multicolumn{3}{@{}l}{\textit{Physical variables and scaling parameters}}\\*[.4em]
$x,M,M_\beta$
& Enclosed mass and total mass; $M=M_\beta$ is fixed along the constructed trajectory.
& \hyperref[en:force]{\ref*{en:force}}, p.~\pageref{en:force}\\[.25em]
$M_{\rm Ch}$
& Chandrasekhar mass in the normalization of the pressure law.
& \hyperref[in:threshold-inequality]{\ref*{in:threshold-inequality}}, p.~\pageref{in:threshold-inequality}\\[.25em]
$t,\tau,s$
& Physical time $t$, reversed time $\tau=-t$, and rescaled time $s$; $\partial_\tau=\lambda^{-3/2}\partial_s$.
& \hyperref[eq:scale]{\ref*{eq:scale}}, p.~\pageref{eq:scale}\\[.25em]
$\lambda,b,e$
& Radius scale, $b=\lambda_s/\lambda=\sqrt{\beta^2+e\lambda}$, and the nonnegative parameter in the scale equation.
& \hyperref[eq:b]{\ref*{eq:b}}, p.~\pageref{eq:b}\\[.25em]
$\lambda_{\rm lead}(t)$
& Leading power-law radius $(3\beta(-t)/2)^{2/3}$; $\lambda$ is the reference scale and $R_b$ is the support radius.
& \hyperref[eq:leading-power-scale]{\ref*{eq:leading-power-scale}}, p.~\pageref{eq:leading-power-scale}\\[.25em]
$T_A$
& Physical time required for the scale to move between $0$ and $A$.
& \hyperref[eq:TA]{\ref*{eq:TA}}, p.~\pageref{eq:TA}\\[.25em]
$\beta,\delta$
& Positive profile parameter and Goldreich--Weber parameter $\delta=-\beta^2/2$.
& \hyperref[mat:seed]{\ref*{mat:seed}}, p.~\pageref{mat:seed}\\[.25em]
$P,h,\rho[R]$
& Exact pressure, enthalpy, and density $(4\pi R^2R_x)^{-1}$.
& \hyperref[mat:eos]{\ref*{mat:eos}}, p.~\pageref{mat:eos}\\[.25em]
$\rho_E,u_E,\rho_L,u_L$
& Eulerian density and radial velocity, and their values along the material shells; $u_L=r_t$.
& \hyperref[in:mass-equation]{p.~\pageref*{in:mass-equation}}\\[.25em]
$F,\mathfrak e,\mathscr E,I_2$
& Internal energy density, specific internal energy, conserved total energy, and $I_2(\delta)=4\pi\int_0^1w_\delta^3z^4\dd z$.
& \hyperref[eq:physical-energy-clock]{\ref*{eq:physical-energy-clock}}, p.~\pageref{eq:physical-energy-clock}\\[.25em]
$\mathcal M[R]$
& Pressure and gravitational force in mass coordinates.
& \hyperref[en:force]{\ref*{en:force}}, p.~\pageref{en:force}\\[.25em]
$\mathcal U[R]$
& Mechanical potential; its first variation is $\mathcal M[R]$.
& \hyperref[en:potential]{\ref*{en:potential}}, p.~\pageref{en:potential}\\[.25em]
\addlinespace[.5em]
\multicolumn{3}{@{}l}{\textit{Coordinates and expansion coefficients}}\\*[.4em]
$z,w_\beta,m_\beta$
& Normalized leading radius, cube root of the rescaled density, and enclosed-mass map $x=m_\beta(z)$.
& \hyperref[mat:seed]{\ref*{mat:seed}}, p.~\pageref{mat:seed}\\[.25em]
$z_c,\xi$
& Fixed endpoint coordinates $z_c=x^{1/3}$ and $\xi=(M-x)^{1/5}$.
& \hyperref[sec:class-and-main]{\ref*{sec:class-and-main}}, p.~\pageref{sec:class-and-main}\\[.25em]
$d,D,y_n$
& $d=M-x$, $D=d/\lambda^4$, and $y_n=D^{1/n}$, $n=5,8$.
& \hyperref[at:edge-coordinates]{\ref*{at:edge-coordinates}}, p.~\pageref{at:edge-coordinates}\\[.25em]
$q,Q$ in the profile
& $q=1-z$ and $Q=D^{1/4}=y_8^2$. Elsewhere $q$ is a derivative order and $Q$ can denote the cutoff power.
& \hyperref[mat:dictionary]{p.~\pageref*{mat:dictionary}}\\[.25em]
$X,Y$ in matching
& $X=(M-x)^{1/4}$, $Y=\lambda/X$; $XY=\lambda$. These variables compare the core and layer expansions on their overlap.
& \hyperref[mat:connection-53]{\ref*{mat:connection-53}}, p.~\pageref{mat:connection-53}\\[.25em]
$L,\ell,p_{\log}$
& $L=\log\lambda$, $\ell=\log X$, and a finite logarithmic degree. The logarithms are independent during coefficient extraction.
& \hyperref[mat:dictionary]{p.~\pageref*{mat:dictionary}}\\[.25em]
$A,R,R_b$
& Matched profile, evolving material radius, and physical support radius. Local uses of $R$ for a trial profile are stated where they occur.
& \hyperref[mat:dictionary]{p.~\pageref*{mat:dictionary}}\\[.25em]
$A$ in the local vacuum equation
& Scalar principal coefficient in the fixed-scale force calculation; distinct from the matched radius.
& \hyperref[loc:vacuum-coefficients]{\ref*{loc:vacuum-coefficients}}, p.~\pageref{loc:vacuum-coefficients}\\[.25em]
$\Phi_n,Z_j$
& Coefficients of the relative core radius and the inward displacement in the layer.
& \hyperref[mat:connection-3]{\ref*{mat:connection-3}}, p.~\pageref{mat:connection-3}\\[.25em]
$Z_0,\rho_0$
& Leading inward displacement and density in the boundary layer.
& \hyperref[mat:connection-10]{\ref*{mat:connection-10}}, p.~\pageref{mat:connection-10}\\[.25em]
$\mathcal Z_\lambda,\mathcal Q_\lambda,\mathcal G_\lambda,\Theta_\lambda$
& Rescaled radius correction and its regular derivative, thickness, and density factors.
& \hyperref[eq:actual-layer-factors]{\ref*{eq:actual-layer-factors}}, p.~\pageref{eq:actual-layer-factors}\\[.25em]
$\varkappa_\delta$
& Boundary slope $-w_\delta'(1)>0$; distinct from the viscosity $\kappa$.
& \hyperref[mat:dictionary]{p.~\pageref*{mat:dictionary}}\\[.25em]
$A_0(D)$
& Layer flux coefficient $16\pi^2\rho_0^2P'(\rho_0)$; distinct from the constant upper bound $A_0$ for the scale.
& \hyperref[mat:dictionary]{p.~\pageref*{mat:dictionary}}\\[.25em]
$E_Q,E_q$
& Differential operators $E_Q=Q\partial_Q$ and $E_q=q\partial_q$.
& \hyperref[mat:finite-endpoint-class]{\ref*{mat:finite-endpoint-class}}, p.~\pageref{mat:finite-endpoint-class}\\[.25em]
$G,Q$ at the vacuum
& Regular endpoint factors. This local quotient notation is distinct from the large-layer coordinate $Q$.
& \hyperref[mat:endpoint-5]{\ref*{mat:endpoint-5}}, p.~\pageref{mat:endpoint-5}\\[.25em]
\addlinespace[.5em]
\multicolumn{3}{@{}l}{\textit{Spaces, cutoffs, and operator domains}}\\*[.4em]
$\chi_c,\chi_i,\chi_5,\chi_8$
& Fixed-mass center and interior cutoffs and the two moving edge cutoffs; their sum is one.
& \hyperref[at:power-partition]{\ref*{at:power-partition}}, p.~\pageref{at:power-partition}\\[.25em]
$U_f$
& Cartesian equivariant lift $f(|X|^3)X/|X|$ of a radial displacement.
& \hyperref[at:section]{\ref*{at:section}}, p.~\pageref{at:section}\\[.25em]
$\mathcal B_n^q,B_5^q,C_8^q$
& Weighted Sobolev spaces using ordinary derivatives and measure $y^{n-1}\dd y$; $B_5^q$ is on the bounded vacuum interval and $C_8^q$ on the longer interval.
& \hyperref[at:atlas]{\ref*{at:atlas}}, p.~\pageref{at:atlas}\\[.25em]
$X^q(s)$
& Completion of compatible smooth functions in the localized weighted norm, including derivatives of the cutoffs.
& \hyperref[at:graph-norm]{\ref*{at:graph-norm}}, p.~\pageref{at:graph-norm}\\[.25em]
$\|\cdot\|_{q,a}$
& Localized weighted derivative norm, with derivatives acting on the cutoffs; no prescribed vacuum trace.
& \hyperref[at:atlas]{\ref*{at:atlas}}, p.~\pageref{at:atlas}\\[.25em]
$\mathcal S_m$
& Strong source norm; the compact components use $g$ and the edge components use $\lambda g$.
& \hyperref[at:source-norm]{\ref*{at:source-norm}}, p.~\pageref{at:source-norm}\\[.25em]
$\mathcal D^m,\mathcal Y^m$
& Fixed-scale little-H\"older domain and source spaces. Write $\mathcal D=\mathcal D^0$, $\mathcal Y=\mathcal Y^0$.
& \hyperref[loc:domain]{\ref*{loc:domain}}, p.~\pageref{loc:domain}\\[.25em]
$H$
& Mass Hilbert space $L^2((0,M),\dd x)$.
& \hyperref[en:force]{\ref*{en:force}}, p.~\pageref{en:force}\\[.25em]
$H_R,L_R$
& $H_R=D\mathcal M[R]$, $L_R=\lambda^3H_R$. The closed graph requires convergence of both a direction and its image under the operator.
& \hyperref[loc:finite-strong-graph]{\ref*{loc:finite-strong-graph}}, p.~\pageref{loc:finite-strong-graph}\\[.25em]
$\mathcal T_A,\mathcal T_R$
& Closed first-order operators defining the pressure forms at the approximate radius and at the evolving radius.
& \hyperref[pr:pressure-operator]{\ref*{pr:pressure-operator}}, p.~\pageref{pr:pressure-operator}\\[.25em]
$Q_R,\mathcal P_A$
& Normalized Hessian $Q_R=\lambda^3D^2\mathcal U[R]$ and positive profile pressure form.
& \hyperref[en:forms]{\ref*{en:forms}}, p.~\pageref{en:forms}\\[.25em]
$V_a,V_a^*$
& Fixed form space on a positive-scale interval and its dual with pivot $H$.
& \hyperref[exlim:inputs]{\ref*{exlim:inputs}}, p.~\pageref{exlim:inputs}\\[.25em]
$F_a^q$
& Compatible Sobolev space with a fixed norm equivalent to the $X^q(s)$ norm on $I_a$.
& \hyperref[exlim:spatial-compactness]{\ref*{exlim:spatial-compactness}}, p.~\pageref{exlim:spatial-compactness}\\[.25em]
\addlinespace[.5em]
\multicolumn{3}{@{}l}{\textit{Time derivatives, recovery, and energy}}\\*[.4em]
$Y,\Psi_j$
& In recovery and energy, $Y=R-A$ and $\Psi_j=\lambda^{-4}\partial_s^jY$. This $Y$ is distinct from the matching coordinate.
& \hyperref[rec:rows]{\ref*{rec:rows}}, p.~\pageref{rec:rows}\\[.25em]
$R_j,A_j$
& Material derivatives $\partial_s^jR$ and $\partial_s^jA$ at fixed mass.
& \hyperref[exlim:rows]{\ref*{exlim:rows}}, p.~\pageref{exlim:rows}\\[.25em]
$\kappa,\nu,\Lambda$
& Viscosity, $\nu=\kappa\lambda^{-3/2}$, and a nonnegative shift constant along the trajectory.
& \hyperref[rec:rows]{\ref*{rec:rows}}, p.~\pageref{rec:rows}\\[.25em]
$Z_j$ in recovery
& Combination $Z_j=\Psi_j+\nu\Psi_{j+1}$ used in viscous recovery; distinct from the layer coefficient $Z_j$.
& \hyperref[rec:relaxation]{\ref*{rec:relaxation}}, p.~\pageref{rec:relaxation}\\[.25em]
$\mathfrak C_j,\mathfrak C_j^{\rm vis}$
& Force and viscous commutators in the differentiated equation.
& \hyperref[rec:force-commutator]{\ref*{rec:force-commutator}}, p.~\pageref{rec:force-commutator}\\[.25em]
$F,F_j$
& Forcing in the perturbation equation and its normalized time derivatives $F_j=\lambda^{-4}\partial_s^j(\lambda^4F)$.
& \hyperref[rec:source]{\ref*{rec:source}}, p.~\pageref{rec:source}\\[.25em]
$\mathcal R_{Q_*}$
& Sum of strong source norms over the mixed derivative range, $Q_*=8$ or $23$.
& \hyperref[rec:source]{\ref*{rec:source}}, p.~\pageref{rec:source}\\[.25em]
$\mathcal X_{Q_*},G_8$
& Sum of the mixed derivative norms for $q\le Q_*$, $j+q\le23$, and the square root of their squared sum for $q\le8$, respectively.
& \hyperref[rec:cones]{\ref*{rec:cones}}, p.~\pageref{rec:cones}\\[.25em]
$B,\eta$
& Norm controlling one spatial derivative and the normalized low-order norm $\eta=\lambda^{-2}G_8$.
& \hyperref[rec:base]{\ref*{rec:base}}, p.~\pageref{rec:base}\\[.25em]
$H^\sharp$
& Normalized running supremum $\sup_{s_*\le\sigma\le s}\lambda(\sigma)^{-p/2}H(\sigma)$. Here $H$ denotes a nonnegative quantity.
& \hyperref[rec:relaxation]{\ref*{rec:relaxation}}, p.~\pageref{rec:relaxation}\\[.25em]
$\varphi_\beta,f^\perp$
& Fixed normalized homologous direction and orthogonal projection in the mass inner product.
& \hyperref[en:mode]{\ref*{en:mode}}, p.~\pageref{en:mode}\\[.25em]
$V(f),Z(f)$
& Norm defined by the pressure form and the mass norm, and the norm with weight $b^2$ on the homologous component.
& \hyperref[en:norms]{\ref*{en:norms}}, p.~\pageref{en:norms}\\[.25em]
$a_f,d_f,G_f$
& Strains $f/R$, $f_x/R_x$, and $2a_f+d_f$.
& \hyperref[en:strain-control]{\ref*{en:strain-control}}, p.~\pageref{en:strain-control}\\[.25em]
$T_\beta,H_*$
& Pullback by the leading profile and common weighted Hilbert space used to compare the pressure forms.
& \hyperref[en:common-seed-pivot]{\ref*{en:common-seed-pivot}}, p.~\pageref{en:common-seed-pivot}\\[.25em]
$\epsilon_A$
& Error for the chosen profile $C_{J,K,\beta}\lambda^{1/2}(1+|\log\lambda|)^{p_{\log}}$.
& \hyperref[en:profile-error]{\ref*{en:profile-error}}, p.~\pageref{en:profile-error}\\[.25em]
$Q_A^{\langle k\rangle}$
& Rescaled time derivative of the Hessian; the exterior factor $\lambda^3$ is held fixed.
& \hyperref[en:covariant-form-definition]{\ref*{en:covariant-form-definition}}, p.~\pageref{en:covariant-form-definition}\\[.25em]
$C_j$
& Remainder after subtracting the linearized force from the differentiated force difference.
& \hyperref[en:commutator-definition]{\ref*{en:commutator-definition}}, p.~\pageref{en:commutator-definition}\\[.25em]
$\omega_j,E$
& Fixed positive weights and the positive quadratic energy.
& \hyperref[en:positive-energy]{\ref*{en:positive-energy}}, p.~\pageref{en:positive-energy}\\[.25em]
$\widetilde E$
& Modified nonlinear energy, including the relative potential and the commutator corrections.
& \hyperref[en:modified-energy]{\ref*{en:modified-energy}}, p.~\pageref{en:modified-energy}\\[.25em]
$\mathfrak A_j$
& Commutator correction added to the energy at positive derivative orders.
& \hyperref[en:primitive-definition]{\ref*{en:primitive-definition}}, p.~\pageref{en:primitive-definition}\\[.25em]
$\mathcal D$ in the energy estimate
& Sum of the shifted forms of the kinetic terms; distinct from the local little-H\"older domain.
& \hyperref[ct:dissipation]{\ref*{ct:dissipation}}, p.~\pageref{ct:dissipation}\\[.25em]
\addlinespace[.5em]
\multicolumn{3}{@{}l}{\textit{Parameter choices and the limits}}\\*[.4em]
$J,p,J^\sharp,K$
& $p=2J+1$, $J^\sharp=J+1713$, and $K=2J+3512$; these integers are fixed before the leading profile is selected.
& \hyperref[ct:orders]{\ref*{ct:orders}}, p.~\pageref{ct:orders}\\[.25em]
$A_*,A_0$
& Upper bounds for the scale in the profile construction and in the evolution; both are chosen after the leading profile and finite orders.
& \hyperref[ct:parameter-admissibility]{\ref*{ct:parameter-admissibility}}, p.~\pageref{ct:parameter-admissibility}\\[.25em]
$\lambda_*,a,I_a$
& Initial scale, positive terminal scale, and $I_a=\{s:a\le\lambda(s)\le A_0\}$.
& \hyperref[exlim:inputs]{\ref*{exlim:inputs}}, p.~\pageref{exlim:inputs}\\[.25em]
$d_0^{\rm edge},d_0$
& Width of the fixed boundary neighborhood in mass coordinates and, in the energy argument, the perturbation amplitude.
& \hyperref[ct:small-constants]{\ref*{ct:small-constants}}, p.~\pageref{ct:small-constants}\\[.25em]
$C_0,C_\beta,C_a$
& Constants uniform in the leading family, for a chosen profile, and at fixed terminal scale, respectively. Dependencies are stated in each estimate.
& \hyperref[exlim:inputs]{\ref*{exlim:inputs}}, p.~\pageref{exlim:inputs}\\[.25em]
\addlinespace[.5em]
\end{longtable}
\endgroup

\end{document}